\documentclass[12pt,a4paper]{book}

\usepackage{latexsym}
\usepackage{tikz} 
\usepackage{multirow}
\usepackage{amsfonts,amsmath,amsthm}
\usepackage{graphicx}
\usepackage{amssymb}
\usepackage{mathrsfs}
\usepackage[top=1.5cm, bottom=1.8cm, left=1.5cm, right=1.5cm]{geometry}
\usepackage{blkarray}
\usepackage[all]{xy}
\usepackage{color}

\usepackage{pdfpages}

\usepackage{fancyhdr}

\fancypagestyle{abstract}{
\fancyhf{}
\fancyhead[R]{{\em R\'esum\'e - Abstract}}
\fancyfoot[C]{\thepage}
}

\fancypagestyle{remerciements}{
\fancyhf{}
\fancyhead[R]{{\em Remerciements}}
\fancyfoot[C]{\thepage}
}

\newtheorem{theorem}{Theorem}[section]
\newtheorem{proposition}[theorem]{Proposition}
\newtheorem{corollary}[theorem]{Corollary}
\newtheorem{definition}[theorem]{Definition}
\newtheorem{lemma}[theorem]{Lemma}
\newtheorem{conjecture}[theorem]{Conjecture}

\theoremstyle{remark}
\newtheorem{remark}[theorem]{Remark}
\newtheorem{exemple}[theorem]{Example}

\newcommand{\fleche}[4]{                     
            \begin{array}{rcl} #1 & \rightarrow & #2 \\   %
                         #3 &\mapsto & #4          %
            \end{array}}

\newcommand{\fonc}[5]{                     
            \begin{array}{crll}#1 :& #2 & \rightarrow & #3 \\   %
                         &#4 &\mapsto & #5          %
            \end{array}}

\newcommand{\exposantGauche}[2]{{\vphantom{#2}}^{#1}#2}

\newcommand\xqed[1]{%
  \leavevmode\unskip\penalty9999 \hbox{}\nobreak\hfill
  \quad\hbox{#1}}
\newcommand{\finEx}{\xqed{$\triangle$}}

\DeclareMathOperator{\End}{End}
\DeclareMathOperator{\SLF}{SLF}
\DeclareMathOperator{\Hom}{Hom}
\DeclareMathOperator{\vect}{vect}
\DeclareMathOperator{\Ob}{Ob}
\DeclareMathOperator{\Proj}{Proj}
\DeclareMathOperator{\Top}{Top}
\DeclareMathOperator{\Soc}{Soc}

\title{Quantification combinatoire, mapping class groups et alg\`ebres skein}
\author{Matthieu Faitg}

\begin{document}

\includepdf{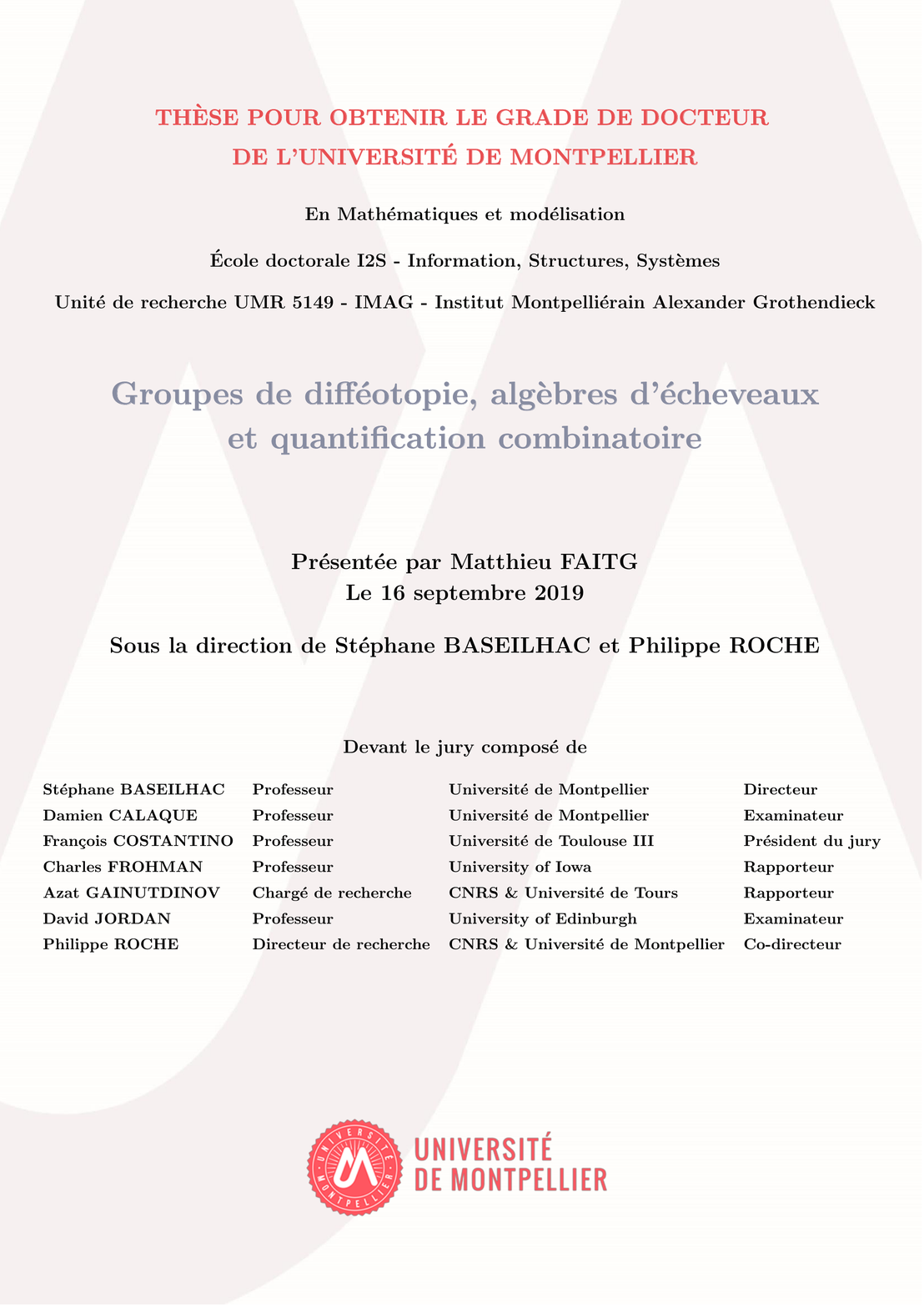}

\vspace*{\stretch{1}}

\thispagestyle{empty}
\begin{center}
{\huge \textbf{Mapping class groups, skein algebras \\and combinatorial quantization} }

%
%
%
\end{center}

\vspace*{\stretch{1}}

\newpage
\thispagestyle{remerciements}
\vspace*{10em}
{\huge \noindent  \textbf{Remerciements}}

\vspace{6em}

\indent Faire une th\`ese n'est pas une mince affaire et y parvenir n\'ecessite le soutien d'un certain nombre de personnes, qu'il m'est agr\'eable de remercier ici.

\medskip

\indent Je souhaite en premier lieu remercier sinc\`erement mes directeurs de th\`ese, St\'ephane Baseilhac et Philippe Roche. Merci de m'avoir accord\'e votre confiance en me proposant ce sujet et merci pour votre disponibilit\'e, votre comp\'etence, vos conseils avis\'es, vos encouragements et votre grande sympathie. De fa\c{c}on g\'en\'erale, merci pour votre accompagnement de qualit\'e pendant ces trois ann\'ees.

\medskip

\indent Merci \`a Charles Frohman et Azat Gainutdinov d'avoir accet\'e d'\^etre les rapporteurs de cette th\`ese, ainsi qu'\`a Damien Calaque, Fran\c{c}ois Costantino et David Jordan d'avoir accept\'e d'en \^etre les examinateurs. Un merci de plus \`a Azat pour ses nombreux commentaires et suggestions tr\`es utiles concernant mes pr\'epublications et ce manuscrit de th\`ese.

\medskip

\indent Un grand merci collectif \`a tous les doctorants de l'IMAG (je ne prendrai pas le risque d'essayer de tous les citer) pour l'ambiance tr\`es amicale et la bonne humeur que vous avez fait r\'egner; vous avez \'et\'e un soutien tr\`es important. Plus g\'en\'eralement je remercie tous les membres de l'IMAG pour leur bienveillance ainsi que tous les personnels administratifs du b\^atiment 9 pour leur efficacit\'e.

\medskip

\indent Je remercie bien \'evidemment ma famille pour leur aide et leurs encouragements pendant toutes ces ann\'ees, malgr\'e la nature hautement myst\'erieuse de cette activit\'e pour un regard ext\'erieur. Enfin, j'ai une pens\'ee particuli\`ere pour mes parents, qui nous ont quitt\'e lors de ces ann\'ees de th\`ese; c'est avec beaucoup d'\'emotion que je leur rend hommage et que je leur t\'emoigne ma reconnaissance pour tout le soutien qu'ils m'ont apport\'e.

\newpage

\vspace*{\stretch{1}}

\thispagestyle{abstract}
\noindent \textbf{R\'esum\'e}. \hspace{2pt} Les alg\`ebres $\mathcal{L}_{g,n}(H)$ ont \'et\'e introduites par Alekseev--Grosse--Schomerus et Buffenoir--Roche au milieu des ann\'ees 1990, dans le cadre de la quantification combinatoire de l'espace de modules des $G$-connexions plates sur la surface $\Sigma_{g,n}$ de genre $g$ avec $n$ disques ouverts enlev\'es. L'alg\`ebre de Hopf $H$, appel\'ee alg\`ebre de jauge, \'etait \`a l'origine le groupe quantique $U_q(\mathfrak{g})$, avec $\mathfrak{g} = \mathrm{Lie}(G)$. Dans cette th\`ese nous appliquons les alg\`ebres $\mathcal{L}_{g,n}(H)$ \`a la topologie en basses dimensions (groupe de diff\'eotopie et alg\`ebres d'\'echeveaux des surfaces), sous l'hypoth\`ese que $H$ est une alg\`ebre de Hopf de dimension finie, factorisable et enrubann\'ee mais pas n\'ecessairement semi-simple, l'exemple phare d'une telle alg\`ebre de Hopf \'etant le groupe quantique restreint $\bar U_q(\mathfrak{sl}_2)$ (o\`u $q$ est une racine $2p$-i\`eme de l'unit\'e). 

\smallskip

\indent D'abord, nous construisons en utilisant $\mathcal{L}_{g,n}(H)$ une repr\'esentation projective des groupes de diff\'eotopie de $\Sigma_{g,0} \!\setminus\! D$ et de $\Sigma_{g,0}$ (o\`u $D$ est un disque ouvert). Nous donnons des formules pour les repr\'esentations d'un ensemble de twists de Dehn qui engendre le groupe de diff\'eotopie; en particulier ces formules nous permettent de montrer que notre repr\'esentation est \'equivalente \`a celle construite par Lyubashenko--Majid et Lyubashenko \textit{via} des m\'ethodes cat\'egoriques. Pour le tore $\Sigma_{1,0}$ avec l'alg\`ebre de jauge $\bar U_q(\mathfrak{sl}_2)$, nous calculons explicitement la repr\'esentation de $\mathrm{SL}_2(\mathbb{Z})$ en utilisant une base convenable de l'espace de repr\'esentation et nous en d\'eterminons la structure.

\smallskip

\indent Ensuite, nous introduisons une description diagrammatique de $\mathcal{L}_{g,n}(H)$ qui nous permet de d\'efinir de fa\c{c}on tr\`es naturelle l'application boucle de Wilson $W$. Cette application associe un \'el\'ement de $\mathcal{L}_{g,n}(H)$ \`a chaque entrelac dans $(\Sigma_{g,n} \!\setminus\! D) \times [0,1]$ qui est parall\'elis\'e, orient\'e et colori\'e par des $H$-modules. Quand l'alg\`ebre de jauge est $H = \bar U_q(\mathfrak{sl}_2)$, nous utilisons $W$ et les repr\'esentations de $\mathcal{L}_{g,n}(H)$ pour construire des repr\'esentations des alg\`ebres d'\'echeveaux $\mathcal{S}_q(\Sigma_{g,n})$. Pour le tore $\Sigma_{1,0}$ nous \'etudions explicitement cette repr\'esentation.

\vspace{5em}

\noindent \textbf{Abstract}. \hspace{2pt} The algebras $\mathcal{L}_{g,n}(H)$ have been introduced by Alekseev--Grosse--Schomerus and Buffenoir--Roche in the middle of the 1990's, in the program of combinatorial quantization of the moduli space of flat $G$-connections over the surface $\Sigma_{g,n}$ of genus $g$ with $n$ open disks removed. The Hopf algebra $H$, called gauge algebra, was originally the quantum group $U_q(\mathfrak{g})$, with $\mathfrak{g} = \mathrm{Lie}(G)$. In this thesis we apply these algebras $\mathcal{L}_{g,n}(H)$ to low-dimensional topology (mapping class groups and skein algebras of surfaces), under the assumption that $H$ is a finite dimensional factorizable ribbon Hopf algebra which is not necessarily semisimple, the guiding example of such a Hopf algebra being the restricted quantum group $\bar U_q(\mathfrak{sl}_2)$ (where $q$ is a $2p$-th root of unity). 

\smallskip

\indent First, we construct from $\mathcal{L}_{g,n}(H)$ a projective representation of the mapping class groups of $\Sigma_{g,0} \!\setminus\! D$ and of $\Sigma_{g,0}$ ($D$ being an open disk). We provide formulas for the representations of Dehn twists generating the mapping class group; in particular these formulas allow us to show that our representation is equivalent to the one constructed by Lyubashenko--Majid and Lyubashenko \textit{via} categorical methods. For the torus $\Sigma_{1,0}$ with the gauge algebra $\bar U_q(\mathfrak{sl}_2)$, we compute explicitly the representation of $\mathrm{SL}_2(\mathbb{Z})$ using a suitable basis of the representation space and we determine the structure of this representation.

\smallskip

\indent Second, we introduce a diagrammatic description of $\mathcal{L}_{g,n}(H)$ which enables us to define in a very natural way the Wilson loop map $W$. This map associates an element of $\mathcal{L}_{g,n}(H)$ to any link in $(\Sigma_{g,n} \!\setminus\! D) \times [0,1]$ which is framed, oriented and colored by $H$-modules. When the gauge algebra is $H = \bar U_q(\mathfrak{sl}_2)$, we use $W$ and the representations of $\mathcal{L}_{g,n}(H)$ to construct representations of the skein algebras $\mathcal{S}_q(\Sigma_{g,n})$. For the torus $\Sigma_{1,0}$ we explicitly study this representation.

\vspace*{\stretch{1}}

\tableofcontents

\chapter{Introduction}

\section{Introduction en fran\c{c}ais}
\indent Soit $\Sigma_{g,n}$ une surface compacte orient\'ee de genre $g$ avec $n$ disques ouverts enlev\'es. ``L'alg\`ebre de graphe'' $\mathcal{L}_{g,n}$ a \'et\'e introduite et \'etudi\'ee par Alekseev \cite{alekseev}, Alekseev--Grosse--Schomerus \cite{AGS, AGS2} et Buffenoir--Roche \cite{BR, BR2} au milieu des ann\'ees 1990, dans le cadre de la quantification combinatoire de l'espace de modules des connexions plates sur $\Sigma_{g,n}$. C'est une alg\`ebre associative (non commutative) d\'efinie par g\'en\'erateurs et relations, les relations \'etant donn\'ees sous une forme matricielle. Le th\`eme principal de cette th\`ese est d'appliquer ces alg\`ebres \`a la construction de repr\'esentations quantiques des groupes de diff\'eotopie et des alg\`ebres d'\'echeveaux des surfaces aux racines de l'unit\'e.

\smallskip

\indent Dans la section \ref{introCombQuantFr} ci-dessous, nous expliquons le contexte et les id\'ees de la quantification combinatoire et la d\'efinition de l'alg\`ebre $\mathcal{L}_{g,n}$. Puis de la section \ref{introPropLgnFr} \`a la section \ref{introRepSkeinFr} nous \'enon\c{c}ons et expliquons nos principaux r\'esultats. Enfin, la section \ref{introPerspectivesFr} contient des conjectures et probl\`emes qui peuvent \^etre le point de d\'epart d'autres travaux.

\subsection{Quantification combinatoire}\label{introCombQuantFr}
%

\indent Nous rappelons rapidement les principaux ingr\'edients de la quantification combinatoire. Soit $G$ un groupe de Lie alg\'ebrique (g\'en\'eralement suppos\'e connexe et simplement connexe, par exemple $G = \mathrm{SL}_2(\mathbb{C})$) et $\Sigma_{g,n}$ une surface compacte orient\'ee de genre $g$ avec $n$ disques ouverts enlev\'es. On consid\`ere l'espace de modules des $G$-connexions plates $\mathcal{M}_{g,n} = \mathcal{A}_f/\mathcal{G}$, o\`u $\mathcal{A} = \Omega^1(\Sigma_{g,n}, \mathfrak{g})$ est identifi\'e avec l'espace de toutes les $G$-connexions, $\mathcal{A}_f$ est le sous-espace des connexions plates, et $\mathcal{G} = C^{\infty}(\Sigma_{g,n},G)$ est le groupe de jauge. Ces objets peuvent \^etre d\'ecrits de fa\c{c}on discr\`ete et combinatoire, en utilisant les holonomies le long des ar\^etes d'un graphe remplissant. Il s'agit d'un graphe orient\'e plong\'e sur $\Sigma_{g,n}$ (ses sommets $v \in V$ sont des points de $\Sigma_{g,n}$ et ses ar\^etes $e \in E$ sont des courbes simples orient\'ees sur $\Sigma_{g,n}$ qui relient deux sommets et qui ne se croisent pas entre elles) tel que $\Sigma_{g,n} \!\setminus\! \Gamma$ est une r\'eunion de disques ouverts. Soit $\mathcal{A}_d = G^E$. Un \'el\'ement de $\mathcal{A}_d$ est appel\'e une connexion discr\`ete; il doit \^etre pens\'e comme la collection $(h_e)_{e \in E}$ des holonomies d'une connexion le long des ar\^etes de $\Gamma$. Si $\gamma = (e_1, \ldots, e_k)$ est un chemin dans $\Gamma$, on d\'efinit l'holonomie discr\`ete d'une connexion discr\`ete $(h_e)_{e \in E}$ le long $\gamma$ comme \'etant le produit $h_{e_1} \ldots h_{e_k}$. Une connexion discr\`ete est dite plate si son holonomie le long de toute face du graphe vaut $1$. Ceci donne l'ensemble $\mathcal{A}_{df} \subset \mathcal{A}_d$ des connexions discr\`etes plates. Enfin, le groupe de jauge $\mathcal{G}$ agit par conjugaison sur l'holonomie le long d'une courbe d'une connexion dans $\mathcal{A}$. Ainsi, nous d\'efinissons le groupe de jauge discret comme \'etant $\mathcal{G}_d = G^V$ et son action sur les connexions discr\`etes est $(h_v)_{v \in V} \cdot (h_e)_{e \in E} = (h_{e^-} h_e h_{e^+}^{-1})$, o\`u $e^-$ est le point de d\'epart de $e$ et $e^+$ est son point d'arriv\'ee. Un r\'esultat connu affirme que $\mathcal{A}_{df}/\mathcal{G}_d \cong \Hom\bigl( \pi_1(\Sigma_{g,n}), G \bigr)/G$ (en principe le quotient est \`a consid\'erer dans le cadre de la th\'eorie g\'eom\'etrique des invariants (quotient GIT), mais ici la discussion est informelle). Donc cette construction est \'equivalente \`a la vari\'et\'e des caract\`eres, qui est un mod\`ele pour $\mathcal{M}_{g,n}$. Pour plus d'informations sur l'espace de modules et sa description combinatoire, une r\'ef\'erence accessible est \cite{labourie}. Cette description est aussi appel\'ee une th\'eorie de jauge discr\`ete, \textit{cf.} \cite{BFKB}.

\smallskip 

\indent L'espace de modules $\mathcal{A}_f/\mathcal{G}$ est muni de la structure de Poisson d'Atiyah--Bott--Goldman \cite{AB, goldman}, c'est-\`a-dire qu'on a un crochet de Poisson sur l'alg\`ebre des fonctions $\mathbb{C}[\mathcal{A}_f/\mathcal{G}] = \mathbb{C}[\mathcal{A}_f]^{\mathcal{G}}$. La structure de Poisson correspondante sur la discr\'etisation $\mathcal{A}_{df}/\mathcal{G}_d$ a \'et\'e d\'ecrite par Fock--Rosly \cite{FockRosly0}; c'est un crochet de Poisson d\'efinit de fa\c{c}on matricielle sur l'alg\`ebre des fonctions $\mathbb{C}[\mathcal{A}_d]$ et qui induit un crochet de Poisson sur $\mathbb{C}[\mathcal{A}_{df}]^{\mathcal{G}_d}$ (o\`u le groupe de jauge agit \`a droite sur les fonctions de fa\c{c}on \'evidente). L'alg\`ebre $\mathcal{L}_{g,n}$ est une quantification de $\mathbb{C}[\mathcal{A}_d]$. Nous n'avons pas besoin de d\'etailler plus ce fait puisque nous ne l'utilisons pas dans cette th\`ese. Dans la suite, nous expliquons simplement l'analogie entre $\mathcal{L}_{g,n}$ et $\mathbb{C}[\mathcal{A}_d]$.

\smallskip

\indent Ici nous utiliserons toujours le graphe $\Gamma = \Gamma_{g,n} \subset \Sigma_{g,n}$ qui a un seul sommet et dont les ar\^etes forment un syst\`eme de g\'en\'erateurs du groupe fondamental :
\[ \Gamma_{g,n} = \bigl( \{\bullet\}, \{b_1, a_1, \ldots, b_g, a_g, m_{g+1}, \ldots, m_{g+n} \} \bigr). \]
Il est repr\'esent\'e ci-dessous :
\begin{center}
\begingroup%
  \makeatletter%
  \providecommand\color[2][]{%
    \errmessage{(Inkscape) Color is used for the text in Inkscape, but the package 'color.sty' is not loaded}%
    \renewcommand\color[2][]{}%
  }%
  \providecommand\transparent[1]{%
    \errmessage{(Inkscape) Transparency is used (non-zero) for the text in Inkscape, but the package 'transparent.sty' is not loaded}%
    \renewcommand\transparent[1]{}%
  }%
  \providecommand\rotatebox[2]{#2}%
  \newcommand*\fsize{\dimexpr\f@size pt\relax}%
  \newcommand*\lineheight[1]{\fontsize{\fsize}{#1\fsize}\selectfont}%
  \ifx\svgwidth\undefined%
    \setlength{\unitlength}{521.43291761bp}%
    \ifx\svgscale\undefined%
      \relax%
    \else%
      \setlength{\unitlength}{\unitlength * \real{\svgscale}}%
    \fi%
  \else%
    \setlength{\unitlength}{\svgwidth}%
  \fi%
  \global\let\svgwidth\undefined%
  \global\let\svgscale\undefined%
  \makeatother%
  \begin{picture}(1,0.16094534)%
    \lineheight{1}%
    \setlength\tabcolsep{0pt}%
    \put(0,0){\includegraphics[width=\unitlength,page=1]{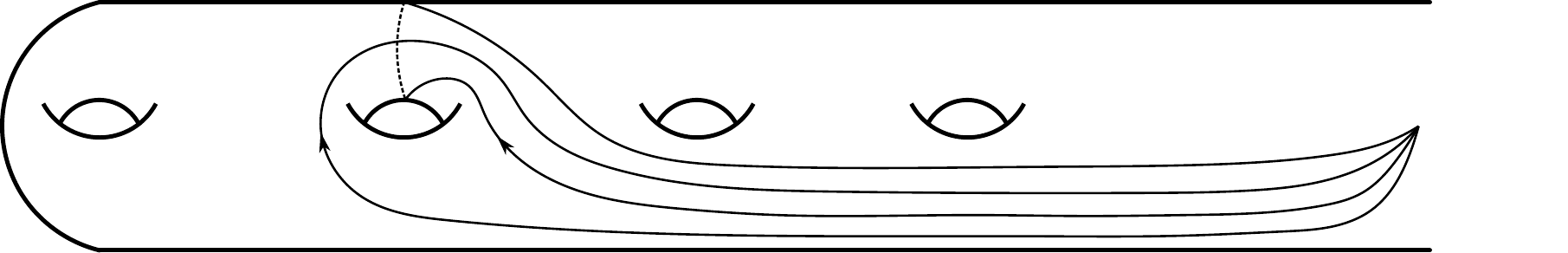}}%
    \put(0.20428392,0.02383234){\color[rgb]{0,0,0}\makebox(0,0)[lt]{\lineheight{1.25}\smash{\begin{tabular}[t]{l}$b_i$\end{tabular}}}}%
    \put(0.32337079,0.03838344){\color[rgb]{0,0,0}\makebox(0,0)[lt]{\lineheight{1.25}\smash{\begin{tabular}[t]{l}$a_i$\end{tabular}}}}%
    \put(0.05743641,0.10283815){\color[rgb]{0,0,0}\makebox(0,0)[lt]{\lineheight{1.25}\smash{\begin{tabular}[t]{l}1\end{tabular}}}}%
    \put(0.2515373,0.05084048){\color[rgb]{0,0,0}\makebox(0,0)[lt]{\lineheight{1.25}\smash{\begin{tabular}[t]{l}$i$\end{tabular}}}}%
    \put(0.42391052,0.10218237){\color[rgb]{0,0,0}\makebox(0,0)[lt]{\lineheight{1.25}\smash{\begin{tabular}[t]{l}$i+1$\end{tabular}}}}%
    \put(0.61165366,0.1065983){\color[rgb]{0,0,0}\makebox(0,0)[lt]{\lineheight{1.25}\smash{\begin{tabular}[t]{l}$g$\end{tabular}}}}%
    \put(0,0){\includegraphics[width=\unitlength,page=2]{surfaceAvecGrapheIntro.pdf}}%
    \put(0.6736706,0.12454969){\color[rgb]{0,0,0}\makebox(0,0)[lt]{\lineheight{1.25}\smash{\begin{tabular}[t]{l}$1$\end{tabular}}}}%
    \put(0.87077249,0.12354791){\color[rgb]{0,0,0}\makebox(0,0)[lt]{\lineheight{1.25}\smash{\begin{tabular}[t]{l}$n$\end{tabular}}}}%
    \put(0.78437644,0.12378273){\color[rgb]{0,0,0}\makebox(0,0)[lt]{\lineheight{1.25}\smash{\begin{tabular}[t]{l}$j+1$\end{tabular}}}}%
    \put(0.69691183,0.08056347){\color[rgb]{0,0,0}\makebox(0,0)[lt]{\lineheight{1.25}\smash{\begin{tabular}[t]{l}$m_{g+j}$\end{tabular}}}}%
    \put(0.74884075,0.14005954){\color[rgb]{0,0,0}\makebox(0,0)[lt]{\lineheight{1.25}\smash{\begin{tabular}[t]{l}$j$\end{tabular}}}}%
    \put(0,0){\includegraphics[width=\unitlength,page=3]{surfaceAvecGrapheIntro.pdf}}%
    \put(0.94474914,0.07473774){\color[rgb]{0,0,0}\makebox(0,0)[lt]{\lineheight{1.25}\smash{\begin{tabular}[t]{l}$D$\end{tabular}}}}%
    \put(0.91106599,0.03741553){\color[rgb]{0,0,0}\makebox(0,0)[lt]{\lineheight{1.25}\smash{\begin{tabular}[t]{l}$c_{g,n}$\end{tabular}}}}%
  \end{picture}%
\endgroup%

\end{center}
On a $\Sigma_{g,n} \!\setminus\! \Gamma_{g,n} \cong D$, o\`u $D$ est un disque ouvert. Ainsi, le voisinage tubulaire ferm\'e de $\Gamma_{g,n}$ est hom\'eomorphe \`a $\Sigma_{g,n} \!\setminus\! D$ :
\begin{center}
\begingroup%
  \makeatletter%
  \providecommand\color[2][]{%
    \errmessage{(Inkscape) Color is used for the text in Inkscape, but the package 'color.sty' is not loaded}%
    \renewcommand\color[2][]{}%
  }%
  \providecommand\transparent[1]{%
    \errmessage{(Inkscape) Transparency is used (non-zero) for the text in Inkscape, but the package 'transparent.sty' is not loaded}%
    \renewcommand\transparent[1]{}%
  }%
  \providecommand\rotatebox[2]{#2}%
  \newcommand*\fsize{\dimexpr\f@size pt\relax}%
  \newcommand*\lineheight[1]{\fontsize{\fsize}{#1\fsize}\selectfont}%
  \ifx\svgwidth\undefined%
    \setlength{\unitlength}{554.41358436bp}%
    \ifx\svgscale\undefined%
      \relax%
    \else%
      \setlength{\unitlength}{\unitlength * \real{\svgscale}}%
    \fi%
  \else%
    \setlength{\unitlength}{\svgwidth}%
  \fi%
  \global\let\svgwidth\undefined%
  \global\let\svgscale\undefined%
  \makeatother%
  \begin{picture}(1,0.20278397)%
    \lineheight{1}%
    \setlength\tabcolsep{0pt}%
    \put(0,0){\includegraphics[width=\unitlength,page=1]{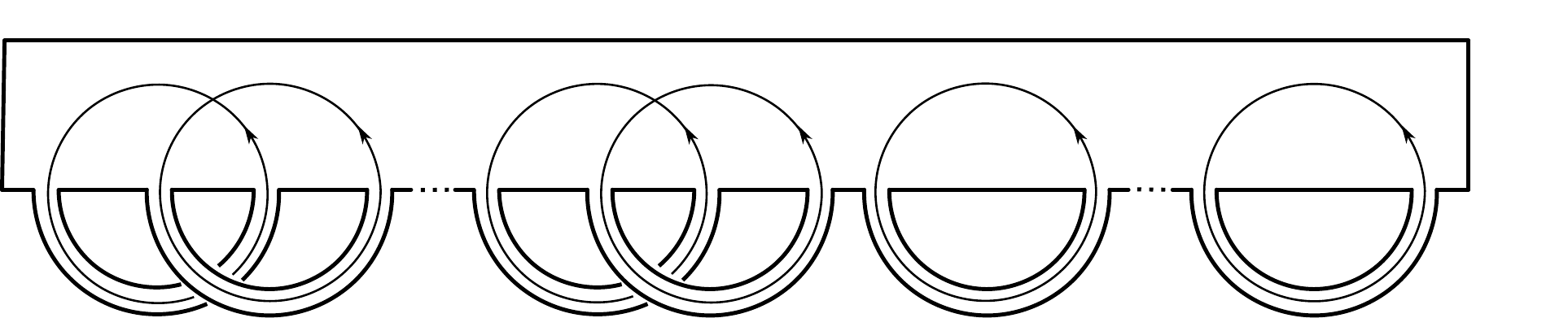}}%
    \put(0.02975383,0.14410272){\color[rgb]{0,0,0}\makebox(0,0)[lt]{\lineheight{1.25}\smash{\begin{tabular}[t]{l}$[b_1]$\end{tabular}}}}%
    \put(0.21139,0.1442029){\color[rgb]{0,0,0}\makebox(0,0)[lt]{\lineheight{1.25}\smash{\begin{tabular}[t]{l}$[a_1]$\end{tabular}}}}%
    \put(0.30990471,0.14420288){\color[rgb]{0,0,0}\makebox(0,0)[lt]{\lineheight{1.25}\smash{\begin{tabular}[t]{l}$[b_g]$\end{tabular}}}}%
    \put(0.49286062,0.1442029){\color[rgb]{0,0,0}\makebox(0,0)[lt]{\lineheight{1.25}\smash{\begin{tabular}[t]{l}$[a_g]$\end{tabular}}}}%
    \put(0.66898082,0.14420289){\color[rgb]{0,0,0}\makebox(0,0)[lt]{\lineheight{1.25}\smash{\begin{tabular}[t]{l}$[m_{g+1}]$\end{tabular}}}}%
    \put(0.87686702,0.14420289){\color[rgb]{0,0,0}\makebox(0,0)[lt]{\lineheight{1.25}\smash{\begin{tabular}[t]{l}$[m_{g+n}]$\end{tabular}}}}%
    \put(0,0){\includegraphics[width=\unitlength,page=2]{grapheEpaissiIntro.pdf}}%
    \put(0.38769684,0.18908002){\color[rgb]{0,0,0}\makebox(0,0)[lt]{\lineheight{1.25}\smash{\begin{tabular}[t]{l}$c_{g,n}$\end{tabular}}}}%
  \end{picture}%
\endgroup%

\end{center}
o\`u $[x]$ d\'enote la classe d'homotopie libre de $x \in \pi_1(\Sigma_{g,n} \!\setminus\! D)$. L'unique face du graphe $\Gamma_{g,n}$ est la courbe induite par la suppression de $D$ :
\[ c_{g,n} = b_1 a_1^{-1} b_1^{-1} a_1 \ldots b_g a_g^{-1} b_g^{-1} a_g m_{g+1} \ldots m_{g+n}. \]
Avec ce choix de graphe, une connexion discr\`ete $A_d \in \mathcal{A}_d$ associe un \'el\'ement de $G$ \`a chaque g\'en\'erateur de $\pi_1(\Sigma_{g,n} \!\setminus\! D)$, et peut donc \^etre identifi\'ee avec une liste d'\'el\'ements de $G$ :
\[ A_d = \bigl( h_{b_1}, h_{a_1}, \ldots, h_{b_g}, h_{a_g}, h_{m_{g+1}}, \ldots, h_{m_{g+n}} \bigr) \in G^{2g+n}. \]
Une connexion discr\`ete plate $A_d \in \mathcal{A}_{df}$ associe un \'el\'ement de $G$ \`a chaque g\'en\'erateur de $\pi_1(\Sigma_{g,n}) = \pi_1(\Sigma_{g,n} \!\setminus\! D)/\langle c_{g,n} \rangle$. Il s'agit d'une liste $\bigl( h_{b_1}, h_{a_1}, \ldots, h_{b_g}, h_{a_g}, h_{m_{g+1}}, \ldots, h_{m_{g+n}} \bigr)$ d'\'el\'ements de $G$ qui v\'erifie que
\begin{equation}\label{flatnessConstraintFr}
\mathrm{Hol}(A_d, c_{g,n}) = h_{b_1} h_{a_1}^{-1} h_{b_1}^{-1} h_{a_1} \ldots h_{b_g} h_{a_g}^{-1} h_{b_g}^{-1} h_{a_g} h_{m_{g+1}} \ldots h_{m_{g+n}} = 1. 
\end{equation}
Le groupe de jauge discret est simplement $\mathcal{G}_d = G$ (puisque $V = \{\bullet\}$). L'action de $h \in G$ sur une connexion discr\`ete se fait par conjugaison :
\begin{equation*}
\begin{split}
&h \cdot \bigl( h_{b_1}, h_{a_1}, \ldots, h_{b_g}, h_{a_g}, h_{m_{g+1}}, \ldots, h_{m_{g+n}} \bigr)\\
= \: &\bigl( h h_{b_1} h^{-1}, h h_{a_1} h^{-1}, \ldots, h h_{b_g} h^{-1}, h h_{a_g} h^{-1}, h h_{m_{g+1}} h^{-1}, \ldots, h h_{m_{g+n}} h^{-1} \bigr). 
\end{split}
\end{equation*}
En d'autres termes :
\[ \begin{array}{ll}
\mathcal{A}_d = \Hom\bigl( \pi_1(\Sigma_{g,n} \!\setminus\! D), G \bigr), & \quad \mathcal{A}_{df} = \Hom\bigl( \pi_1(\Sigma_{g,n}), G \bigr),\\
 \mathcal{A}_d/\mathcal{G}_d = \Hom\bigl( \pi_1(\Sigma_{g,n} \!\setminus\! D), G \bigr)/G, & \quad \mathcal{A}_{df}/\mathcal{G}_d = \Hom\bigl( \pi_1(\Sigma_{g,n}), G \bigr)/G,
\end{array} \]
et on retrouve la vari\'et\'e des caract\`eres.

\smallskip

\indent Il est pertinent pour la suite de d\'ecrire l'alg\`ebre commutative des fonctions $\mathbb{C}[\mathcal{A}_d] = \mathbb{C}[G]^{\otimes (2g+n)}$ en termes de matrices (o\`u $\mathbb{C}[G]$ est l'alg\`ebre des fonctions sur $G$). Soit $V$ une repr\'esentation (de dimension finie) de $G$ avec une base $(v_i)$ et une base duale $(v^j)$. On rappelle que les coefficients matriciels de $V$ dans cette base sont $\overset{V}{T}{^i_j} \in \mathbb{C}[G]$, d\'efinit par $\overset{V}{T}{^i_j}(h) = v^i(h \cdot v_j)$. Ceci donne une matrice $\overset{V}{T}$ qui a ses coefficients dans $\mathbb{C}[G]$. Les coefficients matriciels $\overset{V}{T}{^i_j}$, o\`u $V$ parcourt l'ensemble des $G$-modules de dimension finie, engendre lin\'eairement $\mathbb{C}[G]$. D\'efinissons $\overset{V}{B}(k), \overset{V}{A}(k), \overset{V}{M}(l) \in \mathrm{Mat}_{\dim(V)}\bigl( \mathbb{C}[\mathcal{A}_d] \bigr)$ par
\begin{align*}
\overset{V}{B}(k){^i_j}\bigl( h_{b_1}, h_{a_1}, \ldots, h_{b_g}, h_{a_g}, h_{m_{g+1}}, \ldots, h_{m_{g+n}} \bigr) &= \overset{V}{T}{^i_j}(h_{b_k}),\\
\overset{V}{A}(k){^i_j}\bigl( h_{b_1}, h_{a_1}, \ldots, h_{b_g}, h_{a_g}, h_{m_{g+1}}, \ldots, h_{m_{g+n}} \bigr) &= \overset{V}{T}{^i_j}(h_{a_k}),\\
\overset{V}{M}(l){^i_j}\bigl( h_{b_1}, h_{a_1}, \ldots, h_{b_g}, h_{a_g}, h_{m_{g+1}}, \ldots, h_{m_{g+n}} \bigr) &= \overset{V}{T}{^i_j}(h_{m_l}).\\
\end{align*}
Les coefficients de ces matrices engendrent $\mathbb{C}[\mathcal{A}_d]$ en tant qu'alg\`ebre ($V$ parcourant l'ensemble des $G$-modules de dimension finie). Le groupe de jauge $G$ agit sur $\mathbb{C}[\mathcal{A}_d]$ \`a droite : $(f \cdot h)(x) = f(h\cdot x)$. En termes de matrices, l'action se fait par conjugaison :
\begin{equation}\label{actionGroupeDeJaugeFr}
\forall\, h \in G, \:\:\:\:\: \overset{V}{U}(k) \cdot h = \overset{V}{h} \, \overset{V}{U}(k) \, \overset{V}{h}{^{-1}}
\end{equation}
o\`u $\overset{V}{h} = \overset{V}{T}(h)$ est la repr\'esentation de $h$ sur $V$ et $U$ est $B, A$ ou $M$. Les fonctions invariantes forment une sous-alg\`ebre, $\mathbb{C}[\mathcal{A}_d/G] = \mathbb{C}[\mathcal{A}_d]^G$, et sont appel\'ees observables  (classiques) dans le contexte des th\'eories de jauge discr\`etes. On a ainsi d\'ecrit $\mathbb{C}[\mathcal{A}_d]$ de fa\c{c}on matricielle.

\smallskip

\indent Dans le but de quantifier la structure de Poisson de Fock--Rosly, Alekseev \cite{alekseev}, Alekseev--Grosse--Schomerus \cite{AGS, AGS2} et Buffenoir--Roche \cite{BR, BR2} ont remplac\'e le groupe de Lie $G$ par un groupe quantique $U_q(\mathfrak{g})$, avec $\mathfrak{g} = \mathrm{Lie}(G)$ (mentionnons tout de suite que d'un point de vue purement alg\'ebrique on peut prendre n'importe quelle alg\`ebre de Hopf enrubann\'ee \`a la place de $U_q(\mathfrak{g})$). Ils ont d\'efini une alg\`ebre associative (non-commutative)
\[ \mathcal{L}_{g,n} = \mathbb{C}\bigl\langle \overset{I}{B}(k){^i_j}, \overset{I}{A}(k){^i_j}, \overset{I}{M}(l){^i_j} \: \bigl| \: \text{relations } (\mathcal{R}) \bigr. \bigr\rangle_{I, i, j, k, l} \]
qui est une d\'eformation de $\mathbb{C}[\mathcal{A}_d]$, engendr\'ee par les variables $\overset{I}{B}(k){^i_j}, \overset{I}{A}(k){^i_j}, \overset{I}{M}(l){^i_j}$ o\`u $I$ parcourt maintenant l'ensemble des $U_q(\mathfrak{g})$-modules de dimension finie. Les relations $(\mathcal{R})$ sont donn\'ees sous forme matricielle. Elles font intervenir la $R$-matrice de $U_q(\mathfrak{g})$ et sont con\c{c}ues de sorte que $\mathcal{L}_{g,n}$ est un $U_q(\mathfrak{g})$-module-alg\`ebre \`a droite pour l'action
\[ \forall \, h \in U_q(\mathfrak{g}), \:\:\:\:\: \overset{I}{U}(k) \cdot h = \overset{I}{h'} \, \overset{I}{U}(k) \, \overset{I}{S(h'')} \]
qui est l'analogue de \eqref{actionGroupeDeJaugeFr}, $S$ \'etant l'antipode de $U_q(\mathfrak{g})$, $\Delta(h) = h' \otimes h''$ le coproduit et $U = B, A$ ou $M$. En particulier, on a une sous-alg\`ebre des \'el\'ements invariants $\mathcal{L}^{\mathrm{inv}}_{g,n}$, qui est l'analogue de la sous-alg\`ebre des observables classiques $\mathbb{C}[\mathcal{A}_d]^G$. Notons que $\mathcal{L}_{g,n}$ est vraiment une alg\`ebre quantique de fonctions, dans le sens qu'il est possible d'\'evaluer n'importe quel \'el\'ement de $\mathcal{L}_{g,n}$ sur une connexion discr\`ete; ceci est discut\'e en d\'etail dans les Remarques \ref{F01}, \ref{remarkProductGaugeFields} et dans la section \ref{lienAvecLesLGFT}.

\smallskip

\indent Il est plus difficile de mettre en oeuvre l'analogue quantifi\'e de la contrainte de platitude \eqref{flatnessConstraintFr}, \`a savoir que pour tout $I$ :
\begin{equation}\label{flatnessConstraintMatricesFr}
\overset{I}{C}_{g,n} = \overset{I}{B}(1) \overset{I}{A}(1){^{-1}} \overset{I}{B}(1){^{-1}} \overset{I}{A}(1) \ldots \overset{I}{B}(g) \overset{I}{A}(g){^{-1}} \overset{I}{B}(g){^{-1}} \overset{I}{A}(g) \overset{I}{M}(g+1) \ldots \overset{I}{M}(g+n) = \mathbb{I}_{\dim(I)}.
\end{equation}
En effet, il n'est pas possible de consid\'erer le quotient $\mathcal{L}_{g,n}^{\mathrm{inv}}/\langle (\overset{I}{C}_{g,n})^i_j - \delta^i_j \rangle_{I,i,j}$, puisque les \'el\'ements $(\overset{I}{C}_{g,n})^i_j - \delta^i_j$ ne sont pas invariants. Il n'est pas non plus possible de quotienter toute l'alg\`ebre $\mathcal{L}_{g,n}$ (avant de se restreindre aux \'el\'ements invariants) car l'alg\`ebre qui en r\'esulte peut \^etre r\'eduite \`a $0$\footnote{Par exemple, sous nos hypoth\`eses sur l'alg\`ebre de jauge d\'etaill\'ees ci-dessous, $\mathcal{L}_{1,0}(H)$ est isomorphe \`a une alg\`ebre de matrices sur $\mathbb{C}$ et son seul quotient possible est $0$.}. La mise en oeuvre de la contrainte de platitude, qui donne lieu \`a ``l'alg\`ebre de modules'' (qui est un analogue quantique de $\mathbb{C}[\mathcal{A}_{df}/G]$), diff\`ere selon les auteurs et leurs hypoth\`eses sur l'alg\`ebre de jauge. Par exemple :
\begin{itemize}
\item Dans \cite{AGS, AGS2, AS}, l'alg\`ebre de jauge n'est pas $U_q(\mathfrak{g})$ mais plut\^ot une alg\`ebre de Hopf modulaire (pens\'ee comme une troncation semi-simple de $U_q(\mathfrak{g})$, o\`u $q$ est une racine de l'unit\'e). Ils construisent des ``projecteurs caract\'eristiques'' gr\^ace aux bonnes propri\'et\'es de la $S$-matrice dans le cadre modulaire et utilisent ces projecteurs pour d\'efinir l'alg\`ebre de modules comme une sous-alg\`ebre de $\mathcal{L}_{g,n}^{\mathrm{inv}}$ (c'est-\`a-dire le produit de $\mathcal{L}_{g,n}^{\mathrm{inv}}$ par ces projecteurs). 

\item Dans \cite{BNR}, l'alg\`ebre de jauge est $U_q(\mathfrak{g})$ avec $q$ g\'en\'erique. Ils consid\`erent toutes les matrices $\overset{I}{Y} \in \mathcal{L}_{g,n} \otimes \mathrm{End}_{\mathbb{C}}(I)$ (o\`u $I$ n'est pas fix\'e) qui satisfont $\overset{I}{Y} \cdot h = \overset{I}{h'} \overset{I}{Y} \overset{I}{S(h'')}$ (les produits de $\overset{I}{B}(i), \overset{I}{A}(j), \overset{I}{M}(k)$ sont des exemples imm\'ediats de telles matrices). La trace quantique $\mathrm{tr}\biggl( \overset{I}{g} \overset{I}{Y} \bigl(\overset{I}{C}_{g,n} - \mathbb{I}_{\dim(I)}\bigr) \biggr)$ (o\`u $g$ est l'\'el\'ement pivot) est un \'el\'ement invariant et on peut consid\'erer l'id\'eal $\mathcal{I} \subset \mathcal{L}_{g,n}^{\mathrm{inv}}$ engendr\'e par toutes ces traces quantiques.  Alors ils d\'efinissent l'alg\`ebre de modules comme \'etant $\mathcal{L}_{g,n}^{\mathrm{inv}}/\mathcal{I}$. Cette construction est insuffisante quand l'alg\`ebre de jauge n'est pas semi-simple car dans ce cas il y a des invariants qui ne peuvent pas s'\'ecrire sous la forme $\mathrm{tr}\bigl( \overset{I}{g} \overset{I}{X}\bigr)$.

\item Dans \cite{MW}, l'alg\`ebre de jauge $K$ est de dimension finie et semi-simple pour la construction de l'alg\`ebre de modules. Ils utilisent une coint\'egrale bilat\`ere de $H$ (appel\'ee int\'egrale de Haar dans leur article, et dont l'existence est assur\'ee par ces hypoth\`eses sur $K$) pour construire des projecteurs associ\'es \`a chaque face du graphe $\Gamma$. Ils d\'efinissent alors l'alg\`ebre de modules comme l'image de l'alg\`ebre des observables par ces projecteurs (ceci donne une sous-alg\`ebre de l'alg\`ebre des observables, \'equivalente \`a celle de \cite{AGS, AGS2, AS}). Notons cependant que le formalisme de leur article (qui contient une axiomatisation et une \'etude des th\'eories de jauge avec les alg\`ebres de Hopf) est diff\'erent de celui utilis\'e ici.
\end{itemize}
\noindent Une d\'efinition possible de l'alg\`ebre de modules sous nos hypoth\`eses sur l'alg\`ebre de jauge (qui n'incluent pas la semi-simplicit\'e) est donn\'ee dans la section \ref{introPropLgnFr} ci-dessous.

\smallskip

\indent En plus des articles d\'ej\`a mentionn\'es, ces alg\`ebres de fonctions quantifi\'ees et leurs g\'en\'eralisations apparaissent dans de nombreux travaux. Par exemple : \cite{BFKB2} (d\'efinition d'une comultiplication sur les connexions (duale au produit dans $\mathcal{L}_{g,n}$, voir section \ref{lienAvecLesLGFT}), de l'holonomie et des boucles de Wilson au moyen de multitangles, qui sont des transformations agissant sur le graphe $\Gamma$ et sur les connexions discr\`etes), \cite{BZBJ} (g\'en\'eralisation de  $\mathcal{L}_{g,n}$ dans un cadre cat\'egorique avec l'homologie de factorisation),  \cite{AGPS} (\'etude de $\mathcal{L}_{1,0}$ avec la (super-) alg\`ebre de jauge $\bar U_q\bigl(\mathfrak{gl}(1|1)\bigr)$ et de la repr\'esentation de $\mathrm{SL}_2(\mathbb{Z})$ associ\'ee), \cite{GJS} (vari\'et\'es de caract\`eres quantifi\'ees aux racines de l'unit\'e et d\'efinition de l'alg\`ebre de modules \textit{via} leur proc\'ed\'e de r\'eduction Hamiltonienne quantique).

\smallskip

\indent Comme nous l'avons d\'ej\`a dit, la d\'efinition de  $\mathcal{L}_{g,n}$ est purement alg\'ebrique et n'importe quelle alg\`ebre de Hopf enrubann\'ee $H$ peut jouer le r\^ole de l'alg\`ebre de jauge; l'alg\`ebre de graphe correspondante sera not\'ee $\mathcal{L}_{g,n}(H)$. Dans cette th\`ese, nous supposons de plus que $H$ est de dimension finie et factorisable, mais pas n\'ecessairement semi-simple. L'exemple phare d'une telle alg\`ebre de Hopf sera pour nous le groupe quantique restreint $\bar U_q(\mathfrak{sl}_2)$, not\'e $\bar U_q$ dans la suite; le Chapitre \ref{chapitreUqSl2} est d\'edi\'e aux propri\'et\'es de $\bar U_q$. 

\smallskip

\indent Nous pr\'esentons maintenant nos principaux r\'esultats en d\'etail.

\subsection{Propri\'et\'es de $\mathcal{L}_{g,n}(H)$, impl\'ementation de la contrainte de platitude}\label{introPropLgnFr}

\indent Sous nos hypoth\`eses sur $H$ (dimension finie, factorisable, enrubann\'ee), l'alg\`ebre de boucle $\mathcal{L}_{0,1}(H)$ est isomorphe \`a $H$ (Proposition \ref{IsoRSD}), l'alg\`ebre de anse $\mathcal{L}_{1,0}(H)$ est isomorphe au double d'Heisenberg $\mathcal{H}(\mathcal{O}(H))$ et donc \`a une alg\`ebre de matrices sur $\mathbb{C}$ (Proposition \ref{isoL10Heisenberg}), et l'isomorphisme d'Alekseev reste vrai (Proposition \ref{isoAlekseev}) : 
\[ \mathcal{L}_{0,1}(H) \cong H, \:\:\:\:\:\:\: \mathcal{L}_{1,0}(H) \cong \mathcal{H}(\mathcal{O}(H)) \cong \mathrm{End}_{\mathbb{C}}(H^*), \:\:\:\:\:\:\: \mathcal{L}_{g,n}(H) \cong \mathcal{L}_{1,0}(H)^{\otimes g} \otimes \mathcal{L}_{0,1}(H)^{\otimes n}. \] 
Il s'ensuit que les repr\'esentations ind\'ecomposables de $\mathcal{L}_{g,n}(H)$ sont de la forme
\[ (H^*)^{\otimes g} \otimes I_1 \otimes \ldots \otimes I_n \]
o\`u $H^*$ est l'unique repr\'esentation ind\'ecomposable (et simple) de $\mathcal{H}(\mathcal{O}(H)) \cong \mathrm{End}_{\mathbb{C}}(H^*)$ et $I_1, \ldots, I_n$ sont des repr\'esentations de $H$.

\smallskip

\indent Comme nous ne supposons pas que $H$ est semi-simple, les d\'efinitions de l'alg\`ebre de modules mentionn\'ees plus haut ne peuvent pas \^etre utilis\'ees. A la place, nous allons impl\'ementer la contrainte de platitude au niveau des repr\'esentations : pour chaque repr\'esentation $V$ de $\mathcal{L}_{g,n}(H)$, nous d\'efinissons un sous-espace $\mathrm{Inv}(V)$ par la condition que les matrices $\overset{I}{C}_{g,n}$ agissent trivialement sur $\mathrm{Inv}(V)$ :
\[ \mathrm{Inv}(V) = \left\{ v \in V \, \left| \, \forall \, I, \: \overset{I}{C}_{g,n} \triangleright v = \mathbb{I}_{\dim(I)}v \right.\right\}. \]
Notons que la relation qui d\'efinit $\mathrm{Inv}(V)$ correspond g\'eom\'etriquement \`a recoller le disque $D$ \`a $\Sigma_{g,n} \!\setminus\! D$ (puisque la boucle $c_{g,n}$ qui engendre $\ker\bigl( \pi_1(\Sigma_{g,n} \!\setminus\! D) \to \pi_1(\Sigma_{g,n}) \bigr)$ est d\'etruite).

\smallskip

\indent En fait, le sous-espace $\mathrm{Inv}(V)$ est stable sous l'action de la sous-alg\`ebre des \'el\'ements invariants $\mathcal{L}_{g,n}^{\mathrm{inv}}(H)$ :

\medskip

\noindent \textbf{Th\'eor\`eme \ref{thmInv}.}
{\em 1) Un \'el\'ement $x \in \mathcal{L}_{g,n}(H)$ est invariant sous l'action de l'alg\`ebre de jauge $H$ si, et seulement si, pour tout $H$-module $I$, $\overset{I}{C}_{g,n}x = x\overset{I}{C}_{g,n}$. \\
2) Soit $V$ une repr\'esentation de $\mathcal{L}_{g,n}(H)$. Alors $\mathrm{Inv}(V)$ est stable sous l'action des \'el\'ements invariants et fournit donc une repr\'esentation de $\mathcal{L}^{\mathrm{inv}}_{g,n}(H)$. }

\medskip

\indent Les matrices $\overset{I}{C}_{g,n}$ utilis\'ees dans ce th\'eor\`eme \'etaient d\'ej\`a dans \cite{alekseev} (avec $H = U_q(\mathfrak{g})$, $q$ g\'en\'erique), mais ici nous avons besoin de g\'en\'eraliser et d'adapter la construction des repr\'esentations des \'el\'ements invariants \`a nos hypoth\`eses sur $H$.

\smallskip

\indent Quand $(g,n)=(1,0)$, $\mathcal{L}_{1,0}(H)$ est isomorphe \`a une alg\`ebre de matrices et son unique repr\'esentation ind\'ecomposable (et simple) est $H^*$; dans ce cas, nous avons une repr\'esentation de $\mathcal{L}_{1,0}^{\mathrm{inv}}(H)$ sur $\mathrm{Inv}(H^*) = \mathrm{SLF}(H)$ (Th\'eor\`eme \ref{repInv}), o\`u $\mathrm{SLF}(H)$ est l'alg\`ebre des formes lin\'eaires sym\'etriques sur $H$ :
\[ \mathrm{SLF}(H) = \left\{ \varphi \in H^* \, \left| \: \forall \, x,y \in H, \:\:\: \varphi(xy) = \varphi(yx) \right. \right\}.  \]

\smallskip

\indent La d\'efinition de l'alg\`ebre de modules dans \cite{AGS, AGS2, AS, BNR} requiert de marquer les cercles de bord par des $H$-modules $I_1, \ldots, I_n$. Alors $V = (H^*)^{\otimes g} \otimes I_1 \otimes \ldots \otimes I_n$ est une repr\'esentation de $\mathcal{L}_{g,n}(H)$; soit $\rho_{\mathrm{inv}}^V$ la repr\'esentation associ\'ee de $\mathcal{L}_{g,n}^{\mathrm{inv}}(H)$ sur $\mathrm{Inv}(V)$. L'alg\`ebre $\rho_{\mathrm{inv}}^V\bigl( \mathcal{L}_{g,n}^{\mathrm{inv}}(H) \bigr)$ impl\'emente \`a la fois les contraintes d'invariance et de platitude. C'est donc un analogue quantique de $\mathbb{C}[\mathcal{A}_{df}/G]$ et c'est le candidat naturel pour la d\'efinition de l'alg\`ebre de modules sous nos hypoth\`eses sur $H$ :
\begin{equation}\label{defModuliFr}
\mathcal{M}_{g,n}(H, I_1, \ldots, I_n) \: = \: \rho_{\mathrm{inv}}^{(H^*)^{\otimes g} \otimes I_1 \otimes \ldots \otimes I_n}\bigl( \mathcal{L}_{g,n}^{\mathrm{inv}}(H) \bigr).
\end{equation}
\noindent Cependant, nous n'aurons pas besoin de l'alg\`ebre de modules dans cette th\`ese. \`A la place, l'objet important pour nos desseins est la repr\'esentation $\mathrm{Inv}(V)$ de $\mathcal{L}_{g,n}^{\mathrm{inv}}(H)$ (qui est aussi une repr\'esentation de l'alg\`ebre de modules).

\subsection{Repr\'esentations de groupes de diff\'eotopie}\label{introRepMCGFr}
\indent Le groupe de diff\'eotopie $\mathrm{MCG}(\Sigma_{g,n})$ agit sur $\mathbb{C}\bigl[ \mathrm{Hom}(\pi_1(\Sigma_{g,n}), G)/G \bigr] = \mathbb{C}[\mathcal{A}_{df}/G]$. Dans les Chapitres \ref{chapitreTore} (cas du tore) et \ref{chapitreLgnMCG} (cas g\'en\'eral), nous construisons l'analogue de cette repr\'esentation (pour $n=0$) bas\'e sur $\mathcal{L}_{g,0}(H)$ et nous obtenons une repr\'esentation projective de $\mathrm{MCG}(\Sigma_{g,0} \!\setminus\! D)$ et surtout de $\mathrm{MCG}(\Sigma_{g,0})$. Nous utilisons l'id\'ee, propos\'ee dans \cite{AS}, de remplacer les g\'en\'erateurs $b_i, a_i, m_j$ de $\pi_1(\Sigma_{g,0} \!\setminus\! D)$ par les matrices $\overset{I}{B}(k), \overset{I}{A}(k), \overset{I}{M}(l)$ \`a coefficients dans $\mathcal{L}_{g,0}(H)$. De cette fa\c{c}on, chaque classe de diff\'eotopie $f \in \mathrm{MCG}(\Sigma_{g,0} \!\setminus\! D)$ peut \^etre vue comme un automorphisme de $\mathcal{L}_{g,0}(H)$, not\'e $\widetilde{f}$ et que nous appelons le relev\'e de $f$. En effet, $f$ d\'etermine un automorphisme de $\pi_1(\Sigma_{g,0} \!\setminus\! D)$ et en rempla\c{c}ant les (classes d'homotopie des) boucles par des matrices on d\'efinit $\widetilde{f}$ (\`a une normalisation pr\`es, \textit{cf.} D\'efinition \eqref{liftHomeo}). Gr\^ace aux isomorphismes de la section pr\'ec\'edente
\[ \mathcal{L}_{g,0}(H) \cong \mathcal{L}_{1,0}(H)^{\otimes g} \cong \mathcal{H}(\mathcal{O}(H))^{\otimes g} \cong \mathrm{End}_{\mathbb{C}}\bigl( (H^*)^{\otimes g} \bigr), \]
et on obtient que $\mathcal{L}_{g,0}(H)$ est une alg\`ebre de matrices, son unique repr\'esentation ind\'ecomposable (et simple) \'etant $(H^*)^{\otimes g}$. Il s'ensuit que tout automorphisme de $\mathcal{L}_{g,0}(H)$ est int\'erieur; en particulier, \`a la classe de diff\'eotopie $f$ est associ\'e un \'el\'ement $\widehat{f} \in \mathcal{L}_{g,n}(H)$, unique \`a scalaire pr\`es, et tel que $\widetilde{f}$ est la conjugaison par $\widehat{f}$. Un tel \'el\'ement $\widehat{f}$ est $H$-invariant (Corollaire \eqref{coroFTildeInv}). En repr\'esentant les \'el\'ements $\widehat{f}$ sur $(H^*)^{\otimes g}$ on obtient une repr\'esentation projective de $\mathrm{MCG}(\Sigma_{g,0} \!\setminus\! D)$ et en repr\'esentant les \'el\'ements $\widehat{f}$ sur $\mathrm{Inv}\bigl((H^*)^{\otimes g}\bigr)$ on obtient une repr\'esentation projective de $\mathrm{MCG}(\Sigma_{g,0})$, ce qui correspond au fait que le disque $D$ est ``recoll\'e'' dans $\mathrm{Inv}\bigl((H^*)^{\otimes g}\bigr)$. 

\smallskip

\indent Le cas du tore $\Sigma_{1,0}$ est consid\'er\'e en premier et \`a part car il m\'erite une attention particuli\`ere. Le r\'esultat est \'enonc\'e comme suit :

\medskip

\noindent \textbf{Th\'eor\`eme \ref{repSL2} (cas du tore).} {\em 1) L'affectation
$$ \tau_a \mapsto \rho\!\left(v_A^{-1}\right), \:\:\:\:\: \tau_b \mapsto \rho\!\left(v_B^{-1}\right) $$
o\`u $\rho$ est la repr\'esentation de $\mathcal{L}_{1,0}(H)$ sur $H^*$, d\'efinit une repr\'esentation $\theta_1^D$ de $\mathrm{MCG}(\Sigma_{1,0} \!\setminus\! D)$ sur $H^*$.
\\2) L'affectation
$$ \tau_a \mapsto \rho_{\mathrm{SLF}}\!\left(v_A^{-1}\right), \:\:\:\:\: \tau_b \mapsto \rho_{\mathrm{SLF}}\!\left(v_B^{-1}\right) $$
o\`u $\rho_{\mathrm{SLF}}$ est la repr\'esentation de $\mathcal{L}_{1,0}^{\mathrm{inv}}(H)$ sur $\mathrm{Inv}(H^*) = \mathrm{SLF}(H)$, d\'efinit une repr\'esentation projective $\theta_1$ de $\mathrm{MCG}(\Sigma_{1,0}) = \mathrm{SL}_2(\mathbb{Z})$ sur $\SLF(H)$. Si de plus $S(\psi) = \psi$ pour tout $\psi \in \SLF(H)$, alors ceci d\'efinit en r\'ealit\'e une repr\'esentation projective de $\mathrm{PSL}_2(\mathbb{Z}) = \mathrm{SL}_2(\mathbb{Z})/\{\pm \mathbb{I}_2\}$.}

\medskip

\noindent Les \'el\'ements $v_A^{-1}, v_B^{-1} \in \mathcal{L}_{1,0}^{\mathrm{inv}}(H)$ qui apparaissent dans le th\'eor\`eme sont d\'efinis de la fa\c{c}on suivante. Les coefficients $\bigl( \overset{I}{A}{^i_j}\bigr)_{I,i,j}$ (resp. $\bigl( \overset{I}{B}{^i_j}\bigr)_{I,i,j}$) engendrent une sous-alg\`ebre de $\mathcal{L}_{1,0}(H)$ isomorphe \`a $\mathcal{L}_{0,1}(H)$, qui est lui-m\^eme isomorphe \`a $H$. Donc on a un morphisme $j_A : H \to \mathbb{C}\bigl\langle \overset{I}{A}{^i_j}\bigr\rangle_{I,i,j} \subset \mathcal{L}_{1,0}(H)$ (resp. $j_B : H \to \mathbb{C}\bigl\langle \overset{I}{B}{^i_j}\bigr\rangle_{I,i,j} \subset \mathcal{L}_{1,0}(H)$), et on d\'efinit $v_A^{-1} = j_A(v^{-1})$ (resp. $v_B^{-1} = j_B(v^{-1})$) o\`u $v$ est l'\'el\'ement ruban de $H$. Ces \'el\'ements impl\'ementent les relev\'es $\widetilde{\tau_a}, \widetilde{\tau_b}$ des twists de Dehn $\tau_a, \tau_b$ respectivement (voir Figure \ref{Sigma10}): $v_A^{-1} = \widehat{\tau_a}, v_B^{-1} = \widehat{\tau_b}$ (Proposition \ref{valueHatAlphaBeta}). En combinant les Propositions \ref{propVIntegrale} et \ref{actionAB}, nous obtenons que les repr\'esentations de $\tau_a$ et $\tau_b$ sur $\varphi \in H^*$ sont explicitement donn\'ees par :
\begin{equation}\label{rhoTauATauBIntroFr}
\begin{split}
\theta_1^D(\tau_a)(\varphi) = v_A^{-1} \triangleright \varphi &= \varphi^{v^{-1}},\\
\theta_1^D(\tau_b)(\varphi) = v_B^{-1} \triangleright \varphi &= \mu^l(v)^{-1}\bigl(\mu^l\!\left(g^{-1}v\,?\right) \varphi^v\bigr)^{v^{-1}}
\end{split}
\end{equation}
o\`u pour tout $\beta \in H^*$ et $h \in H$, $\beta^h \in H^*$ est d\'efini par $\beta^h(x) = \beta(hx)$, $\mu^l \in H^*$ est l'int\'egrale \`a gauche de $H$ et $g$ est l'\'el\'ement pivot (la forme lin\'eaire $\mu^l(v)^{-1}\mu^l\!\left(g^{-1}v\,?\right) : x \mapsto \mu^l(v)^{-1}\mu^l\!\left(g^{-1}vx\right)$ est reli\'ee \`a l'inverse de l'\'el\'ement ruban, \textit{cf.} Proposition \ref{propVIntegrale}).

\medskip

\indent Pour une surface de genre quelconque $\Sigma_{g,0}$, le r\'esultat est \'enonc\'e comme suit :

\medskip

\noindent \textbf{Th\'eor\`eme \ref{thmRepMCG} (cas g\'en\'eral).} {\em 1) L'application
\[ \fonc{\theta_g^D}{\mathrm{MCG}(\Sigma_{g,0} \!\setminus\! D)}{\mathrm{GL}\bigl((H^*)^{\otimes g}\bigr)}{f}{\rho(\widehat{f})} \]
o\`u $\rho$ est la repr\'esentation de $\mathcal{L}_{g,0}(H)$ sur $(H^*)^{\otimes g}$, est une repr\'esentation projective.
\\2) L'application
\[ \fonc{\theta_g}{\mathrm{MCG}(\Sigma_{g,0})}{\mathrm{GL}\bigl(\mathrm{Inv}\bigl((H^*)^{\otimes g}\bigr)\bigr)}{f}{\rho_{\mathrm{inv}}(\widehat{f})}\]
o\`u $\rho_{\mathrm{inv}}$ est la repr\'esentation de $\mathcal{L}_{g,0}^{\mathrm{inv}}(H)$ sur $\mathrm{Inv}\bigl((H^*)^{\otimes g}\bigr)$, est une repr\'esentation projective.}

\medskip

\noindent Un r\'esultat similaire a \'et\'e donn\'e dans \cite{AS} sous l'hypoth\`ese que l'alg\`ebre de jauge $H$ est modulaire, l'espace de repr\'esentation \'etant l'alg\`ebre de modules. Ainsi, notre travail fournit une preuve et g\'en\'eralise \`a un cadre non semi-simple leur r\'esultat. Notons qu'on a une repr\'esentation projective car les \'el\'ements $\widehat{f}$ sont d\'efinis \`a un scalaire pr\`es. Notons aussi que $\theta_g$ est juste une restriction de l'espace de repr\'esentation :
\[ \forall \, f \in \mathrm{MCG}(\Sigma_{g,0}), \:\:\:\:\: \theta_g(f) = \theta_g^D(f)_{\bigl| \mathrm{Inv}\bigl( (H^*)^{\otimes g} \bigr)}\bigr.. \]
Le r\'esultat est \'enonc\'e pour $\Sigma_{g,0}$, mais nous discutons son extension \`a une surface g\'en\'erale $\Sigma_{g,n}$ dans la section \ref{CasGeneral}. 

\smallskip

\indent Le relev\'e $\widetilde{\tau_{\gamma}}$ d'un twist de Dehn autour d'une courbe simple $\gamma$ est impl\'ement\'e par conjugaison par l'\'el\'ement $v_{\widetilde{\gamma}}^{-1} \in \mathcal{L}_{g,n}(H)$; en d'autres termes, $\widehat{\tau_{\gamma}} = v_{\widetilde{\gamma}}^{-1}$ (Proposition \ref{propDehnTwist}). Cet \'el\'ement est d\'efini comme ceci. Tout d'abord, exprimons $\gamma$ en fonction des g\'en\'erateurs $b_i, a_i, m_j$ de $\pi_1(\Sigma_{g,n} \!\setminus\! D)$. Puis rempla\c{c}ons $b_i, a_i, m_j$ par les matrices $\overset{I}{B}(i), \overset{I}{A}(i), \overset{I}{M}(j)$ (\`a une normalisation pr\`es par $\overset{I}{v}{^r}$); ceci donne une matrice $\overset{I}{\widetilde{\gamma}}$, appel\'ee le relev\'e de $\gamma$ (Definition \ref{defLiftLoop}). Enfin, les coefficients $\bigl(\overset{I}{\widetilde{\gamma}}{^i_j}\bigr)_{I,i,j}$ satisfont les relations qui d\'efinissent $\mathcal{L}_{0,1}(H)$ (Proposition \ref{propFusionCourbeSimple}), qui est lui-m\^eme isomorphe \`a $H$. Ainsi, nous avons un morphisme $j_{\widetilde{\gamma}} : H \to \mathbb{C}\bigl\langle \overset{I}{\widetilde{\gamma}}{^i_j}\bigr\rangle_{I,i,j} \subset \mathcal{L}_{g,n}(H)$ et nous d\'efinissons $v_{\widetilde{\gamma}}^{-1} = j_{\widetilde{\gamma}}(v^{-1})$, o\`u $v$ est l'\'el\'ement ruban de $H$. Gr\^ace \`a ces \'el\'ements $v_{\widetilde{\gamma}}^{-1}$, nous obtenons des formules pour les repr\'esentations des twists de Dehn autour des courbes ferm\'ees  simples $a_i, b_i, d_i, e_i$ (repr\'esent\'ees dans la Figure \ref{figureCourbesCanoniques}) sur $(H^*)^{\otimes g}$ :

\medskip

\noindent \textbf{Th\'eor\`eme \ref{formulesExplicites}.} 
{\em Soit $\theta_g^D : \mathrm{MCG}(\Sigma_{g,0} \!\setminus\! D) \to \mathrm{PGL}\bigl((H^*)^{\otimes g}\bigr)$ la repr\'esentation projective obtenue dans le Th\'eor\`eme \ref{thmRepMCG}. On a les formules suivantes :
\begin{align*}
\theta_g^D(\tau_{a_i})\bigl(\varphi_1 \otimes \ldots \otimes \varphi_g\bigr) &= \varphi_1 \otimes \ldots \otimes \varphi_{i-1} \otimes \theta_1^D(\tau_a)(\varphi_i) \otimes \varphi_{i+1} \otimes \ldots \otimes \varphi_g, \\
\theta_g^D(\tau_{b_i})\bigl(\varphi_1 \otimes \ldots \otimes \varphi_g\bigr) &= \varphi_1 \otimes \ldots \otimes \varphi_{i-1} \otimes \theta_1^D(\tau_b)(\varphi_i) \otimes \varphi_{i+1} \otimes \ldots \otimes \varphi_g, \\
\theta_g^D(\tau_{d_i})\bigl(\varphi_1 \otimes \ldots \otimes \varphi_g\bigr) &= \varphi_1 \otimes \ldots \otimes \varphi_{i-2} \otimes \varphi_{i-1}\!\left(S^{-1}(a_j)a_k?b_k v''^{-1} b_j\right) \otimes \varphi_i\!\left( S^{-1}(a_l)  S^{-1}(v'^{-1}) a_m ? b_m b_l \right)\\
& \:\:\:\:\: \otimes \varphi_{i+1} \otimes \ldots \otimes \varphi_g, \\
\theta_g^D(\tau_{e_i})\bigl(\varphi_1 \otimes \ldots \otimes \varphi_g\bigr) &= \varphi_1\!\left(S^{-1}\!\left(v^{(2i-2)-1}\right) ? v^{(2i-1)-1} \right) \otimes \ldots \otimes \varphi_{i-1}\!\left(S^{-1}\!\left( v^{(2)-1} \right) ?  v^{(3)-1} \right)\\
&\:\:\:\:\: \otimes \varphi_i\!\left(S^{-1}(a_j)  S^{-1}\!\left(v^{(1)-1}\right) a_k ? b_k b_j \right) \otimes \varphi_{i+1} \otimes \ldots \otimes \varphi_g,
\end{align*}
avec $i \geq 2$ pour les deux derni\`eres \'egalit\'es, $R = a_j \otimes b_j \in H \otimes H$ est la $R$-matrice\footnote{Nous utilisons la sommation implicite sur l'indice dans l'expression de $R$; il ne faut pas confondre les composants $a_j, b_j$ de la $R$-matrice et les boucles $a_i, b_i \in \pi_1(\Sigma_{g,0} \!\setminus\! D)$.}, et les formules pour $\theta_1^{D}(\tau_a), \theta_1^{D}(\tau_b)$ sont donn\'ees dans \eqref{rhoTauATauBIntroFr} ci-dessus.}

\medskip

Pour le tore $\Sigma_{1,0}$ et l'alg\`ebre de jauge $H = \bar U_q = \bar U_q(\mathfrak{sl}_2)$, nous \'etudions explicitement la repr\'esentation projective de $\mathrm{SL}_2(\mathbb{Z})$ sur $\mathrm{SLF}(\bar U_q)$.\footnote{Notons que l'alg\`ebre de Hopf $\bar U_q$ n'est pas tress\'ee. Cependant, l'extension de $\bar U_q$ par une racine carr\'ee de $K$ est enrubann\'ee; de plus, la $R$-matrice et l'\'el\'ement ruban satisfont de bonnes propri\'et\'es qui nous permettent d'appliquer le Th\'eor\`eme \ref{repSL2} avec $H = \bar U_q$. Voir sections \ref{braidedExtension}, \ref{technicalDetails}, \ref{sectionL01Uq}.} Pour ce faire, nous avons besoin d'une base convenable de $\mathrm{SLF}(\bar U_q)$, qui est une alg\`ebre de dimension $3p-1$. Cette base est celle introduite dans \cite{GT} et \cite{arike} et que nous appelons la base GTA; sa d\'efinition est rappel\'ee en d\'etail dans la section \ref{sectionSLF}. Elle contient les caract\`eres $\chi^{\epsilon}_s$ des modules simples $\mathcal{X}^{\epsilon}(s)$, avec $\epsilon \in \{\pm\}$ et $1 \leq s \leq p$; ceci donne $2p$ \'el\'ements. Les $p-1$ formes manquantes, not\'ees $G_s$ ($1 \leq s \leq p-1$), sont construites gr\^ace aux propri\'et\'es des $\bar U_q$-modules projectifs $\mathcal{P}^{\epsilon}(s)$. Une propri\'et\'e importante de cette base pour nos desseins est que ses r\`egles de multiplication, determin\'ees dans le Th\'eor\`eme \ref{ProduitArike} (et ind\'ependamment avant dans \cite{GT}, voir les commentaires au d\'ebut du Chapitre \ref{chapitreUqSl2}), sont simples. Cette particularit\'e nous permet de calculer les formules suivantes :

\medskip

\noindent \textbf{Th\'eor\`eme \ref{actionSL2ZArike}.} 
{\em Soit $\theta_1 : \mathrm{SL}_2(\mathbb{Z}) \to \mathrm{PGL}(\bar U_q^*)$ la repr\'esentation projective obtenue dans le Th\'eor\`eme \ref{repSL2}, avec l'alg\`ebre de jauge $\bar U_q = \bar U_q(\mathfrak{sl}_2)$. Les repr\'esentations des twists de Dehn $\tau_a$ et $\tau_b$ sur la base GTA sont donn\'ees par :
$$ \theta_1(\tau_a)(\chi^{\epsilon}_s) = v^{-1}_{\mathcal{X}^{\epsilon}(s)}\chi^{\epsilon}_s, \:\:\:\:\:\:\: \theta_1(\tau_a)(G_{s'}) = v^{-1}_{\mathcal{X}^+(s')}G_{s'} -  v_{\mathcal{X}^+(s')}^{-1}\hat q\left( \frac{p-s'}{[s']}\chi^+_{s'} - \frac{s'}{[s']}\chi^-_{p-s'} \right) $$
et
\begin{align*}
\theta_1(\tau_b)(\chi^{\epsilon}_s) & = \xi \epsilon(-\epsilon)^{p-1}s q^{-(s^2-1)} \left(\sum_{\ell=1}^{p-1}(-1)^s(-\epsilon)^{p-\ell}\left(q^{\ell s} + q^{-\ell s}\right)\left(\chi^+_{\ell} + \chi^-_{p-\ell}\right) +  \chi^+_p + (-\epsilon)^p(-1)^s\chi^-_p \right)\\  
&\:\:\:+\xi\epsilon (-1)^sq^{-(s^2-1)}\sum_{j=1}^{p-1} (-\epsilon)^{j+1}[j][js]G_j,\\
\theta_1(\tau_b)(G_{s'}) & = \xi (-1)^{s'}q^{-(s'^2-1)}\frac{\hat q p}{[s']}\sum_{j=1}^{p-1}(-1)^{j+1}[j][js']\left(2G_j - \hat q \frac{p-j}{[j]}\chi^+_j + \hat q \frac{j}{[j]}\chi^-_{p-j}\right),
\end{align*}
avec $\epsilon \in \{\pm\}$, $0 \leq s \leq p$, $1 \leq s' \leq p-1$ et $\xi^{-1} = \frac{1-i}{2\sqrt{p}} \frac{\hat q^{p-1}}{[p-1]!} (-1)^p q^{-(p-3)/2}$.}

\medskip

\noindent De ces formules, on d\'eduit la structure de la repr\'esentation :

\medskip

\noindent \textbf{Th\'eor\`eme \ref{thDecRep}.}
{\em Le sous-espace $\mathcal{P} = \vect\!\left(\chi^+_s + \chi^-_{p-s}, \chi^+_p, \chi^-_p\right)_{1 \leq s \leq p-1}$ des caract\`eres des $\bar U_q$-modules projectifs, qui est de dimension $p+1$, est stable sous l'action de $\mathrm{SL}_2(\mathbb{Z})$ calcul\'ee dans le Th\'eor\`eme \ref{actionSL2ZArike}. De plus, il existe une repr\'esentation projective $\mathcal{W}$ de $\mathrm{SL}_2(\mathbb{Z})$, de dimension $p-1$, telle que
$$ \SLF(\bar U_q) = \mathcal{P} \oplus \left(\mathbb{C}^2 \otimes \mathcal{W}\right) $$
o\`u $\mathbb{C}^2$ est la repr\'esentation naturelle de $\mathrm{SL}_2(\mathbb{Z})$ (action par multiplication \`a gauche). Les formules pour l'action sur $\mathcal{W}$ sont donn\'ees dans \eqref{FormulesActionW}.}

\subsection{\'Equivalence avec la repr\'esentation de Lyubashenko}\label{introRepLyubFr}
\indent En utilisant des m\'ethodes cat\'egoriques bas\'ees sur le coend d'une cat\'egorie enrubann\'e, Lyubashenko-Majid \cite{LM} (cas du tore avec une cat\'egorie de modules) et  Lyubashenko \cite{lyu95b, lyu96} (cas g\'en\'eral)  ont construit des repr\'esentations projectives de groupes de diff\'eotopie. Nos hypoth\`eses sur $H$ nous permettent d'appliquer leurs constructions \`a la cat\'egorie enrubann\'ee $\mathrm{mod}_l(H)$, c'est-\`a-dire la cat\'egorie des $H$-modules \`a gauche de dimension finie, et d'obtenir les formules correspondantes. Gr\^ace aux formules de \eqref{rhoTauATauBIntroFr} et du Th\'eor\`eme \ref{formulesExplicites}, nous montrons que ces repr\'esentations sont \'equivalentes \`a celles construites ici :

\medskip

\noindent \textbf{Th\'eor\`eme \ref{EquivalenceLMandSLF} (cas du tore).}
{\em La repr\'esentation projective de $\mathrm{MCG}(\Sigma_{1}) = \mathrm{SL}_2(\mathbb{Z})$ d\'efinie dans le Th\'eor\`eme \ref{repSL2} est \'equivalente \`a celle d\'efinie dans \cite{LM}.}

\medskip

\noindent Pour $H = \bar U_q$, la repr\'esentation de Lyubashenko-Majid de $\mathrm{SL}_2(\mathbb{Z})$ sur $\mathcal{Z}(\bar U_q)$ a \'et\'e \'etudi\'ee explicitement dans \cite{FGST} en relation avec la th\'eorie conforme logarithmique des champs. En particulier, ils ont d\'etermin\'e la structure de la repr\'esentation, et le Th\'eor\`eme \ref{thDecRep} est en parfait accord avec leur r\'esultat.

\medskip

\noindent \textbf{Th\'eor\`eme \ref{thmEquivalenceReps} (cas g\'en\'eral).}
{\em Les repr\'esentations projectives de $\mathrm{MCG}(\Sigma_{g} \!\setminus\! D)$ et $\mathrm{MCG}(\Sigma_{g})$ d\'efinies dans le Th\'eor\`eme \ref{thmRepMCG} sont \'equivalentes \`a celles d\'efinies dans \cite{lyu95b, lyu96}.}

\medskip

\indent Cette \'equivalence est int\'eressante car la construction de la repr\'esentation projective dans le cadre de la quantification combinatoire utilise des techniques diff\'erentes du cadre de Lyubashenko--Majid et Lyubashenko, et est peut-\^etre plus \'el\'ementaire puisque le point de d\'epart est simplement d'imiter l'action du groupe de diff\'eotopie sur le groupe fondamental au niveau de l'alg\`ebre. De plus, bien que les repr\'esentations de groupes de diff\'eotopie r\'esultantes sont \'equivalentes, dans la quantification combinatoire nous avons aussi les alg\`ebres d'observables $\mathcal{L}_{g,n}^{\mathrm{inv}}(H)$ et leurs repr\'esentations; ceci donne lieu \`a des repr\'esentations des alg\`ebres d'\'echeveaux des surfaces (aux racines de l'unit\'e), ce qui est un des sujets du Chapitre \ref{chapitreGraphiqueSkein}.

\subsection{Calcul graphique et th\'eorie d'\'echeveau}\label{introRepSkeinFr}

\indent Dans le Chapitre \ref{chapitreGraphiqueSkein}, nous d\'eveloppons tout d'abord un calcul graphique pour $\mathcal{L}_{g,n}(H)$ et nous reformulons les relations qui d\'efinissent $\mathcal{L}_{g,n}(H)$ en termes de diagrammes. Puis nous utilisons ce calcul graphique pour d\'efinir l'application boucle de Wilson, qui associe un \'el\'ement de $\mathcal{L}_{g,n}(H)$ \`a chaque entrelac parall\'elis\'e, orient\'e et colori\'e (D\'efinition \ref{defDefWilson} et Figure \ref{defWilson}). Notre d\'efinition est compl\`etement naturelle \'etant donn\'e qu'elle est enti\`erement diagrammatique. Elle est \'equivalente mais diff\'erente de celles donn\'es dans \cite{BR2} et \cite{BFKB2}. Il n'est pas difficile de montrer que (entre autres) l'application boucle de Wilson prend ses valeurs dans l'alg\`ebre des  observables et surtout qu'elle est compatible avec le produit en pile de deux entrelacs (ces propri\'et\'es sont aussi dans \cite{BR2} et \cite{BFKB2} bien s\^ur, mais sont d\'emontr\'ees en utilisant leurs formalismes et d\'efinitions respectifs) :

\medskip

\noindent \textbf{Th\'eor\`eme \ref{wilsonStack}.}
{\em L'application boucle de Wilson $W$ est compatible avec le produit en pile :
\[ W(L_1 \ast L_2) = W(L_1)W(L_2). \]}
\noindent Avec notre d\'efinition de $W$, la preuve de ce th\'eor\`eme est enti\`erement diagrammatique (Figure \ref{preuveWilson}).

\smallskip

\indent Quand l'alg\`ebre de jauge est $\bar U_q = \bar U_q(\mathfrak{sl}_2)$, la boucle de Wilson est ind\'ependante de l'orientation de l'entrelac et satisfait la relation d'\'echeveaux du crochet de Kauffman. Ces faits impliquent que nous avons une repr\'esentation de l'alg\`ebre d'\'echeveaux du crochet de Kauffman $\mathcal{S}_q(\Sigma_{g,n} \!\setminus\! D)$ sur n'importe quelle repr\'esentation $V$ de $\mathcal{L}_{g,n}(\bar U_q)$ et si nous nous restreignons \`a $\mathrm{Inv}(V)$, nous obtenons une repr\'esentation de $\mathcal{S}_q(\Sigma_{g,n})$, gr\^ace au fait que le disque $D$ est ``recoll\'e'' dans $\mathrm{Inv}(V)$. La derni\`ere partie de cette affirmation est \'enonc\'ee et prouv\'ee pour $n=0$ dans le th\'eor\`eme ci-dessous, mais est probablement vraie pour tout $g,n$.

\medskip

\noindent \textbf{Th\'eor\`eme \ref{theoRepSkein}.}
{\em 1. Soit $\rho : \mathcal{L}_{g,n}(\bar U_q) \to \mathrm{End}_{\mathbb{C}}(V)$ une repr\'esentation (avec $V = (\bar U_q^*)^{\otimes g} \otimes I_1 \otimes \ldots \otimes I_n$, o\`u $I_1, \ldots, I_n$ sont des repr\'esentations de $\bar U_q$). L'application 
\[ \fleche{\mathcal{S}_q(\Sigma_{g,n} \!\setminus\! D)}{\mathrm{End}_{\mathbb{C}}(V)}{L}{\rho(W(L))} \]
est une repr\'esentation de $\mathcal{S}_q(\Sigma_{g,n} \!\setminus\! D)$.
\\2. Supposons $n=0$ et soit $\rho_{\mathrm{inv}}$ la repr\'esentation de $\mathcal{L}_{g,0}^{\mathrm{inv}}(\bar U_q)$ sur $\mathrm{Inv}\bigl( (\bar U_q^*)^{\otimes g} \bigr)$. L'application
\[ \fleche{\mathcal{S}_q(\Sigma_{g,0})}{\mathrm{End}_{\mathbb{C}}\!\left(\mathrm{Inv}\bigl( (\bar U_q^*)^{\otimes g} \bigr)\right)}{L}{\rho_{\mathrm{inv}}\!\left(W(L^{\mathrm{o}})\right)} \]
est bien d\'efinie et est une repr\'esentation de $\mathcal{S}_q(\Sigma_{g,0})$.}

\medskip

\noindent Dans le th\'eor\`eme, $L^{\mathrm{o}}$ est n'importe quel entrelac dans $(\Sigma_{g,0} \!\setminus\! D) \times [0,1]$ tel que $(j \times \mathrm{id})(L^{\mathrm{o}}) = L$, o\`u $j : \Sigma_{g,0} \!\setminus\! D \to \Sigma_{g,0}$ est l'injection canonique.

\smallskip

\indent Pour le tore $\Sigma_{1,0}$, nous \'etudions explicitement cette repr\'esentation sur $\mathrm{SLF}(\bar U_q)$, en utilisant encore la base GTA et ses r\`egles de multiplication. Il suffit de consid\'erer l'action des boucles $a,b$ puisqu'elles engendrent l'image de l'application boucle de Wilson (Proposition \ref{propWAWB}). Les sous-espaces
\[ \mathcal{P} = \mathrm{vect}\bigl( \chi^+_s + \chi^-_{p-s}, \chi^+_p, \chi^-_p \bigr)_{1 \leq s \leq p-1}, \:\:\:\:\:\: \mathcal{U} = \mathrm{vect}\bigl( \chi^+_s \bigr)_{1 \leq s \leq p-1}, \:\:\:\:\:\: \mathcal{V} = \mathrm{vect}\bigl( G_s \bigr)_{1 \leq s \leq p-1}\]
d\'eterminent la structure de cette repr\'esentation. De plus, nous avons une repr\'esentation naturelle de $\mathcal{S}_q(\Sigma_{1,0})$ sur le module d'\'echeveaux $\mathcal{S}_q(H_1)$, o\`u $H_1$ est un corps \`a anses de genre $1$ (\textit{i.e.}  un anneau \'epaissi). Le module d'\'echeveaux r\'eduit $\mathcal{S}_q^{\mathrm{red}}(H_1)$ est isomorphe \`a un facteur de composition de $\mathrm{SLF}(\bar U_q)$:

\medskip

\noindent \textbf{Propositions \ref{structureRepSqS1} et \ref{propLienSkeinEtL10}.}
{\em $J_1 = \mathcal{P} \subset J_2 = \mathrm{vect}\bigl( \mathcal{P} \cup \mathcal{U} \bigr) \subset J_3 = \mathrm{vect}\bigl( \mathcal{P} \cup \mathcal{U} \cup \mathcal{V}\bigr)$ est une s\'erie de composition de $\mathrm{SLF}(\bar U_q)$ sous l'action de $\mathcal{S}_q(\Sigma_{1,0})$. Cette repr\'esentation est ind\'ecomposable et sa structure est sch\'ematis\'ee par le diagramme suivant :
\[
\xymatrix{
\mathcal{U}  \ar[rd]_{W_{\! B}} & &  \mathcal{V} \ar[ld]^{W_{\! A}, W_{\! B}}\\  
 &\mathcal{P} &
}
\]
De plus, les $\mathcal{S}_q(\Sigma_{1,0})$-modules $\mathcal{S}_q^{\mathrm{red}}(H_1)$ et $\overline{\mathcal{U}} = J_2/J_1$ sont isomorphes.}

\medskip

\noindent Nous conjecturons que la derni\`ere affirmation est vraie en genre quelconque, c'est-\`a-dire que $\mathcal{S}_q^{\mathrm{red}}(H_g)$ est un facteur de composition de $\mathrm{Inv}\bigl( (\bar U_q)^{\otimes g} \bigr)$ sous l'action de $\mathcal{S}_q(\Sigma_{g,0})$ (Conjecture \ref{conjectureSkein}).

\subsection{Perspectives}\label{introPerspectivesFr}
\indent Plusieurs questions et probl\`emes bas\'es sur cette th\`ese peuvent faire l'objet de travaux suppl\'ementaires, en particulier lorsque l'alg\`ebre de jauge est $\bar U_q = \bar U_q(\mathfrak{sl}_2)$.

\smallskip

\indent Le premier groupe de questions concerne la description de l'alg\`ebre des observables $\mathcal{L}_{g,n}^{\mathrm{inv}}(\bar U_q)$. C'est un fait g\'en\'eral (Proposition \ref{propWilsonSimple}) que la valeur de la boucle de Wilson d'une courbe ferm\'ee simple $x \in \pi_1(\Sigma_{g,n} \! \setminus\! D)$ colori\'ee par $I$ est la trace quantique de son relev\'e :
\[ \overset{I}{W}(x) = \mathrm{tr}\bigl( \overset{I}{K^{p+1}} \overset{I}{\widetilde{x}}  \bigr), \]
o\`u le relev\'e $\overset{I}{\widetilde{x}}$ (D\'efinition \ref{defLiftLoop}) est d\'efini en rempla\c{c}ant les g\'en\'erateurs de $\pi_1(\Sigma_{g,n} \! \setminus\! D)$ par des matrices dans l'expression de $x$ \textit{via} la correspondance $\overset{I}{B(i)} \leftrightarrow b_i$, $\overset{I}{A(j)} \leftrightarrow a_j$, $\overset{I}{M(k)} \leftrightarrow m_k$, \`a une normalisation pr\`es. Ces \'el\'ements sont des ``observables semi-simples'' car ils se scindent sur les extensions :
\[ 0 \to I \to V \to J \to 0 \:\: \mathrm{ exacte} \:\:\: \implies \:\:\: \overset{V}{W}(x) = \overset{I}{W}(x) + \overset{J}{W}(x). \]
Dans la section \ref{calculGraphiqueSl2}, tous les entrelacs sont colori\'es par la repr\'esentation fondamentale $\mathcal{X}^+(2)$:
\[ W(x) = \mathrm{tr}\bigl( \overset{\mathcal{X}^+(2)}{K^{p+1}} \overset{\mathcal{X}^+(2)}{\widetilde{x}}  \bigr), \]
C'est suffisant pour retrouver tous les $\overset{\mathcal{X}^{\epsilon}(s)}{W}\!\!(x)$ (o\`u les $\mathcal{X}^{\epsilon}(s)$ sont les $\bar U_q$-modules simples) gr\^ace aux formules (qui sont des cons\'equences de la Proposition \ref{propFusionCourbeSimple} et de \eqref{produitChi}) :
\[ W(x)\overset{\mathcal{X}^{\epsilon}(s)}{W}\!\!(x) = \overset{\mathcal{X}^{\epsilon}(s-1)}{W}\!\!(x) + \overset{\mathcal{X}^{\epsilon}(s+1)}{W}\!\!(x), \:\:\:\:\:\: W(x)\overset{\mathcal{X}^{\epsilon}(p)}{W}\!\!(x) = 2\overset{\mathcal{X}^{\epsilon}(p-1)}{W}\!\!(x) + 2\overset{\mathcal{X}^{-\epsilon}(1)}{W}\!\!(x). \]
En revanche, c'est insuffisant pour retrouver tous les observables. En effet, il y a aussi des observables non semi-simples bas\'es sur les pseudo-traces $G_s$ (voir \eqref{GsAvecTopSoc}) :
\[ V^s(x) = \mathrm{tr}\bigl( \sigma_s \overset{\mathcal{P}^+(s)}{K^{p+1}} \overset{\mathcal{P}^+(s)}{\widetilde{x}} \bigr) + \mathrm{tr}\bigl( \sigma_{p-s} \overset{\mathcal{P}^-(p-s)}{K^{p+1}} \overset{\mathcal{P}^-(p-s)}{\widetilde{x}} \bigr) \:\:\:\:\:\:\: (1 \leq s \leq p-1). \]
Ce $V^s$ est un analogue non semi-simple de l'application boucle de Wilson $W$ et n'est d\'efini que sur les boucles simples pour le moment. Ceci nous am\`ene aux probl\`emes suivants (qui peuvent avoir une solution ou pas) :
\begin{itemize}
\item D\'efinir $V^s(L)$ pour n'importe quel entrelac orient\'e parall\'elis\'e $L \in (\Sigma_{g,n} \!\setminus\! D) \times [0,1]$.
\item D\'eterminer les relations d'\'echeveau satisfaites par l'application $V^s$.
\item Est-ce que la collection des observables $W(x), V^s(x)$ (pour $1 \leq s \leq p-1$ et $x$ une boucle simple)\footnote{Nous pouvons prendre $s=1$ gr\^ace \`a la relation $W(x)V^s(x) = \frac{[s-1]}{[s]}V^{s-1}(x) + \frac{[s+1]}{[s]}V^{s+1}(x)$, qui est due \`a la Proposition \ref{propFusionCourbeSimple} et au Th\'eor\`eme \ref{ProduitArike}.} engendre l'alg\`ebre $\mathcal{L}_{g,n}^{\mathrm{inv}}(\bar U_q)$ ? Plus g\'en\'eralement, d\'ecrire aussi pr\'ecis\'ement que possible $\mathcal{L}_{g,n}^{\mathrm{inv}}(\bar U_q)$.
\end{itemize}

\noindent Le dernier point est probablement tr\`es difficile. \`A la place, on peut recoller le disque $D$ en \'etudiant la repr\'esentation de $\mathcal{L}_{g,n}^{\mathrm{inv}}(\bar U_q)$ sur $\mathrm{Inv}(V)$, o\`u $V = (\bar U_q^*)^{\otimes g} \otimes I_1 \otimes \ldots \otimes I_n$. Ceci pourrait \^etre un premier pas pour comprendre la structure de l'alg\`ebre de modules $\mathcal{M}_{g,n}(H, I_1, \ldots, I_n)$ telle que d\'efinie en \eqref{defModuliFr}. Ces questions sont encore tr\`es difficiles et il vaut mieux se restreindre \`a $n=0$ afin d'\'eviter le choix de $I_1, \ldots, I_n$. La section \ref{sectionConjecture} contient des remarques sur le cas $(g,n) = (1,0)$. La premi\`ere difficult\'e de ce type de questions en genre sup\'erieur est que nous ne connaissons pas de base de $\mathrm{Inv}\bigl( (\bar U_q^*)^{\otimes g} \bigr)$ qui g\'en\'eraliserait la base GTA de $\mathrm{SLF}(\bar U_q) = \mathrm{Inv}(U_q^*)$.

\begin{itemize}
\item D\'eterminer une base ``convenable'' de $\mathrm{Inv}\bigl( (\bar U_q^*)^{\otimes g} \bigr)$.
\item D\'eterminer la structure de la repr\'esentation de $\mathcal{L}_{g,0}^{\mathrm{inv}}(\bar U_q)$ sur $\mathrm{Inv}\bigl( (\bar U_q^*)^{\otimes g} \bigr)$ et d\'eduire des cons\'equences sur la structure de l'alg\`ebre de modules $\mathcal{M}_{g,0}(H)$.
\end{itemize}

\indent La partie semi-simple de $\mathcal{L}_{g,n}^{\mathrm{inv}}(\bar U_q)$ m\'erite cependant une attention particuli\`ere puisque c'est l'image par $W$ (avec tous les entrelacs colori\'es par $\mathcal{X}^+(2)$) de l'alg\`ebre d'\'echeveaux $\mathcal{S}_q\bigl( \Sigma_{g,n} \!\setminus\! D \bigr)$. De plus, $\mathrm{Inv}\bigl( (\bar U_q^*)^{\otimes g} \bigr)$ est une repr\'esentation de $\mathcal{S}_q( \Sigma_{g,0})$. Pour $g=1$, la structure de cette repr\'esentation est d\'etermin\'ee dans la Proposition \ref{structureRepSqS1}, et dans la Proposition \ref{propLienSkeinEtL10} il est montr\'e qu'elle contient la repr\'esentation naturelle de $\mathcal{S}_q( \Sigma_{1,0})$ sur $\mathcal{S}_q^{\mathrm{red}}(H_1)$ en tant que facteur de composition.
\begin{itemize}
\item D\'eterminer la structure de la repr\'esentation de $\mathcal{S}_q( \Sigma_{g,0})$ sur $\mathrm{Inv}\bigl( (\bar U_q^*)^{\otimes g} \bigr)$.
\item Prouver que la repr\'esentation naturelle de $\mathcal{S}_q( \Sigma_{g,0})$ sur $\mathcal{S}_q^{\mathrm{red}}(H_g)$ est un facteur de composition de la repr\'esentation de $\mathcal{S}_q( \Sigma_{g,0})$ sur $\mathrm{Inv}\bigl( (\bar U_q^*)^{\otimes g}\bigr)$ (Conjecture \ref{conjectureSkein}).
\end{itemize}

\smallskip

\indent Un autre probl\`eme est de g\'en\'eraliser le Th\'eor\`eme \ref{thDecRep} en genre sup\'erieur :
\begin{itemize}
\item D\'eterminer la structure de la repr\'esentation projective de $\mathrm{MCG}(\Sigma_{g,0})$ sur $\mathrm{Inv}\bigl( (\bar U_q^*)^{\otimes g} \bigr)$.
\end{itemize}

\smallskip

\indent Enfin, on peut essayer de g\'en\'eraliser les alg\`ebres $\mathcal{L}_{g,n}$ dans un contexte cat\'egorique. Ceci peut avoir plusieurs sens. On sait que $\mathcal{L}_{0,1}(H)$ est un coend (Proposition \ref{L01Coend}); en utilisant ceci, on peut partir d'une cat\'egorie enrubann\'ee avec coend $K$ et tout r\'e\'ecrire de fa\c{c}on cat\'egorique gr\^ace \`a la propri\'et\'e universelle du coend (par exemple, $\mathcal{L}_{1,0}(H)$ serait $K \otimes K$ et son produit serait d\'ecrit par un morphisme qui factorise une certaine famille dinaturelle). Nous pouvons aussi essayer de cat\'egorifier $\mathcal{L}_{g,n}$ (par exemple avec l'alg\`ebre de jauge $U_q(\mathfrak{sl}_2)$, $q$ g\'en\'erique), ce qui signifie qu'on cherche une cat\'egorie mono\"idale $\mathcal{C}_{g,n}$ telle que $K_0(\mathcal{C}_{g,n}) \cong \mathcal{L}_{g,n}$; il serait int\'eressant de voir ce qui joue le r\^ole de l'alg\`ebre des observables dans une telle cat\'egorification.

\section{Introduction in english}
\indent Let $\Sigma_{g,n}$ be a compact oriented surface of genus $g$ with $n$ open disks removed. The ``graph algebra'' $\mathcal{L}_{g,n}$ has been introduced and studied by Alekseev \cite{alekseev}, Alekseev--Grosse--Schomerus \cite{AGS, AGS2} and Buffenoir--Roche \cite{BR, BR2} in the middle of the 1990's, in the program of the combinatorial quantization of the moduli space of flat connections over $\Sigma_{g,n}$. It is an associative (non-commutative) algebra defined by generators and relations, the relations being given in a matrix form. The main theme of this thesis is to apply these algebras to the construction of quantum representations of mapping class groups and of skein algebras of surfaces at roots of unity.

\smallskip

\indent In section \ref{introCombQuant} below we explain the context underlying the combinatorial quantization and the definition of the algebra $\mathcal{L}_{g,n}$. Then from section \ref{introPropLgn} to section \ref{introRepSkein} we state and explain our main results. Finally, the section \ref{introPerspectives} contains conjectures and problems that can serve as a starting point for further work.

\subsection{Combinatorial quantization}\label{introCombQuant}
%

\indent Let us recall quickly the main ingredients of combinatorial quantization. Let $G$ be an algebraic Lie group (generally assumed connected and simply-connected, \textit{e.g.} $G = \mathrm{SL}_2(\mathbb{C})$) and $\Sigma_{g,n}$ be a compact oriented surface of genus $g$ with $n$ open disks removed. We consider the moduli space of flat $G$-connections $\mathcal{M}_{g,n} = \mathcal{A}_f/\mathcal{G}$, where $\mathcal{A} = \Omega^1(\Sigma_{g,n}, \mathfrak{g})$ is identified with the space of all $G$-connections, $\mathcal{A}_f$ is the subspace of flat $G$-connections, and $\mathcal{G} = C^{\infty}(\Sigma_{g,n},G)$ is the gauge group. These objects can be described in a discrete and combinatorial way, using holonomies along the edges of a filling graph $\Gamma = (V,E)$. This is an embedded oriented graph on $\Sigma_{g,n}$ (its vertices $v \in V$ are points of $\Sigma_{g,n}$ and its edges $e \in E$ are simple oriented curves on $\Sigma_{g,n}$ between two vertices which do not intersect pairwise) such that $\Sigma_{g,n} \!\setminus\! \Gamma$ is a union of open disks. Let $\mathcal{A}_d = G^E$. Call an element of $\mathcal{A}_d$ a discrete connection; it is to be thought as the collection $(h_e)_{e \in E}$ of holonomies of a connection along the edges of $\Gamma$. If $\gamma = (e_1, \ldots, e_k)$ is a path in $\Gamma$, we define the discrete holonomy of a discrete connection $(h_e)_{e \in E}$ along $\gamma$ as the product $h_{e_1} \ldots h_{e_k}$. A discrete connection is called flat if its holonomy along any face of the graph is $1$.  This gives a set $\mathcal{A}_{df} \subset \mathcal{A}_d$ of flat discrete connections. Finally, the gauge group $\mathcal{G}$ acts by conjugation on the holonomy along a curve of a connection in $\mathcal{A}$. Hence, we define the discrete gauge group to be $\mathcal{G}_d = G^V$ and its action on discrete connections is $(h_v)_{v \in V} \cdot (h_e)_{e \in E} = (h_{e^-} h_e h_{e^+}^{-1})$, where $e^-$ is the source of $e$ and $e^+$ its target. It is a result that $\mathcal{A}_{df}/\mathcal{G}_d \cong \Hom\bigl( \pi_1(\Sigma_{g,n}), G \bigr)/G$ (usually the quotient is in the sense of geometric invariant theory but our discussion is informal). Hence, this construction is equivalent to the character variety, which is a model for $\mathcal{M}_{g,n}$. For more informations about the moduli space and its combinatorial description, an accessible reference is \cite{labourie}. This description is also called a lattice gauge field theory, see \textit{e.g.} \cite{BFKB}.

\smallskip 

\indent The moduli space $\mathcal{A}_f/\mathcal{G}$ carries the Atiyah--Bott--Goldman Poisson structure \cite{AB, goldman}, namely a Poisson bracket on the algebra of functions $\mathbb{C}[\mathcal{A}_f/\mathcal{G}] = \mathbb{C}[\mathcal{A}_f]^{\mathcal{G}}$. The corresponding Poisson structure on the discretization $\mathcal{A}_{df}/\mathcal{G}_d$ has been described by Fock--Rosly \cite{FockRosly0}; this is a Poisson bracket defined in a matrix way on the algebra of functions $\mathbb{C}[\mathcal{A}_d]$ and which induces a Poison bracket on $\mathbb{C}[\mathcal{A}_{df}]^{\mathcal{G}_d}$ (the gauge group acts on functions on the right in the obvious way). The algebra $\mathcal{L}_{g,n}$ is a quantization of $\mathbb{C}[\mathcal{A}_d]$. We will not need to discuss this fact for the purposes of that thesis. Instead, in the sequel, we simply explain the analogy between $\mathcal{L}_{g,n}$ and $\mathbb{C}[\mathcal{A}_d]$.

\smallskip

\indent Here we always take the graph $\Gamma = \Gamma_{g,n} \subset \Sigma_{g,n}$ with one vertex and whose edges represent a generating system of the fundamental group:
\[ \Gamma_{g,n} = \bigl( \{\bullet\}, \{b_1, a_1, \ldots, b_g, a_g, m_{g+1}, \ldots, m_{g+n} \} \bigr). \]
It is represented below:
\begin{center}

\end{center}
It holds $\Sigma_{g,n} \!\setminus\! \Gamma_{g,n} \cong D$, where $D$ is an open disk. Hence, the closed tubular neighborhood of $\Gamma_{g,n}$ is homeomorphic to $\Sigma_{g,n} \!\setminus\! D$:
\begin{center}

\end{center}
where $[x]$ denotes the free homotopy class of $x \in \pi_1(\Sigma_{g,n} \!\setminus\! D)$. The unique face of the graph $\Gamma_{g,n}$ is the curve induced by the deletion of $D$:
\[ c_{g,n} = b_1 a_1^{-1} b_1^{-1} a_1 \ldots b_g a_g^{-1} b_g^{-1} a_g m_{g+1} \ldots m_{g+n}. \]
With this choice of graph, a discrete connection $A_d \in \mathcal{A}_d$ assigns an element of $G$ to each generator of $\pi_1(\Sigma_{g,n} \!\setminus\! D)$, and can thus be identified with a tuple
\[ A_d = \bigl( h_{b_1}, h_{a_1}, \ldots, h_{b_g}, h_{a_g}, h_{m_{g+1}}, \ldots, h_{m_{g+n}} \bigr) \in G^{2g+n}. \]
A flat discrete connection $A_d \in \mathcal{A}_{df}$ assigns an element of $G$ to each generator of $\pi_1(\Sigma_{g,n}) = \pi_1(\Sigma_{g,n} \!\setminus\! D)/\langle c_{g,n} \rangle$. It is given by a tuple $\bigl( h_{b_1}, h_{a_1}, \ldots, h_{b_g}, h_{a_g}, h_{m_{g+1}}, \ldots, h_{m_{g+n}} \bigr)$ such that
\begin{equation}\label{flatnessConstraint}
\mathrm{Hol}(A_d, c_{g,n}) = h_{b_1} h_{a_1}^{-1} h_{b_1}^{-1} h_{a_1} \ldots h_{b_g} h_{a_g}^{-1} h_{b_g}^{-1} h_{a_g} h_{m_{g+1}} \ldots h_{m_{g+n}} = 1. 
\end{equation}
The discrete gauge group is just $\mathcal{G}_d = G$ (since $V = \{\bullet\}$). The action of $h \in G$ on a discrete connection is by conjugation:
\begin{equation*}
\begin{split}
&h \cdot \bigl( h_{b_1}, h_{a_1}, \ldots, h_{b_g}, h_{a_g}, h_{m_{g+1}}, \ldots, h_{m_{g+n}} \bigr)\\
= \: &\bigl( h h_{b_1} h^{-1}, h h_{a_1} h^{-1}, \ldots, h h_{b_g} h^{-1}, h h_{a_g} h^{-1}, h h_{m_{g+1}} h^{-1}, \ldots, h h_{m_{g+n}} h^{-1} \bigr). 
\end{split}
\end{equation*}
In other words
\[ \begin{array}{ll}
\mathcal{A}_d = \Hom\bigl( \pi_1(\Sigma_{g,n} \!\setminus\! D), G \bigr), & \quad \mathcal{A}_{df} = \Hom\bigl( \pi_1(\Sigma_{g,n}), G \bigr),\\
 \mathcal{A}_d/\mathcal{G}_d = \Hom\bigl( \pi_1(\Sigma_{g,n} \!\setminus\! D), G \bigr)/G, & \quad \mathcal{A}_{df}/\mathcal{G}_d = \Hom\bigl( \pi_1(\Sigma_{g,n}), G \bigr)/G,
\end{array} \]
and we recover the character variety.

\smallskip

\indent For our purposes, it is worthwhile to describe the commutative algebra of functions $\mathbb{C}[\mathcal{A}_d] = \mathbb{C}[G]^{\otimes (2g+n)}$ in terms of matrices (where $\mathbb{C}[G]$ is the algebra of functions on $G$). Let $V$ be a (finite dimensional) representation of $G$ with basis $(v_i)$ and dual basis $(v^j)$. Recall that the matrix coefficients of $V$ in that basis are $\overset{V}{T}{^i_j} \in \mathbb{C}[G]$, defined by $\overset{V}{T}{^i_j}(h) = v^i(h \cdot v_j)$. This gives a matrix $\overset{V}{T}$ with coefficients in $\mathbb{C}[G]$. The matrix coefficients $\overset{V}{T}{^i_j}$, where $V$ runs in the set of finite dimensional $G$-modules, span linearly $\mathbb{C}[G]$. Define $\overset{V}{B}(k), \overset{V}{A}(k), \overset{V}{M}(l) \in \mathrm{Mat}_{\dim(V)}\bigl( \mathbb{C}[\mathcal{A}_d] \bigr)$ by
\begin{align*}
\overset{V}{B}(k){^i_j}\bigl( h_{b_1}, h_{a_1}, \ldots, h_{b_g}, h_{a_g}, h_{m_{g+1}}, \ldots, h_{m_{g+n}} \bigr) &= \overset{V}{T}{^i_j}(h_{b_k}),\\
\overset{V}{A}(k){^i_j}\bigl( h_{b_1}, h_{a_1}, \ldots, h_{b_g}, h_{a_g}, h_{m_{g+1}}, \ldots, h_{m_{g+n}} \bigr) &= \overset{V}{T}{^i_j}(h_{a_k}),\\
\overset{V}{M}(l){^i_j}\bigl( h_{b_1}, h_{a_1}, \ldots, h_{b_g}, h_{a_g}, h_{m_{g+1}}, \ldots, h_{m_{g+n}} \bigr) &= \overset{V}{T}{^i_j}(h_{m_l}).\\
\end{align*}
The coefficients of these matrices span $\mathbb{C}[\mathcal{A}_d]$ as an algebra ($V$ running in the set of finite dimensional $G$-modules). The gauge group $G$ acts on $\mathbb{C}[\mathcal{A}_d]$ on the right: $(f \cdot h)(x) = f(h\cdot x)$. In terms of matrices, the action is by conjugation:
\begin{equation}\label{actionGroupeDeJauge}
\forall\, h \in G, \:\:\:\:\: \overset{V}{U}(k) \cdot h = \overset{V}{h} \, \overset{V}{U}(k) \, \overset{V}{h}{^{-1}}
\end{equation}
where $\overset{V}{h} = \overset{V}{T}(h)$ is the representation of $h$ on $V$ and $U$ is $B, A$ or $M$. The invariant functions form a subalgebra, $\mathbb{C}[\mathcal{A}_d/G] = \mathbb{C}[\mathcal{A}_d]^G$, and are called (classical) observables in the context of lattice gauge field theory. This is the matrix description of $\mathbb{C}[\mathcal{A}_d]$.

\smallskip

\indent In order to quantize the Fock--Rosly Poisson structure, Alekseev \cite{alekseev}, Alekseev--Grosse--Schomerus \cite{AGS, AGS2} and Buffenoir--Roche \cite{BR, BR2} replaced the Lie group $G$ by a quantum group $U_q(\mathfrak{g})$, with $\mathfrak{g} = \mathrm{Lie}(G)$ (but from a purely algebraic point of view we can take any ribbon Hopf algebra instead of $U_q(\mathfrak{g})$). They defined an associative (non-commutative) algebra 
\[ \mathcal{L}_{g,n} = \mathbb{C}\bigl\langle \overset{I}{B}(k){^i_j}, \overset{I}{A}(k){^i_j}, \overset{I}{M}(l){^i_j} \: \bigl| \: \text{relations } (\mathcal{R}) \bigr. \bigr\rangle_{I, i, j, k, l} \]
which is a deformation of $\mathbb{C}[\mathcal{A}_d]$, generated by variables $\overset{I}{B}(k){^i_j}, \overset{I}{A}(k){^i_j}, \overset{I}{M}(l){^i_j}$ where $I$ now runs in the set of finite dimensional $U_q(\mathfrak{g})$-modules. The defining relations $(\mathcal{R})$ are given in a matrix form. They involve the $R$-matrix of $U_q(\mathfrak{g})$ and are designed so that $\mathcal{L}_{g,n}$ is a right $U_q(\mathfrak{g})$-module-algebra for the action
\[ \forall \, h \in U_q(\mathfrak{g}), \:\:\:\:\: \overset{I}{U}(k) \cdot h = \overset{I}{h'} \, \overset{I}{U}(k) \, \overset{I}{S(h'')} \]
which is the analogue of \eqref{actionGroupeDeJauge}, $S$ being the antipode of $U_q(\mathfrak{g})$, $\Delta(h) = h' \otimes h''$ the coproduct and $U = B, A$ ou $M$. In particular, we have a subalgebra of invariant elements $\mathcal{L}^{\mathrm{inv}}_{g,n}$, which is the analogue of the algebra of classical observables $\mathbb{C}[\mathcal{A}_d]^G$. Note that $\mathcal{L}_{g,n}$ really is a quantum algebra of functions, in the sense that it is possible to evaluate any element of $\mathcal{L}_{g,n}$ on a discrete connection; this is discussed in detail in Remarks \ref{F01}, \ref{remarkProductGaugeFields} and in section \ref{lienAvecLesLGFT}.

\smallskip

\indent It is more difficult to implement the quantized analogue of the flatness constraint \eqref{flatnessConstraint}, namely for all $I$:
\begin{equation}\label{flatnessConstraintMatrices}
\overset{I}{C}_{g,n} = \overset{I}{B}(1) \overset{I}{A}(1){^{-1}} \overset{I}{B}(1){^{-1}} \overset{I}{A}(1) \ldots \overset{I}{B}(g) \overset{I}{A}(g){^{-1}} \overset{I}{B}(g){^{-1}} \overset{I}{A}(g) \overset{I}{M}(g+1) \ldots \overset{I}{M}(g+n) = \mathbb{I}_{\dim(I)}.
\end{equation}
Indeed, it is not possible to consider the quotient $\mathcal{L}_{g,n}^{\mathrm{inv}}/\langle (\overset{I}{C}_{g,n})^i_j - \delta^i_j \rangle_{I,i,j}$ since the elements $(\overset{I}{C}_{g,n})^i_j - \delta^i_j$ are not invariant. It is also not possible to quotient the whole algebra $\mathcal{L}_{g,n}$ (before taking the invariant elements) because the resulting algebra may be equal to $0$\footnote{For instance, under our assumptions on the gauge algebra stated below, $\mathcal{L}_{1,0}(H)$ is isomorphic to a matrix algebra over $\mathbb{C}$ and its only possible quotient is $0$.}. The implementation of the flatness constraint, giving rise to the ``moduli algebra'' (which is a quantum analogue of $\mathbb{C}[\mathcal{A}_{df}/G]$), differs depending on the authors and their assumptions on the gauge algebra. For instance:
\begin{itemize}
\item In \cite{AGS, AGS2, AS}, the gauge algebra is not $U_q(\mathfrak{g})$ but rather a modular Hopf algebra (thought of as a semisimple truncation of $U_q(\mathfrak{g})$, where $q$ is a root of unity). They construct ``characteristic projectors'' thanks to the nice properties of the $S$-matrix in the modular setting and use these projectors to define the moduli algebra as a subalgebra of $\mathcal{L}_{g,n}^{\mathrm{inv}}$ (namely the product of $\mathcal{L}_{g,n}^{\mathrm{inv}}$ by these projectors). 

\item In \cite{BNR}, the gauge algebra is $U_q(\mathfrak{g})$ with $q$ generic. They consider all the matrices $\overset{I}{Y} \in \mathcal{L}_{g,n} \otimes \mathrm{End}_{\mathbb{C}}(I)$ (where $I$ is not fixed) satisfying $\overset{I}{Y} \cdot h = \overset{I}{h'} \overset{I}{Y} \overset{I}{S(h'')}$ (products of $\overset{I}{B}(i), \overset{I}{A}(j), \overset{I}{M}(k)$ are obvious examples of such matrices). The quantum trace $\mathrm{tr}\biggl( \overset{I}{g} \overset{I}{Y} \bigl(\overset{I}{C}_{g,n} - \mathbb{I}_{\dim(I)}\bigr) \biggr)$ ($g$ being the pivotal element) is an invariant element and we can consider the ideal $\mathcal{I} \subset \mathcal{L}_{g,n}^{\mathrm{inv}}$ generated by all these quantum traces.  Then they define the moduli algebra as the quotient $\mathcal{L}_{g,n}^{\mathrm{inv}}/\mathcal{I}$. This construction is not sufficient when the gauge algebra is not semisimple because in this case there are invariants which cannot be written as $\mathrm{tr}\bigl( \overset{I}{g} \overset{I}{X}\bigr)$.

\item In \cite{MW}, the gauge algebra $K$ is assumed finite dimensional and semisimple for the construction of the moduli algebra. They use a two-sided cointegral of $H$ (called Haar integral in their paper, and whose existence is guaranteed by these assumptions on $K$) to construct projectors associated to each face of the graph $\Gamma$. Then they define the moduli algebra as the image of the algebra of observables by these projectors (this gives a subalgebra of the algebra of observables, equivalent to the one of \cite{AGS, AGS2, AS}). Note however that the formalism of their paper (which contains an axiomatisation and study of Hopf algebra gauge theory) is different from the one used here.
\end{itemize}
\noindent A possible definition of the moduli algebra under our assumptions on the gauge algebra (which do not include semisimplicity) is given in section \ref{introPropLgn} below.

\smallskip

\indent In addition to the papers already mentionned, these quantized algebras of functions and their generalizations appear in various works. For instance: \cite{BFKB2} (definition of a comultiplication on the connections (dual to the product in $\mathcal{L}_{g,n}$, see section \ref{lienAvecLesLGFT}), of the holonomy and of the Wilson loops by means of multitangles, which are transformations acting on the graph $\Gamma$ and on the discrete connections), \cite{BZBJ} (generalization of $\mathcal{L}_{g,n}$ in a categorical setting using factorization homology),  \cite{AGPS} (study of $\mathcal{L}_{1,0}$ with the gauge (super-) algebra $\bar U_q\bigl(\mathfrak{gl}(1|1)\bigr)$ and of the associated representation of $\mathrm{SL}_2(\mathbb{Z})$), \cite{GJS} (quantized character varieties at roots of unity and definition of the moduli algebra \textit{via} their process of quantum Hamiltonian reduction).

\smallskip

\indent As already said, the definition of $\mathcal{L}_{g,n}$ is purely algebraic and any ribbon Hopf algebra $H$ can play the role of the gauge algebra; the corresponding graph algebra will be denoted $\mathcal{L}_{g,n}(H)$. In that thesis, we assume furthermore  that $H$ is finite dimensional and factorizable, but not necessarily semisimple. For us, the guiding example of such a Hopf algebra will be the restricted quantum group $\bar U_q(\mathfrak{sl}_2)$, denoted $\bar U_q$ in the sequel; Chapter \ref{chapitreUqSl2} is devoted to the properties of $\bar U_q$. 

\smallskip

\indent We now discuss our main results in detail.

\subsection{Properties of $\mathcal{L}_{g,n}(H)$, implementation of the flatness constraint}\label{introPropLgn}

\indent Under our assumptions on $H$ (finite dimensional, factorizable, ribbon), the loop algebra $\mathcal{L}_{0,1}(H)$ is isomorphic to $H$ (Proposition \ref{IsoRSD}), the handle algebra $\mathcal{L}_{1,0}(H)$ is isomorphic to the Heisenberg double $\mathcal{H}(\mathcal{O}(H))$ and thus to a matrix algebra over $\mathbb{C}$ (Proposition \ref{isoL10Heisenberg}), and the Alekseev isomorphism remains valid (Proposition \ref{isoAlekseev}): 
\[ \mathcal{L}_{0,1}(H) \cong H, \:\:\:\:\:\:\: \mathcal{L}_{1,0}(H) \cong \mathcal{H}(\mathcal{O}(H)) \cong \mathrm{End}_{\mathbb{C}}(H^*), \:\:\:\:\:\:\: \mathcal{L}_{g,n}(H) \cong \mathcal{L}_{1,0}(H)^{\otimes g} \otimes \mathcal{L}_{0,1}(H)^{\otimes n}. \] 
It follows that the indecomposable representations of $\mathcal{L}_{g,n}(H)$ have the form
\[ (H^*)^{\otimes g} \otimes I_1 \otimes \ldots \otimes I_n \]
where $H^*$ is the only indecomposable (and simple) representation of $\mathcal{H}(\mathcal{O}(H)) \cong \mathrm{End}_{\mathbb{C}}(H^*)$ and $I_1, \ldots, I_n$ are representations of $H$.

\smallskip

\indent Since we do not assume that $H$ is semisimple, the definitions of the moduli algebra mentionned above cannot be used. Instead, we will implement the flatness constraint \eqref{flatnessConstraintMatrices} at the level of the representations: for each representation $V$ of $\mathcal{L}_{g,n}(H)$, we define a subspace $\mathrm{Inv}(V)$ by the requirement that the matrices $\overset{I}{C}_{g,n}$ act trivially on $\mathrm{Inv}(V)$:
\[ \mathrm{Inv}(V) = \left\{ v \in V \, \left| \, \forall \, I, \: \overset{I}{C}_{g,n} \triangleright v = \mathbb{I}_{\dim(I)}v \right.\right\}. \]
Note that the defining relation of $\mathrm{Inv}(V)$ corresponds geometrically to gluing back the disc $D$ to $\Sigma_{g,n} \!\setminus\! D$ (since the loop $c_{g,n}$ generating $\ker\bigl( \pi_1(\Sigma_{g,n} \!\setminus\! D) \to \pi_1(\Sigma_{g,n}) \bigr)$ is killed).

\smallskip

\indent In fact, the subspace $\mathrm{Inv}(V)$ is stable under the action of the subalgebra of invariant elements $\mathcal{L}_{g,n}^{\mathrm{inv}}(H)$:

\medskip

\noindent \textbf{Theorem \ref{thmInv}.}
{\em 1) An element $x \in \mathcal{L}_{g,n}(H)$ is invariant under the action of the gauge algebra $H$ if, and only if, for every $H$-module $I$, $\overset{I}{C}_{g,n}x = x\overset{I}{C}_{g,n}$. \\
2) Let $V$ be a representation of $\mathcal{L}_{g,n}(H)$. Then $\mathrm{Inv}(V)$ is stable under the action of invariant elements and thus provides a representation of $\mathcal{L}^{\mathrm{inv}}_{g,n}(H)$. }

\medskip

\indent The matrices $\overset{I}{C}_{g,n}$ used in that theorem already appeared in \cite{alekseev} (with $H = U_q(\mathfrak{g})$, $q$ generic), but here we need to generalize and adapt the construction of the representations of the invariant elements to our assumptions on $H$.

\smallskip

\indent When $(g,n)=(1,0)$, $\mathcal{L}_{1,0}(H)$ is isomorphic to a matrix algebra and its unique indecomposable (and simple) representation is $H^*$; in that case, we have a representation of $\mathcal{L}_{1,0}^{\mathrm{inv}}(H)$ on $\mathrm{Inv}(H^*) = \mathrm{SLF}(H)$ (Theorem \ref{repInv}), where $\mathrm{SLF}(H)$ is the algebra of symmetric linear forms on $H$:
\[ \mathrm{SLF}(H) = \left\{ \varphi \in H^* \, \left| \: \forall \, x,y \in H, \:\:\: \varphi(xy) = \varphi(yx) \right. \right\}.  \]

\smallskip

\indent The definition of the moduli algebra in \cite{AGS, AGS2, AS, BNR} requires to label the boundary circles by $H$-modules $I_1, \ldots, I_n$. Then $V = (H^*)^{\otimes g} \otimes I_1 \otimes \ldots \otimes I_n$ is a representation of $\mathcal{L}_{g,n}(H)$; let $\rho_{\mathrm{inv}}^V$ be the corresponding representation of $\mathcal{L}_{g,n}^{\mathrm{inv}}(H)$ on $\mathrm{Inv}(V)$. The algebra $\rho_{\mathrm{inv}}^V\bigl( \mathcal{L}_{g,n}^{\mathrm{inv}}(H) \bigr)$ both implements the invariant and flatness constraints. Hence it is a quantum analogue of $\mathbb{C}[\mathcal{A}_{df}/G]$ and is the natural candidate for the definition of the moduli algebra under our assumptions on $H$:
\begin{equation}\label{defModuli}
\mathcal{M}_{g,n}(H, I_1, \ldots, I_n) \: = \: \rho_{\mathrm{inv}}^{(H^*)^{\otimes g} \otimes I_1 \otimes \ldots \otimes I_n}\bigl( \mathcal{L}_{g,n}^{\mathrm{inv}}(H) \bigr).
\end{equation}
\noindent However, we will not need the moduli algebra in this thesis. The important object for our purposes is instead the representation $\mathrm{Inv}(V)$ of $\mathcal{L}_{g,n}^{\mathrm{inv}}(H)$ (which is also a representation of the moduli algebra).

\subsection{Representations of mapping class groups}\label{introRepMCG}
\indent The mapping class group $\mathrm{MCG}(\Sigma_{g,n})$ acts on $\mathbb{C}\bigl[ \mathrm{Hom}(\pi_1(\Sigma_{g,n}), G)/G \bigr] = \mathbb{C}[\mathcal{A}_{df}/G]$. In Chapters \ref{chapitreTore} (case of the torus) and \ref{chapitreLgnMCG} (general case), we construct the analogue of this representation (for $n=0$) based on $\mathcal{L}_{g,0}(H)$ and we get a projective representation of $\mathrm{MCG}(\Sigma_{g,0} \!\setminus\! D)$ and above all of $\mathrm{MCG}(\Sigma_{g,0})$. We use the idea, proposed in \cite{AS}, of replacing the generators $b_i, a_i, m_j$ of $\pi_1(\Sigma_{g,0} \!\setminus\! D)$ by the matrices $\overset{I}{B}(k), \overset{I}{A}(k), \overset{I}{M}(l)$ with coefficients in $\mathcal{L}_{g,0}(H)$. In that way, each mapping class $f \in \mathrm{MCG}(\Sigma_{g,0} \!\setminus\! D)$ can be seen as an automorphism of $\mathcal{L}_{g,0}(H)$, denoted by $\widetilde{f}$ and which we call the lift of $f$. Indeed, $f$ determines an automorphism of $\pi_1(\Sigma_{g,0} \!\setminus\! D)$ and replacing (homotopy classes of) loops by matrices defines $\widetilde{f}$ (up to some normalization, see Definition \eqref{liftHomeo}). Thanks to the isomorphisms of the previous section
\[ \mathcal{L}_{g,0}(H) \cong \mathcal{L}_{1,0}(H)^{\otimes g} \cong \mathcal{H}(\mathcal{O}(H))^{\otimes g} \cong \mathrm{End}_{\mathbb{C}}\bigl( (H^*)^{\otimes g} \bigr)\]
and hence $\mathcal{L}_{g,0}(H)$ is a matrix algebra, its unique indecomposable (and simple) representation being $(H^*)^{\otimes g}$. It follows that each automorphism of $\mathcal{L}_{g,0}(H)$ is inner; in particular, to the mapping class $f$ is associated an element $\widehat{f} \in \mathcal{L}_{g,n}(H)$, unique up to scalar, and such that $\widetilde{f}$ is the conjugation by $\widehat{f}$. Such an element $\widehat{f}$ is $H$-invariant (Corollary \eqref{coroFTildeInv}). Representing the elements $\widehat{f}$ on $(H^*)^{\otimes g}$ gives a projective representation of $\mathrm{MCG}(\Sigma_{g,0} \!\setminus\! D)$ and representing the elements $\widehat{f}$ on $\mathrm{Inv}\bigl((H^*)^{\otimes g}\bigr)$ gives a projective representation of $\mathrm{MCG}(\Sigma_{g,0})$, which corresponds to the fact that the disk $D$ is ``glued back'' in $\mathrm{Inv}\bigl((H^*)^{\otimes g}\bigr)$. 

\smallskip

\indent The case of the torus $\Sigma_{1,0}$ is considered first because it deserves particular interest. The result is stated as follows:

\medskip

\noindent \textbf{Theorem \ref{repSL2} (case of the torus).} {\em 1) The assignment
$$ \tau_a \mapsto \rho\!\left(v_A^{-1}\right), \:\:\:\:\: \tau_b \mapsto \rho\!\left(v_B^{-1}\right) $$
where $\rho$ is the representation of $\mathcal{L}_{1,0}(H)$ on $H^*$, defines a representation $\theta_1^D$ of $\mathrm{MCG}(\Sigma_{1,0} \!\setminus\! D)$ on $H^*$.
\\2) The assignment
$$ \tau_a \mapsto \rho_{\mathrm{SLF}}\!\left(v_A^{-1}\right), \:\:\:\:\: \tau_b \mapsto \rho_{\mathrm{SLF}}\!\left(v_B^{-1}\right) $$
where $\rho_{\mathrm{SLF}}$ is the representation of $\mathcal{L}_{1,0}^{\mathrm{inv}}(H)$ on $\mathrm{Inv}(H^*) = \mathrm{SLF}(H)$, defines a projective representation $\theta_1$ of $\mathrm{MCG}(\Sigma_{1,0}) = \mathrm{SL}_2(\mathbb{Z})$ on $\SLF(H)$. If moreover $S(\psi) = \psi$ for all $\psi \in \SLF(H)$, then this defines actually a projective representation of $\mathrm{PSL}_2(\mathbb{Z}) = \mathrm{SL}_2(\mathbb{Z})/\{\pm \mathbb{I}_2\}$.}

\medskip

\noindent The elements $v_A^{-1}, v_B^{-1} \in \mathcal{L}_{1,0}^{\mathrm{inv}}(H)$ appearing in the theorem are defined as follows. The coefficients $\bigl( \overset{I}{A}{^i_j}\bigr)_{I,i,j}$ (resp. $\bigl( \overset{I}{B}{^i_j}\bigr)_{I,i,j}$) generate a subalgebra  of $\mathcal{L}_{1,0}(H)$ isomorphic to $\mathcal{L}_{0,1}(H)$, which is itself isomorphic to $H$. Hence we have a morphism $j_A : H \to \mathbb{C}\bigl\langle \overset{I}{A}{^i_j}\bigr\rangle_{I,i,j} \subset \mathcal{L}_{1,0}(H)$ (resp. $j_B : H \to \mathbb{C}\bigl\langle \overset{I}{B}{^i_j}\bigr\rangle_{I,i,j} \subset \mathcal{L}_{1,0}(H)$), and we define $v_A^{-1} = j_A(v^{-1})$ (resp. $v_B^{-1} = j_B(v^{-1})$) where $v$ is the ribbon element of $H$. These elements implements the lifts $\widetilde{\tau_a}, \widetilde{\tau_b}$ of the Dehn twists $\tau_a, \tau_b$ respectively (see Figure \ref{Sigma10}): $v_A^{-1} = \widehat{\tau_a}, v_B^{-1} = \widehat{\tau_b}$ (Proposition \ref{valueHatAlphaBeta}). Combining Propositions \ref{propVIntegrale} and \ref{actionAB}, we get that the representations of $\tau_a$ and $\tau_b$ on $\varphi \in H^*$ are explicitly given by:
\begin{equation}\label{rhoTauATauBIntro}
\begin{split}
\theta_1^D(\tau_a)(\varphi) = v_A^{-1} \triangleright \varphi &= \varphi^{v^{-1}},\\
\theta_1^D(\tau_b)(\varphi) = v_B^{-1} \triangleright \varphi &= \mu^l(v)^{-1}\bigl(\mu^l\!\left(g^{-1}v\,?\right) \varphi^v\bigr)^{v^{-1}}
\end{split}
\end{equation}
where for any $\beta \in H^*$ and $h \in H$, $\beta^h \in H^*$ is defined by $\beta^h(x) = \beta(hx)$, $\mu^l \in H^*$ is the left integral of $H$ and $g$ is the pivotal element (the linear form $\mu^l(v)^{-1}\mu^l\!\left(g^{-1}v\,?\right) : x \mapsto \mu^l(v)^{-1}\mu^l\!\left(g^{-1}vx\right)$ is related to the inverse of the ribbon element, see Proposition \ref{propVIntegrale}).

\medskip

\indent For a surface of arbitrary genus $\Sigma_{g,0}$, the result is stated as follows:

\medskip

\noindent \textbf{Theorem \ref{thmRepMCG} (general case).} {\em 1) The map
\[ \fonc{\theta_g^D}{\mathrm{MCG}(\Sigma_{g,0} \!\setminus\! D)}{\mathrm{GL}\bigl((H^*)^{\otimes g}\bigr)}{f}{\rho(\widehat{f})} \]
where $\rho$ is the representation of $\mathcal{L}_{g,0}(H)$ on $(H^*)^{\otimes g}$, is a projective representation.
\\2) The map
\[ \fonc{\theta_g}{\mathrm{MCG}(\Sigma_{g,0})}{\mathrm{GL}\bigl(\mathrm{Inv}\bigl((H^*)^{\otimes g}\bigr)\bigr)}{f}{\rho_{\mathrm{inv}}(\widehat{f})}\]
where $\rho_{\mathrm{inv}}$ is the representation of $\mathcal{L}_{g,0}^{\mathrm{inv}}(H)$ on $\mathrm{Inv}\bigl((H^*)^{\otimes g}\bigr)$, is a projective representation.}

\medskip

\noindent A similar result was given in \cite{AS} under the assumption that the gauge algebra $H$ is modular, the representation space being the moduli algebra. Thus our work provides a proof and generalizes their result to a non-semisimple setting. Note that we have a projective representation because the elements $\widehat{f}$ are defined only up to scalar. Also note that $\theta_g$ is just a restriction of the representation space:
\[ \forall \, f \in \mathrm{MCG}(\Sigma_{g,0}), \:\:\:\:\: \theta_g(f) = \theta_g^D(f)_{\bigl| \mathrm{Inv}\bigl( (H^*)^{\otimes g} \bigr)}\bigr.. \]
The result is stated for $\Sigma_{g,0}$, but we discuss its extension to a general surface $\Sigma_{g,n}$ in section \ref{CasGeneral}. 

\smallskip

\indent The lift $\widetilde{\tau_{\gamma}}$ of a Dehn twist about a simple closed curve $\gamma$ is implemented by conjugation by an element $v_{\widetilde{\gamma}}^{-1} \in \mathcal{L}_{g,n}(H)$; in other words, $\widehat{\tau_{\gamma}} = v_{\widetilde{\gamma}}^{-1}$ (Proposition \ref{propDehnTwist}). This element is defined as follows. First, express $\gamma$ in terms of the generators $b_i, a_i, m_j$ of $\pi_1(\Sigma_{g,n} \!\setminus\! D)$. Then replace $b_i, a_i, m_j$ by the matrices $\overset{I}{B}(i), \overset{I}{A}(i), \overset{I}{M}(j)$ (up to some normalization by $\overset{I}{v}{^r}$); this gives a matrix $\overset{I}{\widetilde{\gamma}}$, called the lift of $\gamma$ (Definition \ref{defLiftLoop}). Finally, the coefficients $\bigl(\overset{I}{\widetilde{\gamma}}{^i_j}\bigr)_{I,i,j}$ satisfy the defining relations of $\mathcal{L}_{0,1}(H)$ (Proposition \ref{propFusionCourbeSimple}), which is itself isomorphic to $H$. Hence we have a morphism $j_{\widetilde{\gamma}} : H \to \mathbb{C}\bigl\langle \overset{I}{\widetilde{\gamma}}{^i_j}\bigr\rangle_{I,i,j} \subset \mathcal{L}_{g,n}(H)$ and we define $v_{\widetilde{\gamma}}^{-1} = j_{\widetilde{\gamma}}(v^{-1})$, where $v$ is the ribbon element of $H$. Thanks to these elements $v_{\widetilde{\gamma}}^{-1}$, we obtain formulas for the representations of the Dehn twists about the simple closed curves $a_i, b_i, d_i, e_i$ (depicted in Figure \ref{figureCourbesCanoniques}) on $(H^*)^{\otimes g}$:

\medskip

\noindent \textbf{Theorem \ref{formulesExplicites}.} 
{\em Let $\theta_g^D : \mathrm{MCG}(\Sigma_{g,0} \!\setminus\! D) \to \mathrm{PGL}\bigl((H^*)^{\otimes g}\bigr)$ be the projective representation obtained in Theorem \ref{thmRepMCG}. The following formulas hold:
\begin{align*}
\theta_g^D(\tau_{a_i})\bigl(\varphi_1 \otimes \ldots \otimes \varphi_g\bigr) &= \varphi_1 \otimes \ldots \otimes \varphi_{i-1} \otimes \theta_1^D(\tau_a)(\varphi_i) \otimes \varphi_{i+1} \otimes \ldots \otimes \varphi_g, \\
\theta_g^D(\tau_{b_i})\bigl(\varphi_1 \otimes \ldots \otimes \varphi_g\bigr) &= \varphi_1 \otimes \ldots \otimes \varphi_{i-1} \otimes \theta_1^D(\tau_b)(\varphi_i) \otimes \varphi_{i+1} \otimes \ldots \otimes \varphi_g, \\
\theta_g^D(\tau_{d_i})\bigl(\varphi_1 \otimes \ldots \otimes \varphi_g\bigr) &= \varphi_1 \otimes \ldots \otimes \varphi_{i-2} \otimes \varphi_{i-1}\!\left(S^{-1}(a_j)a_k?b_k v''^{-1} b_j\right) \otimes \varphi_i\!\left( S^{-1}(a_l)  S^{-1}(v'^{-1}) a_m ? b_m b_l \right)\\
& \:\:\:\:\: \otimes \varphi_{i+1} \otimes \ldots \otimes \varphi_g, \\
\theta_g^D(\tau_{e_i})\bigl(\varphi_1 \otimes \ldots \otimes \varphi_g\bigr) &= \varphi_1\!\left(S^{-1}\!\left(v^{(2i-2)-1}\right) ? v^{(2i-1)-1} \right) \otimes \ldots \otimes \varphi_{i-1}\!\left(S^{-1}\!\left( v^{(2)-1} \right) ?  v^{(3)-1} \right)\\
&\:\:\:\:\: \otimes \varphi_i\!\left(S^{-1}(a_j)  S^{-1}\!\left(v^{(1)-1}\right) a_k ? b_k b_j \right) \otimes \varphi_{i+1} \otimes \ldots \otimes \varphi_g,
\end{align*}
with $i \geq 2$ for the two last equalities, $R = a_j \otimes b_j \in H \otimes H$ is the $R$-matrix\footnote{We use implicit summation on the index in the expression of $R$; do not confuse the components $a_j, b_j$ of the $R$-matrix and the loops $a_i, b_i \in \pi_1(\Sigma_{g,0} \!\setminus\! D)$.}, and the formulas for $\theta_1^{D}(\tau_a), \theta_1^{D}(\tau_b)$ are given in \eqref{rhoTauATauBIntro} above.}

\medskip

For the torus $\Sigma_{1,0}$ and the gauge algebra $H = \bar U_q = \bar U_q(\mathfrak{sl}_2)$, we explicitly study the projective representation of $\mathrm{SL}_2(\mathbb{Z})$ on $\mathrm{SLF}(\bar U_q)$.\footnote{We mention that the Hopf algebra $\bar U_q$ is not braided. But the extension of $\bar U_q$ by a square root of $K$ is ribbon; moreover, the $R$-matrix and the ribbon element satisfy nice properties which allow us to apply Theorem \ref{repSL2} with $H = \bar U_q$. See sections \ref{braidedExtension}, \ref{technicalDetails}, \ref{sectionL01Uq}.} To do so, we need a suitable basis of $\mathrm{SLF}(\bar U_q)$, which is a $(3p-1)$-dimensional algebra. This basis is the one introduced in \cite{GT} and \cite{arike}, which we call the GTA basis; its definition is recalled in detail in section \ref{sectionSLF}. It contains the characters $\chi^{\epsilon}_s$ of the simple modules $\mathcal{X}^{\epsilon}(s)$, with $\epsilon \in \{\pm\}$ and $1 \leq s \leq p$; this gives $2p$ elements. The missing $p-1$ forms $G_s$ ($1 \leq s \leq p-1$) are constructed thanks to the properties of the $\bar U_q$-projective modules $\mathcal{P}^{\epsilon}(s)$. An important property of this basis for our purposes is that its multiplication rules, determined in Theorem \ref{ProduitArike} (and independently before in \cite{GT}, see the comments at the begining of Chapter \ref{chapitreUqSl2}), are simple. This feature allows us to compute the following formulas:

\medskip

\noindent \textbf{Theorem \ref{actionSL2ZArike}.} 
{\em Let $\theta_1 : \mathrm{SL}_2(\mathbb{Z}) \to \mathrm{PGL}(\bar U_q^*)$ be the projective representation obtained in Theorem \ref{repSL2}, with gauge algebra $\bar U_q = \bar U_q(\mathfrak{sl}_2)$. The representations of the Dehn twists $\tau_a$ and $\tau_b$ on the GTA basis are given by:
$$ \theta_1(\tau_a)(\chi^{\epsilon}_s) = v^{-1}_{\mathcal{X}^{\epsilon}(s)}\chi^{\epsilon}_s, \:\:\:\:\:\:\: \theta_1(\tau_a)(G_{s'}) = v^{-1}_{\mathcal{X}^+(s')}G_{s'} -  v_{\mathcal{X}^+(s')}^{-1}\hat q\left( \frac{p-s'}{[s']}\chi^+_{s'} - \frac{s'}{[s']}\chi^-_{p-s'} \right) $$
and
\begin{align*}
\theta_1(\tau_b)(\chi^{\epsilon}_s) & = \xi \epsilon(-\epsilon)^{p-1}s q^{-(s^2-1)} \left(\sum_{\ell=1}^{p-1}(-1)^s(-\epsilon)^{p-\ell}\left(q^{\ell s} + q^{-\ell s}\right)\left(\chi^+_{\ell} + \chi^-_{p-\ell}\right) +  \chi^+_p + (-\epsilon)^p(-1)^s\chi^-_p \right)\\  
&\:\:\:+\xi\epsilon (-1)^sq^{-(s^2-1)}\sum_{j=1}^{p-1} (-\epsilon)^{j+1}[j][js]G_j,\\
\theta_1(\tau_b)(G_{s'}) & = \xi (-1)^{s'}q^{-(s'^2-1)}\frac{\hat q p}{[s']}\sum_{j=1}^{p-1}(-1)^{j+1}[j][js']\left(2G_j - \hat q \frac{p-j}{[j]}\chi^+_j + \hat q \frac{j}{[j]}\chi^-_{p-j}\right),
\end{align*}
with $\epsilon \in \{\pm\}$, $0 \leq s \leq p$, $1 \leq s' \leq p-1$ and $\xi^{-1} = \frac{1-i}{2\sqrt{p}} \frac{\hat q^{p-1}}{[p-1]!} (-1)^p q^{-(p-3)/2}$.}

\medskip

\noindent From these formulas, we deduce the structure of the representation:

\medskip

\noindent \textbf{Theorem \ref{thDecRep}.}
{\em The $(p+1)$-dimensional subspace $\mathcal{P} = \vect\!\left(\chi^+_s + \chi^-_{p-s}, \chi^+_p, \chi^-_p\right)_{1 \leq s \leq p-1}$ of the characters of the projective $\bar U_q$-modules is stable under the $\mathrm{SL}_2(\mathbb{Z})$-action of Theorem \ref{actionSL2ZArike}. Moreover, there exists a $(p-1)$-dimensional projective representation $\mathcal{W}$ of $\mathrm{SL}_2(\mathbb{Z})$ such that
$$ \SLF(\bar U_q) = \mathcal{P} \oplus \left(\mathbb{C}^2 \otimes \mathcal{W}\right) $$
where $\mathbb{C}^2$ is the natural representation of $\mathrm{SL}_2(\mathbb{Z})$ (action by left multiplication). The formulas for the action on $\mathcal{W}$ are given in \eqref{FormulesActionW}.}

\subsection{Equivalence with the Lyubashenko representation}\label{introRepLyub}
\indent Using categorical methods based on the coend of a ribbon category, Lyubashenko-Majid \cite{LM} (case of the torus and with a category of modules) and  Lyubashenko \cite{lyu95b, lyu96} (general case) constructed projective representations of mapping class groups. Our assumptions on $H$ allow us to apply their constructions to the ribbon category $\mathrm{mod}_l(H)$, namely the category of finite dimensional left $H$-modules, and to obtain the corresponding formulas. Thanks to the formulas of \eqref{rhoTauATauBIntro} and Theorem \ref{formulesExplicites}, we show that these representations are equivalent to those constructed here:

\medskip

\noindent \textbf{Theorem \ref{EquivalenceLMandSLF} (case of the torus).}
{\em The projective representation of $\mathrm{MCG}(\Sigma_{1}) = \mathrm{SL}_2(\mathbb{Z})$ defined in Theorem \ref{repSL2} is equivalent to that defined in \cite{LM}.}

\medskip

\noindent For $H = \bar U_q$, the Lyubashenko-Majid representation of $\mathrm{SL}_2(\mathbb{Z})$ on $\mathcal{Z}(\bar U_q)$ was studied explicitly in \cite{FGST} in relation to logarithmic conformal field theory. In particular, they determined the structure of that representation, and Theorem \ref{thDecRep} is in perfect agreement with their result.

\medskip

\noindent \textbf{Theorem \ref{thmEquivalenceReps} (general case).}
{\em The projective representations of $\mathrm{MCG}(\Sigma_{g} \!\setminus\! D)$ and $\mathrm{MCG}(\Sigma_{g})$ defined in Theorem \ref{thmRepMCG} are equivalent to those defined in \cite{lyu95b, lyu96}.}

\medskip

\indent This equivalence is interesting because the construction of the projective representation of the mapping class group in the combinatorial quantization setting uses different techniques than in the Lyubashenko-Majid and Lyubashenko settings, and is perhaps more elementary since the starting point is simply to mimic the action of the mapping class group on the fundamental group at the level of the algebra. Moreover, even though the resulting representations of mapping class groups are equivalent, in the combinatorial quantization setting we also have the algebras of observables $\mathcal{L}_{g,n}^{\mathrm{inv}}(H)$ and their representations; they give rise to representations of the skein algebras of surfaces (at roots of unity), which is one of the subjects of Chapter \ref{chapitreGraphiqueSkein}.

\subsection{Graphical calculus and skein theory}\label{introRepSkein}

\indent In Chapter \ref{chapitreGraphiqueSkein}, we first develop a graphical calculus for $\mathcal{L}_{g,n}(H)$ and we reformulate the defining relations of $\mathcal{L}_{g,n}(H)$ in terms of diagrams. Then we use this graphical calculus to define the Wilson loop map, which assigns an element of $\mathcal{L}_{g,n}(H)$ to any colored oriented framed link (Definition \ref{defDefWilson} and Figure \ref{defWilson}). Our definition is completely natural since it is entirely diagrammatic. It is equivalent but different from the ones of \cite{BR2} and \cite{BFKB2}. It is not difficult to show that (among other things) the Wilson loop map take values in the algebra of observables and above all that it is compatible with the stack product of two links (these properties are also in \cite{BR2} and \cite{BFKB2} of course, but are proved by using their respective formalisms and definitions):

\medskip

\noindent \textbf{Theorem \ref{wilsonStack}.}
{\em The Wilson loop map $W$ is compatible with the stack product:
\[ W(L_1 \ast L_2) = W(L_1)W(L_2). \]}
\noindent With our definition of $W$, the proof of that theorem is purely diagrammatic (Figure \ref{preuveWilson}).

\smallskip

\indent When the gauge algebra is $\bar U_q = \bar U_q(\mathfrak{sl}_2)$, the Wilson loop map is independent of the orientation of the link and satisfies the Kauffman bracket skein relation. These facts imply that we have a representation of the Kauffman bracket skein algebra $\mathcal{S}_q(\Sigma_{g,n} \!\setminus\! D)$ on any representation $V$ of $\mathcal{L}_{g,n}(\bar U_q)$ and, if we restrict to $\mathrm{Inv}(V)$, we obtain a representation of $\mathcal{S}_q(\Sigma_{g,n})$, due to the fact that the disk $D$ is ``glued back'' in $\mathrm{Inv}(V)$. The last part of the claim is stated and proved for $n=0$ in the theorem below, but is probably true for any $g,n$.

\medskip

\noindent \textbf{Theorem \ref{theoRepSkein}.}
{\em 1. Let $\rho : \mathcal{L}_{g,n}(\bar U_q) \to \mathrm{End}_{\mathbb{C}}(V)$ be a representation (with $V = (\bar U_q^*)^{\otimes g} \otimes I_1 \otimes \ldots \otimes I_n$, where $I_1, \ldots, I_n$ are representations of $\bar U_q$). The map
\[ \fleche{\mathcal{S}_q(\Sigma_{g,n} \!\setminus\! D)}{\mathrm{End}_{\mathbb{C}}(V)}{L}{\rho(W(L))} \]
is a representation of $\mathcal{S}_q(\Sigma_{g,n} \!\setminus\! D)$.
\\2. Assume $n=0$ and let $\rho_{\mathrm{inv}}$ be the representation of $\mathcal{L}_{g,0}^{\mathrm{inv}}(\bar U_q)$ on $\mathrm{Inv}\bigl( (\bar U_q^*)^{\otimes g} \bigr)$. The map
\[ \fleche{\mathcal{S}_q(\Sigma_{g,0})}{\mathrm{End}_{\mathbb{C}}\!\left(\mathrm{Inv}\bigl( (\bar U_q^*)^{\otimes g} \bigr)\right)}{L}{\rho_{\mathrm{inv}}\!\left(W(L^{\mathrm{o}})\right)} \]
is well-defined and is a representation of $\mathcal{S}_q(\Sigma_{g,0})$.}

\medskip

\noindent In the theorem, $L^{\mathrm{o}}$ is any link in $(\Sigma_{g,0} \!\setminus\! D) \times [0,1]$ such that $(j \times \mathrm{id})(L^{\mathrm{o}}) = L$, where $j : \Sigma_{g,0} \!\setminus\! D \to \Sigma_{g,0}$ is the canonical injection.

\smallskip

\indent For the torus $\Sigma_{1,0}$, we explicitly study this representation on $\mathrm{SLF}(\bar U_q)$, using again the GTA basis and its multiplication rules. It suffices to consider the action of the loops $a,b$ since they generate the image of the Wilson loop map (Proposition \ref{propWAWB}). The subspaces 
\[ \mathcal{P} = \mathrm{vect}\bigl( \chi^+_s + \chi^-_{p-s}, \chi^+_p, \chi^-_p \bigr)_{1 \leq s \leq p-1}, \:\:\:\:\:\: \mathcal{U} = \mathrm{vect}\bigl( \chi^+_s \bigr)_{1 \leq s \leq p-1}, \:\:\:\:\:\: \mathcal{V} = \mathrm{vect}\bigl( G_s \bigr)_{1 \leq s \leq p-1}\]
determine the structure of the representation. Moreover, we have a natural representation of $\mathcal{S}_q(\Sigma_{1,0})$ on the skein module $\mathcal{S}_q(H_1)$, where $H_1$ is a genus $1$ handlebody (\textit{i.e.} a thickened annulus). The reduced skein module $\mathcal{S}_q^{\mathrm{red}}(H_1)$ is isomorphic to a composition factor of $\mathrm{SLF}(\bar U_q)$:

\medskip

\noindent \textbf{Propositions \ref{structureRepSqS1} and \ref{propLienSkeinEtL10}.}
{\em $J_1 = \mathcal{P} \subset J_2 = \mathrm{vect}\bigl( \mathcal{P} \cup \mathcal{U} \bigr) \subset J_3 = \mathrm{vect}\bigl( \mathcal{P} \cup \mathcal{U} \cup \mathcal{V}\bigr)$ is a composition series of $\mathrm{SLF}(\bar U_q)$ under the action of $\mathcal{S}_q(\Sigma_{1,0})$. This representation is indecomposable and its structure is schematized by the following diagram:
\[
\xymatrix{
\mathcal{U}  \ar[rd]_{W_{\! B}} & &  \mathcal{V} \ar[ld]^{W_{\! A}, W_{\! B}}\\  
 &\mathcal{P} &
}
\]
Moreover, the $\mathcal{S}_q(\Sigma_{1,0})$-modules $\mathcal{S}_q^{\mathrm{red}}(H_1)$ and $\overline{\mathcal{U}} = J_2/J_1$ are isomorphic.}

\medskip

\noindent We conjecture that the last claim is true in any genus, namely that $\mathcal{S}_q^{\mathrm{red}}(H_g)$ is a composition factor of $\mathrm{Inv}\bigl( (\bar U_q)^{\otimes g} \bigr)$ under the action of $\mathcal{S}_q(\Sigma_{g,0})$ (Conjecture \ref{conjectureSkein}).

\subsection{Perspectives}\label{introPerspectives}
\indent Several questions and problems based upon this thesis can be investigated in further work, especially when the gauge algebra is $\bar U_q = \bar U_q(\mathfrak{sl}_2)$.

\smallskip

\indent The first set of questions is about the description of the algebra of observables $\mathcal{L}_{g,n}^{\mathrm{inv}}(\bar U_q)$. It is a general fact (Proposition \ref{propWilsonSimple}) that the value of the Wilson loop map around a simple loop $x \in \pi_1(\Sigma_{g,n} \! \setminus\! D)$ colored by $I$ is the quantum trace of its lift:
\[ \overset{I}{W}(x) = \mathrm{tr}\bigl( \overset{I}{K^{p+1}} \overset{I}{\widetilde{x}}  \bigr), \]
where the lift $\overset{I}{\widetilde{x}}$ (Definition \ref{defLiftLoop}) is defined by replacing generators of $\pi_1(\Sigma_{g,n} \! \setminus\! D)$ by matrices in the expression of $x$ \textit{via} the correspondence $\overset{I}{B(i)} \leftrightarrow b_i$, $\overset{I}{A(j)} \leftrightarrow a_j$, $\overset{I}{M(k)} \leftrightarrow m_k$, up to some normalization. These elements are ``semisimple observables'' since they split on extensions:
\[ 0 \to I \to V \to J \to 0 \:\: \mathrm{ exact} \:\:\: \implies \:\:\: \overset{V}{W}(x) = \overset{I}{W}(x) + \overset{J}{W}(x). \]
In section \ref{calculGraphiqueSl2}, all the links are colored by the fundamental representation $\mathcal{X}^+(2)$:
\[ W(x) = \mathrm{tr}\bigl( \overset{\mathcal{X}^+(2)}{K^{p+1}} \overset{\mathcal{X}^+(2)}{\widetilde{x}}  \bigr), \]
This is enough to recover all the $\overset{\mathcal{X}^{\epsilon}(s)}{W}\!\!(x)$ (where the $\mathcal{X}^{\epsilon}(s)$ are the simple $\bar U_q$-modules) thanks to the formulas (which follows from Proposition \ref{propFusionCourbeSimple} and \eqref{produitChi}):
\[ W(x)\overset{\mathcal{X}^{\epsilon}(s)}{W}\!\!(x) = \overset{\mathcal{X}^{\epsilon}(s-1)}{W}\!\!(x) + \overset{\mathcal{X}^{\epsilon}(s+1)}{W}\!\!(x), \:\:\:\:\:\: W(x)\overset{\mathcal{X}^{\epsilon}(p)}{W}\!\!(x) = 2\overset{\mathcal{X}^{\epsilon}(p-1)}{W}\!\!(x) + 2\overset{\mathcal{X}^{-\epsilon}(1)}{W}\!\!(x). \]
However, this is insufficient to recover all the observables. Indeed, there are also non-semisimple observables based on the pseudo-traces $G_s$ (see \eqref{GsAvecTopSoc}):
\[ V^s(x) = \mathrm{tr}\bigl( \sigma_s \overset{\mathcal{P}^+(s)}{K^{p+1}} \overset{\mathcal{P}^+(s)}{\widetilde{x}} \bigr) + \mathrm{tr}\bigl( \sigma_{p-s} \overset{\mathcal{P}^-(p-s)}{K^{p+1}} \overset{\mathcal{P}^-(p-s)}{\widetilde{x}} \bigr) \:\:\:\:\:\:\: (1 \leq s \leq p-1). \]
This $V^s$ is a non-semisimple analogue of the Wilson loop map $W$ and is defined only on simple loops for the moment. This leads to the following problems (which may be solvable or not):
\begin{itemize}
\item Define $V^s(L)$ for any oriented framed link $L \in (\Sigma_{g,n} \!\setminus\! D) \times [0,1]$.
\item Determine skein relations satisfied by the map $V^s$.
\item Does the collection of observables $W(x), V^s(x)$ (for $1 \leq s \leq p-1$ and $x$ a simple loop)\footnote{We can take $s=1$ thanks to the relation $W(x)V^s(x) = \frac{[s-1]}{[s]}V^{s-1}(x) + \frac{[s+1]}{[s]}V^{s+1}(x)$, which follows from Proposition \ref{propFusionCourbeSimple} and Theorem \ref{ProduitArike}.} generates the algebra $\mathcal{L}_{g,n}^{\mathrm{inv}}(\bar U_q)$ ? More generally, describe as precisley as possible $\mathcal{L}_{g,n}^{\mathrm{inv}}(\bar U_q)$.
\end{itemize}

\noindent The last item is probably very difficult. Instead, we can glue back the disk $D$ by studying the representation of $\mathcal{L}_{g,n}^{\mathrm{inv}}(\bar U_q)$ on $\mathrm{Inv}(V)$, where $V = (\bar U_q^*)^{\otimes g} \otimes I_1 \otimes \ldots \otimes I_n$. This would be a first step to understand the structure of the moduli algebra $\mathcal{M}_{g,n}(H, I_1, \ldots, I_n)$ defined in \eqref{defModuli}. These questions are still very difficult and it is better to restrict to $n=0$ in order to avoid the choices of $I_1, \ldots, I_n$. Section \ref{sectionConjecture} contains remarks about the case $(g,n) = (1,0)$. The first difficulty for this kind of questions in higher genus is that we do not know a basis of $\mathrm{Inv}\bigl( (\bar U_q^*)^{\otimes g} \bigr)$ which would generalize the GTA basis of $\mathrm{SLF}(\bar U_q) = \mathrm{Inv}(U_q^*)$.

\begin{itemize}
\item Determine a ``suitable basis'' of $\mathrm{Inv}\bigl( (\bar U_q^*)^{\otimes g} \bigr)$.
\item Determine the structure of the representation of $\mathcal{L}_{g,0}^{\mathrm{inv}}(\bar U_q)$ on $\mathrm{Inv}\bigl( (\bar U_q^*)^{\otimes g} \bigr)$ and deduce consequences about the structure of the moduli algebra $\mathcal{M}_{g,0}(H)$. 
\end{itemize}

\indent The semisimple part of $\mathcal{L}_{g,n}^{\mathrm{inv}}(\bar U_q)$ deserves however a special interest since it is the image by $W$ (with all the links colored by $\mathcal{X}^+(2)$) of the Kauffman bracket skein algebra $\mathcal{S}_q\bigl( \Sigma_{g,n} \!\setminus\! D \bigr)$. Moreover, $\mathrm{Inv}\bigl( (\bar U_q^*)^{\otimes g} \bigr)$ is a representation of $\mathcal{S}_q( \Sigma_{g,0})$. For $g=1$, the structure of this representation is determined in Proposition \ref{structureRepSqS1} and it is shown in Proposition \ref{propLienSkeinEtL10} that it contains the natural representation of $\mathcal{S}_q( \Sigma_{1,0})$ on $\mathcal{S}_q^{\mathrm{red}}(H_1)$ as a composition factor.
\begin{itemize}
\item Determine the structure of the representation of $\mathcal{S}_q( \Sigma_{g,0})$ on $\mathrm{Inv}\bigl( (\bar U_q^*)^{\otimes g} \bigr)$.
\item Prove that the natural representation of $\mathcal{S}_q( \Sigma_{g,0})$ on $\mathcal{S}_q^{\mathrm{red}}(H_g)$ is a composition factor of the representation of $\mathcal{S}_q( \Sigma_{g,0})$ on $\mathrm{Inv}\bigl( (\bar U_q^*)^{\otimes g}\bigr)$ (Conjecture \ref{conjectureSkein}).
\end{itemize}

\smallskip

\indent Another problem is to generalize Theorem \ref{thDecRep} in higher genus:
\begin{itemize}
\item Determine the structure of the projective representation of $\mathrm{MCG}(\Sigma_{g,0})$ on $\mathrm{Inv}\bigl( (\bar U_q^*)^{\otimes g} \bigr)$.
\end{itemize}

\smallskip

\indent Finally, we can try to generalize the algebras $\mathcal{L}_{g,n}$ in a categorical setting. This can have differents meanings. We know that $\mathcal{L}_{0,1}(H)$ is a coend (Proposition \ref{L01Coend}); using this, we can start with a ribbon category with coend $K$ and rewrite all in categorical terms using the universal property of the coend (for instance, $\mathcal{L}_{1,0}(H)$ would be $K \otimes K$ and its product would be described as a morphism which factorizes a certain dinatural family). We can also try to categorify $\mathcal{L}_{g,n}$ (for instance with gauge algebra $U_q(\mathfrak{sl}_2)$, $q$ generic), which means finding a monoidal category $\mathcal{C}_{g,n}$ such that $K_0(\mathcal{C}_{g,n}) \cong \mathcal{L}_{g,n}$; it would be interesting to see what plays the role of the algebra of observables in such a categorification.

\chapter{Notations and preliminaries}

In this chapter, we set the notations and collect important facts which will be intensively used in the subsequent chapters. We assume that the reader is familiar with the basic notions about Hopf algebras, knot theory and tensor categories explained for instance in the three first parts of \cite{kassel}.

\section{General notations and conventions}
\indent In order to simplify notations, we will use implicit summations. First, we use Einstein's notation for the computations involving indices: when an index variable appears twice in a litteral expression, one time in upper position and one time in lower position, it implicitly means summation over all the values of the index. Second, we use Sweedler's notation (see \cite[Not. III.1.6]{kassel}) without summation sign for coproducts, that is we write
\begin{equation*}
\begin{split}
&\Delta(x) = x' \otimes x'', \:\:\: \Delta^{(2)}(x) = (\Delta \otimes \mathrm{id}) \circ \Delta(x) = (\mathrm{id} \otimes \Delta) \circ \Delta(x) = x' \otimes x'' \otimes x''', \:\: \ldots,\\
&\Delta^{(n)}(x) = x^{(1)} \otimes \ldots \otimes x^{(n+1)}.
\end{split}
\end{equation*}
We write the universal $R$-matrix as $R = a_i \otimes b_i$ with implicit summation on $i$ and define $R' = b_i \otimes a_i$. We also denote $RR' = X_i \otimes Y_i$, $(RR')^{-1} = \overline{X}_i \otimes \overline{Y}_i$.

\smallskip

\indent The symbol ``?'' will mean a variable in functional constructions. For instance if $H$ is a finite dimensional Hopf algebra and $\varphi, \psi \in H^*$, $a,b \in H$, then for all $x, y \in H$, $\varphi(a?) : x \mapsto \varphi(ax)$, $\varphi(?a) \otimes \psi(b?) : x\otimes y \mapsto \varphi(xa)\psi(by)$ and $\varphi(?a)\psi(b?) : x \mapsto \varphi(x'a)\psi(bx'')$ (thanks to the dual Hopf algebra structure on $H^*$, see section \ref{rappelHopf}). We will often use the notation $\varphi^a$ as a shortand for $\varphi(a?)$.

\smallskip

\indent All the algebras under consideration in this text are finite dimensional $\mathbb{C}$-algebras. If $A$ is a finite dimensional $\mathbb{C}$-algebra, $V$ is a finite dimensional $A$-module and $x \in A$, we denote by $\overset{V}{x} \in \End_{\mathbb{C}}(V)$ the representation of $x$ on the module $V$. Hence, if $(v_i)$ is a basis of $V$, we have 
\begin{equation}\label{repX}
x v_j = \overset{V}{x}{_j^i} v_i.
\end{equation}
More generally, if $X \in A^{\otimes n}$ and if $V_1, \ldots , V_n$ are $A$-modules, we denote by $\overset{V_1 \ldots V_n}{X}$ the representation of $X$ on $V_1 \otimes \ldots \otimes V_n$. 

\smallskip

\indent As in \cite{CR}, we will use the abbreviation PIM for Principal Indecomposable Module. Recall that the PIMs $P_i$ are the indecomposable direct summands of the regular representation of $A$:
\[ _AA = n_1P_1 \oplus \ldots \oplus n_kP_k \]
where the integer $n_i > 0$ is the multiplicity of $P_i$ in that decomposition. Since we assume that $A$ is finite dimensional, it is well-known (see \textit{e.g.} \cite{CR}) that every finite dimensional projective module is a direct sum of PIMs. Hence the projective cover of any finite dimensional module is a direct sum of PIMs, and it follows that any finite dimensional module is a quotient of a direct sum of PIMs. 

\smallskip

\indent We will work only with finite dimensional modules and ``module'' will implicitly mean ``finite dimensional left module''. When we use right modules we explicitly mention it at each time. The {\em socle} of $V$, denoted by $\Soc(V)$ is the largest semi-simple submodule of $V$. The {\em top} of $V$, denoted by $\Top(V)$, is $V/\text{Rad}(V)$, where $\text{Rad}(V)$ is the Jacobson radical of $V$. See \cite[Chap. IV and VIII]{CR} for background material about representation theory.

\smallskip

\indent  For $q \in \mathbb{C}\setminus \{-1,0,1\}$, we define the $q$-integer $[n]$ (with $n \in \mathbb{Z})$ and the $q$-factorial $[m]!$ (with $m \in \mathbb{N})$ by:
$$ [n] = \frac{q^n - q^{-n}}{q-q^{-1}}, \:\:\:\:\: [0]!=1, \:\: [m]! = [1][2] \ldots [m] \: \text{ for } m \geq 1. $$
We will denote $\hat q = q-q^{-1}$ to shorten formulas. In what follows $q$ is a primitive $2p$-root of unity (where $p$ is a fixed integer $\geq 2$), say $q = e^{i\pi/p}$. Observe that in this case $[n] = \frac{\sin(n\pi/p)}{\sin(\pi/p)}$, $[p]=0$ and $[p-n] = [n]$.

\smallskip

\indent As usual $\delta_{s,t}$ or $\delta^s _t$ is the Kronecker symbol and $\mathbb{I}_n$ is the identity matrix of size $n$.

\smallskip

\indent The letter $H$ will always denote a finite dimensional, factorizable, ribbon Hopf algebra (over $\mathbb{C}$). More notations are defined in the next two sections.

\smallskip

\indent We denote by $\Sigma_{g,n}$ the compact oriented surface of genus $g$ with $n$ open disks removed.

\section{Matrices and tensors}\label{matrices}
\indent Let $A$ be a finite dimensional associative algebra. We denote by $\mathrm{Mat}_n(A)$ the algebra of $n \times n$-matrices with coefficients in $A$:
\[ \mathrm{Mat}_n(A) = A \otimes \mathrm{Mat}_n(\mathbb{C}). \]
Every $M \in \mathrm{Mat}_m(A)$ is written as $M = \sum_{i,j} M^i_j \otimes E^i_j$, where $E^i_j$ is the matrix with $1$ at the intersection of the $i$-th row and the $j$-th column and $0$ elsewhere. More generally, every $L \in A \otimes \mathrm{Mat}_{m_1}(\mathbb{C}) \otimes \ldots \otimes \mathrm{Mat}_{m_l}(\mathbb{C})$ can be written as
\[ L = \sum_{i_1,j_1, \ldots, i_l, j_l} L^{i_1 \ldots i_l}_{j_1 \ldots j_l} \otimes  E^{i_1}_{j_1} \otimes \ldots \otimes E^{i_l}_{j_l} \]
and the elements $L^{i_1 \ldots i_l}_{j_1 \ldots j_l} \in A$ are called the coefficients of $L$. If $f : A \to B$ (with $B$ an algebra) is a morphism we define $f(L) \in B \otimes \mathrm{Mat}_{m_1}(\mathbb{C}) \otimes \ldots \otimes \mathrm{Mat}_{m_l}(\mathbb{C})$ component-by-component: 
\[ f(L) = \sum_{i_1,j_1, \ldots, i_l, j_l} f\bigl(L^{i_1 \ldots i_l}_{j_1 \ldots j_l}\bigr) \otimes E^{i_1}_{j_1} \otimes \ldots \otimes E^{i_l}_{j_l} \]
or equivalently $f(L)^{i_1 \ldots i_l}_{j_1 \ldots j_l} = f\bigl(L^{i_1 \ldots i_l}_{j_1 \ldots j_l}\bigr)$.

\smallskip

Let $M \in A \otimes \mathrm{Mat}_m(\mathbb{C})$, $N \in A \otimes \mathrm{Mat}_n(\mathbb{C})$. We embed $M,N$ in $A \otimes \mathrm{Mat}_m(\mathbb{C}) \otimes \mathrm{Mat}_n(\mathbb{C}) = \mathrm{Mat}_{mn}(A)$ by
\[ M_1 = \sum_{i,j} M^i_j \otimes E^i_j \otimes \mathbb{I}_n, \:\:\:\:\:\:\: N_2 = \sum_{i,j} N^i_j \otimes \mathbb{I}_m \otimes E^i_j \]
where $\mathbb{I}_k = \sum_{i=1}^k E^i_i$ is the identity matrix of size $k$. This can also be written as
\[ (M_1)^{ac}_{bd} = M^a_b \delta^c_d, \:\:\:\:\:\:\:\:\: (N_2)^{ac}_{bd} = \delta^a_b N^c_d \]
(where $\delta^i_j$ is the Kronecker symbol), or also
\[ M_1 = M \otimes \mathbb{I}_n, \:\:\:\:\:\:\:\:\: N_2 = \mathbb{I}_m \otimes N \]
where $\otimes$ is the Kronecker product. Note that $M_1N_2$ (resp. $N_2M_1$) contains all the possible products of coefficients of $M$ (resp. of $N$) by coefficients of $N$ (resp. of $M$): $(M_1N_2)^{ik}_{j\ell} = M^i_jN^k_{\ell}$ (resp. $(N_2M_1)^{ik}_{j\ell} = N^k_{\ell}M^i_j$). In particular, $M_1N_2 = N_2M_1$ if and only if the coefficients of $M$ commute with those of $N$.

\begin{exemple}
Let $M, N \in \mathrm{Mat}_2(A)$:
\[ M = 
\begin{pmatrix}
m^1_1 & m^1_2\\
m^2_1 & m^2_2
\end{pmatrix},
\:\:\:\:\:\:\:\:\:
N = 
\begin{pmatrix}
n^1_1 & n^1_2\\
n^2_1 & n^2_2
\end{pmatrix}.
\]
Then
\[ M_1 =
\begin{pmatrix}
m^1_1\mathbb{I}_2 & m^1_2\mathbb{I}_2\\
m^2_1\mathbb{I}_2 & m^2_2\mathbb{I}_2
\end{pmatrix}
=
\begin{pmatrix}
m^1_1 & 0 & m^1_2 & 0\\
0 & m^1_1 & 0 & m^1_2\\
m^2_1 & 0 & m^2_2 & 0\\
0 & m^2_1 & 0 & m^2_2
\end{pmatrix}, \:\:\:\:\:\:\:
N_2 =
\begin{pmatrix}
N & 0\\
0 & N
\end{pmatrix}
=
\begin{pmatrix}
n^1_1 & n^1_2 & 0 & 0\\
n^2_1 & n^2_2 & 0 & 0\\
0 & 0 &n^1_1 & n^1_2\\
0 & 0 & n^2_1 & n^2_2
\end{pmatrix}
\]
and
\[
 M_1N_2 = 
\begin{pmatrix}
m^1_1n^1_1 & m^1_1n^1_2 & m^1_2 n^1_1 & m^1_2n^1_2\\
m^1_1 n^2_1 & m^1_1 n^2_2 & m^1_2 n^2_1 & m^1_2 n^2_2\\
m^2_1 n^1_1 & m^2_1 n^1_2 & m^2_2 n^1_1 & m^2_2n^1_2\\
m^2_1 n^2_1 & m^2_1 n^2_2 & m^2_2 n^2_1 & m^2_2n^2_2
\end{pmatrix}, \:\:\:\:\:\:\:
 N_2 M_1 = 
\begin{pmatrix}
n^1_1m^1_1 & n^1_2m^1_1 & n^1_1m^1_2 & n^1_2m^1_2\\
n^2_1m^1_1 & n^2_2 m^1_1 & n^2_1 m^1_2 & n^2_2 m^1_2\\
n^1_1 m^2_1 & n^1_2 m^2_1 & n^1_1 m^2_2 & n^1_2 m^2_2\\
n^2_1 m^2_1 & n^2_2 m^2_1 & n^2_1 m^2_2 & n^2_2 m^2_2
\end{pmatrix}.
\]
\finEx
\end{exemple}

\indent This is obviously generalized to more general embeddings. For instance, $L \in A \otimes \mathrm{Mat}_{n}(\mathbb{C})^{\otimes 2}$ can be embedded in $A \otimes \mathrm{Mat}_{n}(\mathbb{C})^{\otimes 2}$ in two ways:
\[ (L_{12})^{ac}_{bd} = L^{ac}_{bd}, \:\:\:\:\:\:\:\:\: (L_{21})^{ac}_{bd} = L^{ca}_{db} \]
and can be embedded in $A \otimes \mathrm{Mat}_{n}(\mathbb{C})^{\otimes 3}$ in several ways, \textit{e.g.}
\[ (L_{12})^{ace}_{bdf} = L^{ac}_{bd}\delta^e_f, \:\:\:\:\:\:\:\:\: (L_{13})^{ace}_{bdf} = L^{ae}_{bf}\delta^c_d, \:\:\:\:\:\:\:\:\: (L_{32})^{ace}_{bdf} = L^{ec}_{fd}\delta^a_b\:\: \ldots \]
Note that $\overset{IJ}{(R')}_{12} = \bigl(\overset{JI}{R}\,\bigr)_{21}$, where $R = a_i \otimes b_i \in H^{\otimes 2}$ and $R' = b_i \otimes a_i$.

\smallskip

\indent Recall that we use Einstein's notation for the computations involving indices. For instance if $X, Y \in A \otimes \mathrm{Mat}_{n}(\mathbb{C})$, $L, M \in \mathrm{Mat}_{n}(\mathbb{C})^{\otimes 2}\otimes A$ and $N \in \mathrm{Mat}_{n}(\mathbb{C})^{\otimes 3} \otimes A$, then 
\[
(XY)^a_b = X^a_i Y^i_b, \:\:\:\:\:\:\:\:\: (X_1M_{12})^{ac}_{bd} = X^a_iM^{ic}_{bd}, \:\:\:\:\:\:\:\:\: (L_{32}M_{13}N_{312})^{ace}_{bdf} = L^{ec}_{ij}M^{ai}_{kl}N^{lkj}_{fbd}.
\]
\indent We will extensively use matrices which are labelled by finite dimensional representations of $A$. Such matrices will be denoted by $\overset{V}{M}$, $V$ being a finite dimensional $A$-module. The matrix $\overset{V}{M}$ will be an element of $\mathrm{Mat}_{\dim(V)}(B)$, where $B$ is some algebra, but since we want to record that it is associated with the $A$-module $V$, it is better to consider it as an element of $B \otimes \mathrm{End}_{\mathbb{C}}(V)$.

\indent An important example of such matrices labelled by modules is provided by representation morphisms. Indeed, let $V$ be $A$-module; then by definition we have a map $\overset{V}{T} : A \to \mathrm{End}_{\mathbb{C}}(V)$:
\begin{equation}\label{defMatriceT}
\overset{V}{T}(x) = \overset{V}{x}
\end{equation}
where we denote by $\overset{V}{x}$ the representation of $x \in A$ on $V$. In other words, if $(v_i)$ is a basis of $V$ and $(v^i)$ is its dual basis, then 
\[ \fonc{\overset{V}{T}{^i_j} = v^i(?\cdot v_j)}{A}{\mathbb{C}}{x}{v^i(x\cdot v_j)} \] 
Hence $\overset{V}{T}$ is a matrix with coefficients in $A^*$: $\overset{V}{T} \in A^* \otimes \mathrm{End}_{\mathbb{C}}(V)$. The linear form $\overset{V}{T}{^i_j}$ is called a matrix coefficient associated to the $A$-module $V$. By definition, if $f : V \to W$ is an $A$-morphism it holds 
\begin{equation}\label{naturaliteT}
\overset{W}{T}f = f\overset{V}{T} 
\end{equation}
where we identify $f$ with its matrix. We call this relation the naturality of the (family of) matrices $\overset{V}{T}$. Also note that if $V$ is a submodule or a quotient of $W$, then $\overset{V}{T}$ is a submatrix of $\overset{W}{T}$ and thus the matrix coefficients of $V$ are contained in those of $W$. If $A$ is a Hopf algebra, then $\overset{V}{T}$ is a matrix with coefficients in the dual Hopf algebra $\mathcal{O}(A)$, see \eqref{dualHopf}.

\smallskip

The algebra $A$ being finite dimensional, its dual $A^*$ is generated as a vector space by the matrix coefficients of the PIMs. Indeed, let $(x_1, \ldots, x_n)$ be a basis of $A$ with $x_1 = 1$, let $(x^1, \ldots, x^n) \subset A^*$ be the dual basis and let $_AA$ be the regular representation. It is readily seen that $\overset{_AA}{T}{^i_1}(x_j) = \delta_{i,j}$ and thus $\overset{_AA}{T}{^i_1} = x^i$. Since the PIMs are the direct summands of $_AA$, the claim is proved. 
Note however that the matrix coefficients of the PIMs do not form a basis of $A^*$ in general. Indeed, even if we fix a family $(P_{\alpha})$ representing each isomorphism class of PIMs, it is possible for $P_{\alpha}$ and $P_{\beta}$ to have a composition factor $S$ in common. In this case,  both $\overset{P_{\alpha}}{T}$ and $\overset{P_{\beta}}{T}$ contain $\overset{S}{T}$ as submatrix. This is what happens for $A = \bar U_q(\mathfrak{sl}_2)$, see \eqref{matriceTPim}. In the semi-simple case this phenomenon does not occur.

\smallskip

\indent If $A$ is a Hopf algebra with pivotal element $g$ (see section \ref{rappelHopf}), $V$ is an $A$-module and $\overset{V}{M} \in B \otimes \mathrm{End}_{\mathbb{C}}(V)$ (where $B$ is some algebra), the quantum trace of $\overset{V}{M}$ is
\[ \mathrm{tr}_q(\overset{V}{M}) = \mathrm{tr}(\overset{V}{g}\overset{V}{M}) \in B.  \]

\section{Braided Hopf algebras, factorizability, ribbon element}\label{rappelHopf}
\indent In all this thesis, $H$ is a finite dimensional, factorizable, ribbon Hopf algebra. We recall the meaning of these assumptions and for further use we record some properties of such algebras.

\smallskip

\indent Let $H = (H, \cdot, 1, \Delta, \varepsilon, S, R)$ be a braided Hopf algebra with universal $R$-matrix $R = a_i \otimes b_i$ (see \textit{e.g.} \cite[Chap. VIII]{kassel}). Recall that: 
\begin{align}
& R \, \Delta = \Delta^{\mathrm{op}} \, R, \label{REtCoproduit}\\
&(\Delta \otimes \mathrm{id})(R) = R_{13}R_{23}, \:\:\:\:\: (\mathrm{id} \otimes \Delta)(R) = R_{13}R_{12} \label{deltaR}.\\
&(S \otimes \mathrm{id})(R) = (\mathrm{id} \otimes S^{-1})(R) = R^{-1}, \:\:\:\:\: (S \otimes S)(R) = R \label{propSR}.\\
& R_{12}R_{13}R_{23} = R_{23}R_{13}R_{12} \label{YangBaxter}
\end{align}
with $R_{12} = a_i \otimes b_i \otimes 1, \: R_{13} = a_i \otimes 1 \otimes b_i, \: R_{23} = 1 \otimes a_i \otimes b_i \in H^{\otimes 3}$. The relation (\ref{YangBaxter}) is called the (quantum) Yang-Baxter equation.
\smallskip\\
\indent Consider
\begin{equation*}
\fonc{\Psi}{H^*}{H}{\beta}{(\beta \otimes \mathrm{id})(RR')}
\end{equation*}
where $R' = b_i \otimes a_i$. We say that $H$ is factorizable if $\Psi$ is an isomorphism of vector spaces. Since we assume that $H$ is finite dimensional, we can restrict $\beta$ to be a matrix coefficient $\overset{I}{T}{^i_j}$ of some $H$-module $I$, by the remarks of section \ref{matrices}.
\smallskip\\
\indent Define $R^{(+)} = R$, $R^{(-)} = (R')^{-1}$, and let 
\begin{equation}\label{L}
\overset{I}{L}{^{(\pm)}} = (\overset{I}{T} \otimes \mathrm{id})(R^{(\pm)}) = (\overset{I}{a_i^{(\pm)}})\, b_i^{(\pm)} \in \mathrm{Mat}_{\dim(I)}(H)
\end{equation}
where $R^{(\pm)} = a_i^{(\pm)} \otimes b_i^{(\pm)}$ (note that $\bigl(\overset{I}{a_i^{(\pm)}}\bigr)\, b_i^{(\pm)}$ is the matrix obtained by multiplying each coefficient (which is a scalar) of the matrix $\bigl(\overset{I}{a_i^{(\pm)}}\bigr)$ by the element $b_i \in H$). Recall that $R^{(-)}$ is also a universal $R$-matrix and in particular it satisfies the properties \eqref{REtCoproduit}--\eqref{YangBaxter} above. We use the letters $I,J, \dots$ for modules over $H$. Note that
\[ \overset{I}{L}{^{(\pm)-1}} = S^{-1}(\overset{I}{L}{^{(\pm)}}) = (\overset{I}{a_i^{(\pm)}})\, S^{-1}(b_i^{(\pm)}). \]
If $H$ is factorizable, the coefficients of the matrices $\overset{I}{L} \,\!^{(\pm)}$ generate $H$ as an algebra, since \\$(\overset{I}{T} \otimes \mathrm{id})(RR') = \overset{I}{L}{^{(+)}}\overset{I}{L}{^{(-)-1}}$. These matrices satisfy nice properties which are consequences of (\ref{deltaR}) and (\ref{YangBaxter}):
\begin{equation}\label{propertiesL}
\begin{array}{l}
\overset{I}{L} \,\!^{(\epsilon)}_1\overset{J}{L} \,\!^{(\epsilon)}_2 = \overset{\!\!\!\!\!I\otimes J}{L^{(\epsilon)}_{12}},\\
\overset{IJ}{R}\,\!^{(\epsilon)}_{12} \overset{I}{L} \,\!^{(\epsilon)}_1 \overset{J}{L} \,\!^{(\sigma)}_2 = \overset{J}{L} \,\!^{(\sigma)}_2 \overset{I}{L} \,\!^{(\epsilon)}_1 \overset{IJ}{R}\,\!^{(\epsilon)}_{12} \:\:\:\:\, \forall \, \epsilon, \sigma \in \{\pm\},\\
\overset{IJ}{R}\,\!^{(\epsilon)}_{12} \overset{I}{L} \,\!^{(\sigma)}_1 \overset{J}{L} \,\!^{(\sigma)}_2 = \overset{J}{L} \,\!^{(\sigma)}_2 \overset{I}{L} \,\!^{(\sigma)}_1 \overset{IJ}{R}\,\!^{(\epsilon)}_{12} \:\:\:\: \forall \, \epsilon, \sigma \in \{\pm\},\\
\Delta(\overset{I}{L} \,\!^{(\epsilon)}\,\!^a_b) = \overset{I}{L} \,\!^{(\epsilon)}\,\!^i_b \otimes \overset{I}{L} \,\!^{(\epsilon)}\,\!^a_i, \:\:\:\:\:\:\: \varepsilon(\overset{I}{L}{^{(\epsilon)}}) = \mathbb{I}_{\dim(I)}.
\end{array}
\end{equation}
For instance, here is a proof of the first equality with $\epsilon=+$:
\[ \overset{\!\!\!\!\!I\otimes J}{L^{(+)}_{12}} = \bigl(\overset{I \otimes J}{a_i}\bigr)_{\!12} \, b_i = \bigl(\overset{I}{a_i} \otimes \overset{J}{a_j} \bigr)_{\!12} \, b_i b_j = \bigl(\overset{I}{a_i}\bigr)_{\!1} \bigl(\overset{J}{a_j}\bigr)_{\!2} \, b_i b_j = \bigl(\overset{I}{a_i}\bigr)_{\!1} \, b_i \, \bigl(\overset{J}{a_j}\bigr)_{\!2} \, b_j = \overset{I}{L} \,\!^{(+)}_1\overset{J}{L} \,\!^{(+)}_2 \]
where we used \eqref{deltaR}; note that these are equalities between matrices in $\mathrm{Mat}_{\dim(I)}(\mathbb{C}) \otimes \mathrm{Mat}_{\dim(J)}(\mathbb{C}) \otimes H$ which imply equalities among the coefficients. If the representations $I$ and $J$ are fixed and arbitrary, we will simply write these relations as:
$$ L^{(\epsilon)}_{12} = L^{(\epsilon)}_{1}L^{(\epsilon)}_{2}, \:\:\:\:\: R^{(\epsilon)}_{12} L^{(\epsilon)}_1 L^{(\sigma)}_2 = L^{(\sigma)}_2 L^{(\epsilon)}_1 R^{(\epsilon)}_{12}, \:\:\:\:\: 
R^{(\epsilon)}_{12} L^{(\sigma)}_1 L \,\!^{(\sigma)}_2 = L^{(\sigma)}_2 L^{(\sigma)}_1 R^{(\epsilon)}_{12}, $$
the space $1$ (resp. $2$) corresponding implicitly to the evaluation in the representation $I$ (resp. $J$).

\begin{remark}\label{remarqueGensL}
The set of generators $\bigl(\overset{I}{L}{^{(\pm)}}\bigr)^i_j$ is not at all minimal. However, in practice, due to the fusion relation (the first relation in \eqref{propertiesL}) we can restrict $I$ to belong to a set $\mathcal{G}$ of well-chosen $H$-modules which generate every other module by tensor products, in the sense that every $H$-module is isomorphic to a submodule or a quotient of a tensor product of elements of $\mathcal{G}$. For instance, in the case of $H = \bar U_q(\mathfrak{sl}_2)$, we can restrict $I$ to be the fundamental representation $\mathcal{X}^+(2)$. Also see Remark \ref{remarqueGensT} and section \ref{sectionMatrixCoeffUq}.
\finEx
\end{remark}

\indent Recall that the {\em Drinfeld element} $u$ and its inverse are:
\begin{equation}\label{elementDrinfeld}
u = S(b_i)a_i = b_iS^{-1}(a_i) \:\:\:\: \text{ and } \:\:\:\: u^{-1} = S^{-2}(b_i)a_i = S^{-1}(b_i)S(a_i) = b_iS^2(a_i).
\end{equation}
We assume that $H$ contains a ribbon element $v$. It satisfies:
\begin{equation}\label{ribbon}
v \text{ is central and invertible, } \:\:\:\: \Delta(v) = (R'R)^{-1}v \otimes v, \:\:\:\: S(v) = v, \:\:\:\: \varepsilon(v)=1,  \:\:\: v^2 = uS(u).
\end{equation}
The two last equalities can be deduced easily from the others. A ribbon element is in general not unique. A ribbon Hopf algebra $(H,R,v)$ is a braided Hopf algebra $(H,R)$ together with a ribbon element $v$.
\smallskip\\
\indent We say that $g \in H$ is a {\em pivotal element} if:
\begin{equation}\label{pivot}
\Delta(g) = g \otimes g \:\: \text{ and } \:\: \forall \, x \in H, \:\: S^2(x) = gxg^{-1}.
\end{equation}
Note that $g$ is invertible since it is grouplike: $S(g) = g^{-1}$. A pivotal element is in general not unique. But in a ribbon Hopf algebra $(H,R,v)$ there is a canonical choice:
\begin{equation}\label{pivotCan}
g = uv^{-1}.
\end{equation}
We will always take this canonical pivotal element $g$ in the sequel.

\smallskip

\indent The canonical Hopf algebra structure on $H^*$ is defined by:
\begin{equation}\label{dualHopfDef}
(\varphi \psi)(x) = (\varphi \otimes \psi)(\Delta(x)), \: 1_{H^*} = \varepsilon, \: \Delta(\varphi)(x \otimes y) = \varphi(xy), \: \varepsilon(\varphi) = \varphi(1), \: S(\varphi) = \varphi \circ S.
\end{equation}
with $\varphi, \psi \in H^*$ and $x,y \in H$. When it is endowed with this structure, $H^*$ is called dual Hopf algebra, and denoted $\mathcal{O}(H)$ in the sequel. Recall that $H^*$ is generated as a vector space by the matrix coefficients $\overset{I}{T}{^i_j}$ (see section \ref{matrices}). In terms of these elements, the structure of $\mathcal{O}(H)$ is:
\begin{equation}\label{dualHopf}
\overset{I}{T_1}\overset{J}{T_2} = \overset{I \otimes J}{T}\!\!\!_{12}, \:\: \eta(1) = \overset{\mathbb{C}}{T}, \:\: \Delta(\overset{I}{T^{\, a}_{\, b}}) = \overset{I}{T^{\, a}_{\, i}} \otimes \overset{I}{T^{\, i}_{\, b}}, \:\: \varepsilon(\overset{I}{T}) = \mathbb{I}_{\dim(I)}, \:\: S(\overset{I}{T}) = \overset{I}{T}{^{-1}},
\end{equation}
where $\mathbb{C}$ is the trivial representation. For instance, here are proofs of the first and last equalities (with $h \in H$):
\[ 
\begin{split}
&\bigl(\overset{I \otimes J}{T}\!_{\!\!12}\bigr)(h) = \bigl(\overset{I \otimes J}{h}\bigr)_{12} = \bigl( \overset{I}{h'} \otimes \overset{J}{h''} \bigr)_{12} = \bigl(\overset{I}{h'}\bigr)_1 \, \bigl(\overset{J}{h''}\bigr)_2 = \overset{I}{T_1}(h') \overset{J}{T_2}(h'') = \bigl( \overset{I}{T_1} \overset{J}{T_2} \bigr)(h),\\
&\bigl(S(\overset{I}{T}) \overset{I}{T}\bigr)(h) = S(\overset{I}{T})(h') \overset{I}{T}(h'') = \overset{I}{T}\bigl(S(h')\bigr) \overset{I}{T}(h'') = \overset{I}{T}\bigl( S(h')h'' \bigr) = \varepsilon(h) \overset{I}{T}(1) = \varepsilon(h) \mathbb{I}_{\dim(I)}. 
\end{split}
\]
Note that
\begin{equation}\label{antipodeT}
S(\overset{I}{T}) =\: ^t \overset{\:I^*}{T}
\end{equation}
where $^t$ is the transpose. Indeed, let $(v_i)$ be a basis of $I$ and $(v^j)$ be the dual basis. Denote $v^*_j = v^j$. Then by \eqref{repX}, we have $\langle x v^j, v_i \rangle = \langle x v^*_j, v_i \rangle = \overset{I^*}{x}{_j^k} \langle v^*_k, v_i \rangle = \overset{I^*}{x}{^i_j}$. But using \eqref{actionDual} we also get $\langle x v^j, v_i \rangle = \langle v^j, S(x)v_i \rangle = \overset{I}{S(x)}{_i^k} \langle v^j, v_k \rangle = \overset{I}{S(x)}{_i^j}$. This shows that $\overset{I^*}{x} = {^t}\overset{I}{S(x)}$ as desired.
 Recall the well-known exchange relation
\begin{equation}\label{FRT}
\overset{IJ}{R}_{12} \overset{I}{T_1} \overset{J}{T_2} = \overset{J}{T_2} \overset{I}{T_1} \overset{IJ}{R}_{12},
\end{equation}
which is simply due to \eqref{dualHopf}, \eqref{braidingTwist} below and \eqref{naturaliteT}:
\[ \overset{IJ}{R} \,  \overset{I}{T_1} \, \overset{J}{T_2} = P_{IJ} \, c_{I,J} \overset{I \otimes J}{T} = P_{IJ} \overset{J \otimes I}{T} c_{I,J} = P_{IJ} \, \overset{J}{T_1} \, \overset{I}{T_2} \, P_{IJ} \, \overset{IJ}{R} = \overset{J}{T_2} \, \overset{I}{T_1} \, \overset{IJ}{R} \]
where $P_{IJ}$ is the flip tensor $P_{IJ}(x \otimes y) = y \otimes x$ or $(P_{IJ})^{ac}_{bd} = \delta^a_d \delta^c_b$. As before, if the representations $I$ and $J$ are fixed and arbitrary, we will simply write 
\[ T_{12} = T_1T_2 \:\:\:\: \text{ and } \:\:\:\: R_{12}T_1T_2 = T_2T_1R_{12}. \]

\begin{remark}\label{remarqueGensT}
Exactly as for the matrices $\overset{I}{L}{^{(\pm)}}$, the set of generators $\overset{I}{T}{^i_j}$ is not at all minimal. Due to the fusion relation of \eqref{dualHopf}, we have the same comments that in Remark \ref{remarqueGensL} (see section \ref{sectionMatrixCoeffUq}, where this is discussed in detail for $H = \bar U_q(\mathfrak{sl}_2)$).
\finEx
\end{remark}

\indent 
We denote by $\mathcal{Z}(H)$ the subalgebra of central elements of $H$ and by $\mathrm{SLF}(H) \subset \mathcal{O}(H)$ the subalgebra of symmetric linear forms on $H$:
\begin{equation*}
\mathrm{SLF}(H) = \left\{ \varphi \in H^* \, \left| \: \forall \, x,y \in H, \:\:\: \varphi(xy) = \varphi(yx) \right. \right\}.
\end{equation*}
$\mathrm{SLF}(H)$ is a subalgebra because the coproduct is a morphism of algebras:
\[ \varphi \psi(xy) = \varphi\bigl( (xy)' \bigr) \psi\bigl( (xy)'' \bigr) = \varphi(x' y') \psi(x'' y'') = \varphi(y' x') \psi(y'' x'') = \varphi\bigl( (yx)' \bigr) \psi\bigl( (yx)'' \bigr) = \varphi \psi(yx). \]
Consider the following map, called Drinfeld morphism:
\begin{equation}\label{morphismeDrinfeld}
\fonc{\mathcal{D}}{H^*}{H}{\varphi}{(\varphi \otimes \text{id})\bigl( (g \otimes 1) RR'\bigr) = \varphi(ga_ib_j)b_ia_j}
\end{equation}
where $g$ is the pivotal element \eqref{pivot}.
\begin{lemma}\label{propMorDrinfeld}
If $H$ is factorizable, $\mathcal{D} : H^* \to H$ is an isomorphism of vector spaces and it provides an isomorphism of algebras 
\[ \mathcal{D} : \mathrm{SLF}(H) \overset{\sim}{\longrightarrow} \mathcal{Z}(H). \]
In particular, $\mathrm{SLF}(H) $ is a commutative algebra.
\end{lemma}
\begin{proof}
$\mathcal{D}$ is an isomorphism by assumption since $g$ is invertible. Consider the following actions of $H$ on itself and on $H^*$:
\begin{equation*}
h \diamond a = h'' a S^{-1}(h'), \:\:\:\:\:\:\: h \diamond \varphi = \varphi\bigl( S(h'') ? h' \bigr). 
\end{equation*}
It is easy to see that
\[ a \in \mathcal{Z}(H) \:\, \iff \:\, \forall \, h \in H, \:\: h \diamond a = \varepsilon(h)a, \:\:\:\:\:\:\:\: \varphi \in \mathrm{SLF}(H) \:\, \iff \:\, \forall \, h \in H, \:\: h \diamond \varphi = \varepsilon(h)\varphi. \]
Moreover, using \eqref{REtCoproduit}, we get that $\mathcal{D}$ intertwines these actions:
\begin{align*}
\mathcal{D}\bigl( h \diamond \varphi \bigr) &=\mathcal{D}\bigl( \varphi( S(h'') ? h' )\bigr) = \varphi\bigl(S(h'') g a_i b_j h' \bigr)b_ia_j = \varphi\bigl(S(h^{(4)}) g a_i b_j h''' \bigr) b_i a_j h'' S^{-1}(h')\\
&= \varphi\bigl(S(h^{(4)}) g a_i h'' b_j \bigr) b_i h''' a_j S^{-1}(h') = \varphi\bigl(S(h^{(4)}) g h''' a_i b_j \bigr) h'' b_i a_j S^{-1}(h')\\
& = \varphi\bigl(g a_i b_j \bigr) h'' b_i a_j S^{-1}(h') = h \diamond \mathcal{D}(\varphi)
\end{align*}
Hence, $\mathcal{D}$ brings symmetric linear forms to central elements. To show that it is a morphism of algebras, we use \eqref{dualHopfDef}, \eqref{pivot} and \eqref{deltaR}:
\begin{align*}
\mathcal{D}(\varphi \psi) &= \varphi\psi(g a_i b_j) b_i a_j = \varphi\bigl( g a_i' b_j' \bigr) \psi\bigl( g a_i'' b_j'' \bigr) b_i a_j = \varphi(g a_i b_l) \psi(g a_k b_j ) b_i b_k a_j a_l\\
&=  \varphi(g a_i b_l) b_i \mathcal{D}(\psi) a_l = \varphi(g a_i b_l) b_i a_l \mathcal{D}(\psi) = \mathcal{D}(\varphi)\mathcal{D}(\psi)
\end{align*}
since $\mathcal{D}(\psi)$ is central.
\end{proof}

\smallskip

\indent Recall that a {\em left integral} (resp. {\em right integral}) is a non-zero linear form $\mu^l$ (resp. $\mu^r$) on $H$ which satisfies:
\begin{equation}\label{integrale}
\forall \, x \in H, \:\:\: (\mathrm{id} \otimes \mu^l) \circ \Delta(x) = \mu^l(x)1 \:\:\: \left(\text{resp. } (\mu^r \otimes \mathrm{id}) \circ \Delta(x) = \mu^r(x)1\right). 
\end{equation}
Since $H$ is finite dimensional, this is equivalent to:
\begin{equation}\label{integrale2}
\forall\, \psi \in \mathcal{O}(H), \:\:\: \psi\mu^l = \varepsilon(\psi)\mu^l \:\:\: \left(\text{resp. } \mu^r\psi = \varepsilon(\psi)\mu^r\right).
\end{equation}
It is well-known that left and right integrals always exist if $H$ is finite dimensional. Moreover, they are unique up to scalar. We fix $\mu^l$. Then $\mu^l \circ S^{-1}$ is a right integral, and we choose 
\begin{equation}\label{muR}
\mu^r = \mu^l \circ S^{-1}.
\end{equation}

Integrals will be important for our purposes because they are related to the ribbon element, as explained in the following proposition (the points 2. and 3. are well-known thanks to results of Radford, see \textit{e.g.} \cite{radfordBook}).

\begin{proposition}\label{propVIntegrale} 
Let $\varphi_v, \varphi_{v^{-1}} \in H^*$ defined by:
$$ \varphi_v = \mu^l(v^{-1})^{-1}\mu^l\bigl( g^{-1}v^{-1}\,? \bigr), \:\:\:\:\:\:\: \varphi_{v^{-1}} = \mu^l(v)^{-1}\mu^l\bigl(g^{-1}v\,?\bigr). $$
Then:
$$ \mathcal{D}(\varphi_v) = v, \:\:\:\:\:\:\: \mathcal{D}(\varphi_{v^{-1}}) = v^{-1}. $$
It follows that:
\begin{align}
&1.\:\:\: \varphi_{v}, \varphi_{v^{-1}} \in \SLF(H),\nonumber \\
&2. \:\:\: \forall\, x,y \in H, \:\:\: \mu^r(gxy) = \mu^r(gyx), \:\:\: \mu^l(g^{-1}xy) = \mu^l(g^{-1}yx), \label{integraleShifteSLF}\\
&3. \:\:\: \forall\, x,y \in H, \:\:\: \mu^r(xy) = \mu^r(S^2(y)x), \:\:\: \mu^l(xy) = \mu^l(S^{-2}(y)x). \label{quasiCyclic}
\end{align}
\end{proposition}
\begin{proof}
Consider the following computation, where we use (\ref{ribbon}) and (\ref{integrale}):
\begin{align*}
\mathcal{D}\!\left(\mu^l\!\left( g^{-1}v^{-1} \,? \right)\right) &= \left\langle \mu^l\!\left( g^{-1}v^{-1} \,? \right) \otimes \mathrm{id}, \, (g \otimes 1) RR'\right\rangle = \left\langle \mu^l\!\left( g^{-1}v^{-1} \,? \right) \otimes \mathrm{id}, \, gv(v^{-1})'' \otimes v(v^{-1})' \right\rangle \\
&=  \mu^l\!\left((v^{-1})''\right) v(v^{-1})' = \mu^l(v^{-1})v.
\end{align*}
Since $H$ is factorizable, the map $\mathcal{D}$ is an isomorphism of vector spaces. The left integral $\mu^l$ is non-zero, so $\mu^l\!\left( g^{-1}v^{-1} ? \right)$ is non-zero either. Since $\mathcal{D}$ is an isomorphism, it follows that $\mathcal{D}\!\left(\mu^l\!\left( g^{-1}v^{-1} \,? \right)\right) = \mu^l(v^{-1})v \neq 0$, and thus $\mu^l(v^{-1}) \neq 0$. Hence the formula for $\varphi_v$ is well defined. Moreover, we have the restriction $\mathcal{D} : \SLF(H) \overset{\sim}{\rightarrow} \mathcal{Z}(H)$, so since $v \in \mathcal{Z}(H)$, we get that $\varphi_v \in \SLF(H)$. This allows us to deduce the properties stated about $\mu^l$. Using (\ref{muR}), we obtain the properties 1, 2 and 3 for $\mu^r$. We can now proceed with the computation for $\varphi_{v^{-1}}$:
\begin{align*}
\mathcal{D}\!\left( \mu^l\!\left(g^{-1}v\,?\right) \right) &= \left\langle \mu^l\!\left(g^{-1}v\,?\right) \otimes \mathrm{id}, \, (g \otimes 1) RR'\right\rangle = \mu^l\!\left( va_ib_j \right) b_i a_j\\
&= \mu^r\!\left( v S(b_j)S(a_i) \right) b_i a_j = \mu^r\!\left( vS(a_i)S^{-1}(b_j) \right)b_ia_j \\
&= \left\langle \mu^r \otimes \mathrm{id}, \, (v \otimes 1) (R'R)^{-1} \right\rangle = \mu^r\!\left( v' \right)v''v^{-1} = \mu^r(v)v^{-1} = \mu^l(v)v^{-1}
\end{align*}
where we used (\ref{muR}), the property 3 previously shown and (\ref{ribbon}). We conclude as before.
\end{proof}

\noindent Let us record an useful formula:
\begin{equation}\label{integraleShifte}
\forall \, h \in H, \: \forall \, \varphi \in \mathcal{O}(H), \:\:\: \mu^r(h?) \varphi = \mu^r\!\left(h'?\right) \varphi\!\left(S^{-1}(h'')\right).
\end{equation}
Indeed, thanks to \eqref{integrale}
\begin{align*}
\left\langle \mu^r(h\,?)\psi, \, x \right\rangle &= \mu^r(hx')\psi(x'') = \mu^r(h'x')\psi\!\left(S^{-1}(h''')h''x''\right) = \psi\!\left(S^{-1}(h''')\mu^r(h'x')h''x''\right)\\
&= \psi\!\left(S^{-1}(h'')\,\mu^r\!\left((h'x)'\right)(h'x)''\right) = \mu^r(h'x)\,\psi\!\left(S^{-1}(h'')\right).
\end{align*}
Finally, it holds (see Lemma \ref{lemmaUnimodUnibalance})
\begin{equation}\label{muLmuRgCarre}
\mu^l = \mu^r(g^2 ?).
\end{equation}

\section{Heinsenberg double of $\mathcal{O}(H)$}\label{heisenbergDouble}
\indent The Heisenberg double of $\mathcal{O}(H)$, denoted by $\mathcal{H}(\mathcal{O}(H))$, is the vector space $\mathcal{O}(H) \otimes H$ endowed with the product
\[ (\varphi \otimes h)(\psi \otimes k) = \varphi \, \psi(?h') \otimes h''k \]
where $\psi(?y) \in \mathcal{O}(H)$ is defined by $\psi(?y)(x) = \psi(xy)$ (for $x \in H$) and $\varphi \, \psi(?y)$ is the product of $\varphi$ and $\psi(?y)$ in $\mathcal{O}(H)$: $\bigl\langle \varphi \, \psi(?y), x \bigr\rangle = \varphi(x')\psi(x''y)$. See \textit{e.g.} \cite[4.1.10]{Mon}. Let us identify $\psi \otimes 1 \in \mathcal{H}(\mathcal{O}(H))$ with $\psi \in \mathcal{O}(H)$ and $\varepsilon \otimes h \in \mathcal{H}(\mathcal{O}(H))$ with $h \in H$. Then we have $\varphi \otimes h = (\varphi \otimes 1)(\varepsilon \otimes h) = \varphi h$ and the structure of algebra on $\mathcal{H}(\mathcal{O}(H))$ is entirely defined by the following exchange relation:
\begin{equation}\label{relDefHeisenberg}
h\psi = \psi(?h')h''.
\end{equation}
There is a faithful representation $\triangleright$ of $\mathcal{H}(\mathcal{O}(H))$ on $\mathcal{O}(H)$ (see \cite[Lem. 9.4.2]{Mon}) defined by
\begin{equation}\label{repHO}
\psi \triangleright \varphi = \psi\varphi, \:\:\:\:\: h \triangleright \varphi = \varphi(?h).
\end{equation}
Hence we have an injective morphism $\rho : \mathcal{H}(\mathcal{O}(H)) \to \End_{\mathbb{C}}(H^*)$; by equality of the dimensions, it follows that 
\begin{equation}\label{HeisenbergMat}
\mathcal{H}(\mathcal{O}(H)) \cong \End_{\mathbb{C}}(H^*). 
\end{equation}
In other words, $\mathcal{H}(\mathcal{O}(H))$ is isomorphic to a matrix algebra. In particular, the elements of $\mathcal{H}(\mathcal{O}(H))$ can be defined by their action on $\mathcal{O}(H)$ under $\triangleright$. 

\smallskip
In terms of matrix coefficients, $\mathcal{H}(\mathcal{O}(H))$ is generated as an algebra by $\overset{I}{T}{^i_j}$ and $\overset{I}{L}{^{(\pm)}}^i_j$ for $I$ running in the set of finite dimensional $H$-modules. This set of generators is of course not at all minimal and in practice we can restrict to well-chosen representations, as explained in Remarks \ref{remarqueGensL} and \ref{remarqueGensT}. With these generators, the exchange relation is
\begin{equation}\label{echangeHeisenberg}
\overset{I}{L}\,\!^{(\pm)}_1 \overset{J}{T}_2 = \overset{J}{T}_2 \overset{I}{L}\,\!^{(\pm)}_1 \overset{IJ}{R}\,\!^{(\pm)}_{12}
\end{equation}
Indeed, using \eqref{L}, \eqref{relDefHeisenberg}, \eqref{deltaR} and obvious commutation relations:
\begin{align*}
\overset{I}{L}{^{(\pm)}_1} \overset{J}{T}_2 &= (\overset{I}{a_i^{(\pm)}})_1 \, b_i^{(\pm)} \overset{J}{T}_2 = (\overset{I}{a_i^{(\pm)}})_1 \, \overset{J}{T}(?b_i^{(\pm)}{'})_2 \, b_i^{(\pm)}{''} = \bigl(\overset{I}{a_i^{(\pm)}}\overset{I}{a_j^{(\pm)}}\bigr)_1 \, \overset{J}{T}(?b_j^{(\pm)})_2 \, b_i^{(\pm)}\\
& = (\overset{I}{a_i^{(\pm)}})_1 \, (\overset{I}{a_j^{(\pm)}})_1 \,  \overset{J}{T}_2 \, (\overset{J}{b_j^{(\pm)}})_2 \, b_i^{(\pm)} = \overset{J}{T}_2 \, (\overset{I}{a_i^{(\pm)}})_1 \, b_i^{(\pm)} \, (\overset{I}{a_j^{(\pm)}})_1 \, (\overset{J}{b_j^{(\pm)}})_2 = \overset{J}{T}_2 \overset{I}{L}{^{(\pm)}_1} \overset{IJ}{R}{^{(\pm)}_{12}}
\end{align*}
where $R^{(\pm)} = a_i^{(\pm)} \otimes b_i^{(\pm)}$. Similarly, the representation $\triangleright$ is
\begin{equation}\label{repTriangleHeisenberg}
\overset{I}{T_1}\triangleright \overset{J}{T}_2 = \overset{I \otimes J}{T}\!\!\!_{12}, \:\:\:\:\:\:\: \overset{I}{L}\,\!^{(\pm)}_1 \triangleright \overset{J}{T}_2 = (\overset{I}{a_i^{(\pm)}})_1 \, b_i^{(\pm)} \triangleright \overset{J}{T}_2 = (\overset{I}{a_i^{(\pm)}})_1 \overset{J}{T}_2 (\overset{I}{b_i^{(\pm)}})_2 = \overset{J}{T}_2 \overset{IJ}{R}\,\!^{(\pm)}_{12}.
\end{equation}


\smallskip

\indent For $h \in H$, let $\widetilde{h}\in \mathcal{H}(\mathcal{O}(H))$ be defined by 
\begin{equation}\label{operateursTilde}
\widetilde{h} \triangleright \varphi = \varphi(S^{-1}(h)?).
\end{equation}
It is easy to see that
\begin{equation}\label{echangeTilde}
\forall\, g \in H, \: \forall \, \psi \in \mathcal{O}(H), \:\:\:\: \widetilde{g}\widetilde{h} = \widetilde{gh}, \:\:\: g\widetilde{h} = \widetilde{h}g, \:\:\: \widetilde{h}\psi = \psi\!\left(S^{-1}(h'')?\right)\widetilde{h'}.
\end{equation}
Applying this to the matrices $\overset{I}{L}{^{(\pm)}}$ of generators of $H$, we define
\begin{equation}\label{defLTilde}
\overset{I}{\widetilde{L}}{^{(+)}} = \overset{I}{a_i} \, \widetilde{b_i}, \:\:\: \overset{I}{\widetilde{L}}{^{(-)}} = \overset{I}{S^{-1}(b_i)} \, \widetilde{a_i} \in \mathrm{Mat}_{\dim(I)}\bigl(\mathcal{H}(\mathcal{O}(H))\bigr)
\end{equation}
or equivalently 
\[ \overset{I}{\widetilde{L}}{^{(\pm)}_1} \triangleright \overset{J}{T}_2 = \overset{IJ}{R}\,\!^{(\pm)-1}_{12} \overset{J}{T}_2. \]
Using the standard properties of the $R$-matrix, it is not difficult to show the following relations:
\begin{equation}\label{LTilde}
\begin{array}{l}
\overset{I}{\widetilde{L}}{^{(\epsilon)}_1} \overset{J}{\widetilde{L}}{^{(\epsilon)}_2} = \overset{I\otimes J}{\widetilde{L}}\!\!{^{(\epsilon)}_{12}}, \:\:\:\:\:\:\:\:\:\:\:\:\:\:\:\:\: \overset{I}{\widetilde{L}}{^{(\epsilon)}_1}\overset{J}{L}{^{(\sigma)}_2} =  \overset{J}{L}{^{(\sigma)}_2}\overset{I}{\widetilde{L}}{^{(\epsilon)}_1}, \:\:\:\:\:\:\:\:\:\:\:\:\:\:\:\:\:  \overset{IJ}{R}{^{(\epsilon)}_{12}} \overset{I}{\widetilde{L}}{^{(\epsilon)}_1}\overset{J}{T}_2 = \overset{J}{T}_2\overset{I}{\widetilde{L}}{^{(\epsilon)}_1},\\
\overset{IJ}{R}\,\!^{(\epsilon)}_{12} \overset{I}{\widetilde{L}} \,\!^{(\epsilon)}_1 \overset{J}{\widetilde{L}} \,\!^{(\sigma)}_2 = \overset{J}{\widetilde{L}} \,\!^{(\sigma)}_2 \overset{I}{\widetilde{L}} \,\!^{(\epsilon)}_1 \overset{IJ}{R}\,\!^{(\epsilon)}_{12} \:\:\:\: \forall\, \epsilon, \sigma \in \{\pm\}, \:\:\:\:\:\:\:\:\:  \overset{IJ}{R}\,\!^{(\epsilon)}_{12} \overset{I}{\widetilde{L}} \,\!^{(\sigma)}_1 \overset{J}{\widetilde{L}} \,\!^{(\sigma)}_2 = \overset{J}{\widetilde{L}} \,\!^{(\sigma)}_2 \overset{I}{\widetilde{L}} \,\!^{(\sigma)}_1 \overset{IJ}{R}\,\!^{(\epsilon)}_{12} \:\:\:\: \forall\, \epsilon, \sigma \in \{\pm\}.
\end{array}
\end{equation}
For instance, here is a proof of the fourth equality with $\epsilon = +, \sigma=-$:
\begin{align*}
\overset{IJ}{R}{^{(+)}_{12}} \overset{I}{\widetilde{L}}{^{(+)}_1} \overset{J}{\widetilde{L}}{^{(-)}_2} &= (\overset{I}{a_i})_1 \, (\overset{J}{b_i})_2 \, (\overset{I}{a_j})_1 \, \widetilde{b_j} \, \overset{J}{S^{-1}(b_k)}_2 \, \widetilde{a_k} =
\bigl(\overset{I}{a_i} \overset{I}{a_j}\bigr)_1 \, \bigl( \overset{J}{b_i} \overset{J}{S^{-1}(b_k)} \bigr)_2 \, \widetilde{b_j a_k}\\
&= \bigl(\overset{I}{a_j} \overset{I}{a_i}\bigr)_1 \, \bigl( \overset{J}{S^{-1}(b_k)} \overset{J}{b_i} \bigr)_2 \, \widetilde{a_k b_j} = \overset{J}{S^{-1}(b_k)}_2 \, \widetilde{a_k} \, (\overset{I}{a_j})_1 \, \widetilde{b_j} \, (\overset{I}{a_i})_1 \, (\overset{J}{b_i})_2 = \overset{J}{\widetilde{L}}{^{(-)}_2} \overset{I}{\widetilde{L}}{^{(+)}_1} \overset{IJ}{R}{^{(+)}_{12}}.
\end{align*}
\noindent We simply used a Yang-Baxter relation:
\[ a_i a_j \otimes b_i S^{-1}(b_k) \otimes b_j a_k = R_{12} R_{13} R_{32}^{-1} = R_{32}^{-1} R_{13} R_{12} = a_j a_i \otimes S^{-1}(b_k) b_i \otimes a_k b_j. \]

\section{Category $\mathrm{mod}_l(H)$, Reshetikhin--Turaev functor}\label{modH}
\indent Let $H$ be as above, with $R$-matrix $R = a_i \otimes b_i$, ribbon element $v$ and pivotal element $g = uv^{-1}$. We recall the ribbon structure of $\mathrm{mod}_l(H)$, the category of finite dimensional left $H$-modules. The notations and conventions are those of \cite[Chap. XIV]{kassel}. The objects of $\mathrm{mod}_l(H)$ will be denoted by $I,J...$ as previously. The unit object is $\mathbb{C}$; its $H$-module structure is given by $h \cdot z = \varepsilon(h)z$.

\smallskip

\indent The braiding $c_{I,J}$ and twist $\theta_I$ are defined by 
\begin{equation}\label{braidingTwist}
\fonc{c_{I,J}}{I \otimes J}{J \otimes I}{x \otimes y}{P_{IJ}\overset{IJ}{R}(x \otimes y) = b_iy \otimes a_ix} \:\:\:\:\:\:\:\:\:\:\:\: \fonc{\theta_I}{I}{I}{x}{v^{-1}x} 
\end{equation}
\noindent where $P_{IJ} : I \otimes J \to J \otimes I$ is the flip map $P_{IJ}(x \otimes y) = y \otimes x$. Note that $c_{I,J}^{-1}(y \otimes x) = S(a_i)x \otimes b_iy$.

\smallskip

\indent The dual $I^*$ of the $H$-module $I$ is defined by 
\begin{equation}\label{actionDual}
\forall \, h \in H, \:\:\: \forall \, \varphi \in I^*, \:\:\:\:\: h\varphi = \varphi\!\left(S(h)?\right).
\end{equation}
The duality morphisms are defined by
\[ \fonc{b_I}{\mathbb{C}}{I \otimes I^*}{1}{v_i \otimes v^i} \:\:\:\:\:\:\:\:\:\:\:\:  \fonc{d_I}{I^* \otimes I}{\mathbb{C}}{\varphi \otimes x}{\varphi(x)} \]
and 
\[ \fonc{b'_I}{\mathbb{C}}{I^* \otimes I}{1}{v^i \otimes g^{-1}v_i} \:\:\:\:\:\:\:\:\:\:\:\:  \fonc{d'_I}{I \otimes I^*}{\mathbb{C}}{x \otimes \varphi}{\varphi(gx)}. \]
\indent The name ``ribbon category'' comes from the well-known fact (\cite{RT}, see also \cite{KM} and \cite[XIV.5.1]{kassel}) that there is a tensor functor $F_{\mathrm{RT}} : \mathcal{RG}_H \to \mathrm{mod}_l(H)$, where $\mathcal{RG}_H$ is the category of $\mathrm{mod}_l(H)$-colored ribbon graphs. We call $F_{\mathrm{RT}}$ the Reshetikhin-Turaev functor; it takes the following values:

\begin{center}
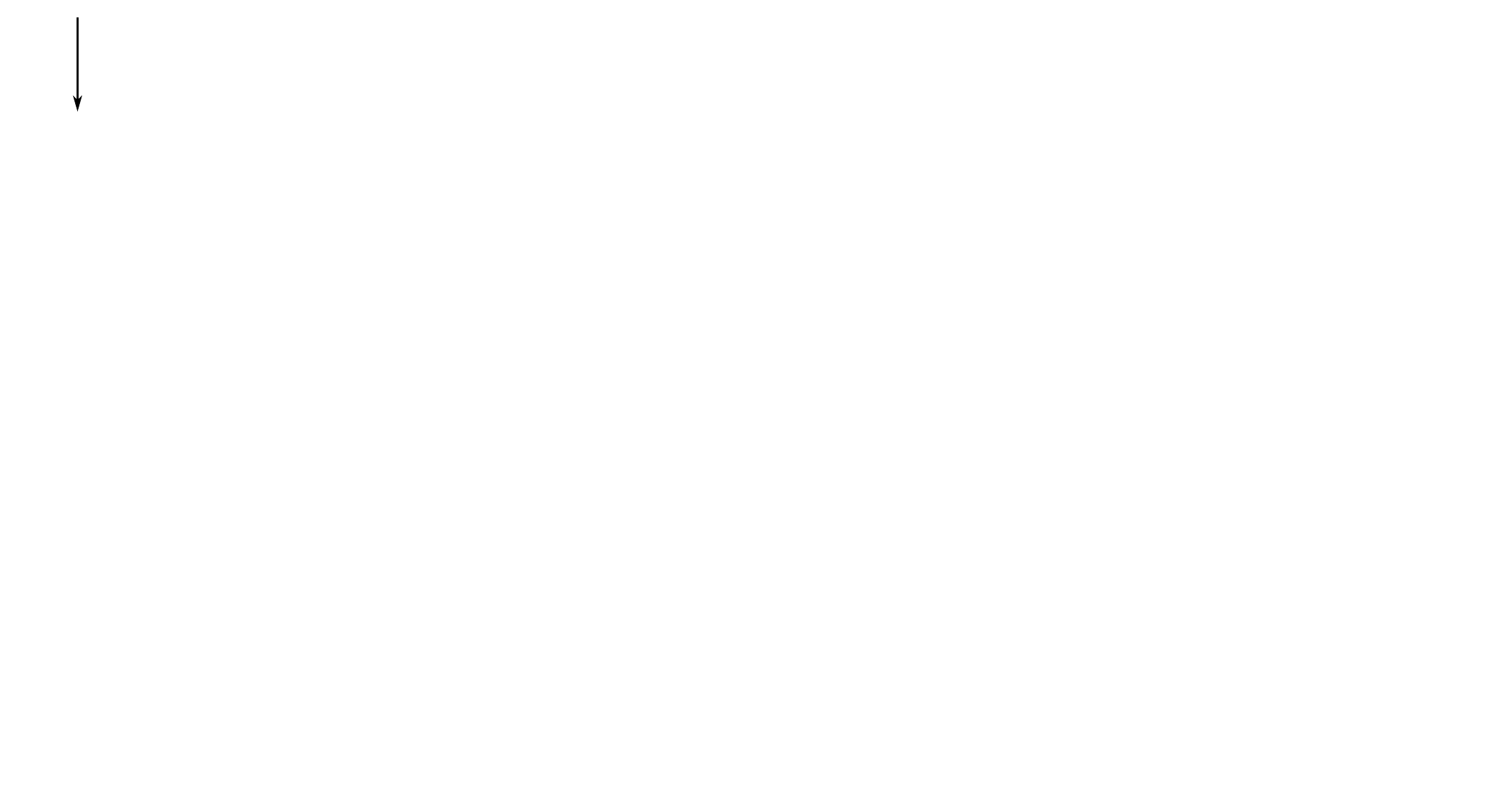
\end{center}

\noindent In the sequel we will identify a colored ribbon graph with its evaluation through $F_{\mathrm{RT}}$. Note that we read diagrams from bottom to top.

\smallskip

\indent We record some facts about duality in a ribbon category and more specifically in $\mathrm{mod}_l(H)$. First, the duality $^*$ is a contravariant functor. Indeed, recall that the transpose of $f : I \to J$ is $f^* : J^* \to I^*$, defined by 
\[ f^* = (\mathrm{id}_{I^*} \otimes d'_J) \circ (\mathrm{id}_{I^*} \otimes f \otimes \mathrm{id}_{J^*}) \circ (b'_I \otimes \mathrm{id}_{J^*}) = (d_J \otimes \mathrm{id}_{I^*}) \circ (\mathrm{id}_{J^*} \otimes f \otimes \mathrm{id}_{I^*}) \circ (\mathrm{id}_{J^*} \otimes b_I). \]
This is represented diagrammaticaly as follows:
\begin{center}
\begingroup%
  \makeatletter%
  \providecommand\color[2][]{%
    \errmessage{(Inkscape) Color is used for the text in Inkscape, but the package 'color.sty' is not loaded}%
    \renewcommand\color[2][]{}%
  }%
  \providecommand\transparent[1]{%
    \errmessage{(Inkscape) Transparency is used (non-zero) for the text in Inkscape, but the package 'transparent.sty' is not loaded}%
    \renewcommand\transparent[1]{}%
  }%
  \providecommand\rotatebox[2]{#2}%
  \newcommand*\fsize{\dimexpr\f@size pt\relax}%
  \newcommand*\lineheight[1]{\fontsize{\fsize}{#1\fsize}\selectfont}%
  \ifx\svgwidth\undefined%
    \setlength{\unitlength}{180.37270373bp}%
    \ifx\svgscale\undefined%
      \relax%
    \else%
      \setlength{\unitlength}{\unitlength * \real{\svgscale}}%
    \fi%
  \else%
    \setlength{\unitlength}{\svgwidth}%
  \fi%
  \global\let\svgwidth\undefined%
  \global\let\svgscale\undefined%
  \makeatother%
  \begin{picture}(1,0.33908107)%
    \lineheight{1}%
    \setlength\tabcolsep{0pt}%
    \put(0,0){\includegraphics[width=\unitlength,page=1]{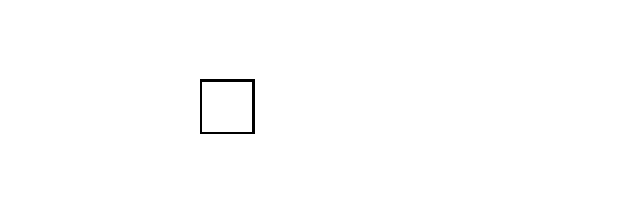}}%
    \put(0.33856131,0.15045779){\color[rgb]{0,0,0}\makebox(0,0)[lt]{\lineheight{1.25}\smash{\begin{tabular}[t]{l}$f$\end{tabular}}}}%
    \put(0.30308749,0.05839667){\color[rgb]{0,0,0}\makebox(0,0)[lt]{\lineheight{1.25}\smash{\begin{tabular}[t]{l}$I$\end{tabular}}}}%
    \put(0.38090322,0.23582409){\color[rgb]{0,0,0}\makebox(0,0)[lt]{\lineheight{1.25}\smash{\begin{tabular}[t]{l}$J$\end{tabular}}}}%
    \put(0,0){\includegraphics[width=\unitlength,page=2]{transpose.pdf}}%
    \put(0.78969154,0.23053261){\color[rgb]{0,0,0}\makebox(0,0)[lt]{\lineheight{1.25}\smash{\begin{tabular}[t]{l}$J$\end{tabular}}}}%
    \put(0.83639601,0.15142058){\color[rgb]{0,0,0}\makebox(0,0)[lt]{\lineheight{1.25}\smash{\begin{tabular}[t]{l}$f$\end{tabular}}}}%
    \put(0.88456691,0.06276635){\color[rgb]{0,0,0}\makebox(0,0)[lt]{\lineheight{1.25}\smash{\begin{tabular}[t]{l}$I$\end{tabular}}}}%
    \put(0.58889506,0.14849809){\color[rgb]{0,0,0}\makebox(0,0)[lt]{\lineheight{1.25}\smash{\begin{tabular}[t]{l}=\end{tabular}}}}%
    \put(-0.00552242,0.1476016){\color[rgb]{0,0,0}\makebox(0,0)[lt]{\lineheight{1.25}\smash{\begin{tabular}[t]{l}$f^* \:\: = $\end{tabular}}}}%
  \end{picture}%
\endgroup%

\end{center}
In $\mathrm{mod}_l(H)$, $f^*$ is simply the usual transpose: $f^*(\varphi) = \varphi \circ f$. The families $(d_I), (d'_I), (b_I), (b'_I)$ satisfy:
\begin{equation}\label{naturaCupCap}
\begin{array}{ll}
d_J \circ (\mathrm{id}_{J^*} \otimes f) = d_I \circ (f^* \otimes \mathrm{id}_I), & d'_J \circ (f \otimes \mathrm{id}_{J^*}) = d'_I \circ (\mathrm{id}_I \otimes f^*),\\
(f \otimes \mathrm{id}_{J^*}) \circ b_I = (\mathrm{id}_J \otimes f^*) \circ b_J, & (\mathrm{id}_{I^*} \otimes f) \circ b'_I = (f^* \otimes \mathrm{id}_J) \circ b'_J.
\end{array}
\end{equation}
This is represented diagrammaticaly as follows:
\begin{center}
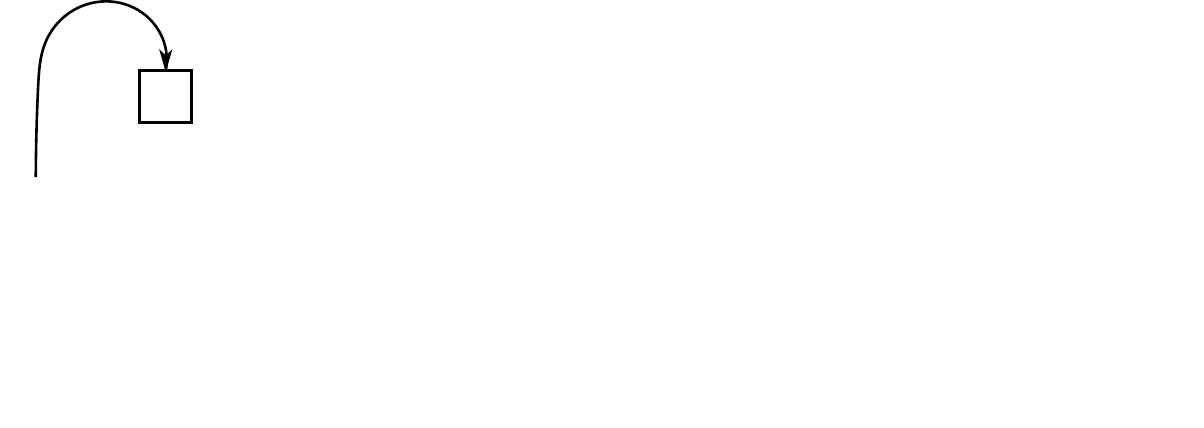
\end{center}

\indent Let
\begin{equation}\label{identificationBidual}
\begin{array}{crll}
e_I : & I^{**} & \overset{\sim}{\rightarrow} & I \\
                  & \langle ?, x \rangle &\mapsto & g^{-1}x
\end{array}
\end{equation}
be the identification with the bidual, where $\langle ?, ? \rangle = d_I : I^* \otimes I \to \mathbb{C}$ is the duality pairing. The morphism $e_I$ and its inverse can be written in terms of duality morphisms:
\begin{align*}
e_I &= (\mathrm{id}_I \otimes d'_{I^*}) \circ (b_I \otimes \mathrm{id}_{I^{**}}) = (d_{I^*} \otimes \mathrm{id}_I) \circ (\mathrm{id}_{I^{**}} \otimes b'_I),\\
e_I^{-1} &= (\mathrm{id}_{I^{**}} \otimes d_I) \circ (b'_{I^*} \otimes \mathrm{id}_I) = (d'_I \otimes \mathrm{id}_{I^{**}}) \circ (\mathrm{id}_I \otimes b_{I^*}).
\end{align*}
This is represented diagrammaticaly as follows:
\begin{center}
\begingroup%
  \makeatletter%
  \providecommand\color[2][]{%
    \errmessage{(Inkscape) Color is used for the text in Inkscape, but the package 'color.sty' is not loaded}%
    \renewcommand\color[2][]{}%
  }%
  \providecommand\transparent[1]{%
    \errmessage{(Inkscape) Transparency is used (non-zero) for the text in Inkscape, but the package 'transparent.sty' is not loaded}%
    \renewcommand\transparent[1]{}%
  }%
  \providecommand\rotatebox[2]{#2}%
  \newcommand*\fsize{\dimexpr\f@size pt\relax}%
  \newcommand*\lineheight[1]{\fontsize{\fsize}{#1\fsize}\selectfont}%
  \ifx\svgwidth\undefined%
    \setlength{\unitlength}{500.50116685bp}%
    \ifx\svgscale\undefined%
      \relax%
    \else%
      \setlength{\unitlength}{\unitlength * \real{\svgscale}}%
    \fi%
  \else%
    \setlength{\unitlength}{\svgwidth}%
  \fi%
  \global\let\svgwidth\undefined%
  \global\let\svgscale\undefined%
  \makeatother%
  \begin{picture}(1,0.15235131)%
    \lineheight{1}%
    \setlength\tabcolsep{0pt}%
    \put(0,0){\includegraphics[width=\unitlength,page=1]{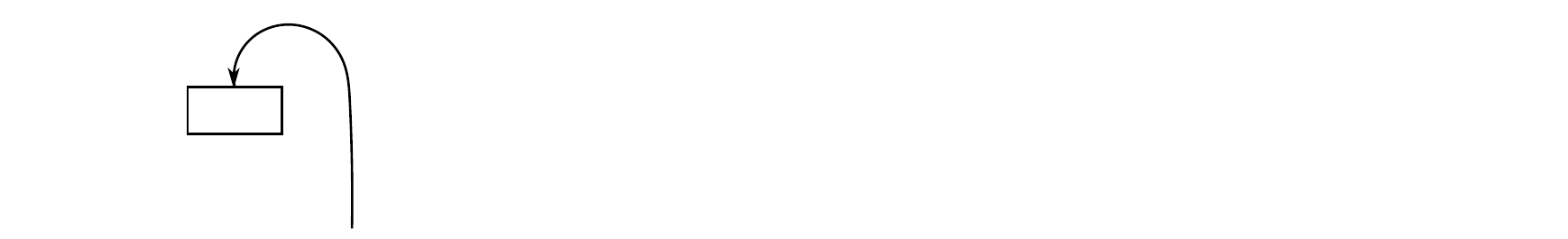}}%
    \put(0.26595352,0.08934378){\color[rgb]{0,0,0}\rotatebox{-180}{\makebox(0,0)[lt]{\lineheight{1.25}\smash{\begin{tabular}[t]{l}=\end{tabular}}}}}%
    \put(0.13423228,0.07499441){\color[rgb]{0,0,0}\makebox(0,0)[lt]{\lineheight{1.25}\smash{\begin{tabular}[t]{l}$\mathrm{id}_{I^*}$\end{tabular}}}}%
    \put(0,0){\includegraphics[width=\unitlength,page=2]{defE.pdf}}%
    \put(0.12576874,0.04306986){\color[rgb]{0,0,0}\makebox(0,0)[lt]{\lineheight{1.25}\smash{\begin{tabular}[t]{l}$I$\end{tabular}}}}%
    \put(0.16104372,0.10354132){\color[rgb]{0,0,0}\makebox(0,0)[lt]{\lineheight{1.25}\smash{\begin{tabular}[t]{l}$I^*$\end{tabular}}}}%
    \put(0,0){\includegraphics[width=\unitlength,page=3]{defE.pdf}}%
    \put(0.33500078,0.09794574){\color[rgb]{0,0,0}\makebox(0,0)[lt]{\lineheight{1.25}\smash{\begin{tabular}[t]{l}$I^*$\end{tabular}}}}%
    \put(0.37686569,0.03558456){\color[rgb]{0,0,0}\makebox(0,0)[lt]{\lineheight{1.25}\smash{\begin{tabular}[t]{l}$I$\end{tabular}}}}%
    \put(0.35297108,0.0671852){\color[rgb]{0,0,0}\makebox(0,0)[lt]{\lineheight{1.25}\smash{\begin{tabular}[t]{l}$\mathrm{id}_{I^*}$\end{tabular}}}}%
    \put(-0.00199019,0.07448421){\color[rgb]{0,0,0}\makebox(0,0)[lt]{\lineheight{1.25}\smash{\begin{tabular}[t]{l}$e_I \:\: = $\end{tabular}}}}%
    \put(0,0){\includegraphics[width=\unitlength,page=4]{defE.pdf}}%
    \put(0.81976793,0.08383206){\color[rgb]{0,0,0}\rotatebox{-180}{\makebox(0,0)[lt]{\lineheight{1.25}\smash{\begin{tabular}[t]{l}=\end{tabular}}}}}%
    \put(0.68804672,0.06948269){\color[rgb]{0,0,0}\makebox(0,0)[lt]{\lineheight{1.25}\smash{\begin{tabular}[t]{l}$\mathrm{id}_{I^*}$\end{tabular}}}}%
    \put(0,0){\includegraphics[width=\unitlength,page=5]{defE.pdf}}%
    \put(0.67643366,0.03755809){\color[rgb]{0,0,0}\makebox(0,0)[lt]{\lineheight{1.25}\smash{\begin{tabular}[t]{l}$I^*$\end{tabular}}}}%
    \put(0.71485815,0.09802961){\color[rgb]{0,0,0}\makebox(0,0)[lt]{\lineheight{1.25}\smash{\begin{tabular}[t]{l}$I$\end{tabular}}}}%
    \put(0,0){\includegraphics[width=\unitlength,page=6]{defE.pdf}}%
    \put(0.89731899,0.08991438){\color[rgb]{0,0,0}\makebox(0,0)[lt]{\lineheight{1.25}\smash{\begin{tabular}[t]{l}$I$\end{tabular}}}}%
    \put(0.93068009,0.03007284){\color[rgb]{0,0,0}\makebox(0,0)[lt]{\lineheight{1.25}\smash{\begin{tabular}[t]{l}$I^*$\end{tabular}}}}%
    \put(0.90678559,0.06167348){\color[rgb]{0,0,0}\makebox(0,0)[lt]{\lineheight{1.25}\smash{\begin{tabular}[t]{l}$\mathrm{id}_{I^*}$\end{tabular}}}}%
    \put(0.54152921,0.07577279){\color[rgb]{0,0,0}\makebox(0,0)[lt]{\lineheight{1.25}\smash{\begin{tabular}[t]{l}$e_I^{-1} \:\: = $\end{tabular}}}}%
  \end{picture}%
\endgroup%

\end{center}
The family of morphisms $(e_I)$ is natural:
\[ f \circ e_I = e_J \circ f^{**} \]
and satisfies the following property:
\begin{equation}\label{propCoherenceBidual}
e_{I^*} = (e_I^{-1})^*.
\end{equation}
This last equality is easy to see diagramatically:
\begin{center}
\begingroup%
  \makeatletter%
  \providecommand\color[2][]{%
    \errmessage{(Inkscape) Color is used for the text in Inkscape, but the package 'color.sty' is not loaded}%
    \renewcommand\color[2][]{}%
  }%
  \providecommand\transparent[1]{%
    \errmessage{(Inkscape) Transparency is used (non-zero) for the text in Inkscape, but the package 'transparent.sty' is not loaded}%
    \renewcommand\transparent[1]{}%
  }%
  \providecommand\rotatebox[2]{#2}%
  \newcommand*\fsize{\dimexpr\f@size pt\relax}%
  \newcommand*\lineheight[1]{\fontsize{\fsize}{#1\fsize}\selectfont}%
  \ifx\svgwidth\undefined%
    \setlength{\unitlength}{397.19234315bp}%
    \ifx\svgscale\undefined%
      \relax%
    \else%
      \setlength{\unitlength}{\unitlength * \real{\svgscale}}%
    \fi%
  \else%
    \setlength{\unitlength}{\svgwidth}%
  \fi%
  \global\let\svgwidth\undefined%
  \global\let\svgscale\undefined%
  \makeatother%
  \begin{picture}(1,0.29214058)%
    \lineheight{1}%
    \setlength\tabcolsep{0pt}%
    \put(0,0){\includegraphics[width=\unitlength,page=1]{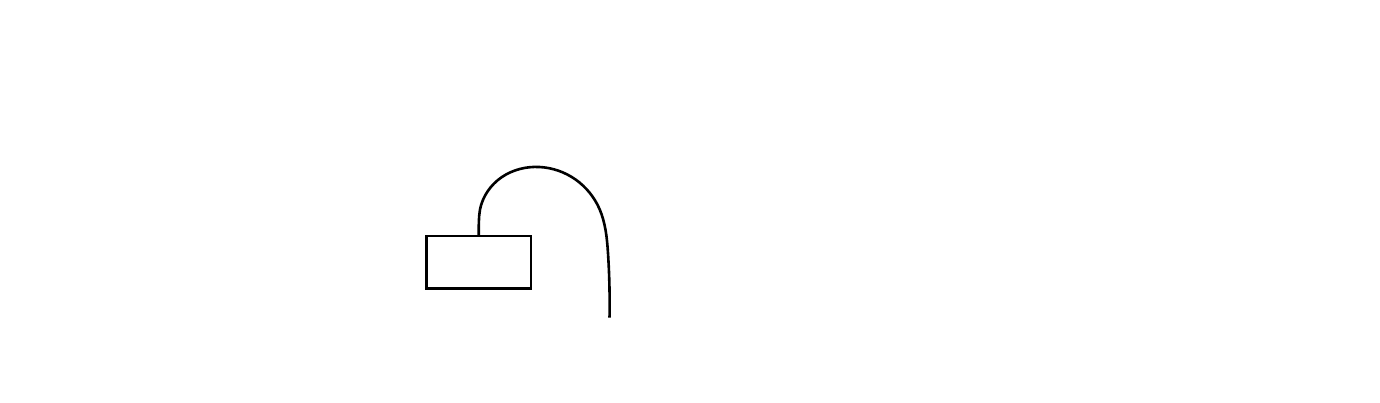}}%
    \put(0.32769222,0.0933049){\color[rgb]{0,0,0}\makebox(0,0)[lt]{\lineheight{1.25}\smash{\begin{tabular}[t]{l}$\mathrm{id}_{I^*}$\end{tabular}}}}%
    \put(0,0){\includegraphics[width=\unitlength,page=2]{propE.pdf}}%
    \put(0.31305858,0.05307678){\color[rgb]{0,0,0}\makebox(0,0)[lt]{\lineheight{1.25}\smash{\begin{tabular}[t]{l}$I^*$\end{tabular}}}}%
    \put(0.36147724,0.1292768){\color[rgb]{0,0,0}\makebox(0,0)[lt]{\lineheight{1.25}\smash{\begin{tabular}[t]{l}$I$\end{tabular}}}}%
    \put(-0.00144899,0.13124576){\color[rgb]{0,0,0}\makebox(0,0)[lt]{\lineheight{1.25}\smash{\begin{tabular}[t]{l}$(e_I^{-1})^* \:\: = $\end{tabular}}}}%
    \put(0,0){\includegraphics[width=\unitlength,page=3]{propE.pdf}}%
    \put(0.22808625,0.19471362){\color[rgb]{0,0,0}\makebox(0,0)[lt]{\lineheight{1.25}\smash{\begin{tabular}[t]{l}$\mathrm{id}_{I^{**}}$\end{tabular}}}}%
    \put(0,0){\includegraphics[width=\unitlength,page=4]{propE.pdf}}%
    \put(0.2030464,0.23432542){\color[rgb]{0,0,0}\makebox(0,0)[lt]{\lineheight{1.25}\smash{\begin{tabular}[t]{l}$I^{**}$\end{tabular}}}}%
    \put(0.59789082,0.14670617){\color[rgb]{0,0,0}\rotatebox{-180}{\makebox(0,0)[lt]{\lineheight{1.25}\smash{\begin{tabular}[t]{l}=\end{tabular}}}}}%
    \put(0,0){\includegraphics[width=\unitlength,page=5]{propE.pdf}}%
    \put(0.67647802,0.15698423){\color[rgb]{0,0,0}\makebox(0,0)[lt]{\lineheight{1.25}\smash{\begin{tabular}[t]{l}$I^{**}$\end{tabular}}}}%
    \put(0.73765087,0.07896435){\color[rgb]{0,0,0}\makebox(0,0)[lt]{\lineheight{1.25}\smash{\begin{tabular}[t]{l}$I^*$\end{tabular}}}}%
    \put(0.699207,0.11838726){\color[rgb]{0,0,0}\makebox(0,0)[lt]{\lineheight{1.25}\smash{\begin{tabular}[t]{l}$\mathrm{id}_{I^{**}}$\end{tabular}}}}%
    \put(0.85413822,0.12351902){\color[rgb]{0,0,0}\makebox(0,0)[lt]{\lineheight{1.25}\smash{\begin{tabular}[t]{l}$ = \:\: e_{I^*}$\end{tabular}}}}%
    \put(0,0){\includegraphics[width=\unitlength,page=6]{propE.pdf}}%
    \put(0.80011841,0.21359619){\color[rgb]{0,0,0}\makebox(0,0)[lt]{\lineheight{1.25}\smash{\begin{tabular}[t]{l}$\mathrm{id}_{I^*}$\end{tabular}}}}%
    \put(0.8339037,0.24956809){\color[rgb]{0,0,0}\makebox(0,0)[lt]{\lineheight{1.25}\smash{\begin{tabular}[t]{l}$I$\end{tabular}}}}%
    \put(0,0){\includegraphics[width=\unitlength,page=7]{propE.pdf}}%
  \end{picture}%
\endgroup%

\end{center}
Finally, note that:
\begin{equation*}\label{dualitePrime}
\begin{array}{ll}
b'_I = (\mathrm{id}_{I^*} \otimes e_I) \circ b_{I^*}, \:\:\: & d'_I = d_{I^*} \circ (e_I^{-1} \otimes \mathrm{id}_{I^*}),\\
b_I = b'_{I^*} \circ (e_I \otimes \mathrm{id}_{I^*}), \:\:\: & d_I = d'_{I^*} \circ (\mathrm{id}_{I^*} \otimes e_I^{-1}),
\end{array}
\end{equation*}
where the second line of equalities follows from the first line thanks to \eqref{naturaCupCap} and \eqref{propCoherenceBidual}. This is represented diagrammaticaly as follows:
\begin{equation}\label{dualitePrime}
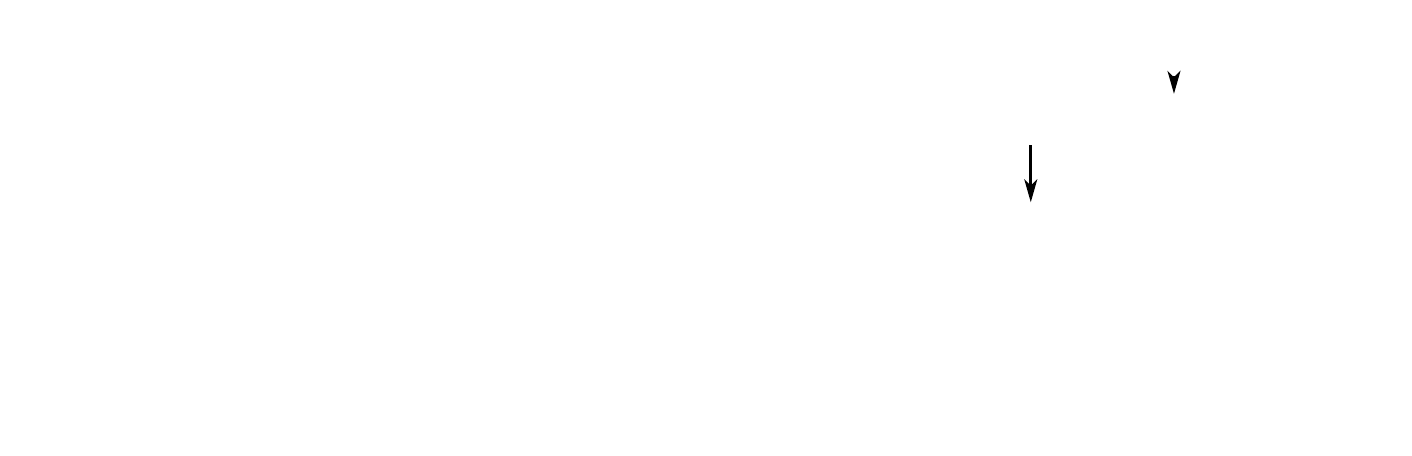
\end{equation}
In particular, we see that it is possible to define the morphisms $d'_I, b'_I$ using $b_I, d_I$ and $e_I$. We will use this remark to define the value of a negatively oriented strand going through an handle in the graphical calculus of Chapter \ref{chapitreGraphiqueSkein} (see \eqref{sensOppose}).

\chapter{The restricted quantum group $\bar U_q(\mathfrak{sl}_2)$}\label{chapitreUqSl2}

\indent The restricted quantum group $\bar U_q(\mathfrak{sl}_2)$ is an important example of a factorizable and ribbon\footnote{These two properties are an abuse of terminology due to the subtlety that the $R$-matrix belongs to an extension of of $\bar U_q(\mathfrak{sl}_2)$ by a square root $K^{1/2}$. But this square root does not appear in the double braiding $RR'$ nor in the ribbon element $v$, see section \ref{braidedExtension}.} finite dimensional Hopf algebra. This algebra is interesting because it is sufficiently simple and well-studied to allow us to carry explicit computations when we use it as a gauge algebra in the combinatorial quantization (see sections \ref{exempleUq} and \ref{calculGraphiqueSl2}). Moreover, its $R$-matrix satisfies the Jones skein relation when evaluated in the fundamental representation. We will use this property to construct representations of skein algebras, thanks to the Wilson loop map, in Chapter \ref{chapitreGraphiqueSkein}.

\smallskip

\indent Note that the aim of this chapter is to collect technical properties which will be used when applying the general constructions of the subsequent chapters to the (important) example of $\bar U_q(\mathfrak{sl}_2)$; the reader can skip this chapter at the first reading and just refer to it from the sections \ref{exempleUq}, \ref{calculGraphiqueSl2}, \ref{sectionRepSkeinAlg}, \ref{repSqSigma1Uq}.

\smallskip

\indent We first recall the definition of $\bar U_q = \bar U_q(\mathfrak{sl}_2)$ and explain some facts about its structure, its representation theory and its matrix coefficients. The material of sections \ref{sectionSimpleProj}, \ref{bimodEtCentre}, \ref{braidedExtension} comes from \cite{FGST} (it is also recalled in detail in \cite{ibanez}).

\smallskip

\indent Then we define the GTA basis of $\SLF(\bar U_q)$ (symmetric linear forms on $\bar U_q$), introduced in \cite{GT} and in \cite{arike}. This basis will be a crucial tool when studying explicitly the representations of $\mathrm{SL}_2(\mathbb{Z})$ and of $\mathcal{S}_q(\Sigma_{1,0})$ on $\SLF(\bar U_q)$ (sections \ref{SL2ZUq} and \ref{repSqSigma1Uq}). A key property of this basis is that its multiplication rules, determined in section \ref{multiplication}, are simple. I mention that such multiplication rules were already given in \cite{GT}, but I was not aware of the existence of their paper when preparing this work (in \cite{GT} they work in the space of $q$-characters $\mathrm{qCh}(\bar U_q)$, which is isomorphic as an algebra to $\mathrm{SLF}(\bar U_q)$ by the shift of the pivotal element; ; also note that they use a normalized version of $G_s$). It turns out that our proofs are different. In \cite{GT}, they use the fact that the multiplication in the canonical basis of $\mathcal{Z}(\bar U_q)$ is very simple. They first express the image of their basis of $\mathrm{qCh}(\bar U_q)$ (which is the GTA basis shifted by the inverse of the pivotal element) through the Radford mapping in the canonical basis of $\mathcal{Z}(\bar U_q)$. This gives a basis of $\mathcal{Z}(\bar U_q)$ called the Radford basis. Then they use the $\mathcal{S}$-transformation of the $\mathrm{SL}_2(\mathbb{Z})$ representation on $\mathcal{Z}(\bar U_q)$ to express the Drinfeld basis (which is the image of their basis of $\mathrm{qCh}(\bar U_q)$ by the Drinfeld map) in the Radford basis. This gives the multiplication rules in the Drinfeld basis. Since the Drinfeld map is an isomorphism of algebras between $\mathrm{qCh}(\bar U_q)$ and $\mathcal{Z}(\bar U_q)$, this gives also the multiplication rules in the GTA basis. Here we directly work in $\mathrm{SLF}(\bar U_q)$. We take advantage of the decomposition rules for tensor products to see that there are not many coefficients to determine, and then we compute these coefficients by using the evaluation on suitable elements of $\bar U_q$.

\smallskip

\indent Section \ref{sectionArike} is a digression on the link between the GTA basis and traces on projective $\bar U_q$-modules. We also compute the decomposition of $\mu^r(K^{p+1}?)$ (right integral shifted by the pivotal element) in the GTA basis (to be used in the proof of Theorem \ref{actionSL2ZArike}). A particular trace on projectives $\bar U_q$-modules is the modified trace computed in \cite{BBGe}, and we observe that $\mu^r(K^{p+1}?)$ is the symmetric linear form corresponding to the modified trace. This last fact is not specific to $\bar U_q$: it has been proved under general assumptions (which cover the case of $\bar U_q$) in \cite{BBG} that there is a unique (up to scalar) modified trace on the ideal of finite-dimensional projective $H$-modules and that the corresponding symmetric linear form is the right integral shifted by the pivotal element.

\smallskip

\indent The material of this chapter is mainly the content of \cite{F}. However, the present chapter contains more details and comments.

\section{Properties of $\bar U_q$}
\indent Let $q = e^{i\pi/p}$ be a primitive root of unity of order $2p$, with $p \geq 2$. Recall that $\bar U_q(\mathfrak{sl}_2)$, the {\em restricted quantum group} associated to $\mathfrak{sl}_2(\mathbb{C})$, is the $\mathbb{C}$-algebra generated by $E, F, K$ together with the relations 
\begin{equation*}
E^p=F^p=0, \:\:\: K^{2p}=1,\:\:\: KE=q^2EK,\:\:\: KF=q^{-2}FK, \:\:\: EF = FE + \frac{K-K^{-1}}{q-q^{-1}}.
\end{equation*}
It will be simply denoted by $\bar U_q$ in the sequel. It is a $2p^3$-dimensional Hopf algebra, with comultiplication $\Delta$, counit $\varepsilon$ and antipode $S$ given by the following formulas:
\begin{equation*}
\begin{array}{lll}
\Delta(E) = 1 \otimes E + E \otimes K, & \Delta(F) = F \otimes 1 + K^{-1} \otimes F, & \Delta(K) = K \otimes K, \\
\varepsilon(E) = 0, & \varepsilon(F) = 0, & \varepsilon(K) = 1,\\
S(E) = -EK^{-1}, & S(F) = -KF, & S(K) = K^{-1}.
\end{array}
\end{equation*}
\indent The monomials $E^mF^nK^{l}$ with $0 \leq m,n \leq p-1, \:0 \leq l \leq 2p-1$, form a basis of $\bar U_q$, usually referred as the PBW-basis. Recall the formula (see \cite[Prop. VII.1.3]{kassel}):
\begin{equation}\label{coproduitMonome}
\Delta(E^mF^nK^{l}) = \sum_{i=0}^m\sum_{j=0}^n q^{i(m-i) + j(n-j) - 2(m-i)(n-j)}{m \brack i}{n \brack j}E^{m-i}F^jK^{l+j-n} \otimes E^iF^{n-j}K^{l+m-i}
\end{equation}
where the $q$-binomial coefficients are defined by ${a \brack b} = \frac{[a]!}{[b]! [a-b]!}$ for $a \geq b$. \\
\indent Since $K$ is annihilated by the polynomial $X^{2p}-1$, which has simple roots over $\mathbb{C}$, the action of $K$ is diagonalizable on each $\bar U_q$-module, and the eigenvalues are $2p$-roots of unity.

\smallskip

\indent The elements $K$ and $K^{p+1}$ both satisfy the properties of a pivotal element, see \eqref{pivot}. In the sequel and as in \cite{FGST}, we always take $K^{p+1}$ for the pivotal element of $\bar U_q$:
\begin{equation}\label{pivotUq}
g = K^{p+1}.
\end{equation}

\indent Due to the Hopf algebra structure on $\bar U_q$, its category of modules is a monoidal category with duals. It is not braided (see \cite{KS}).

\subsection{Simple and projective $\bar U_q$-modules}\label{sectionSimpleProj}
\indent The finite dimensional representations of $\bar U_q$ are classified (\cite{suter} and \cite{GSTF}). Two types of modules are important for our purposes: the simple modules and the projective modules. As in \cite{FGST} (see also \cite{ibanez}), we denote the simple modules by $\mathcal{X}^{\alpha}(s)$, with $\alpha \in \{\pm\}, 1 \leq s \leq p$. $\mathcal{X}^+(2)$ is called the fundamental representation. The modules $\mathcal{X}^{\pm}(p)$ are simple and projective simultaneously. The other indecomposable projective modules are not simple. We denote them by $\mathcal{P}^{\alpha}(s)$ with $\alpha \in \{\pm\}, 1 \leq s \leq p-1$. 
\smallskip\\
\indent The module $\mathcal{X}^{\alpha}(s)$ admits a {\em canonical basis} $\left(v_i\right)_{0 \leq i \leq s-1}$ such that
\begin{equation}\label{BaseSimple}
Kv_i = \alpha q^{s-1-2i}v_i,\: Ev_0=0,\: Ev_i=\alpha[i][s-i]v_{i-1},\: Fv_i=v_{i+1},\: Fv_{s-1}=0.
\end{equation}
The module $\mathcal{P}^{\alpha}(s)$ admits a {\em standard basis} $\left(b_i, x_j, y_k, a_{l}\right)_{\substack{0 \leq i,l \leq s-1 \\ 0 \leq j,k \leq p-s-1}}$ such that

\begin{equation}\label{BaseProjectif}
\begin{array}{lll}
Kb_i=\alpha q^{s-1-2i}b_i, & Eb_i = \alpha[i][s-i]b_{i-1} + a_{i-1}, & Fb_i = b_{i+1}, \\
 & Eb_0 = x_{p-s-1}, & Fb_{s-1}=y_0,\\
Kx_j = -\alpha q^{p-s-1-2j}x_j, & Ex_j = -\alpha[j][p-s-j]x_{j-1}, &  Fx_j = x_{j+1},  \\
 & Ex_0=0, & Fx_{p-s-1}=a_0, \\
Ky_k = -\alpha q^{p-s-1-2k}y_k, & Ey_k = -\alpha[k][p-s-k]y_{k-1}, & Fy_k = y_{k+1}, \\
 & Ey_0 = a_{s-1}, & Fy_{p-s-1} = 0, \\
Ka_{l} = \alpha q^{s-1-2l}a_{l}, & Ea_{l} = \alpha[l][s-l]a_{l-1}, & Fa_{l} = a_{l+1}, \\
 & Ea_0 = 0, & Fa_{s-1}=0.
\end{array}
\end{equation}

Note that such a basis is not unique up to scalar since we can replace $b_i$ by $b_i + \lambda a_i$ (with $\lambda \in \mathbb{C})$ without changing the action. 
\\\indent In terms of composition factors, the structure of $\mathcal{P}^{\alpha}(s)$ can be schematically represented as follows (with the basis vectors corresponding to each factor and the action of $E$ and $F$):
\begin{equation}\label{figureProj}
\xymatrix{
 & \Top\left(\mathcal{P}^{\alpha}(s)\right) \cong \mathcal{X}^{\alpha}(s), (b_i)_{0 \leq i \leq s-1} \ar[ld]^E \ar[rd]^F \ar[dd]^E& \\
(x_j)_{0 \leq j \leq p-s-1}, \mathcal{X}^{-\alpha}(p-s) \!\!\!\!\!\!\!\!\! \ar[rd]^F & & \!\!\!\!\!\!\!\!\! \mathcal{X}^{-\alpha}(p-s) \ar[ld]^E, (y_k)_{0 \leq k \leq p-s-1}\\  
 & \Soc\left(\mathcal{P}^{\alpha}(s)\right) \cong \mathcal{X}^{\alpha}(s), (a_{l})_{0 \leq l \leq s-1} &
}
\end{equation}
If we need to emphasize the module in which we are working, we will use the following notations: $v_i^{\alpha}(s)$ for the canonical basis of $\mathcal{X}^{\alpha}(s)$ and $b_i^{\alpha}(s)$, $x_j^{\alpha}(s)$, $y_k^{\alpha}(s)$, $a_{l}^{\alpha}(s)$ for a  standard basis of $\mathcal{P}^{\alpha}(s)$ (these are the notations used in \cite{arike}).
\smallskip\\
\indent Let us recall the $\bar U_q$-morphisms between these modules. Observe that $\mathcal{X}^{\alpha}(s)$ is $\bar U_q$-generated by $v^{\alpha}_0(s)$ and $\mathcal{P}^{\alpha}(s)$ is $\bar U_q$-generated by $b^{\alpha}_0(s)$, so the images of these vectors suffice to define $\bar U_q$-morphisms. $\mathcal{X}^{\alpha}(s)$ is simple, so by Schur's lemma $\End_{\bar U_q}\left(\mathcal{X}^{\alpha}(s)\right) = \mathbb{C}\text{Id}$. Since 
$$\mathcal{X}^{\alpha}(s) \cong \Top\left(\mathcal{P}^{\alpha}(s)\right) \cong \Soc\left(\mathcal{P}^{\alpha}(s)\right)$$
there exist injection and projection maps defined by:
\begin{equation*}
\begin{array}{lcl}
\mathcal{X}^{\alpha}(s) & \hookrightarrow & \mathcal{P}^{\alpha}(s)\\
v_0^{\alpha}(s) & \mapsto & a_0^{\alpha}(s)
\end{array} \:\:\: \text{ and } \:\:\:\:\,
\begin{array}{lcl}
\mathcal{P}^{\alpha}(s) & \twoheadrightarrow & \mathcal{X}^{\alpha}(s) \\
b_0^{\alpha}(s) & \mapsto & v_0^{\alpha}(s).
\end{array}
\end{equation*}
We have $\End_{\bar U_q}\bigl(\mathcal{P}^{\alpha}(s)\bigr) = \mathbb{C}\text{Id} \oplus \mathbb{C}p^{\alpha}_s$ and $\Hom_{\bar U_q}\bigl(\mathcal{P}^{\alpha}(s), \mathcal{P}^{-\alpha}(p-s)\bigr) = \mathbb{C}P^{\alpha}_s \oplus \mathbb{C}\overline{P}^{\alpha}_s$, where:
\begin{equation}\label{morphismesP}
p^{\alpha}_s\!\left(b_0^{\alpha}(s)\right) = a_0^{\alpha}(s), \:\:\:\:\:\: P^{\alpha}_s\!\left(b_0^{\alpha}(s)\right) = x_0^{-\alpha}(p-s), \:\:\:\:\:\: \overline{P}^{\alpha}_s\!\left(b_0^{\alpha}(s)\right) = y_0^{-\alpha}(p-s).
\end{equation}
The other Hom-spaces involving only simple modules and indecomposable projective modules are null.

\subsection{Structure of the bimodule $_{\bar U_q}\!\left(\bar U_q\right)_{\bar U_q}$ and the center of $\bar U_q$}\label{bimodEtCentre}
\indent Recall that if $M$ is a left module (over any $k$-algebra $A$), then $M^* = \Hom_{\mathbb{C}}(M, k)$ is endowed with a {\em right} $A$-module structure, given by:
\begin{equation*}
\forall \, a \in A, \: \forall\, \varphi \in M^*, \:\: \varphi a = \varphi(a \cdot)
\end{equation*}
where $\cdot$ is the place of the variable. We denote by $R^*(M)$ the so-defined right module. Note that if we define $R^*(f)$ as the transpose of $f$, then $R^*$ becomes a contravariant functor. If $A$ is a Hopf algebra, one must be aware not to confuse $R^*(M)$ with the categorical dual $M^*$, which is a left module on which $A$ acts by:
$$ \forall \, a \in A, \: \forall\, \varphi \in M^*, \:\: a\varphi = \varphi(S(a) \cdot). $$
\begin{lemma}The right $\bar U_q$-module $R^*(\mathcal{X}^{\alpha}(s))$ admits a basis $\left(\bar v_i\right)_{0 \leq i \leq s-1}$ such that
\begin{equation*}
\begin{array}{lllll}
\bar v_i K = \alpha q^{1-s+2i}\bar v_i, & \bar v_i E = \alpha[i][s-i] \bar v_{i-1}, & \bar v_0 E = 0, & v_i F = \bar v_{i+1}, & \bar v_{s-1}F = 0. \\
\end{array}
\end{equation*}
The right $\bar U_q$-module $R^*(\mathcal{P}^{\alpha}(s))$ admits a basis $\left(\bar b_{i}, \bar x_j, \bar y_k, \bar a_l\right)_{\substack{0 \leq i,l \leq s-1 \\ 0 \leq j ,k\leq p-s-1}}$ such that
\begin{equation*}
\begin{array}{lll}
   \bar b_i K = \alpha q^{1-s+2i}\bar b_i, & \bar b_i E = \bar a_{i-1} + \alpha[i][s-i] \bar b_{i-1}, & \bar b_i F = \bar b_{i+1}, \\
      & \bar b_0 E = \bar x_{p-s-1}, & \bar b_{s-1}F = \bar y_0, \\
   \bar x_j K = -\alpha q^{-p+s+1+2j} \bar x_j, & \bar x_j E = -\alpha [j][p-s-j]\bar x_{j-1}, & \bar x_j F = \bar x_{j+1}, \\
      &  \bar x_{0} E = 0, & \bar x_{p-s-1} F = \bar a_0, \\
   \bar y_{k} K = -\alpha q^{-p+s+1+2k} \bar y_k, & \bar y_k E = -\alpha[k][p-s-k] \bar y_{k-1}, & \bar y_k F = \bar y_{k+1}, \\
      &  \bar y_0 E = \bar a_{s-1}, & \bar y_{p-s-1}F = 0, \\
   \bar a_{l} K = \alpha q^{1-s+2l}\bar a_{l}, & \bar a_{l} E = \alpha[l][s-l] \bar a_{l-1}, & \bar a_{l} F = \bar a_{l+1}, \\
      & \bar a_0 E = 0, & \bar a_{s-1}F = 0.
\end{array}
\end{equation*}
\end{lemma}
\noindent Such basis will be termed respectively the {\em canonical basis} and a {\em standard basis} in the sequel.
\begin{proof}
Let $(v^i)_{0 \leq i \leq s-1}$ be the basis dual to the canonical basis given in (\ref{BaseSimple}). Then $\bar v_i = v^{s-1-i}$ gives the desired result. Similarly, let $\left(b^i, x^j, y^k, a^{l}\right)_{\substack{0 \leq i,l \leq s-1 \\ 0 \leq j,k \leq p-s-1}}$ be the basis dual to a standard basis given in (\ref{BaseProjectif}). Then 
$$ \bar b_i = a^{s-1-i}, \:\: \bar x_j = y^{p-s-1-j}, \:\: \bar y_k = x^{p-s-1-k}, \:\: \bar a_{l} = b^{s-1-l} $$
gives the desired result.
\end{proof}

\indent We denote by $_{\bar U_q}\!\left(\bar U_q\right)_{\bar U_q}$ the regular bimodule, where the left and right actions are respectively the left and right multiplication of $\bar U_q$ on itself. Recall that a block of $_{\bar U_q}\!\left(\bar U_q\right)_{\bar U_q}$ is just an indecomposable two-sided ideal (see \cite[Section 55]{CR}). The block decomposition of $\bar U_q$ is (see \cite{FGST})
\begin{equation*}
_{\bar U_q}\!\left(\bar U_q\right)_{\bar U_q} = \bigoplus_{s=0}^{p} Q(s)
\end{equation*}
where the structure of each block $Q(s)$ as a left $\bar U_q$-module is: 
\begin{equation}\label{leftBlock}
\begin{array}{l}
Q(0) \cong p\mathcal{X}^{-}(p), \:\:\:\:\: Q(p) \cong p\mathcal{X}^{+}(p),\\
Q(s) \cong s\mathcal{P}^{+}(s) \oplus (p-s)\mathcal{P}^{-}(p-s) \: \text{ for } 1 \leq s \leq p-1
\end{array}
\end{equation}
and the structure of each block as a right $\bar U_q$-module is:
\begin{equation*}
\begin{array}{l}
Q(0) \cong pR^*\!\left(\mathcal{X}^{-}(p)\right), \:\:\:\:\: Q(p) \cong pR^*\!\left(\mathcal{X}^{+}(p)\right),\\
Q(s) \cong sR^*\!\left(\mathcal{P}^{+}(s)\right) \oplus (p-s)R^*\!\left(\mathcal{P}^{-}(p-s)\right) \: \text{ for } 1 \leq s \leq p-1.
\end{array}
\end{equation*}

The following proposition is a reformulation of \cite[Prop. 4.4.2]{FGST} (see also \cite[Th. II.1.4]{ibanez}). It will be used for the proof of Theorem \ref{MainResult}.
\begin{proposition}\label{baseBloc}
For $1 \leq s \leq p-1$, the block $Q(s)$ admits a basis $$\left(B^{++}_{ab}(s), X^{-+}_{cd}(s), Y^{-+}_{ef}(s), A^{++}_{gh}(s), B^{--}_{ij}(s), X^{+-}_{kl}(s), Y^{+-}_{mn}(s), A^{--}_{or}(s) \right)$$
with $0 \leq a,b,d,f,g,h,k,m \leq s-1, \:\: 0 \leq c, e, i,j,l,n,o,r \leq p-s-1$, such that 
\begin{enumerate}
\item $\forall\, 0 \leq j \leq s-1, \:\: \left(B^{++}_{ij}(s), X^{-+}_{kj}(s), Y^{-+}_{lj}(s), A^{++}_{mj}(s)\right)_{\substack{0 \leq i,m \leq s-1 \\ 0 \leq k,l \leq p-s-1}}$ is a standard basis of $\mathcal{P}^{+}(s)$ for the left action.
\item $\forall \, 0 \leq j \leq p-s-1, \:\: \left(B^{--}_{ij}(s), X^{+-}_{kj}(s), Y^{+-}_{lj}(s), A^{--}_{mj}(s)\right)_{\substack{0 \leq k,l \leq s-1 \\ 0 \leq i,m \leq p-s-1}}$ is a standard basis of $\mathcal{P}^{-}(p-s)$ for the left action.
\item $\forall\, 0 \leq i \leq s-1, \:\: \left(B^{++}_{ij}(s), X^{+-}_{ik}(s), Y^{+-}_{il}(s), A^{++}_{im}(s)\right)_{\substack{0 \leq j,m \leq s-1 \\ 0 \leq k,l \leq p-s-1}}$ is a standard basis of $R^*\left(\mathcal{P}^{+}(s)\right)$ for the right action.
\item $\forall\, 0 \leq i \leq p-s-1, \:\: \left(B^{--}_{ij}(s), X^{-+}_{ik}(s), Y^{-+}_{il}(s), A^{--}_{im}(s)\right)_{\substack{0 \leq k,l \leq s-1 \\ 0 \leq j,m \leq p-s-1}}$ is a standard basis of $R^*\left(\mathcal{P}^{-}(p-s)\right)$ for the right action.
\end{enumerate}
The block $Q(0)$ admits a basis $\left(A^{--}_{ij}(0)\right)_{0 \leq i,j \leq p-1}$ such that
\begin{enumerate}
\item $\forall \, 0 \leq j \leq p-1, \:\: \left(A^{--}_{ij}(0)\right)_{0 \leq i \leq p-1}$ is a standard basis of $\mathcal{X}^{-}(p)$ for the left action.
\item $\forall \, 0 \leq i \leq p-1, \:\: \left(A^{--}_{ij}(0)\right)_{0 \leq j \leq p-1}$ is a standard basis of $R^*\left(\mathcal{X}^{-}(p)\right)$ for the right action.
\end{enumerate}
The block $Q(p)$ admits a basis $\left(A^{++}_{ij}(p)\right)_{0 \leq i,j \leq p-1}$ such that
\begin{enumerate}
\item $\forall \, 0 \leq j \leq p-1, \:\: \left(A^{++}_{ij}(p)\right)_{0 \leq i \leq p-1}$ is a standard basis of $\mathcal{X}^{+}(p)$ for the left action.
\item $\forall \, 0 \leq i \leq p-1, \:\: \left(A^{++}_{ij}(p)\right)_{0 \leq j \leq p-1}$ is a standard basis of $R^*\left(\mathcal{X}^{+}(p)\right)$ for the right action.
\end{enumerate}
\end{proposition}

\noindent As in \cite{FGST}, the structure of $Q(s)$ in terms of composition factors can be schematically represented as follows (each vertex represents a composition factor and is labelled by the basis vectors of this factor):
{\footnotesize
$$
\xymatrix{
        & \left(B^{++}_{ab}(s)\right) \ar[ld]^E \ar[rd]^F \ar[dd]^E &      &        &  \left(B^{--}_{ij}(s)\right) \ar[ld]^E \ar[rd]^F \ar[dd]^E&    \\
\left(X^{-+}_{cd}(s)\right) \ar[dr]^F &    &  (Y^{-+}_{ef}(s)) \ar[dl]^E  &  \left(X^{+-}_{cd}(s)\right) \ar[dr]^F  &  &  \left(Y^{+-}_{mn}(s) \ar[dl]^E\right) \\
        &  (A^{++}_{gh}(s)) &      &        &  \left(A^{--}_{or}(s)\right) & 
}
$$
}
for the left action, and
{\footnotesize
$$
\xymatrix{
        & \left(B^{++}_{ab}(s)\right) \ar[rrd]^E \ar[rrrrd]^F \ar[dd]^E&      &        &  \left(B^{--}_{ij}(s)\right) \ar[lllld]^E \ar[lld]^F \ar[dd]^E&    \\
\left(X^{-+}_{cd}(s)\right) \ar[rrrrd]^F &    &  (Y^{-+}_{ef}(s)) \ar[rrd]^E  &  \left(X^{+-}_{cd}(s)\right) \ar[lld]^F &  &  \left(Y^{+-}_{mn}(s)\right) \ar[lllld]^E\\
        &  (A^{++}_{gh}(s)) &      &        &  \left(A^{--}_{or}(s)\right) & 
}
$$
}
for the right action.
\smallskip\\
\indent The knowledge of the structure of the bimodule $_{\bar U_q}\!\left(\bar U_q\right)_{\bar U_q}$ allows us to determine the center of $\bar U_q$. Indeed, each central element determines a bimodule endomorphism and conversely. Recall from \cite{FGST} that $\mathcal{Z}(\bar U_q)$ is a $(3p-1)$-dimensional algebra with basis elements $e_s \:\:(0 \leq s \leq p)$ and $w^{\pm}_t \:\: (1 \leq t \leq p-1)$. The element $e_s$ is just the unit of the block $Q(s)$, thus by (\ref{leftBlock}) and (\ref{figureProj}) the action of $e_s$ on the simple and the projective modules is given by
\begin{equation}\label{actionEs}
\begin{array}{llll}
\text{For } s=0, \:\: &e_0 v_0^+(t) = 0, & e_0 v_0^-(t) =\delta_{t,p}v_0^-(p), & e_0 b_0^{\pm}(t) = 0, \\
\text{For }1 \leq s \leq p-1, \:\: &e_s v_0^{+}(t) = \delta_{s,t}v_0^+(s), & e_s v_0^-(t) = \delta_{p-s,t}v_0^{-}(p-s), & \\
 & e_s b_0^+(t) = \delta_{s,t}b_0^+(s), & e_s b_0^-(t) = \delta_{t,p-s}b_0^-(p-s), & \\
\text{For } s=p, \:\: &e_p v_0^+(t) = \delta_{t,p}v_0^+(p), & e_p v_0^-(t) =0, & e_p b_0^{\pm}(t) = 0
\end{array}
\end{equation}
while for the elements $w^{\pm}_s$:
\begin{equation}\label{actionWs}
\begin{array}{lll}
w^+_s v_0^{\pm}(t) = 0, & w^+_s b_0^+(t) = \delta_{s,t}a_0^+(s), & w^+_s b_0^-(t) = 0,\\
w^-_s v_0^{\pm}(t) = 0, & w^-_s b_0^+(t) = 0, & w^-_s b_0^-(t) = \delta_{t,p-s}a_0^-(p-s).\\
\end{array}
\end{equation}
Observe that 
$$\overset{\mathcal{P}^+(s)}{w_s^+} = p^+_s, \:\:\:\: \overset{\mathcal{P}^-(p-s)}{w_s^-} = p^-_{p-s}.$$
The action of the central elements on $\mathcal{P}^{\alpha}(s)$ is enough to recover their action on every module, using projective covers. From these formulas, we deduce the multiplication rules of these elements:
\begin{equation}\label{produitCentre}
e_se_t = \delta_{s,t}e_s, \:\:\: e_sw^{\pm}_t = \delta_{s,t}w^{\pm}_s, \:\:\: w^{\pm}_sw^{\pm}_t = 0.
\end{equation}
\noindent Let us mention that the idempotents $e_s$ are not primitive: there exists primitive orthogonal idempotents $e_{s,i}$ such that $e_s = \sum_{i}e_{s,i}$, see \cite{arike}.

\begin{definition}\label{defBaseCanoCentre}
The basis $\{ e_s, w^{\pm}_t \}$ ($0 \leq s \leq p$, $1 \leq t \leq p-1$) will be called the canonical basis of $\mathcal{Z}(\bar U_q)$.
\end{definition}

\indent If $z$ is a central element and $S$ is a simple module, we know by Schur lemma that $\overset{S}{z} = z_S \mathrm{id}_S$ for some scalar $z_S \in \mathbb{C}$. For a simple $\bar U_q$-module, we see thanks to \eqref{actionEs} and \eqref{actionWs} that the scalars $z_{\mathcal{X}^{\alpha}(s)}$ satisfy a symmetry property:
\begin{equation}\label{symetrieCentraux}
\forall \, 1 \leq s \leq p-1, \:\:\:\:z_{\mathcal{X}^+(s)} = z_{\mathcal{X}^-(p-s)}.
\end{equation}
\noindent We will sometimes use the convention $z_{\mathcal{X}^+(0)} = z_{\mathcal{X}^-(p)}$ to unify formulas.

\smallskip

\indent An important and useful central element of $\bar U_q$ is the Casimir element $C$, defined by
\begin{equation}\label{defCasimir}
C = FE + \frac{qK + q^{-1}K^{-1}}{(q - q^{-1})^2}
\end{equation}
Since $Cb_0^{\alpha}(s) = a_0^{\alpha}(s) + \alpha \frac{q^s + q^{-s}}{(q-q^{-1})^2}b_0^{\alpha}(s)$, we get that the expression of $C$ in the canonical basis of $\mathcal{Z}(\bar U_q)$ is
\begin{equation}\label{casimirBaseCan}
C = \sum_{j=0}^p \frac{q^j + q^{-j}}{ (q - q^{-1})^2 } e_j + \sum_{k = 1}^{p-1} (w^+_k + w^-_k).
\end{equation}
Moreover, thanks to \cite[Formula (D.7)]{FGST}, we know that the subalgebra of $\mathcal{Z}(\bar U_q)$ generated by $C$ is
\begin{equation}\label{engendreParCasimir}
\mathbb{C}\langle C \rangle = \mathrm{vect}\bigl( e_s, w_t^+ + w_t^- \bigr)_{0 \leq s \leq p, \: 1 \leq t \leq p-1}.
\end{equation}
\noindent This means that any of the elements $e_s, w_t^+ + w_t^-$ can be written as a polynomial in $C$.

\subsection{The braided extension of $\bar U_q$}\label{braidedExtension}
\indent The Hopf algebra $\bar U_q$ is not braided. Indeed, for $p>2$, one can find (see \cite{KS}) two $\bar U_q$-modules $V, W$ such that $V \otimes W$ and $W \otimes V$ are not isomorphic, which immediately implies that $\bar U_q$ cannot contain a $R$-matrix; the remaining case $p=2$ is considered in \cite{GR} where it is shown directly that there is no $R$-matrix either. However, $\bar U_q$ is very close to be braided, since its extension by a square root of $K$ is braided. Let $\bar U_q^{1/2}$ be this extension; the universal $R$-matrix $R \in \bar U_q^{1/2} \otimes \bar U_q^{1/2}$ is given by
\begin{equation}\label{RMatriceUq}
R = q^{H \otimes H/2} \sum_{m=0}^{p-1}\frac{\hat q^m}{[m]!}q^{m(m-1)/2} E^m \otimes F^m, \: \text{ with } \:\: q^{H \otimes H/2} = \frac{1}{4p}\sum_{n,j=0}^{4p-1}q^{-nj/2} K^{n/2} \otimes K^{j/2} 
\end{equation}
where $q^{1/2}$ is a fixed square root of $q$. We use the notation $q^{H \otimes H/2}$ because $q^{H \otimes H/2} v \otimes w = q^{ab/2}$ if $K^{1/2}v = q^{a/2} v$ and $K^{1/2}w = q^{b/2}w$; also recall the notation $\hat q = q - q^{-1}$. 

\smallskip

Even if $R \not\in \bar U_q \otimes \bar U_q$, the $R$-matrix satisfies the important property that $RR' \in \bar U_q \otimes \bar U_q$. Its value is
\begin{equation}\label{RRPrime}
RR' = \frac{1}{2p} \sum_{m,n = 0}^{p-1}\sum_{i,j=0}^{2p-1} \frac{\hat q^{m+n}}{[m]![n]!}q^{\frac{m(m-1)}{2} + \frac{n(n-1)}{2} - m^2 +m(i-j) - ij } E^mK^iF^n \otimes F^m K^j E^n.
\end{equation}
with $\hat q = q - q^{-1}$. From this expression, we see that the map
\begin{equation*}
\fonc{\Psi}{\bar U_q^*}{\bar U_q}{\beta}{(\beta \otimes \mathrm{id})(RR')}
\end{equation*}
is an isomorphism of vector spaces. Thus by abuse of terminology we will say that $\bar U_q$ is factorizable (this is an abuse of terminology since the usual definition of factorizability requires braiding (\textit{i.e.} existence of a $R$-matrix) and $\bar U_q$ is not braided). Note that the extension $\bar U_q^{1/2}$ is not factorizable since $K^{1/2}$ does not appear in the expression of $RR'$.

\smallskip

The Drinfeld element $u = S(b_i)a_i$ (with $R = a_i \otimes b_i$) also belongs to $\bar U_q$ (the square root of $K$ does not appear in its expression). Moreover, $\bar U_q$ contains two possible ribbon elements, namely elements $v \in \bar U_q$ which satisfy \eqref{ribbon}. Here we take
\[ v = \frac{1 - i}{2 \sqrt{p}} \sum_{m=0}^{p-1}\sum_{j=0}^{2p-1} \frac{\hat q^m}{[m]!} q^{-\frac{m}{2} -mj +\frac{(j+p+1)^2}{2}} F^mK^jE^m.\]
With this choice of $v$, it holds (see \eqref{pivotCan})
\[ g = uv^{-1} = K^{p+1}. \]
The choice of the other possible ribbon element would have led to $g = uv^{-1} = K$, but from the begining we have decided to take $K^{p+1}$ as pivotal element, which forces the choice of $v$. The element $v$ is central and invertible; its expression in the canonical basis of $\mathcal{Z}(\bar U_q)$ is
\begin{equation}\label{rubanCentre}
\begin{split}
v &= \sum_{s=0}^p v_{\mathcal{X}^+(s)} e_s + \hat q \sum_{s=1}^{p-1} v_{\mathcal{X}^+(s)}\left( \frac{p-s}{[s]}w^+_s - \frac{s}{[s]}w^-_s \right),\\
v^{-1} &= \sum_{s=0}^p v_{\mathcal{X}^+(s)}^{-1} e_s - \hat q \sum_{s=1}^{p-1} v_{\mathcal{X}^+(s)}^{-1}\left( \frac{p-s}{[s]}w^+_s - \frac{s}{[s]}w^-_s \right).
\end{split}
\end{equation}
with $\hat q = q - q^{-1}$. The scalar $v_{\mathcal{X}^{\alpha}(s)}$ is defined by $\overset{\mathcal{X}^{\alpha}(s)}{v} = v_{\mathcal{X}^{\alpha}(s)}\mathrm{id}$, its value is
\begin{equation}\label{valueVRep}
v_{\mathcal{X}^+(s)} = v_{\mathcal{X}^-(p-s)} = (-1)^{s-1}q^{\frac{-(s^2-1)}{2}}.
\end{equation}
and $v_{\mathcal{X}^+(0)}$ is just a notation for $v_{\mathcal{X}^-(p)}$ used to unify the formula.

\subsection{Matrix coefficients for $\bar U_q$}\label{sectionMatrixCoeffUq}
\indent This section is an example and will not be used in the rest of the text. The aim is to explain how one can restrict to well-chosen representations when he deals with algebras defined by means of matrix coefficients (which is the case of all the algebras considered in this text). Here we will exhibit a (well-known) minimal set of matrix coefficients which generate $\mathcal{O}(\bar U_q)$ and give a presentation by generators and relations of this algebra based on these matrix coefficients. The reasoning is more interesting than the result because the method presented here can be applied to the algebra $\mathcal{L}_{0,1}(\bar U_q)$ and then to $\mathcal{L}_{g,n}(\bar U_q)$ (see section \ref{sectionL01Uq}).

\smallskip The decomposition formulas for tensor products of simple modules and projective modules are given in \cite{KS}. First, since
\[ \mathcal{X}^+(2) \otimes \mathcal{X}^+(s) \cong \mathcal{X}^+(s-1) \oplus \mathcal{X}^+(s+1) \]
for $2 \leq s \leq p-1$, we see that $\mathcal{X}^+(p) = \mathcal{P}^+(p)$ is a direct summand of $\mathcal{X}^+(2)^{\otimes (p-1)}$. Second, since
\[\mathcal{X}^+(2) \otimes \mathcal{X}^+(p) \cong \mathcal{P}^+(p-1), \:\:\:\:\:\: \mathcal{X}^+(2) \otimes \mathcal{P}^+(s) \cong \mathcal{P}^+(s-1) \oplus \mathcal{P}^+(s+1), \:\:\:\:\:\: \mathcal{X}^+(2) \otimes \mathcal{P}^+(1) \cong 2\mathcal{X}^-(p) \oplus \mathcal{P}^+(2) \]
for $2 \leq s \leq p-1$, we see that $\mathcal{P}^+(j)$ is a direct summand of $\mathcal{X}^+(2)^{\otimes 2p-j-1}$, and that $\mathcal{X}^-(p) = \mathcal{P}^-(p)$ is a direct summand of $\mathcal{X}^+(2)^{\otimes 2p-1}$. Finally, since
\[\mathcal{X}^+(2) \otimes \mathcal{X}^-(p) \cong \mathcal{P}^-(p-1), \:\:\:\:\:\: \mathcal{X}^+(2) \otimes \mathcal{P}^-(s) \cong \mathcal{P}^-(s-1) \oplus \mathcal{P}^-(s+1) \]
for $2 \leq s \leq p-1$, we see that $\mathcal{P}^-(j)$ is a direct summand of $\mathcal{X}^+(2)^{\otimes 3p-j-1}$. It follows that every PIM $P$ is a direct summand of some tensor power $\mathcal{X}^+(2)^{\otimes n}$, and that $\overset{P}{T}$ is a submatrix of $\overset{\mathcal{X}^+(2)^{\otimes n}}{T}$. Hence the coefficients $\bigl(\overset{\mathcal{X}^+(2)^{\otimes n}}{T}\bigr)^{i_1, \ldots, i_n}_{j_1, \ldots, j_n}$ linearly span $\mathcal{O}(\bar U_q)$. But thanks to the fusion relation \eqref{dualHopf} of the matrices $\overset{I}{T}$, it holds:
\begin{equation}\label{TX2tensN}
\bigl(\overset{\mathcal{X}^+(2)^{\otimes n}}{T}\bigr)^{i_1 \ldots i_n}_{j_1 \ldots j_n} = \bigl(\overset{\mathcal{X}^+(2)}{T}\bigr)^{i_1}_{j_1} \ldots \bigl(\overset{\mathcal{X}^+(2)}{T}\bigr)^{i_n}_{j_n}
\end{equation}
and this shows that the matrix coefficients $\bigl(\overset{\mathcal{X}^+(2)}{T}\bigr)^{i}_{j}$ generate $\mathcal{O}(\bar U_q)$ as an algebra.

\smallskip

\indent We denote
\[ \overset{\mathcal{X}^+(2)}{T}=
\begin{pmatrix}
a & b\\
c & d
\end{pmatrix}.
\]
\noindent Let us seek relations between these generators. All the relations are implied by the existence of certain well-chosen morphisms. First, one has the exchange relation \eqref{FRT} (which comes from the existence of the braiding isomorphism $c : \mathcal{X}^+(2)^{\otimes 2} \to \mathcal{X}^+(2)^{\otimes 2}$). The $R$-matrix \eqref{RMatriceUq} evaluated in $\mathcal{X}^+(2)^{\otimes 2}$ is
\[ 
q^{-1/2}
\begin{pmatrix}
q & 0 & 0 & 0\\
0 & 1 & \hat q & 0\\
0 & 0 & 1 & 0\\
0 & 0 & 0 & q
\end{pmatrix}
\]
(with $\hat q =q - q^{-1}$) and relation \eqref{FRT} is equivalent to
\[ ba=qab, \:\:\: db=qbd, \:\:\:, ca=qac,  \:\:\: dc=qcd, \:\:\: bc=cb, \:\:\: ad-da=(q^{-1}-q)bc. \]
Second, since $\mathcal{X}^+(2)^{\otimes 2} \cong \mathcal{X}^+(1) \oplus \mathcal{X}^+(3)$, there exists a unique (up to scalar) morphism \\$\Phi : \mathbb{C} = \mathcal{X}^+(1) \to \mathcal{X}^+(2)^{\otimes 2}$; it is given by $ \Phi(1) = q v_0 \otimes v_1 - v_1 \otimes v_0$.
By naturality \eqref{naturaliteT} and fusion \eqref{dualHopf}, we have $\overset{\mathcal{X}^+(2)}{T_1}\overset{\mathcal{X}^+(2)}{T_2}\Phi = \overset{\mathcal{X}^+(2)^{\otimes 2}}{T}\Phi  = \Phi \overset{\mathbb{C}}{T} = \Phi$.
This gives just one new relation, the quantum determinant:
\begin{equation}\label{quantumDet}
ad-q^{-1}bc=1. 
\end{equation}
Next, since $\mathcal{P}^+(p-1)$ is a direct summand of $\mathcal{X}^+(2)^{\otimes p}$, there exists\footnote{Such a morphism is far from unique. Indeed, using the decomposition rules recalled above, one can show that $\mathcal{X}^+(2)^{\otimes p} \cong W \oplus (p-2)\mathcal{X}^+(p-1) \oplus \mathcal{P}^+(p-1)$, and thanks to the description of Hom-spaces in section \ref{sectionSimpleProj}, we get that $\dim\bigl( \mathrm{Hom}_{\bar U_q}(\mathcal{P}^+(p-1), \mathcal{X}^+(2)^{\otimes p} \bigr) = p$.} an injection $f : \mathcal{P}^+(p-1) \to \mathcal{X}^+(2)^{\otimes p}$. For instance, one can check that the assignment
\[ f\bigl(b_0^+(p-1)\bigr) = v_0^{\otimes (p-1)} \otimes v_1 \]
does the job, and we have $f\bigl(x_0^+(p-1)\bigr) = v_0^{\otimes p}$, $f\bigl(y_0^+(p-1)\bigr) = \lambda v_1^{\otimes p}$ for some $\lambda$ ($\neq 0$ since $f$ is injective). Endowing the tensor basis of $\mathcal{X}^+(2)^{\otimes p}$ with the lexicographic order and using \eqref{TX2tensN} above and \eqref{matriceTPim} below, we get that the matrices under consideration have the following shapes:
\begin{center}
$f=
\begin{blockarray}{ccccc}
 (b_i) & x_0 & y_0 & (a_j)  \\
  \begin{block}{(cccc)c}
    \mathbf{0} & 1 & 0 & \mathbf{0} & v_0^{\otimes p}  \\
  \ast & \mathbf{0} & \mathbf{0} & \ast & \vdots\\
  \mathbf{0}  & 0 & \lambda & \mathbf{0}  & v_1^{\otimes p}  \\
  \end{block}
\end{blockarray}$
~~~~~~~~
$\overset{\mathcal{X}^+(2)^{\otimes p}}{T} = 
\begin{blockarray}{cccc}
 v_0^{\otimes p} & \ldots & v_1^{\otimes p}  \\
  \begin{block}{(ccc)c}
    a^p & \ast & b^p & v_0^{\otimes p}  \\
   \ast  & \ast & \ast & \vdots\\
  c^p & \ast & d^p & v_1^{\otimes p}\\
  \end{block}
\end{blockarray}$
~~~~~~~~
$\overset{\mathcal{P}^+(p-1)}{T} = 
\begin{blockarray}{ccccc}
 (b_i) & x_0 & y_0 & (a_j)  \\
  \begin{block}{(cccc)c}
    \overset{\mathcal{X}^+(p-1)}{T} & \mathbf{0} & \mathbf{0} & \mathbf{0} & (b_i)  \\
   A^+_{p-1}  & \overset{\mathcal{X}^-(1)}{T} & \mathbf{0} & \mathbf{0} & x_0\\
  B^+_{p-1} & \mathbf{0} & \overset{\mathcal{X}^-(1)}{T} &  \mathbf{0} & y_0\\
  H^+_{p-1}  & D^+_{p-1} & C^+_{p-1} & \overset{\mathcal{X}^+(p-1)}{T}  & (a_j)  \\
  \end{block}
\end{blockarray}$
\end{center}
To obtain the blocks of $\mathbf{0}$'s in the matrix of $f$, just compare the weights of the elements. The relation $f\overset{\mathcal{P}^+(p-1)}{T} = \overset{\mathcal{X}^+(2)^{\otimes p}}{T}f$ implies $b^p = c^p = 0$. Finally, $\mathcal{P}^-(p-1)$ is a direct summand of $\mathcal{X}^+(2)^{\otimes 2p}$, and it is not difficult to check that the assignments
\[ f_1\bigl( b_0^-(p-1) \bigr) = v_0^{\otimes (2p-1)} \otimes v_1, \:\:\:\:\:\:\:\:\:\:\: f_2\bigl( b_{p-2}^-(p-1) \bigr) = v_0 \otimes v_1^{\otimes (2p-1)} \]
define injective morphisms $f_1, f_2 : \mathcal{P}^-(p-1) \to \mathcal{X}^+(2)^{\otimes 2p}$, and we have $f_1\bigl(x_0^-(p-1)\bigr) = v_0^{\otimes 2p}$, $f_2\bigl( y_0^-(p-1) \bigr) = v_1^{\otimes 2p}$. One can compute as above that the relation $f_1\overset{\mathcal{P}^-(p-1)}{T} = \overset{\mathcal{X}^+(2)^{\otimes 2p}}{T}f_1$ (resp. $f_2\overset{\mathcal{P}^-(p-1)}{T} = \overset{\mathcal{X}^+(2)^{\otimes 2p}}{T}f_2$) implies $a^{2p}=1$ (resp. $d^{2p}=1$), where $1 = \varepsilon = \overset{\mathcal{X}^+(1)}{T}$. We then arrive to the following Proposition, which is certainly well-known.

\begin{proposition}
The algebra $\mathcal{O}(\bar U_q)$ admits the following presentation:
\[ \left\langle b,c,d \: \left| \: db=qbd, \:\: dc=qcd, \:\: bc=cb, \:\: b^p = c^p=0, \:\: d^{2p}=1 \right. \right\rangle \]
\end{proposition}
\begin{proof}
Note first that the generator $a$ is not required. Indeed since $d$ is invertible, we have $a = d^{-1} + q^{-1}bcd^{-1}$. Let $A$ be the algebra defined by this presentation. Then by the computations above, we have a surjection $p : A \to \mathcal{O}(\bar U_q)$  and thus $\dim\bigl(\mathcal{O}(\bar U_q)\bigr) \leq \dim(A)$. But it is clear that any element of $A$ is a linear combination of monomials $b^ic^jd^k$ with $0 \leq i,j \leq p-1$, $0 \leq k \leq 2p-1$. Hence $\dim(A) \leq 2p^3 = \dim\bigl(\mathcal{O}(\bar U_q)\bigr)$ and $p$ is an isomorphism.
\end{proof}

\section{Symmetric linear forms and the GTA basis}\label{sectionSLF}
\indent Let
$$ \SLF(\bar U_q) = \left\{\varphi \in \bar U_q^* \, \vert \, \forall\, x,y \in \bar U_q, \:\: \varphi(xy) = \varphi(yx)\right\}. $$
From the general comments of section \ref{rappelHopf}, $\SLF(\bar U_q)$ is a subalgebra of $\mathcal{O}(\bar U_q)$. $\mathcal{O}(\bar U_q)$ is more precisely a Hopf algebra, but $\SLF(\bar U_q)$ is not a sub-coalgebra of $\mathcal{O}(\bar U_q)$, see Remark \ref{remarkCogebreSLF} below; it is however stable by the antipode $S$.

\smallskip

\indent Since $\bar U_q$ is factorizable (in the generalized sense of section \ref{braidedExtension}), we know thanks to Lemma \ref{propMorDrinfeld} that
\[ \dim\bigl( \mathrm{SLF}(\bar U_q) \bigr) = \dim\bigl( \mathcal{Z}(\bar U_q) \bigr) = 3p-1. \]

\indent An interesting basis of $\SLF(\bar U_q)$ was found by Gainutdinov and Tipunin in \cite{GT} and by Arike in \cite{arike}. To be precise, a basis of the space $\mathrm{qCh}(\bar U_q)$ of $q$-characters is constructed in \cite{GT}, but the shift by the pivotal element $g=K^{p+1}$ provides an isomorphism
$$\mathrm{qCh}(\bar U_q) \overset{\sim}{\rightarrow} \mathrm{SLF}(\bar U_q), \:\:\:\:\psi \mapsto\psi(g\,\cdot).$$
\indent This basis is built from the simple and the projective modules. First, define $2p$ linear forms\footnote{The correspondence of notations with \cite{arike} is: $T^+_s = \chi^+_s$, $T^-_s = \chi^-_{p-s}$. The letter $T$ is here reserved for the matrices $\overset{V}{T}$ described above.} $\chi^{\alpha}_s$, $\alpha \in \{\pm\}, 1 \leq s \leq p$, by:
\begin{equation}\label{defChis}
\chi^{\alpha}_s = \text{tr}(\overset{\mathcal{X}^{\alpha}(s)}{T}).
\end{equation}
They are obviously symmetric. Observe that $\chi^+_1 = \varepsilon$ is the unit for the algebra structure on $\SLF(\bar U_q)$ described above. To construct the $p-1$ missing linear forms, observe with the help of (\ref{figureProj}) that the matrix of the action on $\mathcal{P}^{\alpha}(s)$ has the following block form in a standard basis:
\begin{equation}\label{matriceTPim}
\overset{\mathcal{P}^{\alpha}(s)}{T} =
\begin{blockarray}{ccccc}
(b_i) & (x_j) & (y_k) & (a_{l}) &\\
\begin{block}{(cccc)c}
\overset{\mathcal{X}^{\alpha}(s)}{T} & \mathbf{0} & \mathbf{0} & \mathbf{0} & (b_i) \\
A^{\alpha}_s & \overset{\mathcal{X}^{-\alpha}(p-s)}{T} & \mathbf{0} & \mathbf{0} & (x_j)\\
B^{\alpha}_s & \mathbf{0} & \overset{\mathcal{X}^{-\alpha}(p-s)}{T} & \mathbf{0} & (y_k)\\
H^{\alpha}_s & D^{\alpha}_s & C^{\alpha}_s & \overset{\mathcal{X}^{\alpha}(s)}{T} & (a_{l}).\\
\end{block}
\end{blockarray}
\end{equation}
It is not difficult to see that these matrices satisfy the following symmetries:
\begin{equation*}
A^-_{p-s} = C^+_s, \:\:\: B^-_{p-s} = D^+_s, \:\:\: D^-_{p-s} = B^+_s, \:\:\: C^-_{p-s} = A^+_s .
\end{equation*}
By computing the matrices $\overset{\mathcal{P}^{+}(s)}{(xy)} = \overset{\mathcal{P}^{+}(s)}{x}\overset{\mathcal{P}^{+}(s)}{y}$ and $\overset{\mathcal{P}^{-}(p-s)}{(xy)} = \overset{\mathcal{P}^{-}(p-s)}{x}\overset{\mathcal{P}^{-}(p-s)}{y}$, these symmetries allow us to see that the linear form $G_s \: (1 \leq s \leq p-1)$ defined by
\begin{equation}\label{defGs}
G_s = \text{tr}(H^+_s) + \text{tr}(H^-_{p-s})
\end{equation}
is a symmetric linear form. This can also be written as
\begin{equation}\label{GsAvecTopSoc}
G_s = \mathrm{tr}\bigl( \sigma_s \overset{\mathcal{P}^+(s)}{T} \bigr) + \mathrm{tr}\bigl( \sigma_{p-s} \overset{\mathcal{P}^-(p-s)}{T} \bigr) 
\end{equation}
where $\sigma_j : \mathcal{P}^{\alpha}(s) \to \mathcal{P}^{\alpha}(s)$ is the linear map (which is not a $\bar U_q$-morphism) sending $\mathrm{Soc}\bigl(\mathcal{P}^{\alpha}(s)\bigr)$ to $\mathrm{Top}\bigl(\mathcal{P}^{\alpha}(s)\bigr)$ (see \eqref{figureProj}):
\[
\sigma_s =
\begin{blockarray}{ccccc}
(b_i) & (x_j) & (y_k) & (a_{l}) &\\
\begin{block}{(cccc)c}
\mathbf{0} & \mathbf{0} & \mathbf{0} & \mathbb{I}_j & (b_i) \\
\mathbf{0} & \mathbf{0} & \mathbf{0} & \mathbf{0} & (x_j)\\
\mathbf{0} & \mathbf{0} & \mathbf{0} & \mathbf{0} & (y_k)\\
\mathbf{0} & \mathbf{0} & \mathbf{0} & \mathbf{0} & (a_{l}).\\
\end{block}
\end{blockarray}
\]

\indent It is instructive for our purposes to see a proof that these symmetric linear forms are linearly independent. Let us begin by introducing important elements for $0 \leq n \leq p-1$ (they are discrete Fourier transforms of $(K^{l})_{0 \leq l \leq 2p-1}$): 
\begin{equation*}
\Phi^{\alpha}_n = \frac{1}{2p}\sum_{l=0}^{2p-1}\left(\alpha q^{-n}\right)^{l}K^{l}.
\end{equation*}
The following easy lemma shows that these elements allow one to select vectors which have a given weight, and this turns out to be very useful. 
\begin{lemma}\label{lemmaPoids}
1) Let $M$ be a left $\bar U_q$-module, and let $m^+_i(s)$ be a vector of weight $q^{s-1-2i}$, $m^-_i(p-s)$ be a vector of weight $-q^{(p-s)-1-2i} = q^{-s-1-2i}$, $m^-_i(s)$ be a vector of weight $-q^{s-1-2i}$, $m^+_i(p-s)$ be a vector of weight $q^{(p-s)-1-2i} = -q^{-s-1-2i}$. Then: 
\begin{align*}
&\Phi^{+}_{s-1}m_i^+(s) = \delta_{i,0}m^+_0(s), \:\:\: \Phi^{+}_{s-1}m_i^-(p-s) = 0,\\
&\Phi^{-}_{s-1}m_i^-(s) =  \delta_{i,0}m^-_0(s), \:\:\: \Phi^{-}_{s-1}m_i^+(p-s) = 0.
\end{align*}
2) Let $N$ be a right $\bar U_q$-module, and let $n^+_i(s)$ be a vector of weight $q^{1-s+2i}$, $n^-_i(p-s)$ be a vector of weight $-q^{1-(p-s)+2i} = q^{1+s+2i}$, $n^-_i(s)$ be a vector of weight $-q^{1-s+2i}$, $n^+_i(p-s)$ be a vector of weight $q^{1-(p-s)+2i} = -q^{1+s+2i}$. Then:
\begin{align*}
&n_i^+(s)\Phi^{+}_{s-1} = \delta_{i,s-1}n^+_{s-1}(s),\:\:\:n_i^-(p-s)\Phi^{+}_{s-1} = 0,\\
&n_i^-(s)\Phi^{-}_{s-1} = \delta_{i,s-1}n^-_{s-1}(s) ,\:\:\:n_i^+(p-s)\Phi^{-}_{s-1} = 0.
\end{align*}
\end{lemma}
\begin{proof}
It follows from easy computations with sums of roots of unity.
\end{proof}

We can now state the key observation.

\begin{proposition}\label{propArike}
Let
$$\varphi = \sum_{s=1}^{p} \left(\lambda^{+}_s\chi^{+}_s + \lambda^{-}_s\chi^{-}_s\right) + \sum_{s'=1}^{p-1}\mu_{s'} G_{s'} \in \text{\em SLF}\left(\bar U_q\right).$$
Then: 
$$ \lambda^+_s = \varphi\left(\Phi^+_{s-1}e_s\right),\:\: \lambda^-_s = \varphi\left(\Phi^-_{s-1}e_{p-s}\right), \:\: \mu_{s'} = \frac{\varphi\left(w^+_{s'}\right)}{s'} = \frac{\varphi(w_{s'}^-)}{p-s'}. $$
\end{proposition}
\begin{proof}
It is a corollary of (\ref{actionEs}) and (\ref{actionWs}). Indeed, we have:
\begin{equation}\label{actionCentreMatricesBlocs}
\begin{array}{l l l}
\overset{\mathcal{X}^+(s)}{T}(e_t) = \delta_{s,t}I_s, \:\:\: & \overset{\mathcal{X}^+(s)}{T}(w^{\pm}_t) = 0, \:\:\: & \overset{\mathcal{X}^-(s)}{T}(e_t) = \delta_{s, p-t}I_s,\\
\overset{\mathcal{X}^-(s)}{T}(w^{\pm}_t) = 0,& H^{\pm}_s(e_t) = 0, & H^+_s(w^+_s) = \delta_{s,t}I_{s}, \\
H^+_s(w^-_t) = 0, & H^-_{p-s}(w^+_t) = 0, & H^-_{p-s}(w^-_t) = \delta_{s,t}I_{p-s}.
\end{array}
\end{equation}
This gives the formula for $\mu_s$. The formulas for $\lambda^{\pm}_s$ follow from this and Lemma \ref{lemmaPoids}.
\end{proof}

\noindent If we have $\sum_{s=1}^{p} \left(\lambda^{+}_s\chi^{+}_s + \lambda^{-}_s\chi^{-}_s\right) + \sum_{s'=1}^{p-1}\mu_{s'} G_{s'} = 0$, we can evaluate the left-hand side on the elements appearing in Proposition \ref{propArike} to get that all the coefficients are equal to $0$. Thus we have a free family of cardinal $3p-1$, hence a basis of $\SLF(\bar U_q)$.

\begin{theorem}\label{thAri}
The symmetric linear forms $\chi^{\pm}_s \: (1 \leq s \leq p)$ and $G_{s'} \: (1 \leq s' \leq p-1)$ form a basis of $\SLF(\bar U_q)$.
\end{theorem}

\begin{definition}
The basis of Theorem \ref{thAri} will be called the {\em GTA basis} (for Gainutdinov, Tipunin, Arike).
\end{definition}

\begin{remark}\label{remarkCogebreSLF}
Let $\varphi \in \SLF(\bar U_q)$. It is easy to see that $\varphi(K^jE^nF^m) = 0$ if $n \neq m$. From this we deduce that $\SLF(\bar U_q)$ is not a sub-coalgebra of $\bar U_q^*$. Indeed, write $\Delta(\chi^+_2) = \sum_i \varphi_i \otimes \psi_i$, and assume that $\varphi_i, \psi_i \in \SLF(\bar U_q)$. Then $ 1 = \chi^+_2(EF) = \sum_i \varphi_i(E)\psi_i(F) = 0 $, a contradiction.
\finEx
\end{remark}

\begin{remark}
If we choose a basis of $\mathcal{Z}(\bar U_q)$, then its dual basis can not be entirely contained in $\SLF(\bar U_q)$. Indeed, let $\varphi = \sum_{s=0}^p\lambda^{\pm}_s \chi^{\pm}_s  + \sum_{s=1}^{p-1}\mu_s G_s \in \SLF(\bar U_q)$. Then $\varphi(w^+_s) = s\mu_s, \varphi(w^-_s) = (p-s)\mu_s$, and we see that there does not exist $\varphi \in \SLF(\bar U_q)$ such that $\varphi(w^+_s) = 1, \: \varphi(w^-_s) = 0$. Hence, $\SLF(\bar U_q) \subset \bar U_q^*$ is not the dual of $\mathcal{Z}(\bar U_q) \subset \bar U_q$.
\finEx
\end{remark}

\indent We have an obvious action of $\mathcal{Z}(\bar U_q)$ on $\mathrm{SLF}(\bar U_q)$:
\begin{equation}\label{defActionCentreSLF}
\forall\, z \in \mathcal{Z}(\bar U_q), \:\: \forall \, \varphi \in \mathrm{SLF}(\bar U_q), \:\:\:\: \varphi^z = \varphi(z?).
\end{equation}
For further use, we record the values of this action with the canonical basis of $\mathcal{Z}(\bar U_q)$ and the GTA basis of $\mathrm{SLF}(\bar U_q)$:
\begin{equation}\label{actionCentreSLF}
\begin{array}{llll}
(\chi^+_s)^{e_t} = \delta_{s,t}\chi^+_s, & (\chi^-_s)^{e_t} = \delta_{p-s,t}\chi^-_s, & G_{s}^{\, e_t} = \delta_{s,t}G_{s} & \\
(\chi^+_s)^{w^{\pm}_{t}} = 0, & (\chi^-_s)^{w^{\pm}_{t}} = 0, & G_{s}^{\, w^+_{t}} = \delta_{s,t}\chi^+_{s}, & G_{s}^{\, w^-_{t}} = \delta_{s,t}\chi^-_{p-s}.
\end{array}
\end{equation}
\noindent This is jut a consequence of \eqref{actionCentreMatricesBlocs}.

\section{Traces on projective $\bar U_q$-modules and the GTA basis}\label{sectionArike}
\indent The material of this section is independent of the rest of the text and will not be used elsewhere (except for the equality \eqref{traceModifieeGTA}).
\subsection{Correspondence between traces and symmetric linear forms}\label{GenCor}
\indent Let $A$ be a finite dimensional $k$-algebra. We have an anti-isomorphism of algebras:
$$A \to\End_{A}(A), \:\:\: a \mapsto \rho_a \text{ defined by } \rho_a(x) = xa. $$
Observe that the right action of $A$ naturally appears.
Let $t$ be a trace on $A$, that is, an element of $\SLF(\End_A(A))$. Then:
$$t(\rho_{ab}) = t(\rho_b \circ \rho_a) = t(\rho_a \circ \rho_b) = t(\rho_{ba}).$$
So we get an isomorphism of vector spaces
\begin{equation*}
\fleche{\left\{\text{Traces on } \End_{A}(A)\right\} = \SLF\left(\End_{A}(A)\right)}{\SLF(A)}{t}{\varphi^t \text{ defined by } \varphi^t(a) = t(\rho_a).}
\end{equation*}
whose inverse is:
\begin{equation*}
\fleche{\SLF(A)}{\left\{\text{Traces on } \End_{A}(A)\right\} = \SLF\left(\End_{A}(A)\right)}{\varphi}{t^{\varphi} \text{ defined by } t^{\varphi}(\rho_a) = \varphi(a).}
\end{equation*}
In the case of $A=\bar U_q$, we can express $\varphi^t$ in the GTA basis, which will be the object of the next section.
\medskip\\
\indent Let $\Proj_A$ be the full subcategory of the category of finite dimensional $A$-modules whose objects are the projective $A$-modules.
\begin{definition}
A {\em trace} on $\Proj_A$ is a family of linear maps $t = \left(t_U : \End_A(U) \to k\right)_{U \in \Proj_A}$ such that
$$
\forall \, f \in \Hom_A(U,V), \: \forall\, g \in \Hom_A(V,U), \:\:  t_V(g \circ f) = t_U(f \circ g).
$$
We denote by $\mathcal{T}_{\Proj_A}$ the vector space of traces on $\Proj_A$.
\end{definition}
This cyclic property of traces on $\Proj_A$ is one of the axioms of the so-called modified traces, defined for instance in \cite{GKP}. Note that this definition could be restated in the following way (and could be generalized to other abelian full subcategories than $\Proj_A$).
\begin{lemma}\label{lemmaCoherence}
Let $t=(t_U : \End_A(U) \to k)_{U \in \Proj_A}$ be a family of linear maps. Then $t$ is a trace on $\text{\em Proj}_A$ if and only if:
\begin{itemize}
\item $\forall\, f,g \in \End_A(U), \:\:t_U(g \circ f) = t_U(f \circ g)$,
\item $t_{U \oplus V}(f) = t_{U}(p_U \circ f \circ i_U) + t_V(p_V \circ f \circ i_V)$, where $p_U, p_V$ are the canonical projection maps and $i_U, i_V$ are the canonical injection maps.
\end{itemize}
\end{lemma} 
\begin{proof}
If $t$ is a trace and $f \in \End_A(U \oplus V)$, we have:
$$ t_{U \oplus V}(f) =  t_{U \oplus V}((i_Up_U + i_Vp_V)f) = t_U(p_Ufi_U) + t_V(p_Vfi_V). $$
Conversely, let $f : U \to V$, $g : V \to U$. Define $F = i_V f p_U, G = i_U g p_V$. Then $FG = i_V f g p_V$ and $GF = i_U g f p_U$. We have $p_U GF i_U = gf$, $p_V GF i_V = 0$, $p_U FG i_U = 0$, $p_V FG i_V = fg$, thus:
$$ t_{V}(fg) = t_{U \oplus V}(FG) = t_{U \oplus V}(GF) = t_U(gf). $$
This shows the equivalence.
\end{proof}

\indent Now, consider:
\begin{equation*}
\begin{array}{crcccl}
\Pi_A \: : & \mathcal{T}_{\Proj_A} & \to & \SLF(\End_A(A)) & \overset{\sim}{\to} & \SLF(A)\\
&t=(t_U)_{U \in \Ob(\text{Proj}_A)} & \mapsto & t_A & \mapsto & \varphi^t \text{ defined by } \varphi^t(a) = t_A(\rho_a).
\end{array}
\end{equation*}

\begin{theorem}\label{corTracesSLF}
The map $\Pi_A$ is an isomorphism. In other words, $t_A$ entirely characterizes $t = (t_U)$.
\end{theorem}
\begin{proof}
For all the facts concerning PIMs (Principal Indecomposable Modules) and idempotents in finite dimensional $k$-algebras, we refer to \cite[Chap. VIII]{CR}. We first show that $\Pi_A$ is surjective. Let: 
$$ 1 = e_1 + \ldots + e_n $$
be a decomposition of the unit into primitive orthogonal idempotents ($e_ie_j = \delta_{i,j}e_i$). Then the PIMs of $A$ are isomorphic to the left ideals $Ae_i$ (possibly with multiplicity). We have isomorphisms of vector spaces:
$$ \Hom_A(Ae_i, Ae_j) \overset{\sim}{\longrightarrow} e_iAe_j,\:\:\:\: f \mapsto f(e_i). $$
For every $\varphi \in \SLF(A)$, define $t^{\varphi}_{Ae_i}$ by:
$$ t^{\varphi}_{Ae_i}(f) = \varphi(f(e_i)). $$
Let $f : Ae_i \to Ae_j$, $g : Ae_j \to Ae_i$, and put $f(e_i) = e_ia_fe_j$, $g(e_j) = e_ja_ge_i$. Then using the idempotence of the $e_i$'s and the symmetry of $\varphi$ we get:
$$ t^{\varphi}_{Ae_i}(g \circ f) = \varphi(g \circ f(e_i)) = \varphi\left((e_ia_fe_j)(e_ja_ge_i)\right) = \varphi\left((e_ja_ge_i)(e_ia_fe_j)\right) = \varphi(f \circ g(e_j)) =  t^{\varphi}_{Ae_j}(f \circ g). $$

We know that every projective module is isomorphic to a direct sum of PIMs, so we extend $t^{\varphi}$ to $\text{Proj}_A$ by the following formula:
$$ t_{\bigoplus_{l} A_{l}}(f) = \sum_{l}t_{Ae_{l}}(i_{l} \circ f \circ p_{l}) $$
where $p_{l}$ and $i_{l}$ are the canonical injection and projection maps. By Lemma \ref{lemmaCoherence}, this defines a trace on $\Proj_A$.
We then show that $\Pi_A(t^{\varphi}) = \varphi$, proving surjectivity:
\begin{align*}
\Pi_A(t^{\varphi})(a) &= t^{\varphi}_A(\rho_a) = \sum_{j=1}^n t^{\varphi}_{Ae_j}\left(p_{j} \circ {\rho_a} \circ i_j\right) = \sum_{j=1}^n \varphi\left(p_{j} \circ {\rho_a}(e_j)\right) = \sum_{j=1}^n \varphi\left(p_{Ae_j}(e_ja)\right) \\&= \sum_{j, k = 1}^n \varphi\left(p_{Ae_j}(e_jae_k)\right) = \sum_{j=1}^n \varphi\left(e_jae_j\right) = \sum_{j=1}^n \varphi\left(ae_j\right) = \varphi(a).
\end{align*}
Note that we used that the $e_j$'s are idempotents and that $a = \sum_{j=1}^n ae_j$. We now show injectivity. Assume that $\Pi_A(t)=0$. Then:
$$\forall \, a \in A, \:\: t_A(\rho_a) = \sum_{j=1}^n t_{Ae_j}(p_{j} \circ {\rho_a} \circ i_j) = 0. $$
Let $f : Ae_j \to Ae_j$, with $f(e_{j}) = e_{j}a_fe_{j}$. Since $\rho_{f(e_{j})}(e_{l}) = \delta_{j,l}e_{j}a_fe_{j}$, we have $p_{j} \circ \rho_{f(e_{j})} \circ i_{j} = f$ and $p_{l} \circ \rho_{f(e_{j})} \circ i_{l}=0$ if $l \neq j$. Hence:
$$ t_{Ae_{j}}(f) = t_A(\rho_{f(e_{j})}) = 0. $$
Then $t_{Ae_{j}} = 0$ for each $j$, so that $t=0$.
\end{proof}

\subsection{Link with the GTA basis}
\indent We leave the general case and focus on $A=\bar U_q$. The following theorem expresses $\Pi_{\bar U_q}$ in the GTA basis.
\begin{theorem}\label{MainResult}
Let $t=(t_U)_{U \in\Proj_{\bar U_q}}$ be a trace on $\Proj_{\bar U_q}$. Then:
\[
\Pi_{\bar U_q}(t) = t_{\mathcal{X}^+(p)}(\text{\em Id})\chi^+_p + t_{\mathcal{X}^-(p)}(\text{\em Id})\chi^-_p + \sum_{s=1}^{p-1} \left(t_{\mathcal{P}^+(s)}(\text{\em Id})\chi^+_s + t_{\mathcal{P}^-(s)}(\text{\em Id})\chi^-_s + t_{\mathcal{P}^+(s)}(p^+_s)G_s\right).
\]
\end{theorem}
\begin{proof}
First of all, we write the decomposition of the left regular representation of $\bar U_q$, assigning an index to the multiple factors: 
$$ \bar U_q = \bigoplus_{s=1}^{p-1}\left(\bigoplus_{j=0}^{s-1}\mathcal{P}_j^+(s) \oplus \mathcal{P}_j^-(s)\right) \oplus \bigoplus_{j=0}^{p-1}\mathcal{X}_j^+(p) \oplus \mathcal{X}_j^-(p). $$
Thus, since $t$ is a trace: 
\begin{align*}
t_{\bar U_q}(\rho_a) = &\sum_{s=1}^{p-1}\left(\sum_{j=0}^{s-1}t_{\mathcal{P}_j^+(s)}\left(p_{\mathcal{P}_j^+(s)} \circ \rho_a \circ i_{\mathcal{P}_j^+(s)}\right) + t_{\mathcal{P}_j^-(s)}\left(p_{\mathcal{P}_j^-(s)} \circ \rho_a \circ i_{\mathcal{P}_j^-(s)}\right)\right)\\
& + \sum_{j=0}^{p-1}t_{\mathcal{X}_j^+(p)}\left(p_{\mathcal{X}_j^+(p)} \circ \rho_a \circ i_{\mathcal{X}_j^+(p)}\right) + t_{\mathcal{X}_j^-(p)}\left(p_{\mathcal{X}_j^-(p)} \circ \rho_a \circ i_{\mathcal{X}_j^-(p)}\right).
\end{align*}
Consider the following composite maps for $1 \leq s \leq p-1$ (note that the blocks appear because $\rho_a$ is the right multiplication by $a$): 
\begin{align*}
&h^+_{s,j,a} : \mathcal{P}^{+}(s) \overset{I^+_{s,j}}{\longrightarrow} \mathcal{P}^+_j(s) \overset{i_{\mathcal{P}^+_j(s)}}{\longrightarrow} Q(s) \overset{\rho_a}{\longrightarrow} Q(s) \overset{p_{\mathcal{P}_j^+(s)}}{\longrightarrow} \mathcal{P}^+_j(s)  \overset{\left(I^+_{s,j}\right)^{-1}}{\longrightarrow} \mathcal{P}^+(s),\\
&h^-_{s,j,a} : \mathcal{P}^{-}(s) \overset{I^-_{s,j}}{\longrightarrow} \mathcal{P}^-_j(s) \overset{i_{\mathcal{P}^-_j(s)}}{\longrightarrow} Q(p-s) \overset{\rho_a}{\longrightarrow} Q(p-s) \overset{p_{\mathcal{P}_j^-(s)}}{\longrightarrow} \mathcal{P}^-_j(s)  \overset{\left(I^-_{s,j}\right)^{-1}}{\longrightarrow} \mathcal{P}^-(s),
\end{align*}
where $I_{s,j}^+$ and $I^-_{s,j}$ are the isomorphisms defined by (see Proposition \ref{baseBloc}):
\begin{align*}
&I_{s,j}^+(b^+_i(s)) = B^{++}_{ij}(s),\: I_{s,j}^+(x^+_i(s))=X^{-+}_{ij}(s),\: I_{s,j}^+(y^+_i(s))=Y^{-+}_{ij}(s),\: I_{s,j}^+(a^+_i(s)) = A^{++}_{ij}(s),\\
&I_{s,j}^-(b^-_i(s)) = B^{--}_{ij}(p-s),\: I_{s,j}^-(x^-_i(s))=X^{+-}_{ij}(p-s),\: I_{s,j}^-(y^-_i(s))=Y^{+-}_{ij}(p-s),\\
& I_{s,j}^-(a^-_i(s))=A^{--}_{ij}(p-s).
\end{align*}
For $s = p$, consider:
\begin{align*}
&h^+_{p,j,a} : \mathcal{X}^{+}(p) \overset{I^+_{p,j}}{\longrightarrow} \mathcal{X}^+_j(p) \overset{i_{\mathcal{X}^+_j(p)}}{\longrightarrow} Q(p) \overset{\rho_a}{\longrightarrow} Q(p) \overset{p_{\mathcal{X}_j^+(p)}}{\longrightarrow} \mathcal{X}^+_j(p)  \overset{\left(I^+_{p,j}\right)^{-1}}{\longrightarrow} \mathcal{X}^+(p),\\
&h^-_{p,j,a} : \mathcal{X}^{-}(p) \overset{I^-_{p,j}}{\longrightarrow} \mathcal{X}^-_j(p) \overset{i_{\mathcal{X}^-_j(p)}}{\longrightarrow} Q(0) \overset{\rho_a}{\longrightarrow} Q(0) \overset{p_{\mathcal{X}_j^-(p)}}{\longrightarrow} \mathcal{X}^-_j(p)  \overset{\left(I^-_{p,j}\right)^{-1}}{\longrightarrow} \mathcal{X}^-(p)
\end{align*}
where $I^+_{p,j}$ and $I^-_{p,j}$ are the isomorphisms defined by (see Proposition \ref{baseBloc}):
$$ I^+_{p,j}(v^+_i(p)) = A^{++}_{ij}(p) \:\:\: \text{ and } \:\:\:  I^-_{p,j}(v^-_i(p)) = A^{--}_{ij}(0). $$
Then for $1 \leq s \leq p-1$:
$$t_{\mathcal{P}_j^{\alpha}(s)}\!\left(p_{\mathcal{P}_j^{\alpha}(s)} \circ \rho_a \circ i_{\mathcal{P}_j^{\alpha}(s)}\right) = t_{\mathcal{P}^{\alpha}(s)}\!\left(h^{\alpha}_{s,j,a}\right)$$
and for $s = p$:
$$t_{\mathcal{X}_j^{\alpha}(p)}\!\left(p_{\mathcal{X}_j^{\alpha}(p)} \circ \rho_a \circ i_{\mathcal{X}_j^{\alpha}(p)} \right) = t_{\mathcal{X}^{\alpha}(p)}\!\left(h^{\alpha}_{p,j,a}\right).$$

\noindent We must determine the endomorphism $h^{\alpha}_{s,j,a}$ when $a$ is replaced by the elements given in Proposition \ref{propArike}. Using (\ref{actionWs}), we get:
$$ \forall\, s' \neq s, \forall \, j, \:\: h^{\pm}_{s',j,w^+_s} = 0 \:\:\:\text{ and } \:\:\: h^{-}_{s,j,w^+_s} = 0 $$
and:
$$\forall\, j, \:\: h^{+}_{s,j,w^+_s} = p^+_s.$$
Since this does not depend on $j$ and since the block $Q(s)$ contains $s$ copies of $\mathcal{P}^+(s)$, we find that $t_{\bar U_q}(\rho_{w^+_s}) = st_{\mathcal{P}^+(s)}(p^+_s)$. So by Proposition \ref{propArike}, the coefficient of $G_s$ is $t_{\mathcal{P}^+(s)}(p^+_s)$.

Next, assume that $1 \leq s \leq p-1$, and let us compute $h^{\alpha}_{s',j,\Phi^+_{s-1}e_s}$. By (\ref{actionEs}), we see that 
$$ \forall\, s' \not\in\{s, p-s\}, \forall \, j, \:\: h^{\pm}_{s',j,\Phi^+_{s-1}e_s} = 0 \:\:\: \text{ and } \:\:\: \forall\, j, \:\: h^-_{s,j,\Phi^+_{s-1}e_s}=0,\, h^+_{p-s,j,\Phi^+_{s-1}e_s}=0. $$
Then, Proposition \ref{baseBloc} together with Lemma \ref{lemmaPoids} gives:
$$ \forall \, j, \:\:  h^-_{p-s,j,\Phi^+_{s-1}e_s}=0 \:\:\:\text{ and } \:\:\:  \forall \, 0 \leq j \leq s-2, \:\: h^+_{s,j,\Phi^+_{s-1}e_s}=0  \:\:\: \text{ and } \:\:\: h^+_{s,s-1,\Phi^+_{s-1}e_s}=\text{Id}. $$
It follows that $t_{\bar U_q}\left(\rho_{\Phi^+_{s-1}e_s}\right) = t_{\mathcal{P}^+(s)}(\text{Id})$. So by Proposition \ref{propArike}, the coefficient of $\chi^+_s$ is $t_{\mathcal{P}^+(s)}(\text{Id})$.\\

\noindent We now consider $h^{\pm}_{s',j,\Phi^-_{s-1}e_{p-s}}$. This time, (\ref{actionEs}) shows that
$$ \forall\, s' \not\in\{s, p-s\}, \forall \, j, \:\: h^{\pm}_{s',j,\Phi^-_{s-1}e_{p-s}} = 0 \:\:\: \text{ and } \:\:\: \forall\, j, \:\: h^-_{p-s,j,\Phi^-_{s-1}e_{p-s}}=0,\, h^+_{s,j,\Phi^-_{s-1}e_{p-s}}=0. $$
Then, Proposition \ref{baseBloc} together with Lemma \ref{lemmaPoids} gives:
$$ \forall \, j, \:\:  h^+_{p-s,j,\Phi^-_{s-1}e_{p-s}}=0 \:\:\: \text{ and } \:\:\: \forall \, 0 \leq j \leq s-2, \:\: h^-_{s,j,\Phi^-_{s-1}e_{p-s}}=0  \:\:\: \text{ and } \:\:\: h^-_{s,s-1,\Phi^-_{s-1}e_{p-s}}=\text{Id}. $$
It follows that $t_{\bar U_q}\left(\rho_{\Phi^-_{s-1}e_{p-s}}\right) = t_{\mathcal{P}^-(s)}(\text{Id})$. So by Proposition \ref{propArike}, the coefficient of $\chi^-_s$ is $t_{\mathcal{P}^-(s)}(\text{Id})$.\\
Finally, in the case where $s=p$:
$$ \forall\, s' \neq p, \forall\, j, \:\: h^{\pm}_{s', j, \Phi^+_{p-1}e_{p}} = 0 \:\:\: \text{ and } \:\:\: h^-_{p,j,\Phi^+_{p-1}e_p} = 0.$$
Then, Proposition \ref{baseBloc} together with Lemma \ref{lemmaPoids} gives:
$$ \forall \, 0 \leq j \leq p-2, \:\: h^+_{p,j,\Phi^+_{p-1}e_p}=0  \:\:\:\text{ and }\:\:\: h^+_{p,p-1,\Phi^+_{p-1}e_p}=\text{Id}. $$
It follows that $t_{\bar U_q}\left(\rho_{\Phi^+_{p-1}e_p}\right) = t_{\mathcal{X}^+(p)}(\text{Id})$. So by Proposition \ref{propArike}, the coefficient of $\chi^+_p$ is $t_{\mathcal{X}^+(p)}(\text{Id})$
One similarly gets the coefficient of $\chi^-_p$.
\end{proof}

\indent By Proposition \ref{propArike}, the coefficient of $G_s$ is also given by: $\frac{1}{p-s}t_{\bar U_q}(\rho_{w^-_s})$. Taking back the notations of the proof above, we see using (\ref{actionWs}) that
$$ \forall\, s' \neq p-s, \forall \, j, \:\: h^{\pm}_{s',j,w^-_s} = 0 \:\:\: \text{ and } \:\:\: h^{+}_{p-s,j,w^-_s} = 0 $$
and:
$$\forall\, j, \:\: h^{-}_{p-s,j,w^-_s} = p^-_{p-s}.$$
Since this does not depend on $j$ and since the block $Q(s)$ contains $p-s$ copies of $\mathcal{P}^-(p-s)$, we find that $t_{\bar U_q}(\rho_{w^-_s}) = (p-s)t_{\mathcal{P}^-(p-s)}(p^-_{p-s})$. So by Proposition \ref{propArike}, the coefficient of $G_s$ is $t_{\mathcal{P}^-(p-s)}(p^-_{p-s})$. We thus have: 
\begin{equation}\label{tracePplusMoins}
t_{\mathcal{P}^-(p-s)}(p^-_{p-s}) = t_{\mathcal{P}^+(s)}(p^+_{s}).
\end{equation}
Note that there is an elementary way to see this. Indeed, the morphisms $P^+_s$ and $\bar P^-_{p-s}$ defined in (\ref{morphismesP}) satisfy:
$$ \bar P^-_{p-s} \circ P^+_s = p^+_s, \:\:\: P^+_s \circ \bar P^-_{p-s} = p^-_{p-s}. $$
Hence, we recover (\ref{tracePplusMoins}) by property of the traces. From this, we deduce the following corollary.

\begin{corollary}
Let
$$\varphi = \sum_{s=1}^{p} \left(\lambda^{+}_s\chi^{+}_s + \lambda^{-}_s\chi^{-}_s\right) + \sum_{s'=1}^{p-1}\mu_{s'} G_{s'} \in \text{\em SLF}\left(\bar U_q\right).$$
Then the trace $t^{\varphi} = \Pi_{\bar U_q}^{-1}(\varphi)$ associated to $\varphi$ is given by:
$$ t^{\varphi}_{\mathcal{X}^{\pm}(p)}(\text{\em Id}) = \lambda^{\pm}_p,\:\:\: t^{\varphi}_{\mathcal{P}^{\pm}(s)}(\text{\em Id}) = \lambda^{\pm}_s, \:\:\: t^{\varphi}_{\mathcal{P}^+(s')}(p^+_{s'})=t^{\varphi}_{\mathcal{P}^-(p-s')}(p^-_{p-s'})=\mu_{s'}. $$
\end{corollary}

\subsection{Symmetric linear form corresponding to the modified trace on $\Proj_{\bar U_q}$}\label{ModTrace}
\indent Let $H$ be a finite dimensional Hopf algebra. Let us recall that a {\em modified trace} $\mathsf{t}$ on $\Proj_H$ is a trace which satisfies the additional property that for $U \in \Proj_{H}$, for each $H$-module $V$ and for $f \in \End_{H}(U \otimes V)$  we have:
\begin{equation*}
\mathsf{t}_{U \otimes V}(f) = \mathsf{t}_U(\mathrm{tr}_R(f))
\end{equation*}
where $\mathrm{tr}_R = \mathrm{Id} \otimes \mathrm{tr}_q$ is the right partial quantum trace (see \cite[(3.2.2)]{GKP}). These modified traces are actively studied, having for motivation the construction of invariants in low dimensional topology. We refer to \cite{GKP} for the general theory in a categorical framework which encapsulates the case of $\Proj_H$.
\medskip\\
\indent In \cite{BBGe}, it is shown that there exists a unique up to scalar modified trace $\mathsf{t} = (\mathsf{t}_U)$ on $\Proj_{\bar U_q}$. Uniqueness comes from the fact that $\mathcal{X}^+(p)$ is both a simple and a projective module. The values of this trace are given by:
\begin{equation*}
\begin{array}{lll}
\mathsf{t}_{\mathcal{X}^+(p)}(\text{Id}) = (-1)^{p-1}, & \mathsf{t}_{\mathcal{X}^-(p)}(\text{Id}) = 1, & \mathsf{t}_{\mathcal{P}^+(s)}(\text{Id}) = (-1)^s(q^s+q^{-s}),\\
 \mathsf{t}_{\mathcal{P}^-(s)}(\text{Id}) = (-1)^{p-s-1}(q^s+q^{-s}), & \mathsf{t}_{\mathcal{P}^+(s)}(p^+_s) = (-1)^s[s]^2 & \mathsf{t}_{\mathcal{P}^-(s)}(p^-_s) = \mathsf{t}_{\mathcal{P}^+(p-s)}(p^+_{p-s}).
\end{array}
\end{equation*}
\indent Let $H$ be a finite dimensional unimodular pivotal Hopf algebra with pivotal element $g$ and let $\mu^r \in H^*$ be a right integral on $H$, which means that
\begin{equation*}
\forall \, x \in H, \:\: (\mu^r \otimes \text{Id})(\Delta(x)) = \mu(x)1.
\end{equation*}
Recall from \eqref{integraleShifteSLF} that $\mu^r(g \cdot)$ is a symmetric linear form. In the recent paper \cite{BBG}, it is shown that modified traces on $\Proj_H$ are unique up to scalar, and that the corresponding symmetric linear forms are scalar multiples of $\mu^r(g\cdot)$. Here, we show how Theorem \ref{MainResult} and computations made in \cite{GT} (see also \cite{arike}) and \cite{FGST} quickly allow us to recover this result in the case of $H = \bar U_q$. First, recall that right integrals $\mu^r_{\zeta}$ of $\bar U_q$ are given by:
\begin{equation*}
\mu^r_{\zeta}(F^mE^nK^j) = \zeta \delta_{m, p-1}\delta_{n,p-1}\delta_{j,p+1},
\end{equation*}
where $\zeta$ is an arbitrary scalar. Hence:
\begin{equation*}
\mu^r_{\zeta}(K^{p+1} F^mE^nK^j) = \zeta \delta_{m, p-1}\delta_{n,p-1}\delta_{j,0}.
\end{equation*}
Using formulas given in \cite{GT} (see also \cite{arike}\footnote{In notations of \cite{arike}, we have $e_s = \sum_{t=1}^se^+(s,t) + \sum_{u=1}^{p-s}e^-(p-s,u)$.}), we have ($1 \leq s \leq p-1$):
\begin{align*}
e_0 &= \frac{(-1)^{p-1}}{2p[p-1]!^2}\sum_{t=0}^{p-1}\sum_{l=0}^{2p-1}q^{-(-2t-1)l}F^{p-1}E^{p-1}K^l + (\text{terms of lower degree in } E \text{ and } F),\\
e_s &= \alpha_s\sum_{t=0}^{p-1}\sum_{l=0}^{2p-1}q^{-(s-2t-1)l}F^{p-1}E^{p-1}K^l + (\text{terms of lower degree in } E \text{ and } F),\\
e_p &= \frac{1}{2p[p-1]!^2}\sum_{t=0}^{p-1}\sum_{l=0}^{2p-1}q^{-(p-2t-1)l}F^{p-1}E^{p-1}K^l + (\text{terms of lower degree in } E \text{ and } F),
\end{align*}
where $\alpha_s$ is given in the last page of \cite{arike} as:
$$ \alpha_s = -\frac{(-1)^{p-s-1}}{2p[p-s-1]!^2[s-1]!^2}\left(\sum_{l=1}^{s-1}\frac{1}{[l][s-l]} - \sum_{l=1}^{p-s-1}\frac{1}{[l][p-s-l]}\right). $$
In order to simplify this, it is observed in \cite[Proof of Proposition 2]{murakami}, that
$$ \sum_{l=1}^{s-1}\frac{1}{[l][s-l]} - \sum_{l=1}^{p-s-1}\frac{1}{[l][p-s-l]} = \frac{-(q^s + q^{-s})}{[s]^{2}}. $$
So, since:
$$ [p-s-1]!^2[s-1]!^2 = \frac{[p-1]!^2}{[s]^2}, $$
we get:
\begin{equation*}
\alpha_s = \frac{(-1)^{p-s-1}}{2p[p-1]!^2}(q^s+q^{-s}).
\end{equation*}
Using formulas given in \cite{FGST} (see also \cite[Prop. II.3.19]{ibanez}), we have:
\begin{align*}
w_s^+ &= \frac{(-1)^{p-s-1}}{2p[p-1]!^2}[s]^2sF^{p-1}E^{p-1} + (\text{other monomials}),\\
w_s^- &= \frac{(-1)^{p-s-1}}{2p[p-1]!^2}[s]^2 (p-s)F^{p-1}E^{p-1} + (\text{other monomials}).
\end{align*}
We now use Proposition \ref{propArike} to get the coefficients of $\mu^r_{\zeta}(K^{p+1}\cdot)$ in the GTA basis. For instance:
\begin{align*}
\frac{\mu^r_{\zeta}(K^{p+1}w^+_s)}{s} &= \zeta \frac{(-1)^{p-s-1}}{2p[p-1]!^2}[s]^2,\\
\mu^r_{\zeta}(K^{p+1}\Phi^+_{s-1}e_s) &= \frac{\alpha_s}{2p}\mu^r_{\zeta}\left(K^{p+1}F^{p-1}E^{p-1}\sum_{t=0}^{p-1}\sum_{l,j=0}^{2p-1} q^{-(s-1)(l+j)+2tl}K^{l+j}\right)\\
& = \zeta\frac{\alpha_s}{2p}\sum_{t=0}^{p-1}\sum_{l=0}^{2p-1} q^{2tl} = \zeta\alpha_s.
\end{align*}
Choose the normalization factor to be $\zeta = (-1)^{p-1}2p[p-1]!^2$, and let $\mu^r$ be the so-normalized integral. Then:
\begin{equation}\label{traceModifieeGTA}
\begin{split}
\mu^r(K^{p+1}\cdot) = (-1)^{p-1}\chi^+_p + \chi^-_p + \sum_{s=1}^{p-1} & \left( (-1)^s(q^s+q^{-s})\chi^+_s + (-1)^{p-s-1}(q^s+q^{-s})\chi^-_s \right. \\
& \left. + \: (-1)^s[s]^2G_s \right).
\end{split}
\end{equation}
By Theorem \ref{MainResult}, we recover $\Pi_{\bar U_q}(\mathsf{t}) = \mu^r(K^{p+1}\cdot)$.

\section{Multiplication rules in the GTA basis}\label{multiplication}
\indent We mentioned in section \ref{sectionSLF} that $\SLF(\bar U_q)$ is a commutative algebra. In this section, we address the problem of the decomposition in the GTA basis of the product of two elements in this basis. The resulting formulas are surprisingly simple.
\smallskip\\
\indent Let us start by recalling some facts. For every $\bar U_q$-module $V$, we define the character of $V$ as (see (\ref{defMatriceT}) for the definition of $\overset{V}{T}$):
$$ \chi^V = \text{tr}(\overset{V}{T}). $$
This splits on extensions:
\begin{equation*}
0 \rightarrow V \rightarrow M \rightarrow W \rightarrow 0 \: \text{ exact} \:\:\: \implies \:\:\: \chi^M = \chi^V + \chi^W.
\end{equation*}
Due to the fact that $\bar U_q$ is finite dimensional, every finite dimensional $\bar U_q$-module has a composition series (\textit{i.e.} is constructed by successive extensions by simple modules). It follows that every $\chi^V$ can be written as a linear combination of the $\chi^{\alpha}_s = \chi^{\mathcal{X}^{\alpha}(s)}$. Moreover, we see by definition of the product on $\bar U_q^*$ that
\begin{equation}\label{Ttens}
\overset{V \otimes W}{T} = \overset{V}{T}_1 \overset{W}{T}_2
\end{equation}
where $\overset{V}{T}_1 = \overset{V}{T} \otimes I_{\dim(W)}$ and $\overset{W}{T}_2 = I_{\dim(V)} \otimes \overset{W}{T}$. Thus $\chi^{V \otimes W} = \chi^V\chi^W$. Hence multiplying two $\chi$'s is equivalent to tensoring two simples modules and finding the decomposition into simple factors. This means that
\begin{equation*}
\vect(\chi^{\alpha}_s)_{\alpha \in \{\pm\}, 1 \leq s \leq p} \overset{\sim}{\rightarrow} \mathfrak{G}(\bar U_q) \otimes_{\mathbb{Z}} \mathbb{C}, \:\:\: \chi^I \mapsto [I]
\end{equation*}
where $\mathfrak{G}(\bar U_q)$ is the Grothendieck ring of $\bar U_q$. By \cite{FGST}, we know the structure of $\mathfrak{G}(\bar U_q)$. Recall the decomposition formulas (with $2 \leq s \leq p-1$):
\begin{equation*}
\mathcal{X}^-(1) \otimes \mathcal{X}^{\alpha}(s) \cong \mathcal{X}^{-\alpha}(s), \:\:\: \mathcal{X}^+(2) \otimes \mathcal{X}^{\alpha}(s) \cong \mathcal{X}^{\alpha}(s-1) \oplus \mathcal{X}^{\alpha}(s+1), \:\:\: \mathcal{X}^+(2) \otimes \mathcal{X}^{\alpha}(p) \cong \mathcal{P}^{\alpha}(p-1) 
\end{equation*}
so that
\begin{equation}\label{produitChi}
\chi^-_1\chi^{\alpha}_s = \chi^{-\alpha}_s, \:\:\: \chi^+_2 \chi^{\alpha}_s = \chi^{\alpha}_{s-1} + \chi^{\alpha}_{s+1}, \:\:\: \chi^{+}_2\chi^{\alpha}_p = 2\chi^{\alpha}_{p-1} + 2\chi^{-\alpha}_1.
\end{equation}
We see in particular that $\chi^+_2$ generates the subalgebra $\vect(\chi^{\alpha}_s)_{\alpha \in \{\pm\}, 1 \leq s \leq p}$. The $\chi^{\alpha}_s$ are expressed as Chebyschev polynomials of $\chi^+_2$, see \cite[section 3.3]{FGST} for details.

\begin{theorem}\label{ProduitArike}
The multiplication rules in the GTA basis are entirely determined by (\ref{produitChi}) and by the following formulas:
\begin{align}
&\chi^+_2G_1 = [2]G_2, \label{ChiG1}\\
&\chi^+_2G_s = \frac{[s-1]}{[s]}G_{s-1} + \frac{[s+1]}{[s]}G_{s+1} \:\text{ for } 2 \leq s \leq p-2,\label{ChiG}\\
&\chi^+_2G_{p-1} = [2]G_{p-2}, \label{ChiGpMoins1}\\
&\chi^-_1G_s = -G_{p-s} \:\text{ for all } s,\label{ChiMoinsGs}\\
&G_sG_t = 0 \:\text{ for all } s, t.\label{GG}
\end{align}
\end{theorem}
Before giving the proof, let us deduce a few consequences.
\begin{corollary}\label{coroVermaIdeal} For all $1 \leq s \leq p-1$ we have:
$$ G_s = \frac{1}{[s]}\chi^+_sG_1, \:\:\:\: \chi^+_pG_1 = 0. $$
It follows that $(\chi^+_s + \chi^-_{p-s})G_t = 0$, and that $\mathcal{P} = \vect(\chi^+_s + \chi^-_{p-s}, \chi^+_p, \chi^-_p)_{1 \leq s \leq p-1}$ is an ideal of $\SLF(\bar U_q)$.
\end{corollary}
\begin{proof}[Proof of Corollary \ref{coroVermaIdeal}]
The formulas for $\chi^+_sG_1$ are proved by induction using $\chi^+_{s+1}= \chi^+_s\chi^+_2 - \chi^+_{s-1}$ together with formula (\ref{ChiG}). We deduce:
$$ (\chi^+_s + \chi^-_{p-s})G_t = \frac{\chi^+_t}{[t]}(\chi^+_sG_1 + \chi^-_{p-s}G_1) = \frac{\chi^+_t}{[t]}([s]G_s + [s]\chi^-_1G_{p-s}) = 0 .$$
It is straightforward that $\mathcal{P}$ is stable by multiplication by $\chi^+_2$, so it is an ideal.
\end{proof}
\begin{remark}
We have $\chi^{\mathcal{P}^{\alpha}(s)} = 2\left(\chi^{\alpha}_s + \chi^{-\alpha}_{p-s}\right)$ for $1 \leq s \leq p-1$. Thus $\mathcal{P}$ is generated by characters of the projective modules. It is well-known that if $H$ is a finite dimensional Hopf algebra, then the full subcategory of finite dimensional projective $H$-modules is a tensor ideal. Thus we can deduce without any computation that $\mathcal{P}$ is stable under the multiplication by every $\chi^{I}$.
\finEx
\end{remark}

\indent We now proceed with the proof of the theorem. Observe that we cannot apply Proposition \ref{propArike} to show it since we do not know expressions of $\Delta(e_s)$ and $\Delta(w^{\pm}_s)$ which are easy to evaluate in the GTA basis. Recall (\cite{KS}, see also \cite{ibanez}) the following fusion rules:
\begin{align}
&\mathcal{X}^-(1) \otimes \mathcal{P}^{\alpha}(s) \cong \mathcal{P}^{-\alpha}(s)  \:\:\: \text{ for all } s, \label{chi1MoinsTensPs}\\
&\mathcal{X}^+(2) \otimes \mathcal{P}^{\alpha}(1) \cong 2\mathcal{X}^{-\alpha}(p) \oplus \mathcal{P}^{\alpha}(2),\label{chi2TensP1}\\
&\mathcal{X}^+(2) \otimes \mathcal{P}^{\alpha}(s) \cong \mathcal{P}^{\alpha}(s-1) \oplus \mathcal{P}^{\alpha}(s+1) \:\:\: \text{ for } 2 \leq s \leq p-1,\label{chi2TensPs}\\
&\mathcal{X}^+(2) \otimes \mathcal{P}^{\alpha}(p-1) \cong 2\mathcal{X}^{\alpha}(p) \oplus \mathcal{P}^{\alpha}(p-2)\label{chi2TensPpMoins1}.
\end{align}
They imply the following key lemma.
\begin{lemma}\label{lemmaProdChiGs}
There exist scalars $\gamma_s, \beta_s, \lambda_s, \eta_s, \delta_s$ such that
\begin{align*}
&\chi^+_2G_s = \beta_s G_{s-1} + \gamma_s G_{s+1} + \lambda_s\!\left( \chi^+_{s-1} + \chi^-_{p-s+1} - \chi^+_{s+1} - \chi^-_{p-s-1} \right) \:\:\:\: (\text{for } 2 \leq s \leq p-2),\\
&\chi^+_2G_1 = \gamma_1 G_{2} + \lambda_1\!\left(\chi^-_p - \chi^+_2 - \chi^-_{p-2} \right), \:\:\:\: \chi^+_2G_{p-1} = \beta_{p-1} G_{p-2} + \lambda_{p-1}\!\left( \chi^+_{p-2} + \chi^-_{2} - \chi^+_p \right),\\
&\chi^-_1G_s = \eta_s G_{p-s} + \delta_s\!\left( \chi^+_{p-s} + \chi^-_s \right).
\end{align*}
\end{lemma}
\begin{proof}
Let us fix $2 \leq s \leq p-2$; by (\ref{defChis}), (\ref{defGs}), (\ref{Ttens}) and (\ref{chi2TensPs}) we have:
\begin{align*}
\chi^+_2G_s &\in \vect\!\left(\overset{\mathcal{X}^{+}(2)}{T_{ij}}\cdot\overset{\mathcal{P}^{+}(s)}{T_{kl}}, \overset{\mathcal{X}^{+}(2)}{T_{ij}}\cdot\overset{\mathcal{P}^{-}(p-s)}{T_{kl}}\right)_{ijkl} = \vect\!\left(\overset{\mathcal{X}^{+}(2) \otimes \mathcal{P}^+(s)}{T_{ijkl}}, \overset{\mathcal{X}^{+}(2) \otimes \mathcal{P}^-(p-s)}{T_{ijkl}}\right)_{ijkl}\\
&= \vect\!\left(\overset{\mathcal{P}^+(s-1)}{T_{ij}}, \overset{\mathcal{P}^+(s+1)}{T_{ij}}, \overset{\mathcal{P}^-(p-s+1)}{T_{ij}}, \overset{\mathcal{P}^-(p-s-1)}{T_{ij}}\right)_{ij}
\end{align*}
where $\overset{V}{T}\!_{ij}$ is the matrix element at the $i$-th row and $j$-th column of the representation matrix $\overset{V}{T}$ and $\overset{V \otimes W}{T}\!\!\!\!_{ijkl}$ is the matrix element at the $(i,j)$-th row and $(k,l)$-th column of the representation matrix $\overset{V \otimes W}{T}$. Hence, since $\chi^+_2G_s$ is symmetric, it is necessarily of the form 
$$
\chi^+_2G_s = \beta_s G_{s-1} + \gamma_s G_{s+1} + z_1 \chi^+_{s-1} + z_2 \chi^+_{s+1} + z_3 \chi^-_{p-s+1} + z_4 \chi^-_{p-s-1}.
$$
Evaluating this equality on $K$ and $K^2$, we find (since $G_t(K^{l}) = 0$ for all $t$ and $l$):
\[ [s-1](z_1 - z_3) + [s+1](z_2 - z_4) = 0, \:\:\:\: [s-1]_{q^2}(z_1 - z_3) + [s+1]_{q^2}(z_2 - z_4) = 0, \]
with $[n]_{q^2} = \frac{q^{2n} - q^{-2n}}{q^2 - q^{-2}}$. The determinant of this linear system with unknowns $z_1 - z_3, z_2 - z_4$ is $\frac{2\sin( (s-1)\pi/p ) \sin( (s+1)\pi/p )}{\sin( \pi/p )\sin( 2\pi/p )}\left( \cos( (s+1)\pi/p ) - \cos( (s-1)\pi/p ) \right) \neq 0$.  Hence $z_1 = z_3$, $z_2 = z_4$. Moreover, evaluating the above equality on $1$, we find $p(z_1 + z_2) = 0$. Letting $\lambda_s = z_1$, the result follows. The other formulas are obtained in a similar way using (\ref{chi1MoinsTensPs}), (\ref{chi2TensP1}) and (\ref{chi2TensPpMoins1}).
\end{proof}

\indent We will use the Casimir element $C$ (see \eqref{defCasimir}) to make computations easier. Observe that, due to \eqref{casimirBaseCan} and \eqref{actionCentreSLF}, we have
\begin{equation}\label{precCasimir}
\forall \, x \in \bar U_q, \:\:\: \chi^{\alpha}_s(Cx) = \alpha c_s \chi^{\alpha}_s(x), \:\:\: G_s(Cx) = c_sG_s(x) + (\chi^+_s + \chi^-_{p-s})(x). 
\end{equation}
where $\displaystyle c_s = \frac{q^s + q^{-s}}{ (q - q^{-1})^2 }$. Then by induction we get $G_s(C^n) = npc_s^{n-1}$ for $n \geq 1$. We will also denote $\displaystyle c_K = \frac{qK + q^{-1}K^{-1}}{(q - q^{-1})^2}$.

\begin{proof}[Proof of Theorem \ref{ProduitArike}]
~\\
\indent \textbullet \hspace{1pt} \textit{Formula \eqref{ChiG}.} \hspace{2pt} We first evaluate the corresponding formula of Lemma \ref{lemmaProdChiGs} on $FE$. It holds $G_t(FE) = G_t(C) = p$, $(\chi^+_t + \chi^-_{p-t})(FE) = (\chi^+_t + \chi^-_{p-t})(C) = p c_t$ for all $t$ and $\chi^+_2G_s(FE) = \chi^+_2(K^{-1})G_s(FE) = [2]p$. Thus we get:
\begin{equation}\label{eqLin1}
\beta_s + \gamma_s  + (c_{s-1} - c_{s+1}) \lambda_s = \beta_s + \gamma_s  - [s] \lambda_s = [2].
\end{equation}
Next, we evaluate the formula of Lemma \ref{lemmaProdChiGs} on $(FE)^2$. On the one hand,
\begin{align*}
(\chi^+_2 G_s)\!\left( (FE)^2 \right) &= \chi^+_2(K^{-2}) G_s\!\left( (FE)^2 \right) = \chi^+_2(K^{-2}) G_s\!\left( C^2 - 2 C c_K + c_K^2 \right)\\
& = \chi^+_2(K^{-2}) G_s(C^2) = 2p(q^2 + q^{-2})c_s.
\end{align*}
For the first equality, we used that $\varphi(E^iF^jK^{l}) = \delta_{i,j}\varphi(E^iF^iK^{l})$ for all $\varphi \in \mathrm{SLF}(\bar U_q)$, that $G_s(K^{l}) = 0$ and that $G_s(FEK^{l}) = 0$ for $1 \leq l \leq p-1$. The third equality is due to \eqref{precCasimir} and to the fact that $(\chi^+_s + \chi^-_{p-s})(K^{l}) = 0$ for $1 \leq l \leq p-1$. 
On the other hand, using again the Casimir element,
\begin{align*}
&\beta_s G_{s-1}\!\left( (FE)^2 \right) + \gamma_s G_{s+1}\!\left( (FE)^2 \right) + \lambda_s\!\left( \chi^+_{s-1} + \chi^-_{p-s+1} - \chi^+_{s+1} - \chi^-_{p-s-1} \right)\!\left( (FE)^2 \right)\\
&= \beta_s G_{s-1}\!\left( C^2 \right) + \gamma_s G_{s+1}\!\left( C^2 \right) + \lambda_s\!\left( \chi^+_{s-1} + \chi^-_{p-s+1} - \chi^+_{s+1} - \chi^-_{p-s-1} \right)\!\left( C^2 \right)\\
&= 2pc_{s-1} \beta_s + 2pc_{s+1} \gamma_s + p(c_{s-1}^2 - c_{s+1}^2) \lambda_s.
\end{align*}
Since $c_{s-1}^2 - c_{s+1}^2 = -(q + q^{-1})c_s[s]$, we get
\begin{equation}\label{eqLin2}
2c_{s-1} \beta_s + 2c_{s+1} \gamma_s -(q + q^{-1})c_s[s]\lambda_s = 2(q^2 + q^{-2})c_s.
\end{equation}
In order to get a third linear equation between $\beta_s$, $\gamma_s$ and $\lambda_s$, we use evaluation on $E^{p-1}F^{p-1}$. This has the advantage to annihilate all the $\chi^{\alpha}_t$ appearing in the formula of Lemma \ref{lemmaProdChiGs}. First:
\begin{equation}\label{actionP}
\begin{split}
E^{p-1}F^{p-1}b_0^{\alpha}(s) &= E^{p-1}y^{\alpha}_{p-s-1}(s) = (-\alpha)^{p-s-1}[p-s-1]!^2E^sy_0^{\alpha}(s)\\
&= (-\alpha)^{p-s-1}\alpha^{s-1}[p-s-1]!^2[s-1]!^2a_0^{\alpha}(s) \\
&= (-\alpha)^{p-s-1}\alpha^{s-1}\frac{[p-1]!^2}{[s]^2}a_0^{\alpha}(s)
\end{split}
\end{equation}
and $E^{p-1}F^{p-1}$ annihilates all the other basis vectors. Hence:
$$
G_s(E^{p-1}F^{p-1}) = 2(-1)^{p-s-1}\frac{[p-1]!^2}{[s]^2}.
$$
Next by (\ref{coproduitMonome}), we have:
$$
\chi^+_2 \otimes \text{Id}\left(\Delta(E^{p-1}F^{p-1})\right) = -[2]E^{p-1}F^{p-1} - q^2E^{p-2}F^{p-2}K.
$$
As in (\ref{actionP}), we find:
\begin{align*}
E^{p-2}F^{p-2}Kb_0^{\alpha}(s) &= (-\alpha)^{p-s}\alpha^sq^{s-1}\frac{[p-1]!^2}{[s+1][s]^2}a_0^{\alpha}(s),\\
E^{p-2}F^{p-2}Kb_1^{\alpha}(s) &= (-\alpha)^{p-s-1}\alpha^{s-1}q^{s-3}\frac{[p-1]!^2}{[s-1][s]^2} a_1^{\alpha}(s)
\end{align*}
and all the others basis vectors are annihilated. Hence:
$$
G_s(E^{p-2}F^{p-2}K) = 2(-1)^{p-s-1}\frac{[p-1]!^2}{[s]^2}\frac{q^{-2}[2]}{[s-1][s+1]}.
$$
We obtain:
$$ \chi^+_2 \otimes G_s\left(\Delta(E^{p-1}F^{p-1})\right) = 2(-1)^{p-s}[p-1]!^2\frac{[2]}{[s-1][s+1]} $$
and thus:
\begin{equation}\label{eqLin3}
\frac{\beta_s}{[s-1]^2} +  \frac{\gamma_s}{[s+1]^2} = \frac{[2]}{[s-1][s+1]}.
\end{equation}
As a result, we have a linear system \eqref{eqLin1}--\eqref{eqLin2}--\eqref{eqLin3} between $\beta_s$, $\gamma_s$ and $\lambda_s$. It is easy to check that $\beta_s = \frac{[s-1]}{[s]}, \gamma_s = \frac{[s+1]}{[s]}, \lambda_s = 0$ is a solution. Moreover this solution is unique. Indeed, a straightforward computation reveals that
\[
\det
\left(
\begin{array}{ccc}
1 & 1 & -[s]\\
2 c_{s-1} & 2c_{s+1} & -(q + q^{-1})c_s[s]\\
\frac{1}{[s-1]^2} & \frac{1}{[s+1]^2} & 0
\end{array}
\right)
= \frac{[s]^2}{[s-1]^2} + \frac{[s]^2}{[s+1]^2} > 0.
\]

\indent \textbullet \hspace{1pt} \textit{Formulas \eqref{ChiG1} and \eqref{ChiMoinsGs}.} \hspace{2pt} Evaluating as above the corresponding formulas of Lemma \ref{lemmaProdChiGs} on $FE$ and $(FE)^2$, one gets linear systems with non-zero determinants.  It is then easy to see that $\beta_1 = [2], \lambda_1 = 0$ and $\eta_s = -1, \delta_s = 0$ are the unique solutions of each of these two systems.
\medskip\\
\indent \textbullet \hspace{1pt} \textit{Formula \eqref{ChiGpMoins1}.} \hspace{2pt} It can be deduced from the formulas already shown:
\[ \chi^+_2 G_{p-1} = -\chi^+_2 \chi^-_1 G_1 = -\chi^-_1 [2] G_2 = [2] G_{p-2}. \]
\indent \textbullet \hspace{1pt} \textit{Formula \eqref{GG}.} \hspace{2pt} Recall the isomorphism of algebras $\mathcal{D}$ defined in (\ref{morphismeDrinfeld}). Taking into account that $\varphi(K^iF^mE^n)=0$ if $n \neq m$ for any $\varphi \in \SLF(\bar U_q)$ and that $G_s(K^i) = 0$ for all $i$, and making use of the expression of $RR'$ given in \eqref{RRPrime}, we get:
\begin{align*}
\mathcal{D}(G_s) &= \sum_{n=0}^{p-1}\sum_{j=0}^{2p-1}\left(\sum_{i=0}^{2p-1} \frac{(q-q^{-1})^n}{[n]!^2}q^{n(j-i-1)-ij}G_s(K^{p+i+1}E^nF^n)\right)K^jF^nE^n\\
&= \sum_{n=1}^{p-1}\sum_{j=0}^{2p-1}\lambda_{j,n}K^jF^nE^n
\end{align*}
for some coefficients $\lambda_{j,n}$ (observe that $n \geq 1$). From this it follows that for all $\alpha \in \{\pm\}$ and $1 \leq r \leq p-1$: $ \mathcal{D}(G_s)b_0^{\alpha}(r) \in \mathbb{C}a_0^{\alpha}(r)$. By \eqref{actionEs}, we deduce that $\mathcal{D}(G_s) \in \vect(w^{\pm}_r)_{1 \leq r \leq p-1}$ for all $s$. Thus $\mathcal{D}(G_sG_t) = 0$, thanks to \eqref{produitCentre}.
\end{proof}

\chapter{$\mathcal{L}_{0,1}(H)$, $\mathcal{L}_{1,0}(H)$ and representation of the modular group $\mathrm{SL}_2(\mathbb{Z})$}\label{chapitreTore}
\indent In this chapter, we focus on the surfaces
\[ \Sigma_{0,1}^{\mathrm{o}} = \Sigma_{0,1} \setminus D, \:\:\:\:\:\:\: \Sigma_{1,0}^{\mathrm{o}} = \Sigma_{1,0} \setminus D \]
where $D$ is an embedded open disk, and on the associated algebras $\mathcal{L}_{0,1}(H)$ and $\mathcal{L}_{1,0}(H)$. These algebras deserve a particular interest and are a necessary preliminary step because they are the building blocks of $\mathcal{L}_{g,n}(H)$ (we will see what does this means in Definition \ref{definitionLgn}). 

\smallskip

\begin{figure}[!h]
\centering
\begingroup%
  \makeatletter%
  \providecommand\color[2][]{%
    \errmessage{(Inkscape) Color is used for the text in Inkscape, but the package 'color.sty' is not loaded}%
    \renewcommand\color[2][]{}%
  }%
  \providecommand\transparent[1]{%
    \errmessage{(Inkscape) Transparency is used (non-zero) for the text in Inkscape, but the package 'transparent.sty' is not loaded}%
    \renewcommand\transparent[1]{}%
  }%
  \providecommand\rotatebox[2]{#2}%
  \newcommand*\fsize{\dimexpr\f@size pt\relax}%
  \newcommand*\lineheight[1]{\fontsize{\fsize}{#1\fsize}\selectfont}%
  \ifx\svgwidth\undefined%
    \setlength{\unitlength}{161.1839629bp}%
    \ifx\svgscale\undefined%
      \relax%
    \else%
      \setlength{\unitlength}{\unitlength * \real{\svgscale}}%
    \fi%
  \else%
    \setlength{\unitlength}{\svgwidth}%
  \fi%
  \global\let\svgwidth\undefined%
  \global\let\svgscale\undefined%
  \makeatother%
  \begin{picture}(1,0.78508371)%
    \lineheight{1}%
    \setlength\tabcolsep{0pt}%
    \put(0.68539575,0.64953473){\color[rgb]{0,0,0}\makebox(0,0)[lt]{\lineheight{1.25}\smash{\begin{tabular}[t]{l}$m$\end{tabular}}}}%
    \put(0,0){\includegraphics[width=\unitlength,page=1]{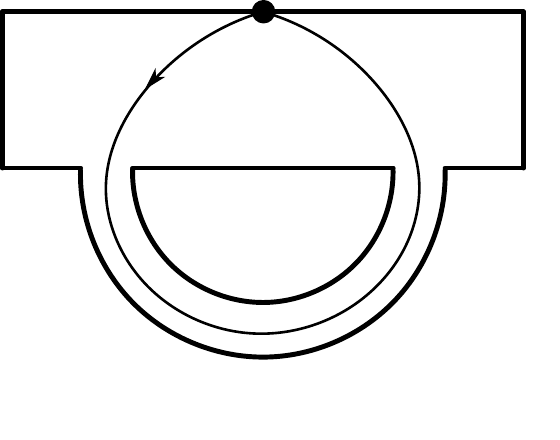}}%
    \put(0.44208621,0.010118){\color[rgb]{0,0,0}\makebox(0,0)[lt]{\lineheight{1.25}\smash{\begin{tabular}[t]{l}$\overset{I}{M}$\end{tabular}}}}%
  \end{picture}%
\endgroup%

\caption{Surface $\Sigma_{0,1}^{\mathrm{o}}$ with basepoint, canonical loop and matrices of generators.}
\label{Sigma01}
\end{figure}

\begin{figure}[!h]
\centering
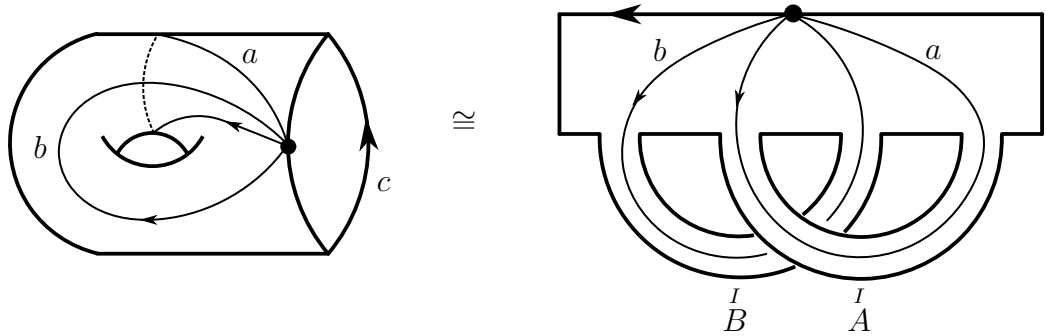
\caption{Two views of the surface $\Sigma_{1,0}^{\mathrm{o}}$, with basepoint, canonical loops and matrices of generators.}
\label{Sigma10}
\end{figure}

\indent Figures \ref{Sigma01} and \ref{Sigma10} are the pictures that one should always keep in mind. We see $\Sigma_{0,1}^{\mathrm{o}}$ and $\Sigma_{1,0}^{\mathrm{o}}$ as thickenings (\textit{i.e.} tubular neighborhoods) of the embedded oriented graphs $\Gamma_{0,1} = \bigl( \{\bullet\}, \{m\} \bigr)$ and $\Gamma_{1,0} = \bigl( \{\bullet\}, \{b,a\} \bigr)$ whose vertex is the basepoint and whose edges are the generators of the fundamental group represented below. To get the second view from the first in Figure \ref{Sigma10}\footnote{Compared to the Figure 1 of \cite{Fai18}, we have done a $180^{\circ}$-rotation around the horizontal axis of $\mathbb{R}^3$, in order to have the handles at the bottom of the Figure. The reason of this change comes from the definition of the graphical calculus and the Wilson loop map in Chapter \ref{chapitreGraphiqueSkein}.}, retract to a tubular neighborhood of the loops $b$ and $a$. Note that with this choice of generators, the boundary loop of $\Sigma_{1,0}^{\mathrm{o}}$ is expressed as 
\begin{equation}\label{boundaryLoop10}
c = ba^{-1}b^{-1}a.
\end{equation}

To each generating loop, or equivalently to each handle, is associated a family of matrices, indexed by the $H$-modules and whose coefficients are generators of the algebra. The defining relations are given in Definitions \ref{defL01} and \ref{definitionL10}, following \cite{AGS, BR, AGS2}. The difference with these original papers is that we directly restrict to the canonical graphs $\Gamma_{0,1}$ and $\Gamma_{1,0}$ and that we do not use Clebsch-Gordan operators to write the fusion relation; indeed, these operators have good properties in the semisimple case only and computations are simpler without using them. Also, the relations may vary from one paper to another due to different choice of $R$-matrices ($R'^{-1}$ instead or $R$) or different labellings of the matrices associated to loops.

\smallskip

\indent The main results of this chapter are
\begin{itemize}
\item The representation of the algebra of invariant elements $\mathcal{L}_{1,0}^{\mathrm{inv}}(H)$ on $\mathrm{SLF}(H)$ (Theorem \ref{repInv}). Note that the matrices $\overset{I}{C}$ (which correspond to the boundary of $\Sigma_{1,0}^{\mathrm{o}}$, see \eqref{boundaryLoop10}) used for the proof of that Theorem already appeared in \cite{alekseev} (with $H = U_q(\mathfrak{g})$), but here we need to generalize and adapt the construction of the representation to our assumptions on $H$.
\item The construction of a projective representation of $\mathrm{SL}_2(\mathbb{Z})$ (mapping class group of the torus $\Sigma_{1,0}$) on $\mathrm{SLF}(H)$ (Theorem \ref{repSL2}). This complete and generalize to a non-semisimple setting the idea of Alekseev--Schomerus \cite{AS}.
\item The equivalence of the representation of $\mathrm{SL}_2(\mathbb{Z})$ with the one constructed by Lyubashenko--Majid in \cite{LM} (Theorem \ref{EquivalenceLMandSLF}).
\item The explicit description of the representation of $\mathrm{SL}_2(\mathbb{Z})$ when $H = \bar U_q(\mathfrak{sl}_2)$ (Theorems \ref{actionSL2ZArike} and \ref{thDecRep}).
\end{itemize}

\indent The material of this chapter is mainly the content of \cite{Fai18}. However here we give more details and comments.

\section{The loop algebra $\mathcal{L}_{0,1}(H)$}\label{sectionLoopL01}
\indent We assume that $H$ is a finite dimensional factorizable Hopf algebra. The ribbon assumption is not needed in this section.
\subsection{Definition of $\mathcal{L}_{0,1}(H)$ and $H$-module-algebra structure}
\indent Let $\mathrm{T}(H^*)$ be the tensor algebra of $H^*$, which by definition is linearly spanned by all the formal products $\varphi_1 \cdots \varphi_n$ (with $n \geq 0$ and $\varphi_i \in H^*$) modulo the obvious multilinear relations. 
There is a canonical injection $j : H^* \to \mathrm{T}(H^*)$ and we define $\overset{I}{M} = j(\overset{I}{T})$, \textit{i.e.} $\overset{I}{M}{^a_b} = j\bigl( \overset{I}{T}{^a_b} \bigr)$.

\begin{definition}\label{defL01}
The loop algebra $\mathcal{L}_{0,1}(H)$ is the quotient of $\mathrm{T}(H^*)$ by the following fusion relations:
$$ \overset{I \otimes J}{M}\!_{12} = \overset{I}{M}_1\overset{IJ}{(R')}_{12}\overset{J}{M}_2\overset{IJ}{(R'^{-1})}_{12} $$
for all finite dimensional $H$-modules $I,J$.
\end{definition}
\noindent The right hand-side of the fusion relation in Definition \ref{defL01} above is the one of \cite[Def 1]{BNR}; the one of \cite[Def 12]{AGS2} and \cite[eq (3.11)]{AS} is different, due to different choices of the action of $H$ on $\mathcal{L}_{0,1}(H)$ and to particular normalization of Clebsch-Gordan operators. Moreover, in the papers \cite{AGS2, AS, BNR} the matrix $\overset{I \otimes J}{M}$ did not appeared, instead it was always decomposed as a sum of the matrices of the irreducible direct summands of $I \otimes J$ thanks to Clebsch-Gordan operators, which is relevant in the semisimple case only; in \cite{S}, $I$ and $J$ are restricted to the regular representation and the matrix $\overset{I \otimes J}{M}$ is denoted by $\Delta_a(M)$ (the relation is again different due to different choices). In the semisimple setting, the algebras resulting from each of these definitions are isomorphic.

\smallskip

\indent Note that the fusion relation is a relation between matrices in $\mathcal{L}_{0,1}(H) \otimes \mathrm{Mat}_{\dim(I)}(\mathbb{C}) \otimes \mathrm{Mat}_{\dim(J)}(\mathbb{C})$ (for all finite dimensional $I,J$) which implies relations among elements of $\mathcal{L}_{0,1}(H)$ (the coefficients of these matrices). Explicitly, in terms of matrix coefficients it is written as
\[ \forall \, I,J,a,b,c,d, \:\:\:\:\: \overset{I \otimes J}{M}{^{ac}_{bd}} = \overset{I}{M}{^a_i} \overset{IJ}{(R')}{^{ic}_{jk}} \overset{J}{M}{^k_l} \overset{IJ}{(R'^{-1})}{^{jl}_{bd}}, \]
see the definition of the subscripts $1$ and $2$ in section \ref{matrices}. If the two representations $I$ and $J$ are fixed and arbitrary, we can simply write
\begin{equation}\label{fusionL01}
M_{12} = M_1\,R_{21}\,M_2\,R_{21}^{-1}.
\end{equation}
Moreover, one can check that $\overset{(I \otimes J)\otimes K}{M} = \overset{I \otimes (J \otimes K)}{M}$ holds thanks to the Yang-Baxter equation. 

\smallskip

Note that if $f : I \to J$ is a $H$-morphism it holds
\begin{equation}\label{naturaliteL01}
\overset{J}{M} f = f \overset{I}{M}
\end{equation}
where we identify $f$ with its matrix (indeed, by \eqref{naturaliteT}, $f\overset{I}{M} = f \bar j\bigl(\overset{I}{T}\bigr) = \bar j\bigl(f\overset{I}{T}\bigr) = \bar j\bigl(\overset{J}{T}f\bigr) = \bar j\bigl(\overset{J}{T}\bigr) f = \overset{J}{M}f $ where $\bar j$ is the linear map $H^* \to \mathrm{T}(H^*) \to \mathcal{L}_{0,1}(H)$). We call this relation the naturality of the (family of) matrices $\overset{I}{M}$.

\begin{remark}\label{remarkRestrictionCoeffL01}
If $I$ is a submodule or a quotient of $J$, then $\overset{I}{M}$ is a submatrix of $\overset{J}{M}$ thanks to \eqref{naturaliteL01}. Let $\mathcal{G}$ be a set of $H$-modules which generate $\mathrm{mod}_l(H)$ by tensor products, in the sense that every $H$-module is isomorphic to a submodule or a quotient of a tensor product of elements of $\mathcal{G}$. Then by the fusion relation, we see that every matrix coefficient $\overset{I}{M}{^i_j}$ is a linear combination of matrix coefficients of the modules in $\mathcal{G}$. In practice, we restrict ourselves to such a well-chosen set $\mathcal{G}$ to obtain presentations of $\mathcal{L}_{0,1}(H)$ with no many generators. For instance, if $H = \bar U_q(\mathfrak{sl}_2)$, we take $\mathcal{G} = \bigl\{ \mathcal{X}^+(2) \bigr\}$, see section \ref{sectionL01Uq}.
\finEx
\end{remark}

\indent We have an useful analogue of relation \eqref{FRT}.

\begin{proposition}\label{EqRef}
The following exchange relation holds in $\mathcal{L}_{0,1}(H)$:
$$ \overset{IJ}{R}_{12}\overset{I}{M}_1\overset{IJ}{(R')}_{12}\overset{J}{M}_2 = \overset{J}{M}_2\overset{IJ}{R}_{12}\overset{I}{M}_1\overset{IJ}{(R')}_{12}. $$
\end{proposition}
\noindent This relation is called the {\em reflection equation}. It can be written in a shortened way if the representations $I$ and $J$ are fixed and arbitrary:
\begin{equation}\label{reflection}
R_{12}M_1R_{21}M_2 = M_2R_{12}M_1R_{21}.
\end{equation}
\begin{proof}
We have the braiding isomorphism $c_{I,J} = P\overset{IJ}{R} : I \otimes J \to J \otimes I$ where $P$ is the flip tensor $P(x \otimes y) = y \otimes x$. Hence:
\begin{align*}
\overset{IJ}{R}_{12} \, \overset{I}{M}_1 \, \overset{IJ}{(R')}_{12} \, \overset{J}{M}_2 \, \overset{IJ}{(R')}{_{12}^{-1}} &= \overset{IJ}{R} \, \overset{I \otimes J}{M} = P \, c_{I,J} \, \overset{I \otimes J}{M} = P \, \overset{J \otimes I}{M} \, c_{I,J} = P_{12} \, \overset{J}{M}_1 \, \overset{JI}{(R')}_{12} \, \overset{I}{M}_2 \, \overset{JI}{(R')}{_{12}^{-1}} \, P_{12} \, \overset{IJ}{R}_{12}\\
&= \overset{J}{M}_2 \, \overset{JI}{(R')}_{21} \, \overset{I}{M}_1 \, \overset{JI}{(R')}{_{21}^{-1}} \, \overset{IJ}{R}_{12} = \overset{J}{M}_2 \, \overset{IJ}{R}_{12} \, \overset{I}{M}_1.
\end{align*}
We simply used the fusion relation, the naturality \eqref{naturaliteL01} and $P_{12} \overset{JI}{(R')}_{12} P_{12} = \overset{JI}{(R')}_{21} = \overset{IJ}{R}_{12}$.
\end{proof}

Consider the following right action $\cdot$ of $H$ on $\mathcal{L}_{0,1}(H)$, which is the analogue of the right action of the gauge group on the functions:
\begin{equation}\label{actionHsurL01}
\overset{I}{M} \cdot h = \overset{I}{h'} \overset{I}{M} \overset{I}{S(h'')}.
\end{equation}
As in \cite{BR}, one can equivalently work with the corresponding left coaction $\Omega : \mathcal{L}_{0,1}(H) \to \mathcal{O}(H) \otimes \mathcal{L}_{0,1}(H)$ defined by  
\[ \Omega(\overset{I}{M}{^a_b}) = \overset{I}{T}{^a_i}S(\overset{I}{T}{^j_b}) \otimes \overset{I}{M}{^i_j} \]
so that we recover $\cdot$ by evaluation: $x \cdot h = (\langle ?, h \rangle \otimes \mathrm{id}) \circ \Omega(x)$. If we view $\mathcal{O}(H)$ and $\mathcal{L}_{0,1}(H)$ as subalgebras of $\mathcal{O}(H) \otimes \mathcal{L}_{0,1}(H)$ in the canonical way, then $\Omega$ is simply written as $\Omega(\overset{I}{M}) = \overset{I}{T}\overset{I}{M}S(\overset{I}{T})$. 

\begin{proposition}\label{L01moduleAlgebra}
The right action $\cdot$ is a $H$-module-algebra structure on $\mathcal{L}_{0,1}(H)$. Equivalently, $\Omega$ is a left $\mathcal{O}(H)$-comodule-algebra structure on $\mathcal{L}_{0,1}(H)$. 
\end{proposition}
\begin{proof}
One must show for instance that $\Omega$ is a morphism of algebras (\textit{i.e.} that it preserves the fusion relation), as in \cite{BR}. With the shortened notation explained before, the computation is as follows:
\[\begin{array}{lll}
\Omega(M)_{12} &= T_{12}\,M_{12} \, S(T_{12}) &\:\:\: \text{(definition)}\\
&= T_1\,T_2\,M_1\,R_{21}\,M_2\,R^{-1}_{21}\,S(T)_2\,S(T)_1 &\:\:\: \text{(eq. (\ref{dualHopf}) and (\ref{fusionL01}))}\\
&= T_1\,M_1\,T_2\,R_{21}\,M_2\,R^{-1}_{21}\,S(T)_2\,S(T)_1 &\:\:\: \text{(commuting elements in tensor product algebra)}\\
&= T_1\,M_1\,T_2\,R_{21}\,M_2\,S(T)_1\,S(T)_2\,R_{21}^{-1} &\:\:\: \text{(eq. (\ref{FRT}))}\\
&= T_1\,M_1\,T_2\,R_{21}\,S(T)_1\,M_2\,S(T)_2\,R_{21}^{-1} &\:\:\: \text{(commuting elements in tensor product algebra)}\\
&= T_1\,M_1\,S(T)_1\,R_{21}\,T_2\,M_2\,S(T)_2\,R_{21}^{-1} &\:\:\: \text{(eq. (\ref{FRT}))}\\
&= \Omega(M)_1\,R_{21}\,\Omega(M)_2\,R_{21}^{-1} &\:\:\: \text{(definition)}.\\
\end{array}\]

\end{proof}

\indent We say that an element $x \in \mathcal{L}_{0,1}(H)$ is invariant if for all $h \in H$, $x \cdot h = \varepsilon(h)x$ (or equivalently, $\Omega(x) = \varepsilon \otimes x$) and we denote by $\mathcal{L}_{0,1}^{\text{inv}}(H)$ the subalgebra of invariant elements of $\mathcal{L}_{0,1}(H)$ (also called ``observables'').

\begin{exemple}\label{exempleCaractereW}
For any representation $I$, the element 
\begin{equation}\label{defWexplicite}
\overset{I}{W} = \mathrm{tr}_q\bigl(\overset{I}{M}\bigr) = \mathrm{tr}\bigl(\overset{I}{g}\overset{I}{M}\bigr)
\end{equation}
is invariant:
\[ \mathrm{tr}(\overset{I}{g}\overset{I}{M}) \cdot h = \mathrm{tr}(\overset{I}{g} \overset{I}{h'} \overset{I}{M} \overset{I}{S(h'')}) = \mathrm{tr}(\overset{I}{g} \overset{I}{S^{-1}(h'')}\overset{I}{h'} \overset{I}{M}) = \varepsilon(h)\mathrm{tr}(\overset{I}{g}\overset{I}{M}). \]
However, this splits on extensions:
\[ 0 \rightarrow I \rightarrow K \rightarrow J \rightarrow 0 \: \text{ exact} \:\:\: \implies \:\:\: \overset{K}{W} = \overset{I}{W} + \overset{J}{W} \]
and in general (when $H$ is non-semisimple), the span of the $\overset{I}{W}$'s is strictly smaller than the subalgebra of invariant elements.
\finEx
\end{exemple}

\subsection{Isomorphism $\mathcal{L}_{0,1}(H) \cong H$}\label{isoL01H}
Recall the right adjoint action of $H$ on itself defined by $a \cdot h = S(h')ah''$ with $a,h \in H$, whose invariant elements are the central elements of $H$.

\begin{proposition}\label{reshetikhin}
If we endow $H$ with the right adjoint action, the following map is a morphism of (right) $H$-module-algebras:
$$\fonc{\Psi_{0,1}}{\mathcal{L}_{0,1}(H)}{H}{\overset{I}{M}}{(\overset{I}{T} \otimes \mathrm{id})(RR') = \overset{I}{L}{^{(+)}}\overset{I}{L}{^{(-)-1}}.}$$
In particular, $\Psi_{0,1}$ brings invariant elements to central elements.
\end{proposition}
\begin{proof}
Using the relations of (\ref{propertiesL})\footnote{More precisely, we use relations easily implied by \eqref{propertiesL}. For instance for the second equality, we used the relation $L^{(-)-1}_1 R_{21} L^{(+)}_2 = L^{(+)}_2 R_{21} L^{(-)-1}_1$, which is obtained as follows: exchanging $I$ and $J$ in the second relation of \eqref{propertiesL} we have $\overset{JI}{R}_{12} \overset{J}{L}{_1^{(+)}} \overset{I}{L}{^{(-)}_2} = \overset{I}{L}{^{(-)}_2} \overset{J}{L}{_1^{(+)}} \overset{JI}{R}_{12}$ and then applying the flip map $\fonc{P}{H \otimes \mathrm{End}_{\mathbb{C}}(J) \otimes \mathrm{End}_{\mathbb{C}}(I)}{H \otimes \mathrm{End}_{\mathbb{C}}(I) \otimes \mathrm{End}_{\mathbb{C}}(J)}{x \otimes Y \otimes Z}{x \otimes Z \otimes Y}$ we get $\overset{JI}{R}_{21} \overset{J}{L}{_2^{(+)}} \overset{I}{L}{^{(-)}_1} = \overset{I}{L}{^{(-)}_1} \overset{J}{L}{_2^{(+)}} \overset{JI}{R}_{21}$ which is written $R_{21} L_2^{(+)} L_1^{(-)} = L_1^{(-)} L_2^{(+)} R_{21}$ in the shortened notation. Recall that in the shortened notation the index $1$ (resp. $2$) implicitly means evaluation in a representation $I$ (resp. $J$) (thus $R_{21}$ means $\overset{JI}{R}_{21} = \bigl(\overset{IJ}{R'}\bigr)_{12}$). More generally, any permutation of the indices in defining relations is allowed when one does computations with the shortened notation. We will no longer give such details in subsequent computations.}, we check that $\Psi_{0,1}$ preserves the relation of Definition \ref{defL01}:
\begin{align*}
\Psi_{0,1}(M)_1R_{21}\Psi_{0,1}(M)_2R^{-1}_{21} &= L^{(+)}_1L^{(-)-1}_1 R_{21} L^{(+)}_2L^{(-)-1}_2R_{21}^{-1} = L^{(+)}_1  L^{(+)}_2 R_{21} L^{(-)-1}_1 L^{(-)-1}_2R_{21}^{-1}\\
&= L^{(+)}_1  L^{(+)}_2 L^{(-)-1}_2 L^{(-)-1}_1 = L^{(+)}_{12}L^{(-)-1}_{12} = \Psi_{0,1}(M)_{12}.
\end{align*}
For the $H$-linearity:
\begin{align*}
\Psi_{0,1}(\overset{I}{h'}\overset{I}{M}\overset{I}{S(h'')}) &= (\overset{I}{T} \otimes \mathrm{id})\!\left(h' \otimes 1 \, RR' \, S(h'') \otimes 1\right)\\
&= (\overset{I}{T} \otimes \mathrm{id})\!\left(h' \otimes 1 \, RR' \, S(h)''' \otimes S(h)''h''''\right)\\
&= (\overset{I}{T} \otimes \mathrm{id})\!\left(h'S(h)''' \otimes S(h)''\, RR' \, 1 \otimes h''''\right)\\
&= (\overset{I}{T} \otimes \mathrm{id})\!\left(1 \otimes S(h')\, RR' \, 1 \otimes h''\right) = S(h') \Psi_{0,1}(\overset{I}{M})h''.
\end{align*}
We used the basic properties of $S$ and the fact that $\Delta^{\mathrm{op}}R = R\Delta$, with $\Delta^{\mathrm{op}}(h) = h'' \otimes h'$.
\end{proof}
\noindent We call $\Psi_{0,1}$ the {\em Reshetikhin -- Semenov-Tian-Shansky -- Drinfeld morphism} (RSD morphism for short) \cite{RS, drinfeld}. The difference with the morphism $\Psi$ of section \ref{rappelHopf} is that the source spaces are different.

\smallskip

\indent Write $\mathrm{T}(H^*) = \bigoplus_{ n \in \mathbb{N}} \mathrm{T}_n(H^*)$, where $\mathrm{T}_n(H^*)$ is the subspace generated by all the products $\psi_1 \cdots \psi_n$, with $\psi_i \in H^*$ for each $i$.

\begin{lemma}\label{lemmaDegre}
Each element of $\mathrm{T}(H^*)$ is equivalent modulo the fusion relation of $\mathcal{L}_{0,1}(H)$ to an element of $\mathrm{T}_1(H^*)$. It follows that $\dim\!\left(\mathcal{L}_{0,1}(H)\right) \leq \dim\!\left(H^*\right)$.
\end{lemma}
\begin{proof}
It suffices to show that the product of two elements of $T_1(H^*)$ is equivalent to a linear combination of elements of $\mathrm{T}_1(H^*)$, and the result follows by induction. We can restrict to matrix coefficients since they linearly span $H^*$. 
If we write $R = a_i \otimes b_i$, then the fusion relation is rewritten as:
$$ \overset{I \otimes J}{M}\!\!_{12} (\overset{IJ}{R'})_{12} = (\overset{J}{a_i})_2 \overset{I}{M}\!_1 \overset{J}{M}\!_2 (\overset{I}{b_i})_1 $$
Using $a_ja_i \otimes b_iS(b_j) = 1 \otimes 1$, we get:
\begin{equation}\label{eqFusionInverse}
\overset{I}{M}\!_1 \overset{J}{M}\!_2 = \overset{J}{(a_i)}_2\overset{I \otimes J}{M}\!\!_{12} (\overset{IJ}{R'})_{12} \overset{J}{S(b_i)}_1
\end{equation}
and this give the result since $\overset{I}{M}_1 \overset{J}{M}_2$ contains all the possible products between the coefficients of $\overset{I}{M}$ and those of $\overset{J}{M}$.
\end{proof}

\begin{proposition}\label{IsoRSD}
Recall that we assume that $H$ is a finite dimensional factorizable Hopf algebra. Then the RSD morphism $\Psi_{0,1}$ gives an isomorphism of $H$-module-algebras $\mathcal{L}_{0,1}(H) \cong H$. It follows that $\mathcal{L}_{0,1}^{\mathrm{inv}}(H) \cong \mathcal{Z}(H)$.
\end{proposition}
\begin{proof}
Since $H$ is factorizable, $\Psi_{0,1}$ is surjective. Hence $\dim\!\left(\mathcal{L}_{0,1}(H)\right) \geq \dim(H)$. But by Lemma \ref{lemmaDegre}, $\dim\!\left(\mathcal{L}_{0,1}(H)\right) \leq \dim\!\left(H^*\right) = \dim(H)$. Thus $\dim\!\left(\mathcal{L}_{0,1}(H)\right) = \dim\!\left(H\right)$.
\end{proof}

\noindent Let us point out obvious consequences. First, the matrices $\overset{I}{M}$ are invertible since $RR'$ is invertible. Second, this theorem allows us to identify $\mathcal{L}_{0,1}(H)$ with $H$ via $\overset{I}{M} = \overset{I}{L} \,\!^{(+)}\overset{I}{L} \,\!^{(-)-1}$, where the matrices $L^{(\pm)}$ are defined in $(\ref{L})$. We will always work with this identification in the sequel. 

\begin{remark}\label{F01}
Due to Proposition \ref{IsoRSD}, there is an isomorphism of vector spaces $f : \mathcal{L}_{0,1}(H) \to H^*$ given by $\overset{I}{M}{^i_j} \mapsto \overset{I}{T}{^i_j}$. We define a $H$-module-algebra structure on $H^*$, denoted by $\mathcal{F}_{0,1}(H)$ and with product $\ast$, by requiring $f$ to be an isomorphism of (right) $H$-module-algebras. The right $H$-action is
\begin{equation*}
\varphi \cdot h = \varphi\!\left( h' ? S(h'') \right).
\end{equation*}
Using \eqref{eqFusionInverse}, \eqref{dualHopf} and obvious commutation relations we have
\begin{align*}
\overset{I}{T}{^{\alpha}_{\beta}} \ast \overset{I}{T}{^{\gamma}_{\delta}} &= f\!\left(\overset{I}{M}\!_1 \overset{J}{M}\!_2\right)^{\alpha \gamma}_{\beta \delta} = f\!\left(\overset{J}{(a_i)}_2\overset{I \otimes J}{M}\!\!_{12} (\overset{IJ}{R'})_{12} \overset{I}{S(b_i)}_1\right)^{\alpha \gamma}_{\beta \delta} = \left(\overset{J}{(a_i)}_2\overset{I \otimes J}{T}\!\!\!_{12} \, \overset{I}{(b_j)}_1 \overset{I}{S(b_i)}_1 \overset{J}{(a_j)}_2 \right)^{\alpha \gamma}_{\beta \delta}\\
& =  \left( \bigl( \overset{I}{T} \overset{I}{b_j} \overset{I}{S(b_i)} \bigr)_1 \bigl( \overset{J}{a_i} \overset{J}{T}  \overset{J}{a_j} \bigr)_2 \right)^{\alpha \gamma}_{\beta \delta} = \overset{I}{T}{^{\alpha}_{\beta}}\!\left( ? b_j S(b_i) \right) \overset{I}{T}{^{\gamma}_{\delta}}\!\left( a_i ? a_j\right).
\end{align*}
In other words,
\begin{equation}\label{produitL01Explicite}
\varphi \ast \psi = \varphi\!\left( ? b_j S(b_i) \right) \psi\!\left(a_i ? a_j\right), \:\:\: x_m \mapsto \varphi\!\left( x_m' b_j S(b_i) \right) \psi\!\left(a_i x_m'' a_j\right).
\end{equation}
This is the product of the functions $\varphi, \psi \in \mathcal{F}_{0,1}(H)$ and its evaluation on the discrete connection which assigns $x_m \in H$ to the loop $m$, see Figure \ref{Sigma01AvecConnexion} and the Introduction.
\begin{figure}[!h]
\centering
\begingroup%
  \makeatletter%
  \providecommand\color[2][]{%
    \errmessage{(Inkscape) Color is used for the text in Inkscape, but the package 'color.sty' is not loaded}%
    \renewcommand\color[2][]{}%
  }%
  \providecommand\transparent[1]{%
    \errmessage{(Inkscape) Transparency is used (non-zero) for the text in Inkscape, but the package 'transparent.sty' is not loaded}%
    \renewcommand\transparent[1]{}%
  }%
  \providecommand\rotatebox[2]{#2}%
  \newcommand*\fsize{\dimexpr\f@size pt\relax}%
  \newcommand*\lineheight[1]{\fontsize{\fsize}{#1\fsize}\selectfont}%
  \ifx\svgwidth\undefined%
    \setlength{\unitlength}{151.4375063bp}%
    \ifx\svgscale\undefined%
      \relax%
    \else%
      \setlength{\unitlength}{\unitlength * \real{\svgscale}}%
    \fi%
  \else%
    \setlength{\unitlength}{\svgwidth}%
  \fi%
  \global\let\svgwidth\undefined%
  \global\let\svgscale\undefined%
  \makeatother%
  \begin{picture}(1,0.76417815)%
    \lineheight{1}%
    \setlength\tabcolsep{0pt}%
    \put(0.72950757,0.61990526){\color[rgb]{0,0,0}\makebox(0,0)[lt]{\lineheight{1.25}\smash{\begin{tabular}[t]{l}$m$\end{tabular}}}}%
    \put(0,0){\includegraphics[width=\unitlength,page=1]{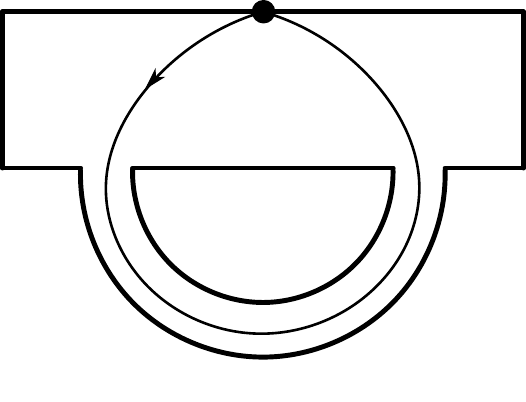}}%
    \put(0.46187537,0.01557341){\color[rgb]{0,0,0}\makebox(0,0)[lt]{\lineheight{1.25}\smash{\begin{tabular}[t]{l}$x_m$\end{tabular}}}}%
  \end{picture}%
\endgroup%

\caption{A discrete connection $x_m \in H$ on $\Sigma_{0,1}^{\mathrm{o}}$.}
\label{Sigma01AvecConnexion}
\end{figure}
\finEx
\end{remark}

\indent Recall that we denote by $\SLF(H)$ the space of symmetric linear forms on $H$:
$$ \SLF(H) = \left\{ \psi \in H^* \left\vert\,  \forall\, x,y \in H, \:\: \psi(xy) = \psi(yx) \right.\right\} $$
and that $\SLF(H)$ is obviously a subalgebra of $\mathcal{O}(H)$. Also recall from section \ref{rappelHopf} the Drinfeld morphism
\begin{equation*}
\fonc{\mathcal{D}}{H^*}{H}{\psi}{(\psi \otimes \text{id})\bigl( (g \otimes 1) RR'\bigr) = \psi(ga_ib_j)b_ia_j}.
\end{equation*}
By Lemma \ref{propMorDrinfeld}, $\mathcal{D}$ induces an isomorphism of algebras between $\SLF(H)$ and $\mathcal{Z}(H)$. Hence $\SLF(H) \cong \mathcal{Z}(H) = \mathcal{L}^{\mathrm{inv}}_{0,1}(H)$.

\begin{exemple}\label{exempleWChi}
Consider the invariant elements of Example \ref{exempleCaractereW}; then it holds by definition 
\begin{equation}\label{DWChi}
\mathcal{D}\bigl(\chi^I\bigr) = \Psi_{0,1}\bigl(\overset{I}{W}\bigr),
\end{equation}
where $\chi^I = \mathrm{tr}\bigl( \overset{I}{T} \bigr)$ is the character of the representation $I$. Due to \eqref{dualHopf}, we have $\chi^{I \otimes J} = \chi^I \chi^J$ and hence the same fusion rule applies to the observables $\overset{I}{W}$:
\[ \overset{I \otimes J}{W} = \overset{I}{W} \overset{J}{W}. \]
This implies that the span of the $\overset{I}{W}$'s is a subalgebra of $\mathcal{L}_{0,1}^{\mathrm{inv}}(H)$ (in general strictly smaller, see for instance \eqref{SLFcentrauxUq}).
\finEx
\end{exemple}

\indent Let us fix a notation. Every $\psi \in H^*$ can be written as $\psi = \sum_{i,j,I}\lambda_{ij}^I\overset{I}{T}\,\!^i_j$ with $\lambda_{ij}^I \in \mathbb{C}$. In order to avoid the indices, define for each $I$ a matrix $\Lambda_I \in \mathrm{Mat}_{\dim(I)}(\mathbb{C})$ by $(\Lambda_I)^i_j = \lambda^I_{ji}$. Then $\psi$ can be expressed as:
$$ \psi = \sum_I \mathrm{tr}(\Lambda_I\overset{I}{T}). $$
We record these observations as a lemma.

\begin{lemma}\label{writingInvariants}
Every $x \in \mathcal{L}_{0,1}(H)$ can be expressed as:
$$ x = \sum_I \mathrm{tr}(\Lambda_I\overset{I}{g}\overset{I}{M}) $$
such that $\mathcal{D}^{-1}(x) = \sum_I \mathrm{tr}(\Lambda_I\overset{I}{T})$. Moreover, if $x \in \mathcal{L}^{\mathrm{inv}}_{0,1}(H)$, then $\mathcal{D}^{-1}(x) \in \SLF(H)$.
\end{lemma}

\begin{remark}\label{remarkWriting}
Let us stress that, due to non-semi-simplicity, this way of writing elements of $\mathcal{L}_{0,1}(H)$ and of $\SLF(H)$ is in general not unique, see the comments in section \ref{matrices}.
\finEx
\end{remark}

\section{The handle algebra $\mathcal{L}_{1,0}(H)$}\label{sectionHandleL10}
\indent We assume that $H$ is a finite dimensional factorizable ribbon Hopf algebra. Note however that the ribbon assumption is not needed in sections \ref{defL10} and \ref{sectionIsoL10Heisenberg}.

\subsection{Definition of $\mathcal{L}_{1,0}(H)$ and $H$-module-algebra structure}\label{defL10}
\indent Consider the free product $\mathcal{L}_{0,1}(H) \ast \mathcal{L}_{0,1}(H)$, and let $j_1$ (resp. $j_2$) be the canonical injection in the first (resp. second) copy of $\mathcal{L}_{0,1}(H)$. We define $\overset{I}{B} = j_1(\overset{I}{M})$ and $\overset{I}{A} = j_2(\overset{I}{M})$.

\begin{definition}\label{definitionL10}
The handle algebra $\mathcal{L}_{1,0}(H)$ is the quotient of $\mathcal{L}_{0,1}(H) \ast \mathcal{L}_{0,1}(H)$ by the following exchange relations:
$$ \overset{IJ}{R}_{12}\overset{I}{B}_1\overset{IJ}{(R')}_{12}\overset{J}{A}_2 = \overset{J}{A}_2\overset{IJ}{R}_{12}\overset{I}{B}_1\overset{IJ}{R}{^{-1}_{12}} $$
for all finite dimensional $H$-modules $I, J$.
\end{definition}
\noindent The exchange relation above is the same as in \cite[Def 1]{BNR} except that $A$ and $B$ are switched; the one of \cite[Def 12]{AGS2} and \cite[eq (3.14)]{AS} is different, due to a different choice of the action of $H$ on $\mathcal{L}_{1,0}(H)$. In the semisimple setting, the algebras resulting from each of these definitions are isomorphic. 

\smallskip

\indent The exchange relation is a relation between matrices in $\mathcal{L}_{1,0}(H) \otimes \mathrm{Mat}_{\dim(I)}(\mathbb{C}) \otimes \mathrm{Mat}_{\dim(J)}(\mathbb{C})$ (for all finite dimensional $I,J$) which implies relations among elements of $\mathcal{L}_{1,0}(H)$ (the coefficients of these matrices), namely
\[ \forall \, I,J,a,b,c,d, \:\:\:\:\: \overset{IJ}{R}{^{ac}_{ij}} \, \overset{I}{B}{^i_k } \, (\overset{IJ}{R'})^{kj}_{bl} \, \overset{J}{A}{^l_d} = \overset{J}{A}{^c_i} \, \overset{IJ}{R}{^{ai}_{jk}} \, \overset{I}{B}{^j_l} \, (\overset{IJ}{R}{^{-1}})^{lk}_{bd}. \]
Like the other relations before, the $\mathcal{L}_{1,0}(H)$-exchange relation can be written more simply as:
\begin{equation}\label{echangeL10}
R_{12}\,B_1\,R_{21}\,A_2 = A_2\,R_{12}\,B_1\,R_{12}^{-1}.
\end{equation}
By \eqref{naturaliteL01}, if $f : I \to J$ is a morphism it holds
\begin{equation}\label{naturaliteL10}
\overset{J}{B} f = f \overset{I}{B}, \:\:\:\:\:\: \overset{J}{A} f = f \overset{I}{A}
\end{equation}
where we identify $f$ with its matrix. We call this relation the naturality of the (families of) matrices $\overset{I}{B}, \overset{I}{A}$. Also note that the content of Remark \ref{remarkRestrictionCoeffL01} also applies to $\mathcal{L}_{1,0}(H)$: in practice, we can restrict to a set $\mathcal{G}$ of well-chosen $H$-modules and when we write $\overset{I}{B}$ and $\overset{I}{A}$, we can assume that $I \in \mathcal{G}$. For instance, if $H = \bar U_q(\mathfrak{sl}_2)$, we take $\mathcal{G} = \bigl\{ \mathcal{X}^+(2) \bigr\}$, see section \ref{sectionL01Uq}.

\smallskip

\indent Similarly to $\mathcal{L}_{0,1}(H)$, consider the following right action of $H$ on $\mathcal{L}_{1,0}(H)$, which is the analogue of the action of the gauge group on the functions:
\begin{equation}\label{actionHsurL10}
\overset{I}{B} \cdot h = \overset{I}{h'} \overset{I}{B} \overset{I}{S(h'')}, \:\:\:\:\: \overset{I}{A} \cdot h = \overset{I}{h'} \overset{I}{A} \overset{I}{S(h'')}.
\end{equation}
As above and like in \cite{BR}, it is equivalent to work with the corresponding left coaction $\Omega : \mathcal{L}_{1,0}(H) \to \mathcal{O}(H) \otimes \mathcal{L}_{1,0}(H)$ defined by  
\[ \Omega(\overset{I}{B}) = \overset{I}{T}\overset{I}{B}S(\overset{I}{T}), \:\:\:\:\: \Omega(\overset{I}{A}) = \overset{I}{T}\overset{I}{A}S(\overset{I}{T}). \]

\begin{proposition}\label{L10moduleAlgebra}
The right action $\cdot$ is a $H$-module-algebra structure on $\mathcal{L}_{1,0}(H)$. Equivalently, $\Omega$ is a left $\mathcal{O}(H)$-comodule-algebra structure on $\mathcal{L}_{1,0}(H)$. 
\end{proposition}
\begin{proof}
One must show that $\Omega$ is an algebra morphism, as in \cite{BR}. This amounts to check that $\Omega$ is compatible with the exchange relation, which is similar to the proof of Proposition \ref{L01moduleAlgebra} and is left to the reader. 
\end{proof}

\indent We denote by $\mathcal{L}_{1,0}^{\text{inv}}(H)$ the subalgebra of invariant elements of $\mathcal{L}_{1,0}(H)$ (also called ``observables''). For instance, the elements 
\begin{equation}\label{FamilleInvL10}
\mathrm{tr}_{12}\!\left(\overset{I \otimes J}{g}\!\!\!_{12}\Phi\overset{I}{A}_1\overset{IJ}{(R')}_{12}\overset{J}{B}_2\overset{IJ}{R}_{12}\right)
\end{equation}
with $\Phi \in \End_H(I \otimes J)$ and $\mathrm{tr}_{12} = \mathrm{tr} \otimes \mathrm{tr}$, are invariant.

\smallskip

\noindent \textbf{Notation.}~~ Let $\overset{I}{N} = \overset{I}{v}{^m} \overset{I}{N}{^{n_1}_1} \ldots \overset{I}{N}{^{n_l}_l} \in \mathrm{Mat}_{\dim(I)}\!\left(\mathcal{L}_{1,0}(H)\right)$, where $m, n_i \in \mathbb{Z}$  and each $N_i$ is $A$ or $B$. By definition of the right action of $H$ on $\mathcal{L}_{1,0}(H)$, we have a morphism of $H$-modules 
\begin{equation*}
\fonc{j_N}{\mathcal{L}_{0,1}(H)}{\mathcal{L}_{1,0}(H)}{\overset{I}{M}}{\overset{I}{N}}.
\end{equation*}
Let $x \in \mathcal{L}_{0,1}(H)$, we introduce the notation
\begin{equation}\label{notationImbed}
x_N = j_N(x).
\end{equation}
Since we identify $\mathcal{L}_{0,1}(H)$ with $H$ through $\Psi_{0,1}$ we also use this notation when $x \in H$: $x_N = \Psi_{0,1}^{-1}(x)_N$. Note that if $x \in \mathcal{L}_{0,1}^{\mathrm{inv}}(H) \cong \mathcal{Z}(H)$, then $x_N \in \mathcal{L}_{1,0}^{\mathrm{inv}}(H)$. The following lemma is an obvious fact.
\begin{lemma}\label{injectionFusion10}
If $N$ satisfies the fusion relation of $\mathcal{L}_{0,1}(H)$, $\overset{I \otimes J}{N}\!\!_{12} = \overset{I}{N}(i)_1\,\overset{IJ}{(R')}_{12}\,\overset{J}{N}(i)_2\, \overset{IJ}{(R')}{_{12}^{-1}}$, then $j_N$ is a morphism of $H$-module-algebras: $(xy)_N = x_N y_N$.
\end{lemma}
Note that we allow $\overset{I}{v}{^m}$ in the formula of $\overset{I}{N}$ due to the fusion relation. Indeed, a suitable product of matrices $\overset{I}{A}{^{\pm 1}}, \overset{I}{B}{^{\pm 1}}$ satisfies the fusion relation when it is correctly normalized by some power of $v$, see \textit{e.g.} \eqref{normaFusion1}, \eqref{normaFusion2} and Proposition \ref{autoAlphaBeta}.
\begin{exemple}
Taking back the elements introduced in Example \ref{exempleCaractereW}, we have
\begin{equation}\label{WAWBWBA}
\overset{I}{W}_{\!\! A} = \mathrm{tr}_q\bigl( \overset{I}{A} \bigr), \:\:\:\:\: \overset{I}{W}_{\!\! B} = \mathrm{tr}_q\bigl( \overset{I}{B} \bigr), \:\:\:\:\: \overset{I}{W}_{\!\! vB^{-1}A} = \mathrm{tr}_q\bigl( \overset{I}{v}\overset{I}{B}{^{-1}} \overset{I}{A} \bigr), \ldots
\end{equation}
These are invariant elements of $\mathcal{L}_{1,0}(H)$.
\finEx
\end{exemple}

\begin{remark}
Recall from remark \ref{remarkWriting} that the matrix coefficients do not form a basis of $\mathcal{L}_{0,1}(H)$. They just linearly span this space. However, the maps $j_{wA^{m_1}B^{n_1} \ldots A^{m_k}B^{n_k}}$ are well-defined. Indeed, first observe that 
\begin{align*}
j_B &:  \mathcal{L}_{0,1}(H) \overset{j_1}{\hookrightarrow} \mathcal{L}_{0,1}(H) \ast \mathcal{L}_{0,1}(H) \overset{\pi}{\longrightarrow\!\!\!\!\!\to} \mathcal{L}_{1,0}(H)\\
j_A &:  \mathcal{L}_{0,1}(H) \overset{j_2}{\hookrightarrow} \mathcal{L}_{0,1}(H) \ast \mathcal{L}_{0,1}(H) \overset{\pi}{\longrightarrow\!\!\!\!\!\to} \mathcal{L}_{1,0}(H)
\end{align*}
are well-defined. Let us show for instance that the map $j_{A^{-1}B^{-1}A}$ is well-defined. Assume that $\lambda^a_b\overset{I}{T}{^b_a} = 0$. Applying the coproduct in $\mathcal{O}(H)$ twice and tensoring with $\mathrm{id}_H$, we get:
$$ \lambda^a_b\,\overset{I}{T}{^b_k} \otimes \mathrm{id}_H \otimes \overset{I}{T}{^k_l} \otimes \mathrm{id}_H \otimes \overset{I}{T}{^l_a} \otimes \mathrm{id}_H = 0. $$
We evaluate this on $\left(RR'\right)^{-1} \otimes \left(RR'\right)^{-1} \otimes RR'$:
$$ \lambda^a_b\,(\overset{I}{M}{^{-1}}){^b_k} \otimes (\overset{I}{M}{^{-1}}){^k_l} \otimes \overset{I}{M}{^l_a} = 0.$$
Finally, we apply the map $j_A \otimes j_B \otimes j_A$ and multiplication in $\mathcal{L}_{1,0}(H)$:
$$ \lambda^a_b\,(\overset{I}{A}{^{-1}}\overset{I}{B}{^{-1}}\overset{I}{A}){^b_a} = 0 $$
as desired. A similar proof can be used to show that all the other maps defined by means of matrix coefficients (like $\Psi_{1,0}$ or $\alpha, \beta$ below etc..) are well-defined.
\finEx
\end{remark}

\subsection{Isomorphism $\mathcal{L}_{1,0}(H) \cong \mathcal{H}(\mathcal{O}(H))$}\label{sectionIsoL10Heisenberg}
Recall that the definition and properties of the Heisenberg double $\mathcal{H}(\mathcal{O}(H))$ are summarized in section \ref{heisenbergDouble}.
\begin{proposition}
The following map is a morphism of algebras:
$$\begin{array}{crll}
\Psi_{1,0} :& \mathcal{L}_{1,0}(H) & \rightarrow & \mathcal{H}(\mathcal{O}(H)) \\
                         & \overset{I}{B} &\mapsto &  \overset{I}{L}\,\!^{(+)}\overset{I}{T}\overset{I}{L}\,\!^{(-)-1}.\\
                         & \overset{I}{A} &\mapsto & \overset{I}{L}\,\!^{(+)}\overset{I}{L}\,\!^{(-)-1}\\
\end{array}$$
\end{proposition}
\begin{proof}
We have to check that the fusion and exchange relations are compatible with $\Psi_{1,0}$. Observe that the restriction of $\Psi_{1,0}$ to the first copy of $\mathcal{L}_{0,1}(H) \subset \mathcal{L}_{1,0}(H)$ is just the RSD morphism $\Psi_{0,1}$, thus $\Psi_{1,0}$ is compatible with the fusion relation over $A$. For the fusion relation over $B$, we have:
$$\begin{array}{lll}
\Psi_{1,0}(B)_{12} &= L^{(+)}_{12} \, T_{12} \, L^{(-)-1}_{12} &\:\:\: \text{(definition)}\\
&= L^{(+)}_1 \, L^{(+)}_2 \, T_1 \, T_2 \, L^{(-)-1}_2 \, L^{(-)-1}_1 &\:\:\: \text{(eq. (\ref{propertiesL}) and (\ref{dualHopf}))}\\
&= L^{(+)}_1 \, T_1 \,  L^{(+)}_2 \, R_{21} \, T_2 \, L^{(-)-1}_2 \, L^{(-)-1}_1 &\:\:\: \text{(Lemma \ref{echangeHeisenberg})}\\
&=  L^{(+)}_1 \, T_1 \,  L^{(+)}_2 \, R_{21} \, T_2 \, R_{21} \, L^{(-)-1}_1 \, L^{(-)-1}_2 R_{21}^{-1} &\:\:\: \text{(eq. (\ref{propertiesL}))}\\
&= L^{(+)}_1 \, T_1 \,  L^{(+)}_2 \, R_{21} \, L^{(-)-1}_1 \,T_2 \,  L^{(-)-1}_2 R_{21}^{-1} &\:\:\: \text{(Lemma \ref{echangeHeisenberg})}\\
&= L^{(+)}_1 \, T_1 \,  L^{(-)-1}_1 \, R_{21} \, L^{(+)}_2 \,T_2 \,  L^{(-)-1}_2 R_{21}^{-1} &\:\:\: \text{(eq. (\ref{propertiesL}))}\\
&= \Psi_{1,0}(B)_1\,R_{21}\,\Psi_{1,0}(B)_2\,R_{21}^{-1} &\:\:\: \text{(definition)}.\\
\end{array}$$
The same kind of computation allows one to show that $\Psi_{1,0}$ is compatible with the $\mathcal{L}_{1,0}$-exchange relation.
\end{proof}

\indent We wish to show that $\Psi_{1,0}$ is an isomorphism.

\begin{lemma}\label{lemmaDimL10}
Every element in $\mathcal{L}_{1,0}(H)$ can be written as $\sum_i (x_i)_B (y_i)_A$ with $x_i, y_i \in \mathcal{L}_{0,1}(H)$. It follows that $\dim\!\left(\mathcal{L}_{1,0}(H)\right) \leq \dim\!\left(\mathcal{L}_{0,1}(H)\right)^2 = \dim(H)^2$.
\end{lemma}
\begin{proof}
This is the same proof as in Lemma \ref{lemmaDegre}. It suffices to show that an element like $y_A x_B$ can be expressed as $\sum_i (x_i)_B (y_i)_A$. 
The exchange relation can be rewritten as:
\begin{equation}\label{echangeL10Inverse}
\overset{I}{A}_1\overset{J}{B}_2 = \overset{J}{(a_i)}_2\overset{IJ}{(R')}_{12}\overset{J}{B}_2\overset{IJ}{R}_{12}\overset{I}{A}_1\overset{IJ}{(R')}_{12}\overset{I}{S(b_i)}_1.
\end{equation}
and the result follows since $\overset{I}{A}_1 \overset{J}{B}_2$ contains all the possible products between the coefficients of $\overset{I}{A}$ and those of $\overset{J}{B}$.
\end{proof}

\begin{proposition}\label{isoL10Heisenberg}
Recall that we assume that $H$ is a finite dimensional factorizable Hopf algebra. $\Psi_{1,0}$ gives an isomorphism of algebras $\mathcal{L}_{1,0}(H) \cong \mathcal{H}(\mathcal{O}(H))$. It follows that $\mathcal{L}_{1,0}(H)$ is a matrix algebra: $\mathcal{L}_{1,0}(H) \cong \mathrm{Mat}_{\dim(H)}(\mathbb{C})$ and in particular has trivial center.
\end{proposition}
\begin{proof}
Observe that $\Psi_{1,0} \circ j_A = i_H \circ \Psi_{0,1}$ where $i_H : H \to \mathcal{H}(\mathcal{O}(H))$ is the canonical inclusion. Since $\Psi_{0,1}$ is an isomorphism, there exist matrices $\overset{I}{A}\,\!^{(\pm)}$ such that 
$$ \Psi_{1,0}(\overset{I}{A}\,\!^{(\pm)}) = \overset{I}{L}\,\!^{(\pm)} \in \mathrm{Mat}_{\dim(I)}(\mathcal{H}(\mathcal{O}(H))). $$
Moreover, we have:
$$ \Psi_{1,0}(\overset{I}{A}\,\!^{(+)-1}\overset{I}{B}\overset{I}{A}\,\!^{(-)}) = \overset{I}{T} \in \mathrm{Mat}_{\dim(I)}(\mathcal{H}(\mathcal{O}(H))). $$
Thus $\Psi_{1,0}$ is surjective, and hence $\dim\!\left(\mathcal{L}_{1,0}(H)\right) \geq \dim\!\left(\mathcal{H}(\mathcal{O}(H))\right) = \dim(H)^2$. This together with Lemma \ref{lemmaDimL10} gives $\dim\!\left(\mathcal{L}_{1,0}(H)\right) = \dim\!\left(\mathcal{H}(\mathcal{O}(H))\right)$. The last claim is a general fact, see (\ref{HeisenbergMat}).
\end{proof}


\begin{remark}\label{remarkProductGaugeFields}
Due to Proposition \ref{isoL10Heisenberg}, there is an isomorphism of vector spaces $f : \mathcal{L}_{1,0}(H) \to H^* \otimes H^*$ given by $\overset{I}{B}{^i_j} \overset{J}{A}{^k_l} \mapsto \overset{I}{T}{^i_j} \otimes \overset{J}{T}{^k_l}$. We define a $H$-module-algebra structure on $H^* \otimes H^*$, denoted by $\mathcal{F}_{1,0}(H)$ and with product $\ast$, by requiring $f$ to be an isomorphism of (right) $H$-module-algebras. The right $H$-action is
\[ (\varphi \otimes \psi) \cdot h = \varphi\!\left( h' ? S(h'') \right) \otimes \psi\!\left( h''' ? S(h^{(4)}) \right). \]
For $\varphi \in H^*$, let $\varphi_b = \varphi \otimes \varepsilon$ and $\varphi_a = \varepsilon \otimes \varphi$. It is clear that
\begin{equation}\label{propAstPhiAB}
\varphi_b \ast \psi_a = \varphi \otimes \psi, \:\:\:\:\:\: \varphi_a \ast \psi_a = (\varphi \ast \psi)_a, \:\:\:\:\:\: \varphi_b \ast \psi_b = (\varphi \ast \psi)_b
\end{equation}
(recall the algebra $\mathcal{F}_{0,1}(H)$ defined in Remark \ref{F01}). Moreover, using \eqref{echangeL10Inverse}:
\begin{align*}
(\overset{I}{T}{^{\alpha}_{\beta}})_a \ast  (\overset{J}{T}{^{\gamma}_{\delta}})_b &= f\!\left(\overset{I}{A}_1\overset{J}{B}_2\right)^{\alpha \gamma}_{\beta \delta} = f\!\left(\overset{J}{(a_i)}_2\overset{IJ}{(R')}_{12}\overset{J}{B}_2\overset{IJ}{R}_{12}\overset{I}{A}_1\overset{IJ}{(R')}_{12}\overset{I}{S(b_i)}_1\right)^{\alpha \gamma}_{\beta \delta}\\
&= f\!\left( \overset{J}{a_i} \overset{J}{a_j} \overset{J}{B} \overset{J}{b_k} \overset{J}{a_l} \right)^{\gamma}_{\delta} f\!\left( \overset{I}{b_j} \overset{I}{a_k} \overset{I}{A} \overset{I}{b_l} \overset{I}{S(b_i)} \right)^{\alpha}_{\beta} = \overset{J}{T}{^{\gamma}_{\delta}}\!\left( a_i a_j ? b_k a_l \right) \otimes \overset{I}{T}{^{\alpha}_{\beta}}\!\left( b_j a_k ? b_l S(b_i) \right)
\end{align*}
and it follows that
\begin{equation}\label{produitGaugeFields10}
\varphi_a \ast \psi_b = \psi(a_i a_j ? b_k a_l)_b \ast \varphi(b_j a_k ? b_l S(b_i))_a.
\end{equation}
Combining \eqref{propAstPhiAB} and \eqref{produitGaugeFields10}, we obtain the general formula:
\begin{equation*}
\begin{split} 
&(\varphi_1 \otimes \psi_1) \ast (\varphi_2 \otimes \psi_2) = \varphi_1\bigl(? b_m S(b_n)\bigr) \varphi_2\bigl(a_i a_j a_n ? a_m b_k a_l \bigr) \otimes \psi_1\bigl( b_j a_k ? b_o S(b_p) b_l S(b_i) \bigr) \psi_2\bigl(a_p ? a_o\bigr), \\
&x_b \otimes x_a \mapsto \varphi_1\bigl(x_b' b_m S(b_n)\bigr) \varphi_2\bigl( a_i a_j a_n x_b'' a_m b_k a_l \bigr) \otimes \psi_1\bigl( b_j a_k x_a' b_o S(b_p) b_l S(b_i) \bigr) \psi_2\bigl(a_p x_a'' a_o\bigr).
\end{split}
\end{equation*}
This is the product of the functions $\varphi_1 \otimes \psi_1, \varphi_2 \otimes \psi_2 \in \mathcal{F}_{1,0}(H)$ and its evaluation on the discrete connection which assigns $x_b$ to the loop $b$ and $x_a$ to the loop $a$, see Figure \ref{Sigma10AvecConnexion} and the Introduction.
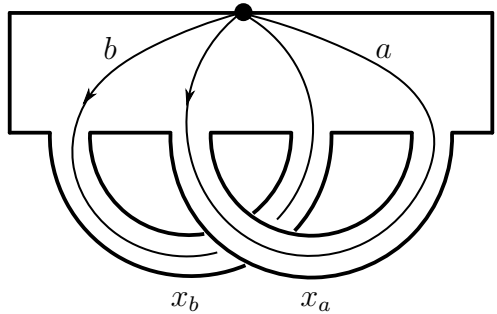
\begin{figure}[!h]
\centering
\begingroup%
  \makeatletter%
  \providecommand\color[2][]{%
    \errmessage{(Inkscape) Color is used for the text in Inkscape, but the package 'color.sty' is not loaded}%
    \renewcommand\color[2][]{}%
  }%
  \providecommand\transparent[1]{%
    \errmessage{(Inkscape) Transparency is used (non-zero) for the text in Inkscape, but the package 'transparent.sty' is not loaded}%
    \renewcommand\transparent[1]{}%
  }%
  \providecommand\rotatebox[2]{#2}%
  \newcommand*\fsize{\dimexpr\f@size pt\relax}%
  \newcommand*\lineheight[1]{\fontsize{\fsize}{#1\fsize}\selectfont}%
  \ifx\svgwidth\undefined%
    \setlength{\unitlength}{181.43750252bp}%
    \ifx\svgscale\undefined%
      \relax%
    \else%
      \setlength{\unitlength}{\unitlength * \real{\svgscale}}%
    \fi%
  \else%
    \setlength{\unitlength}{\svgwidth}%
  \fi%
  \global\let\svgwidth\undefined%
  \global\let\svgscale\undefined%
  \makeatother%
  \begin{picture}(1,0.64063813)%
    \lineheight{1}%
    \setlength\tabcolsep{0pt}%
    \put(0.19668408,0.52418621){\color[rgb]{0,0,0}\makebox(0,0)[lt]{\lineheight{1.25}\smash{\begin{tabular}[t]{l}$b$\end{tabular}}}}%
    \put(0.7556463,0.52712437){\color[rgb]{0,0,0}\makebox(0,0)[lt]{\lineheight{1.25}\smash{\begin{tabular}[t]{l}$a$\end{tabular}}}}%
    \put(0,0){\includegraphics[width=\unitlength,page=1]{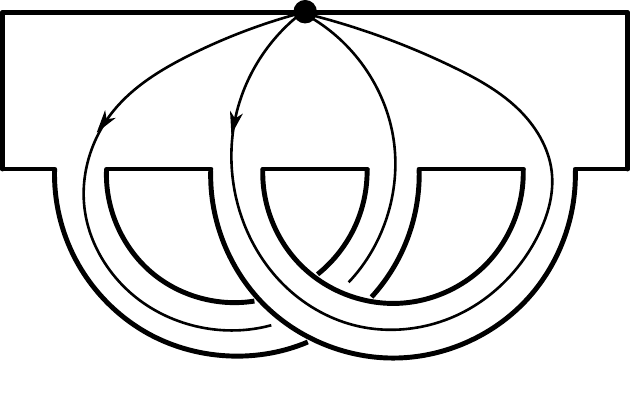}}%
    \put(0.3352267,0.01315116){\color[rgb]{0,0,0}\makebox(0,0)[lt]{\lineheight{1.25}\smash{\begin{tabular}[t]{l}$x_b$\end{tabular}}}}%
    \put(0.60209499,0.01299841){\color[rgb]{0,0,0}\makebox(0,0)[lt]{\lineheight{1.25}\smash{\begin{tabular}[t]{l}$x_a$\end{tabular}}}}%
  \end{picture}%
\endgroup%

\caption{A discrete connection $x_b \otimes x_a \in H^{\otimes 2}$ on $\Sigma_{1,0}^{\mathrm{o}}$.}
\label{Sigma10AvecConnexion}
\end{figure}
\finEx
\end{remark}

\subsection{Representation of $\mathcal{L}_{1,0}^{\mathrm{inv}}(H)$ on $\SLF(H)$}\label{sectionRepInvariants}
\indent In this section we construct representations of the subalgebra of invariants $\mathcal{L}^{\mathrm{inv}}_{1,0}(H)$. This will be extended to any $g, n$ in the next chapter.

\smallskip

\indent Recall from (\ref{repHO}) that there is a faithful representation $\triangleright$ of $\mathcal{H}(\mathcal{O}(H))$ on $\mathcal{O}(H)$. Using the isomorphism $\Psi_{1,0}$, we get a representation of $\mathcal{L}_{1,0}(H)$ on $\mathcal{O}(H)$, still denoted $\triangleright$:
\begin{equation}\label{defActionL10SurDual}
\forall\, x \in  \mathcal{L}_{1,0}(H), \: \forall\, \psi \in \mathcal{O}(H), \:\: x \triangleright \psi = \Psi_{1,0}(x) \triangleright \psi. 
\end{equation}
Using \eqref{echangeHeisenberg}, it is easy to get:
\begin{equation}\label{actionL10}
\overset{I}{A}_1 \triangleright \overset{J}{T}\!_2 = \overset{J}{T}\!_2 \overset{IJ}{(RR')}_{12} \:\:\:\: \text{ and } \:\:\:\: \overset{I}{B}_1 \triangleright \overset{J}{T}\!_2 = \overset{I}{(a_i)}_1\overset{I \otimes J}{T}\!\!\!_{12} \overset{I \otimes J}{(b_i)}_{12}\overset{IJ}{(R')}_{12} = \overset{I}{(a_ia_j)}_1\overset{I \otimes J}{T}\!\!\!_{12} \overset{I}{(b_j)}_{1}\overset{J}{(b_i)}_{2}\overset{IJ}{(R')}_{12}
\end{equation}
where as usual $R = a_i \otimes b_i$ and the last equality is obtained using (\ref{deltaR}).

\smallskip

Let us define matrices (see section \ref{heisenbergDouble} for the definition of $\overset{I}{\widetilde{L}}{^{(\pm)}}$)
\begin{equation}\label{matriceCL10}
\overset{I}{C} = \overset{I}{v}{^{2}}\overset{I}{B}\overset{I}{A}{^{-1}}\overset{I}{B}{^{-1}}\overset{I}{A}, \:\:\:\:\: \overset{I}{C}{^{(\pm)}} = \Psi_{1,0}^{-1}(\overset{I}{L}{^{(\pm)}}\overset{I}{\widetilde{L}}{^{(\pm)}})
\end{equation}
Observe that geometrically, $\overset{I}{C}$ corresponds to the boundary of the surface $\Sigma_{1,0}^{\mathrm{o}}$, see \eqref{boundaryLoop10} and Figure \ref{Sigma10}.
\begin{lemma}\label{decGauss}
It holds:
\[ \overset{I}{C} = \overset{I}{C}{^{(+)}}\overset{I}{C}{^{(-)-1}}. \]
Moreover, the matrices $\overset{I}{C}$ satisfy the fusion relation of $\mathcal{L}_{0,1}(H)$:
\[ \overset{I \otimes J}{C}\!\!_{12} = \overset{I}{C}_1 \, \overset{IJ}{(R')}_{12} \, \overset{J}{C}_2 \, \overset{IJ}{(R')}{^{-1}_{12}}. \]
\end{lemma}
\begin{proof}
We have
\[ \Psi_{1,0}\!\left(\overset{I}{v}{^{2}}\overset{I}{B}\overset{I}{A}{^{-1}}\overset{I}{B}{^{-1}}\overset{I}{A}\right) = \overset{I}{L}{^{(+)}} \left(\overset{I}{v}{^{2}} \overset{I}{T}\overset{I}{L}{^{(+)-1}}\overset{I}{L}{^{(-)}} S(\overset{I}{T}) \right)\overset{I}{L}{^{(-)-1}}. \]
Let us simplify the middle term. It is equal to:
\begin{align*}
\overset{I}{v}{^{2}} \overset{I}{T} \overset{I}{S(a_i)} \overset{I}{S^{-1}(b_j)} b_i a_j S(\overset{I}{T}) &= \overset{I}{v}{^{2}} \overset{I}{T} \overset{I}{S(a_i)} \overset{I}{S^{-1}(b_j)} \overset{I}{S(a'_j)} \overset{I}{S(b'_i)} S(\overset{I}{T}) b''_i a''_j \\
&= \overset{I}{v}{^{2}}\overset{I}{T} \overset{I}{S(a_ia_k)}\overset{I}{S^{-1}(b_j b_{\ell})} \overset{I}{S(a_j)} \overset{I}{S(b_k)} S(\overset{I}{T}) b_ia_{\ell} = \overset{I}{T} \overset{I}{S(b_{\ell}a_i)} S(\overset{I}{T}) a_{\ell} b_i.
\end{align*}
The first equality is the exchange relation \eqref{relDefHeisenberg} in $\mathcal{H}(\mathcal{O}(H))$ and the second follows from the properties of the $R$-matrix. The third equality is obtained as follows: denoting $m : H \otimes H \to H$ the multiplication, we can write 
\begin{align*}
&v^2 S(a_ia_k) S^{-1}(b_{\ell}) S^{-1}(b_j) S(a_j) S(b_k) \otimes b_ia_{\ell} = v S(a_ia_k) S^{-1}(b_k b_{\ell}) g^{-1} \otimes b_ia_{\ell}\\
= \:&v S(a_k a_i) S^{-1}(b_{\ell} b_k) g^{-1} \otimes a_{\ell} b_i = v S(a_i) S\!\left(S^{-2}(b_k)a_k\right) S^{-1}(b_{\ell})g^{-1} \otimes a_{\ell} b_i = S(a_i) S(b_{\ell}) \otimes a_{\ell}b_i.
\end{align*}
For the second equality, we used the Yang-Baxter relation $R_{13} R_{12} R_{32} = R_{32} R_{12} R_{13}$; the others equalities follows from \eqref{elementDrinfeld} and the standard properties for $g$ and $v$. Now, we have:
\[ \overset{I}{T}_1 \overset{I}{S(b_{\ell} a_i)}_1 S(\overset{I}{T})_1 a_{\ell} b_i \triangleright \overset{J}{T}_2 = \overset{I}{T}_1 \overset{I}{S(b_{\ell} a_i)}_1 S(\overset{I}{T})_1 \overset{J}{T}_2 \overset{J}{(a_{\ell} b_i)}_2
= \overset{I}{S(b_{\ell} a_i)}_1 \overset{J}{(a_{\ell} b_i)}_2 \overset{J}{T}_2 = \overset{I}{(a_i b_{\ell})}_1 \widetilde{b_i} \widetilde{a_{\ell}} \triangleright \overset{J}{T}_2. \]
For the second equality, we used that for any $h \in H$:
\[ \langle \overset{I}{S(b_{\ell} a_i)}_1 S(\overset{I}{T})_1 \overset{J}{T}_2 \overset{J}{(a_{\ell} b_i)}_2, h \rangle = \overset{I}{S(h'b_{\ell} a_i)}_1 \overset{J}{(h''a_{\ell} b_i)}_2 = \overset{I}{S(b_{\ell} a_ih')}_1 \overset{J}{(a_{\ell} b_ih'')}_2 = \langle S(\overset{I}{T})_1 \overset{I}{S(b_{\ell} a_i)}_1  \overset{J}{(a_{\ell} b_i)}_2 \overset{J}{T}_2, h \rangle. \]
Since $\triangleright$ is faithful, we finally get 
\[ \overset{I}{v}{^{2}} \overset{I}{T} \overset{I}{S(a_i)} \overset{I}{S^{-1}(b_j)} b_i a_j S(\overset{I}{T}) = \overset{I}{(a_{\ell}b_i)} \widetilde{b_{\ell}} \widetilde{a_i} = \overset{I}{\widetilde{L}}{^{(+)}}\overset{I}{\widetilde{L}}{^{(-)-1}}.\]
Hence 
\[ \Psi_{1,0}(\overset{I}{C}) = \overset{I}{L}{^{(+)}}\overset{I}{\widetilde{L}}{^{(+)}}(\overset{I}{L}{^{(-)}}\overset{I}{\widetilde{L}}{^{(-)}})^{-1} 
= \Psi_{1,0}(\overset{I}{C}{^{(+)}}\overset{I}{C}{^{(-)-1}}) \]
as desired. Now, consider the morphism of algebras 
\[\fonc{f}{H}{\mathcal{H}(\mathcal{O}(H))}{h}{\widetilde{h'}h''}\]
and observe that 
\begin{equation}\label{expressionC10}
\Psi_{1,0}( \overset{I}{C}{^{(+)}}) = \overset{I}{a_i}\:\widetilde{b_i'} b_i'', \:\:\:\:\:\:\:\: \Psi_{1,0}( \overset{I}{C}{^{(-)}} ) = \overset{I}{S^{-1}(b_i)} \: \widetilde{a_i'} a_i'', \:\:\:\:\:\:\:\: \Psi_{1,0}( \overset{I}{C} ) = \overset{I}{X_i}\: \widetilde{Y_i'} Y_i'',
\end{equation}
where $X_i \otimes Y_i = RR'$. Then we have a morphism
\[\begin{array}{ccccccc}
\mathcal{L}_{0,1}(H) &\overset{\Psi_{0,1}}{\longrightarrow}& H &\overset{f}{\longrightarrow}& \mathcal{H}(\mathcal{O}(H)) &\overset{\Psi_{1,0}^{-1}}{\longrightarrow}& \mathcal{L}_{1,0}(H)\\
\overset{I}{M} &\longmapsto& \overset{I}{X_i}Y_i &\longmapsto& \overset{I}{X_i}\: \widetilde{Y_i'} Y_i'' &\longmapsto& \overset{I}{C}.
\end{array}\]
It follows that $\overset{I}{C}$ satisfies the fusion relation. Alternatively, one can write $\Psi_{1,0}(\overset{I}{C}) = \overset{I}{L}{^{(+)}}\overset{I}{\widetilde{L}}{^{(+)}}\overset{I}{\widetilde{L}}{^{(-)-1}}\overset{I}{L}{^{(-)-1}}$ and check the fusion relation directly using \eqref{propertiesL} and \eqref{LTilde}.
\end{proof}

\noindent Thanks to Lemmas \ref{injectionFusion10} and \ref{decGauss}, we have a morphism
\[ \fonc{j_C}{\mathcal{L}_{0,1}(H)}{\mathcal{L}_{1,0}(H)}{\overset{I}{M}}{\overset{I}{C}}. \]
The previous proof shows that $j_C = \Psi_{1,0}^{-1} \circ f \circ \Psi_{0,1}$. Moreover, the algebra generated by the coefficients $\overset{I}{C^{(\pm)}}{^i_j}$ equals the vector space generated by the coefficients $\overset{I}{C}{^i_j}$:
\[ \mathbb{C}\langle \overset{I}{C^{(\pm)}}{^i_j} \rangle_{I,i,j} = \mathrm{vect}(\overset{I}{C}{^i_j})_{I,i,j}. \]
Indeed, since $H$ is factorisable, it is generated as an algebra by the coefficients $\overset{I}{L}{^{(\pm)}}^i_j = (\overset{I}{a_l^{(\pm)}})^i_j \, b_l^{(\pm)}$ and as a vector space by the coefficients $(\overset{I}{X_l})^i_j \, Y_l$. The claim follows from $\overset{I}{C^{(\pm)}}{^i_j} = j_C(\overset{I}{a_l^{(\pm)}}{^i_j} \, b_l^{(\pm)})$ and $\overset{I}{C}{^i_j} = j_C((\overset{I}{X_l})^i_j \, Y_l)$.

\begin{lemma}\label{conjugaisonC10}
It holds
\[ \overset{I}{C}{^{(\pm)}_1} \, \overset{J}{U}_2 \, \overset{I}{C}{^{(\pm)-1}_1} = \overset{IJ}{R}{^{(\pm)-1}_{12}} \, \overset{J}{U}_2 \, \overset{IJ}{R}{^{(\pm)}_{12}} \]
where $U$ is $A$ or $B$.
\end{lemma}
\begin{proof}
We use the isomorphism $\Psi_{1,0}$ together with relations \eqref{propertiesL}, \eqref{echangeHeisenberg} and \eqref{LTilde}:
\begin{align*}
\Psi_{1,0}\!\left(C^{(\pm)}_1 A_2 C^{(\pm)-1}_1\right) &= L^{(\pm)}_1 \widetilde{L}^{(\pm)}_1 L^{(+)}_2 L^{(-)-1}_2 \widetilde{L}^{(\pm)-1}_1 L^{(\pm)-1}_1 = L^{(\pm)}_1 L^{(+)}_2 L^{(-)-1}_2 L^{(\pm)-1}_1\\
& = R^{(\pm)-1}_{12} L^{(+)}_2 L^{(\pm)}_1 R^{(\pm)}_{12} L^{(-)-1}_2 L^{(\pm)-1}_1 = R^{(\pm)-1}_{12} L^{(+)}_2 L^{(-)-1}_2 R^{(\pm)}_{12}\\
& = \Psi_{1,0}\!\left(R^{(\pm)-1}_{12} A_2 R^{(\pm)}_{12}\right)
\end{align*}
and
\begin{align*}
\Psi_{1,0}\!\left(C^{(\pm)}_1 B_2 C^{(\pm)-1}_1\right) &= L^{(\pm)}_1 \widetilde{L}^{(\pm)}_1 L^{(+)}_2 T_2 L^{(-)-1}_2 \widetilde{L}^{(\pm)-1}_1 L^{(\pm)-1}_1 = L^{(\pm)}_1 L^{(+)}_2 \widetilde{L}^{(\pm)}_1 T_2 \widetilde{L}^{(\pm)-1}_1 L^{(-)-1}_2 L^{(\pm)-1}_1\\
&= L^{(\pm)}_1 L^{(+)}_2 R^{(\pm)-1}_{12} T_2  L^{(-)-1}_2 L^{(\pm)-1}_1 = R^{(\pm)-1}_{12} L^{(+)}_2 L^{(\pm)}_1 T_2  L^{(-)-1}_2 L^{(\pm)-1}_1\\
& = R^{(\pm)-1}_{12} L^{(+)}_2 T_2 L^{(\pm)}_1 R^{(\pm)}_{12} L^{(-)-1}_2 L^{(\pm)-1}_1 = R^{(\pm)-1}_{12} L^{(+)}_2 T_2  L^{(-)-1}_2 R^{(\pm)}_{12}\\
& = \Psi_{1,0}\!\left(R^{(\pm)-1}_{12} B_2 R^{(\pm)}_{12}\right).
\end{align*}
The subscript $1$ (resp. $2$) implicitly means evaluation in a representation $I$ (resp. $J$).
\end{proof}

Thanks to Lemma \ref{injectionFusion10}, we can define an action of $H$ (identified with $\mathcal{L}_{0,1}(H)$ through $\Psi_{0,1}$) on $H^*$ by
\begin{equation}\label{actionHTriangle}
h \cdot \varphi = h_C \triangleright \varphi. 
\end{equation}

\begin{lemma}\label{CSLF}
The action \eqref{actionHTriangle} of $H$ on $H^*$ is
\[ h \cdot \varphi = \varphi(S^{-1}(h') ? h''). \]
It follows that we have the equivalence
\begin{equation*}
\varphi \in \mathrm{SLF}(H) \:\:\:\: \iff \:\:\:\: \overset{I}{C} \triangleright \varphi = \mathbb{I}_{\dim(I)}\varphi.
\end{equation*}
\end{lemma}
\begin{proof}
Since $H$ is factorizable, we can assume that $h = (\overset{I}{X_i})^a_b \, Y_i$, where $RR' = X_i \otimes Y_i$. Due to \eqref{defActionL10SurDual} and \eqref{expressionC10}, we obtain
\[ h \cdot \varphi = (\overset{I}{X_i})^a_b \, (Y_i)_C \triangleright \varphi = \overset{I}{C}{^a_b} \triangleright \varphi = (\overset{I}{X_i})^a_b\, \widetilde{Y_i'} Y_i'' \triangleright \varphi = \varphi\!\left((\overset{I}{X_i})^a_b S^{-1}(Y_i') ? Y_i''  \right) = \varphi\!\left( S^{-1}(h') ? h'' \right) \]
as desired. Next, it is easy to see that $\varphi \in \mathrm{SLF}(H)$ if, and only if, $h \cdot \varphi = \varepsilon(h)\varphi$ for all $h \in H$. Applying this to $h = (\overset{I}{X_i})^a_b \, Y_i$ and using that $X_i \varepsilon(Y_i) = 1$, we find that $\varphi \in \mathrm{SLF}(H)$ if, and only if, $\overset{I}{C}{^a_b} \triangleright \varphi = \delta^a_b \varphi$ for all $I,a,b$.
\end{proof}
\noindent It follows from this lemma that, in the case of the torus, $\mathrm{SLF}(H)$ implements the flatness constraint \eqref{flatnessConstraintMatrices} discussed in the Introduction.
\begin{theorem}\label{repInv}
1) An element $x \in \mathcal{L}_{1,0}(H)$ is invariant under the action of $H$ (or equivalently under the coaction $\Omega$ of $\mathcal{O}(H)$) if, and only if, for every $H$-module $I$, $\overset{I}{C}x = x\overset{I}{C}$. \\
2) The restriction of $\triangleright$ to $\mathcal{L}_{1,0}^{\text{\em inv}}(H)$ leaves the subspace $\SLF(H) \subset H^*$ stable:
$$ \forall\, x \in  \mathcal{L}^{\text{\em inv}}_{1,0}(H), \: \forall\, \psi \in \SLF(H), \:\: x \triangleright \psi \in \SLF(H). $$
Hence, we have a representation of $\mathcal{L}_{1,0}^{\text{\em inv}}(H)$ on $\SLF(H)$. We denote it $\rho_{\mathrm{SLF}}$.
\end{theorem}
\begin{proof}
1) Letting $U$ be $A$ or $B$, $R^{(\pm)} = a_i^{(\pm)} \otimes b_i^{(\pm)}$ and using Lemma \ref{conjugaisonC10}, we have that the right action $ \cdot$ of $H$ on $\mathcal{L}_{1,0}(H)$ satisfies:
\begin{align*}
\overset{J}{U}_2 \cdot  \overset{I}{L}{^{(\pm)-1}_1} &= \overset{J}{U}_2 \cdot S^{-1}(b_i^{(\pm)}) \overset{I}{(a_i^{(\pm)})}_1 = \overset{J}{S^{-1}(b_i^{(\pm)}{''})}_2 \overset{J}{U}_2 \overset{J}{(b_i^{(\pm)}{'})}_2 \overset{I}{(a_i^{(\pm)})}_1\\
& = \overset{J}{S^{-1}(b_i^{(\pm)})}_2 \overset{J}{U}_2 \overset{J}{(b_j^{(\pm)})}_2 \overset{I}{(a_i^{(\pm)}a_j^{(\pm)})}_1 = \overset{IJ}{R}{^{(\pm)-1}_{12}} \, \overset{J}{U}_2 \, \overset{IJ}{R}{^{(\pm)}_{12}} = \overset{I}{C}{^{(\pm)}_1} \, \overset{J}{U}_2 \, \overset{I}{C}{^{(\pm)-1}_1}.
\end{align*}
We have thus shown that
\[ (\overset{J}{U})^c_d \cdot  S^{-1}(\overset{I}{L}{^{(\pm)}})^a_b = (\overset{I}{C}{^{(\pm)}})^a_i \, (\overset{J}{U})^c_d \, (\overset{I}{C}{^{(\pm)-1}})^i_b \]
or in other words
\[ \forall \,x \in \mathcal{L}_{1,0}(H), \:\:\: x \cdot  S^{-1}(\overset{I}{L}{^{(\pm)}}) = \overset{I}{C}{^{(\pm)}} \, x \, \overset{I}{C}{^{(\pm)-1}} \]
Since $H$ is factorizable, the elements $S^{-1}(\overset{I}{L}{^{(\pm)}})^a_b$ generate $H$ as an algebra. Hence the previous equation shows that $x$ is an invariant element if, and only if, it commutes with the cofficients of the matrices $\overset{I}{C}{^{(\pm)}}$. As remarked above, the algebra generated by the coefficients $\overset{I}{C}{^{(\pm)}}{^i_j}$ equals the algebra generated by the coefficients $\overset{I}{C}{^i_j}$. Hence, an element is invariant if, and only if, it commutes with the coefficients of the matrices $\overset{I}{C}$.
\\2) Let $x \in \mathcal{L}_{1,0}^{\mathrm{inv}}(H)$ and $\varphi \in \mathrm{SLF}(H)$, then
\[ \overset{I}{C} \triangleright (x \triangleright \varphi) = (\overset{I}{C}x) \triangleright \varphi = (x\overset{I}{C}) \triangleright \varphi = x \triangleright (\overset{I}{C} \triangleright \varphi) = \mathbb{I}_{\dim(I)}(x \triangleright \varphi) \]
and it follows that $x \triangleright \varphi \in \mathrm{SLF}(H)$ thanks to Lemma \ref{CSLF}.
\end{proof}

\indent We now need to determine explicit formulas for the representation of particular types of invariant elements that will appear in the proof of the modular identities in section \ref{sectionModulaire}. If $\psi \in H^*$ and $a \in H$, we define:
$$ \psi^a = \psi(a?) $$
where $\psi(a?) : x \mapsto \psi(ax)$. This defines a right representation of $H$ on $H^*$. Obviously, if $z \in  \mathcal{Z}(H)$ and $\psi \in \SLF(H)$ then $\psi^z \in \SLF(H)$.
\\\indent Recall that $z_A = j_A(z)$ (resp. $z_B = j_B(z)$) is the image of $z \in \mathcal{L}_{0,1}(H)$ by the map $j_A(\overset{I}{M}) = \overset{I}{A}$ (resp. $j_B(\overset{I}{M}) = \overset{I}{B}$). See (\ref{notationImbed}) for the general definition.
\begin{proposition}\label{actionAB}
Let $z \in \mathcal{L}_{0,1}^{\mathrm{inv}}(H) = \mathcal{Z}(H)$ and let $\psi \in \mathrm{SLF}(H)$. Then:
$$ z_A \triangleright \psi = \psi^z \:\:\:\: \text{ and } \:\:\:\: z_B \triangleright \psi = \left(\mathcal{D}^{-1}(z)\psi^v\right)^{v^{-1}} $$
where $\mathcal{D}$ is the isomorphism defined in (\ref{morphismeDrinfeld}).
\end{proposition}
\begin{proof}
The first relation is obvious.  
For the second formula, we write $z_B = \sum_I \mathrm{tr}(\Lambda_I\overset{I}{g}\overset{I}{B})$ with $\mathcal{D}^{-1}(z) = \sum_I \mathrm{tr}(\Lambda_I\overset{I}{T}) \in \SLF(H)$ by Lemma \ref{writingInvariants}. We also write $\psi = \sum_J \mathrm{tr}(\Theta_J\overset{J}{T})$. Then, using (\ref{actionL10}):
\begin{align*}
z_B \triangleright \psi &= \sum_{I,J} \mathrm{tr}_{12}\!\left((\Lambda_I)_1(\Theta_J)_2\,\overset{I}{g}_1\overset{I}{B}_1 \triangleright \overset{J}{T}_2\right)\\
&=\sum_{I,J} \mathrm{tr}_{12}\!\left((\Lambda_I)_1(\Theta_J)_2\,\overset{I}{g}_1\overset{I}{(a_ia_j)}_1\overset{I \otimes J}{T}\!\!\!_{12} \overset{I}{(b_j)}_{1}\overset{J}{(b_i)}_{2}\overset{IJ}{(R')}_{12}\right)\\
&= \mathcal{D}^{-1}(z)\!\left( g a_i a_j ? b_j b_k \right) \psi\!\left( ?b_ia_k \right) = \mathcal{D}^{-1}(z)\!\left( ? b_j b_k  S^2(a_ia_j) g\right) \psi\!\left( ?b_ia_k \right)
\end{align*}
with $\mathrm{tr}_{12} = \mathrm{tr} \otimes \mathrm{tr}$, $R = a_i \otimes b_i$. Thanks to the Yang-Baxter equation, we have:
$$ b_j b_k \otimes a_ia_j \otimes b_ia_k = R_{23} R_{21} R_{31} = R_{31} R_{21} R_{23} = b_ib_j \otimes a_ja_k \otimes a_ib_k. $$
It follows that:
\begin{align*}
z_B \triangleright \psi &= \mathcal{D}^{-1}(z)\!\left( ? b_i b_j  S^2(a_j a_k) g\right) \psi\!\left( ? a_ib_k \right)\\ &= \mathcal{D}^{-1}(z)\!\left(? v^{-1} b_i a_k\right) \psi\!\left( ? a_ib_k \right) = \mathcal{D}^{-1}(z)\!\left(? (v^{-1})'\right) \psi\!\left( ? v (v^{-1})'' \right)
\end{align*}
where we used \eqref{elementDrinfeld}, \eqref{pivotCan} and \eqref{ribbon}. Hence for $x \in H$:
$$  \left(z_B \triangleright \psi\right)\!(x) = \mathcal{D}^{-1}(z)\!\left((v^{-1})'x'\right) \psi\!\left(v(v^{-1})''x''\right) = \left(\mathcal{D}^{-1}(z)\,\psi^v\right)\!(v^{-1}x) = \left(\mathcal{D}^{-1}(z)\,\psi^v\right)^{v^{-1}}\!(x) $$
as desired.
\end{proof}

\begin{lemma}\label{lemmaBMoinsUn}
Let $z \in \mathcal{L}_{0,1}^{\mathrm{inv}}(H) = \mathcal{Z}(H)$ and let $\psi \in \SLF(H)$. Then:
$$ z_{B^{-1}} \triangleright \psi = \left( S\!\left(\mathcal{D}^{-1}(z)\right) \psi^v \right)^{v^{-1}}. $$
It follows that if $S(\psi) = \psi$ for all $\psi \in \SLF(H)$, then $\rho_{\mathrm{SLF}}(z_{B^{-1}}) = \rho_{\mathrm{SLF}}(z_B)$.
\end{lemma}
\begin{proof}
This proof is quite similar to that of the previous proposition. Using the fact that $\Psi_{1,0}(\overset{I}{B}\,\!^{-1}) = \overset{I}{L}\,\!^{(-)} S(\overset{I}{T}) \overset{I}{L}\,\!^{(+)-1} $ together with Lemma \ref{echangeHeisenberg} and formulas (\ref{antipodeT}), (\ref{propSR}) and (\ref{deltaR}), it is not too difficult to show that
$$ \overset{I}{B}\,\!^{-1}_1 \triangleright \overset{J}{T}\!_2 = \exposantGauche{t \otimes \mathrm{id}\!\!}{\left( \overset{I^*}{(a_i)}_1 \overset{I^* \otimes J}{T}\!\!\!\!_{12} \,\overset{I^*}{(a_jS^{-2}(b_jb_k))}_1\overset{J}{(a_kb_i)}_2 \right)} $$
where $^{t \otimes \mathrm{id}}$ means transpose on the first tensorand.  Write $z_{B^{-1}} = \sum_I \mathrm{tr}(\Lambda_I\overset{I}{g}\overset{I}{B}\,\!^{-1})$ with $\mathcal{D}^{-1}(z) = \sum_I \mathrm{tr}(\Lambda_I\overset{I}{T}) \in \SLF(H)$ by Lemma \ref{writingInvariants}, and $\psi = \sum_J \mathrm{tr}(\Theta_J\overset{J}{T})$. Observe using (\ref{antipodeT}) that:
$$ S\!\left(\mathcal{D}^{-1}(z)\right) = \sum_I \mathrm{tr}\!\left( \Lambda_I S(\overset{I}{T})\right) = \sum_I \mathrm{tr}\!\left( \exposantGauche{t}{\Lambda_I} \overset{I^*}{T}\right). $$
Using the fact that $S(g) = g^{-1}$ and (\ref{antipodeT}), we thus get:
\begin{align*}
z_{B^{-1}} \triangleright \psi &= \sum_{I,J} \mathrm{tr}_{12}\!\left( (\Lambda_I\overset{I}{g})_1(\Theta_J)_2\, \exposantGauche{t \otimes \mathrm{id}\!\!}{\left( \overset{I^*}{(a_i)}_1 \overset{I^* \otimes J}{T}\!\!\!\!_{12} \,\overset{I^*}{(a_jS^{-2}(b_jb_k))}_1\overset{J}{(a_kb_i)}_2 \right)} \right)\\
&= \sum_{I,J} \mathrm{tr}_{12}\!\left( (\exposantGauche{t}{\Lambda_I})_1(\Theta_J)_2\, \overset{I^*}{(a_i)}_1 \overset{I^* \otimes J}{T}\!\!\!\!_{12} \,\overset{I^*}{(a_jS^{-2}(b_jb_k)g^{-1})}_1\overset{J}{(a_kb_i)}_2 \right)\\
&= S\!\left(\mathcal{D}^{-1}(z)\right)\!\left( a_i\,?\,a_jS^{-2}(b_jb_k)g^{-1} \right) \, \psi\!\left(?\,a_kb_i\right)\\
&= S\!\left(\mathcal{D}^{-1}(z)\right)\!\left( ?\,a_jS^{-2}(b_jb_k)g^{-1}a_i \right) \, \psi\!\left(?\,a_kb_i\right) = S\!\left(\mathcal{D}^{-1}(z)\right)\!\left( ?(v^{-1})'\right) \, \psi\!\left(?\,v (v^{-1})''\right).
\end{align*}
For the last equality we used 
\eqref{elementDrinfeld}, \eqref{pivotCan} and \eqref{ribbon}. Hence we get as in the previous proof $ z_{B^{-1}} \triangleright \psi =  \left(S\!\left(\mathcal{D}^{-1}(z)\right)\psi^v\right)^{v^{-1}}$.
\end{proof}

\section{Projective representation of $\mathrm{SL}_2(\mathbb{Z})$}\label{sectionModulaire}
\indent As previously, $H$ is a finite dimensional factorizable ribbon Hopf algebra.

\subsection{Mapping class group of the torus}
\indent First, recall the general definition (see \cite{FM}).
\begin{definition}\label{defMCG}
If $S$ is a compact oriented surface, we denote by $\mathrm{MCG}(S)$ its mapping class group, that is the group of isotopy classes of orientation preserving homeomorphisms of $S$ which fix the boundary pointwise.
\end{definition}
\indent For this chapter we focus on the torus $\Sigma_{1,0} = S^1 \times S^1$. Let $\Sigma_{1,0}^{\mathrm{o}} = \Sigma_{1,0} \setminus D$, where $D$ is an embedded open disk. The surface $\Sigma_{1,0}^{\mathrm{o}}$ together with the canonical curves $a$ and $b$ are represented in Figure \ref{Sigma10}. The groups $\mathrm{MCG}(\Sigma_{1,0}^{\mathrm{o}})$ and $\mathrm{MCG}(\Sigma_{1,0})$ are generated by the Dehn twists $\tau_a, \tau_b$ about the free homotopy classes of the curves $a$ and $b$.  It is well-known (see \cite{FM}) that
\begin{align*}
\mathrm{MCG}(\Sigma_{1,0}^{\mathrm{o}}) &= B_3 = \left\langle \tau_a, \tau_b \, \left\vert \, \tau_a\tau_b\tau_a = \tau_b\tau_a\tau_b \right. \right\rangle,\\
\mathrm{MCG}(\Sigma_{1,0}) &= \mathrm{SL}_2(\mathbb{Z}) = \left\langle \tau_a, \tau_b \, \left\vert \, \tau_a\tau_b\tau_a = \tau_b\tau_a\tau_b, \: (\tau_a\tau_b)^6=1 \right. \right\rangle. 
\end{align*}
This presentation is not the usual one of $\mathrm{SL}_2(\mathbb{Z})$, which is:
\[ \mathrm{SL}_2(\mathbb{Z}) = \left\langle s, t \, \left\vert \, (st)^3=s^2, \: s^4=1 \right. \right\rangle. \]
The link between the two presentations is $s = \tau_a^{-1}\tau_b^{-1}\tau_a^{-1}$, $t = \tau_a$.
\smallskip\\
\indent Recall that if we have two simple closed curves $\gamma, x$ then $\tau_{\gamma}(x)$ is obtained as follows: at each intersection point between $x$ and $\gamma$, resolve the intersection by plugging a copy of $\gamma$ into $x$, in such a way that $x$ turns left into the copy of $\gamma$:
\begin{equation}\label{figureDehnTwist}
\begingroup%
  \makeatletter%
  \providecommand\color[2][]{%
    \errmessage{(Inkscape) Color is used for the text in Inkscape, but the package 'color.sty' is not loaded}%
    \renewcommand\color[2][]{}%
  }%
  \providecommand\transparent[1]{%
    \errmessage{(Inkscape) Transparency is used (non-zero) for the text in Inkscape, but the package 'transparent.sty' is not loaded}%
    \renewcommand\transparent[1]{}%
  }%
  \providecommand\rotatebox[2]{#2}%
  \newcommand*\fsize{\dimexpr\f@size pt\relax}%
  \newcommand*\lineheight[1]{\fontsize{\fsize}{#1\fsize}\selectfont}%
  \ifx\svgwidth\undefined%
    \setlength{\unitlength}{244.2039792bp}%
    \ifx\svgscale\undefined%
      \relax%
    \else%
      \setlength{\unitlength}{\unitlength * \real{\svgscale}}%
    \fi%
  \else%
    \setlength{\unitlength}{\svgwidth}%
  \fi%
  \global\let\svgwidth\undefined%
  \global\let\svgscale\undefined%
  \makeatother%
  \begin{picture}(1,0.25652361)%
    \lineheight{1}%
    \setlength\tabcolsep{0pt}%
    \put(0,0){\includegraphics[width=\unitlength,page=1]{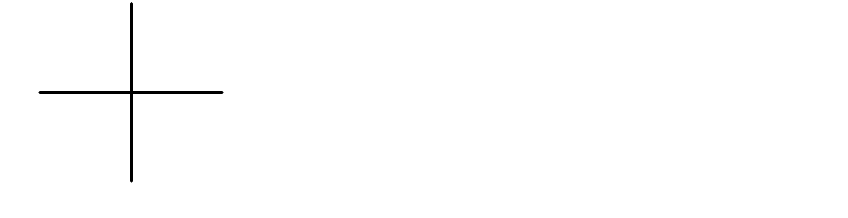}}%
    \put(0.14589021,0.0040971){\color[rgb]{0,0,0}\makebox(0,0)[lt]{\lineheight{1.25}\smash{\begin{tabular}[t]{l}$x$\end{tabular}}}}%
    \put(-0.0023757,0.13938583){\color[rgb]{0,0,0}\makebox(0,0)[lt]{\lineheight{1.25}\smash{\begin{tabular}[t]{l}$\gamma$\end{tabular}}}}%
    \put(0,0){\includegraphics[width=\unitlength,page=2]{dehnTwist.pdf}}%
    \put(0.79666278,0.14063171){\color[rgb]{0,0,0}\makebox(0,0)[lt]{\lineheight{1.25}\smash{\begin{tabular}[t]{l}$\tau_{\gamma}(x)$\end{tabular}}}}%
    \put(0,0){\includegraphics[width=\unitlength,page=3]{dehnTwist.pdf}}%
  \end{picture}%
\endgroup%

\end{equation}
Since $\mathrm{MCG}(\Sigma_{1,0}^{\mathrm{o}})$ fixes the boundary, it fixes the basepoint (see Figure \ref{Sigma10}) and hence we have an action of $\mathrm{MCG}(\Sigma_{1,0}^{\mathrm{o}})$ on $\pi_1(\Sigma_{1,0}^{\mathrm{o}})$. The actions of the Dehn twists $\tau_a$ and $\tau_b$ are given by:
\begin{equation}\label{actionMCG}
\tau_a(a) = a, \:\:\: \tau_a(b) = ba \:\:\: \text{ and } \:\:\: \tau_b(a) = b^{-1}a, \:\:\: \tau_b(b) = b.
\end{equation}
For instance, the action of $\tau_a$ on $b$ is depicted by:
\begin{center}
\begingroup%
  \makeatletter%
  \providecommand\color[2][]{%
    \errmessage{(Inkscape) Color is used for the text in Inkscape, but the package 'color.sty' is not loaded}%
    \renewcommand\color[2][]{}%
  }%
  \providecommand\transparent[1]{%
    \errmessage{(Inkscape) Transparency is used (non-zero) for the text in Inkscape, but the package 'transparent.sty' is not loaded}%
    \renewcommand\transparent[1]{}%
  }%
  \providecommand\rotatebox[2]{#2}%
  \newcommand*\fsize{\dimexpr\f@size pt\relax}%
  \newcommand*\lineheight[1]{\fontsize{\fsize}{#1\fsize}\selectfont}%
  \ifx\svgwidth\undefined%
    \setlength{\unitlength}{436.43750923bp}%
    \ifx\svgscale\undefined%
      \relax%
    \else%
      \setlength{\unitlength}{\unitlength * \real{\svgscale}}%
    \fi%
  \else%
    \setlength{\unitlength}{\svgwidth}%
  \fi%
  \global\let\svgwidth\undefined%
  \global\let\svgscale\undefined%
  \makeatother%
  \begin{picture}(1,0.2207242)%
    \lineheight{1}%
    \setlength\tabcolsep{0pt}%
    \put(0.06939571,0.17946065){\color[rgb]{0,0,0}\makebox(0,0)[lt]{\lineheight{1.25}\smash{\begin{tabular}[t]{l}$b$\end{tabular}}}}%
    \put(0.33639752,0.17192598){\color[rgb]{0,0,0}\makebox(0,0)[lt]{\lineheight{1.25}\smash{\begin{tabular}[t]{l}$[a]$\end{tabular}}}}%
    \put(0,0){\includegraphics[width=\unitlength,page=1]{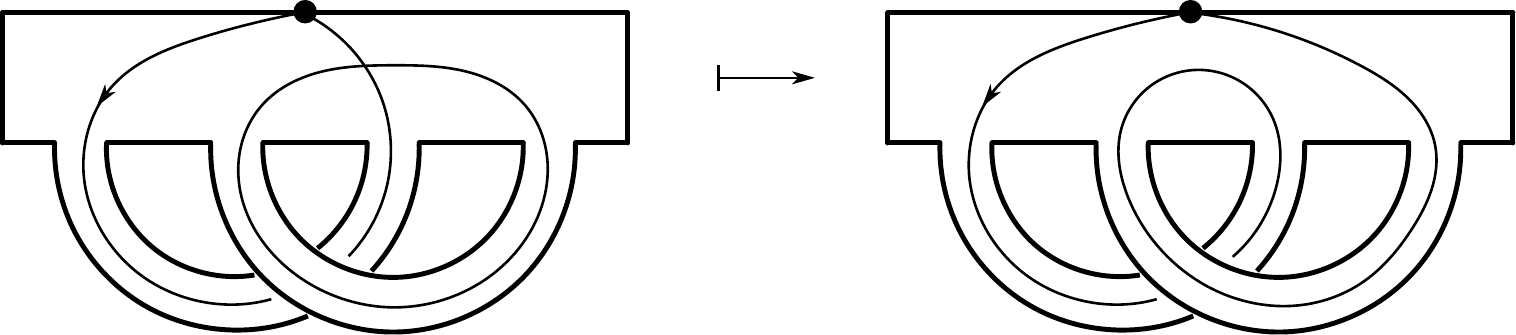}}%
    \put(0.49180015,0.17744841){\color[rgb]{0,0,0}\makebox(0,0)[lt]{\lineheight{1.25}\smash{\begin{tabular}[t]{l}$\tau_a$\end{tabular}}}}%
    \put(0.64657771,0.18199624){\color[rgb]{0,0,0}\makebox(0,0)[lt]{\lineheight{1.25}\smash{\begin{tabular}[t]{l}$ba$\end{tabular}}}}%
    \put(0,0){\includegraphics[width=\unitlength,page=2]{twistBA.pdf}}%
  \end{picture}%
\endgroup%

\end{center}

\subsection{Automorphisms $\widetilde{\tau}_a$ and $\widetilde{\tau}_b$}

\indent The fundamental idea, proposed in \cite{AS} and \cite{AS2}, is to lift the action of the Dehn twists of $\mathrm{MCG}(\Sigma_{g,n} \setminus D)$ on $\pi_1(\Sigma_{g,n}\setminus D)$ at the level of the algebra $\mathcal{L}_{g,n}(H)$. Let us be more precise. In $\pi_1(\Sigma_{1,0}^{\mathrm{o}})$ we have the two canonical curves $a$ and $b$, while in $\mathcal{L}_{1,0}(H)$ we have the matrices $\overset{I}{A}$ and $\overset{I}{B}$. Using (\ref{actionMCG}), let us try to define two morphisms $f_a, f_b : \mathcal{L}_{1,0}(H) \to \mathcal{L}_{1,0}(H)$ by the same formulas:
$$\begin{array}{ll}
f_a(\overset{I}{A}) = \overset{I}{A}, & f_a(\overset{I}{B}) = \overset{I}{B}\overset{I}{A} \\
f_b(\overset{I}{A}) = \overset{I}{B}\,\!^{-1}\overset{I}{A}, & f_b(\overset{I}{B}) = \overset{I}{B}. \\
\end{array}$$
Let us see the behavior of these mappings under the fusion and exchange relations. For the exchange relation, no problem arises:
$$\begin{array}{lll}
R_{12}\,f_a(B)_1\,R_{21}\,f_a(A)_2 &= R_{12}\,B_1\,A_1\,R_{21} \, A_2 & \text{(definition)}\\
&= R_{12}\,B_1\,R_{21}\,A_2\,R_{12}\,A_1\,R_{12}^{-1} & \text{(eq. (\ref{reflection}))}\\
&= A_2 \, R_{12} \, B_1 \, A_1 \, R_{12}^{-1} & \text{(eq. (\ref{echangeL10}))}\\
&=f_a(A)_2 \, R_{12} \, f_a(B)_1 \, R_{12}^{-1} & \text{(definition)}\\
\end{array}$$
and a similar computation holds for $f_b$. The fusion relation is almost satisfied:
\begin{equation}\label{normaFusion1}
\begin{array}{lll}
f_a(B)_{12}&= B_{12}\,A_{12} & \text{(definition)}\\
&= \Delta(v)_{12} \, B_{12} \, v_1^{-1}\,v_2^{-1}\,R_{21}\,R_{12}\,A_{12} & \text{(trick)}\\
&= \Delta(v)_{12} \, B_1 \, R_{21} B_2 \, R_{21}^{-1} \, v_1^{-1} \, v_2^{-1} \, R_{21} \, R_{12} \, A_1 \, R_{21} \, A_2 \, R_{21}^{-1} & \text{(eq. (\ref{fusionL01}))}\\
&= \Delta(v)_{12} \, v_1^{-1} \, v_2^{-1} \, B_1 \, R_{21} B_2 \, R_{12} \, A_1 \, R_{21} \, A_2 \, R_{21}^{-1} & \text{($v$ is central)}\\
&= \Delta(v)_{12} \, v_1^{-1} \, v_2^{-1} \, B_1\,A_1 \, R_{21} \, B_2 \, A_2 \, R_{21}^{-1} & \text{(eq. (\ref{echangeL10}))}\\
&= \Delta(v)_{12} \, v_1^{-1} \, v_2^{-1} \, f_a(B)_1 \, R_{21} \, f_a(B)_2 \, R_{21}^{-1} & \text{(definition)}\\
\end{array}
\end{equation}
and we get similarly:
\begin{equation}\label{normaFusion2}
f_b(A)_{12} = B^{-1}_{12}\, A_{12} = \Delta(v^{-1})_{12} \, v_1 \, v_2 \, f_b(A)_1 \, R_{21} \, f_b(A)_2 \, R_{21}^{-1}.
\end{equation}
From this we conclude that the elements $\overset{I}{v}\,\!^{-1}\overset{I}{B}\overset{I}{A}$ and $\overset{I}{v}\overset{I}{B}\,\!^{-1}\overset{I}{A}$ satisfy the relation (\ref{fusionL01}). Since $v$ is central, we see that the exchange relation still holds with these elements. We thus have found the morphisms which lift $\tau_a$ and $\tau_b$; we denote them by $\widetilde{\tau}_a$ and $\widetilde{\tau}_b$ respectively (these morphisms appeared first in \cite[Lem. 2]{AS} and \cite[eqs (4.1), (4.2)]{AS2}).

\begin{proposition}\label{autoAlphaBeta}
We have two automorphisms $\widetilde{\tau}_a, \widetilde{\tau}_b$ of $\mathcal{L}_{1,0}(H)$ defined by:
$$\begin{array}{ll}
\widetilde{\tau}_a(\overset{I}{A}) = \overset{I}{A}, & \widetilde{\tau}_a(\overset{I}{B}) = \overset{I}{v}\,\!^{-1}\overset{I}{B}\overset{I}{A} \\
\widetilde{\tau}_b(\overset{I}{A}) = \overset{I}{v} \overset{I}{B}\,\!^{-1}\overset{I}{A}, & \widetilde{\tau}_b(\overset{I}{B}) = \overset{I}{B}. \\
\end{array}$$

Moreover, these automorphisms are inner: there exist $\widehat{\tau}_a, \widehat{\tau}_b \in \mathcal{L}_{1,0}(H)$ unique up to scalar such that
$$ \forall \, x \in \mathcal{L}_{1,0}(H), \:\:\: \widetilde{\tau}_a(x) = \widehat{\tau}_a x \widehat{\tau}_a^{-1}, \:\:\: \widetilde{\tau}_b(x) = \widehat{\tau}_b x \widehat{\tau}_b^{-1}. $$
\end{proposition}
\begin{proof}
By Proposition \ref{isoL10Heisenberg}, $\mathcal{L}_{1,0}(H)$ is a matrix algebra. Hence, by the Skolem-Noether theorem, every automorphism of $\mathcal{L}_{1,0}(H)$ is inner.
\end{proof}

A natural question is then to find explicitly the elements $\widehat{\tau}_a, \widehat{\tau}_b$. The answer is amazingly simple (it has been given in \cite[eq (9.7)]{AS} for the modular case; there they express these elements as linear combinations of traces which form a basis in the modular case only). Recall the notation \eqref{notationImbed}.

\begin{proposition}\label{valueHatAlphaBeta}
Up to scalar, $ \widehat{\tau}_a = v_A^{-1} \in \mathcal{L}_{1,0}^{\mathrm{inv}}(H)$ and $\widehat{\tau}_b = v_B^{-1} \in \mathcal{L}_{1,0}^{\mathrm{inv}}(H)$.
\end{proposition}
\begin{proof}
We must show that:
\[
v^{-1}_A\overset{I}{A} = \overset{I}{A} v^{-1}_A, \:\:\: v^{-1}_A \overset{I}{B} = \overset{I}{v}\,\!^{-1} \overset{I}{B}\overset{I}{A} v^{-1}_A \:\:\: \text{ and } \:\:\: v^{-1}_B \overset{I}{A} = \overset{I}{v} \overset{I}{B}\,\!^{-1}\overset{I}{A}v^{-1}_B, \:\:\:  v^{-1}_B \overset{I}{B} = \overset{I}{B} v^{-1}_B.
\]
It is obvious that $v_A^{-1}$ (resp. $v_B^{-1}$) commutes with the matrices $\overset{I}{A}$ (resp. $\overset{I}{B}$) since it is central in $j_A(\mathcal{L}_{0,1}(H))$ (resp. in $j_B(\mathcal{L}_{0,1}(H))$). Let us show the other commutation relation for $v_A^{-1}$. We use the isomorphism $\Psi_{1,0}$. Observe that $\Psi_{1,0}(x_A) = x$ for all $x \in H$. Hence, using the exchange relation of Definition \ref{relDefHeisenberg} and (\ref{ribbon}), we have:
\begin{align*}
\Psi_{1,0}(v_A^{-1}\overset{I}{B}) & = \overset{I}{L}\,\!^{(+)} v^{-1} \overset{I}{T} \overset{I}{L}\,\!^{(-)-1} = \overset{I}{L}\,\!^{(+)} \overset{I}{T}(? v'^{-1}) v''^{-1} \overset{I}{L}\,\!^{(-)-1} = \overset{I}{L}\,\!^{(+)} \overset{I}{T} \, \overset{I}{(v')}{^{-1}} v''^{-1} \overset{I}{L}\,\!^{(-)-1}\\
& = \overset{I}{L}\,\!^{(+)} \overset{I}{T} \, \overset{I}{v}{^{-1}} \overset{I}{b_i} \overset{I}{a_j} \, v^{-1} a_i b_j \overset{I}{L}\,\!^{(-)-1} = \overset{I}{v}{^{-1}} \overset{I}{L}\,\!^{(+)} \overset{I}{T} \, \overset{I}{b_i}a_i \, \overset{I}{a_j}b_j \overset{I}{L}\,\!^{(-)-1} v^{-1}\\
& = \overset{I}{v}{^{-1}} \overset{I}{L}\,\!^{(+)} \overset{I}{T} \overset{I}{L}\,\!^{(-)-1} \, \overset{I}{L}\,\!^{(+)} \overset{I}{L}\,\!^{(-)-1} \, v^{-1} = \Psi_{1,0}(\overset{I}{v}\,\!^{-1} \overset{I}{B}\overset{I}{A} v^{-1}_A)
\end{align*}
as desired. We now apply the morphism $\widetilde{\tau}_a^{-1} \circ \widetilde{\tau}_b^{-1}$ to the equality $v_A^{-1}\overset{I}{B} = \overset{I}{v}\,\!^{-1} \overset{I}{B}\overset{I}{A} v_A^{-1}$:
$$ \widetilde{\tau}_a^{-1} \circ \widetilde{\tau}_b^{-1}(v_A^{-1}\overset{I}{B}) = \overset{I}{v} v_B^{-1} \overset{I}{B}\overset{I}{A}{^{-1}} = \widetilde{\tau}_a^{-1} \circ \widetilde{\tau}_b^{-1}(\overset{I}{v}^{-1}\overset{I}{B}\overset{I}{A}v_A^{-1}) = \overset{I}{B} \overset{I}{A}{^{-1}} \overset{I}{B} v_B^{-1}. $$
Using that $v_B$ and $\overset{I}{B}$ commute, we easily get the desired equality.
\end{proof}


\subsection{Projective representation of $\mathrm{SL}_2(\mathbb{Z})$ on $\SLF(H)$}\label{sectionProjRepSL2Z}
\indent Observe that 
\[ \tau_a\tau_b\tau_a = \tau_b\tau_a\tau_b, \:\:\:\:\:\:\: \left(\tau_a\tau_b\right)^6 \neq \mathrm{id} \:\:\:\:\:\:\:\:\: \text{ in } \pi_1(\Sigma_{1,0}^{\mathrm{o}}). \]
Since $\tau_a(c) = \tau_b(c) = c$ where $c = ba^{-1}b^{-1}a$ is the boundary loop induced by the deletion of the open disk $D$, $\tau_a$ and $\tau_b$ are well-defined in $\pi_1(\Sigma_{1,0}) = \pi_1(\Sigma_{1,0}^{\mathrm{o}})/\langle c \rangle$, and we have
\[ \tau_a\tau_b\tau_a = \tau_b\tau_a\tau_b, \:\:\:\:\:\:\: \left(\tau_a\tau_b\right)^6 = \mathrm{id} \:\:\:\:\:\:\:\:\: \text{ in } \pi_1(\Sigma_{1,0}). \]
Recall that $\mathcal{L}_{1,0}(H)$ is associated to $\Sigma_{1,0}^{\mathrm{o}}$, and it is easy to check that
\[ \widetilde{\tau}_a \widetilde{\tau}_b \widetilde{\tau}_a = \widetilde{\tau}_b \widetilde{\tau}_a \widetilde{\tau}_b, \:\:\:\:\:\:\: \left(\widetilde{\tau}_a \widetilde{\tau}_b\right)^6 \neq \mathrm{id} \:\:\:\:\:\:\:\:\: \text{ in } \mathcal{L}_{1,0}(H).\]
It follows from Propositions \ref{autoAlphaBeta} and \ref{valueHatAlphaBeta} that
\[ v_A^{-1} v_B^{-1} v_A^{-1} \sim v_B^{-1} v_A^{-1} v_B^{-1}, \:\:\:\:\:\:\: \left( v_A^{-1} v_B^{-1} \right)^6 \not\sim 1 \:\:\:\:\:\:\:\:\: \text{ in } \mathcal{L}_{1,0}(H) \]
where $\sim$ means equality up to scalar (we will see that $\sim$ is actually $=$ for the braid relation). Hence, if we want a representation of $\mathrm{MCG}(\Sigma_{1,0})$ based on the elements $v_A^{-1}$ and $v_B^{-1}$, we have to glue back the disc $D$. Recall that the matrices $\overset{I}{C}$ corresponding to the boundary circle are killed on $\mathrm{SLF}(H)$ (Lemma \ref{CSLF}). Hence, it is natural to think that representing $v_A^{-1}$ and $v_B^{-1}$ on $\mathrm{SLF}(H)$ (see Theorem \ref{repInv}) will provide a projective representation of $\mathrm{MCG}(\Sigma_{1,0})$:
\[ \rho_{\mathrm{SLF}}\!\left( v_A^{-1}\,v_B^{-1}\,v_A^{-1} \right) \sim \rho_{\mathrm{SLF}}\!\left( v_B^{-1}\,v_A^{-1}\,v_B^{-1} \right), \:\:\:\:\:\:\: \rho_{\mathrm{SLF}}\!\left( v_A^{-1}\,v_B^{-1} \right)^6 \sim 1. \]
We will show that this indeed holds.

\smallskip


\indent Recall from Proposition \ref{propVIntegrale} the symmetric linear forms 
\[ \varphi_v = \mu^l(v^{-1})^{-1}\mu^l\bigl( g^{-1}v^{-1}\,? \bigr), \:\:\:\:\:\:\: \varphi_{v^{-1}} = \mu^l(v)^{-1}\mu^l\bigl(g^{-1}v\,?\bigr). \]
satisfying $\mathcal{D}(\varphi_{v^{\pm 1}}) = v^{\pm 1}$. Due to the fact that $\varphi_{v^{-1}} = \varphi_v^{-1}$ (since $\mathcal{D}$ is an isomorphism of algebras), we see that
\begin{equation}\label{phiVcarre}
\varphi_{v^{-1}}\varphi_{v^{-1}}^{v^{-2}} = \frac{\mu^l(v^{-1})}{\mu^l(v)}\varepsilon.
\end{equation}
where $\beta^h = \beta(h?)$ for all $\beta \in H^*, h \in H$. By Proposition \ref{actionAB}, the actions of $v_A^{-1}$ and $v_B^{-1}$ on $\SLF(H)$ are:
\begin{equation}\label{actionsV}
\forall \, \psi \in \SLF(H), \:\:\: v_A^{-1} \triangleright \psi = \psi^{v^{-1}} = \psi(v^{-1}?) \:\:\: \text{ and } \:\:\: v_B^{-1} \triangleright \psi = \left(\varphi_{v^{-1}}\psi^v\right)^{v^{-1}}.
\end{equation}

\begin{lemma}\label{phiPhiV}
$\varphi_{v^{-1}} \varphi_{v^{-1}}^{v^{-1}} = \varphi_{v^{-1}}^{v^{-1}}$.
\end{lemma}
\begin{proof}
For $x \in H$:
\begin{align*}
\left\langle \varphi_{v^{-1}} \varphi_{v^{-1}}^{v^{-1}},\, x \right\rangle &= \mu^l\!\left(v\right)^{-2} \mu^l\!\left(v g^{-1}x'\right)\mu^l\!\left(g^{-1}x''\right) = \mu^l\!\left(v\right)^{-2} \left\langle \mu^l\!\left(v ?\right)\mu^l, \, g^{-1}x \right\rangle\\
&= \mu^l\!\left(v\right)^{-1} \mu^l(g^{-1}x) = \varphi_{v^{-1}}^{v^{-1}}(x).
\end{align*}
We simply used (\ref{integrale2}).
\end{proof}

\indent This lemma has an important consequence.

\begin{proposition}\label{braidV}
The following braid relation holds in $\mathcal{L}_{1,0}(H)$:
$$ v_A^{-1} \, v_B^{-1} \, v_A^{-1} = v_B^{-1} \, v_A^{-1} \, v_B^{-1}. $$
\end{proposition}
\begin{proof}
The morphisms $\widetilde{\tau}_a$ and $\widetilde{\tau}_b$ satisfy the braid relation $\widetilde{\tau}_a\widetilde{\tau}_b\widetilde{\tau}_a = \widetilde{\tau}_b\widetilde{\tau}_a\widetilde{\tau}_b$. Hence by Proposition \ref{valueHatAlphaBeta} and since $\mathcal{Z}(\mathcal{L}_{1,0}(H)) \cong \mathbb{C}$, we have: $\lambda v_A^{-1} \, v_B^{-1} \, v_A^{-1} = v_B^{-1} \, v_A^{-1} \, v_B^{-1}$ for some $\lambda \in \mathbb{C}$. We evaluate on the counit:
\begin{align*}
\lambda v_A^{-1} \, v_B^{-1} \, v_A^{-1} \triangleright \varepsilon &= \lambda v_A^{-1} \, v_B^{-1} \triangleright \varepsilon = \lambda v_A^{-1} \triangleright \varphi_{v^{-1}}^{v^{-1}} = \lambda \varphi_{v^{-1}}^{v^{-2}}\\
v_B^{-1} \, v_A^{-1} \, v_B^{-1} \triangleright \varepsilon &= v_B^{-1} \, v_A^{-1} \triangleright \varphi_{v^{-1}}^{v^{-1}} = v_B \triangleright \varphi_{v^{-1}}^{v^{-2}} = \left(\varphi_{v^{-1}}\varphi_{v^{-1}}^{v^{-1}}\right)^{v^{-1}} = \left(\varphi_{v^{-1}}^{v^{-1}}\right)^{v^{-1}} = \varphi_{v^{-1}}^{v^{-2}}.
\end{align*}
We used $\varepsilon(v?) = \varepsilon(v)\varepsilon = \varepsilon$ and Lemma \ref{phiPhiV}. It follows that $\lambda=1$.
\end{proof}

\indent Consider $\widehat{\omega} = v_A^{-1} \, v_B^{-1} \, v_A^{-1} \in \mathcal{L}_{1,0}(H)$, which implements the automorphism $\omega = \widetilde{\tau}_a\widetilde{\tau}_b\widetilde{\tau}_a$: $\omega(x) = \widehat{\omega} \, x \, \widehat{\omega}^{-1}$.
The key observation is the following lemma.

\begin{lemma}\label{keylemmaOmega}
For all $\psi \in \SLF(H)$:
$$ \widehat{\omega}^2 \triangleright \psi = \frac{\mu^l(v^{-1})}{\mu^l(v)}S(\psi). $$
\end{lemma}
\begin{proof}
First, we show the formula for $\psi = \varepsilon$:
$$ \widehat{\omega}^2 \triangleright \varepsilon = (v_A^{-1} \, v_B^{-1} \, v_A^{-1})^2 \triangleright \varepsilon = v_A^{-1} \, v_B^{-1} \, v_A^{-1} \triangleright \varphi_{v^{-1}}^{v^{-2}} = v_A^{-1} \triangleright  \bigl(\varphi_{v^{-1}}\varphi_{v^{-1}}^{v^{-2}}\bigr)^{v^{-1}} = \frac{\mu^l(v^{-1})}{\mu^l(v)} v_A^{-1} \triangleright \varepsilon^{v^{-1}} = \frac{\mu^l(v^{-1})}{\mu^l(v)} \varepsilon $$
where we used (\ref{actionsV}) and (\ref{phiVcarre}). Second, note that $\omega(\overset{I}{A}) = \overset{I}{v}{^2} \overset{I}{A}{^{-1}}\overset{I}{B}{^{-1}}\overset{I}{A} = \overset{I}{B}{^{-1}}\overset{I}{C}$ and $\omega(\overset{I}{B}) = \overset{I}{A}$. It follows that, for every $z \in \mathcal{Z}(H) = \mathcal{L}^{\mathrm{inv}}_{0,1}(H)$, $\omega^2(z_B) = \omega(z_A) = z_{B^{-1}C}$ and thus
\[ \hat \omega^2 z_B = z_{B^{-1}C} \hat \omega^2. \]
Observe by Proposition \ref{actionAB} that for every $\psi \in \mathrm{SLF}(H)$ we have $\psi = \mathcal{D}\!\left(\psi^v\right)_B \triangleright \varepsilon$. Hence, we get
\begin{align*}
\widehat{\omega}^2 \triangleright \psi &= \widehat{\omega}^2 \mathcal{D}\!\left(\psi^v\right)_B \triangleright \varepsilon = \mathcal{D}\!\left(\psi^v\right)_{B^{-1}C} \widehat{\omega}^2 \triangleright \varepsilon = \frac{\mu^l(v^{-1})}{\mu^l(v)}  \mathcal{D}\!\left(\psi^v\right)_{B^{-1}C} \triangleright \varepsilon = \frac{\mu^l(v^{-1})}{\mu^l(v)}  \mathcal{D}\!\left(\psi^v\right)_{B^{-1}} \triangleright \varepsilon\\
& = \frac{\mu^l(v^{-1})}{\mu^l(v)} S\!\left(\psi^v\right)^{v^{-1}} = \frac{\mu^l(v^{-1})}{\mu^l(v)}S(\psi).
\end{align*}
We simply used Lemmas \ref{CSLF} and \ref{lemmaBMoinsUn}. Also recall that if $\varphi \in \mathrm{SLF}(H)$ then 
\[ z_{B^{-1}C} \triangleright \varphi = \sum_I \mathrm{tr}\!\left(\Lambda_I \overset{I}{g} \overset{I}{B}{^{-1}} \overset{I}{C} \triangleright \varphi \right) = \sum_I \mathrm{tr}\!\left(\Lambda_I \overset{I}{g} \overset{I}{B}{^{-1}} \triangleright \varphi \right) = z_{B^{-1}} \triangleright \varphi, \]
with the notation of Lemma \ref{writingInvariants}.
\end{proof}

\indent Recall that the group $\mathrm{PSL}_2(\mathbb{Z}) = \mathrm{SL}_2(\mathbb{Z})/\{\pm \mathbb{I}_2\}$ admits the following presentations:
$$ \mathrm{PSL}_2(\mathbb{Z}) = \left\langle \tau_a, \tau_b \, \left\vert \, \tau_a\tau_b\tau_a = \tau_b\tau_a\tau_b, \: (\tau_a\tau_b)^3=1 \right. \right\rangle = \left\langle s, t \, \left\vert \, (st)^3=1, \: s^2=1 \right. \right\rangle. $$
We denote by $\rho$ the representation of $\mathcal{L}_{0,1}(H)$ on $H^*$ and by $\rho_{\mathrm{SLF}}$ the representation of $\mathcal{L}_{0,1}^{\mathrm{inv}}(H)$ on $\mathrm{SLF}(H)$.

\begin{theorem}\label{repSL2}
Recall that we assume that $H$ is a finite dimensional factorizable ribbon Hopf algebra.\\
1) The assignment
$$ \tau_a \mapsto \rho\!\left(v_A^{-1}\right), \:\:\:\:\: \tau_b \mapsto \rho\!\left(v_B^{-1}\right) $$
defines a representation $\theta_1^{\mathrm{o}}$ of $\mathrm{MCG}(\Sigma_{1,0}^{\mathrm{o}}) = B_3$ on $H^*$.
\\2) The assignment
$$ \tau_a \mapsto \rho_{\mathrm{SLF}}\!\left(v_A^{-1}\right), \:\:\:\:\: \tau_b \mapsto \rho_{\mathrm{SLF}}\!\left(v_B^{-1}\right) $$
defines a projective representation $\theta_1$ of $\mathrm{MCG}(\Sigma_{1,0}) = \mathrm{SL}_2(\mathbb{Z})$ on $\SLF(H)$. If moreover $S(\psi) = \psi$ for all $\psi \in \SLF(H)$, then this defines actually a projective representation of $\mathrm{PSL}_2(\mathbb{Z})$.
\end{theorem}
\begin{proof}
By Proposition \ref{braidV}, we know that the braid relation is satisfied in $\mathcal{L}_{1,0}(H)$, thus the first claim holds. By Lemma \ref{keylemmaOmega}, we have:
$$ (v_A^{-1} v_B^{-1})^3 \triangleright \psi = \widehat{\omega}^2 \triangleright \psi = \frac{\mu^l(v^{-1})}{\mu^l(v)}S(\psi). $$
If $S_{\vert \mathrm{SLF}(H)} = \mathrm{id}$, then
$$ \rho_{\mathrm{SLF}}\!\left(v_A^{-1} \,v_B^{-1} \right)^3 = \frac{\mu^l(v^{-1})}{\mu^l(v)}\mathrm{id}. $$
Otherwise, 
$$ (v_A^{-1} v_B^{-1})^6 \triangleright \psi = \frac{\mu^l(v^{-1})}{\mu^l(v)}\widehat{\omega}^2 \triangleright S(\psi) = \frac{\mu^l(v^{-1})^2}{\mu^l(v)^2}S^{2}(\psi) = \frac{\mu^l(v^{-1})^2}{\mu^l(v)^2}\psi(g \, ? \, g^{-1}) = \frac{\mu^l(v^{-1})^2}{\mu^l(v)^2}\psi. $$
\end{proof}

\noindent Observe that the quantity $\frac{\mu^l(v^{-1})}{\mu^l(v)}$ does not depend on the choice of $\mu^l$ since it is unique up to scalar.

%

\subsection{Equivalence with the Lyubashenko-Majid representation}\label{equivalenceAvecLM}
\indent Recall that $H$ is a finite dimensional factorizable ribbon Hopf algebra. Under this assumption, two operators $\mathcal{S}, \mathcal{T} : H \to H$ are defined in \cite{LM}:
$$ \mathcal{S}(x) = \left(\mathrm{id} \otimes \mu^l\right)\!\left(R^{-1}(1 \otimes x)R'^{-1}\right), \:\:\:\:\:\:\: \mathcal{T}(x) = v^{-1}x. $$
It is shown that they are invertible and satisfy $(\mathcal{S}\mathcal{T})^3 = \lambda\mathcal{S}^2$, $\mathcal{S}^2 = S^{-1}$, with $\lambda \in \mathbb{C}\!\setminus\!\{0\}$. 
We warn the reader that in \cite{LM}, they consider the {\em inverse} of the ribbon element (see the bottom of the third page of their paper). That is why there is $v^{-1}$ in the formula for $\mathcal{T}$.
\smallskip\\
\indent Now we introduce two maps. The first is
$$\fonc{\boldsymbol{\chi}}{H^*}{H}{\beta}{(\beta \otimes \mathrm{id})(R'R)}$$
while the second is
$$\fonc{\boldsymbol{\gamma}}{H}{H^*}{x}{\mu^r(S(x)\,?).}$$
The map $\boldsymbol{\chi}$ is a slight variant of the map $\Psi$ of section \ref{rappelHopf} and is called Drinfeld morphism in \cite{FGST}. The map $\boldsymbol{\gamma}$ is denoted $\widehat{\boldsymbol{\phi}}{^{-1}}$ in \cite{FGST} and is the inverse of the Radford map ($\widehat{\boldsymbol{\phi}}(\varphi) = \varphi(c')c''$, where $c$ is the two-sided cointegral of $H$; see \cite{radford, radfordBook}). Consider the space of left $q$-characters:
$$ \mathrm{Ch}^l(H) = \left\{ \beta \in H^* \left\vert \, \forall\, x,y \in H, \:\:\beta(xy) = \beta\!\left(S^2(y)x\right) \right. \right\}. $$
These maps satisfy the following restrictions:
$$ \boldsymbol{\chi} \: : \: \mathrm{Ch}^l(H) \longrightarrow \mathcal{Z}(H), \:\:\:\:\: \boldsymbol{\gamma} \: : \: \mathcal{Z}(H) \longrightarrow \mathrm{Ch}^l(H). $$
This is due to the fact (observed by Drinfeld \cite{drinfeld} and Radford \cite{radford} respectively) that they intertwine the adjoint and the coadjoint actions (for the first the computation is analogous to that of the proof of Proposition \ref{reshetikhin}, while the second is immediate by Proposition \ref{propVIntegrale}).
\smallskip\\
\indent It is not too difficult to show (see e.g. \cite[Remark IV.1.2]{ibanez}) that
$$ \forall \, z \in \mathcal{Z}(H), \:\: \mathcal{S}(z) = \boldsymbol{\chi} \circ \boldsymbol{\gamma}(z). $$
It follows that $\mathcal{Z}(H)$ is stable under $\mathcal{S}$ and $\mathcal{T}$. But since $S^2$ is inner, we have $\mathcal{S}^4(z) = S^{-2}(z) = z$ for each $z \in \mathcal{Z}(H)$. Thus there exists a projective representation $\rho_{\mathrm{LM}}$ of $\mathrm{SL}_2(\mathbb{Z})$ on $\mathcal{Z}(H)$, defined by:
$$ \rho_{\mathrm{LM}}(s) = \mathcal{S}_{\vert \mathcal{Z}(H)}, \:\:\:\:\:\:\: \rho_{\mathrm{LM}}(t) = \mathcal{T}_{\vert \mathcal{Z}(H)}. $$

\smallskip

\indent The left $q$-characters are nothing more than shifted symmetric linear forms. More precisely, we have an isomorphism of algebras:
$$ \fonc{\left(g^{-1}\right)^*}{\SLF(H)}{\mathrm{Ch}^l(H)}{\psi}{\psi(g^{-1}?).} $$
Let us define shifted versions of $\boldsymbol{\chi}$ and of $\boldsymbol{\gamma}$:
$$ \boldsymbol{\chi}_{g^{-1}} = \boldsymbol{\chi} \circ \left(g^{-1}\right)^* \: : \: \SLF(H) \overset{\cong}{\longrightarrow} \mathcal{Z}(H), \:\:\:\:\: \boldsymbol{\gamma}_g = g^* \circ \boldsymbol{\gamma} \: : \: \mathcal{Z}(H) \overset{\cong}{\longrightarrow} \SLF(H). $$
The equality $\mathcal{S} = \boldsymbol{\chi}_{g^{-1}} \circ \boldsymbol{\gamma}_g$ still holds, but we have now $\SLF(H)$ instead of $\mathrm{Ch}^l(H)$.

\smallskip

\indent We will need the following relation between left and right integrals.
\begin{lemma}\label{lemmaUnimodUnibalance}
Under our assumptions $H$ is unibalanced, which means that $\mu^l = \mu^r(g^2\,?)$.
\end{lemma}
\begin{proof}
The terminology ``unibalanced'' is picked from \cite{BBG}, where some facts about integrals and cointegrals are recalled. Recall (see \textit{e.g.} \cite[Prop. 8.10.10]{EGNO}) that a finite dimensional factorizable Hopf algebra is unimodular, which means that there exists $c \in \mathcal{Z}(H)$, called two-sided cointegral, such that $xc = \varepsilon(x)c$ for all $x \in H$. Let $a \in H$ be the comodulus of $\mu^r$: $\psi\mu^r = \psi(a)\mu^r$ for all $\psi \in \mathcal{O}(H)$ (see \textit{e.g.} \cite[eq. 4.9]{BBG}). By a result of Drinfeld (see \cite[Prop. 10.1.14]{Mon}, but be aware that in this book the notations and conventions for $a$ and $g$ are different from those we use), we know that: 
$$ uS(u)^{-1} = a(\mathfrak{a} \otimes \mathrm{id})(R) $$
where $\mathfrak{a} \in H^*$ is the modulus of the left cointegral $c^l$ of $H$. Here, since $c = c^l$ is two-sided, we have $\mathfrak{a} = \varepsilon$. Thus $g^2 = u^2v^{-2} = uS(u)^{-1} = a$ by (\ref{ribbon}) and (\ref{pivotCan}). We deduce that 
$$\mu^l = \mu^r \circ S = \mu^r(a?) = \mu^r(g^2?)$$
where the second equality is \cite[Prop. 4.7]{BBG}.
\end{proof}

\begin{lemma}\label{lemmaPourEquivalence}
It holds:
$$ \rho_{\mathrm{SLF}}\!\left(v_A^2\,v_B\right) = \mu^l(v^{-1})^{-1} \, \boldsymbol{\gamma}_g \circ \boldsymbol{\chi}_{g^{-1}}. $$
\end{lemma}
\begin{proof}
We compute each side of the equality. On the one hand:
$$
\boldsymbol{\gamma}_g \circ \boldsymbol{\chi}_{g^{-1}}(\psi) = \boldsymbol{\gamma}_g\!\left((\psi \otimes \mathrm{id})\!\left(g^{-1}(v^{-1})'v \otimes (v^{-1})''v\right)\right) = \psi\!\left(g^{-1}(v^{-1})'v\right) \mu^r\!\left(gS\!\left((v^{-1})''\right)v\,?\right)
$$
whereas on the other hand:
\begin{align*}
v_A^2\,v_B \triangleright \psi &= (\varphi_v\psi^v)^v = \mu^l(v^{-1})^{-1} \left( \mu^r\!\left(gv^{-1}\,?\right)\psi^v \right)^v\\
&= \mu^l(v^{-1})^{-1}\left[\mu^r\!\left(g(v^{-1})'\,?\right)\psi\!\left( vg^{-1}S^{-1}\!\left((v^{-1})''\right) \right)\right]^v\\
&= \mu^l(v^{-1})^{-1}\mu^r\!\left(gvS\!\left((v^{-1})''\right)\,?\right)\psi\!\left( vg^{-1}(v^{-1})' \right)
\end{align*}
as desired. We used the formulas (which are analogous to \eqref{actionsV})
\[ v_A \triangleright \psi = \psi^{v} = \psi(v?), \:\:\:\:\:\:\: v_B \triangleright \psi = \mu^l(v^{-1})^{-1} \bigl( \mu^l\bigl( g^{-1}v^{-1}\,? \bigr) \psi^v\bigr)^{v^{-1}}\]
together with Lemma \ref{lemmaUnimodUnibalance}, the property \eqref{integraleShifte} and the equality $(v^{-1})' \otimes S^{-1}\!\left((v^{-1})''\right) = S\!\left((v^{-1})''\right) \otimes (v^{-1})'$ which is due to $S(v)=v$.
\end{proof}

\indent The link between the two presentations of $\mathrm{SL}_2(\mathbb{Z})$ is $s = \tau_a^{-1}\tau_b^{-1}\tau_a^{-1}$, $t = \tau_a$. Hence we define two operators $\mathcal{S}', \mathcal{T}' : \SLF(H) \to \SLF(H)$ by:
$$ \mathcal{S}' = \theta_1(s) = \rho_{\mathrm{SLF}}(v_Av_Bv_A), \:\:\:\:\: \mathcal{T}' = \theta_1(t) = \rho_{\mathrm{SLF}}(v_A^{-1}). $$

\begin{theorem}\label{EquivalenceLMandSLF}
Recall that we assume that $H$ is a finite dimensional factorizable ribbon Hopf algebra. 
Then the projective representation $\theta_1$ of Theorem \ref{repSL2} is equivalent to $\rho_{\mathrm{LM}}$.
\end{theorem}
\begin{proof}
Consider the following isomorphism of vector spaces:
$$ \fonc{f = \rho_{\mathrm{SLF}}(v_A^{-1}) \circ \boldsymbol{\gamma}_g}{\mathcal{Z}(H)}{\SLF(H)}{z}{\boldsymbol{\gamma}_g(z)^{v^{-1}} = \mu^r\!\left(gv^{-1}S(z)\,?\right).} $$
By Lemma \ref{lemmaPourEquivalence}:
$$  \mathcal{S}' = \mu^l(v^{-1})^{-1}\rho_{\mathrm{SLF}}(v_A^{-1}) \circ \boldsymbol{\gamma}_g \circ \boldsymbol{\chi}_{g^{-1}} \circ \rho_{\mathrm{SLF}}(v_A). $$
Thus:
$$ f \circ \mathcal{S} = \rho_{\mathrm{SLF}}(v_A^{-1}) \circ \boldsymbol{\gamma}_g \circ \boldsymbol{\chi}_{g^{-1}} \circ \boldsymbol{\gamma}_g = \mu^l(v^{-1})\mathcal{S}' \circ f. $$
Next,
$$ f \circ \mathcal{T}(z) = f(v^{-1}z) = \boldsymbol{\gamma}_g(z)^{v^{-2}} = \rho_{\mathrm{SLF}}(v_A^{-1})\!\left(\boldsymbol{\gamma}_g(z)^{v^{-1}}\right) = \mathcal{T}' \circ f(z) $$
Then $f$ is an intertwiner of projective representations.
\end{proof}

\section{The case of $H = \bar U_q(\mathfrak{sl}_2)$}\label{exempleUq}

\indent Let $q$ be a primitive root of unity of order $2p$, with $p \geq 2$. We now work in some detail the case of $H = \bar U_q = \bar U_q(\mathfrak{sl}_2)$, the restricted quantum group associated to $\mathfrak{sl}_2(\mathbb{C})$. We take back all the notations and results from Chapter \ref{chapitreUqSl2}. In particular, to explicitly describe the representation of $\mathrm{SL}_2(\mathbb{Z})$, we will use the GTA basis of $\SLF(\bar U_q)$ introduced in section \ref{sectionSLF}.

\subsection{Technical details}\label{technicalDetails}
\indent In principle, since $\bar U_q$ is not braided (see section \ref{braidedExtension}), it is not clear that the previous definitions and results remain valid. In practice, the universal $R$-matrix simply belongs to the extension $\bar U_q^{1/2}$ of $\bar U_q$ by a square root of $K$, and although some computations occur in the extension, the final result always belongs to $\bar U_q$. The important things are that the $M$-matrix $RR'$ belongs to $\bar U_q^{\otimes 2}$, that the ribbon element $v$ belongs to $\bar U_q$ and that $\bar U_q$ is factorizable (recall that this last claim is an abuse of terminology since $\bar U_q$ is not braided, see details in section \ref{braidedExtension}).

\smallskip

\indent In order to define $\mathcal{L}_{0,1}(\bar U_q)$ and $\mathcal{L}_{1,0}(\bar U_q)$ we introduce some terminology. Let $I$ be a $\bar U_q^{1/2}$-module. Since $\bar U_q \subset \bar U_q^{1/2}$, $I$ determines a $\bar U_q$-module, which we denote $I_{\vert \bar U_q}$. We say that a $\bar U_q$-module $J$ is {\em liftable} if there exists a $\bar U_q^{1/2}$-module $\widetilde{J}$ such that $\widetilde{J}_{\vert \bar U_q} = J$. Not every $\bar U_q$-module is liftable. Indeed, if it was the case, this would imply that $V \otimes W \cong W \otimes V$ (since $\bar U_q^{1/2}$ is braided) for every $\bar U_q$-modules $V,W$, which is false: a counter-example is given in \cite{KS}. However, the simple modules and the PIMs are liftable, which is enough for us. It suffices to define the action of $K^{1/2}$ on these modules. Take back the notations of section \ref{sectionSimpleProj}  for the canonical basis of modules. For the simple module $\mathcal{X}^{\epsilon}(s)$ ($\epsilon \in \{\pm\})$, there are two choices for $\epsilon^{1/2}$, and so the two possible lifts are defined by 
$$ K^{1/2}v_j = \epsilon^{1/2}q^{(s-1-2j)/2}v_j $$
and the action of $E$ and $F$ is unchanged. Similarly, the two possible lifts of the PIM $\mathcal{P}^{\epsilon}(s)$ are defined by
$$
\begin{array}{ll}
K^{1/2}b_0 = \epsilon^{1/2}q^{(s-1)/2}b_0, & K^{1/2}x_0 = \left(\epsilon^{1/2}q^{p/2}\right)q^{(p-s-1)/2}x_0 \\
K^{1/2}y_0 = \left(-\epsilon^{1/2}q^{p/2}\right)q^{(p-s-1)/2}y_0, & K^{1/2}a_0 = \epsilon^{1/2}q^{(s-1)/2}a_0
\end{array}
$$
and the action of $E$ and $F$ is unchanged. 
\smallskip\\
\indent Let $\widetilde{\mathbb{C}}^-$ be the $1$-dimensional $\overline{U}\!_q{^{\!1/2}}$-module with basis $v$ defined by $Ev = Fv = 0$, $K^{1/2}v=-v$ (which is a lift of $\mathcal{X}^+(1) = \mathbb{C}$). If $\widetilde{I}$ is a lift of a simple module or a PIM $I$, then we have seen that the only possible lift of $I$ are $\widetilde{I}^+ = \widetilde{I}$ and $\widetilde{I}^- = \widetilde{I} \otimes \widetilde{\mathbb{C}}^-$. Moreover, using \eqref{deltaR}, we get equalities which will be used in the next section:
\begin{equation}\label{LiftR}
\begin{array}{ll}
\bigl(\overset{\widetilde{I}^- \widetilde{J}}{R}\,\bigr)_{\! 12} = \bigl(\overset{\widetilde{I}^+ \widetilde{J}}{R}\bigr)_{\! 12}\, \bigl(\overset{\widetilde{J}}{K^p}\bigr)_2, & \bigl(\overset{\widetilde{I} \, \widetilde{J}^-}{\!R\:}\bigr)_{12} = \bigl(\overset{\widetilde{I} \, \widetilde{J}^+}{\!R\:}\bigr)_{\! 12} \, \bigl(\overset{\widetilde{I}}{K^p}\bigr)_1\\
\bigl(\overset{\widetilde{I}^-\, \widetilde{J}}{R'}\bigr)_{\! 12} = \bigl(\overset{\widetilde{I}^+ \, \widetilde{J}}{R'}\bigr)_{\! 12}\, \bigl(\overset{\widetilde{J}}{K^p}\bigr)_2, & \bigl(\overset{\widetilde{I}\, \widetilde{J}^-}{R'}\bigr)_{\! 12} = \bigl(\overset{\widetilde{I}\widetilde{J}^+}{R'}\bigr)_{\! 12}\, \bigl(\overset{\widetilde{I}}{K^p}\bigr)_1.
\end{array}
\end{equation}

\subsection{$\mathcal{L}_{0,1}(\bar U_q)$ and $\mathcal{L}_{1,0}(\bar U_q)$}\label{sectionL01Uq}
\indent We define $\mathcal{L}_{0,1}(\bar U_q)$ as the quotient of $\mathrm{T}(\bar U_q^*)$ by the fusion relation
$$ \overset{I \otimes J}{M}\!_{12} = \overset{I}{M}_1\overset{\widetilde{I}\widetilde{J}}{(R')}_{12}\overset{J}{M}_2\overset{\widetilde{I}\widetilde{J}}{(R'^{-1})}_{12} $$
where $I, J$ are simple modules or PIMs and $\widetilde{I}, \widetilde{J}$ are lifts of $I$ and $J$. From (\ref{LiftR}) and the fact that $K^p$ is central, we see that this does not depend on the choice of $\widetilde{I}$ and $\widetilde{J}$. As we saw in section \ref{matrices}, the matrix coefficients of the PIMs linearly span $\mathcal{L}_{0,1}(H)$, thus we can restrict to them in the definition. However, it is important to have the matrices associated to simple modules and more precisely to $\mathcal{X}^+(2)$, as we shall see below. All the results of section \ref{sectionLoopL01} remain true for $\mathcal{L}_{0,1}(\bar U_q)$. In particular, $\Psi_{0,1}$ is an isomorphism since $\bar U_q$ is factorizable.
\smallskip\\
\indent We now describe $\mathcal{L}_{0,1}(\bar U_q)$ by generators and relations. Let
$$ M = \overset{\mathcal{X}^+(2)}{M} = 
\begin{pmatrix}
a&b\\
c&d
\end{pmatrix}
\:\:\:\:\:\: \text{ and } \:\:\:\:\:\: \widetilde{R} = \overset{\widetilde{\mathcal{X}}^+(2)\widetilde{\mathcal{X}}^+(2)}{R} = q^{-1/2}
\begin{pmatrix}
q & 0 & 0 & 0\\
0 & 1 & \hat q & 0\\
0 & 0 & 1 & 0\\
0 & 0 & 0 & q
\end{pmatrix}
$$
where $\widetilde{\mathcal{X}}^+(2)$ is the lift of $\mathcal{X}^+(2)$ defined by $K^{1/2}v_0 = q^{1/2}v_0$. Using the same reasoning as in section \ref{sectionMatrixCoeffUq}, we know that the coefficients of $\overset{\mathcal{X}^+(2)^{\otimes n}}{M}$ ($n \in \mathbb{N}$) linearly span $\mathcal{L}_{0,1}(\bar U_q)$. But thanks to the fusion relation, each such coefficient is a linear combination of products of the elements $a,b,c,d$:
\[ \bigl(\overset{\mathcal{X}^+(2)^{\otimes n}}{M}\bigr)^{i_1, \ldots, i_n}_{j_1, \ldots, j_n} = \lambda({^{i_1, \ldots, i_n}_{j_1, \ldots, j_n}})_{k_1, \ldots, k_n}^{l_1, \ldots, l_n} \bigl(\overset{\mathcal{X}^+(2)}{M}\bigr){^{k_1}_{l_1}} \ldots \bigl(\overset{\mathcal{X}^+(2)}{M}\bigr){^{k_n}_{l_n}}, \]
which is the analogue of \eqref{TX2tensN}. It follows that $a,b,c,d$ generate $\mathcal{L}_{0,1}(\bar U_q)$. Let us determine relations between these elements. First, we have the reflection equation, which comes from the existence of the braiding morphism $c : \mathcal{X}^+(2)^{\otimes 2} \to \mathcal{X}^+(2)^{\otimes 2}$: 
\[ \widetilde{R}_{12} M_1 \widetilde{R}_{21} M_2 = M_2 \widetilde{R}_{12} M_1 \widetilde{R}_{21}. \]

\noindent This equation is equivalent to the following exchange relations:
$$
\begin{array}{ccc}
da=ad, & db=q^2bd, & dc=q^{-2}cd, \\
ba = ab+q^{-1}\hat qbd, & cb = bc+q^{-1}\hat q(da-d^2), & ca = ac-q^{-1}\hat q dc.
\end{array}
$$
with $\hat q = q - q^{-1}$. Second, since $\mathcal{X}^+(2)^{\otimes 2} \cong \mathcal{X}^+(1) \oplus \mathcal{X}^+(3)$\footnote{This decomposition does not hold if $p=2$: in that case, we have $\mathcal{X}^+(2)^{\otimes 2} \cong \mathcal{P}^+(1)$. But there is still the morphism $\Psi : \mathbb{C} \to \mathcal{X}^+(2)^{\otimes 2}$ which corresponds to sending $\mathbb{C} = \mathcal{X}^+(1)$ in $\mathrm{Soc}\bigl(\mathcal{P}^+(1)\bigr)$.}, there exists a unique (up to scalar) morphism $\Phi : \mathbb{C} = \mathcal{X}^+(1) \to \mathcal{X}^+(2)^{\otimes 2}$. It is easily computed:
$$ \Phi(1) = q v_0 \otimes v_1 - v_1 \otimes v_0. $$
By fusion and naturality \eqref{naturaliteL01}, we have
$$ M_1\widetilde{R}_{21} M_2 \widetilde{R}_{21}^{-1}\Phi = \overset{\mathcal{X}^+(2)^{\otimes 2}}{M}\!\!\!\!\!\!_{12}\Phi = \Phi\overset{\mathbb{C}}{M} = \Phi. $$
This gives just one new relation, which is the analogue of the quantum determinant \eqref{quantumDet}:
$$ ad-q^2bc=1. $$ 
Finally, let us compute the RSD isomorphism on $M$:
$$ \Psi_{0,1}
\begin{pmatrix}
a & b\\
c & d
\end{pmatrix}
= \overset{\widetilde{\mathcal{X}}^+(2)}{L^{(+)}}\overset{\widetilde{\mathcal{X}}^+(2)}{L^{(-)-1}} = 
\begin{pmatrix}
K^{1/2} & \hat q K^{1/2}F\\
0 & K^{-1/2}
\end{pmatrix}
\begin{pmatrix}
K^{1/2} & 0\\
\hat q K^{-1/2}E & K^{-1/2}
\end{pmatrix}
=
\begin{pmatrix}
K+q^{-1}\hat q^2FE & q^{-1}\hat q F\\
\hat q K^{-1}E & K^{-1}
\end{pmatrix}.
$$
We deduce the relations $b^p = c^p = 0$ and $d^{2p}=1$ from the defining relations of $\bar U_q$. 

\begin{theorem}\label{theoPresL01Uq}
The algebra $\mathcal{L}_{0,1}(\bar U_q)$ admits the following presentation:
$$\left\langle 
a,b,c,d\,
\left\vert
\begin{array}{ccc}
da=ad, & db=q^2bd, & dc=q^{-2}cd \\
ba = ab+q^{-1}\hat qbd, & cb = bc+q^{-1}\hat q(da-d^2), & ca = ac-q^{-1}\hat q dc\\
ad-q^2bc = 1, & b^p=c^p=0, & d^{2p}=1
\end{array}
\right.
\right\rangle.$$
A basis is given by the monomials $b^ic^jd^k$ with $0 \leq i,j \leq p-1$, $0 \leq k \leq 2p-1$.
\end{theorem}
\begin{proof}
Let $A$ be the algebra defined by this presentation. It is readily seen that $a = d^{-1} + q^2bcd^{-1}$ and that the monomials $b^ic^jd^k$ with $0 \leq i,j \leq p-1$, $0 \leq k \leq 2p-1$ linearly span $A$. Thus $\dim(A) \leq 2p^3$. But we know that $2p^3 = \dim(\bar U_q) = \dim\!\left(\mathcal{L}_{0,1}(\bar U_q)\right)$ since the monomials $E^iF^jK^{\ell}$ with $0 \leq i,j \leq p-1$, $0 \leq k \leq 2p-1$ form the PBW basis of $\bar U_q$. It follows that $\dim(A) \leq \dim\!\left(\mathcal{L}_{0,1}(\bar U_q)\right)$. Since these relations are satisfied in $\mathcal{L}_{0,1}(\bar U_q)$, there exists a surjection $p : A \to \mathcal{L}_{0,1}(\bar U_q)$. Thus $\dim(A) \geq \dim\!\left(\mathcal{L}_{0,1}(\bar U_q)\right)$, and the theorem is proved.
\end{proof}

\begin{remark}
It is possible to get the relations $b^p=c^p=0, d^{2p}=1$ by fusion and naturality, as it was done in section \ref{sectionMatrixCoeffUq}. The big difference is that the fusion relation of $\mathcal{L}_{0,1}(\bar U_q)$ is more complicated than the fusion relation of $\mathcal{O}(\bar U_q)$ and consequently the matrices $\overset{\mathcal{X}^+(2)^{\otimes n}}{M}$ are much more complicated than the matrices $\overset{\mathcal{X}^+(2)^{\otimes n}}{T}$. First one must show by induction \textit{via} tedious matrix reasonings that
\begin{center}
 $\overset{\mathcal{X}^+(2)^{\otimes n}}{M} = 
\begin{blockarray}{cccc}
 v_0^{\otimes p} & \ldots & v_1^{\otimes p}  \\
  \begin{block}{(ccc)c}
    \ast & \ast & \lambda_n b^n & v_0^{\otimes p}  \\
   \ast  & \ast & \ast & \vdots\\
  \mu_n c^n & \ast & d^n & v_1^{\otimes p}\\
  \end{block}
\end{blockarray}$ 
\end{center}
where $\lambda_n, \mu_n$ are non-zero scalars. Then one can take back the reasoning of section \ref{sectionMatrixCoeffUq}: the morphism $f : \mathcal{P}^+(p-1) \to \mathcal{X}^+(2)^{\otimes p}$ will imply $b^p=c^p=0$ and the morphism $f_2$ will imply $d^{2p}=1$. In contrast, the relation $a^{2p}=1$ is not true in $\mathcal{L}_{0,1}(\bar U_q)$.
\finEx
\end{remark}

\begin{remark}Theorem \ref{theoPresL01Uq} indicates that $\mathcal{L}_{0,1}(\bar U_q)$ is a restricted version (\textit{i.e.} a finite dimensional quotient by monomial central elements) of $\mathcal{L}_{0,1}(U_q)^{\mathrm{spe}}$, the specialization at our root of unity $q$ of the algebra $\mathcal{L}_{0,1}(U_q)$. A complete study of the algebra $\mathcal{L}_{0,1}(U_q)^{\mathrm{spe}}$ will appear in \cite{BaR}. 
\finEx
\end{remark}

\indent Applying the isomorphism of algebras $\mathcal{D}$ defined in (\ref{morphismeDrinfeld}) to the GTA basis of $\SLF(\bar U_q)$ (defined in section \ref{sectionSLF}), we get a basis of $\mathcal{Z}(\bar U_q) = \mathcal{L}_{0,1}^{\mathrm{inv}}(\bar U_q)$. We introduce notations for these basis elements\footnote{The elements $\boldsymbol{\chi}(s)$ defined in \cite{GT} correspond to $[s]V^s$ here}:
\begin{equation}\label{SLFcentrauxUq}
\overset{\mathcal{X}^{\epsilon}(s)}{W} = \mathcal{D}\!\left( \chi^{\epsilon}_s \right), \:\:\:\:\:\: V^{s'} = \mathcal{D}\!\left( G_{s'} \right)
\end{equation}
with $1 \leq s \leq p$, $\epsilon \in \{\pm\}$ and $1 \leq s' \leq p-1$. They satisfy the same multiplication rules than the elements of the GTA basis, see Theorem \ref{ProduitArike}. Unwinding the definitions, this reads 
\begin{equation}\label{formuleWV}
\begin{split}
\overset{\mathcal{X}^{\epsilon}(s)}{W} &= (\chi^{\epsilon}_s \otimes \mathrm{id})\bigl( (K^{p+1} \otimes 1) RR' \bigr) \cong \mathrm{tr}\bigl(\overset{\mathcal{X}^{\epsilon}(s)}{K^{p+1}}\overset{\mathcal{X}^{\epsilon}(s)}{M}\bigr) \\
V^{s} &= (G_s \otimes \mathrm{id})\bigl( (K^{p+1} \otimes 1) RR' \bigr) \cong \mathrm{tr}\bigl(\sigma_s \, \overset{\mathcal{P}^+(s)}{K^{p+1}}\overset{\mathcal{P}^+(s)}{M}\bigr) + \mathrm{tr}\bigl(\sigma_{p-s} \overset{\mathcal{P}^-(p-s)}{K^{p+1}}\overset{\mathcal{P}^-(p-s)}{M}\bigr)
\end{split}
\end{equation}
since we choose $K^{p+1}$ as pivotal element and where $\cong$ is the identification $\mathcal{L}_{0,1}(\bar U_q) = \bar U_q$ \textit{via} $\Psi_{0,1}$ (recall that usually we use this identification without mention). In particular, the expression \eqref{RRPrime} of $RR'$ allows us to compute that
\begin{equation}\label{WCasimir}
 \overset{\mathcal{X}^+(2)}{W} = -qa - q^{-1}d = \mathcal{D}(\chi^+_2) = -\hat q^2FE - qK - q^{-1}K^{-1} = -\hat q^2 C
\end{equation}
where $C$ is the Casimir element \eqref{defCasimir}.

\medskip

\indent Similarly, we define $\mathcal{L}_{1,0}(\bar U_q)$ as the quotient of $\mathcal{L}_{0,1}(\bar U_q) \ast \mathcal{L}_{0,1}(\bar U_q)$ by the exchange relations:
$$ \overset{\widetilde{I}\widetilde{J}}{R}_{12}\overset{I}{B}_1\overset{\widetilde{I}\widetilde{J}}{(R')}_{12}\overset{J}{A}_2 = \overset{J}{A}_2\overset{\widetilde{I}\widetilde{J}}{R}_{12}\overset{I}{B}_1\overset{\widetilde{I}\widetilde{J}}{(R^{-1})}_{12} $$
where $I, J$ are simple modules or PIMs and $\widetilde{I}, \widetilde{J}$ are liftings of $I$ and $J$. From (\ref{LiftR}), we see again that this does not depend on the choice of $\widetilde{I}$ and $\widetilde{J}$. The coefficients of $\overset{\mathcal{X}^+(2)}{A}$ and of $\overset{\mathcal{X}^+(2)}{B}$:
\[
\overset{\mathcal{X}^+(2)}{A} = 
\begin{pmatrix}
a_1 & b_1\\
c_1 & d_1
\end{pmatrix},\:\:\:\:
\overset{\mathcal{X}^+(2)}{B} = 
\begin{pmatrix}
a_2 & b_2\\
c_2 & d_2
\end{pmatrix}
\]
generate $\mathcal{L}_{1,0}(\bar U_q)$. Using the commutation relations of the Heisenberg double, it is easy to show that $\Psi_{1,0}$ indeed takes values in $\mathcal{H}(\mathcal{O}(\bar U_q))$ (the square root of $K$ does not appear). In order to obtain a presentation of $\mathcal{L}_{1,0}(\bar U_q)$, one can again restrict to $I = J = \mathcal{X}^+(2)$ and write down the corresponding exchange relations together with the relations coming from the fact that the variables $a_i, b_i, c_i, d_i$ ($i=1,2$) generate a copy of $\mathcal{L}_{0,1}(\bar U_q)$. We do not give this presentation of $\mathcal{L}_{1,0}(\bar U_q)$ since it is quite cumbersome and we will not use it in this work. Let us just mention that the monomials
\begin{equation}\label{PBWL10}
b_1^i c_1^j d_1^k b_2^l c_2^m d_2^n, \:\:\:\:\: 0 \leq i,j,l,m \leq p-1, \:\: 0 \leq k, n \leq 2p-1
\end{equation}
form a basis. Indeed, they are a generating set thanks to the exchange relations and to the restriction relations $b_1^p = c_1^p = b_2^p = c_2^p = 0, d_1^{2p} = d_2^{2p} = 1$; moreover, the number of such monomials is $4p^6 = \dim\bigl(\mathcal{L}_{1,0}(\bar U_q)\bigr)$. The generators $a_1, a_2$ do not appear in the monomial basis because they can be expressed as $a_1 = d_1^{-1} + q^2b_1c_1d^{-1}_1$, $a_2 = d_2^{-1} + q^2b_2c_2d^{-1}_2$.

\smallskip

\indent In view of the next section, let us precise that by definition (see \eqref{notationImbed} and \eqref{formuleWV}) we have
\begin{equation}\label{traceWV}
\overset{\mathcal{X}^{\epsilon}(s)}{W}\!\!\!\!_X = \mathrm{tr}\bigl(\overset{\mathcal{X}^{\epsilon}(s)}{K^{p+1}}\overset{\mathcal{X}^{\epsilon}(s)}{X}\bigr), \:\:\:\:\:\:\:\:\:\: V^{s}_X = \mathrm{tr}\bigl(\sigma_s \, \overset{\mathcal{P}^+(s)}{K^{p+1}}\overset{\mathcal{P}^+(s)}{X}\bigr) + \mathrm{tr}\bigl(\sigma_{p-s} \overset{\mathcal{P}^-(p-s)}{K^{p+1}}\overset{\mathcal{P}^-(p-s)}{X}\bigr),
\end{equation}
where $X$ is any product of the matrices $B, A$ with some normalization by $v^r$, for instance $X= A$, $X = B$, $X = v B^{-1} A$.

\subsection{Explicit description of the $\mathrm{SL}_2(\mathbb{Z})$-projective representation}\label{SL2ZUq}
\indent Note that it can be shown directly that $\bar U_q$ is unimodular and unibalanced, see for instance \cite[Cor. II.2.8]{ibanez} (also note that in \cite{BBG} it is shown that all the simply laced restricted quantum groups at roots of unity are unibalanced). 

\begin{proposition}\label{antipodeSLF}
For all $z \in \mathcal{Z}(\bar U_q)$, $S(z) =z$ and for all $\psi \in \SLF(\bar U_q)$, $S(\psi) = \psi$. It follows that in the case of $\bar U_q$, $\rho_{\SLF}$ is in fact a projective representation of $\mathrm{PSL}_2(\mathbb{Z})$.
\end{proposition}
\begin{proof}
By \cite[Appendix D]{FGST}, the canonical central elements are expressed as $e_s = P_s(C)$, $w^{\pm}_s = \pi^{\pm}_sQ_s(C)$ where $P_s$ and $Q_s$ are polynomials, $C$ is the Casimir element \eqref{defCasimir} and $\pi_s^{\pm}$ are discrete Fourier transforms of $(K^j)_{0 \leq j \leq 2p-1}$. It is easy to check that $S(C) = C$ and that $S(\pi^{\pm}_s) = \pi^{\pm}_s$, thus $S(e_s) = e_s$ and $S(w^{\pm}_s) = w^{\pm}_s$. Next, let $\psi \in \SLF(\bar U_q)$. Since $\boldsymbol{\gamma}_g$ is an isomorphism, we can write $\psi = \boldsymbol{\gamma}_g(z) = \mu^r(gS(z)\,?)$ with $z \in \mathcal{Z}(\bar U_q)$. Then:
$$ S(\psi) = S\!\left(\mu^r(gz\,?)\right) = \mu^r \circ S(?zg^{-1}) = \mu^l(g^{-1}z?) = \mu^r(gz?) = \psi. $$
We used that $S(z)=z$, Proposition \ref{propVIntegrale} and the fact that $\bar U_q$ is unibalanced.
\end{proof}

\indent We want to determine the action of $\theta_1(\tau_a) = \rho_{\mathrm{SLF}}(v_A^{-1})$ and $\theta_1(\tau_b) = \rho_{\mathrm{SLF}}(v_B^{-1})$ on the GTA basis. For this, we will need the expression \eqref{rubanCentre} of $v^{\pm 1}$ in the canonical basis of $\mathcal{Z}(\bar U_q)$, the formulas \eqref{actionCentreSLF} for the action \eqref{defActionCentreSLF} of $\mathcal{Z}(\bar U_q)$ on $\SLF(\bar U_q)$, the multiplication rules in the GTA basis (Theorem \ref{ProduitArike}) and the following lemma (in which we use the notation \eqref{notationImbed}).
\begin{lemma}\label{actionBmoinsUnA}
Let $z \in \mathcal{L}_{0,1}^{\mathrm{inv}}(H) = \mathcal{Z}(H)$ and let $\psi \in \SLF(H)$. Then:
$$ z_{vB^{-1}A} \triangleright \psi = S\!\left(\mathcal{D}^{-1}(z)\right) \psi. $$
\end{lemma}
\begin{proof}
The proof is analogous to those of the two similar results in section \ref{sectionRepInvariants} and is thus left to the reader. Note that this lemma is not specific to the case of $\bar U_q$.
\end{proof}

\begin{theorem}\label{actionSL2ZArike}
Let $\theta_1 : \mathrm{SL}_2(\mathbb{Z}) \to \mathrm{PGL}(\bar U_q^*)$ be the projective representation obtained in Theorem \ref{repSL2}, with gauge algebra $\bar U_q = \bar U_q(\mathfrak{sl}_2)$. The representations of the Dehn twists $\tau_a$ and $\tau_b$ on the GTA basis are given by:
$$ \theta_1(\tau_a)(\chi^{\epsilon}_s) = v^{-1}_{\mathcal{X}^{\epsilon}(s)}\chi^{\epsilon}_s, \:\:\:\:\:\:\: \theta_1(\tau_a)(G_{s'}) = v^{-1}_{\mathcal{X}^+(s')}G_{s'} -  v_{\mathcal{X}^+(s')}^{-1}\hat q\left( \frac{p-s'}{[s']}\chi^+_{s'} - \frac{s'}{[s']}\chi^-_{p-s'} \right) $$
and
\begin{align*}
\theta_1(\tau_b)(\chi^{\epsilon}_s) & = \xi \epsilon(-\epsilon)^{p-1}s q^{-(s^2-1)} \left(\sum_{\ell=1}^{p-1}(-1)^s(-\epsilon)^{p-\ell}\left(q^{\ell s} + q^{-\ell s}\right)\left(\chi^+_{\ell} + \chi^-_{p-\ell}\right) +  \chi^+_p + (-\epsilon)^p(-1)^s\chi^-_p \right)\\  
&\:\:\:+\xi\epsilon (-1)^sq^{-(s^2-1)}\sum_{j=1}^{p-1} (-\epsilon)^{j+1}[j][js]G_j,\\
\theta_1(\tau_b)(G_{s'}) & = \xi (-1)^{s'}q^{-(s'^2-1)}\frac{\hat q p}{[s']}\sum_{j=1}^{p-1}(-1)^{j+1}[j][js']\left(2G_j - \hat q \frac{p-j}{[j]}\chi^+_j + \hat q \frac{j}{[j]}\chi^-_{p-j}\right),
\end{align*}
with $\epsilon \in \{\pm\}$, $0 \leq s \leq p$, $1 \leq s' \leq p-1$ and $\xi^{-1} = \frac{1-i}{2\sqrt{p}} \frac{\hat q^{p-1}}{[p-1]!} (-1)^p q^{-(p-3)/2}$.
\end{theorem}
\begin{proof}
The formulas for $\theta_1(\tau_a) = \rho_{\mathrm{SLF}}(v_A^{-1})$ are easily deduced from Proposition \ref{actionAB}, (\ref{rubanCentre}) and (\ref{actionCentreSLF}). Computing the action of $\theta_1(\tau_b) = \rho_{\mathrm{SLF}}(v_B^{-1})$ is more difficult. 
\noindent We will use the commutation relations of $v_B^{-1}$ with the $A,B$-matrices, namely
\begin{equation}\label{commBeta}
v_B^{-1}\overset{I}{A} = \overset{I}{v} \overset{I}{B}{^{-1}}\overset{I}{A}v_B^{-1}, \:\:\:\:\: v_B^{-1}\overset{I}{B} = \overset{I}{B}v_B^{-1}
\end{equation}
to compute the action of $v_B^{-1}$ by induction. The multiplication rules of the GTA basis (Theorem \ref{ProduitArike}) will be used several times. Let us denote
$$ v_B^{-1} \triangleright \chi^{\epsilon}_s = \sum_{\sigma \in \{\pm\}}\sum_{\ell=1}^p \lambda^{\sigma}_{\ell}(\epsilon, s) \chi^{\sigma}_{\ell} + \sum_{j=1}^{p-1}\delta_j(\epsilon, s)G_j. $$
Taking the quantum trace of relation (\ref{commBeta}) provides $ v_B^{-1}\,\overset{\mathcal{X}^+(2)}{W}\!\!\!\!\!_A = \overset{\mathcal{X}^+(2)}{W}\!\!\!\!_{vB{^{-1}}A}\,v_B^{-1}$ (recall \eqref{traceWV}). 
On the one hand, we obtain by (\ref{WCasimir}), \eqref{casimirBaseCan} and \eqref{actionCentreSLF}:
$$ v_B^{-1}\,\overset{\mathcal{X}^+(2)}{W}\!\!\!\!\!_A \triangleright \chi^{\epsilon}_s = v_B^{-1} \triangleright  \chi^{\epsilon}_s(-\hat q^2C\, ?) = \sum_{\substack{\ell=1\\ \sigma \in \{\pm\}}}^p -\epsilon (q^s+q^{-s}) \lambda^{\sigma}_{\ell}(\epsilon, s) \chi^{\sigma}_{\ell} + \sum_{j=1}^{p-1}-\epsilon(q^s + q^{-s})\delta_j(\epsilon, s)G_j. $$
On the other hand, we use Lemma \ref{actionBmoinsUnA} together with \eqref{SLFcentrauxUq} and the multiplication rules:
\begin{align*}
& \overset{\mathcal{X}^+(2)}{W}\!\!\!\!_{vB{^{-1}}A}v_B^{-1} \triangleright \chi^{\epsilon}_s = \sum_{\sigma \in \{\pm\}}\sum_{\ell=1}^p \lambda^{\sigma}_{\ell}(\epsilon, s) \chi^+_2\chi^{\sigma}_{\ell} + \sum_{j=1}^{p-1}\delta_j(\epsilon, s)\chi^+_2G_j\\
&= \sum_{\sigma \in \{\pm\}} \left(\lambda^{\sigma}_2(\epsilon,s) + 2\lambda^{-\sigma}_p(\epsilon,s)\right)\!\chi^{\sigma}_1 + \sum_{\ell=2}^{p-2} \left(\lambda^{\sigma}_{s-1}(\epsilon,s) + \lambda^{\sigma}_{s+1}(\epsilon,s)\right)\!\chi^{\sigma}_{\ell} + \left(\lambda^{\sigma}_{p-2}(\epsilon,s) + 2\lambda^{\sigma}_p(\epsilon,s)\right)\!\chi^{\sigma}_{p-1}\\
&\:\:\:\:\:\:\:\:\: + \lambda^{\sigma}_{p-1}(\epsilon,s)\chi^{\sigma}_p + \frac{\delta_2(\epsilon, s)}{[2]}G_1 + \sum_{j=2}^{p-2}[j]\left(\frac{\delta_{j-1}(\epsilon, s)}{[j-1]} + \frac{\delta_{j+1}(\epsilon, s)}{[j+1]}\right)\!G_j + \frac{\delta_{p-2}(\epsilon, s)}{[2]}G_{p-1}.
\end{align*}
This gives recurrence equations between the coefficients which are easily solved:
\begin{align*}
v_B^{-1} \triangleright \chi^{\epsilon}_s &= \lambda(\epsilon, s) \left( \sum_{\ell=1}^{p-1}(-1)^s(-\epsilon)^{p-\ell}\left(q^{\ell s} + q^{-\ell s}\right)\left(\chi^+_{\ell} + \chi^-_{p-\ell}\right) +  \chi^+_p + (-\epsilon)^p(-1)^s\chi^-_p \right)\\  
&\:\:\:+\delta(\epsilon,s)\sum_{j=1}^{p-1} (-\epsilon)^{j+1}\frac{[j][js]}{[s]}G_j.
\end{align*}
The coefficients $\lambda(\epsilon, s) = \lambda^+_p(\epsilon, s)$ and $\delta(\epsilon,s) = \delta_1(\epsilon,s)$ are still unknown. In order to compute them by induction, we use the relation $v_B^{-1}\,\overset{\mathcal{X}^+(2)}{W}\!\!\!\!\!_B = \overset{\mathcal{X}^+(2)}{W}\!\!\!\!\!_B\,v_B^{-1}$, which is another consequence of (\ref{commBeta}) (and of \eqref{traceWV}). Before, note that
$$ \overset{\mathcal{X}^+(2)}{W}\!\!\!\!\!_B \triangleright \chi^{\epsilon}_s = \left(\chi^+_2(\chi^{\epsilon}_s)^v\right)^{v^{-1}} = \frac{v_{\mathcal{X}^{\epsilon}(s)}}{v_{\mathcal{X}^{\epsilon}(s-1)}}\chi^{\epsilon}_{s-1} + \frac{v_{\mathcal{X}^{\epsilon}(s)}}{v_{\mathcal{X}^{\epsilon}(s+1)}}\chi^{\epsilon}_{s+1} = -\epsilon q^{-s+\frac{1}{2}}\chi^{\epsilon}_{s-1} - \epsilon q^{s+\frac{1}{2}}\chi^{\epsilon}_{s+1}. $$
with $1 \leq s \leq p-1$ and the convention that $\chi^{\pm}_0=0$. It follows that
\begin{equation}\label{eqInductionBeta}
v_B^{-1} \triangleright \chi^{\epsilon}_{s+1} = -\epsilon q^{-s-\frac{1}{2}}\overset{\mathcal{X}^+(2)}{W}\!\!\!\!\!_B \triangleright \left(v_B^{-1} \triangleright \chi^{\epsilon}_{s}\right) - q^{-2s}v_B^{-1} \triangleright \chi^{\epsilon}_{s-1}.
\end{equation}
Due to (\ref{rubanCentre}), (\ref{actionCentreSLF}) and the multiplication rules, we have
\begin{align*}
&\overset{\mathcal{X}^+(2)}{W}\!\!\!\!\!_B \triangleright \left(v_B^{-1} \triangleright \chi^{\epsilon}_{s}\right) = \left( \chi^+_2 \left(v_B^{-1} \triangleright \chi^{\epsilon}_{s}\right)^v\right)^{v^{-1}} 
= \frac{v_{\mathcal{X}^+(p-1)}}{v_{\mathcal{X}^+(p)}}\left( \lambda^+_{p-1}(\epsilon,s) + \hat q \delta_{p-1}(\epsilon, s)\right)\!\chi^+_{p} + \frac{v_{\mathcal{X}^+(2)}}{[2]}\delta_2(\epsilon, s)G_1 + \ldots\\
\end{align*}
where the dots ($\ldots$) mean the remaining of the linear combination in the GTA basis. After replacing by the values found previously and inserting in relation (\ref{eqInductionBeta}), this yields
\begin{align*}
&\lambda(\epsilon,s+1)\chi^+_p + \delta(\epsilon, s+1)G_1 + \ldots \\
&= \left( q^{-(s+1)}(q^s+q^{-s})\lambda(\epsilon,s) - q^{-2s}\lambda(\epsilon,s-1) + (-\epsilon)^{p-1}(-1)^{s-1}\hat qq^{-(s+1)}\delta(\epsilon,s) \right)\!\chi^+_p\\
&\:\:\:\:\: + \left( -q^{-(s+2)}(q^s + q^{-s})\delta(\epsilon,s) - q^{-2s}\delta(\epsilon,s-1) \right)\!G_1 +\ldots.
\end{align*}
These are recurrence equations. It just remains to determine the first values $\lambda(\epsilon, 1)$, $\delta(\epsilon,1)$. Observe that, since $\bar U_q$ is unibalanced:
\begin{equation}\label{relXi}
v_B^{-1} \triangleright \chi^+_1 = \left(\varphi_{v^{-1}} (\chi^+_1)^v\right)^{v^{-1}} = \mu^l(v)^{-1}\mu^l\!\left(K^{p-1}v\,?\right)^{v^{-1}} = \mu^l(v)^{-1}\mu^r(K^{p+1}?). 
\end{equation}
In \eqref{traceModifieeGTA}, the decomposition of $\mu^r(K^{p+1}?)$ in the GTA is given (when $\mu^r$ is suitably normalized). Thanks to this, we obtain
$$ v_B^{-1} \triangleright \chi^+_1 = \lambda(+, 1)\chi^+_p + \delta(+,1)G_1 + \ldots = \xi(-1)^{p-1}\chi^+_p - \xi G_1 + \ldots $$
and
\begin{align*}
v_B^{-1} \triangleright \chi^-_1 &= v_{\mathcal{X}^-(1)} v_B^{-1}\overset{\mathcal{X}^-(1)}{W}\!\!\!\!\!_B \triangleright \chi^+_1 = v_{\mathcal{X}^-(1)} \overset{\mathcal{X}^-(1)}{W}\!\!\!\!\!_B\triangleright v_B^{-1} \triangleright \chi^+_1 = v_{\mathcal{X}^-(1)}\left( \chi^-_1 \left(v_B^{-1} \triangleright \chi^+_1\right)^v \right)^{v^{-1}}\\
&= -\xi\chi^+_p + \xi G_1 + \ldots = \lambda(-,1)\chi^+_p + \delta(-,1)G_1 + \ldots.
\end{align*}
The scalar $\xi$ does not depend on the choice of $\mu^r$ thanks to the factor $\mu^l(v)^{-1} = \mu^l \circ S(v)^{-1} = \mu^r(v)^{-1}$ in \eqref{relXi}. Using the formulas \cite{FGST} for $\mu^r$ and $v$ in the PBW basis to compute $\mu^r(v)$ gives the value of $\xi$. We are now in position to solve the recurrence equations. It is easy to check that the solutions are
$$ \delta(\epsilon, s) = \xi\epsilon (-1)^sq^{-(s^2-1)}[s], \:\:\:\:\:\:\: \lambda(\epsilon, s) = \xi\epsilon(-\epsilon)^{p-1}s q^{-(s^2-1)}. $$
We now proceed with the proof of the formula for $G_{s'}$. Relation (\ref{commBeta}) implies $v_B^{-1} V_B^1 = V_B^1 v_B^{-1}$ (see \eqref{traceWV}). By (\ref{rubanCentre}), (\ref{actionCentreSLF}) and the multiplication rules, we have on the one hand:
$$ v_B^{-1} V_B^1 \triangleright \chi^+_s = v_B^{-1} \triangleright \left(G_1(\chi^+_s)^v \right)^{v^{-1}} = [s] v_B^{-1} \triangleright G_s - \hat q(p-s) v_B^{-1} \triangleright \chi^+_s + \hat q s v_B^{-1} \triangleright \chi^-_{p-s} $$
whereas on the other hand:
$$ V_B^1 v_B^{-1} \triangleright \chi^+_s = \left(G_1 (v_B^{-1} \triangleright \chi^+_s)^v\right)^{v^{-1}} = \hat q p \sum_{j=1}^{p-1} \delta_j(+,s) \left(G_j - \hat q \frac{p-j}{[j]} \chi^+_j + \hat q\frac{j}{[j]}\chi^-_{p-j}\right). $$
Equalizing both sides and inserting the previously found values, we obtain the desired formula.
\end{proof}

\begin{remark}
The guiding principle of the previous computations was that the mutiplication of two symmetric linear forms in the GTA basis is easy when one of them is $\chi^+_2$, $\chi^-_1$ or $G_1$ (see Theorem \ref{ProduitArike}), and that all the formulas can be derived from $v_B^{-1} \triangleright \chi^+_1$ using only such products.
\finEx
\end{remark}
\smallskip
\indent Recall that the standard representation $\mathbb{C}^2$ of $\mathrm{SL}_2(\mathbb{Z}) = \mathrm{MCG}(\Sigma_{1,0})$ is defined by
$$ \tau_a \mapsto 
\begin{pmatrix}
1 & 0\\
-1 & 1
\end{pmatrix}, \:\:\:\:\:\:
\tau_b \mapsto 
\begin{pmatrix}
1 & 1\\
0 & 1
\end{pmatrix}.
$$

\begin{lemma}\label{lemmaPropRepSL2Z}
Let $V$ be a (projective) representation of $\mathrm{SL}_2(\mathbb{Z})$ which admits a basis $\left(x_s, y_s\right)$ such that
$$
\begin{array}{ll}
\tau_a\,x_s = \sum_{\ell} a_{\ell}(s)x_{\ell}, & \tau_b\,x_s = \sum_{\ell} b_{\ell}(s)(x_{\ell} + y_{\ell})\\
\tau_a\,y_s = \sum_{\ell} a_{\ell}(s)(y_{\ell} - x_{\ell}), & \tau_b\,y_s = \sum_{\ell} b_{\ell}(s)y_{\ell}.
\end{array}
$$
Then there exists a (projective) representation $W$ of $\mathrm{SL}_2(\mathbb{Z})$ such that $V \cong \mathbb{C}^2 \otimes W$. More precisely, $W$ admits a basis $(w_s)$ such that 
$$\begin{array}{ll}
\tau_a\,w_s = \sum_{\ell} a_{\ell}(s)w_{\ell}, & \tau_b\,v_s = \sum_{\ell} b_{\ell}(s) w_{\ell}. 
\end{array}$$
\end{lemma}
\begin{proof}
It is easy to check that the formulas for $\tau_a \, w_s$ and $\tau_b \, w_s$ indeed define a $\mathrm{SL}_2(\mathbb{Z})$-representation on $W$. Let $(e_1, e_2)$ be the canonical basis of $\mathbb{C}^2$. Then 
$$ e_1 \otimes w_s \mapsto y_s, \:\:\:\:\:\: e_2 \otimes w_s \mapsto x_s $$
is an isomorphism which intertwines the $\mathrm{SL}_2(\mathbb{Z})$-action.
\end{proof}

\begin{theorem}\label{thDecRep}
The $(p+1)$-dimensional subspace $\mathcal{P} = \vect\!\left(\chi^+_s + \chi^-_{p-s}, \chi^+_p, \chi^-_p\right)_{1 \leq s \leq p-1}$ is stable under the $\mathrm{SL}_2(\mathbb{Z})$-action of Theorem \ref{actionSL2ZArike}. Moreover, there exists a $(p-1)$-dimensional projective representation $\mathcal{W}$ of $\mathrm{SL}_2(\mathbb{Z})$ such that
$$ \SLF(\bar U_q) = \mathcal{P} \oplus \left(\mathbb{C}^2 \otimes \mathcal{W}\right). $$
\end{theorem}
\begin{proof}
By Corollary \ref{coroVermaIdeal}, $\mathcal{P}$ is an ideal of $\SLF(\bar U_q)$. It is easy to see that $\mathcal{P}$ is moreover stable under the action (\ref{actionCentreSLF}) of $\mathcal{Z}(\bar U_q)$. Thus we deduce without any computation that $\mathcal{P}$ is $\mathrm{SL}_2(\mathbb{Z})$-stable. Next, in view of the formulas in Theorem \ref{actionSL2ZArike}, it is natural to define
$$ x_s = \hat q \frac{p-s}{[s]}\chi^+_s - \hat q \frac{s}{[s]}\chi^-_{p-s}, \:\:\:\:\:\:\: y_s = G_s - x_s. $$
Then:
\begin{center}
\begin{tabular}{ll}
$\displaystyle \theta_1(\tau_a)(x_s) = v_{\mathcal{X}^+(s)}^{-1} x_s$, & $\displaystyle \theta_1(\tau_b)(x_s) = \xi (-1)^{s}q^{-(s^2-1)}\frac{\hat q p}{[s]}\sum_{j=1}^{p-1}(-1)^{j+1}[j][js]\!\left(x_j + y_j\right)$\\
$\displaystyle \theta_1(\tau_a)(y_s) = v_{\mathcal{X}^+(s)}^{-1}(y_s - x_s) $, & $\displaystyle \theta_1(\tau_b)(y_s) = \xi (-1)^{s}q^{-(s^2-1)}\frac{\hat q p}{[s]}\sum_{j=1}^{p-1}(-1)^{j+1}[j][js] y_j$.
\end{tabular}
\end{center}
The result follows from Lemma \ref{lemmaPropRepSL2Z}.
\end{proof}

\noindent We precise that, explicitly, the projective representation $\mathcal{W}$ has a basis $(w_s)_{1 \leq s \leq p-1}$ such that
\begin{equation}\label{FormulesActionW}
 \tau_a \, w_s = v_{\mathcal{X}^+(s)}^{-1} w_s, \:\:\:\:\: \tau_b \, w_s = \xi (-1)^{s}q^{-(s^2-1)}\frac{\hat q p}{[s]}\sum_{j=1}^{p-1}(-1)^{j+1}[j][js] w_j.
\end{equation}

\noindent The structure of the Lyubashenko-Majid representation on $\mathcal{Z}(\bar U_q)$, which by Theorem \ref{EquivalenceLMandSLF} is equivalent to the one constructed here, was described in \cite{FGST} in relation to logarithmic conformal field theory.  Theorem \ref{thDecRep} is in perfect agreement with their result.

\begin{remark}
The subspace $\mathcal{P}$ appearing in Theorem \ref{thDecRep} is spanned as a vector space by the characters of the projective $\bar U_q$-modules. Indeed, since the characters split on extensions, we have $\chi^{\mathcal{P}^+(s)} = \chi^{\mathcal{P}^-(s)} = 2(\chi^+_s + \chi^-_{p-s})$ for $1 \leq s \leq p-1$ and the simple projective modules $\mathcal{X}^{\pm}(p) = \mathcal{P}^{\pm}(p)$ give $\chi^{\pm}_p$.
\end{remark}

\subsection{A conjecture about the representation of $\mathcal{L}_{1,0}^{\mathrm{inv}}(\bar U_q)$ on $\SLF(\bar U_q)$}\label{sectionConjecture}
\indent Another natural (but harder) question is to determine the structure of $\SLF(\bar U_q)$ under the action of $\mathcal{L}_{1,0}^{\mathrm{inv}}(\bar U_q)$\footnote{A weakened version of this problem will be solved in section \ref{structRepSqS1}, where we determine the structure of $\SLF(\bar U_q)$ under the action of the subalgebra $\mathbb{C}\langle \overset{\mathcal{X}^+(2)}{W_{\! A}}, \overset{\mathcal{X}^+(2)}{W_{\! B}} \rangle \subsetneq \mathcal{L}_{1,0}^{\mathrm{inv}}(\bar U_q)$; this subalgebra generated by $\overset{\mathcal{X}^+(2)}{W_{\! A}}$ and $\overset{\mathcal{X}^+(2)}{W_{\! B}}$ is the image of the skein algebra of the torus $\mathcal{S}_q(\Sigma_1)$ by the Wilson loop map $W$, which will be defined in Chapter \ref{chapitreGraphiqueSkein}.}. As mentioned in the proof of Theorem \ref{thDecRep}, the subspace $\mathcal{P}= \vect\!\left(\chi^+_s + \chi^-_{p-s}, \chi^{\pm}_p\right)_{1 \leq s \leq p-1}$ is quite ``stable''. We propose the following conjecture.
\begin{conjecture}\label{conjSLFL10}
$\mathcal{P}$ is a $\mathcal{L}_{1,0}^{\mathrm{inv}}(\bar U_q)$-submodule of $\SLF(\bar U_q)$.
\end{conjecture}
In order to prove this conjecture one needs to find a basis or a generating set of $\mathcal{L}_{1,0}^{\mathrm{inv}}(\bar U_q)$, and then to show that $\mathcal{P}$ is stable under the action of the basis elements (or of the generating elements). Both tasks are difficult. 
\smallskip\\
\indent Recall that since $\mathcal{P}$ is an ideal of $\SLF(\bar U_q)$ which is stable under the action \eqref{actionCentreSLF} of $\mathcal{Z}(\bar U_q)$, it follows from the formulas of Proposition \ref{actionAB} and Lemma \ref{actionBmoinsUnA} that $\mathcal{P}$ is stable under the representations of $z_{A}$, $z_{B}$ and $z_{B^{-1}A}$ for all $z \in \mathcal{Z}(\bar U_q) = \mathcal{L}^{\mathrm{inv}}_{0,1}(\bar U_q)$. Also recall the wide family of invariants given in \eqref{FamilleInvL10}; we can try to test the conjecture with them. A long computation (which is not specific to $\bar U_q$) shows that
$$ \mathrm{tr}_{12}\!\left(\overset{I \otimes J}{g}\!\!\!_{12}\Phi_{12}\overset{I}{A}_1\overset{IJ}{(R')}_{12}\overset{J}{B}_2\overset{IJ}{R}_{12}\right) \triangleright \chi^K = v_J \mathrm{tr}_{13}\!\left( \overset{I \otimes K}{T}\!\!\!\!_{13} \, \overset{I \otimes K}{v}{^{-1}}\!\!\!\!\!\!\!\!\!_{13} \:\:\:\, s_{IJ,K}(\Phi)_{13} \right) $$
where $\chi^K$ is the character of $K$, $\overset{J}{v} = v_J\mathrm{id}$ (note that we may assume that $I,J,K$ are simple modules) and
$$ s_{IJ,K}(\Phi) = \mathrm{tr}_{2}\!\left( \overset{J}{g}_2 \overset{JK}{R}_{23}\Phi_{12}\overset{JK}{(R')}_{23} \right). $$
Proving that $\mathcal{P}$ is stable under the action of these invariants amounts to show symmetry properties between $s_{IJ, \mathcal{X}^{+}(s)}$ and $s_{IJ, \mathcal{X}^{-}(p-s)}$ for all simple $\bar U_q$-modules $I,J$. We have checked that it is true if $\Phi = \mathrm{id}_{I \otimes J}$ (in this case $ s_{IJ,K}(\mathrm{id}_{I \otimes J}) = s_{J,K} \mathrm{id}_{I \otimes K}$, where $s_{J,K}$ is the usual $S$-matrix) for all simple modules $I,J$, and also that it holds for $I = J = \mathcal{X}^+(2)$ with every $\Phi$.

\begin{proposition}\label{propStructureSLFSousL10inv}
1) $\SLF(\bar U_q)$ is indecomposable as a $\mathcal{L}_{1,0}^{\mathrm{inv}}(\bar U_q)$-module.\\
2) Assume that Conjecture \ref{conjSLFL10} holds. Then the $\mathcal{L}_{1,0}^{\mathrm{inv}}(\bar U_q)$-modules $\mathcal{P}$ and $\SLF(\bar U_q) / \mathcal{P}$ are simple. It follows that $\SLF(\bar U_q)$ has length $2$ as a $\mathcal{L}_{1,0}^{\mathrm{inv}}(\bar U_q)$-module.
\end{proposition}
\begin{proof}
These are basically consequences of (\ref{actionCentreSLF}) and of the multiplication rules in the GTA basis (Theorem \ref{ProduitArike}). To avoid particular cases, let $\chi^{\epsilon}_0 = 0$, $\chi^{\epsilon}_{p+1} = \chi^{-\epsilon}_1$, $\chi^{\epsilon}_{-1} = \chi^{-\epsilon}_{p-1}$ and $e_{-1} = e_{p+1} = 0$. 
\\\noindent 1) Observe that $\SLF(\bar U_q)$ is generated by $\chi^+_1 = \varepsilon$ as a $\mathcal{L}_{1,0}^{\mathrm{inv}}(\bar U_q)$-module: this is a general fact which follows immediately from Lemma \ref{actionBmoinsUnA}. Explicitly:
\[ \overset{\mathcal{X}^{\epsilon}(s)}{W}\!\!\!\!_{vB^{-1}A} \triangleright \chi^+_1 = \chi^{\epsilon}_s\chi^+_1 = \chi^{\epsilon}_s, \:\:\:\:\:\:\:\:\: V^s_{vB^{-1}A} \triangleright \chi^+_1 = G_s \chi^+_1 = G_s. \]
Write $\SLF(\bar U_q) = U_1 \oplus U_2$. At least one of the two subspaces $U_1, U_2$ necessarily contains an element of the form $u = G_1 + \sum_{i \neq 1} \lambda_i G_i + \sum_{j, \epsilon} \eta^{\epsilon}_j \chi^{\epsilon}_j$; assume that it is $U_1$. Then $(w_1^+)_A \triangleright u = \chi^+_1 \in U_1$ thanks to \eqref{actionCentreSLF}. It follows that $U_1 = \SLF(\bar U_q)$ and $U_2 = \{0\}$, as desired.
\\\noindent 2) Let $0 \neq U \subset \mathcal{P}$ be a submodule, and let $v = \sum_{j=0}^p \lambda_{j}(\chi^+_j + \chi^-_{p-j}) \in U$ with $\lambda_s \neq 0$ for some $s$. Then using Proposition \ref{actionAB} and \eqref{actionCentreSLF}, we get $(e_s)_A \triangleright v = \lambda_{s}(\chi^+_s + \chi^-_{p-s})$, and thus $\chi^+_s + \chi^-_{p-s} \in U$. Apply $\overset{\mathcal{X}^+(2)}{W}\!\!\!\!_{vB^{-1}A}$ (we use Lemma \ref{actionBmoinsUnA} and Proposition \ref{antipodeSLF}):
$$ \overset{\mathcal{X}^+(2)}{W}\!\!\!\!_{vB^{-1}A} \triangleright (\chi^+_s + \chi^-_{p-s}) = \chi^+_2(\chi^+_s + \chi^-_{p-s}) = (\chi^+_{s-1} + \chi^-_{p-s+1}) + (\chi^+_{s+1} + \chi^-_{p-s-1}). $$
Hence:
\begin{align*}
&(e_{s-1})_A\overset{\mathcal{X}^+(2)}{W}\!\!\!\!_{vB^{-1}A} \triangleright (\chi^+_s + \chi^-_{p-s}) = \chi^+_{s-1} + \chi^-_{p-s+1}\\
&(e_{s+1})_A\overset{\mathcal{X}^+(2)}{W}\!\!\!\!_{vB^{-1}A} \triangleright (\chi^+_s + \chi^-_{p-s}) = \chi^+_{s+1} + \chi^-_{p-s-1}.
\end{align*}
It follows that $\chi^+_{s-1} + \chi^-_{p-s+1}, \chi^+_{s+1} + \chi^-_{p-s-1} \in U$. Continuing like this, one gets step by step that all the basis vectors belong to $U$, hence $U = \mathcal{P}$. 
\\Next, let $\overline{G}_s$ and $\overline{\chi}^+_s$ be the classes of $G_s$ and $\chi^+_s$ modulo $\mathcal{P}$ (with $\overline{\chi}^+_0 = \overline{\chi}^+_p = 0$). Let $0 \neq U \subset \SLF(\bar U_q) / \mathcal{P}$ be a submodule and $w = \sum_{j = 1}^{p-1} \nu_{j}\overline{G}_{j} + \sigma_{j}\overline{\chi}^+_{j} \in U$ be non-zero. If all the $\nu_j$ are $0$, then there exists $\sigma_s \neq 0$ and $(e_s)_A \triangleright w = \sigma_s\overline{\chi}^+_s \in U$. If one of the $\nu_j$, say $\nu_s$, is non-zero, then $(w^+_s)_A \triangleright w = \nu_s \overline{\chi}^+_s  \in U$. In both cases we get $\overline{\chi}^+_s \in U$. Now we proceed as previously:
$$ (e_{s-1})_A\overset{\mathcal{X}^+(2)}{W}\!\!\!\!_{vB^{-1}A} \triangleright \overline{\chi}^+_s = \overline{\chi}^+_{s-1}, \:\:\:\:\:\:\: (e_{s+1})_A\overset{\mathcal{X}^+(2)}{W}\!\!\!\!_{vB^{-1}A} \triangleright \overline{\chi}^+_s = \overline{\chi}^+_{s+1}. $$
Thus we get step by step that $\overline{\chi}^+_{j} \in U$ for all $j$. Apply $V^1_{vB^{-1}A}$ and use Corollary \ref{coroVermaIdeal}:
$$ V_{vB^{-1}A}^1 \triangleright \overline{\chi}^+_{j} = G_1\chi^+_{j} + \mathcal{P} = [j]G_j + \mathcal{P}. $$
It follows that $\overline{G}_j \in U$ for all $j$, and thus $U = \SLF(\bar U_q) / \mathcal{P}$ as desired.
\end{proof}


\begin{remark}
As suggested to me by A. Gainutdinov, we can extend Conjecture \ref{conjSLFL10} to any finite dimensional, factorizable, ribbon Hopf algebra $H$. Let $\mathcal{P}_H = \mathrm{vect}\bigl( \chi^P \bigr)_{P \in \mathrm{Proj}(H)} \subset \mathrm{SLF}(H)$, where $\mathrm{Proj}(H) \subset \mathrm{mod}_l(H)$ is the subcategory (in fact, the ideal) of finite dimensional projective $H$-modules and $\chi^P = \mathrm{tr}\bigl( \overset{P}{T} \bigr)$ is the character of $P$ (we can restrict $P$ to be some PIM).

\medskip

\noindent \textbf{Generalized Conjecture \ref{conjSLFL10}}. {\em $\mathcal{P}_H$ is a $\mathcal{L}_{1,0}^{\mathrm{inv}}(H)$-submodule of $\SLF(H)$.}

\medskip

\noindent Note that it is known that $\mathcal{P}_H$ is stable under the Lyubashenko-Majid $\mathrm{SL}_2(\mathbb{Z})$-action (\cite{CW}, also see \cite{GR2}). Hence, by Theorem \ref{EquivalenceLMandSLF}, $\mathcal{P}_H$ is stable under the action of $v_A^{-1}$ and $v_B^{-1}$.
\end{remark}

\chapter{$\mathcal{L}_{g,n}(H)$ and projective representations of mapping class groups}\label{chapitreLgnMCG}

\begin{figure}[h]
\centering
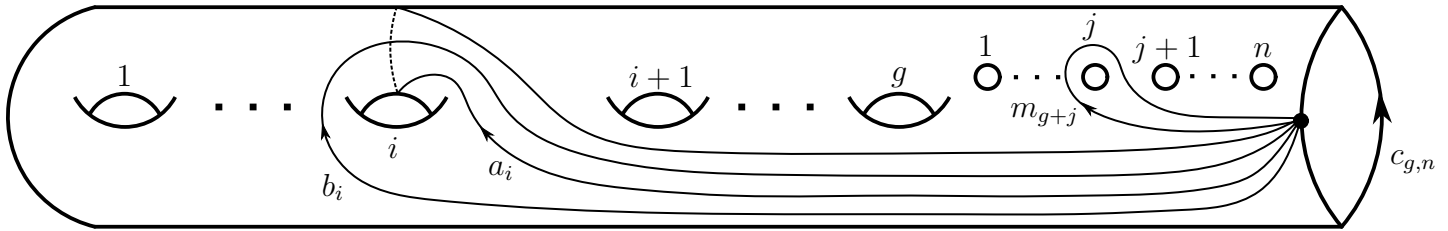
\caption{Surface $\Sigma_{g,n}^{\mathrm{o}}$ with basepoint ($\bullet$) and generators for $\pi_1(\Sigma_{g,n}^{\mathrm{o}})$.}
\label{surfaceAvecCourbes}
\end{figure}

Let $\Sigma_{g,n}$ be the compact oriented surface of genus $g$ with $n$ open disks removed. Let $D \subset \Sigma_{g,n}$ be an open disk, then we define $\Sigma_{g,n}^{\mathrm{o}} = \Sigma_{g,n} \setminus D$. Of course, $\Sigma_{g,n}^{\mathrm{o}} = \Sigma_{g, n+1}$, but the boundary circle $c_{g,n}$ induced by the deletion of $D$ plays a particular role since we put a basepoint on it, see Figure \ref{surfaceAvecCourbes}. In this chapter, we consider the algebra $\mathcal{L}_{g,n}(H)$ associated to $\Sigma_{g,n}^{\mathrm{o}}$ (as everywhere in this thesis, $H$ denotes a finite dimensional, factorizable, ribbon Hopf algebra). 

\smallskip 

\begin{figure}[h]
\centering
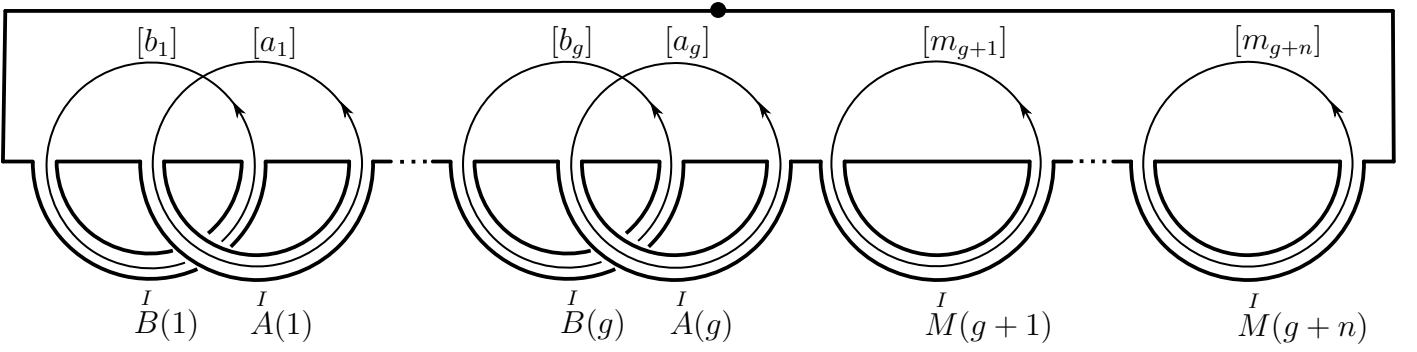
\caption{$\Sigma_{g,n}^{\mathrm{o}}$ represented as a thickened graph and matrices of generators of $\mathcal{L}_{g,n}(H)$.}
\label{surfaceAvecMatrices}
\end{figure}

Figure \ref{surfaceAvecMatrices} is the picture that one should always keep in mind\footnote{Compared to the Figure 1 of \cite{Fai18c}, we have done a $180^{\circ}$-rotation around the horizontal axis of $\mathbb{R}^3$, in order to have the handles at the bottom of the Figure. The reason of this change comes from the definition of the graphical calculus and the Wilson loop map in Chapter \ref{chapitreGraphiqueSkein}.}. In this picture, we see $\Sigma_{g,n}^{\mathrm{o}}$ as a thickening (\textit{i.e.} tubular neighborhood) of the embedded oriented graph 
\begin{equation*}
\Gamma_{g,n} = \bigl( \{\bullet\},  \{ b_1, a_1, \ldots, b_g, a_g, m_{g+1}, \ldots,m_{g+n}\}\bigr)
\end{equation*}
where the loops $b_i, a_i, m_j$ generating the free group $\pi_1(\Sigma_{g,n}^{\mathrm{o}})$ are represented in Figure \ref{surfaceAvecCourbes}. Note that with these generators, the loop $c_{g,n}$ induced by the deletion of the disc $D$ is expressed as
\begin{equation}\label{boundaryCirclePi1}
c_{g,n} = b_1 a_1^{-1} b_1^{-1} a_1 \ldots b_g a_g^{-1} b_g^{-1} a_g m_{g+1} \ldots m_{g+n}.
\end{equation}

To each generating loop, or equivalently to each handle, is associated a family of matrices, indexed by the $H$-modules and whose coefficients are generators of the algebra $\mathcal{L}_{g,n}(H)$. The defining relations are given in \eqref{PresentationLgn}, following \cite{AGS, BR, AGS2} (modulo the same remarks that in the introduction of Chapter \ref{chapitreTore}). We define $\mathcal{L}_{g,n}(H)$ as a braided tensor product, as in \cite{AS2}. This has the advantage to show immediately that $\mathcal{L}_{g,n}(H)$ is a $H$-module-algebra and to emphasize the role of the two building blocks of the theory, namely $\mathcal{L}_{0,1}(H)$ and $\mathcal{L}_{1,0}(H)$.

\smallskip

\indent The main results of this chapter are
\begin{itemize}
\item The construction of a representation $\mathrm{Inv}(V)$ of the algebra of invariant elements $\mathcal{L}_{g,n}^{\mathrm{inv}}(H)$ for any representation $V$ of $\mathcal{L}_{g,n}(H)$ (Theorem \ref{thmInv}). Note that the matrices $\overset{I}{C}_{g,n}$ (which correspond to the loop $c_{g,n}$, see \eqref{boundaryLoop10}) used for the proof of that theorem already appeared in \cite{alekseev} (with $H = U_q(\mathfrak{g})$), but here we need to generalize and adapt the construction of the representation to our assumptions on $H$.
\item The construction of a projective representation of the mapping class group of $\Sigma_g = \Sigma_{g,0}$ on $\mathrm{Inv}\bigl( (H^*)^{\otimes g} \bigr)$ (Theorem \ref{thmRepMCG}). This implements and generalizes to a non-semisimple setting an idea of Alekseev--Schomerus \cite[Section 9]{AS}.
\item The explicit formulas for the representation of the Dehn twists about the circles depicted in Figure \ref{figureCourbesCanoniques} (Theorem \ref{formulesExplicites}).
\item The equivalence of the representation of Theorem \ref{thmRepMCG} with the one constructed by Lyubashenko using categorical techniques in \cite{lyu95b, lyu96} (Theorem \ref{thmEquivalenceReps}).
\end{itemize}

Most of the material presented in this chapter is the content of \cite{Fai18c}. Here we added section \ref{lienAvecLesLGFT} to explain how $\mathcal{L}_{g,n}(H)$ is related to the work of \cite{BFKB, BFKB2}. We also added section \ref{sectionNorma} in which we define the normalization of a simple closed curve and this allows us to define the lift of any simple loop in section \ref{sectionLifts} (the difference with \cite{Fai18c} is explained in Remark \ref{remarkLiftArticle}).


\section{Definition and properties of $\mathcal{L}_{g,n}(H)$}

\subsection{Braided tensor product and definition of $\mathcal{L}_{g,n}(H)$}\label{defLgn}
Let $\mathrm{mod}_r(H)$ be the category of finite dimensional right $H$-modules (or, equivalently, of finite dimensional left $\mathcal{O}(H)$-comodules). The braiding in $\mathrm{mod}_r(H)$ is given by:
\[ \fonc{c_{I,J}}{I \otimes J}{J \otimes I}{v \otimes w}{w \cdot a_i \otimes v \cdot b_i} \]
with $R = a_i \otimes b_i$. Let $(A, m_A, 1_A)$ and $(B, m_B, 1_B)$ be two algebras in $\mathrm{mod}_r(H)$ (that is, $H$-module-algebras), and define:
\begin{align*}
&m_{A \widetilde{\otimes} B} = (m_A \otimes m_B) \circ (\mathrm{id}_A \otimes c_{B,A} \otimes \mathrm{id}_B) \: : \: (A \otimes B) \otimes (A \otimes B) \to A \otimes B,\\  
&1_{A \widetilde{\otimes} B} = 1_A \otimes 1_B \: : \: \mathbb{C} \to A \otimes B.
\end{align*}
This endows $A \otimes B$ with a structure of algebra in $\mathrm{mod}_r(H)$, denoted $A\,\widetilde{\otimes}\,B$ and called braided tensor product of $A$ and $B$ (see \cite[Lemma 9.2.12]{majid}).  Note that $\widetilde{\otimes}$ is associative.
\\\indent There are two canonical algebra embeddings $j_A, j_B : A, B \hookrightarrow A\,\widetilde{\otimes}\,B$ respectively defined by $j_A(x) = x \otimes 1_B$, $j_B(y) = 1_A \otimes y$. We identify $x \in A$ (resp. $y \in B$) with $j_A(x) \in A\,\widetilde{\otimes}\,B$ (resp. $j_B(y)$). Under these identifications, the multiplication rule in $A\,\widetilde{\otimes}\,B$ is entirely given by:
\begin{equation}\label{braidedProduct}
\forall \, x \in A, \forall \, y \in B, \:\: yx = (x \cdot a_i)(y \cdot b_i) .
\end{equation}
\indent Since $\mathcal{L}_{0,1}(H)$ and $\mathcal{L}_{1,0}(H)$ are algebras in $\mathrm{mod}_r(H)$, we can apply the braided tensor product to them.
\begin{definition}\label{definitionLgn}
$\mathcal{L}_{g,n}(H)$ is the $H$-module-algebra $\mathcal{L}_{1,0}(H)^{\widetilde{\otimes} g} \, \widetilde{\otimes} \, \mathcal{L}_{0,1}(H)^{\widetilde{\otimes} n}$.
\end{definition}
\noindent It is essential to keep in mind that the $H$-module-algebra $\mathcal{L}_{g,n}(H)$ is associated with the surface $\Sigma_{g,n}^{\mathrm{o}} = \Sigma_{g,n} \setminus D$; in order to make this precise we now define the matrices introduced in Figure \ref{surfaceAvecMatrices}. There are canonical algebra embeddings $j_i : \mathcal{L}_{1,0}(H) \hookrightarrow \mathcal{L}_{g,n}(H)$ for $1 \leq i \leq g$ and $j_i : \mathcal{L}_{0,1}(H) \hookrightarrow \mathcal{L}_{g,n}(H)$ for $g+1 \leq i \leq g+n$, given by $j_i(x) = 1^{\otimes i-1} \otimes x \otimes 1^{\otimes g+n-i}$. Define 
\[ \overset{I}{A}(i) = j_i(\overset{I}{A}), \:\:\: \overset{I}{B}(i) = j_i(\overset{I}{B}) \: \text{ for } 1 \leq i \leq g \:\: \text{ and } \:\: \overset{I}{M}(i) = j_i(\overset{I}{M}) \: \text{ for } g+1 \leq i \leq g+n. \]
\noindent The right action of $H$ on $\mathcal{L}_{g,n}(H)$ (or equivalently the left coaction $\Omega$ of $\mathcal{O}(H)$) is of course
\begin{equation}\label{actionHLgn}
\overset{I}{U}(i) \cdot h = \overset{I}{h'} \overset{I}{U}(i) \overset{I}{S(h'')} \:\:\:\:\:\:\:\: \text{(or equivalently } \Omega\bigl(\overset{I}{U}(i)\bigr) = \overset{I}{T} \overset{I}{U}(i) S(\overset{I}{T}) \text{ )}
\end{equation}
where $U$ is $B,A$ or $M$. By \eqref{naturaliteL01} and \eqref{naturaliteL10}, if $f : I \to J$ is a morphism of $H$-modules it holds
\begin{equation}\label{naturaliteLgn}
\overset{J}{B}(i) f = f \overset{I}{B(i)}, \:\:\:\:\:\: \overset{J}{A}(i) f = f \overset{I}{A}(i), \:\:\:\:\:\: \overset{J}{M}(j) f = f \overset{I}{M}(j)
\end{equation}
for all $i,j$, where we identify $f$ with its matrix. We call this relation the naturality of the (families of) matrices $\overset{I}{B}(i), \overset{I}{A}(i), \overset{I}{M}(j)$. 

\smallskip

\indent Relation \eqref{braidedProduct} indicates that $\mathcal{L}_{g,n}(H)$ is an exchange algebra. Let us write the exchange relations in a matrix form. Let $U,V$ be $B$ or $A$ or $M$. 
Then, by definition of the right action and by \eqref{braidedProduct}:
\begin{equation}\label{eqEchangeLgn}
\overset{J}{V}(j)_2 \overset{I}{U}(i)_1 = (\overset{I}{a'_k})_1 \overset{I}{U}(i)_1 \overset{I}{S(a''_k)}_1 (\overset{J}{b'_k})_2 \overset{J}{V}(j)_2 \overset{J}{S(b''_k)}_2 
= (\overset{I}{a_l})_1 \overset{IJ}{R}_{12} \overset{I}{U}(i)_1 \overset{IJ}{R}{^{-1}_{12}}\overset{J}{V}(j)_2 \overset{IJ}{R}_{12} \overset{J}{S(b_l)}_2
\end{equation}
where for the second equality we applied properties of the $R$-matrix and obvious commutation relations in $\mathrm{End}_{\mathbb{C}}(I) \otimes \mathrm{End}_{\mathbb{C}}(J) \otimes \mathcal{L}_{g,n}(H)$. Using that $a_m a_l \otimes S(b_l) b_m = 1 \otimes 1$ together with obvious commutation relations, we obtain the desired exchange relation:
\begin{equation*}
\overset{IJ}{R}_{12} \overset{I}{U}(i)_1 \overset{IJ}{R}{^{-1}_{12}}\overset{J}{V}(j)_2 = \overset{J}{V}(j)_2 \overset{IJ}{R}_{12} \overset{I}{U}(i)_1 \overset{IJ}{R}{^{-1}_{12}}.
\end{equation*}
To sum up, the presentation of $\mathcal{L}_{g,n}(H)$ by generators and relations is:
\begin{equation}\label{PresentationLgn}
  \left\{
      \begin{aligned}
& \overset{I \otimes J}{U}\!\!(i)_{12} = \overset{I}{U}(i)_1\,(\overset{IJ}{R'})_{12}\,\overset{J}{U}(i)_2\, (\overset{IJ}{R'}){_{12}^{-1}}& &\text{ for } 1 \leq i \leq g+n\\
& \overset{IJ}{R}_{12}\,\overset{I}{U}(i)_1\, \overset{IJ}{R}{^{-1}_{12}} \,\overset{J}{V}(j)_2 = \overset{J}{V}(j)_2\, \overset{IJ}{R}_{12}\,\overset{I}{U}(i)_1\, \overset{IJ}{R}{^{-1}_{12}}& &\text{ for } 1 \leq i < j \leq g+n\\
& \overset{IJ}{R}_{12}\,\overset{I}{B}(i)_1\, (\overset{IJ}{R'})_{12}\,\overset{J}{A}(i)_2 = \overset{J}{A}(i)_2\, \overset{IJ}{R}_{12}\,\overset{I}{B}(i)_1\, \overset{IJ}{R}{_{12}^{-1}} & &\text{ for } 1 \leq i \leq g
\end{aligned}
    \right.
\end{equation}
where $U(i)$ (resp. $V(i)$) is $A(i)$ or $B(i)$ if $1 \leq i \leq g$ and is $M(i)$ if $g+1 \leq i \leq g+n$. These are relations between matrices in $\mathcal{L}_{g,n}(H) \otimes \mathrm{Mat}_{\dim(I)}(\mathbb{C}) \otimes \mathrm{Mat}_{\dim(J)}(\mathbb{C})$ (for all finite dimensional $I,J$) which imply relations among elements of $\mathcal{L}_{g,n}(H)$ (the coefficients of these matrices). Such a presentation was first introduced in \cite{alekseev} and \cite{AGS}. Recall that the first line of relations is the $\mathcal{L}_{0,1}(H)$-fusion relation on each loop, the second line is the exchange relation of the braided tensor product and the third line is the $\mathcal{L}_{1,0}(H)$-exchange-relation. These are the same relations as in \cite{BNR}, except that $A(i)$ and $B(i)$ are switched for all $1 \leq i \leq g$; the ones of \cite[Def 12]{AGS2} and \cite[eqs (3.11)--(3.21)]{AS} are different, due to a different choice of the action of $H$ on $\mathcal{L}_{1,0}(H)$, but yield the same algebra when $H$ is semisimple. Thanks to these relations, we see that a generic element in $\mathcal{L}_{g,n}(H)$ is a linear combination of elements of the form
\[ \overset{I_1}{B}(1)^{i_1}_{j_1} \overset{I_2}{A}(1)^{i_2}_{j_2} \ldots \overset{I_{2g-1}}{B}\!\!\!(g)^{i_{2g-1}}_{j_{2g-1}} \overset{I_{2g}}{A}(g)^{i_{2g}}_{j_{2g}}  \overset{I_{2g+1}}{M}\!\!(g+1)^{i_{2g+1}}_{j_{2g+1}} \ldots \overset{I_{2g+n}}{M}\!\!(g+n)^{i_{2g+n}}_{j_{2g+n}}. \]
Note that the content of Remark \ref{remarkRestrictionCoeffL01} also applies to $\mathcal{L}_{g,n}(H)$: in practice, we can assume that the representations labelling the matrices belong to a set $\mathcal{G}$ of well-chosen $H$-modules.
For instance, if $H = \bar U_q(\mathfrak{sl}_2)$, we take $\mathcal{G} = \bigl\{ \mathcal{X}^+(2) \bigr\}$. 

\noindent \textbf{Notation.}~~ Let $\overset{I}{N} = \overset{I}{v}{^m} \overset{I}{N}{^{n_1}_1} \ldots \overset{I}{N}{^{n_l}_l} \in \mathrm{Mat}_{\dim(I)}\!\left(\mathcal{L}_{g,n}(H)\right)$, where $m, n_i \in \mathbb{Z}$  and each $N_i$ is one of the $A(j), B(j), M(k)$ for some $j$ or $k$. By definition of the right action on $\mathcal{L}_{g,n}(H)$, we have a morphism of $H$-modules 
\begin{equation}\label{morphismeEmbed}
\fonc{j_N}{\mathcal{L}_{0,1}(H)}{\mathcal{L}_{g,n}(H)}{\overset{I}{M}}{\overset{I}{N}}.
\end{equation}
Let $x \in \mathcal{L}_{0,1}(H)$, then we denote 
\begin{equation}\label{notationEmbed}
x_N = j_N(x).
\end{equation}
Since we identify $\mathcal{L}_{0,1}(H)$ with $H$ we also use this notation when $x \in H$. Note that if $x \in \mathcal{L}_{0,1}^{\mathrm{inv}}(H) \cong \mathcal{Z}(H)$, then $x_N \in \mathcal{L}_{g,n}^{\mathrm{inv}}(H)$. The following lemma is an obvious fact.
\begin{lemma}\label{injectionFusion}
If $N$ satisfies the fusion relation of $\mathcal{L}_{0,1}(H)$, $\overset{I \otimes J}{N}\!\!_{12} = \overset{I}{N}(i)_1\,\overset{IJ}{(R')}_{12}\,\overset{J}{N}(i)_2\, \overset{IJ}{(R')}{_{12}^{-1}}$, then $j_N$ is a morphism of $H$-module-algebras: $(xy)_N = x_N y_N$.
\end{lemma}
\noindent See e.g. \eqref{actionHsurFormes} for an application of this lemma.

\subsection{The Alekseev isomorphism}
\indent Consider the tensor product algebra $\mathcal{L}_{1,0}(H)^{\otimes g} \otimes \mathcal{L}_{0,1}(H)^{\otimes n}$. We have canonical algebra embeddings $j_i : \mathcal{L}_{1,0}(H) \hookrightarrow \mathcal{L}_{1,0}(H)^{\otimes g} \otimes \mathcal{L}_{0,1}(H)^{\otimes n}$ for $1 \leq i \leq g$ and $j_i : \mathcal{L}_{0,1}(H) \hookrightarrow \mathcal{L}_{1,0}(H)^{\otimes g} \otimes \mathcal{L}_{0,1}(H)^{\otimes n}$ for $g+1 \leq i \leq g+n$, defined by $j_i(x) = 1^{\otimes i-1} \otimes x \otimes 1^{\otimes g+n-i}$. Define $\overset{I}{\underline{A}}(i) = j_i(\overset{I}{A})$, $\overset{I}{\underline{B}}(i) = j_i(\overset{I}{B})$ for $1 \leq i \leq g$ and $\overset{I}{\underline{M}}(i) = j_i(\overset{I}{M})$ for $g+1 \leq i \leq g+n$. We underline these matrices to avoid confusion with prior matrices having coefficients in $\mathcal{L}_{g,n}(H)$. By definition, the exchange relation between copies in $\mathcal{L}_{1,0}(H)^{\otimes g} \otimes \mathcal{L}_{0,1}(H)^{\otimes n}$ is simply
\[ \overset{I}{\underline{U}}(i)_1\,\overset{J}{\underline{V}}(j)_2 = \overset{J}{\underline{V}}(j)_2\,\overset{I}{\underline{U}}(i)_1 \]
where $i \neq j$, $\underline{U}(i), \underline{V}(i)$ is $\underline{A}(i)$ or $\underline{B}(i)$ if $1 \leq i \leq g$ and is $\underline{M}(i)$ if $g+1 \leq i \leq g+n$.
\smallskip\\
\indent The next result is due to Alekseev (see \cite{alekseev}). Consider the matrices $\overset{I}{M}{^{(-)}} = \Psi_{0,1}^{-1}(\overset{I}{L}{^{(-)}})$ and $\overset{I}{C}{^{(-)}} = \Psi_{1,0}^{-1}(\overset{I}{L}{^{(-)}}\overset{I}{\widetilde{L}}{^{(-)}})$ (recall \eqref{L} and \eqref{defLTilde}). Let
\begin{equation}\label{matricesAlekseev}
\begin{array}{ll}
\overset{I}{\Lambda}_1 = \mathbb{I}_{\dim(I)}, & \:\: \overset{I}{\Lambda}_i = \overset{I}{\underline{C}}{^{(-)}}(1) \ldots \overset{I}{\underline{C}}{^{(-)}}(i-1)  \:\: \text{ for } 2 \leq i \leq g+1,\\
\overset{I}{\Gamma}_{g+1} = \overset{I}{\Lambda}_{g+1},  & \:\: \overset{I}{\Gamma}_i = \overset{I}{\Lambda}_{g+1} \overset{I}{\underline{M}}{^{(-)}}(g+1) \ldots \overset{I}{\underline{M}}{^{(-)}}(i-1) \:\: \text{ for } g+2 \leq i \leq g+n.
\end{array}
\end{equation}
be matrices with coefficients in $\mathcal{L}_{1,0}(H)^{\otimes g} \otimes \mathcal{L}_{0,1}(H)^{\otimes n}$ (with $\mathbb{I}_s$ the identity matrix of size $s$).

\begin{proposition}\label{isoAlekseev}
The map
\[\begin{array}{crll}\alpha_{g,n} :& \mathcal{L}_{g,n}(H) = \mathcal{L}_{1,0}(H)^{\widetilde{\otimes} g} \, \widetilde{\otimes} \, \mathcal{L}_{0,1}(H)^{\widetilde{\otimes} n} & \rightarrow & \mathcal{L}_{1,0}(H)^{\otimes g} \otimes \mathcal{L}_{0,1}(H)^{\otimes n} \\
                         & \overset{I}{A}(i) &\mapsto &  \overset{I}{\Lambda}_i\,\overset{I}{\underline{A}}(i)\,\overset{I}{\Lambda}{_i^{-1}} \:\: \text{ for } 1 \leq i \leq g \\
                         & \overset{I}{B}(i) &\mapsto &  \overset{I}{\Lambda}_i\,\overset{I}{\underline{B}}(i)\,\overset{I}{\Lambda}{_i^{-1}} \:\: \text{ for } 1 \leq i \leq g \\
                         & \overset{I}{M}(i) &\mapsto &  \overset{I}{\Gamma}_i\,\overset{I}{\underline{M}}(i)\,\overset{I}{\Gamma}{_i^{-1}} \:\: \text{ for } g+1 \leq i \leq g+n \\
\end{array}\]
is an isomorphism of algebras, which we call the Alekseev isomorphism.
\end{proposition}
\begin{proof}
In order to show that it is a morphism of algebras, one must check using various exchange relations that the defining relations \eqref{PresentationLgn} of $\mathcal{L}_{g,n}(H)$ are preserved under $\alpha_{g,n}$. This is a straightforward but tedious task and we will not give the details. Let us prove that $\alpha_{g,n}$ is bijective. We first show that $\alpha_{g,0}$ is surjective for all $g$ by induction. For $g=1$, $\alpha_{1,0}$ is the identity. For $g \geq 2$, we embed $\mathcal{L}_{g-1,0}(H)$ in $\mathcal{L}_{g,0}(H)$ in an obvious way by $\overset{I}{A}(i) \mapsto \overset{I}{A}(i)$ and $\overset{I}{B}(i) \mapsto \overset{I}{B}(i)$ for $1 \leq i \leq g-1$. Then the restriction of $\alpha_{g,0}$ to $\mathcal{L}_{g-1,0}(H)$ is $\alpha_{g-1,0}$, and by induction we assume that $\alpha_{g-1,0}(\mathcal{L}_{g-1,0}(H)) = \mathcal{L}_{1,0}(H)^{\otimes g-1}$. Since $\overset{I}{\Lambda}_i \in \mathrm{Mat}_{\dim(I)}\!\left( \mathcal{L}_{1,0}(H)^{\otimes i-1} \otimes \mathbb{C}^{\otimes g+1-i} \right)$, there exists matrices $\overset{I}{\mathcal{N}}_i$ $(1 \leq i \leq g)$ such that $\alpha_{g,0}(\overset{I}{\mathcal{N}}_i) = \overset{I}{\Lambda}_i$. Then $\alpha_{g,0}(\overset{I}{\mathcal{N}}{^{-1}_{i}} \overset{I}{U}(i) \overset{I}{\mathcal{N}}_{i}) = \overset{I}{\underline{U}}(i)$, with $U = A$ or $B$ and $\alpha_{g,0}$ is surjective. Similarly, for $g$ fixed and $n \geq 1$, we can embed $\mathcal{L}_{g,n-1}(H)$ into $\mathcal{L}_{g,n}(H)$ and reproduce the same reasoning. Hence $\alpha_{g,n}$ is surjective for all $g,n$. Since the domain and the range of $\alpha_{g,n}$ have the same dimension, it is an isomorphism.
\end{proof}

\indent We generalize the isomorphisms $\Psi_{0,1}$ and $\Psi_{1,0}$ by
\begin{equation}\label{isoPsi}
\Psi_{g,n} = \left(\Psi_{1,0}^{\otimes g} \otimes \Psi_{0,1}^{\otimes n}\right) \circ \alpha_{g,n} \: : \: \mathcal{L}_{g,n}(H) \overset{\sim}{\rightarrow} \mathcal{H}(\mathcal{O}(H))^{\otimes g} \otimes H^{\otimes n}.
\end{equation}
In particular $\mathcal{L}_{g,0}(H)$ is a matrix algebra, since $\mathcal{H}(\mathcal{O}(H))$ is. 
\smallskip\\
\indent Thanks to $\Psi_{g,n}$, the representation theory of $\mathcal{L}_{g,n}(H)$ is entirely determined by the representation theory of $H$. Indeed, the only indecomposable (and simple) representation of $\mathcal{H}(\mathcal{O}(H)) \cong \End_{\mathbb{C}}(H^*)$ is $H^*$, thus it follows that the indecomposable representations of $\mathcal{L}_{g,n}(H)$ are of the form
\[ (H^*)^{\otimes g} \otimes I_1 \otimes \ldots \otimes I_n \]
where $I_1, \ldots, I_n$ are indecomposable representations of $H$. We will denote the action of $\mathcal{L}_{g,n}(H)$ on $(H^*)^{\otimes g} \otimes I_1 \otimes \ldots \otimes I_n$ by $\triangleright$, namely:
\begin{equation}\label{actionTriangle}
x \triangleright (\varphi_1 \otimes \ldots \otimes \varphi_g \otimes v_1 \otimes \ldots \otimes v_n) = \Psi_{g,n}(x) \cdot (\varphi_1 \otimes \ldots \otimes \varphi_g \otimes v_1 \otimes \ldots \otimes v_n)
\end{equation}
for $x \in \mathcal{L}_{g,n}(H)$, where $\cdot$ is the action component-by-component of $\Psi_{g,n}(x)$ on $(H^*)^{\otimes g} \otimes I_1 \otimes \ldots \otimes I_n$.

\subsection{$\mathcal{L}_{g,n}(H)$ as an algebra of functions and LGFT}\label{lienAvecLesLGFT}
\indent In this section, we discuss the fact that $\mathcal{L}_{g,n}(H)$ is the algebra of gauge fields of a lattice gauge field theory (LGFT) as defined in \cite{BFKB, BFKB2}; this is independent of the rest of the text and will not be used elsewhere. 

\smallskip

First, we describe $\mathcal{L}_{g,n}(H)$ as an algebra of functions. Let $E_{g,n} = \{b_1, a_1, \ldots, b_g, a_g, m_{g+1}, \ldots, m_{g+n}\}$\footnote{We endow $E_{g,n}$ with the total order $b_1 < a_1 < \ldots < b_g < a_g < m_{g+1} < \ldots < m_{g+n}$.} and let $\mathcal{F}_{g,n}(H)$ be the vector space $\bigotimes_{e \in E_{g,n}} H^*_e$, where $H_e^*$ is a copy of $H^*$ labelled by $e$. For $e \in E_{g,n}$ and $\varphi \in H^*$, define an element $\varphi_e \in \mathcal{F}_{g,n}(H)$ by
\[ \varphi_e\!\left( \bigotimes_{a \in E} x_a \right) = \varphi(x_e) \prod_{a \in E_{g,n}\setminus \{e\}} \varepsilon(x_a). \]
Consider the linear map $f : \mathcal{L}_{g,n}(H) \to \mathcal{F}_{g,n}(H)$ defined by
\[ \overset{I_1}{B}(1)^{i_1}_{j_1} \overset{I_2}{A}(1)^{i_2}_{j_2} \ldots \overset{I_{2g-1}}{B}\!\!\!(g)^{i_{2g-1}}_{j_{2g-1}} \overset{I_{2g}}{A}(g)^{i_{2g}}_{j_{2g}}  \overset{I_{2g+1}}{M}\!\!(g+1)^{i_{2g+1}}_{j_{2g+1}} \ldots \overset{I_{2g+n}}{M}\!\!(g+n)^{i_{2g+n}}_{j_{2g+n}} \mapsto \overset{I_1}{T}{^{i_1}_{j_1}} \otimes \ldots \otimes \overset{I_{2g+n}}{T}{^{i_{2g+n}}_{j_{2g+n}}}. \]
Thanks to the Alekseev isomorphism, $\dim\!\left(\mathcal{L}_{g,n}(H)\right) = \dim\!\left( \mathcal{F}_{g,n}(H) \right)$, and thus $f$ is an isomorphism of vector spaces. We define a structure of right $H$-module-algebra (with product denoted by $\ast$) on $\mathcal{F}_{g,n}(H)$ by requiring $f$ to be an isomorphism of right $H$-module-algebras.

\begin{proposition}
The right $H$-module-algebra $\mathcal{F}_{g,n}(H)$ is generated by the elements $\varphi_e$, with $\varphi \in H^*$ and $e \in E_{g,n}$. These elements satisfy
\begin{enumerate}
\item $\varphi_e \cdot h = \varphi\bigl(h' ? S(h'')\bigr)_e$ (with $h \in H$),
\item $(\varphi_1)_{b_1} \ast (\varphi_2)_{a_1} \ast \ldots \ast (\varphi_{2g-1})_{b_g} \ast (\varphi_{2g})_{a_g} \ast (\varphi_{2g+1})_{m_{g+1}} \ast \ldots \ast (\varphi_{2g+n})_{m_{g+n}} = \varphi_1 \otimes \ldots \otimes \varphi_{2g+n}$,
\item $\varphi_{v_{\beta}} \ast \psi_{u_{\alpha}} = \psi\bigl( s_i s_j ? S(s_k) s_l \bigr)_{u_{\alpha}} \ast \varphi\bigl( t_j t_k ? t_l S(t_i) \bigr)_{v_{\beta}}$ where $u,v$ are $a$ or $b$ or $m$ and $\alpha < \beta$,
\item $\varphi_{a_i} \ast \psi_{b_i} = \psi\bigl( s_i s_j ? t_k s_l \bigr)_{b_i} \ast \varphi\bigl( t_j s_k ? t_l S(t_i) \bigr)_{a_i}$ for all $1 \leq i \leq g$,
\item $\varphi_e \ast \psi_e = \bigl[\varphi\bigl( ? t_j S(t_i) \bigr) \psi\bigl(s_i ? s_j\bigr)\bigr]_e$ for all $e \in E_{g,n}$,
\end{enumerate}
where exceptionally we denote $R = s_i \otimes t_i$ to avoid confusion between the usual notation $R = a_i \otimes b_i$ and the loops $a_i$ and $b_i$. These formulas allow one to compute the product of any two elements in $\mathcal{F}_{g,n}(H)$.
\end{proposition}
\begin{proof}
Since $\overset{\mathbb{C}}{U}(i) = 1$ for any $U = A, B, M$ and $1 \leq i \leq 2g+n$, we have
\[ f(\overset{I}{B}(i){^j_k}) = (\overset{I}{T}{^j_k})_{b_i}, \:\:\:\:\:\:\: f(\overset{I}{A}(i){^j_k}) = (\overset{I}{T}{^j_k})_{a_i}, \:\:\:\:\:\:\: f(\overset{I}{M}(i){^j_k}) = (\overset{I}{T}{^j_k})_{m_i}. \]
This implies the first claim and the first and second equalities (recall that any $\varphi \in H^*$ is a linear combination of matrix coefficients $\overset{I}{T}{^i_j}$). For the third equality, \eqref{eqEchangeLgn} gives
\begin{align*}
(\overset{J}{T}{^a_b})_{v_{\beta}} \ast (\overset{I}{T}{^c_d})_{u_{\alpha}} &= f\!\left(\overset{J}{V}(\beta)_1 \overset{I}{U}(\alpha)_2\right)^{ac}_{bd} = f\!\left((\overset{I}{s_i})_2 \overset{IJ}{R}_{21} \overset{I}{U}(\alpha)_2 \overset{IJ}{R}{^{-1}_{21}} \overset{J}{V}(\beta)_1 \overset{IJ}{R}_{21} \overset{J}{S(t_i)}_1\right)^{ac}_{bd}\\
&= f\!\left(\overset{I}{s_i} \overset{I}{s_j}\overset{I}{U}(\alpha) \overset{I}{S(s_k)} \overset{I}{s_l} \right)^c_d  f\!\left(\overset{J}{t_j} \overset{J}{t_k} \overset{J}{V}(\beta) \overset{J}{t_l} \overset{J}{S(t_i)}\right)^a _b = \overset{I}{T}{^c_d}\!\left( s_i s_j ? S(s_k) s_l \right)_{u_{\alpha}} \ast \overset{J}{T}{^a_b}\!\left( t_j t_k ? t_l S(t_i) \right)_{v_{\beta}}
\end{align*}
as desired. The fourth and firth equalities are \eqref{produitL01Explicite} and \eqref{produitGaugeFields10}.
\end{proof}

\begin{figure}[!h]
\centering
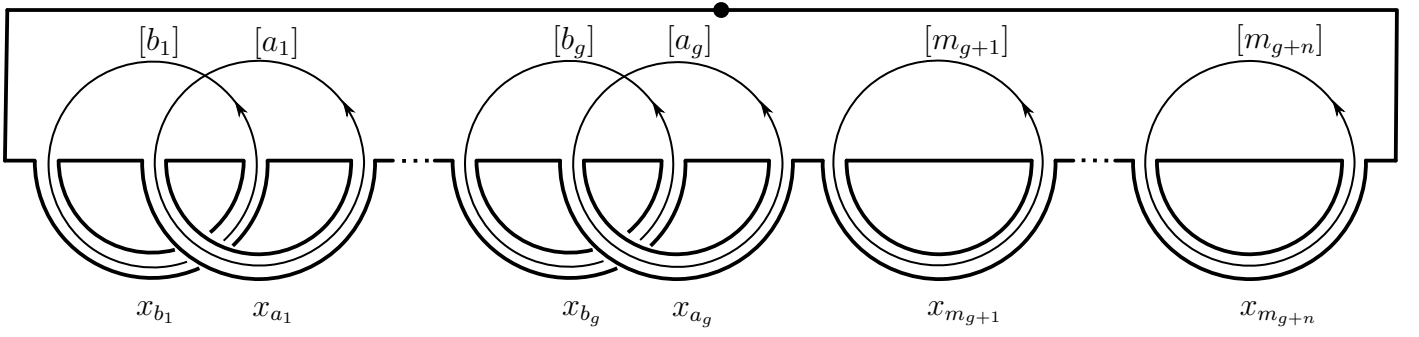
\caption{Surface $\Sigma_{g,n}^{\mathrm{o}}$ viewed as a thickening of $\Gamma_{g,n}$, canonical loops and a discrete connection.}
\label{surfaceGN}
\end{figure}

\indent Now, let $\Gamma = (V,E)$ be a filling graph of $\Sigma_{g,n}$ (an embedded oriented graph such that $\Sigma_{g,n} \setminus \Gamma$ is a union of open discs). Recall (see \cite{BFKB, BFKB2} for the precise definitions) that a lattice gauge field theory on $\Gamma$ consists of 
\begin{itemize}
\item[-] a space of (discrete) connections $\mathcal{A}_{\Gamma} = \bigotimes_{e \in E} H_e$,
\item[-] a space of gauge fields $\mathbb{C}[\mathcal{A}_{\Gamma}] = \bigotimes_{e \in E} H^*_e$ (functions on $\mathcal{A}_{\Gamma}$),
\item[-] a gauge algebra $\mathcal{G}_{\Gamma} = \bigotimes_{v \in V} H_v.$
\end{itemize}
In \cite{BFKB, BFKB2}, an action of $\mathcal{G}_{\Gamma}$ on $\mathcal{A}_{\Gamma}$ as well as a $H$-equivariant comultiplication \\$\nabla : \mathcal{A}_{\Gamma} \to \mathcal{A}_{\Gamma}^{\otimes 2}$ are defined. Dualizing this, we get a structure of (right) $H$-module-algebra (with product denoted by $\star$) on $\mathbb{C}[\mathcal{A}_{\Gamma}]$:
\[ \forall \, h \in H, \:\:\: \forall \, \varphi, \psi \in \mathbb{C}[\mathcal{A}_{\Gamma}], \:\:\:\: \varphi\cdot h = \varphi(h\cdot ?), \:\:\:\:\: \varphi \star \psi = (\varphi \otimes \psi) \circ \nabla.\]
Here we take the most natural graph
\[ \Gamma_{g,n} = \left(V = \{\bullet\}, \: E_{g,n} = \{b_1, a_1, \ldots, b_g, a_g, m_{g+1}, \ldots, m_{g+n}\} \right), \]
(see Figure \ref{surfaceGN}) and we denote $\mathcal{A}_{\Gamma_{g,n}} = \mathcal{A}_{g,n}$, $\mathcal{G}_{\Gamma_{g,n}} = H$. With this choice, a discrete connection is
\[ h_{b_1} \otimes h_{a_1} \otimes \ldots \otimes h_{b_g} \otimes h_{a_g} \otimes h_{m_{g+1}} \otimes \ldots \otimes h_{m_{g+n}} \in \mathcal{A}_{g,n}. \]
 We will consider the examples of $\Sigma_{0,1}$ and $\Sigma_{1,0}$, write down the $H$-action and the product in $\mathbb{C}[\mathcal{A}_{g,n}]$ and observe that it is isomorphic to $\mathcal{F}_{g,n}(H)$ (and hence to $\mathcal{L}_{g,n}(H)$). We use the rules given in \cite{BFKB} without explanation.

\smallskip

\indent For $\Sigma_{0,1}$, we have the graph at the left below (with ciliation and cyclic order, see \cite{BFKB} and the references therein):
\begin{center}
\begingroup%
  \makeatletter%
  \providecommand\color[2][]{%
    \errmessage{(Inkscape) Color is used for the text in Inkscape, but the package 'color.sty' is not loaded}%
    \renewcommand\color[2][]{}%
  }%
  \providecommand\transparent[1]{%
    \errmessage{(Inkscape) Transparency is used (non-zero) for the text in Inkscape, but the package 'transparent.sty' is not loaded}%
    \renewcommand\transparent[1]{}%
  }%
  \providecommand\rotatebox[2]{#2}%
  \newcommand*\fsize{\dimexpr\f@size pt\relax}%
  \newcommand*\lineheight[1]{\fontsize{\fsize}{#1\fsize}\selectfont}%
  \ifx\svgwidth\undefined%
    \setlength{\unitlength}{466.57606161bp}%
    \ifx\svgscale\undefined%
      \relax%
    \else%
      \setlength{\unitlength}{\unitlength * \real{\svgscale}}%
    \fi%
  \else%
    \setlength{\unitlength}{\svgwidth}%
  \fi%
  \global\let\svgwidth\undefined%
  \global\let\svgscale\undefined%
  \makeatother%
  \begin{picture}(1,0.17208037)%
    \lineheight{1}%
    \setlength\tabcolsep{0pt}%
    \put(0,0){\includegraphics[width=\unitlength,page=1]{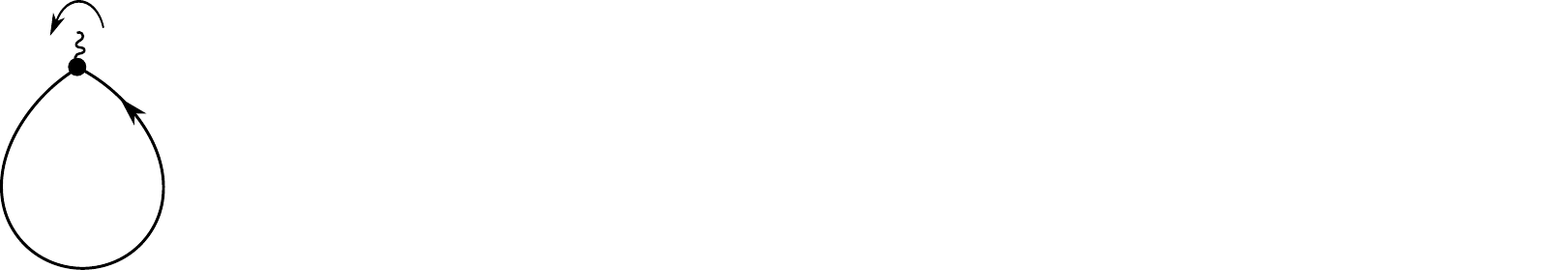}}%
    \put(0.10706903,0.02772659){\color[rgb]{0,0,0}\makebox(0,0)[lt]{\lineheight{1.25}\smash{\begin{tabular}[t]{l}$m$\end{tabular}}}}%
    \put(0,0){\includegraphics[width=\unitlength,page=2]{graphe01.pdf}}%
    \put(0.37793382,0.03639499){\color[rgb]{0,0,0}\makebox(0,0)[lt]{\lineheight{1.25}\smash{\begin{tabular}[t]{l}$x_m^{\text{out}}$\end{tabular}}}}%
    \put(0,0){\includegraphics[width=\unitlength,page=3]{graphe01.pdf}}%
    \put(0.50257185,0.03791977){\color[rgb]{0,0,0}\makebox(0,0)[lt]{\lineheight{1.25}\smash{\begin{tabular}[t]{l}$x_m^{\text{in}}$\end{tabular}}}}%
    \put(0.30192884,0.09313787){\color[rgb]{0,0,0}\makebox(0,0)[lt]{\lineheight{1.25}\smash{\begin{tabular}[t]{l}$\mathrm{v} \:\:\: = $\end{tabular}}}}%
    \put(0,0){\includegraphics[width=\unitlength,page=4]{graphe01.pdf}}%
    \put(0.78845785,0.03662345){\color[rgb]{0,0,0}\makebox(0,0)[lt]{\lineheight{1.25}\smash{\begin{tabular}[t]{l}$h'x_m^{\text{out}}$\end{tabular}}}}%
    \put(0,0){\includegraphics[width=\unitlength,page=5]{graphe01.pdf}}%
    \put(0.91242012,0.03814824){\color[rgb]{0,0,0}\makebox(0,0)[lt]{\lineheight{1.25}\smash{\begin{tabular}[t]{l}$x_m^{\text{in}}S(h'')$\end{tabular}}}}%
    \put(0.68039684,0.10107873){\color[rgb]{0,0,0}\makebox(0,0)[lt]{\lineheight{1.25}\smash{\begin{tabular}[t]{l}$h\cdot \mathrm{v} \:\:\: =$\end{tabular}}}}%
  \end{picture}%
\endgroup%

\end{center}
A connection is an assignment of an element of $H$ to each edge (the holonomy of that edge), so here a connection is simply $x_m$. To compute the action of $H$ on $x_m$, we determine the action on the vertex $\mathrm{v}$ as represented above. Then, gluing $x_m^{\text{out}} \sim x_m^{\text{in}}$, we get $h \cdot x_m = h' x_m S(h'')$. Hence the right action on $\varphi \in \mathbb{C}[\mathcal{A}_{0,1}]$ is $\varphi \cdot h = \varphi(h' ? S(h''))$. To compute $\nabla$, we consider the operator $F_{\mathrm{v}} : H^{\otimes 2} \to H^{\otimes 4}$ associated to the vertex $\mathrm{v}$:
\begin{center}
\begingroup%
  \makeatletter%
  \providecommand\color[2][]{%
    \errmessage{(Inkscape) Color is used for the text in Inkscape, but the package 'color.sty' is not loaded}%
    \renewcommand\color[2][]{}%
  }%
  \providecommand\transparent[1]{%
    \errmessage{(Inkscape) Transparency is used (non-zero) for the text in Inkscape, but the package 'transparent.sty' is not loaded}%
    \renewcommand\transparent[1]{}%
  }%
  \providecommand\rotatebox[2]{#2}%
  \newcommand*\fsize{\dimexpr\f@size pt\relax}%
  \newcommand*\lineheight[1]{\fontsize{\fsize}{#1\fsize}\selectfont}%
  \ifx\svgwidth\undefined%
    \setlength{\unitlength}{202.18042433bp}%
    \ifx\svgscale\undefined%
      \relax%
    \else%
      \setlength{\unitlength}{\unitlength * \real{\svgscale}}%
    \fi%
  \else%
    \setlength{\unitlength}{\svgwidth}%
  \fi%
  \global\let\svgwidth\undefined%
  \global\let\svgscale\undefined%
  \makeatother%
  \begin{picture}(1,0.42417658)%
    \lineheight{1}%
    \setlength\tabcolsep{0pt}%
    \put(0,0){\includegraphics[width=\unitlength,page=1]{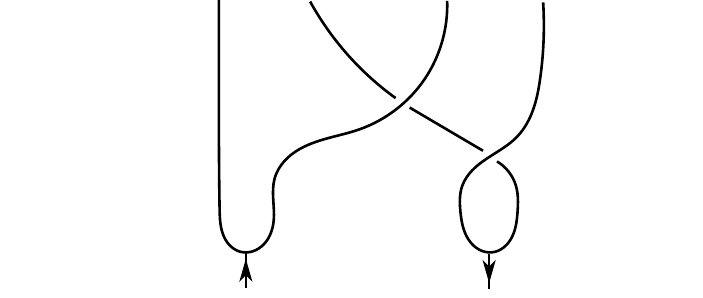}}%
    \put(0.36323657,0.00537845){\color[rgb]{0,0,0}\makebox(0,0)[lt]{\lineheight{1.25}\smash{\begin{tabular}[t]{l}$x_m^{\text{out}}$\end{tabular}}}}%
    \put(0.71152412,0.00692097){\color[rgb]{0,0,0}\makebox(0,0)[lt]{\lineheight{1.25}\smash{\begin{tabular}[t]{l}$x_m^{\text{in}}$\end{tabular}}}}%
    \put(0.149162,0.10505505){\color[rgb]{0,0,0}\makebox(0,0)[lt]{\lineheight{1.25}\smash{\begin{tabular}[t]{l}$(x_m^{\text{out}})'$\end{tabular}}}}%
    \put(0.40836038,0.15136552){\color[rgb]{0,0,0}\makebox(0,0)[lt]{\lineheight{1.25}\smash{\begin{tabular}[t]{l}$(x_m^{\text{out}})''$\end{tabular}}}}%
    \put(0.50477529,0.0818998){\color[rgb]{0,0,0}\makebox(0,0)[lt]{\lineheight{1.25}\smash{\begin{tabular}[t]{l}$(x_m^{\text{in}})''$\end{tabular}}}}%
    \put(0.75090979,0.10505505){\color[rgb]{0,0,0}\makebox(0,0)[lt]{\lineheight{1.25}\smash{\begin{tabular}[t]{l}$(x_m^{\text{in}})'$\end{tabular}}}}%
    \put(0,0){\includegraphics[width=\unitlength,page=2]{tangle01.pdf}}%
    \put(-0.00336664,0.26172376){\color[rgb]{0,0,0}\makebox(0,0)[lt]{\lineheight{1.25}\smash{\begin{tabular}[t]{l}$F_{\mathrm{v}} \:\:\: =$\end{tabular}}}}%
  \end{picture}%
\endgroup%

\end{center}
Evaluating the tangle\footnote{We point out a misprint in \cite[Fig. 11]{BFKB}: the second and third crossings are inverted.} yields $F_{\mathrm{v}}\!\left( x_m^{\text{out}} \otimes x_m^{\text{in}} \right) = (x_m^{\text{out}})' \otimes (x_m^{\text{in}})'b_iS(b_j) \otimes a_j (x_m^{\text{out}})'' \otimes (x_m^{\text{in}})'' a_i$. Then, gluing $x_m^{\text{out}} \sim x_m^{\text{in}}$, we get $\nabla(x_m) = x_m' b_i S(b_j) \otimes a_j x_m'' a_i$ and thus
\[ \varphi \star \psi = \varphi\!\left(? b_i S(b_j)\right) \psi\!\left( a_j ? a_i \right). \]
We see that $\mathbb{C}[\mathcal{A}_{0,1}] = \mathcal{F}_{0,1}(H)$ (see Remark \ref{F01}).

\smallskip

\indent For $\Sigma_{1,0}$, we have the graph at the left below:
\begin{center}
\begingroup%
  \makeatletter%
  \providecommand\color[2][]{%
    \errmessage{(Inkscape) Color is used for the text in Inkscape, but the package 'color.sty' is not loaded}%
    \renewcommand\color[2][]{}%
  }%
  \providecommand\transparent[1]{%
    \errmessage{(Inkscape) Transparency is used (non-zero) for the text in Inkscape, but the package 'transparent.sty' is not loaded}%
    \renewcommand\transparent[1]{}%
  }%
  \providecommand\rotatebox[2]{#2}%
  \newcommand*\fsize{\dimexpr\f@size pt\relax}%
  \newcommand*\lineheight[1]{\fontsize{\fsize}{#1\fsize}\selectfont}%
  \ifx\svgwidth\undefined%
    \setlength{\unitlength}{502.07366917bp}%
    \ifx\svgscale\undefined%
      \relax%
    \else%
      \setlength{\unitlength}{\unitlength * \real{\svgscale}}%
    \fi%
  \else%
    \setlength{\unitlength}{\svgwidth}%
  \fi%
  \global\let\svgwidth\undefined%
  \global\let\svgscale\undefined%
  \makeatother%
  \begin{picture}(1,0.17819056)%
    \lineheight{1}%
    \setlength\tabcolsep{0pt}%
    \put(0,0){\includegraphics[width=\unitlength,page=1]{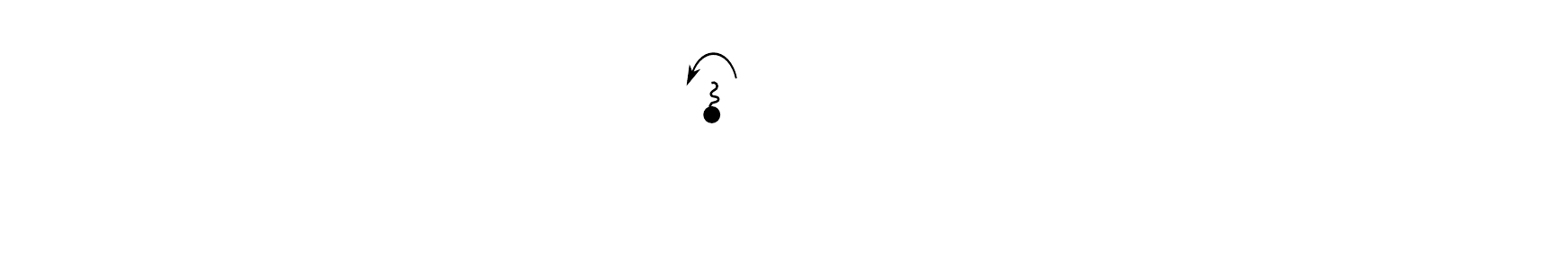}}%
    \put(0.40299427,0.02075827){\color[rgb]{0,0,0}\makebox(0,0)[lt]{\lineheight{1.25}\smash{\begin{tabular}[t]{l}$x_a^{\text{out}}$\end{tabular}}}}%
    \put(0,0){\includegraphics[width=\unitlength,page=2]{graphe10.pdf}}%
    \put(0.51882014,0.02217525){\color[rgb]{0,0,0}\makebox(0,0)[lt]{\lineheight{1.25}\smash{\begin{tabular}[t]{l}$x_b^{\text{in}}$\end{tabular}}}}%
    \put(0.28593706,0.07487458){\color[rgb]{0,0,0}\makebox(0,0)[lt]{\lineheight{1.25}\smash{\begin{tabular}[t]{l}$\mathrm{v} \:\:\: = $\end{tabular}}}}%
    \put(0,0){\includegraphics[width=\unitlength,page=3]{graphe10.pdf}}%
    \put(0.1899408,0.06206373){\color[rgb]{0,0,0}\makebox(0,0)[lt]{\lineheight{1.25}\smash{\begin{tabular}[t]{l}$a$\end{tabular}}}}%
    \put(0,0){\includegraphics[width=\unitlength,page=4]{graphe10.pdf}}%
    \put(-0.00117013,0.06495378){\color[rgb]{0,0,0}\makebox(0,0)[lt]{\lineheight{1.25}\smash{\begin{tabular}[t]{l}$b$\end{tabular}}}}%
    \put(0,0){\includegraphics[width=\unitlength,page=5]{graphe10.pdf}}%
    \put(0.34363077,0.13489043){\color[rgb]{0,0,0}\makebox(0,0)[lt]{\lineheight{1.25}\smash{\begin{tabular}[t]{l}$x_b^{\text{out}}$\end{tabular}}}}%
    \put(0.54166326,0.1344464){\color[rgb]{0,0,0}\makebox(0,0)[lt]{\lineheight{1.25}\smash{\begin{tabular}[t]{l}$x_a^{\text{in}}$\end{tabular}}}}%
    \put(0,0){\includegraphics[width=\unitlength,page=6]{graphe10.pdf}}%
    \put(0.76285791,0.01906371){\color[rgb]{0,0,0}\makebox(0,0)[lt]{\lineheight{1.25}\smash{\begin{tabular}[t]{l}$h''x_a^{\text{out}}$\end{tabular}}}}%
    \put(0,0){\includegraphics[width=\unitlength,page=7]{graphe10.pdf}}%
    \put(0.87868383,0.02048072){\color[rgb]{0,0,0}\makebox(0,0)[lt]{\lineheight{1.25}\smash{\begin{tabular}[t]{l}$x_b^{\text{in}}S(h''')$\end{tabular}}}}%
    \put(0.63324205,0.07255209){\color[rgb]{0,0,0}\makebox(0,0)[lt]{\lineheight{1.25}\smash{\begin{tabular}[t]{l}$h \cdot \mathrm{v} \:\:\: = $\end{tabular}}}}%
    \put(0,0){\includegraphics[width=\unitlength,page=8]{graphe10.pdf}}%
    \put(0.70349446,0.13319588){\color[rgb]{0,0,0}\makebox(0,0)[lt]{\lineheight{1.25}\smash{\begin{tabular}[t]{l}$h'x_b^{\text{out}}$\end{tabular}}}}%
    \put(0.9015269,0.13275185){\color[rgb]{0,0,0}\makebox(0,0)[lt]{\lineheight{1.25}\smash{\begin{tabular}[t]{l}$x_a^{\text{in}}S(h^{(4)})$\end{tabular}}}}%
  \end{picture}%
\endgroup%

\end{center}
Here a connection is $x_b \otimes x_a \in H^{\otimes 2}$. Gluing $x_b^{\text{out}} \sim x_b^{\text{in}}$ and $x_a^{\text{out}} \sim x_a^{\text{in}}$ in the action on the vertex, we get that the left action of $H$ on $\mathcal{A}_{1,0}$ is $h \cdot (x_b \otimes x_a) = h' x_b S(h''') \otimes h'' x_a S(h^{(4)})$. Hence the right action of $H$ on $\varphi \otimes \psi \in \mathbb{C}[\mathcal{A}_{1,0}]$ is 
\[ (\varphi \otimes \psi) \cdot h = \varphi\!\left( h' ? S(h''') \right) \otimes \psi\!\left( h'' ? S(h^{(4)}) \right) \]
To compute $\nabla$, we consider the operator $F_{\mathrm{v}} : H^{\otimes 4} \to H^{\otimes 8}$ associated to the vertex $\mathrm{v}$:
\begin{center}
\begingroup%
  \makeatletter%
  \providecommand\color[2][]{%
    \errmessage{(Inkscape) Color is used for the text in Inkscape, but the package 'color.sty' is not loaded}%
    \renewcommand\color[2][]{}%
  }%
  \providecommand\transparent[1]{%
    \errmessage{(Inkscape) Transparency is used (non-zero) for the text in Inkscape, but the package 'transparent.sty' is not loaded}%
    \renewcommand\transparent[1]{}%
  }%
  \providecommand\rotatebox[2]{#2}%
  \newcommand*\fsize{\dimexpr\f@size pt\relax}%
  \newcommand*\lineheight[1]{\fontsize{\fsize}{#1\fsize}\selectfont}%
  \ifx\svgwidth\undefined%
    \setlength{\unitlength}{349.00606149bp}%
    \ifx\svgscale\undefined%
      \relax%
    \else%
      \setlength{\unitlength}{\unitlength * \real{\svgscale}}%
    \fi%
  \else%
    \setlength{\unitlength}{\svgwidth}%
  \fi%
  \global\let\svgwidth\undefined%
  \global\let\svgscale\undefined%
  \makeatother%
  \begin{picture}(1,0.33846952)%
    \lineheight{1}%
    \setlength\tabcolsep{0pt}%
    \put(0,0){\includegraphics[width=\unitlength,page=1]{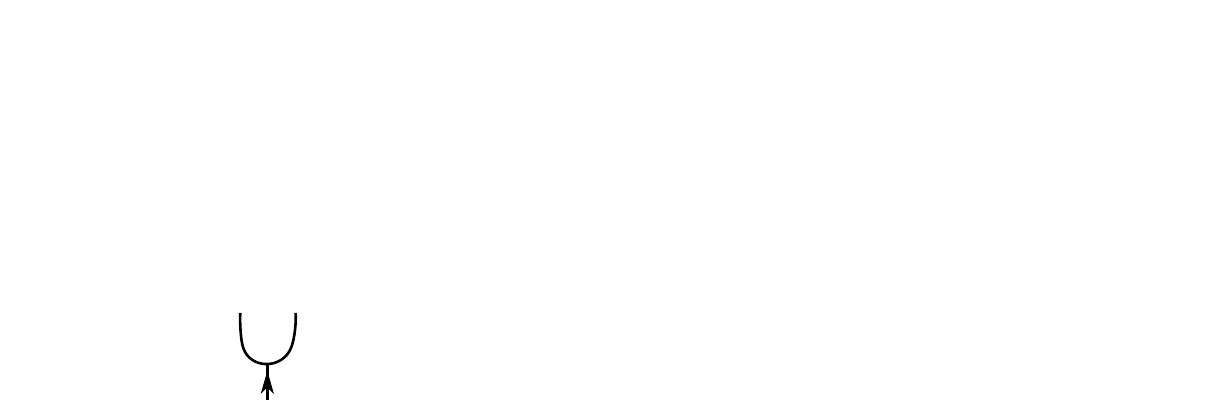}}%
    \put(0.22974126,0.00311582){\color[rgb]{0,0,0}\makebox(0,0)[lt]{\lineheight{1.25}\smash{\begin{tabular}[t]{l}$x_b^{\text{out}}$\end{tabular}}}}%
    \put(0,0){\includegraphics[width=\unitlength,page=2]{tangle10.pdf}}%
    \put(0.43389247,0.00311582){\color[rgb]{0,0,0}\makebox(0,0)[lt]{\lineheight{1.25}\smash{\begin{tabular}[t]{l}$x_a^{\text{out}}$\end{tabular}}}}%
    \put(0,0){\includegraphics[width=\unitlength,page=3]{tangle10.pdf}}%
    \put(0.63804369,0.00311576){\color[rgb]{0,0,0}\makebox(0,0)[lt]{\lineheight{1.25}\smash{\begin{tabular}[t]{l}$x_b^{\text{in}}$\end{tabular}}}}%
    \put(0,0){\includegraphics[width=\unitlength,page=4]{tangle10.pdf}}%
    \put(0.84219485,0.00311576){\color[rgb]{0,0,0}\makebox(0,0)[lt]{\lineheight{1.25}\smash{\begin{tabular}[t]{l}$x_a^{\text{in}}$\end{tabular}}}}%
    \put(0,0){\includegraphics[width=\unitlength,page=5]{tangle10.pdf}}%
    \put(0.10081922,0.05096964){\color[rgb]{0,0,0}\makebox(0,0)[lt]{\lineheight{1.25}\smash{\begin{tabular}[t]{l}$(x_b^{\text{out}})'$\end{tabular}}}}%
    \put(0.25258024,0.08936157){\color[rgb]{0,0,0}\makebox(0,0)[lt]{\lineheight{1.25}\smash{\begin{tabular}[t]{l}$(x_b^{\text{out}})''$\end{tabular}}}}%
    \put(0.30542561,0.04961462){\color[rgb]{0,0,0}\makebox(0,0)[lt]{\lineheight{1.25}\smash{\begin{tabular}[t]{l}$(x_a^{\text{out}})'$\end{tabular}}}}%
    \put(0.46847836,0.08123152){\color[rgb]{0,0,0}\makebox(0,0)[lt]{\lineheight{1.25}\smash{\begin{tabular}[t]{l}$(x_a^{\text{out}})''$\end{tabular}}}}%
    \put(0.5119547,0.04916295){\color[rgb]{0,0,0}\makebox(0,0)[lt]{\lineheight{1.25}\smash{\begin{tabular}[t]{l}$(x_b^{\text{in}})''$\end{tabular}}}}%
    \put(0.65843644,0.08123152){\color[rgb]{0,0,0}\makebox(0,0)[lt]{\lineheight{1.25}\smash{\begin{tabular}[t]{l}$(x_b^{\text{in}})'$\end{tabular}}}}%
    \put(0.72347688,0.04961462){\color[rgb]{0,0,0}\makebox(0,0)[lt]{\lineheight{1.25}\smash{\begin{tabular}[t]{l}$(x_a^{\text{in}})''$\end{tabular}}}}%
    \put(0.8604748,0.04961462){\color[rgb]{0,0,0}\makebox(0,0)[lt]{\lineheight{1.25}\smash{\begin{tabular}[t]{l}$(x_a^{\text{in}})'$\end{tabular}}}}%
    \put(-0.00232153,0.1838989){\color[rgb]{0,0,0}\makebox(0,0)[lt]{\lineheight{1.25}\smash{\begin{tabular}[t]{l}$F_{\mathrm{v}} \:\:\: = $\end{tabular}}}}%
  \end{picture}%
\endgroup%

\end{center}
Evaluating the tangle yields
\begin{align*}
F_{\mathrm{v}}\!\left( x_b^{\text{out}} \otimes x_a^{\text{out}} \otimes x_b^{\text{in}} \otimes x_a^{\text{in}} \right) =& \: (x_b^{\text{out}})' \otimes b_k(x_a^{\text{out}})' \otimes (x_b^{\text{in}})' b_i S(b_l) S(b_n) \otimes (x_a^{\text{in}})' b_j b_m S(b_o) S(b_p)\\
& \otimes a_p a_n a_k (x_b^{\text{out}})'' \otimes a_o a_l (x_a^{\text{out}})'' \otimes (x_b^{\text{in}})'' a_i a_m \otimes (x_a^{\text{in}})'' a_j. 
\end{align*}
Then, gluing $x_b^{\text{out}} \sim x_b^{\text{in}}$ and $x_a^{\text{out}} \sim x_a^{\text{in}}$, we get
\[ \nabla(x_b \otimes x_a) = x_b' b_i S(b_l) S(b_n) \otimes b_k x_a' b_j b_m S(b_o) S(b_p) \otimes a_p a_n a_k x_b'' a_i a_m \otimes a_o a_l x_a'' a_j \]
and thus
\[ (\varphi_1 \otimes \psi_1) \star (\varphi_2 \otimes \psi_2) = \varphi_1\!\left( ? b_i S(b_l) S(b_n) \right) \varphi_2\!\left( a_p a_n a_k ? a_i a_m \right) \otimes \psi_1\!\left( b_k ? b_j b_m S(b_o) S(b_p) \right) \psi_2\!\left( a_o a_l ? a_j \right) \]
Recall that we defined $\varphi_b = \varphi \otimes \varepsilon$ and $\varphi_a = \varepsilon \otimes \varphi$. We have
\[ \varphi_b \star \psi_a = \varphi\!\left( ? S(b_i) \right) \otimes \psi\!\left( a_i ? \right) = \varphi\!\left( ? S(b_i) \right)_b \ast \psi\!\left( a_i ? \right)_a \]
and this suggests to consider
\[ \fonc{\eta}{\mathbb{C}[\mathcal{A}_{1,0}]}{\mathcal{F}_{1,0}(H)}{\varphi \otimes \psi}{\varphi\!\left( ? b_i \right) \otimes \psi\!\left( a_i ? \right)}. \]
A straightforward computation using \eqref{deltaR} and \eqref{YangBaxter} shows that $\eta$ is indeed an isomorphism of $H$-modules-algebras (see Remark \ref{remarkProductGaugeFields}). The $H$-module-algebras $\mathbb{C}[\mathcal{A}_{g,n}]$ and $\mathcal{F}_{g,n}(H)$ are isomorphic for any $g,n$. The example of $\Sigma_{1,0}^{\mathrm{o}}$ is generalized as follows: for each crossing in the graph $\Gamma_{g,n}$ use an $R$-matrix as above to define the isomorphism.

\section{Representation of $\mathcal{L}^{\mathrm{inv}}_{g,n}(H)$}\label{RepInvariants}
\indent An element $x \in \mathcal{L}_{g,n}(H)$ is called {\em invariant} if $x \cdot h = \varepsilon(h)x$ for all $h \in H$, or equivalently, if $\Omega(x) = \varepsilon \otimes x$. Such elements are also called ``observables''. In this section we construct representations of the subalgebra of invariant elements $\mathcal{L}^{\mathrm{inv}}_{g,n}(H)$. 

\smallskip

Recall from section \ref{sectionRepInvariants} that the matrices 
\[ \overset{I}{C} = \overset{I}{v}{^{2}}\overset{I}{B}\overset{I}{A}{^{-1}}\overset{I}{B}{^{-1}}\overset{I}{A}, \:\:\:\: \overset{I}{C}{^{(\pm)}} = \Psi_{1,0}^{-1}(\overset{I}{L}{^{(\pm)}}\overset{I}{\widetilde{L}}{^{(\pm)}}) \:\:\: \in \mathrm{Mat}_{\dim(I)}\bigl( \mathcal{L}_{1,0}(H) \bigr) \]
satisfy the decomposition $\overset{I}{C} = \overset{I}{C}{^{(+)}} \overset{I}{C}{^{(-)-1}}$ and allow for a simple characterization of the invariant elements. We generalize this to any $g,n$. For $i \leq g$, let $\overset{I}{C}(i)$ be the embedding of $\overset{I}{C}$ previously defined on the $i$-th copy of $\mathcal{L}_{1,0}(H)$ in $\mathcal{L}_{g,n}(H)$: $\overset{I}{C}(i) = \overset{I}{v}{^{2}}\overset{I}{B}(i) \overset{I}{A}{^{-1}}(i) \overset{I}{B}{^{-1}}(i) \overset{I}{A}(i)$.
\begin{definition} 
$\overset{I}{C}_{g,n} = \overset{I}{C}(1) \ldots \overset{I}{C}(g) \overset{I}{M}(g+1) \ldots \overset{I}{M}(g+n) \in \mathrm{Mat}_{\dim(I)}\!\left(\mathcal{L}_{g,n}(H)\right)$.
\end{definition}
\noindent In particular $\overset{I}{C}_{1,0} = \overset{I}{C}$. Geometrically (see \eqref{boundaryCirclePi1} and Figure \ref{surfaceAvecMatrices}), for each $I$ the matrix $\overset{I}{C}_{g,n}$ corresponds to the boundary circle induced by the removal of the disk $D$ in $\Sigma_{g,n}$.
\smallskip\\
\indent There is a decomposition analogous to Lemma \ref{decGauss}. Indeed, let
\[ \overset{I}{C}{_{g,n}^{(\pm)}} = \alpha_{g,n}^{-1}\!\left(\overset{I}{\underline{C}}{^{(\pm)}}(1) \ldots \overset{I}{\underline{C}}{^{(\pm)}}(g) \overset{I}{\underline{M}}{^{(\pm)}}(g+1) \ldots \overset{I}{\underline{M}}{^{(\pm)}}(g+n)\right) \in \mathrm{Mat}_{\dim(I)}\!\left(\mathcal{L}_{g,n}(H)\right), \]
where $\overset{I}{C}{^{(\pm)}} = \Psi_{1,0}^{-1}(\overset{I}{L}{^{(\pm)}}\overset{I}{\widetilde{L}}{^{(\pm)}})$ and $\overset{I}{M}{^{(\pm)}} = \Psi_{0,1}^{-1}(\overset{I}{L}{^{(\pm)}})$ (recall \eqref{L} and \eqref{defLTilde}).
\begin{proposition}\label{propMatriceCgn}
The following equality holds in $\mathcal{L}_{g,n}(H)$:
\[ \overset{I}{C}_{g,n} = \overset{I}{C}{^{(+)}_{g,n}}\,\overset{I}{C}{^{(-)-1}_{g,n}}. \]
Moreover, the matrices $\overset{I}{C}_{g,n}$ satisfy the fusion relation of $\mathcal{L}_{0,1}(H)$:
\[ (\overset{I \otimes J}{C}\!\!_{g,n})_{12} = (\overset{I}{C}_{g,n})_1 \, \overset{IJ}{(R')}_{12} \, (\overset{J}{C}_{g,n})_2 \, \overset{IJ}{(R'^{-1})}_{12}. \]
\end{proposition}
\begin{proof}
The first claim is a simple consequence of the definition of $\alpha_{g,n}$ and of Lemma \ref{decGauss}. The fusion relation is a consequence of a more general fact which is easy to show, namely: if $i_1 < \ldots < i_k$ and if $\overset{I}{X^1}(i_1), \ldots, \overset{I}{X^k}(i_k)$ are matrices satisfying the fusion relation of $\mathcal{L}_{0,1}(H)$, then their product $\overset{I}{X^1}(i_1) \ldots \overset{I}{X^k}(i_k)$ also satisfies the fusion relation of $\mathcal{L}_{0,1}(H)$.
\end{proof}

\indent The image of these matrices have simple expressions in $\mathcal{H}(\mathcal{O}(H))^{\otimes g} \otimes H^{\otimes n}$:
\begin{lemma}\label{expressionM}
It holds
\begin{align*}
\Psi_{g,n}( \overset{I}{C}{_{g,n}^{(+)}} ) &= \overset{I}{a_i}\:\widetilde{b_i^{(2g -1 + n)}} b_i^{(2g + n)} \otimes \ldots \otimes \widetilde{b_i^{(1+n)}}b_i^{(2+n)} \otimes b_i^{(n)} \otimes \ldots \otimes b^{(1)}_i\\
\Psi_{g,n}( \overset{I}{C}{_{g,n}^{(-)}} ) &= \overset{I}{S^{-1}(b_i)} \: \widetilde{a_i^{(2g -1 + n)}} a_i^{(2g + n)} \otimes \ldots \otimes \widetilde{a_i^{(1+n)}}a_i^{(2+n)} \otimes a_i^{(n)} \otimes \ldots \otimes a^{(1)}_i\\
\Psi_{g,n}( \overset{I}{C}_{g,n} ) &= \overset{I}{X_i}\: \widetilde{Y_i^{(2g -1 + n)}} Y_i^{(2g + n)} \otimes \ldots \otimes \widetilde{Y_i^{(1+n)}}Y_i^{(2+n)} \otimes Y_i^{(n)} \otimes \ldots \otimes Y^{(1)}_i
\end{align*}
where $X_i \otimes Y_i = RR'$ and the superscripts mean iterated coproduct.
\end{lemma}
\begin{proof}
As an immediate consequence of quasitriangularity, we have for all $n \geq 2$ 
\[ (\mathrm{id} \otimes \Delta^{(n-1)})(R) = a_i \otimes b^{(1)}_i \otimes \ldots \otimes b^{(n)}_i = a_{i_1} \ldots a_{i_n} \otimes b_{i_n} \otimes \ldots \otimes b_{i_1}. \]
with implicit summation on $i_1, \ldots, i_n $. It follows that
\begin{align*}
\Psi_{g,n}( \overset{I}{C}{_{g,n}^{(+)}} ) &= \overset{I}{L}{^{(+)}}(1)\overset{I}{\widetilde{L}}{^{(+)}}(1) \ldots \overset{I}{L}{^{(+)}}(g)\overset{I}{\widetilde{L}}{^{(+)}}(g) \overset{I}{L}{^{(+)}}(g+1) \ldots \overset{I}{L}{^{(+)}}(g+n)\\
&= \overset{I}{a_{i_1}} \ldots \overset{I}{a_{i_{2g+n}}} \: \widetilde{b_{i_2}} b_{i_1} \otimes \ldots \otimes \widetilde{b_{i_{2g}}} b_{i_{2g-1}} \otimes b_{i_{2g+1}} \otimes \ldots \otimes b_{i_{2g+n}}\\
&= \overset{I}{a_i}\:\widetilde{b_i^{(2g -1 + n)}} b_i^{(2g + n)} \otimes \ldots \otimes \widetilde{b_i^{(1+n)}}b_i^{(2+n)} \otimes b_i^{(n)} \otimes \ldots \otimes b^{(1)}_i
\end{align*}
as desired. The second is shown similarly since $R'^{-1}$ is also an universal $R$-matrix. The third is an immediate consequence.
\end{proof}

\noindent Lemma \ref{expressionM} indicates that the algebra generated by the coefficients $\overset{I}{C^{(\pm)}}{^i_j}$ equals the vector space generated by the coefficients $\overset{I}{C}{^i_j}$:
\begin{equation}\label{CoeffsCCPlusCMoins}
\mathbb{C}\langle \overset{I}{C^{(\pm)}}{^i_j} \rangle_{I,i,j} = \mathrm{vect}(\overset{I}{C}{^i_j})_{I,i,j}. 
\end{equation}
\indent The matrices $\overset{I}{C}_{g,n}$ satisfying the fusion relation of $\mathcal{L}_{0,1}(H)$, we can apply Lemma \ref{injectionFusion} and define a representation of $H$ on $V = (H^*)^{\otimes g} \otimes I_1 \otimes \ldots \otimes I_n$ by
\begin{equation}\label{actionHsurFormes}
h \cdot v = h_{C_{g,n}} \triangleright v.
\end{equation}
Recall that
\[ \begin{array}{ccccc}
H & \overset{\Psi_{0,1}^{-1}}{\longrightarrow} & \mathcal{L}_{0,1}(H) & \overset{j_{C_{g,n}}}{\longrightarrow} & \mathcal{L}_{g,n}(H)\\
\bigl(\overset{I}{X}_i\bigr) Y_i & \longmapsto & \overset{I}{M} & \longmapsto & \overset{I}{C}_{g,n}
\end{array}
\]
Since $H$ is factorizable, each $h \in H$ is a linear combination of coefficients of the matrices $\overset{I}{L}{^{(+)}}\overset{I}{L}{^{(-)-1}} = \bigl(\overset{I}{X}_i\bigr) Y_i$ and thus $h_{C_{g,n}} = j_{C_{g,n}} \circ \Psi_{0,1}^{-1}(h)$ is a linear combination of coefficients of the matrices $\overset{I}{C}_{g,n}$. It follows from Lemma \ref{expressionM} that the action \eqref{actionHsurFormes} is explicitly given by
\begin{equation}\label{actionH}
\begin{split}
&h \cdot \varphi_1 \, \otimes \, \ldots \, \otimes \, \varphi_g \, \otimes v_1 \otimes \, \ldots \, \otimes v_n\\
&= \varphi_1\!\left(S^{-1}\!\left(h^{(2g -1 + n)}\right) ? h^{(2g + n)}\right)\,  \otimes \ldots \otimes \, \varphi_g\!\left(S^{-1}\!\left(h^{(1+n)}\right) ? h^{(2+n)}\right)
\otimes \,h^{(n)}v_1 \otimes \ldots \otimes \, h^{(1)}v_n. 
\end{split}
\end{equation}

\indent As in the case of $\mathcal{L}_{1,0}(H)$, the matrices $\overset{I}{C}_{g,n}$ allow one to give a simple characterization of the invariant elements of $\mathcal{L}_{g,n}(H)$ and to construct representations of them. We begin with a technical lemma.
\begin{lemma}\label{conjugaisonM}
It holds
\[ (\overset{I}{C}{^{(\pm)}_{g,n}})_1 \, \overset{J}{U}(i)_2 \, (\overset{I}{C}{^{(\pm)}_{g,n}})^{-1}_1 \:\:= \overset{IJ}{R}{^{(\pm)-1}_{12}} \, \overset{J}{U}(i)_2 \, \overset{IJ}{R}{^{(\pm)}_{12}} \]
where $U$ is $A$, $B$ or $M$.
\end{lemma}
\begin{proof}
The case $(g,n) = (1,0)$ is Lemma \ref{conjugaisonC10}. The case $(g, n) = (0,1)$ is easy with \eqref{propertiesL}. Similarly, thanks to \eqref{propertiesL}, \eqref{echangeHeisenberg} and \eqref{LTilde}, we obtain
\[ \overset{IJ}{R}{^{(\pm)}_{12}} \, \overset{I}{C}{^{(\pm)}_{1}} \, \overset{J}{C}{^{(-)}_{2}} = \overset{J}{C}{^{(-)}_{2}} \, \overset{I}{C}{^{(\pm)}_{1}} \, \overset{IJ}{R}{^{(\pm)}_{12}}, 
\:\:\:\: \overset{IJ}{R}{^{(\pm)}_{12}} \, \overset{I}{M}{^{(\pm)}_{1}} \, \overset{J}{M}{^{(-)}_{2}} = \overset{J}{M}{^{(-)}_{2}} \, \overset{I}{M}{^{(\pm)}_{1}} \, \overset{IJ}{R}{^{(\pm)}_{12}}. \]
Using these preliminary facts, we can carry out the general computation. For instance, for $i \leq g$ 
{\footnotesize\begin{align*}
&\alpha_{g,n}\!\left( (\overset{I}{C}{^{(\pm)}_{g,n}})_1 \, \overset{J}{U}(i)_2 \, (\overset{I}{C}{^{(\pm)}_{g,n}})^{-1}_1\:\, \right) \\
& = \overset{I}{\underline{C}}{^{(\pm)}}(1)_1 \, \ldots \, \overset{I}{\underline{C}}{^{(\pm)}}(i)_1 \, \overset{J}{\underline{C}}{^{(-)}}(1)_2 \, \ldots \, \overset{J}{\underline{C}}{^{(-)}}(i-1)_2 \,
\overset{J}{\underline{U}}(i)_2 \, \overset{J}{\underline{C}}{^{(-)}}(i-1)^{-1}_2 \, 
 \ldots \, \overset{J}{\underline{C}}{^{(-)}}(1)^{-1}_2 \, \overset{I}{\underline{C}}{^{(\pm)}}(i)^{-1}_1 \, \ldots \, \overset{I}{\underline{C}}{^{(\pm)}}(1)^{-1}_1\\
& = \overset{I}{\underline{C}}{^{(\pm)}}(1)_1 \, \overset{J}{\underline{C}}{^{(-)}}(1)_2 \, \ldots \, \overset{I}{\underline{C}}{^{(\pm)}}(i-1)_1 \, \overset{J}{\underline{C}}{^{(-)}}(i-1)_2 \, 
\overset{I}{\underline{C}}{^{(\pm)}}(i)_1 \, \overset{J}{\underline{U}}(i)_2 \, 
 \overset{I}{\underline{C}}{^{(\pm)}}(i)^{-1}_1 \, \overset{J}{\underline{C}}{^{(-)}}(i-1)^{-1}_2 \, \overset{I}{\underline{C}}{^{(\pm)}}(i-1)^{-1}_1 \, \ldots\\
& \:\:\:\:\:\: \overset{J}{\underline{C}}{^{(-)}}(1)^{-1}_2 \, \overset{I}{\underline{C}}{^{(\pm)}}(1)^{-1}_1\\
& = \overset{I}{\underline{C}}{^{(\pm)}}(1)_1 \, \overset{J}{\underline{C}}{^{(-)}}(1)_2 \, \ldots \overset{I}{\underline{C}}{^{(\pm)}}(i-1)_1 \, \overset{J}{\underline{C}}{^{(-)}}(i-1)_2 \, 
\overset{IJ}{R}{^{(\pm)-1}_{12}} \, \overset{J}{\underline{U}}(i)_2 \, 
\overset{IJ}{R}{^{(\pm)}_{12}} \, \overset{J}{\underline{C}}{^{(-)}}(i-1)^{-1}_2 \, \overset{I}{\underline{C}}{^{(\pm)}}(i-1)^{-1}_1 \, \ldots \\
& \:\:\:\:\:\: \overset{J}{\underline{C}}{^{(-)}}(1)^{-1}_2 \, \overset{I}{\underline{C}}{^{(\pm)}}(1)^{-1}_1\\ 
& =\overset{IJ}{R}{^{(\pm)-1}_{12}}  \, \overset{J}{\underline{C}}{^{(-)}}(1)_2 \, \overset{I}{\underline{C}}{^{(\pm)}}(1)_1 \, \ldots \overset{J}{\underline{C}}{^{(-)}}(i-1)_2 \, 
\overset{I}{\underline{C}}{^{(\pm)}}(i-1)_1 \, \overset{J}{\underline{U}}(i)_2 \, 
 \overset{I}{\underline{C}}{^{(\pm)}}(i-1)^{-1}_1 \, \overset{J}{\underline{C}}{^{(-)}}(i-1)^{-1}_2 \, \ldots \\
& \:\:\:\:\:\: \overset{I}{\underline{C}}{^{(\pm)}}(1)^{-1}_1 \, \overset{J}{\underline{C}}{^{(-)}}(1)^{-1}_2 \, \overset{IJ}{R}{^{(\pm)}_{12}}\\ 
& =\overset{IJ}{R}{^{(\pm)-1}_{12}} \, \overset{J}{\underline{C}}{^{(-)}}(1)_2 \, \ldots \overset{J}{\underline{C}}{^{(-)}}(i-1)_2 \,
 \overset{J}{\underline{U}}(i)_2 \,  \overset{J}{\underline{C}}{^{(-)}}(i-1)^{-1}_2 \, \ldots \overset{J}{\underline{C}}{^{(-)}}(1)^{-1}_2 \, \overset{IJ}{R}{^{(\pm)}_{12}}
= \alpha_{g,n}\!\left(  \overset{IJ}{R}{^{(\pm)-1}_{12}} \, \overset{J}{U}(i)_2 \, \overset{IJ}{R}{^{(\pm)}_{12}} \right).
\end{align*}}
The case $i > g$ is treated in a similar way.
\end{proof}

\indent For $(V, \triangleright)$ a representation of $\mathcal{L}_{g,n}(H)$, let
\begin{equation}\label{defInvRep}
\mathrm{Inv}(V) = \left\{ v \in V \, \left| \, \forall \, I, \: \overset{I}{C}_{g,n} \triangleright v = \mathbb{I}_{\dim(I)}v \right.\right\} = \left\{ v \in V \, \left| \, \forall \, h \in H, \:\: h \cdot v = \varepsilon(h)v \right.\right\} 
\end{equation}
where $\mathbb{I}_k$ is the identity matrix of size $k$, and the action $\cdot$ of $H$ on $V$ is defined in \eqref{actionHsurFormes} and \eqref{actionH}. This subspace $\mathrm{Inv}(V)$ implements the flatness constraint \eqref{flatnessConstraintMatrices} discussed in the Introduction.

\begin{remark}
For $(g,n) = (1,0)$, $\mathrm{Inv}(H^*) = \mathrm{SLF}(H)$ thanks to Lemma \ref{CSLF}.
\finEx
\end{remark}

\begin{theorem}\label{thmInv}
1) An element $x \in \mathcal{L}_{g,n}(H)$ is invariant under the action of $H$ (or equivalently under the coaction $\Omega$ of $\mathcal{O}(H)$) if, and only if, for every $H$-module $I$, $\overset{I}{C}_{g,n}x = x\overset{I}{C}_{g,n}$. \\
2) Let $V$ be a representation of $\mathcal{L}_{g,n}(H)$. Then $\mathrm{Inv}(V)$ is stable under the action of invariant elements and thus provides a representation of $\mathcal{L}^{\mathrm{inv}}_{g,n}(H)$.
\end{theorem}
\begin{proof}
1) Letting $U$ be $A(i), B(i)$ or $M(j)$, $R^{(\pm)} = a_i^{(\pm)} \otimes b_i^{(\pm)}$ and using Lemma \ref{conjugaisonM}, we have that the right action $ \cdot$ of $H$ on $\mathcal{L}_{1,0}(H)$ satisfies:
\begin{align*}
\overset{J}{U}_2 \cdot  \overset{I}{L}{^{(\pm)-1}_1} &= \overset{J}{U}_2 \cdot S^{-1}(b_i^{(\pm)}) \overset{I}{(a_i^{(\pm)})}_1 = \overset{J}{S^{-1}(b_i^{(\pm)}{''})}_2 \overset{J}{U}_2 \overset{J}{(b_i^{(\pm)}{'})}_2 \overset{I}{(a_i^{(\pm)})}_1\\
& = \overset{J}{S^{-1}(b_i^{(\pm)})}_2 \overset{J}{U}_2 \overset{J}{(b_j^{(\pm)})}_2 \overset{I}{(a_i^{(\pm)}a_j^{(\pm)})}_1 = \overset{IJ}{R}{^{(\pm)-1}_{12}} \, \overset{J}{U}_2 \, \overset{IJ}{R}{^{(\pm)}_{12}} = (\overset{I}{C}{^{(\pm)}_{g,n}})_1 \, \overset{J}{U}_2 \, (\overset{I}{C}{^{(\pm)-1}_{g,n}})_1.
\end{align*}
We have thus shown that
\[ (\overset{J}{U})^c_d \cdot  S^{-1}(\overset{I}{L}{^{(\pm)}})^a_b = (\overset{I}{C}{^{(\pm)}_{g,n}})^a_i \, (\overset{J}{U})^c_d \, (\overset{I}{C}{^{(\pm)-1}_{g,n}})^i_b \]
or in other words
\[ \forall \,x \in \mathcal{L}_{g,n}(H), \:\:\: x \cdot  S^{-1}(\overset{I}{L}{^{(\pm)}}) = \overset{I}{C}{^{(\pm)}_{g,n}} \, x \, \overset{I}{C}{^{(\pm)-1}_{g,n}} \]
Since $H$ is factorizable, the elements $S^{-1}(\overset{I}{L}{^{(\pm)}})^a_b$ generate $H$ as an algebra. Hence the previous equation shows that $x$ is an invariant element if, and only if, it commutes with the cofficients of the matrices $\overset{I}{C}{^{(\pm)}_{g,n}}$. As remarked in \eqref{CoeffsCCPlusCMoins}, the algebra generated by the coefficients $(\overset{I}{C}{^{(\pm)}_{g,n}}){^i_j}$ equals the algebra generated by the coefficients $(\overset{I}{C}_{g,n})^i_j$. Hence, an element is invariant if, and only if, it commutes with the coefficients of the matrices $\overset{I}{C}_{g,n}$.
\\2) Let $x \in \mathcal{L}_{g,n}^{\mathrm{inv}}(H)$ and $v \in \mathrm{Inv}(v)$, then
\[ \overset{I}{C}_{g,n} \triangleright (x \triangleright v) = (\overset{I}{C}_{g,n}x) \triangleright v = (x\overset{I}{C}_{g,n}) \triangleright v = x \triangleright (\overset{I}{C}_{g,n} \triangleright v) = \mathbb{I}_{\dim(I)}(x \triangleright v) \]
and thus $x \triangleright \varphi \in \mathrm{Inv}(V)$ by definition.
\end{proof}

\section{Projective representation of the mapping class group}
\indent Recall that the mapping class group $\mathrm{MCG}(\Sigma_{g,n})$ is the group of all isotopy classes of orientation-preserving homeomorphisms which fix the boundary pointwise. 
%
%

\smallskip

For simplicity we will mainly consider the case of $\Sigma_g$ ($n = 0$)\footnote{except in section \ref{sectionNorma} and at the beggining of section \ref{sectionLifts} were we deal with the general case ($n \geq 0$).}. The particular features in this case are that the presentation of the mapping class group is easier and that the associated algebra $\mathcal{L}_{g,0}(H) \cong \mathcal{H}(\mathcal{O}(H))^{\otimes g}$ is isomorphic to a matrix algebra. We discuss the extension of the construction to the case of $n>0$ in subsection \ref{CasGeneral}.

\subsection{Mapping class group of $\Sigma_g$}\label{MCGSigmaG}
\indent We begin with some terminology. A curve on $\Sigma_g^{\mathrm{o}}$ is called {\em simple} if it does not contain self-crossings (up to free homotopy). A simple closed curve (not necessarily oriented) on a surface up to free homotopy is simply called a {\em circle}. Elements of $\pi_1(\Sigma_g^{\mathrm{o}})$ (oriented based curves up to homotopy) are called {\em loops}. We say that a loop is {\em simple} if it does not contain self-crossing (up to homotopy). If $\gamma \in \pi_1(\Sigma_g^{\mathrm{o}})$, we denote by $[\gamma]$ the free homotopy class of $\gamma$. For $\alpha$ a circle, recall that we denote by $\tau_{\alpha}$ the Dehn twist about $\alpha$ (see \cite{FM}). If $\gamma \in \pi_1(\Sigma_g^{\mathrm{o}})$, then $\tau_{\gamma}$ is a shortand for $\tau_{[\gamma]}$, thus defined as follows: consider a circle $\gamma'$ freely homotopic to $\gamma$ and which does not intersect the boundary circle $c_g = \partial\!\left( \Sigma_g^{\mathrm{o}} \right)$, then $\tau_{\gamma} = \tau_{\gamma'}$. Of course all these notions make sense for $\Sigma_{g,n}$ as well.

\smallskip

\begin{figure}[!h]
\centering
\begingroup%
  \makeatletter%
  \providecommand\color[2][]{%
    \errmessage{(Inkscape) Color is used for the text in Inkscape, but the package 'color.sty' is not loaded}%
    \renewcommand\color[2][]{}%
  }%
  \providecommand\transparent[1]{%
    \errmessage{(Inkscape) Transparency is used (non-zero) for the text in Inkscape, but the package 'transparent.sty' is not loaded}%
    \renewcommand\transparent[1]{}%
  }%
  \providecommand\rotatebox[2]{#2}%
  \newcommand*\fsize{\dimexpr\f@size pt\relax}%
  \newcommand*\lineheight[1]{\fontsize{\fsize}{#1\fsize}\selectfont}%
  \ifx\svgwidth\undefined%
    \setlength{\unitlength}{542.66669247bp}%
    \ifx\svgscale\undefined%
      \relax%
    \else%
      \setlength{\unitlength}{\unitlength * \real{\svgscale}}%
    \fi%
  \else%
    \setlength{\unitlength}{\svgwidth}%
  \fi%
  \global\let\svgwidth\undefined%
  \global\let\svgscale\undefined%
  \makeatother%
  \begin{picture}(1,0.15466737)%
    \lineheight{1}%
    \setlength\tabcolsep{0pt}%
    \put(0,0){\includegraphics[width=\unitlength,page=1]{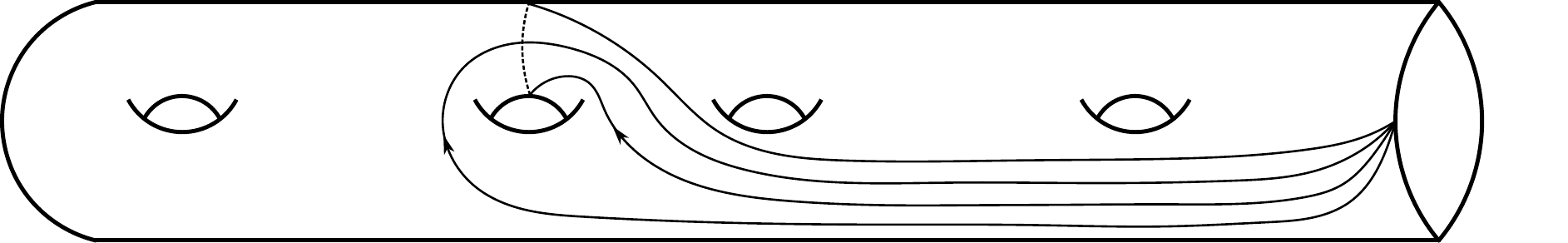}}%
    \put(0.28092784,0.02813831){\color[rgb]{0,0,0}\makebox(0,0)[lt]{\lineheight{1.25}\smash{\begin{tabular}[t]{l}$b_i$\end{tabular}}}}%
    \put(0.39128823,0.04212004){\color[rgb]{0,0,0}\makebox(0,0)[lt]{\lineheight{1.25}\smash{\begin{tabular}[t]{l}$a_i$\end{tabular}}}}%
    \put(0.11047155,0.09882403){\color[rgb]{0,0,0}\makebox(0,0)[lt]{\lineheight{1.25}\smash{\begin{tabular}[t]{l}1\end{tabular}}}}%
    \put(0.33152915,0.04886096){\color[rgb]{0,0,0}\makebox(0,0)[lt]{\lineheight{1.25}\smash{\begin{tabular}[t]{l}$i$\end{tabular}}}}%
    \put(0.46951632,0.09819392){\color[rgb]{0,0,0}\makebox(0,0)[lt]{\lineheight{1.25}\smash{\begin{tabular}[t]{l}$i-1$\end{tabular}}}}%
    \put(0.71901647,0.10243706){\color[rgb]{0,0,0}\makebox(0,0)[lt]{\lineheight{1.25}\smash{\begin{tabular}[t]{l}$g$\end{tabular}}}}%
    \put(0,0){\includegraphics[width=\unitlength,page=2]{surfaceAvecCourbes.pdf}}%
    \put(0.9518385,0.04924261){\color[rgb]{0,0,0}\makebox(0,0)[lt]{\lineheight{1.25}\smash{\begin{tabular}[t]{l}$c_g$\end{tabular}}}}%
  \end{picture}%
\endgroup%

\caption{Surface $\Sigma_{g}^{\mathrm{o}}$ with basepoint ($\bullet$), generators for $\pi_1(\Sigma_g^{\mathrm{o}})$ and boundary circle $c_g$.}
\label{figureSurface}
\end{figure}

\indent Recall that we take the loops $b_i, a_i$ ($1 \leq i \leq g$) represented in Figure \ref{figureSurface} as generators for the free group $\pi_1(\Sigma_g^{\mathrm{o}})$. With these generators, the boundary circle $c_g$ has the following expression: 
\[ c_g = b_1 a_1^{-1} b_1^{-1} a_1 \ldots b_g a_g^{-1} b_g^{-1} a_g. \]

\begin{figure}[!h]
\centering
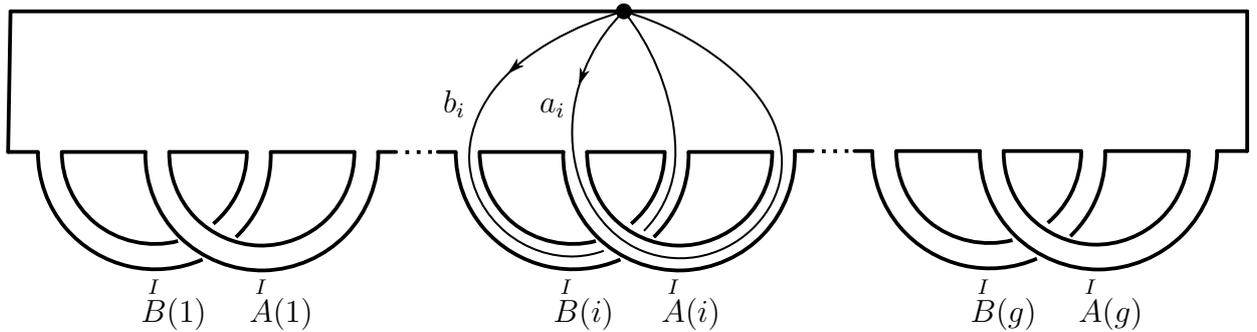
\caption{Surface $\Sigma_g^{\mathrm{o}}$ viewed as a thickening of the graph $\Gamma_{g,0}$, and matrices of generators of $\mathcal{L}_{g,0}(H)$.}
\label{figureSurfaceRuban}
\end{figure}

Retracting $\Sigma_g^{\mathrm{o}}$ to a tubular neighborhood of the loops $b_i$ and $a_i$, we get Figure \ref{figureSurfaceRuban}. In other words, $\Sigma_g^{\mathrm{o}}$ is homeomorphic to the thickening of the embedded oriented graph $\Gamma_{g,0}$ with vertex $\bullet$ and with edges $b_i$ and $a_i$.

\begin{figure}[!h]
\centering
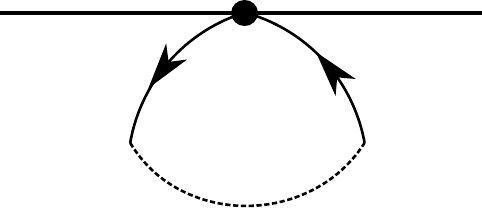
\caption{A positively oriented loop near the basepoint fixed in Figure \ref{surfaceAvecMatrices} or \ref{figureSurfaceRuban}.}
\label{figureCourbeOriente}
\end{figure}

\indent An important notion for the sequel is that of a positively oriented simple loop. 
\begin{definition}\label{defPositivelyOriented}
We say that a loop in $\pi_1(\Sigma_{g}^{\mathrm{o}})$ (or more generally in $\pi_1(\Sigma_{g,n}^{\mathrm{o}})$) is positively oriented if its orientation is counterclockwise, as indicated in Figure \ref{figureCourbeOriente}\footnote{Compared to the Figure 5 of \cite{Fai18c}, we have done a $180^{\circ}$-rotation around the horizontal axis of $\mathbb{R}^3$, in order to have the handles of $\Sigma_{g,n}^{\mathrm{o}}$ at the bottom of the Figure. The reason of this change comes from the definition of the graphical calculus and the Wilson loop map in Chapter \ref{chapitreGraphiqueSkein}.}; we say that it is negatively oriented if it is not positively oriented.
\end{definition}
\noindent Note that it is possible to have two simple loops $x_+, x_- \in \pi_1(\Sigma_{g,n}^{\mathrm{o}})$ such that $x_+$ is positively oriented, $x_-$ is negatively oriented and $[x_+] = [x_-]$. For instance in $\Sigma_1^{\mathrm{o}}$ take $x_+ = b a^{-1} b^{-1}$ and $x_- = a^{-1}$. Thus it makes no sense to say that a circle is positively or negatively oriented. Also recall that a simple loop or a circle is {\em non-separating} if it does not cut the surface into two connected components, otherwise it is called {\em separating}. All these properties (simple, non-separating, positively oriented) are preserved under the action of Dehn twists on the loops, hence they are preserved under the action of $\mathrm{MCG}(\Sigma_g^{\mathrm{o}})$ on the loops.

\smallskip

\begin{figure}[h]
\centering
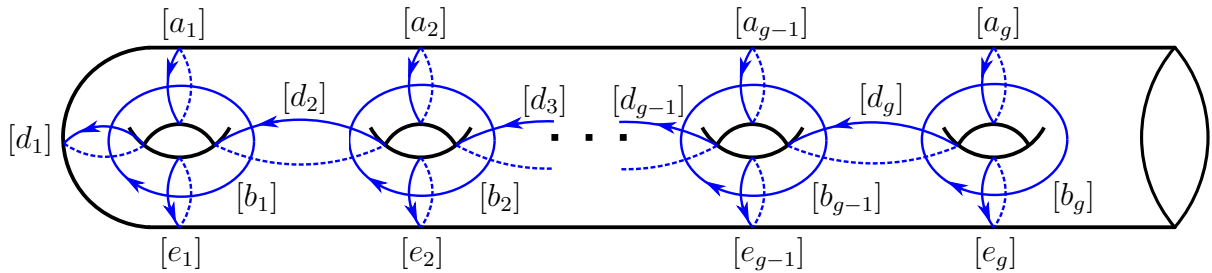
\caption{A canonical set of non-separating curves on the surface $\Sigma_g^{\mathrm{o}}$.}
\label{figureCourbesCanoniques}
\end{figure}

\indent In addition to the generating loops $b_i$ and $a_i$, we define the following loops on $\Sigma_g^{\mathrm{o}}$:
\begin{equation}\label{courbesGenHumphries}
\begin{split}
& d_1 = b_1 a_1^{-1} b_1^{-1}, \:\:\: d_i = a_{i-1} b_{i} a_{i}^{-1} b_{i}^{-1} \:\:\: \text{ for } 2 \leq i \leq g, \\
& e_1 = b_1 a_1^{-1} b_1^{-1}, \:\:\: e_i = b_1 a_1^{-1} b_1^{-1} a_1 \ldots b_{i-1} a_{i-1}^{-1} b_{i-1}^{-1} a_{i-1} b_i a_i^{-1} b_i^{-1} \:\:\: \text{ for } 2 \leq i \leq g,\\
&s_i = b_1 a_1^{-1} b_1^{-1} a_1 \ldots b_i a_i^{-1} b_i^{-1} a_i  \:\:\: \text{ for } 1 \leq i \leq g.
\end{split}
\end{equation}
The loops $b_i, a_i, d_i, e_i$ are simple, non-separating and positively oriented; their free homotopy class are depicted in Figure \ref{figureCourbesCanoniques}.

\begin{figure}[h]
\centering
\begingroup%
  \makeatletter%
  \providecommand\color[2][]{%
    \errmessage{(Inkscape) Color is used for the text in Inkscape, but the package 'color.sty' is not loaded}%
    \renewcommand\color[2][]{}%
  }%
  \providecommand\transparent[1]{%
    \errmessage{(Inkscape) Transparency is used (non-zero) for the text in Inkscape, but the package 'transparent.sty' is not loaded}%
    \renewcommand\transparent[1]{}%
  }%
  \providecommand\rotatebox[2]{#2}%
  \newcommand*\fsize{\dimexpr\f@size pt\relax}%
  \newcommand*\lineheight[1]{\fontsize{\fsize}{#1\fsize}\selectfont}%
  \ifx\svgwidth\undefined%
    \setlength{\unitlength}{428.33999733bp}%
    \ifx\svgscale\undefined%
      \relax%
    \else%
      \setlength{\unitlength}{\unitlength * \real{\svgscale}}%
    \fi%
  \else%
    \setlength{\unitlength}{\svgwidth}%
  \fi%
  \global\let\svgwidth\undefined%
  \global\let\svgscale\undefined%
  \makeatother%
  \begin{picture}(1,0.19148461)%
    \lineheight{1}%
    \setlength\tabcolsep{0pt}%
    \put(0,0){\includegraphics[width=\unitlength,page=1]{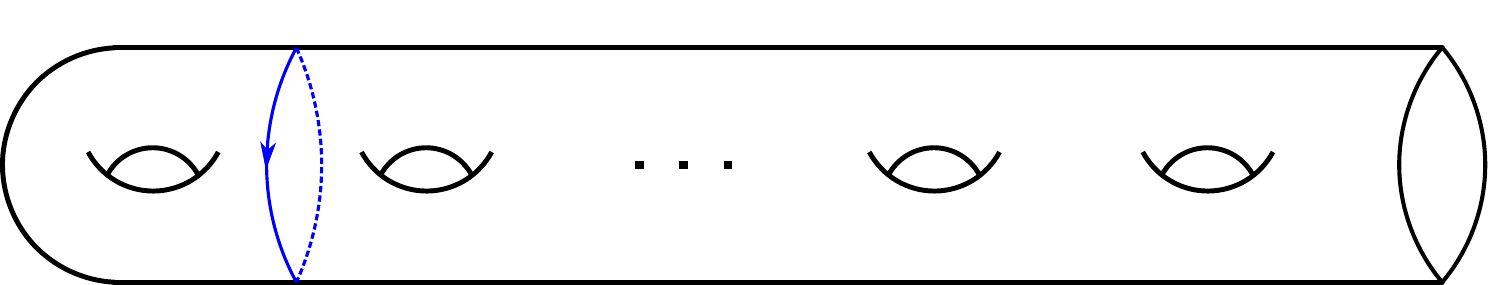}}%
    \put(0.17850201,0.17244639){\color[rgb]{0,0,0}\makebox(0,0)[lt]{\lineheight{1.25}\smash{\begin{tabular}[t]{l}$[s_1]$\end{tabular}}}}%
    \put(0,0){\includegraphics[width=\unitlength,page=2]{courbesSeparantes.pdf}}%
    \put(0.36428899,0.17121673){\color[rgb]{0,0,0}\makebox(0,0)[lt]{\lineheight{1.25}\smash{\begin{tabular}[t]{l}$[s_2]$\end{tabular}}}}%
    \put(0,0){\includegraphics[width=\unitlength,page=3]{courbesSeparantes.pdf}}%
    \put(0.5199364,0.17244638){\color[rgb]{0,0,0}\makebox(0,0)[lt]{\lineheight{1.25}\smash{\begin{tabular}[t]{l}$[s_{g-2}]$\end{tabular}}}}%
    \put(0,0){\includegraphics[width=\unitlength,page=4]{courbesSeparantes.pdf}}%
    \put(0.69503097,0.17244638){\color[rgb]{0,0,0}\makebox(0,0)[lt]{\lineheight{1.25}\smash{\begin{tabular}[t]{l}$[s_{g-1}]$\end{tabular}}}}%
    \put(0,0){\includegraphics[width=\unitlength,page=5]{courbesSeparantes.pdf}}%
    \put(0.88081788,0.17121673){\color[rgb]{0,0,0}\makebox(0,0)[lt]{\lineheight{1.25}\smash{\begin{tabular}[t]{l}$[s_g]$\end{tabular}}}}%
  \end{picture}%
\endgroup%

\caption{Canonical separating curves on the surface $\Sigma_g^{\mathrm{o}}$.}
\label{courbesSeparantes}
\end{figure}

The loops $s_i$ are simple, separating and positively oriented; their free homotopy class are depicted in Figure \ref{courbesSeparantes}. Note that $s_g = c_g$.

The Dehn twists $\tau_{e_2}, \tau_{b_i}, \tau_{d_i}$ are called the Humphries generators. There exists presentations of $\mathrm{MCG}(\Sigma_g)$ and $\mathrm{MCG}(\Sigma_g^{\mathrm{o}})$ due to Wajnryb \cite{wajnryb} (also see \cite[Sect. 5.2.1]{FM}):   $\mathrm{MCG}(\Sigma_g^{\mathrm{o}})$ is generated by the Humphries generators together with four families of relations called disjointness relations, braid relations, 3-chain relation and lantern relation, see \cite[Theorem 5.3]{FM} (the correspondence of notations with \cite[Figure 5.7]{FM} is $c_0 = [e_2]$, $c_{2j} = [b_j]$, $c_{2j-1} = [d_j]$). The presentation of $\mathrm{MCG}(\Sigma_g)$ is obtained as the quotient of $\mathrm{MCG}(\Sigma_g^{\mathrm{o}})$ by the hyperelliptic relation:
\begin{equation}\label{hyperelliptic}
\left( \tau_{b_g} \tau_{d_g} \ldots \tau_{b_1} \tau_{d_1} \tau_{d_1} \tau_{b_1} \ldots \tau_{d_g} \tau_{b_g} \right) w = w  \left( \tau_{b_g} \tau_{d_g} \ldots \tau_{b_1} \tau_{d_1} \tau_{d_1} \tau_{b_1} \ldots \tau_{d_g} \tau_{b_g} \right)
\end{equation}
where $w$ is any word in the Humphries generators which equals $\tau_{a_g}$.

\smallskip

\indent Since the mapping class group fixes the boundary pointwise, the basepoint is fixed and we can consider the action of $\mathrm{MCG}(\Sigma_{g}^{\mathrm{o}})$ on $\pi_1(\Sigma_{g}^{\mathrm{o}})$. The actions of the Humphries generators on the fundamental group are easily computed (see \eqref{figureDehnTwist}). We just indicate the non-trivial actions:
\begin{equation}\label{actionPi1}
\begin{split}
& \tau_{e_2}(a_1) = e_2^{-1} a_1 e_2, \:\: \tau_{e_2}(b_1) = e_2^{-1} b_1 e_2, \:\: \tau_{e_2}(b_2) = e_2^{-1} b_2, \\
& \tau_{b_i}(a_i) = b^{-1}_i a_i, \\
& \tau_{a_1}(b_1) = b_1 a_1 \:\:\:\: (\text{note that } [a_1] = [d_1^{-1}]),\\
& \tau_{d_i}(a_{i-1}) = d_i^{-1} a_{i-1} d_i, \:\: \tau_{d_i}(b_{i-1}) = b_{i-1} d_i, \:\: \tau_{d_i}(b_i) = d_i^{-1} b_i \:\:\: \text{(with } i \geq 2\text{)}.
\end{split}
\end{equation}

\subsection{Normalization of simple closed curves}\label{sectionNorma}
\indent In this section, we associate an integer to any oriented circle and to any simple loop in $\pi_1(\Sigma_{g,n}^{\mathrm{o}})$ (we define this integer for the general case $n \geq 0$). To define these quantities, we use the view of $\Sigma_{g,n}^{\mathrm{o}}$ depicted in Figure \ref{surfaceAvecMatrices}. Note that it is not the assignment defined in \cite[Section 9]{AS} (it gives different values). This comes from the fact that their normalizations by powers of $v$ of particular product of matrices (see section \ref{sectionLifts}) differ from ours. This is maybe due to their normalization of Clebsch-Gordan operators (which are not used here since they are defined in the semi-simple case only).

\smallskip

\indent Let $\gamma \subset \Sigma_{g,n}^{\mathrm{o}}$ be an oriented circle (which is not isotopic to a point). Using isotopy, we may assume that the handles of Figure \ref{surfaceAvecMatrices} contain only bunches of parallel strands. Since $\gamma$ is simple, the ``rectangle'' in Figure \ref{surfaceAvecMatrices} can contain only vertical strands, caps $\cap$ and cups $\cup$. We choose the following preferred direction on $\Sigma_{g,n}^{\mathrm{o}}$:
\begin{center}
\begingroup%
  \makeatletter%
  \providecommand\color[2][]{%
    \errmessage{(Inkscape) Color is used for the text in Inkscape, but the package 'color.sty' is not loaded}%
    \renewcommand\color[2][]{}%
  }%
  \providecommand\transparent[1]{%
    \errmessage{(Inkscape) Transparency is used (non-zero) for the text in Inkscape, but the package 'transparent.sty' is not loaded}%
    \renewcommand\transparent[1]{}%
  }%
  \providecommand\rotatebox[2]{#2}%
  \newcommand*\fsize{\dimexpr\f@size pt\relax}%
  \newcommand*\lineheight[1]{\fontsize{\fsize}{#1\fsize}\selectfont}%
  \ifx\svgwidth\undefined%
    \setlength{\unitlength}{520.1446052bp}%
    \ifx\svgscale\undefined%
      \relax%
    \else%
      \setlength{\unitlength}{\unitlength * \real{\svgscale}}%
    \fi%
  \else%
    \setlength{\unitlength}{\svgwidth}%
  \fi%
  \global\let\svgwidth\undefined%
  \global\let\svgscale\undefined%
  \makeatother%
  \begin{picture}(1,0.16612622)%
    \lineheight{1}%
    \setlength\tabcolsep{0pt}%
    \put(0,0){\includegraphics[width=\unitlength,page=1]{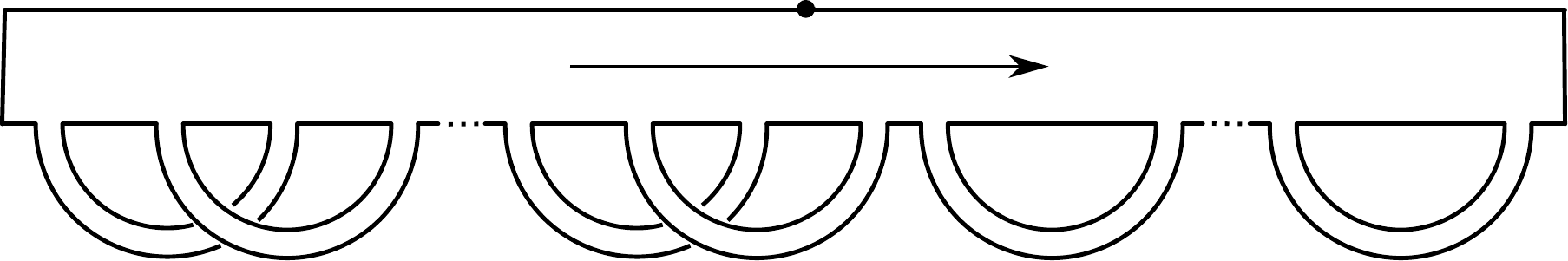}}%
  \end{picture}%
\endgroup%

\end{center}
We define $N_{\cup}(\gamma)$ as the number of cups and strands-in-handle which run against this preferred direction; similarly, we define $N_{\cap}(\gamma)$ as the number of caps which run against this preferred direction. More precisely, $N_{\cup}(\gamma)$ (resp. $N_{\cap}(\gamma)$) is the number of pieces of $\gamma$ which look like the following:
\begin{center}
\begingroup%
  \makeatletter%
  \providecommand\color[2][]{%
    \errmessage{(Inkscape) Color is used for the text in Inkscape, but the package 'color.sty' is not loaded}%
    \renewcommand\color[2][]{}%
  }%
  \providecommand\transparent[1]{%
    \errmessage{(Inkscape) Transparency is used (non-zero) for the text in Inkscape, but the package 'transparent.sty' is not loaded}%
    \renewcommand\transparent[1]{}%
  }%
  \providecommand\rotatebox[2]{#2}%
  \newcommand*\fsize{\dimexpr\f@size pt\relax}%
  \newcommand*\lineheight[1]{\fontsize{\fsize}{#1\fsize}\selectfont}%
  \ifx\svgwidth\undefined%
    \setlength{\unitlength}{462.37670637bp}%
    \ifx\svgscale\undefined%
      \relax%
    \else%
      \setlength{\unitlength}{\unitlength * \real{\svgscale}}%
    \fi%
  \else%
    \setlength{\unitlength}{\svgwidth}%
  \fi%
  \global\let\svgwidth\undefined%
  \global\let\svgscale\undefined%
  \makeatother%
  \begin{picture}(1,0.25490063)%
    \lineheight{1}%
    \setlength\tabcolsep{0pt}%
    \put(0,0){\includegraphics[width=\unitlength,page=1]{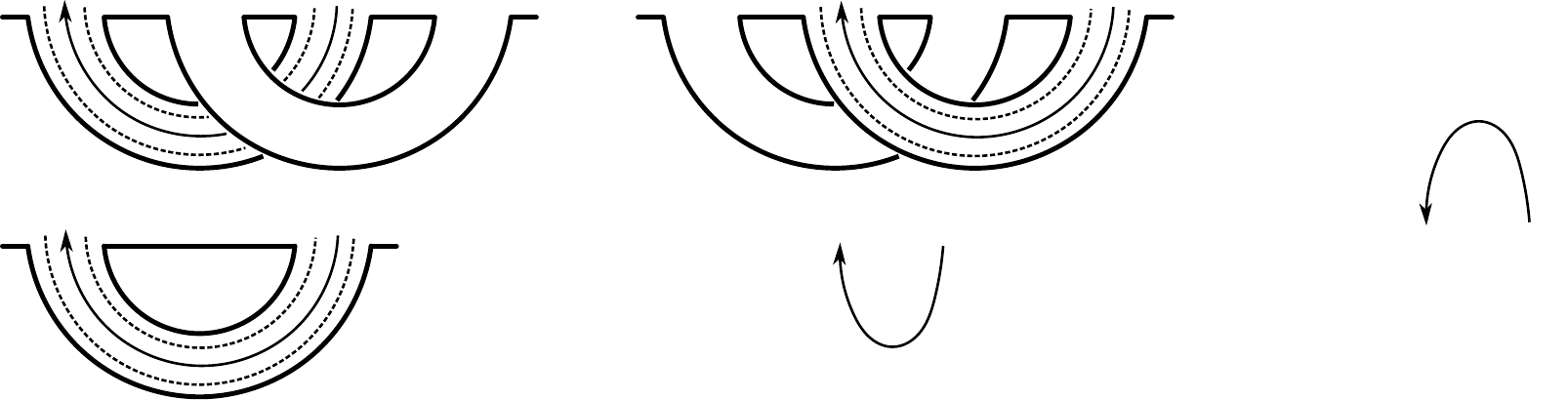}}%
    \put(0.8323205,0.13468499){\color[rgb]{0,0,0}\makebox(0,0)[lt]{\lineheight{1.25}\smash{\begin{tabular}[t]{l}(resp. \end{tabular}}}}%
    \put(0.98986216,0.13490337){\color[rgb]{0,0,0}\makebox(0,0)[lt]{\lineheight{1.25}\smash{\begin{tabular}[t]{l}) \end{tabular}}}}%
  \end{picture}%
\endgroup%

\end{center}
If $x \in \pi_1\!\left( \Sigma_{g,n}^{\mathrm{o}} \right)$ is a non-trivial loop, we define $N_{\cup}(x)$ and $N_{\cap}(x)$ by the same formula (we stress that the junction of the loop at  the basepoint is not considered as a cap).

\begin{definition}
Let $\eta$ be an oriented circle $\subset \Sigma_{g,n}^{\mathrm{o}}$ or a loop $\in \pi_1\!\left( \Sigma_{g,n}^{\mathrm{o}} \right)$. The normalization of $\eta$ is $N(\eta) = N_{\cup}(\eta) - N_{\cap}(\eta) \in \mathbb{Z}$.
\end{definition}
It is clear that $N(\eta)$ does not depend of the homotopy class of $\eta$. Note that if $x \in \pi_1\!\left( \Sigma_{g,n}^{\mathrm{o}} \right)$, we have 
\begin{equation}\label{NLoopNFree}
N(x) =
\begin{cases}
N([x]) + 1 & \text{if } x \text{ is positively oriented}\\
N([x]) & \text{if } x \text{ is negatively oriented}
\end{cases}
\end{equation}
Also observe that if $\gamma$ is an oriented circle, then 
\begin{equation}\label{NInvFree}
N(\gamma^{-1}) = -N(\gamma)
\end{equation}
where $\gamma^{-1}$ is $\gamma$ with the opposite orientation.
\begin{exemple}\label{exempleN}
Consider the oriented circle $\gamma$ and the loop $x$ depicted below:
\begin{center}
\begingroup%
  \makeatletter%
  \providecommand\color[2][]{%
    \errmessage{(Inkscape) Color is used for the text in Inkscape, but the package 'color.sty' is not loaded}%
    \renewcommand\color[2][]{}%
  }%
  \providecommand\transparent[1]{%
    \errmessage{(Inkscape) Transparency is used (non-zero) for the text in Inkscape, but the package 'transparent.sty' is not loaded}%
    \renewcommand\transparent[1]{}%
  }%
  \providecommand\rotatebox[2]{#2}%
  \newcommand*\fsize{\dimexpr\f@size pt\relax}%
  \newcommand*\lineheight[1]{\fontsize{\fsize}{#1\fsize}\selectfont}%
  \ifx\svgwidth\undefined%
    \setlength{\unitlength}{380.32771534bp}%
    \ifx\svgscale\undefined%
      \relax%
    \else%
      \setlength{\unitlength}{\unitlength * \real{\svgscale}}%
    \fi%
  \else%
    \setlength{\unitlength}{\svgwidth}%
  \fi%
  \global\let\svgwidth\undefined%
  \global\let\svgscale\undefined%
  \makeatother%
  \begin{picture}(1,0.21549293)%
    \lineheight{1}%
    \setlength\tabcolsep{0pt}%
    \put(0,0){\includegraphics[width=\unitlength,page=1]{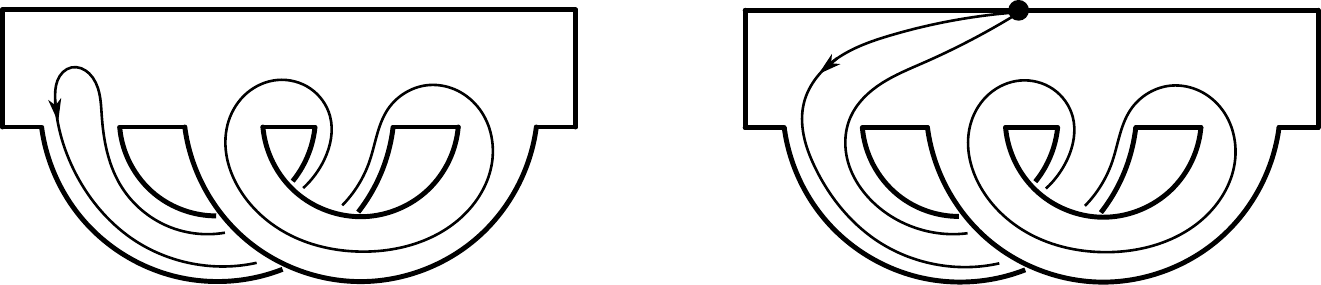}}%
    \put(0.87734232,0.16473637){\color[rgb]{0,0,0}\makebox(0,0)[lt]{\lineheight{1.25}\smash{\begin{tabular}[t]{l}$x$\end{tabular}}}}%
    \put(0.315683,0.16237054){\color[rgb]{0,0,0}\makebox(0,0)[lt]{\lineheight{1.25}\smash{\begin{tabular}[t]{l}$\gamma$\end{tabular}}}}%
  \end{picture}%
\endgroup%

\end{center}
We have $N_{\cup}(\gamma) = 2$ and $N_{\cap}(\gamma) = 1$, thus $N(\gamma) = 1$. For $x = ba^{-1}b^{-1}$, we have $N_{\cup}(x) = 2$ and $N_{\cap}(x) = 0$, thus $N(x)=2$.
\finEx
\end{exemple}

We can define $N$ in a different manner. Consider the following figure:
\begin{center}
\begingroup%
  \makeatletter%
  \providecommand\color[2][]{%
    \errmessage{(Inkscape) Color is used for the text in Inkscape, but the package 'color.sty' is not loaded}%
    \renewcommand\color[2][]{}%
  }%
  \providecommand\transparent[1]{%
    \errmessage{(Inkscape) Transparency is used (non-zero) for the text in Inkscape, but the package 'transparent.sty' is not loaded}%
    \renewcommand\transparent[1]{}%
  }%
  \providecommand\rotatebox[2]{#2}%
  \newcommand*\fsize{\dimexpr\f@size pt\relax}%
  \newcommand*\lineheight[1]{\fontsize{\fsize}{#1\fsize}\selectfont}%
  \ifx\svgwidth\undefined%
    \setlength{\unitlength}{520.14459159bp}%
    \ifx\svgscale\undefined%
      \relax%
    \else%
      \setlength{\unitlength}{\unitlength * \real{\svgscale}}%
    \fi%
  \else%
    \setlength{\unitlength}{\svgwidth}%
  \fi%
  \global\let\svgwidth\undefined%
  \global\let\svgscale\undefined%
  \makeatother%
  \begin{picture}(1,0.18996805)%
    \lineheight{1}%
    \setlength\tabcolsep{0pt}%
    \put(0,0){\includegraphics[width=\unitlength,page=1]{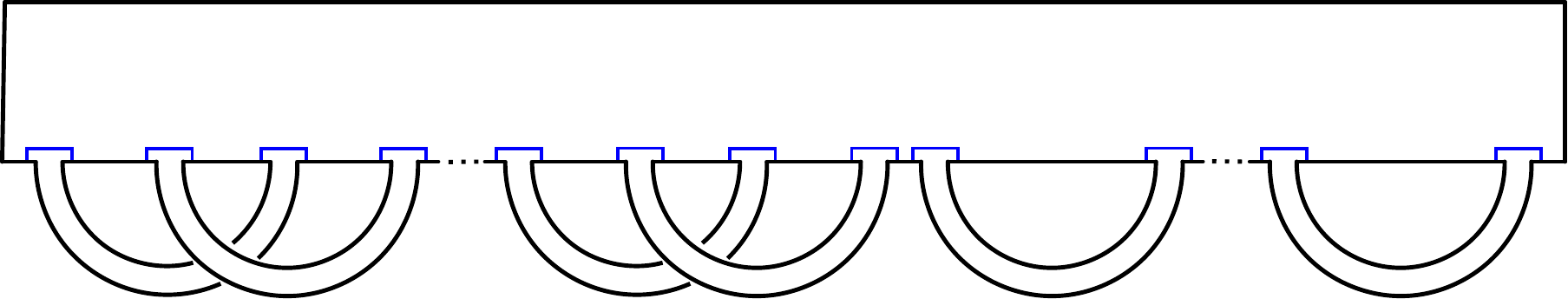}}%
    \put(0.02732971,0.10032977){\color[rgb]{0,0,0}\makebox(0,0)[lt]{\lineheight{1.25}\smash{\begin{tabular}[t]{l}$1$\end{tabular}}}}%
    \put(0.10419083,0.10000991){\color[rgb]{0,0,0}\makebox(0,0)[lt]{\lineheight{1.25}\smash{\begin{tabular}[t]{l}$2$\end{tabular}}}}%
    \put(0.1762455,0.10179344){\color[rgb]{0,0,0}\makebox(0,0)[lt]{\lineheight{1.25}\smash{\begin{tabular}[t]{l}$3$\end{tabular}}}}%
    \put(0.2543641,0.10215016){\color[rgb]{0,0,0}\makebox(0,0)[lt]{\lineheight{1.25}\smash{\begin{tabular}[t]{l}$4$\end{tabular}}}}%
    \put(0.31016076,0.10311298){\color[rgb]{0,0,0}\makebox(0,0)[lt]{\lineheight{1.25}\smash{\begin{tabular}[t]{l}$4g-3$\end{tabular}}}}%
    \put(0.38658569,0.1038837){\color[rgb]{0,0,0}\makebox(0,0)[lt]{\lineheight{1.25}\smash{\begin{tabular}[t]{l}$4g-2$\end{tabular}}}}%
    \put(0.45847941,0.10363427){\color[rgb]{0,0,0}\makebox(0,0)[lt]{\lineheight{1.25}\smash{\begin{tabular}[t]{l}$4g-1$\end{tabular}}}}%
    \put(0.54169642,0.10368793){\color[rgb]{0,0,0}\makebox(0,0)[lt]{\lineheight{1.25}\smash{\begin{tabular}[t]{l}$4g$\end{tabular}}}}%
    \put(0.5795809,0.10425428){\color[rgb]{0,0,0}\makebox(0,0)[lt]{\lineheight{1.25}\smash{\begin{tabular}[t]{l}$4g+1$\end{tabular}}}}%
    \put(0.72237872,0.10382574){\color[rgb]{0,0,0}\makebox(0,0)[lt]{\lineheight{1.25}\smash{\begin{tabular}[t]{l}$4g+2$\end{tabular}}}}%
    \put(0.79525542,0.10412884){\color[rgb]{0,0,0}\makebox(0,0)[lt]{\lineheight{1.25}\smash{\begin{tabular}[t]{l}$4g+2n-1$\end{tabular}}}}%
    \put(0.9214874,0.10447928){\color[rgb]{0,0,0}\makebox(0,0)[lt]{\lineheight{1.25}\smash{\begin{tabular}[t]{l}$4g+2n$\end{tabular}}}}%
  \end{picture}%
\endgroup%

\end{center}
The extremities of each handle are endowed with the blue lines, which we call ``gates'', numbered from $1$ to $4g+2n$. First, let $x \in \pi_1\!\left( \Sigma_{g,n}^{\mathrm{o}} \right)$ be a simple loop. Starting from the basepoint and following $x$ along its orientation, we meet a first gate numbered $g_1$, then a second gate numbered $g_2$ and so on. This provides a sequence $g(x) = (g_1, \ldots, g_{2k})$ and we have
\[ N(x) = \sum_{i=1}^{k} \delta\!\left( g_{2i-1} > g_{2i} \right) - \sum_{i=1}^{k-1} \delta\!\left( g_{2i} \geq g_{2i+1} \right), \]
where $\delta(a > b)$ (resp. $\delta(a \geq b)$) is $1$ if $a>b$ (resp. $a \geq b$) and $0$ otherwise. Now, if $\gamma$ is a circle, we do not have a canonical starting point. Instead, choose a point where $\gamma$ meets one of the gates to {\em enter} in a handle, numbered $g_1$, and follow $\gamma$ along its orientation. This gives a sequence as previously except that we meet $g_1$ two times, at the begining and at the end. In other words, $g(\gamma) = (g_1, \ldots, g_{2k}, g_{2k+1} = g_1)$, and we have
\[ N(\gamma) = \sum_{i=1}^{k} \delta\!\left( g_{2i-1} > g_{2i} \right) - \sum_{i=1}^{k} \delta\!\left( g_{2i} \geq g_{2i+1} \right). \]
It is clear that this quantity does not depend on the choice of the starting point.
\begin{exemple}
Take back the cases of Example \ref{exempleN}. For the circle $\gamma$, we put the starting point at the left of the gate $1$. Then we obtain $g(\gamma) = \left( 1, \: 3, \: 4, \: 2, \: 3, \: 1, \: 1 \right)$ and we recover $N(\gamma)=1$. We might as well have started from the left of the gate $3$ or from the gate $4$. For the loop $x$, $g(x) = \left( 1, \: 3, \: 4, \: 2, \: 3, \: 1 \right)$ and we recover $N(x) = 2$.
\finEx
\end{exemple}

Let $\omega$ be the algebraic intersection form of simple curves (simple loops or oriented circles); recall that $\omega(\alpha, \beta)$ is the sum of the indices of the intersection points of $\alpha$ and $\beta$, as follows:
\begin{center}
\begingroup%
  \makeatletter%
  \providecommand\color[2][]{%
    \errmessage{(Inkscape) Color is used for the text in Inkscape, but the package 'color.sty' is not loaded}%
    \renewcommand\color[2][]{}%
  }%
  \providecommand\transparent[1]{%
    \errmessage{(Inkscape) Transparency is used (non-zero) for the text in Inkscape, but the package 'transparent.sty' is not loaded}%
    \renewcommand\transparent[1]{}%
  }%
  \providecommand\rotatebox[2]{#2}%
  \newcommand*\fsize{\dimexpr\f@size pt\relax}%
  \newcommand*\lineheight[1]{\fontsize{\fsize}{#1\fsize}\selectfont}%
  \ifx\svgwidth\undefined%
    \setlength{\unitlength}{429.87985528bp}%
    \ifx\svgscale\undefined%
      \relax%
    \else%
      \setlength{\unitlength}{\unitlength * \real{\svgscale}}%
    \fi%
  \else%
    \setlength{\unitlength}{\svgwidth}%
  \fi%
  \global\let\svgwidth\undefined%
  \global\let\svgscale\undefined%
  \makeatother%
  \begin{picture}(1,0.11975845)%
    \lineheight{1}%
    \setlength\tabcolsep{0pt}%
    \put(0,0){\includegraphics[width=\unitlength,page=1]{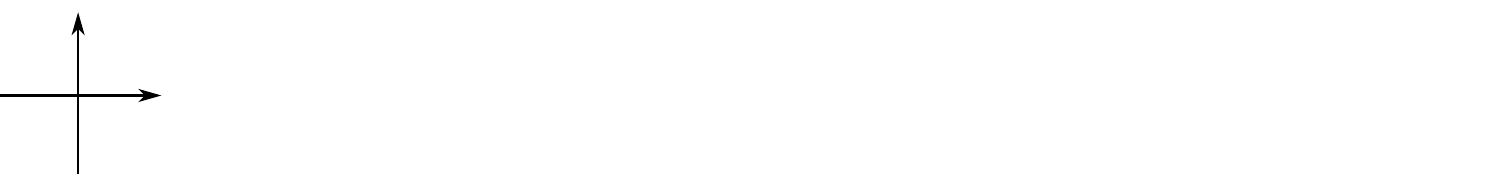}}%
    \put(0.08833594,0.0044284){\color[rgb]{0,0,0}\makebox(0,0)[lt]{\lineheight{1.25}\smash{\begin{tabular}[t]{l}$+1$\end{tabular}}}}%
    \put(0.06077433,0.1091705){\color[rgb]{0,0,0}\makebox(0,0)[lt]{\lineheight{1.25}\smash{\begin{tabular}[t]{l}$\beta$\end{tabular}}}}%
    \put(0.10261499,0.06644169){\color[rgb]{0,0,0}\makebox(0,0)[lt]{\lineheight{1.25}\smash{\begin{tabular}[t]{l}$\alpha$\end{tabular}}}}%
    \put(0,0){\includegraphics[width=\unitlength,page=2]{formeIntersection.pdf}}%
    \put(0.36696513,0.00520627){\color[rgb]{0,0,0}\makebox(0,0)[lt]{\lineheight{1.25}\smash{\begin{tabular}[t]{l}$-1$\end{tabular}}}}%
    \put(0.3399221,0.1091705){\color[rgb]{0,0,0}\makebox(0,0)[lt]{\lineheight{1.25}\smash{\begin{tabular}[t]{l}$\beta$\end{tabular}}}}%
    \put(0.38176276,0.06644169){\color[rgb]{0,0,0}\makebox(0,0)[lt]{\lineheight{1.25}\smash{\begin{tabular}[t]{l}$\alpha$\end{tabular}}}}%
    \put(0,0){\includegraphics[width=\unitlength,page=3]{formeIntersection.pdf}}%
    \put(0.64611288,0.00520628){\color[rgb]{0,0,0}\makebox(0,0)[lt]{\lineheight{1.25}\smash{\begin{tabular}[t]{l}$-1$\end{tabular}}}}%
    \put(0.61906988,0.1091705){\color[rgb]{0,0,0}\makebox(0,0)[lt]{\lineheight{1.25}\smash{\begin{tabular}[t]{l}$\beta$\end{tabular}}}}%
    \put(0.66091054,0.06644169){\color[rgb]{0,0,0}\makebox(0,0)[lt]{\lineheight{1.25}\smash{\begin{tabular}[t]{l}$\alpha$\end{tabular}}}}%
    \put(0,0){\includegraphics[width=\unitlength,page=4]{formeIntersection.pdf}}%
    \put(0.92526081,0.00520628){\color[rgb]{0,0,0}\makebox(0,0)[lt]{\lineheight{1.25}\smash{\begin{tabular}[t]{l}$+1$\end{tabular}}}}%
    \put(0.89821755,0.10917051){\color[rgb]{0,0,0}\makebox(0,0)[lt]{\lineheight{1.25}\smash{\begin{tabular}[t]{l}$\beta$\end{tabular}}}}%
    \put(0.94005837,0.0664417){\color[rgb]{0,0,0}\makebox(0,0)[lt]{\lineheight{1.25}\smash{\begin{tabular}[t]{l}$\alpha$\end{tabular}}}}%
  \end{picture}%
\endgroup%

\end{center}

\begin{lemma}\label{normImageTwist}
Let $\gamma$ be a circle endowed with an arbitrary orientation and let $x \in \pi_1\!\left(\Sigma_{g,n}^{\mathrm{o}}\right)$ be a simple loop; it holds
\[ N\!\left( \tau_{\gamma}(x) \right) = N(x) + \omega(x,\gamma) N(\gamma). \]
This formula is also true if $x$ is an oriented circle.
\end{lemma}
\begin{proof}
Let $p_1, \ldots, p_k$ be the intersection points between $x$ and $\gamma$. By the method for computing the action of a Dehn twist on a simple closed curve (see \eqref{figureDehnTwist}), it is clear that $N\!\left( \tau_{\gamma}(x) \right) = N(x) + \epsilon_1 N(\gamma) + \ldots + \epsilon_k N(\gamma)$, with $\epsilon_i = \pm 1$. Moreover, it is not difficult to check that each $\epsilon_i$ indeed is the index of the intersection point $p_i$, as defined above. For instance, assume that the intersection looks like in the figure below
\begin{center}
\begingroup%
  \makeatletter%
  \providecommand\color[2][]{%
    \errmessage{(Inkscape) Color is used for the text in Inkscape, but the package 'color.sty' is not loaded}%
    \renewcommand\color[2][]{}%
  }%
  \providecommand\transparent[1]{%
    \errmessage{(Inkscape) Transparency is used (non-zero) for the text in Inkscape, but the package 'transparent.sty' is not loaded}%
    \renewcommand\transparent[1]{}%
  }%
  \providecommand\rotatebox[2]{#2}%
  \newcommand*\fsize{\dimexpr\f@size pt\relax}%
  \newcommand*\lineheight[1]{\fontsize{\fsize}{#1\fsize}\selectfont}%
  \ifx\svgwidth\undefined%
    \setlength{\unitlength}{244.20398846bp}%
    \ifx\svgscale\undefined%
      \relax%
    \else%
      \setlength{\unitlength}{\unitlength * \real{\svgscale}}%
    \fi%
  \else%
    \setlength{\unitlength}{\svgwidth}%
  \fi%
  \global\let\svgwidth\undefined%
  \global\let\svgscale\undefined%
  \makeatother%
  \begin{picture}(1,0.25652359)%
    \lineheight{1}%
    \setlength\tabcolsep{0pt}%
    \put(0,0){\includegraphics[width=\unitlength,page=1]{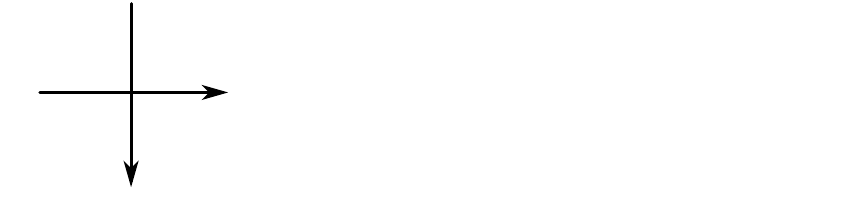}}%
    \put(0.14589018,0.0040971){\color[rgb]{0,0,0}\makebox(0,0)[lt]{\lineheight{1.25}\smash{\begin{tabular}[t]{l}$x$\end{tabular}}}}%
    \put(-0.0023757,0.13938584){\color[rgb]{0,0,0}\makebox(0,0)[lt]{\lineheight{1.25}\smash{\begin{tabular}[t]{l}$\gamma$\end{tabular}}}}%
    \put(0,0){\includegraphics[width=\unitlength,page=2]{twistOriente.pdf}}%
    \put(0.79666278,0.14063172){\color[rgb]{0,0,0}\makebox(0,0)[lt]{\lineheight{1.25}\smash{\begin{tabular}[t]{l}$\tau_{\gamma}(x)$\end{tabular}}}}%
    \put(0,0){\includegraphics[width=\unitlength,page=3]{twistOriente.pdf}}%
  \end{picture}%
\endgroup%

\end{center}
We see that in the resulting curve we follow a copy of $\gamma$ in the sense of its orientation. Hence, the contribution of this operation to the final result is $+N(\gamma)$. And indeed we have $\omega(x,\gamma) = -\omega(\gamma, x) = +1$.
\end{proof}

\indent Now, let $\pi_1^v\!\left( \Sigma_{g,n}^{\mathrm{o}} \right)$ be the group $\mathbb{Z} \times \pi_1\!\left( \Sigma_{g,n}^{\mathrm{o}} \right)$. The elements of $\pi_1^v\!\left( \Sigma_{g,n}^{\mathrm{o}} \right)$ are of the form $v^nx$, where $x \in \pi_1\!\left( \Sigma_{g,n}^{\mathrm{o}} \right)$ and $v$ is a formal element commuting with $\pi_1\!\left( \Sigma_{g,n}^{\mathrm{o}} \right)$. 

\smallskip

\indent For the remaining of this section, we restrict to $n=0$ (see however the discussion in section \ref{CasGeneral}). We define group automorphisms $\tau^v_{e_2}, \tau^v_{b_i}, \tau^v_{d_i} : \pi_1^v\!\left(\Sigma_g^{\mathrm{o}}\right) \to \pi_1^v\!\left(\Sigma_g^{\mathrm{o}}\right)$ which are normalized versions of \eqref{actionPi1}:
\begin{equation}\label{actionPi1V}
\begin{split}
& \tau^v_{e_2}(v) =v, \:\:\:\: \tau^v_{b_i}(v) = v, \:\:\:\: \tau^v_{a_1}(v) = v, \:\:\:\: \tau^v_{d_i}(v) = v,\\
& \tau^v_{e_2}(a_1) = e_2^{-1} a_1 e_2, \:\: \tau^v_{e_2}(b_1) = e_2^{-1} b_1 e_2, \:\: \tau^v_{e_2}(b_2) = v^{-3}e_2^{-1} b_2, \\
& \tau^v_{b_i}(a_i) = vb^{-1}_i a_i, \\
& \tau^v_{a_1}(b_1) = v^{-1}b_1 a_1 \:\:\:\: (\text{recall that } [a_1] = [d_1^{-1}]),\\
& \tau^v_{d_i}(a_{i-1}) = d_i^{-1} a_{i-1} d_i, \:\: \tau^v_{d_i}(b_{i-1}) = vb_{i-1} d_i, \:\: \tau^v_{d_i}(b_i) = v^{-1}d_i^{-1} b_i \:\:\: \text{(with } i \geq 2\text{)}.
\end{split}
\end{equation}
We just indicate the values on the generators of $\pi_1\!\left( \Sigma_{g}^{\mathrm{o}} \right)$ which are not fixed. Note that we have defined these group automorphisms in the following way:
\begin{equation}\label{normTauV}
\tau^v_{y}(u_i) = v^{N(\tau_y(u_i))} \tau_y(u_i) = v^{\omega(u_i, [y])N([y])} \tau_y(u_i)
\end{equation}
where $u_i$ is $b_i$ or $a_i$ and $y$ is $e_2, b_i$ or $d_i$\footnote{We recall from section \ref{MCGSigmaG} that if $y \in \pi_1\!\left( \Sigma_g^{\mathrm{o}} \right)$, then $\tau_y$ is a shortand for $\tau_{[y]}$ and in this case $y$ must be considered as a circle.}.

\begin{proposition}\label{propNormalisation}
1) The assignment
\[ \tau_{e_2} \mapsto \tau^v_{e_2}, \:\:\: \tau_{b_i} \mapsto \tau^v_{b_i}, \:\:\: \tau_{d_i} \mapsto \tau^v_{d_i} \]
extends to a morphism of groups $\mathrm{MCG}(\Sigma_g^{\mathrm{o}}) \to \mathrm{Aut}\!\left(\pi_1^v\!\left( \Sigma_{g}^{\mathrm{o}} \right)\right)$.
\\2) For $f \in \mathrm{MCG}(\Sigma_g^{\mathrm{o}})$, denote by $f^v : \pi_1^v\!\left( \Sigma_{g}^{\mathrm{o}} \right) \to \pi_1^v\!\left( \Sigma_{g}^{\mathrm{o}} \right)$ the image of $f$ by this morphism. Then for all simple loop $x \in \pi_1\!\left( \Sigma_{g}^{\mathrm{o}} \right)$, it holds
\begin{equation}\label{fvEtNSeMarientBien}
f^v(v^{N(x)}x) = v^{N(f(x))}f(x). 
\end{equation}
\end{proposition}
\begin{proof}
Thanks to Lemma \ref{normImageTwist} and to \eqref{normTauV}, we see that
\begin{equation}\label{NEtComposition}
\tau^v_{y_1} \circ \tau_{y_2}^v(v^{N(x)}x) = v^{N(\tau_{y_1} \circ \tau_{y_2}(x))} \tau_{y_1} \circ \tau_{y_2}(x)
\end{equation}
for all $x \in \pi_1\!\left( \Sigma_g^{\mathrm{o}} \right)$, where $y_1, y_2$ are $e_2, b_i$ or $d_i$. Let $\tau_{y_1}^{\epsilon_1} \ldots \tau_{y_k}^{\epsilon_k} = \mathrm{id}$ be a relation in $\mathrm{MCG}(\Sigma_g^{\mathrm{o}})$ (where each $y_j$ is $e_2, b_i$ or $d_i$ and $\epsilon_j \in \mathbb{Z}$), then we have $\tau_{y_1}^{\epsilon_1} \ldots \tau_{y_k}^{\epsilon_k}(u_i) = u_i$ for each $u_i = a_i$ or $b_i$. It follows that $(\tau_{y_1}^v)^{\epsilon_1} \ldots (\tau_{y_k}^v)^{\epsilon_k}(u_i) = v^{N(u_i)} u_i = u_i$ and thus $(\tau_{y_1}^v)^{\epsilon_1} \ldots (\tau_{y_k}^v)^{\epsilon_k}(u_i) = \mathrm{id}$. Alternatively, one can check tediously that the assignment preserve the Wajnryb relations. The second claim follows from \eqref{NEtComposition} and the fact that $\tau_{e_2}, \tau_{b_i}, \tau_{d_i}$ generate $\mathrm{MCG}(\Sigma_g^{\mathrm{o}})$.
\end{proof}

\subsection{Lifting simple loops and mapping classes to $\mathcal{L}_{g,0}(H)$}\label{sectionLifts}
\indent We define the lift of a simple loop in the general case ($n \geq 0$). In the group $\pi_1^v(\Sigma_{g,n}^{\mathrm{o}})$ defined in the previous section, we have the loops $b_i, a_i, m_j$ and the formal variable $v$, while in $\mathcal{L}_{g,n}(H)$ we have the matrices $\overset{I}{B}(i), \overset{I}{A}(i), \overset{I}{M}(j)$ and $\overset{I}{v}$. Hence for each $H$-module $I$, we have an evaluation map $\mathrm{ev}_I : \pi_1^v(\Sigma_{g,n}^{\mathrm{o}}) \to \mathcal{L}_{g,n}(H) \otimes \mathrm{End}_{\mathbb{C}}(I)$ defined by
\[ b_i \mapsto \overset{I}{B}(i), \:\:\: a_i \mapsto \overset{I}{A}(i), \:\:\: m_j \mapsto \overset{I}{M}(j), \:\:\: v \mapsto \overset{I}{v}, \:\:\: \mathrm{ev}_I(xy) = \mathrm{ev}_I(x) \mathrm{ev}_I(y). \]
This observation together with the normalization introduced in section \ref{sectionNorma} will allow us to define the {\em lift} of a simple loop for any $g,n$ and the {\em lift} of a homeomorphism for $n=0$ (see however the discussion in section \ref{CasGeneral}). 

\begin{definition}\label{defLiftLoop}
Let $x \in \pi_1(\Sigma_{g,n}^{\mathrm{o}})$ be a positively oriented simple loop. The lift of $x$ (in the representation $I$) is 
\[ \overset{I}{\widetilde{x}} = \mathrm{ev}_I\!\left( v^{N(x)}x \right). \]
Let $x \in \pi_1(\Sigma_{g,n}^{\mathrm{o}})$ be a negatively oriented simple loop. The lift of $x$ (in the representation $I$) is 
\[ \overset{I}{\widetilde{x}} = \biggl(\overset{I}{\widetilde{x^{-1}}}\biggr)^{-1}. \]
\end{definition}
\noindent In Remark \ref{pourquoiDefLift}, we will see why we must distinguish the positively oriented case from the negatively oriented case in the definition of the lift. Of course, $\overset{I}{\widetilde{b_i}} = \overset{I}{B}(i), \overset{I}{\widetilde{a_i}} = \overset{I}{A}(i), \overset{I}{\widetilde{m_j}} = \overset{I}{M}(j)$ since $N(b_i) = N(a_i) = N(m_j) = 0$. 

\smallskip

\indent Until now, we restrict to $n=0$. For the loops of \eqref{courbesGenHumphries}, we have ($1 \leq i \leq g$):
\begin{equation*}
\begin{split}
& \overset{I}{\widetilde{d_i}} = \mathrm{ev}_I\!\left(v^2 d_i\right) = \overset{I}{v}{^2}\overset{I}{A}(i-1) \overset{I}{B}(i) \overset{I}{A}(i){^{-1}} \overset{I}{B}(i){^{-1}}, \\
&  \overset{I}{\widetilde{e_i}} = \mathrm{ev}_I\!\left(v^{2i} e_i\right) = \overset{I}{v}{^{2i}} \overset{I}{B}(1) \overset{I}{A}(1){^{-1}} \overset{I}{B}(1){^{-1}} \overset{I}{A}(1) \ldots \overset{I}{B}(i-1)\overset{I}{A}(i-1){^{-1}}\overset{I}{B}(i-1){^{-1}}\overset{I}{A}(i-1) \overset{I}{B}(i) \overset{I}{A}(i){^{-1}}\overset{I}{B}(i){^{-1}},\\
&\overset{I}{\widetilde{s_i}} = \mathrm{ev}_I\!\left( v^{2i} s_i \right) = \overset{I}{v}{^{2i}} \overset{I}{B}(1) \overset{I}{A}(1){^{-1}} \overset{I}{B}(1){^{-1}} \overset{I}{A}(1) \ldots \overset{I}{B}(i) \overset{I}{A}(i){^{-1}}\overset{I}{B}(i){^{-1}} \overset{I}{A}(i).
\end{split}
\end{equation*}
We note that these lifts satisfy the $\mathcal{L}_{0,1}(H)$-fusion relation:
\begin{equation}\label{liftHumphriesFusion}
\overset{I \otimes J}{\widetilde{d_i}} = \overset{I}{\widetilde{d_i}} \overset{IJ}{(R')}_{12} \overset{J}{\widetilde{d_i}} (\overset{IJ}{R'})^{-1}_{12}, \:\:\:\:\:\:\:\: \overset{I \otimes J}{\widetilde{e_i}} = \overset{I}{\widetilde{e_i}} \overset{IJ}{(R')}_{12} \overset{J}{\widetilde{e_i}} (\overset{IJ}{R'})^{-1}_{12}, \:\:\:\:\:\:\:\: \overset{I \otimes J}{\widetilde{s_i}} = \overset{I}{\widetilde{s_i}} \overset{IJ}{(R')}_{12} \overset{J}{\widetilde{s_i}} (\overset{IJ}{R'})^{-1}_{12}.
\end{equation}
To check this easily, observe that $\overset{I}{\widetilde{d_i}} = \overset{I}{A}(i-1) \overset{I}{C}_{1,0} \overset{I}{A}(i){^{-1}}$, $\overset{I}{\widetilde{e_i}} = \overset{I}{C}_{i,0} \overset{I}{A}(i){^{-1}}$ and use Lemma \ref{decGauss} and relations \eqref{PresentationLgn} to write the fusion and reorder the matrices, which is a straightforward computation.

\smallskip

Now, we define maps $\widetilde{\tau_{e_2}}, \widetilde{\tau_{b_i}}, \widetilde{\tau_{d_i}} : \mathcal{L}_{g,0}(H) \to \mathcal{L}_{g,0}(H)$ which lift the action \eqref{actionPi1} of the Humphries generators on the fundamental group by the following formulas (recall \eqref{actionPi1V}):
\begin{equation}\label{defLiftHumphries}
\widetilde{\tau_{e_2}}\!\left( \mathrm{ev}_I(u_j) \right) = \mathrm{ev}_I\!\left( \tau_{e_2}^v(u_j) \right), \:\:\:\: \widetilde{\tau_{b_i}}\!\left( \mathrm{ev}_I(u_j) \right) = \mathrm{ev}_I\!\left( \tau_{b_i}^v(u_j) \right), \:\:\:\: \widetilde{\tau_{d_i}}\!\left( \mathrm{ev}_I(u_j) \right) = \mathrm{ev}_I\!\left( \tau_{d_i}^v(u_j) \right) 
\end{equation}
where $u_j$ is $a_j$ or $b_j$. Thanks to \eqref{fvEtNSeMarientBien} and the fact that $N(u_j)=0$, this can also be written as
\[ \widetilde{\tau_{e_2}}(\widetilde{u_j}) = \widetilde{\tau_{e_2}(u_j)}, \:\:\:\:\: \widetilde{\tau_{b_i}}(\widetilde{u_j}) = \widetilde{\tau_{b_i}(u_j)}, \:\:\:\:\: \widetilde{\tau_{d_i}}(\widetilde{u_j}) = \widetilde{\tau_{d_i}(u_j)}. \]
More explicitly:
\begin{equation}\label{courbesDeviennentMatrices}
\begin{split}
& \widetilde{\tau_{e_2}}(\overset{I}{A}(1)) = \overset{I}{\widetilde{e_2}}{^{-1}} \overset{I}{A}(1) \overset{I}{\widetilde{e_2}}, \:\:\:\: \widetilde{\tau_{e_2}}(\overset{I}{B}(1)) = \overset{I}{\widetilde{e_2}}{^{-1}} \overset{I}{B}(1) \overset{I}{\widetilde{e_2}}, \:\:\:\: \widetilde{\tau_{e_2}}(\overset{I}{B}(2)) = \overset{I}{v}\overset{I}{\widetilde{e_2}}{^{-1}} \overset{I}{B}(2), \\
& \widetilde{\tau_{b_i}}(\overset{I}{A}(i)) = \overset{I}{v} \overset{I}{B}(i){^{-1}} \overset{I}{A}(i), \\
& \widetilde{\tau_{a_1}}(\overset{I}{B}(1)) = \overset{I}{v}{^{-1}} \overset{I}{B}(1) \overset{I}{A}(1) \:\:\:\: (\text{recall that } [a_1] = [d_1^{-1}]),\\
& \widetilde{\tau_{d_j}}(\overset{I}{A}(j-1)) = \overset{I}{\widetilde{d_j}}{^{-1}} \overset{I}{A}(j-1) \overset{I}{\widetilde{d_j}}, \:\:\:\: \widetilde{\tau_{d_j}}(\overset{I}{B}(j-1)) = \overset{I}{v}{^{-1}}\overset{I}{B}(j-1) \overset{I}{\widetilde{d_j}}, \:\:\:\: \widetilde{\tau_{d_j}}(\overset{I}{B}(j)) = \overset{I}{v}\overset{I}{\widetilde{d_j}}{^{-1}}\overset{I}{B}(j),
\end{split}
\end{equation}
for $j \geq 2$, and the other matrices are fixed.

\begin{proposition}\label{liftHumphries} 1) The maps $\widetilde{\tau_{e_2}}, \widetilde{\tau_{b_i}}, \widetilde{\tau_{d_i}}$ are automorphisms of $\mathcal{L}_{g,0}(H)$.
\\2) The assignment
\[ \tau_{e_2} \mapsto \widetilde{\tau_{e_2}}, \:\:\: \tau_{b_i} \mapsto \widetilde{\tau_{b_i}}, \:\:\: \tau_{d_i} \mapsto \widetilde{\tau_{d_i}} \]
extends to a morphism of groups $\mathrm{MCG}(\Sigma_g^{\mathrm{o}}) \to \mathrm{Aut}(\mathcal{L}_{g,0}(H))$.
\end{proposition}
\begin{proof}
1) We have to check that these maps are compatible with the defining relations \eqref{PresentationLgn}. This relies on straightforward but tedious computations. For instance, let us show that $\widetilde{\tau_{d_j}}(\overset{I}{B}(j-1))$ satisfies the fusion relation. First, it is easy to establish the following exchange relation:
\[ \overset{IJ}{R'} \, \overset{J}{B}(j-1)_2 \, \overset{IJ}{R} \, (\overset{I}{\widetilde{d_j}})_1 \, \overset{IJ}{R'} = (\overset{I}{\widetilde{d_j}})_1 \, \overset{IJ}{R'} \, \overset{J}{B}(j-1)_2. \]
Hence, using \eqref{liftHumphriesFusion},
\begin{align*}
\overset{I \otimes J}{v^{-1}} \, \overset{I \otimes J}{B}\!(j-1) \, \overset{I \otimes J}{\widetilde{d_j}} &= \overset{I \otimes J}{B}\!(j-1) \, \overset{I \otimes J}{v^{-1}} \, \overset{I \otimes J}{\widetilde{d_j}} = \overset{I \otimes J}{B}\!(j-1) \, \overset{IJ}{R'} \, \overset{IJ}{R} \, \overset{I}{v}{_1^{-1}} \, \overset{J}{v}{_2^{-1}} \, \overset{I \otimes J}{\widetilde{d_j}}\\
&= \overset{I}{v}{_1^{-1}} \, \overset{J}{v}{_2^{-1}} \, \overset{I}{B}(j-1)_1 \, \overset{IJ}{R'} \, \overset{J}{B}(j-1)_2 \, \overset{IJ}{R} \, (\overset{I}{\widetilde{d_j}})_1 \, \overset{IJ}{R'}  \, (\overset{J}{\widetilde{d_j}})_2 \, \overset{IJ}{R'}{^{-1}}\\
&= \biggl(\overset{I}{v}{^{-1}} \, \overset{I}{B}(j-1) \, (\overset{I}{\widetilde{d_j}}) \biggr)_{\!\!1} \, \overset{IJ}{R'} \, \biggl( \overset{J}{v}{^{-1}} \, \overset{J}{B}(j-1) \, (\overset{J}{\widetilde{d_j}}) \biggr)_{\!\!2} \, \overset{IJ}{R'}{^{-1}}.
\end{align*}
This computation reveals the role of the power of $v$ which appears in $\widetilde{\tau_{d_j}}(\overset{I}{B}(j-1))$: it replaces $R'^{-1}$ by $R$ and allows us to apply the previously established exchange relation. We used \eqref{ribbon} and the fact that $\overset{K}{v} \, \overset{K}{U}(k) = \overset{K}{U}(k) \, \overset{K}{v}$ where $U$ is $B$ or $A$. Note that the normalizations by powers of $v$ have no importance when one checks the compatibility with the other defining relations of $\mathcal{L}_{g,0}(H)$, they are only used for the fusion relation.
\\\noindent 2) This is obvious thanks to Proposition \ref{propNormalisation} and \eqref{defLiftHumphries}.
\end{proof}

\begin{definition}\label{liftHomeo}
The lift of an element $f \in \mathrm{MCG}(\Sigma_g^{\mathrm{o}})$, denoted by $\widetilde{f}$, is its image by the morphism of Proposition \ref{liftHumphries}.
\end{definition}
\noindent Due to \eqref{defLiftHumphries}, it holds $\widetilde{f} \circ \mathrm{ev}_I = \mathrm{ev}_I \circ  f^v$. In other words,  $f^v$ and $\widetilde{f}$ are formally identical. Moreover, we have the following lemma which is an expected fact.
\begin{lemma}\label{pasSurprenant}
If $x \in \pi_1\!\left( \Sigma_{g}^{\mathrm{o}} \right)$ is a simple loop it holds
\[ \widetilde{f}\bigl(\overset{I}{\widetilde{x}}\bigr) = \overset{I}{\widetilde{f(x)}}. \]
\end{lemma}
\begin{proof}
If $x$ is positively oriented, so is $f(x)$ and we have
\[ \widetilde{f}\bigl(\overset{I}{\widetilde{x}}\bigr) = \widetilde{f}\bigl( \mathrm{ev}_I( v^{N(x)}x) \bigr) = \mathrm{ev}_I\bigl( f^v(v^{N(x)}x) \bigr) = \mathrm{ev}_I\bigl( v^{N(f(x))} f(x) \bigr) = \overset{I}{\widetilde{f(x)}} \]
thanks to \eqref{propNormalisation} and \eqref{fvEtNSeMarientBien}. If $x$ is negatively oriented, so is $f(x)$ and we have
\[ \widetilde{f}\bigl( \widetilde{x} \bigr) = \widetilde{f}\bigl( \widetilde{x^{-1}}^{-1} \bigr) = \widetilde{f}\bigl( \widetilde{x^{-1}} \bigr)^{-1} = \widetilde{f(x^{-1})}^{-1} = \widetilde{f(x)^{-1}}^{-1} = \widetilde{f(x)}. \qedhere \]
\end{proof}

\indent If $\gamma_1, \gamma_2$ are circles on a surface which have the same topological type\footnote{Two circles $\gamma_1, \gamma_2$ on a surface $S$ are said to have the same topological type if the cut surfaces $S_{\gamma_1}, S_{\gamma_2}$ are the same (up to homeomorphism), see \cite[p. 38]{FM}.}, there exists a homeomorphism $f$ such that $f(\gamma_1)$ is freely homotopic to $\gamma_2$ (see e.g. \cite[Sect. 1.3.1]{FM}). Here we need to consider fixed-point homotopies. We say that two simple loops $x_1, x_2 \in \pi_1\!\left( \Sigma_g^{\mathrm{o}} \right)$ have the same topological type if the circles $[x_1], [x_2]$ have the same topological type.

\begin{lemma}\label{transfoLoop}
Let $x_1, x_2$ be positively oriented simple loops in $\pi_1(\Sigma_g^{\mathrm{o}})$ which have the same topological type, then there exists $f \in \mathrm{MCG}(\Sigma_g^{\mathrm{o}})$ such that $f(x_1) = x_2$ in $\pi_1(\Sigma_g^{\mathrm{o}})$.
\end{lemma}
\begin{proof}
As mentionned, we already know that there exists $\eta \in \mathrm{MCG}(\Sigma_g^{\mathrm{o}})$ such that $\eta(x_1) = x_2' = \alpha^{\varepsilon} x_2^{\pm 1} \alpha^{-\varepsilon}$ in $\pi_1(\Sigma_g^{\mathrm{o}})$ for some loop $\alpha$ and some $\varepsilon \in \{\pm 1\}$. $x'_2$ is positively oriented, non-separating and simple since $x_1$ is, and thus we can assume that $\alpha$ is simple and does not intersect $x_2$ (except at the basepoint). There are six possible configurations for the loops $\alpha$ and $x_2$ in a neighbourhood of the basepoint:

\begin{center}
\begin{tabular}{l l l}
\begin{tikzpicture}
\begin{scope}[yscale=-1,xscale=1]
\draw [line width=0.8pt] (1,1)-- (5,1);
\draw [line width=0.8pt, color=red] (3,1) to [bend left=10] (1.9983888861853318,1.332201709197648);
\draw [line width=0.8pt, color=red, ->] (1.9983888861853318,1.332201709197648) to [bend left=30] (1.7552819764080068,2.26861350982142);
\draw [line width=0.8pt,dash pattern=on 1pt off 1pt, color=red] (1.7552819764080068,2.26861350982142) to [bend left=50] (4.3484223473661405,2.2596095502000373);
\draw [line width=0.8pt, color=red] (3,1) to [bend right=10] (4.096311477967434,1.332201709197648);
\draw [line width=0.8pt, color=red] (4.096311477967434,1.332201709197648) to [bend right=30] (4.3484223473661405,2.2596095502000373);
\draw [line width=0.8pt, color=blue, ->] (3,1) to [bend left=30] (2.4575908268758346,2.103680132957095);
\draw [line width=0.8pt, color=blue] (3,1) to [bend right=30] (3.5290620218203412,2.103680132957095);
\draw [line width=0.8pt,dash pattern=on 1pt off 1pt, color=blue] (2.4575908268758346,2.103680132957095) to [bend left=50] (3.5290620218203412,2.103680132957095);
\draw [fill=black] (3,1) circle (2.5pt);
\draw[color=blue] (3.5290620218203412+0.3,2.103680132957095) node{$x_2$};
\draw[color=red] (4.3484223473661405+0.3,2.2596095502000373) node{$\alpha$};
\draw (0.3, 1) node{1.};
\end{scope}
\end{tikzpicture}
&
~~~~\begin{tikzpicture}
\begin{scope}[yscale=-1,xscale=1]
\draw [line width=0.8pt] (1,1)-- (5,1);
\draw [line width=0.8pt, color=blue] (3,1) to [bend left=10] (1.9983888861853318,1.332201709197648);
\draw [line width=0.8pt, color=blue, ->] (1.9983888861853318,1.332201709197648) to [bend left=30] (1.7552819764080068,2.26861350982142);
\draw [line width=0.8pt,dash pattern=on 1pt off 1pt, color=blue] (1.7552819764080068,2.26861350982142) to [bend left=50] (4.3484223473661405,2.2596095502000373);
\draw [line width=0.8pt, color=blue] (3,1) to [bend right=10] (4.096311477967434,1.332201709197648);
\draw [line width=0.8pt, color=blue] (4.096311477967434,1.332201709197648) to [bend right=30] (4.3484223473661405,2.2596095502000373);
\draw [line width=0.8pt, color=red, ->] (3,1) to [bend left=30] (2.4575908268758346,2.103680132957095);
\draw [line width=0.8pt, color=red] (3,1) to [bend right=30] (3.5290620218203412,2.103680132957095);
\draw [line width=0.8pt,dash pattern=on 1pt off 1pt, color=red] (2.4575908268758346,2.103680132957095) to [bend left=50] (3.5290620218203412,2.103680132957095);
\draw [fill=black] (3,1) circle (2.5pt);
\draw[color=red] (3.5290620218203412+0.3,2.103680132957095) node{$\alpha$};
\draw[color=blue] (4.3484223473661405+0.4,2.2596095502000373) node{$x_2$};
\draw (0.3, 1) node{2.};
\end{scope}
\end{tikzpicture}
&
~~~~\begin{tikzpicture}
\begin{scope}[yscale=-1.15,xscale=1]
\draw [line width=0.8pt] (1,1)-- (5,1);
\draw [line width=0.8pt, color=red, ->] (3,1) to [bend left=30] (1.26,2.18);
\draw [line width=0.8pt, color=red] (3,1) to [bend right=20] (2.38,2.54);
\draw [line width=0.8pt, color=red, dash pattern=on 1pt off 1pt] (1.26,2.18) to [bend left=50] (2.38,2.54);
\draw [line width=0.8pt, color=blue, ->] (3,1) to [bend left=20] (3.54,2.56);
\draw [line width=0.8pt,dash pattern=on 1pt off 1pt, color=blue] (3.54,2.56) to [bend left=50] (4.74,1.94);
\draw [line width=0.8pt, color=blue] (3,1) to [bend right=30] (4.74,1.94);
\draw [fill=black] (3,1) circle (2.5pt);
\draw [color=red] (1,2.18) node{$\alpha$};
\draw[color=blue] (5.1,1.94) node{$x_2$};
\draw (0.3, 1) node{3.};
\end{scope}
\end{tikzpicture}
\\
\begin{tikzpicture}
\begin{scope}[yscale=-1,xscale=1]
\draw [line width=0.8pt] (1,1)-- (5,1);
\draw [line width=0.8pt, color=blue, ->] (3,1) to [bend left=30] (1.26,2.18);
\draw [line width=0.8pt, color=blue] (3,1) to [bend right=20] (2.38,2.54);
\draw [line width=0.8pt, color=blue, dash pattern=on 1pt off 1pt] (1.26,2.18) to [bend left=50] (2.38,2.54);
\draw [line width=0.8pt, color=red, ->] (3,1) to [bend left=20] (3.54,2.56);
\draw [line width=0.8pt,dash pattern=on 1pt off 1pt, color=red] (3.54,2.56) to [bend left=50] (4.74,1.94);
\draw [line width=0.8pt, color=red] (3,1) to [bend right=30] (4.74,1.94);
\draw [fill=black] (3,1) circle (2.5pt);
\draw [color=blue] (1-0.1,2.18) node{$x_2$};
\draw[color=red] (5,1.94) node{$\alpha$};
\draw (0.3, 1) node{4.};
\end{scope}
\end{tikzpicture}
&
~~~~\begin{tikzpicture}
\begin{scope}[yscale=-1,xscale=1]
\draw [line width=0.8pt] (1,1)-- (5,1);
\draw [line width=0.8pt, color=blue, ->] (3,1) to [bend left=50] (1.5121750666306817,2.3046293483069493+0.2); 
\draw [line width=0.8pt,dash pattern=on 1pt off 1pt, color=blue] (1.5121750666306817,2.3046293483069493+0.2) to [bend left=20] (2.3315353921764808,2.5837520965698046+0.2);
\draw [line width=0.8pt,dash pattern=on 1pt off 1pt, color=blue] (2.3315353921764808,2.5837520965698046+0.2) to [bend left=20] (3.060856121508456,2.23259767133589+0.2);
\draw [line width=0.8pt, color=blue] (3.060856121508456,2.23259767133589+0.2) to [bend left=30] (3,1);
\draw [line width=0.8pt, color=red, ->] (3,1) to [bend left=30] (2.3585472710406283,2.0705263981510065+0.2);
\draw [line width=0.8pt,dash pattern=on 1pt off 1pt, color=red] (2.3585472710406283,2.0705263981510065+0.2) to [bend left=20] (3.141891758100898,2.5207243792201277+0.2);
\draw [line width=0.8pt,dash pattern=on 1pt off 1pt, color=red] (3.141891758100898,2.5207243792201277+0.2) to [bend left=30] (3.943244164403932+0.3,1.8994511653447406+0.2);
\draw [line width=0.8pt, color=red] (3.943244164403932+0.3,1.8994511653447406+0.2) to [bend left=40] (3,1);
\draw [fill=black] (3,1) circle (2.5pt);
\draw[color=blue] (1.5121750666306817-0.4,2.3046293483069493+0.2) node{$x_2$};
\draw[color=red] (3.943244164403932+0.6,1.8994511653447406+0.2) node{$\alpha$};
\draw (0.3, 1) node{5.};
\end{scope}
\end{tikzpicture}
&
~~~~\begin{tikzpicture}
\begin{scope}[yscale=-1,xscale=1]
\draw [line width=0.8pt] (1,1)-- (5,1);
\draw [line width=0.8pt, color=red, ->] (3,1) to [bend left=50] (1.5121750666306817,2.3046293483069493+0.2); 
\draw [line width=0.8pt,dash pattern=on 1pt off 1pt, color=red] (1.5121750666306817,2.3046293483069493+0.2) to [bend left=20] (2.3315353921764808,2.5837520965698046+0.2);
\draw [line width=0.8pt,dash pattern=on 1pt off 1pt, color=red] (2.3315353921764808,2.5837520965698046+0.2) to [bend left=20] (3.060856121508456,2.23259767133589+0.2);
\draw [line width=0.8pt, color=red] (3.060856121508456,2.23259767133589+0.2) to [bend left=30] (3,1);
\draw [line width=0.8pt, color=blue, ->] (3,1) to [bend left=30] (2.3585472710406283,2.0705263981510065+0.2);
\draw [line width=0.8pt,dash pattern=on 1pt off 1pt, color=blue] (2.3585472710406283,2.0705263981510065+0.2) to [bend left=20] (3.141891758100898,2.5207243792201277+0.2);
\draw [line width=0.8pt,dash pattern=on 1pt off 1pt, color=blue] (3.141891758100898,2.5207243792201277+0.2) to [bend left=30] (3.943244164403932+0.3,1.8994511653447406+0.2);
\draw [line width=0.8pt, color=blue] (3.943244164403932+0.3,1.8994511653447406+0.2) to [bend left=40] (3,1);
\draw [fill=black] (3,1) circle (2.5pt);
\draw[color=red] (1.5121750666306817-0.4,2.3046293483069493+0.2) node{$\alpha$};
\draw[color=blue] (3.943244164403932+0.6,1.8994511653447406+0.2) node{$x_2$};
\draw (0.3, 1) node{6.};
\end{scope}
\end{tikzpicture}
\end{tabular}
\end{center}
In case 1, $x_2' = \alpha x_2 \alpha^{-1}$, and then $\tau_{\alpha}(x'_2) = \alpha^{-1} x'_2 \alpha = x_2$. Case 2 is impossible because none of the four possible loops $\alpha^{\varepsilon} x_2^{\pm 1} \alpha^{-\varepsilon}$ is simple. In case 3, $x'_2 = \alpha x_2 \alpha^{-1}$. For $\beta = \alpha x_2$, we have $\tau_{\beta}(\alpha) = \beta^{-1}\alpha \beta$, $\tau_{\beta}(x_2) = \beta^{-1} x_2 \beta$, and thus $\tau_{\beta}(x'_2) = x_2$. In case 4, $x'_2 = \alpha^{-1} x_2 \alpha$. For $\delta = x_2 \alpha$, we get similarly to case 3 that $\tau_{\delta}^{-1}(x'_2) = \delta x'_2 \delta^{-1} = x_2$. In case 5, $x'_2 = \alpha^{-1} x_2^{-1} \alpha$. Observe that $\tau_{\alpha}(x_2) = x_2 \alpha$, $\tau_{x_2}(\alpha) = x_2^{-1}\alpha$, and then
\[ \tau_{\alpha}^{-1} \tau_{x_2}^{-2} \tau_{\alpha}^{-1}(\alpha^{-1} x^{-1}_2 \alpha) = \tau_{\alpha}^{-1} \tau_{x_2}^{-2}(x^{-1}_2 \alpha) = \tau_{\alpha}^{-1}(x_2 \alpha) =  x_2.\]
In case 6, $x'_2 = \alpha x_2^{-1} \alpha^{-1}$, and we get similarly to case 5 that $\tau_{\alpha} \tau_{x_2}^{2} \tau_{\alpha}(x'_2) = x_2$.
\end{proof}

\begin{exemple}\label{exempleTransfoLoopA1}
 We have
\[
\tau_{b_i} \tau_{a_i}(b_i) = a_i, \:\:\:\:\:
\tau_{d_i}^{-1} \tau_{b_{i-1}}^{-1}(d_i) = b_{i-1}, \:\:\:\:\:
\tau_{y_i}^{-1} \tau_{a_i}^{-1} \tau_{b_{i-1}}^{-1} \tau_{y_i}^{-1}(a_i) = b_{i-1}, \:\:\:\:\:
\tau_{y_2}^{-1}\tau_{b_1}^{-1}\tau_{e_2}\tau_{y_2}(e_2) = b_1
\]
where $y_i=a_{i-1}b_i$. This allows to transform any of the loops $a_i, b_i, d_i, e_2$ into $a_1$.
\finEx
\end{exemple}

\smallskip

\begin{remark}\label{remarkLiftArticle}
In \cite{Fai18c}, Lemma \ref{transfoLoop} was the starting point to define the lift of non-separating positively oriented simple loops. More precisely, we first defined the lifts of the Humphries generators and checked that they satisfy the Wajnryb relations. This defines the lift of every $f \in \mathrm{MCG}\!\left( \Sigma_g^{\mathrm{o}} \right)$. Then we declared that the lift of $a_1$ (in the representation $I$) is $\overset{I}{A}(1)$ and that the lift of a non-separating positively oriented simple loop $x \in \pi_1(\Sigma_g^{\mathrm{o}})$ is $\widetilde{f}(\overset{I}{A}(1))$ where $f$ is such that $x = f(a_1)$. For this, it was necessary to show that $f(a_1) = g(a_1)$ implies $\widetilde{f}(\overset{I}{A}(1)) = \widetilde{g}(\overset{I}{A}(1))$, which was done in \cite[Lemma 5.5]{Fai18c}. This had the advantage to be shorter than the construction presented here since it does not require to define the normalization $N$. However, it is less general because it is not adapted to non-separating loops.
\finEx
\end{remark}

\smallskip

\indent Lemma \ref{transfoLoop} has the following important consequence.

\begin{proposition}\label{propFusionCourbeSimple}
Let $x \in \pi_1(\Sigma_g^{\mathrm{o}})$ be a positively oriented simple loop. Then the lift of $x$ satisfies the fusion relation of $\mathcal{L}_{0,1}(H)$:
\[ \overset{I \otimes J}{\widetilde{x}} = \overset{I}{\widetilde{x}}_1 \, \overset{IJ}{R'} \, \overset{J}{\widetilde{x}}_2 \, \overset{IJ}{R'}{^{-1}}. \]
It follows that there exists a morphism of $H$-module-algebras $j_{\widetilde{x}} : \mathcal{L}_{0,1}(H) \to \mathcal{L}_{g,0}(H)$ given by $\overset{I}{M} \mapsto \overset{I}{\widetilde{x}}$.
\end{proposition}
\begin{proof}
There are $g+1$ possible topological types for loops in $\Sigma_g^{\mathrm{o}}$. The more simple positively oriented loops representing each topological type are $a_1, s_1, \ldots, s_g$ (whose cut surfaces are respectively $\Sigma_{g,3}, \Sigma_{1,1} \sqcup \Sigma_{g-1,2}, \ldots, \Sigma_{g, 1} \sqcup \Sigma_{0,2}$, see Figures \ref{figureCourbesCanoniques} and \ref{courbesSeparantes}). These particular loops satisfy the fusion relation. This is obviously true for $\widetilde{a_1} = A(1)$ . For $\widetilde{s_i}$, observe that $\overset{I}{\widetilde{s_i}} = j(\overset{I}{C}_{i,0})$, where $\overset{I}{C}_{i,0}$ is defined is defined in section \ref{RepInvariants} and $j : \mathcal{L}_{i,0}(H) \to \mathcal{L}_{g,0}(H)$ is the obvious embedding. Since by Proposition \ref{propMatriceCgn}, $C_{i,0}$ satifies the fusion relation, so does $\widetilde{s_i}$. By Lemma \ref{transfoLoop}, there exists $f\in \mathrm{MCG}\!\left( \Sigma_g^{\mathrm{o}} \right)$ such that $f(x)$ is $a_1$ or $s_1 \ldots$ or $s_g$ and hence by Lemma \ref{pasSurprenant}: 
\[ \overset{I \otimes J}{\widetilde{x}} = \overset{I \otimes J}{\widetilde{f^{-1}(f(x))}} = \widetilde{f}^{-1}\biggl( \overset{I \otimes J}{\widetilde{f(x)}} \biggr) = \widetilde{f}^{-1}\biggl( \overset{I}{\widetilde{f(x)}}_1 \, \overset{IJ}{R'} \, \overset{J}{\widetilde{f(x)}}_2 \, \overset{IJ}{R'}{^{-1}} \biggr) = \overset{I}{\widetilde{x}}_1 \, \overset{IJ}{R'} \, \overset{J}{\widetilde{x}}_2 \, \overset{IJ}{R'}{^{-1}}, \]
as desired.
\end{proof}

\smallskip

\indent Recall that $\mathcal{L}_{g,0}(H) \cong \mathrm{End}_{\mathbb{C}}\!\left((H^*)^{\otimes g}\right)$ is a matrix algebra. By the Skolem-Noether theorem, every automorphism of $\mathcal{L}_{g,0}(H)$ is inner. Hence to each $f \in \mathrm{MCG}(\Sigma_g^{\mathrm{o}})$ is associated an element $\widehat{f} \in \mathcal{L}_{g,0}(H)$, unique up to scalar, such that
\begin{equation}\label{conjugaison}
\forall\, x \in \mathcal{L}_{g,0}(H), \:\:\: \widetilde{f}(x) = \widehat{f}x\widehat{f}^{-1}.
\end{equation}
\indent We will determine the elements $\widehat{\tau_{\gamma}}$ associated to Dehn twists about non-separating circles and use this to show that $\widehat{f} \in \mathcal{L}_{g,0}^{\mathrm{inv}}(H)$.

\begin{lemma}\label{lemmaA1}
We have $\widehat{\tau_{a_1}} = v_{A(1)}^{-1}$. In other words:
\[ \forall\, x \in \mathcal{L}_{g,0}(H), \:\: \widetilde{\tau_{a_1}}(x) = v_{A(1)}^{-1} \, x \, v_{A(1)}. \]
\end{lemma}
\begin{proof}
We have $v_{A(1)}^{-1} \overset{I}{A}(1) = \overset{I}{A}(1) v_{A(1)}^{-1} = \widetilde{\tau_{a_1}}(\overset{I}{A}(1))v_{A(1)}^{-1}$. Indeed, since $v^{-1}$ is central in $H$, $v_{A(1)}^{-1}$ is central in the subalgebra generated by the coefficients of the matrices $\overset{I}{A}(1)$. Next, let $j_1 : \mathcal{H}(\mathcal{O}(H)) \to \mathcal{H}(\mathcal{O}(H))^{\otimes g}$ be the canonical embedding on the first copy. Observe that for all $x \in H$, $\Psi_{g,0}(x_{A(1)}) = j_1(x)$. Then:
\begin{align*}
\Psi_{g,0}\bigl(v_{A(1)}^{-1} \overset{I}{B}(1)\bigr) &= j_1\bigl(v^{-1} \overset{I}{L}{^{(+)}} \, \overset{I}{T} \, \overset{I}{L}{^{(-)-1}}\bigr) = j_1\bigl(\overset{I}{L}{^{(+)}} \, \overset{I}{T} \, \overset{I}{(v'^{-1})} \, v''^{-1} \,\overset{I}{L}{^{(-)-1}}\bigr)\\
& = j_1\bigl(\overset{I}{L}{^{(+)}} \, \overset{I}{T} \, \overset{I}{v}{^{-1}}\overset{I}{b_i} \overset{I}{a_j} a_i b_j \, v^{-1} \,\overset{I}{L}{^{(-)-1}}\bigr) = j_1\bigl(\overset{I}{v}{^{-1}} \overset{I}{L}{^{(+)}} \, \overset{I}{T} \, \overset{I}{L}{^{(-)-1}} \, \overset{I}{L}{^{(+)}} \,\overset{I}{L}{^{(-)-1}} v^{-1}\bigr)\\
& = \Psi_{g,0}\bigl(\overset{I}{v}{^{-1}} \overset{I}{B}(1) \overset{I}{A}(1) v_{A(1)}^{-1}\bigr) =  \Psi_{g,0}\bigl( \widetilde{\tau_{a_1}}(\overset{I}{B}(1)) v_{A(1)}^{-1}\bigr).
\end{align*}
We used the exchange relation \eqref{relDefHeisenberg} of $\mathcal{H}(\mathcal{O}(H))$ together with \eqref{ribbon} and the definition of the matrices $\overset{I}{L}{^{(\pm)}}$. Finally, recall the matrices \eqref{matricesAlekseev} which occur in the definition of the Alekseev isomorphism. The same argument as in the proof of Lemma \ref{expressionM} shows that
\[ \Psi_{1,0}^{\otimes g}(\overset{I}{\Lambda}_i) = \overset{I}{S^{-1}(b_{\ell})} \: \widetilde{a_{\ell}^{(2i -1)}} a_{\ell}^{(2i)} \otimes \ldots \otimes \widetilde{a_{\ell}^{(1)}}b_{\ell}^{(2)}. \]
From this we see that $j_1(v^{-1})$ commutes with $\Psi_{1,0}^{\otimes g}(\overset{I}{\Lambda}_i)$. Eventually it follows that $\Psi_{g,0}(v^{-1}_{A(1)})$ commutes with $\Psi_{g,0}\bigl(\overset{I}{U}(i)\bigr) = \Psi_{g,0}\biggl(\widetilde{\tau_{a_1}}\bigl(\overset{I}{U}(i)\bigr)\biggr)$, where $U$ is $A$ or $B$.
\end{proof}

\noindent Recall the notation \eqref{notationEmbed}. If $\gamma$ is a simple loop, $\widetilde{\gamma}$ satisfies the fusion relation and thus $(v^{-1})_{\widetilde{\gamma}} = (v_{\widetilde{\gamma}})^{-1}$.
\begin{proposition}\label{propDehnTwist}
For any non-separating circle $\gamma$ on $\Sigma_g^{\mathrm{o}}$, we have $\widehat{\tau_{\gamma}} = v_{\widetilde{\gamma}}^{-1}$. In other words:
\[ \forall\, x \in \mathcal{L}_{g,0}(H), \:\: \widetilde{\tau_{\gamma}}(x) = v_{\widetilde{\gamma}}^{-1} \, x \, v_{\widetilde{\gamma}}. \]
If $\gamma, \delta \in \pi_1(\Sigma_g^{\mathrm{o}})$ are positively oriented non-separating simple loops such that $[\gamma] = [\delta]$, then $v_{\widetilde{\gamma}}$ is proportional to $v_{\widetilde{\delta}}$.
\end{proposition}
\begin{proof}
We represent the circle $[\gamma]$ by a positively oriented, non-separating simple loop $\gamma \in \pi_1(\Sigma_g^{\mathrm{o}})$. Let $f \in \mathrm{MCG}(\Sigma_g^{\mathrm{o}})$ be such that $f(a_1) = \gamma$, then
\[ \widetilde{\tau_{\gamma}} = \widetilde{\tau_{f(a_1)}} = \widetilde{f \tau_{a_1}f^{-1}} = \widetilde{f} \, \widetilde{\tau_{a_1}} \, \widetilde{f}^{-1}. \]
Hence, by Lemma \ref{lemmaA1},
\[ \forall\, x \in \mathcal{L}_{g,0}(H), \:\:\: \widetilde{\tau_{\gamma}}\!\left(\widetilde{f}(x)\right) = \widetilde{f}\!\left(\widetilde{\tau_{a_1}}(x)\right) = \widetilde{f}\!\left(v_{A(1)}^{-1} x v_{A(1)}\right) = v_{\widetilde{\gamma}}^{-1} \widetilde{f}(x) v_{\widetilde{\gamma}}. \]
Replacing $x$ by $\widetilde{f}^{-1}(x)$, we get the result. The second claim follows from a similar reasoning together with the fact that $\tau_{\gamma}$ depends only of the free homotopy class of $\gamma$.
\end{proof}

\noindent An analogous result in the modular setting has been given in \cite[eq (9.7)]{AS}. The notation $v_{\widetilde{\gamma}}^{-1}$ does not appear in their work; instead, they express this element as a linear combination of traces which form a basis in the modular case only.

\begin{corollary}\label{coroFTildeInv}
For all $f \in \mathrm{MCG}(\Sigma_g^{\mathrm{o}})$, it holds $\widehat{f} \in \mathcal{L}_{g,0}^{\mathrm{inv}}(H)$.
\end{corollary}
\begin{proof}
Let $\gamma$ be a positively oriented, non-separating simple loop. Then $\widetilde{\gamma}$ satisfies the fusion relation of $\mathcal{L}_{0,1}(H)$, and thus $j_{\widetilde{\gamma}}$ is a morphism of $H$-module-algebras (Lemma \ref{injectionFusion}). Hence, since $v^{-1} \in \mathcal{Z}(H) = \mathcal{L}_{0,1}^{\mathrm{inv}}(H)$, we have $v_{\widetilde{\gamma}}^{-1} \in \mathcal{L}_{g,0}^{\mathrm{inv}}(H)$. In particular, the statement is true for the Humphries generators thanks to Proposition \ref{propDehnTwist} and thus for any $f$.
\end{proof}

\subsection{Representation of the mapping class group}\label{sectionRepMCG}
\indent The only additional fact needed is the following lemma.
\begin{lemma}\label{lemmavA1Moins1}
It holds: $v_{A(g)}^{-1} = v_{A(g)^{-1}}^{-1}$.
\end{lemma}
\begin{proof}
Denote as usual $X_i \otimes Y_i = RR'$, $\overline{X}_i \otimes \overline{Y}_i = (RR')^{-1}$ and let $\mu^l$ be the left integral on $H$ (unique up to scalar). We have:
\[ \mu^l(vX_i)Y_i = \mu^l(v)v^{-1} = \mu^l(v\overline{X}_i)\overline{Y}_i. \]
The first equality is shown in the proof of Proposition \ref{propVIntegrale} while the second is easy using \eqref{ribbon} and \eqref{integrale}. Let us write $\mu^l(v)^{-1}\mu^l(v?) = \sum_{i,j,I} c_{I,i}^j \overset{I}{T}{^i_j}$ with $c_{I,i}^j \in \mathbb{C}$. Then, using the identification $\overset{I}{M} = (\overset{I}{X_i})Y_i$ between $\mathcal{L}_{0,1}(H)$ and $H$, the fact that $\overset{I}{M}{^{-1}} = (\overset{I}{\overline{X}_i})\overline{Y}_i$ and the equalities above, we get
\[ v_{A(g)}^{-1} = j_{A(g)}\!\left( \sum_{i,j,I} c_{I,i}^j \overset{I}{M}{^i_j} \right) = j_{A(g)}\!\left( \sum_{i,j,I} c_{I,i}^j (\overset{I}{M}{^{-1}})^i_j \right) = j_{A(g)^{-1}}\!\left( \sum_{i,j,I} c_{I,i}^j \overset{I}{M}{^i_j} \right) = v_{A(g)^{-1}}^{-1} \]
where the morphisms $j_{\bullet}$ are defined at the end of subsection \ref{defLgn}. We used that $j_{A(g)}$ is a morphism of algebras (see Lemma \ref{injectionFusion}).
\end{proof}

It is clear that the lemma holds for the lift of any positively oriented, non-separating simple loop, but we do not need this.
\smallskip\\  
\indent Recall that we have a representation of $\mathcal{L}_{g,0}(H)$ on $(H^*)^{\otimes g}$, let us denote it $\rho$. We also have the associated representation of $\mathcal{L}^{\mathrm{inv}}_{g,0}(H)$ on $\mathrm{Inv}\!\left((H^*)^{\otimes g}\right)$, let us denote it $\rho_{\mathrm{inv}}$. Also recall that the elements $\widehat{f}$ are defined in \eqref{conjugaison}. We can now state the representation of the mapping class groups $\mathrm{MCG}(\Sigma_g^{\mathrm{o}})$ and $\mathrm{MCG}(\Sigma_g)$. An analogous result was announced in \cite{AS} under the assumption that the gauge algebra is modular.

\begin{theorem}\label{thmRepMCG}
1) The map
\[ \fonc{\theta_g^{\mathrm{o}}}{\mathrm{MCG}(\Sigma_g^{\mathrm{o}})}{\mathrm{GL}\bigl((H^*)^{\otimes g}\bigr)}{f}{\rho(\widehat{f})} \]
is a projective representation.
\\2) The map
\[ \fonc{\theta_g}{\mathrm{MCG}(\Sigma_g)}{\mathrm{GL}\bigl(\mathrm{Inv}\bigl((H^*)^{\otimes g}\bigr)\bigr)}{f}{\rho_{\mathrm{inv}}(\widehat{f})}\]
is a projective representation.
\end{theorem}
\begin{proof}
1) This is an immediate consequence of Proposition \ref{liftHumphries}.
\\2) We must show that the hyperelliptic relation \eqref{hyperelliptic} is projectively satisfied. The word $w$ can be constructed as follows: take $f \in \mathrm{MCG}(\Sigma_g^{\mathrm{o}})$ such that $f(a_1) = a_g$ and express it as a word in the Humphries generators $f = \tau_{\gamma_1} \ldots \tau_{\gamma_n}$. Then $\tau_{a_g} = f \tau_{a_1} f^{-1}$ and $w = \tau_{\gamma_1} \ldots \tau_{\gamma_n} \tau_{a_1} \tau_{\gamma_n}^{-1} \ldots \tau_{\gamma_1}^{-1}$. The automorphism $\widetilde{\tau_{a_g}}$ is implemented by conjugation by $\widehat{f} v_{A(1)}^{-1} \widehat{f}^{-1}$ and also by conjugation by $v_{A(g)}^{-1}$ (Proposition \ref{propDehnTwist}). Hence, $\widehat{f} v_{A(1)}^{-1} \widehat{f}^{-1} \sim v_{A(g)}^{-1}$, where $\sim$ means proportional. Now, let $h = \tau_{b_g} \tau_{d_g} \ldots \tau_{b_1} \tau_{d_1} \tau_{d_1} \tau_{b_1} \ldots \tau_{d_g} \tau_{b_g}$. A computation gives $\widetilde{h}(\overset{I}{A}(g)) = \overset{I}{A}(g){^{-1}} \overset{I}{C}_{g,0}$. Thus
\[ \widehat{h} \widehat{f} v_{A(1)}^{-1} \widehat{f}^{-1} \widehat{h}^{-1} \sim \widehat{h} v_{A(g)}^{-1} \widehat{h}^{-1} = \widetilde{h}(v_{A(g)}^{-1}) = v_{A(g)^{-1}C_{g,0}}^{-1}. \]
By definition of $\mathrm{Inv}\!\left((H^*)^{\otimes g}\right)$ and Lemma \ref{lemmavA1Moins1}, we have 
\[ \rho_{\mathrm{inv}}(v_{A(g)^{-1}C_{g,0}}^{-1}) = \rho_{\mathrm{inv}}(v_{A(g)^{-1}}^{-1}) = \rho_{\mathrm{inv}}(v_{A(g)}^{-1}). \]
It follows that $\rho_{\mathrm{inv}}\!\left(\widehat{h} \left(\widehat{f} v_{A(1)}^{-1} \widehat{f}^{-1}\right) \widehat{h}^{-1}\right) \sim \rho_{\mathrm{inv}}\!\left(\widehat{f} v_{A(1)}^{-1} \widehat{f}^{-1}\right)$. This shows that the map is well-defined since $\mathrm{MCG}(\Sigma_g)$ is the quotient of $\mathrm{MCG}(\Sigma_g^{\mathrm{o}})$ by the hyperelliptic relation and that it is a projective representation.
\end{proof}

\indent Note that, as in the case of the torus (see discussion at the beginning of section \ref{sectionProjRepSL2Z}), the representation $\mathrm{Inv}\bigl( (H^*)^{\otimes g} \bigr)$ allowed us to glue back the disc $D$, and hence to obtain a representation of $\mathrm{MCG}(\Sigma_g)$ (and not just of $\mathrm{MCG}(\Sigma_g^{\mathrm{o}})$).

\begin{remark}
I think undoubtedly that, when $H$ is modular, the projective representation of Theorem \ref{thmRepMCG} is equivalent to the one of \cite[Th. 28]{AS}. Let $\mathcal{M}_{g,0}^{\mathrm{AS}}(H)$ be the moduli algebra of \cite{AS} (whose construction holds under the assumption that $H$ is modular). In \cite[Th. 28]{AS}, each Dehn twist is mapped to an element of $\mathcal{M}^{\mathrm{AS}}_{g,0}(H)$: $\fleche{\mathrm{MCG}(\Sigma_g)}{\mathcal{M}^{\mathrm{AS}}_{g,0}(H)}{\tau_{\gamma}}{\hat h(\gamma)}$ and these elements $\hat h(\gamma)$ (which implement the lift of $\tau_{\gamma}$ by conjugation in the moduli algebra) are claimed to satisfy the relations of $\mathrm{MCG}(\Sigma_g)$ up to scalar. Recall that in \eqref{defModuli} we proposed to take the image of the representation $\rho_{\mathrm{inv}}$ as a generalization $\mathcal{M}^{\mathrm{gen}}_{g,0}(H)$ of the moduli algebra $\mathcal{M}^{\mathrm{AS}}_{g,0}(H)$ in the non-modular setting. With this definition, we have a map $\fleche{\mathrm{MCG}(\Sigma_g)}{\mathcal{M}^{\mathrm{gen}}_{g,0}(H)}{\tau_{\gamma}}{\rho_{\mathrm{inv}}(v_{\widetilde{\gamma}}^{-1})}$ which satisfies the relations of $\mathrm{MCG}(\Sigma_g)$ up to scalar. If the algebras $\mathcal{M}^{\mathrm{AS}}_{g,0}(H)$ and $\mathcal{M}^{\mathrm{gen}}_{g,0}(H)$ are isomorphic when $H$ is modular, then the representations of $\mathrm{MCG}(\Sigma_g)$ are equivalent in this case. This remark can of course be generalized to any $\Sigma_{g,n}$.
\finEx
\end{remark}

\subsection{Discussion for the case $n > 0$}\label{CasGeneral}
\indent Let us consider the general case $n>0$, see Figure \ref{surfaceAvecCourbes}. Denote $\Sigma_{g,n}^{\mathrm{o}} = \Sigma_{g,n} \setminus D$, where $D$ is an embedded open disk.

\smallskip

\indent The only difference is that, in general, $\mathcal{L}_{g,n}(H)$ is not a matrix algebra and we cannot claim directly the existence and unicity up to scalar of the elements $\widehat{f}$. Nevertheless, we now propose (without proofs) a program to extend the previous construction which should not be difficult to apply.

\smallskip

\indent The first task is to define the lifting homomorphism $\mathrm{MCG}(\Sigma_{g,n}^{\mathrm{o}}) \to \mathrm{Aut}\bigl(\mathcal{L}_{g,n}(H)\bigr)$ generalizing Definition \ref{liftHomeo}:

\smallskip

\indent \textbullet ~ Consider a generating set $g_1, \ldots, g_k$ of $\mathrm{MCG}(\Sigma^{\mathrm{o}}_{g,n})$ and compute the action of these generators on $\pi_1(\Sigma_{g,n}^{\mathrm{o}})$.

\smallskip

\indent \textbullet ~ As in \eqref{actionPi1V}, define the morphisms $g_i^v \in \mathrm{Aut}\bigl( \pi_1^v(\Sigma_{g,n}^{\mathrm{o}}) \bigr)$ by the formula
\[ \forall x \in \{ b_1, a_1, \ldots, b_g, a_g, m_{g+1}, \ldots, m_{g+n} \}, \:\:\: g_i^v(x) = v^{N(g_i(x))}g_i(x). \]
It is clear that Proposition \ref{propNormalisation} is not at all specific to the case $n=0$ and remains true for any $n$. In particular, the assignment $g_i \mapsto g_i^v$ extends to a morphism of groups $\mathrm{MCG}(\Sigma_{g,n}^{\mathrm{o}}) \to \mathrm{Aut}\bigl( \pi_1^v(\Sigma_{g,n}^{\mathrm{o}}) \bigr)$.

\smallskip

\indent  \textbullet ~ The lift of a simple loop is defined for any $g,n$, see Definition \ref{defLiftLoop}. Using this, define the lifts $\widetilde{g_1}, \ldots, \widetilde{g_k}$ of the generators $g_1, \ldots, g_k$ by the same formula as in \ref{defLiftHumphries}:
\[ \forall x \in \{ b_1, a_1, \ldots, b_g, a_g, m_{g+1}, \ldots, m_{g+n} \}, \:\:\: \widetilde{g_i}\bigl( \mathrm{ev}_I(x) \bigr) = \mathrm{ev}_I\bigl( g_i^v(x) \bigr) \]
which can also be written as
\[ \forall x \in \{ b_1, a_1, \ldots, b_g, a_g, m_{g+1}, \ldots, m_{g+n} \}, \:\:\: \widetilde{g_i}\bigl(\overset{I}{\widetilde{x}}\bigr) = \overset{I}{\widetilde{g_i(x)}}. \]
Then one must show that each $\widetilde{g_i}$ preserves the defining relations of $\mathcal{L}_{g,n}(H)$, and it will follow that the assignment $g_i \mapsto \widetilde{g_i}$ extends to a morphism of groups $\mathrm{MCG}(\Sigma_{g,n}^{\mathrm{o}}) \to \mathrm{Aut}\bigl(\mathcal{L}_{g,n}(H)\bigr)$, which to a mapping class $f$ associates its {\em lift} $\widetilde{f}$.

\medskip

\indent Now, we will associate a (not unique) element $\widehat{f}$ to each mapping class $f$ such that $\widetilde{f}(x) = \widehat{f} x \widehat{f}^{-1}$ for all $x \in \mathcal{L}_{g,n}(H)$. We assume $g > 1$ because in this case there is a  generating set of $\mathrm{MCG}(\Sigma^{\mathrm{o}}_{g,n})$ which consists only of Dehn twists, see \cite[Figure 4.10]{FM}; hence $g_i = \tau_{c_i}$ for each $i$, where $c_i \in \pi_1(\Sigma_{g,n}^{\mathrm{o}})$ is a non-separating, positively oriented simple loop.

\smallskip

\indent \textbullet ~ Lemma \ref{transfoLoop} still holds. In particular, since the loop $c_i$ is non-separating for each $i$, it has the same topological type than $a_1$ and there exists $f_i \in \mathrm{MCG}(\Sigma^{\mathrm{o}}_{g,n})$ such that $f_i(a_1) = c_i$. Note that it is possible to find the $f_i$'s explicitly as in Example \ref{exempleTransfoLoopA1}, without invoking Lemma \ref{transfoLoop}, and this is sufficient for the sequel.

\smallskip

\indent \textbullet ~ Lemma \ref{lemmaA1} remains true for $n > 0$. By reproducing the proof of Proposition \ref{propDehnTwist} with the $f_i$'s, we get that $\widetilde{\tau}_{c_i}(x) = v_{\widetilde{c_i}}^{-1} x v_{\widetilde{c_i}}$. Since the $\tau_{c_i}$ are a generating set, it follows that for each $f \in \mathrm{MCG}(\Sigma_{g,n}^{\mathrm{o}})$, there exists an invertible element $\widehat{f} \in \mathcal{L}^{\mathrm{inv}}_{g,n}(H)$ such that $\widetilde{f}(x) = \widehat{f} x \widehat{f}^{-1}$. The elements $\widehat{f}$ are unique only up to an invertible central invariant element of $\mathcal{L}_{g,n}(H)$; in particular $\widehat{fg} = z \widehat{f}\widehat{g}$ for some $z \in \mathcal{Z}\bigl( \mathcal{L}_{g,n}(H) \bigr)$.

\smallskip

\indent \textbullet ~ Since $\mathcal{Z}\bigl(\mathcal{H}({\mathcal{O}}(H))^{\otimes g}\bigr) \cong \mathbb{C}$, we have for all $z \in \mathcal{Z}\!\left(\mathcal{L}_{g,n}(H)\right)$: 
\begin{equation}\label{actionCentralRepLgn}
\Psi_{g,n}(z) = 1 \otimes \ldots \otimes 1 \otimes z_1 \otimes \ldots \otimes z_n \in \mathcal{H}({\mathcal{O}}(H))^{\otimes g} \otimes H^{\otimes n}
\end{equation}
with $z_i \in \mathcal{Z}(H)$. Let 
\[ \rho^{g, S_1, \ldots, S_n} : \mathcal{L}_{g,n}(H) \to \mathrm{End}_{\mathbb{C}}\bigl( V(g, S_1, \ldots, S_n) \bigr) \]
be a representation, with $V(g, S_1, \ldots, S_n) = (H^*)^{\otimes g} \otimes S_1 \otimes \ldots \otimes S_n$ and where $S_1, \ldots, S_n$ are {\em simple} $H$-modules. Due to \eqref{actionCentralRepLgn} and to Schur lemma, we see that the representation of $z \in \mathcal{Z}\bigl(\mathcal{L}_{g,n}(H)\bigr)$ on $V(g, S_1, \ldots, S_n)$ is a scalar $\lambda_z \mathrm{id}_{V(g, S_1, \ldots, S_n)}$. Hence, the element $\rho^{g, S_1, \ldots, S_n}(\widehat{f})$ is unique up to scalar. Define
\[ \fonc{ \theta_{g, S_1, \ldots, S_n}^{\mathrm{o}} }{ \mathrm{MCG}(\Sigma_{g,n}^{\mathrm{o}}) }{ \mathrm{GL}\bigl( V(g, S_1, \ldots, S_n) \bigr) }{ f } { \rho^{g, S_1, \ldots, S_n}(\widehat{f}) }. \]
This is a projective representation:
\[ \theta_{g, S_1, \ldots, S_n}^{\mathrm{o}}(fg) \sim \rho^{g, S_1, \ldots, S_n}(\widehat{fg}) = \rho^{g, S_1, \ldots, S_n}(z \widehat{f}\widehat{g}) \sim  \rho^{g, S_1, \ldots, S_n}(\widehat{f})  \rho^{g, S_1, \ldots, S_n}(\widehat{g}) \sim \theta_{g, S_1, \ldots, S_n}^{\mathrm{o}}(f) \theta_{g, S_1, \ldots, S_n}^{\mathrm{o}}(g), \]
where $\sim$ means equality up to scalar. Hence, $\theta_{g, S_1, \ldots, S_n}^{\mathrm{o}}$ generalizes the projective representation $\theta_g^{\mathrm{o}}$ of the first part of Theorem \ref{thmRepMCG} (when $S_1, \ldots, S_n$ are simple $H$-modules). Let $\rho^{g, S_1, \ldots, S_n}_{\mathrm{inv}}$ be the representation of $\mathcal{L}_{g,n}^{\mathrm{inv}}(H)$ on $\mathrm{Inv}\bigl( V(g, S_1, \ldots, S_n) \bigr)$. The statement generalizing the second part of Theorem \ref{thmRepMCG} is the following:

\smallskip

\noindent \textbf{Statement (to be proved).} {\em Let $S_1, \ldots, S_n$ be simple $H$-modules. The map
\[ \fonc{ \theta_{g, S_1, \ldots, S_n} }{ \mathrm{MCG}(\Sigma_{g,n}) }{ \mathrm{GL}\bigl( \mathrm{Inv}\bigl( V(g, S_1, \ldots, S_n) \bigr) \bigr) }{ f } { \rho^{g, S_1, \ldots, S_n}_{\mathrm{inv}}(\widehat{f}) } \]
is a projective representation.} 

\smallskip

\noindent To prove this, one must use a presentation of $\mathrm{MCG}(\Sigma_{g,n})$ based on the generators $\tau_{c_1}, \ldots, \tau_{c_k}$ and check that the relations between these generators hold in $\mathrm{Inv}\bigl( V(g, S_1, \ldots, S_n) \bigr)$. Note that the relations of $\mathrm{MCG}(\Sigma_{g,n})$ which already hold in $\mathrm{MCG}(\Sigma_{g,n}^{\mathrm{o}})$ are automatically satisfied; thus it is relevant to use a presentation of $\mathrm{MCG}(\Sigma_{g,n})$ which is a quotient of $\mathrm{MCG}(\Sigma_{g,n}^{\mathrm{o}})$ by some extra relations. For instance, recall that in the case $n=0$ we had to prove only the validity of one relation in $\mathrm{Inv}\bigl( (H^*)^{\otimes g} \bigr)$ (the hyperelliptic relation) because $\mathrm{MCG}(\Sigma_g)$ was the quotient of $\mathrm{MCG}(\Sigma_g^{\mathrm{o}})$ by this relation.

\subsection{Explicit formulas for the representation of some Dehn twists}\label{sectionFormulesExplicites}
\indent We will compute explicitly the representation on $(H^*)^{\otimes g}$ of the Dehn twists $\tau_{\gamma}$, where the curves $\gamma$ are represented in Figure \ref{figureCourbesCanoniques}. Thanks to Proposition \ref{propDehnTwist}, this amounts to compute the action of $v_{\widetilde{\gamma}}^{-1}$ on $(H^*)^{\otimes g}$.
\smallskip\\
\indent We recall that the action $\triangleright$ of $\mathcal{L}_{g,0}(H)$ on $(H^*)^{\otimes g}$ is defined using $\Psi_{g,0}$ in \eqref{actionTriangle} and that we denote the associated representation by $\rho$. Also recall the definition of the elements $\widetilde{h}$ in \eqref{operateursTilde} and the notation  $RR' = X_i \otimes Y_i$. Note that
\begin{equation}\label{coproduitRR}
X_i \otimes Y_i' \otimes Y_i'' = a_j X_i b_k \otimes Y_i \otimes b_j a_k.
\end{equation}

\indent Recall from Chapter \ref{chapitreTore} the representation $\theta_1^{\mathrm{o}} : \mathrm{MCG}(\Sigma_1^{\mathrm{o}}) \to \mathrm{GL}(H^*)$ given by the action of the elements $v_A^{-1}, v_B^{-1} \in \mathcal{L}_{1,0}(H)$  on $H^*$:
\begin{equation}\label{actionvAvB}
\begin{split}
\theta_1^{\mathrm{o}}(\tau_a)(\varphi) = v_A^{-1} \triangleright \varphi &= \varphi^{v^{-1}},\\
\theta_1^{\mathrm{o}}(\tau_b)(\varphi) = v_B^{-1} \triangleright \varphi &= \mu^l(v)^{-1}\left(\mu^l\!\left(g^{-1}v\,?\right) \varphi^v\right)^{v^{-1}}
\end{split}
\end{equation}
where $\beta^h = \beta(h?) \in H^*$ for any $\beta \in H^*, h \in H$ and $\mu^l$ is the left integral on $H$.
\smallskip\\
We will need the following generalization of Lemma \ref{keylemmaOmega} (in which we restricted to $\varphi \in \mathrm{SLF}(H)$).

\begin{lemma}\label{lemmaActionTresseAB}
For all $\varphi \in H^*$:
\begin{align*}
\left( v_{A}^{-1} v_{B}^{-1} v_{A}^{-1} \right)^2 \triangleright \varphi &= \frac{\mu^l(v^{-1})}{\mu^l(v)} \varphi\!\left( S^{-1}(a_i) g^{-1}v^{-1} S(?) b_i \right)\\
\left( v_{A}^{-1} v_{B}^{-1} v_{A}^{-1} \right)^{-2} \triangleright \varphi &= \frac{\mu^l(v)}{\mu^l(v^{-1})} \varphi\!\left( b_j S^{-1}(?) a_j g^{-1}v \right)
\end{align*}
\end{lemma}
\begin{proof}
Write $\varphi = \sum_{I,i,j} \Phi_{I,i}^j \overset{I}{T}{_j^i} = \sum_I \mathrm{tr}\!\left( \Phi_I \overset{I}{T} \right)$ with $\Phi_{I,i}^j \in \mathbb{C}$ and let $z(\varphi) = \sum_I \mathrm{tr}\!\left( \overset{I}{b_i} \Phi_I \overset{I}{S^{-1}(a_i)} \overset{I}{M} \right)$ $\in \mathcal{L}_{0,1}(H)$. Then $z(\varphi)_B \triangleright \varepsilon = \varphi$ (where $z(\varphi)_B = j_B(z(\varphi))$, see notation at the end of section \ref{defLgn}), and $\varepsilon$ is the counit of $H$). Indeed
\[ z(\varphi)_B \triangleright \varepsilon 
= \sum_I \mathrm{tr}\!\left( \overset{I}{b_i} \Phi_I \overset{I}{S^{-1}(a_i)} \overset{I}{L}{^{(+)}} \overset{I}{T} \overset{I}{L}{^{(-)-1}} \triangleright \varepsilon \right)
= \sum_I \mathrm{tr}\!\left( \overset{I}{b_i} \Phi_I \overset{I}{S^{-1}(a_i)} \overset{I}{a_j} \overset{I}{T} \overset{I}{b_j}  \right) = \sum_I \mathrm{tr}\!\left(\Phi_I  \overset{I}{T} \right) = \varphi.\]
We simply used \eqref{actionTriangle}, \eqref{repTriangleHeisenberg}, the cyclicity of the trace and the equality $S^{-1}(a_i)a_j \otimes b_j b_i = 1 \otimes 1$. Observe that 
\[ \left(\widetilde{\tau}_a \widetilde{\tau}_b \widetilde{\tau}_a\right)^2 (\overset{I}{B}) = \overset{I}{v}{^2}\overset{I}{A}{^{-1}}\overset{I}{B}{^{-1}}\overset{I}{A} = \overset{I}{B}{^{-1}} \overset{I}{C} \]
where $\overset{I}{C} = \overset{I}{C}_{1,0}$ is defined in \eqref{matriceCL10}. Hence:
\begin{align*}
\left( v_{A}^{-1} v_{B}^{-1} v_{A}^{-1} \right)^2 \triangleright \varphi &= \left( v_{A}^{-1} v_{B}^{-1} v_{A}^{-1} \right)^2 z(\varphi)_{B} \triangleright \varepsilon = z(\varphi)_{B^{-1}C} \left( v_{A}^{-1} v_{B}^{-1} v_{A}^{-1} \right)^2 \triangleright \varepsilon\\
&= \frac{\mu^l(v^{-1})}{\mu^l(v)} z(\varphi)_{B^{-1}C} \triangleright \varepsilon = \frac{\mu^l(v^{-1})}{\mu^l(v)} z(\varphi)_{B^{-1}} \triangleright \varepsilon.
\end{align*}
We used Proposition \ref{propDehnTwist}, the formula of Lemma \ref{keylemmaOmega} applied to $\varepsilon$, and the fact that $\overset{I}{C} \triangleright \varepsilon = \mathbb{I}_{\dim(I)}\varepsilon$ (which follows from \ref{actionH}). Now we compute
\begin{align*}
z(\varphi)_{B^{-1}} \triangleright \varepsilon &= \sum_I \mathrm{tr}\!\left( \overset{I}{b_i} \Phi_I \overset{I}{S^{-1}(a_i)} \overset{I}{L}{^{(-)}} S(\overset{I}{T}) \overset{I}{L}{^{(+)-1}} \triangleright \varepsilon \right)
= \sum_I \mathrm{tr}\!\left( \overset{I}{b_i} \Phi_I \overset{I}{S^{-1}(a_i)} \overset{I}{S^{-1}(b_j)}a_j \triangleright S(\overset{I}{T}) \right)\\
&= \sum_I \mathrm{tr}\!\left( \overset{I}{b_i} \Phi_I \overset{I}{S^{-1}(a_i)} \overset{I}{S^{-1}(b_j)} \overset{I}{S(a_j)} S(\overset{I}{T}) \right) = \sum_I \mathrm{tr}\!\left( \Phi_I \overset{I}{S^{-1}(a_i)} \overset{I}{g}{^{-1}} \overset{I}{v}{^{-1}} S(\overset{I}{T}) \overset{I}{b_i} \right)\\
&= \varphi\!\left( S^{-1}(a_i) g^{-1}v^{-1} S(?) b_i \right).
\end{align*}
We used \eqref{repHO} and \eqref{elementDrinfeld}. The second formula is easily checked.
\end{proof}

\begin{theorem}\label{formulesExplicites}
Let $\theta_g^{\mathrm{o}} : \mathrm{MCG}(\Sigma_g^{\mathrm{o}}) \to \mathrm{PGL}\bigl((H^*)^{\otimes g}\bigr)$ be the projective representation obtained in Theorem \ref{thmRepMCG}. The following formulas hold:
\begin{align*}
\theta_g^{\mathrm{o}}(\tau_{a_i})\bigl(\varphi_1 \otimes \ldots \otimes \varphi_g\bigr) &= \varphi_1 \otimes \ldots \otimes \varphi_{i-1} \otimes \theta_1^{\mathrm{o}}(\tau_a)(\varphi_i) \otimes \varphi_{i+1} \otimes \ldots \otimes \varphi_g, \\
\theta_g^{\mathrm{o}}(\tau_{b_i})\bigl(\varphi_1 \otimes \ldots \otimes \varphi_g\bigr) &= \varphi_1 \otimes \ldots \otimes \varphi_{i-1} \otimes \theta_1^{\mathrm{o}}(\tau_b)(\varphi_i) \otimes \varphi_{i+1} \otimes \ldots \otimes \varphi_g, \\
\theta_g^{\mathrm{o}}(\tau_{d_i})\bigl(\varphi_1 \otimes \ldots \otimes \varphi_g\bigr) &= \varphi_1 \otimes \ldots \otimes \varphi_{i-2} \otimes \varphi_{i-1}\!\left(S^{-1}(a_j)a_k?b_k v''^{-1} b_j\right) \otimes \varphi_i\!\left( S^{-1}(a_l)  S^{-1}(v'^{-1}) a_m ? b_m b_l \right)\\
& \:\:\:\:\: \otimes \varphi_{i+1} \otimes \ldots \otimes \varphi_g, \\
\theta_g^{\mathrm{o}}(\tau_{e_i})\bigl(\varphi_1 \otimes \ldots \otimes \varphi_g\bigr) &= \varphi_1\!\left(S^{-1}\!\left(v^{(2i-2)-1}\right) ? v^{(2i-1)-1} \right) \otimes \ldots \otimes \varphi_{i-1}\!\left(S^{-1}\!\left( v^{(2)-1} \right) ?  v^{(3)-1} \right)\\
&\:\:\:\:\: \otimes \varphi_i\!\left(S^{-1}(a_j)  S^{-1}\!\left(v^{(1)-1}\right) a_k ? b_k b_j \right) \otimes \varphi_{i+1} \otimes \ldots \otimes \varphi_g,
\end{align*}
with $i \geq 2$ for the two last formulas, $R = a_j \otimes b_j$ is the $R$-matrix\footnote{Do not confuse the components $a_j, b_j$ of the $R$-matrix and the loops $a_i, b_i \in \pi_1(\Sigma_{g,0} \!\setminus\! D)$.} and the formulas for $\theta_1^{\mathrm{o}}(\tau_a), \theta_1^{\mathrm{o}}(\tau_b)$ are recalled in \eqref{actionvAvB} above.
\end{theorem}

\begin{proof} First, it is useful to record that
\begin{equation}\label{formuleLambda}
\begin{split}
\Psi_{1,0}^{\otimes g}(\overset{I}{\Lambda}_i) = \Psi_{1,0}^{\otimes g}\!\left( \overset{I}{\underline{C}}{^{(-)}}(1) \ldots \overset{I}{\underline{C}}{^{(-)}}(i-1) \right) &= \overset{I}{S^{-1}(b_j)} \: \widetilde{a_j^{(2i-3)}} a_j^{(2i-2)}  \otimes  \ldots  \otimes  \widetilde{a_j^{(1)}}a_j^{(2)}  \otimes  1^{\otimes g-i+1} \\
\Psi_{1,0}^{\otimes g}\!\left( \overset{I}{\underline{C}}{^{(+)}}(1) \ldots \overset{I}{\underline{C}}{^{(+)}}(i-1) \right) &= \overset{I}{a_j} \: \widetilde{b_j^{(2i-3)}} b_j^{(2i-2)}  \otimes  \ldots  \otimes  \widetilde{b_j^{(1)}}b_j^{(2)}  \otimes  1^{\otimes g-i+1}
\end{split}
\end{equation}
where the matrix $\overset{I}{\Lambda}_k$ is defined in \eqref{matricesAlekseev}. The proof is a simple computation analogous to that of Lemma \ref{expressionM}. Second, recall from the proof of Lemma \ref{lemmavA1Moins1} that 
\begin{equation}\label{VIntegrale}
\mu^l(v)^{-1}\mu^l(vX_i)Y_i = v^{-1}.
\end{equation}
We will write $\mu^l(v)^{-1}\mu^l(v?) = \sum_{I} \mathrm{tr}\!\left(c_I \overset{I}{T}\right)$. Then $v^{-1} = \sum_{I} \mathrm{tr}\!\left(c_I \overset{I}{M}\right)$ under the identification $\mathcal{L}_{0,1}(H) = H$.
\smallskip\\
\indent \textbullet ~ {\em Proof of the formula for the action of $v_{A(i)}^{-1}$.} By definition and by \eqref{formuleLambda}, we have
\begin{align*}
&\Psi_{g,0}(v_{A(i)}^{-1}) = \sum_{I} \mathrm{tr}\!\left(c_I \Psi_{1,0}^{\otimes g}(\overset{I}{\Lambda}_i \, \overset{I}{\underline{A}}(i) \, \overset{I}{\Lambda}_i{^{-1}})\right)\\
& = \sum_I \mathrm{tr}\!\left(c_I \overset{I}{S^{-1}(b_j)} \overset{I}{X_k} \overset{I}{b_l}\right)  \widetilde{a_j^{(2i-3)}}\widetilde{a_l^{(2i-3)}} a_j^{(2i-2)} a_l^{(2i-2)}  \otimes  \ldots  \otimes  \widetilde{a_j^{(1)}}\widetilde{a_l^{(1)}} a_j^{(2)}a_l^{(2)}  \otimes  Y_k  \otimes  1^{\otimes g-i}\\
&=\mu^l(v)^{-1}\mu^l\!\left( v S^{-1}(b_j) X_k b_l\right)  \widetilde{a_j^{(2i-3)}}\widetilde{a_l^{(2i-3)}} a_j^{(2i-2)} a_l^{(2i-2)}  \otimes  \ldots  \otimes  \widetilde{a_j^{(1)}}\widetilde{a_l^{(1)}} a_j^{(2)}a_l^{(2)}  \otimes  Y_k  \otimes  1^{\otimes g-i}\\
&= \mu^l(v)^{-1}\mu^l\!\left(v S^{-1}\!\left(b_jS^{-1}(b_l)\right) X_k\right)  \widetilde{a_j^{(2i-3)}}\widetilde{a_l^{(2i-3)}} a_j^{(2i-2)} a_l^{(2i-2)}  \otimes  \ldots  \otimes  \widetilde{a_j^{(1)}}\widetilde{a_l^{(1)}} a_j^{(2)}a_l^{(2)}  \otimes  Y_k  \otimes  1^{\otimes g-i}\\
&= \mu^l(v)^{-1}\mu^l\!\left( v X_k\right)  1^{\otimes i-1}  \otimes  Y_k  \otimes  1^{\otimes g-i} = 1^{\otimes i-1}  \otimes  v^{-1}  \otimes  1^{\otimes g-i}
\end{align*}
and the formula follows. We used \eqref{quasiCyclic}, the formula $R^{-1} = a_l \otimes S^{-1}(b_l)$ and \eqref{VIntegrale}.
\smallskip\\
\indent \textbullet ~ {\em Proof of the formula for the action of $v_{B(i)}^{-1}$.} This the same proof as for $v_{A(i)}^{-1}$ (the conjugation by $\overset{I}{\Lambda}_i$ vanishes thanks to \eqref{quasiCyclic}).
\smallskip\\
\indent \textbullet ~ {\em Proof of the formula for the action of $v_{\widetilde{d_i}}^{-1}$, $i \geq 2$.} We first compute the action of $\overset{I}{A}(i-1)\overset{I}{A}(i)$. We have
\begin{align*}
\Psi_{g,0}\!\left( \overset{I}{A}(i-1)\overset{I}{A}(i) \right) &= \Psi_{1,0}^{\otimes g}\!\left(\overset{I}{\Lambda}_{i-1} \, \overset{I}{\underline{A}}(i-1) \, \overset{I}{\Lambda}{_{i-1}^{-1}} \, \overset{I}{\Lambda}_{i} \, \overset{I}{\underline{A}}(i) \, \overset{I}{\Lambda}{_{i}^{-1}}\right)\\
&= \Psi_{1,0}^{\otimes g}\!\left(\overset{I}{\Lambda}_{i-1} \, \overset{I}{\underline{A}}(i-1) \, \overset{I}{\underline{C}}{^{(-)}}(i-1) \, \overset{I}{\underline{A}}(i) \, \overset{I}{\underline{C}}{^{(-)}}(i-1)^{-1} \, \overset{I}{\Lambda}{_{i-1}^{-1}}\right).
\end{align*}
Hence:
\begin{align*}
&\Psi_{g,0}\!\left( v_{A(i-1)A(i)}^{-1} \right)
= \sum_{I} \mathrm{tr}\!\left(c_I \Psi_{1,0}^{\otimes g}\!\left(\overset{I}{\Lambda}_{i-1} \, \overset{I}{\underline{A}}(i-1) \, \overset{I}{\underline{C}}{^{(-)}}(i-1) \, \overset{I}{\underline{A}}(i) \, \overset{I}{\underline{C}}{^{(-)}}(i-1)^{-1} \, \overset{I}{\Lambda}{_{i-1}^{-1}}\right)\right)\\
&= \mu^l(v)^{-1}\mu^l\!\left( v S^{-1}(b_j) X_k S^{-1}(b_l) X_m b_n b_o \right)  \widetilde{a_j^{(2i-5)}}\widetilde{a_o^{(2i-5)}} a_j^{(2i-4)} a_o^{(2i-4)} \, \otimes \, \ldots \, \otimes \, \widetilde{a_j^{(1)}}\widetilde{a_o^{(1)}} a_j^{(2)}a_o^{(2)}\\
& \:\:\:\:\:\:\:\:\:\:\:\:\:\:\:\:\:\:\:\:\:\:\:\:\:\:\:\:\:\:\:\:\:\:\:\:\:\:\:\:\:\:\:\:\:\:\:\:\:\:\:\:\:\:\:\:\:\:\:\:\:\:\:\:\:\:\:\:\:\:\:\:\:\:\:\otimes \widetilde{a_l'}\widetilde{a_n'} Y_k a_l'' a_n'' \otimes Y_m \otimes 1^{\otimes g-i}\\
&= \mu^l(v)^{-1}\mu^l\!\left( v X_k S^{-1}(b_l) X_m b_n \right)  1^{\otimes i-2} \otimes \widetilde{a_l'}\widetilde{a_n'} Y_k a_l'' a_n'' \otimes Y_m \otimes 1^{\otimes g-i}\\
&= \mu^l(v)^{-1}\mu^l\!\left( v a_k S^{-1}(b_l) X_m b_n \right)  1^{\otimes i-2} \otimes \widetilde{a_l}\widetilde{a_n'} b_k a_n'' \otimes Y_m \otimes 1^{\otimes g-i}\\
\end{align*}
We used \eqref{quasiCyclic} and the fact that $X_k S^{-1}(b_l) \otimes a_l' \otimes Y_k a_l'' = X_k S^{-1}(b_p)S^{-1}(b_l) \otimes a_l \otimes Y_k a_p = a_kS^{-1}(b_l) \otimes a_l \otimes b_k$. We can assume without loss of generality that $g=2$ and $i=2$, since the action is ``local''. Moreover, this can be simplified. Let $F : H^* \to H^*$ be the linear map defined by
\[ F(\varphi) = \varphi\!\left(a_j ? b_j\right), \:\:\:\:\: F^{-1}(\varphi) = \varphi\!\left(S^{-1}(a_j) ? b_j\right). \]
We compute:
\begin{align*}
&(F^{-1} \otimes \mathrm{id}) \circ \rho\!\left( v_{A(1)A(2)}^{-1} \right) \circ (F \otimes \mathrm{id})(\varphi \otimes \psi)\\
&= \mu^l(v)^{-1}\mu^l\!\left( v a_k b_l X_m b_n \right)   \varphi\!\left( a_j S^{-1}(a_n')a_l S^{-1}(a_o) ? b_o b_k a_n''b_j \right) \otimes \psi(?Y_m)\\
&= \mu^l(v)^{-1}\mu^l\!\left(v a_k b_l X_m b_n \right)   \varphi\!\left( S^{-1}(a_n'') a_j a_l S^{-1}(a_o) ? b_o b_k b_j a_n' \right) \otimes \psi(?Y_m) = (\star).
\end{align*}
We used the formula $R\Delta = \Delta^{\mathrm{op}}R$. Now, we have a Yang-Baxter identity 
\[ a_k b_l \otimes a_j a_l \otimes b_k b_j = R_{13} R_{23} R_{21} = R_{21} R_{23} R_{13} = b_l a_k \otimes a_l a_j \otimes b_j b_k \]
which allows us to continue the computation:
\begin{align*}
(\star) &= \mu^l(v)^{-1}\mu^l\!\left(v b_l a_k X_m b_n \right)   \varphi\!\left( S^{-1}(a_n'') a_l a_j S^{-1}(a_o) ? b_o b_j b_k a_n' \right) \otimes \psi(?Y_m)\\
&= \mu^l(v)^{-1}\mu^l\!\left(v S^{-2}(b_p) b_l a_k X_m b_n \right)   \varphi\!\left( S^{-1}(a_p) a_l ?  b_k a_n \right) \otimes \psi(?Y_m)\\
&= \mu^l(v)^{-1}\mu^l\!\left(v a_k X_m b_n \right)   \varphi\!\left( ?  b_k a_n \right) \otimes \psi(?Y_m)\\
&= \mu^l(v)^{-1}\mu^l\!\left(v X_m \right)   \varphi\!\left( ?  Y_m'' \right) \otimes \psi(? Y_m') = \varphi\!\left( ?  v''^{-1} \right) \otimes \psi(? v'^{-1}).
\end{align*}
We used basic properties of the $R$-matrix and relations \eqref{coproduitRR}, \eqref{VIntegrale}. We have thus shown that
\[ v_{A(1)A(2)}^{-1} \triangleright \varphi \otimes \psi = \varphi\!\left( S^{-1}(a_j)a_k ? b_k v''^{-1} b_j \right) \otimes \psi(?v'^{-1}). \]
Recall that $\overset{I}{\widetilde{d_2}} = \overset{I}{v}{^2} \overset{I}{A}(1)\overset{I}{B}(2)\overset{I}{A}(2){^{-1}} \overset{I}{B}(2){^{-1}}$. Hence $\left(\widetilde{\tau}_{a_2} \widetilde{\tau}_{b_2} \widetilde{\tau}_{a_2} \right)^{-2}(\overset{I}{A}(1)\overset{I}{A}(2)) = \overset{I}{\widetilde{d_2}}$. It follows that $\left(\widetilde{\tau}_{a_2} \widetilde{\tau}_{b_2} \widetilde{\tau}_{a_2} \right)^{-2}(v_{A(1)A(2)}^{-1}) = v_{\widetilde{d_2}}^{-1}$, and thus by Proposition \ref{propDehnTwist} and Lemma \ref{lemmaActionTresseAB}:
\begin{align*}
v_{\widetilde{d_2}}^{-1} \triangleright \varphi \otimes \psi &= \left( v_{A(2)}^{-1} v_{B(2)}^{-1} v_{A(2)}^{-1} \right)^{-2} v_{A(1)A(2)}^{-1} \left( v_{A(2)}^{-1} v_{B(2)}^{-1} v_{A(2)}^{-1} \right)^2 \triangleright \varphi \otimes \psi\\
&= \left( v_{A(2)}^{-1} v_{B(2)}^{-1} v_{A(2)}^{-1} \right)^{-2} \triangleright \varphi\!\left( S^{-1}(a_j)a_k ? b_k v''^{-1} b_j \right) \otimes \psi\!\left( S^{-1}(a_l) g^{-1} v^{-1}  S(v'^{-1}) S(?) b_l \right)\\
&= \varphi\!\left( S^{-1}(a_j)a_k ? b_k v''^{-1} b_j \right) \otimes \psi\!\left( S^{-1}(a_l) g^{-1} v^{-1}  S(v'^{-1}) S\!\left(b_m S^{-1}(?) a_m g^{-1}v\right) b_l \right)\\
&= \varphi\!\left( S^{-1}(a_j)a_k ? b_k v''^{-1} b_j \right) \otimes \psi\!\left( S^{-1}(a_l)  S^{-1}(v'^{-1}) a_m ? b_mb_l \right)
\end{align*}
which is the announced formula.
\smallskip\\
\indent \textbullet ~ {\em Proof of the formula for the action of $v_{\widetilde{e_i}}^{-1}$, $i \geq 2$.} We first compute the action of $\overset{I}{C}(1) \ldots \overset{I}{C}(i-1) \overset{I}{A}(i)$. We have
\begin{align*}
\Psi_{g,0}\!\left(\overset{I}{C}(1) \ldots \overset{I}{C}(i-1) \overset{I}{A}(i) \right) &= \Psi_{1,0}^{\otimes g}\!\left( \overset{I}{\underline{C}}{^{(+)}}(1) \ldots \overset{I}{\underline{C}}{^{(+)}}(i-1) \left(\overset{I}{\underline{C}}{^{(-)}}(1) \ldots \overset{I}{\underline{C}}{^{(-)}}(i-1)\right)^{-1} \overset{I} \Lambda_i \overset{I}{\underline{A}}(i) \Lambda_i^{-1} \right)\\
&= \Psi_{1,0}^{\otimes g}\!\left( \overset{I}{\underline{C}}{^{(+)}}(1) \ldots \overset{I}{\underline{C}}{^{(+)}}(i-1) \overset{I}{\underline{A}}(i) \left(\overset{I}{\underline{C}}{^{(-)}}(1) \ldots \overset{I}{\underline{C}}{^{(-)}}(i-1)\right)^{-1} \right)\\
&= \overset{I}{a_j} \overset{I}{X_k} \overset{I}{b_l} \, \widetilde{b_j^{(2i-3)}} \widetilde{a_l^{(2i-3)}} b_j^{(2i-2)} a_l^{(2i-2)}  \otimes  \ldots  \otimes  \widetilde{b_j^{(1)}} \widetilde{a_l^{(1)}} b_j^{(2)} a_l^{(2)}  \otimes  Y_k \otimes 1^{\otimes g-i}\\
&= \overset{I}{X_k} \, \widetilde{Y_k^{(2i-2)}}  Y_k^{(2i-1)} \otimes \ldots \otimes \widetilde{Y_k^{(2)}}  Y_k^{(3)}   \otimes  Y_k^{(1)} \otimes 1^{\otimes g-i}
\end{align*}
thanks to \eqref{formuleLambda} and \eqref{coproduitRR}. Hence, by \eqref{VIntegrale}:
\begin{align*}
\Psi_{g,0}\!\left(v_{C(1) \ldots C(i-1) A(i)}^{-1} \right) &= \mu^l(v)^{-1}\mu^l(v X_k)\, \widetilde{Y_k^{(2i-2)}}  Y_k^{(2i-1)} \otimes \ldots \otimes \widetilde{Y_k^{(2)}}  Y_k^{(3)}   \otimes  Y_k^{(1)} \otimes 1^{\otimes g-i}\\
&= \widetilde{v^{(2i-2)-1}}  v^{(2i-1)-1} \otimes \ldots \otimes \widetilde{v^{(2)-1}}  v^{(3)-1}   \otimes  v^{(1)-1} \otimes 1^{\otimes g-i},
\end{align*}
which means that
\begin{align*}
v_{C(1) \ldots C(i-1) A(i)}^{-1} \triangleright \left(\varphi_1 \otimes \ldots \otimes \varphi_g\right) &= \varphi_1\!\left(S^{-1}\!\left(v^{(2i-2)-1}\right) ? v^{(2i-1)-1} \right) \otimes \ldots \otimes \varphi_{i-1}\!\left(S^{-1}\!\left( v^{(2)-1} \right) ?  v^{(3)-1} \right)\\
&\:\:\:\:\: \otimes \varphi_i\!\left( ? v^{(1)-1} \right) \otimes \varphi_{i+1} \ldots \otimes \varphi_g. 
\end{align*}
Recall that $\overset{I}{\widetilde{e_i}} = \overset{I}{v}{^2} \overset{I}{C}(1) \ldots \overset{I}{C}(i-1)\overset{I}{B}(i)\overset{I}{A}(i){^{-1}} \overset{I}{B}(i){^{-1}}$. Hence $\left(\widetilde{\tau}_{a_i} \widetilde{\tau}_{b_i} \widetilde{\tau}_{a_i} \right)^{-2}(\overset{I}{C}(1) \ldots \overset{I}{C}(i-1)\overset{I}{A}(i)) = \overset{I}{\widetilde{e_i}}$. As previously, it follows from Proposition \ref{propDehnTwist} and Lemma \ref{lemmaActionTresseAB} that
\begin{align*}
v_{\widetilde{e_i}}^{-1} \triangleright \left( \varphi_1 \otimes \ldots \otimes \varphi_g \right) &= \left( v_{A(i)}^{-1} v_{B(i)}^{-1} v_{A(i)}^{-1} \right)^{-2} v_{C(1) \ldots C(i-1)A(i)}^{-1} \left( v_{A(i)}^{-1} v_{B(i)}^{-1} v_{A(i)}^{-1} \right)^2 \triangleright \left( \varphi_1 \otimes \ldots \otimes \varphi_g \right) \\
&= \varphi_1\!\left(S^{-1}\!\left(v^{(2i-2)-1}\right) ? v^{(2i-1)-1} \right) \otimes \ldots \otimes \varphi_{i-1}\!\left(S^{-1}\!\left( v^{(2)-1} \right) ?  v^{(3)-1} \right)\\
&\:\:\:\:\: \otimes \varphi_i\!\left(S^{-1}(a_j)  S^{-1}\!\left(v^{(1)-1}\right) a_k ? b_k b_j \right) \otimes \varphi_{i+1} \ldots \otimes \varphi_g,
\end{align*}
which is the announced formula.
\end{proof}

\section{Equivalence with the Lyubashenko representation}
\indent In a series of papers \cite{lyu95a, lyu95b, lyu96}, V. Lyubashenko has constructed projective representations of $\mathrm{MCG}(\Sigma_{g,n})$ by categorical techniques based on the coend of a ribbon category. Our assumptions on $H$ allow to apply his construction to $\mathrm{mod}_l(H)$, the ribbon category of finite dimensional left $H$-modules. Here we will show that these two representations are equivalent. For the case of the torus, we have already shown in Chapter \ref{chapitreTore} that the projective representation of $\mathrm{MCG}(\Sigma_{1,0})$ obtained thanks to $\mathcal{L}_{1,0}(H)$ is equivalent to the Lyubashenko-Majid representation \cite{LM}. For works based on the Lyubashenko representation, see \textit{e.g.} \cite{FSS1, FSS2}.

\subsection{The Lyubashenko representation for $\mathrm{mod}_l(H)$}
\indent Let us first quickly recall the Lyubashenko representation in the general framework of a ribbon category $\mathcal{C}$ satisfying some assumptions (see \cite{lyu95b}). 
\smallskip\\
\indent Let $K = \int^X \! X^* \otimes X$ be the coend of the functor $F : \mathcal{C}^{\mathrm{op}} \times \mathcal{C} \to \mathcal{C}$, $F(X,Y) = X^* \otimes Y$ and let $i_X : X^* \otimes X \to K$ be the associated dinatural transformation (see \cite[IX.6]{ML}). Thanks to the universal property of the coend $K$, Lyubashenko defined several morphisms; we will need some of them which we recall now. The first is an algebra structure $K \otimes K \to K$ (also see \cite{majid93}). Consider the following family of morphisms (for each $X,Y \in \mathcal{C}$)
\begin{equation}\label{dinatProduit}
\begin{split}
d_{X,Y} : X^* \otimes X  \otimes Y^* \otimes Y \xrightarrow{\mathrm{id}_{X^*} \otimes c_{X,Y^*} \otimes \mathrm{id}_Y} X^* \otimes Y^*  \otimes X \otimes Y  &\xrightarrow{\mathrm{id}_{X^*} \otimes \mathrm{id}_{Y^*} \otimes c_{X,Y}} X^* \otimes Y^*  \otimes Y \otimes X\\
&\:\:\:\:\:\:\:\:\: \xrightarrow{\sim} (Y \otimes X)^*  \otimes Y \otimes X \xrightarrow{i_{Y \otimes X}} K.
\end{split}
\end{equation}
Since the family $d_{X,Y}$ is dinatural in $X$ and $Y$, it exists a unique $m_K : K \otimes K \to K$ such that $d_{X,Y} = m_K \circ (i_X \otimes i_Y)$, which is in fact an associative product on $K$. Actually, $K$ is endowed with a Hopf algebra structure whose structure morphisms are similarly defined using the universal property, but we do not need this here.
\smallskip\\
\indent Next, consider the following families of morphisms
\begin{equation}\label{dinatRep}
\begin{split}
&\alpha_X : X^* \otimes X \xrightarrow{\theta_{X^*} \otimes \mathrm{id}_X} X^* \otimes X \xrightarrow{i_X} K,\\
&\beta_{X,Y} : X^* \otimes X \otimes Y^* \otimes Y \xrightarrow{\mathrm{id}_{X^*} \otimes \left(c_{Y^*, X} \circ c_{X, Y^*}\right) \otimes \mathrm{id}_Y} X^* \otimes X \otimes Y^* \otimes Y \xrightarrow{i_X \otimes i_Y} K \otimes K,\\
&\gamma_X^Y : X^* \otimes X \otimes Y \xrightarrow{\mathrm{id}_{X^*} \otimes \left(c_{Y, X} \circ c_{X, Y}\right)} X^* \otimes X \otimes Y \xrightarrow{i_X \otimes \mathrm{id}_Y} K \otimes Y.
\end{split}
\end{equation}
\indent The families $\alpha_X$ and $\gamma_X^Y$ (with $Y$ fixed) are dinatural in $X$, and the family $\beta_{X,Y}$ is dinatural in $X,Y$. Hence by universality of $K$ (and also thanks to the Fubini theorem for multiple coends and to the universality of $K \otimes K$ and of $K \otimes Y$), there exist unique morphisms $\mathcal{T} : K \to K, \mathcal{O} : K \otimes K \to K \otimes K, \mathcal{Q}_Y : K \otimes Y \to K \otimes Y$ such that
\begin{equation}\label{morphismesRep}
\alpha_X = \mathcal{T} \circ i_X, \:\:\:\:\: \beta_{X,Y} = \mathcal{O} \circ (i_X \otimes i_Y), \:\:\:\:\: \gamma_X^Y = \mathcal{Q}_Y \circ (i_X \otimes \mathrm{id}_Y).
\end{equation}

\indent Finally, the morphism $\mathcal{S} : K \to K$ is defined by $\mathcal{S} = (\varepsilon_K \otimes \mathrm{id}_K) \circ \mathcal{O} \circ (\mathrm{id}_K \otimes \Lambda_K)$, where $\varepsilon_K$ is the counit on $K$ and $\Lambda_K$ is the two-sided cointegral. 
\smallskip\\
\indent Let $X$ be any object of $\mathcal{C}$ and $V_X = \mathrm{Hom}_{\mathcal{C}}(X, K^{\otimes g})$. The Lyubashenko projective representation $Z_X : \mathrm{MCG}(\Sigma_g^{\mathrm{o}}) \to \mathrm{PGL}(V_X)$ \cite[Section 4.4]{lyu95b} takes the following values on $f \in V_X$:
\begin{equation}\label{repLyubashenkoC}
\begin{split}
Z_X(\tau_{a_i})(f) &= \bigl(\mathrm{id}_K^{\otimes g-i} \otimes \mathcal{T} \otimes \mathrm{id}_K^{\otimes i-1}\bigr) \circ f,\\
Z_X(\tau_{b_i})(f) &= \bigl(\mathrm{id}_K^{\otimes g-i} \otimes (\mathcal{S}^{-1} \circ \mathcal{T} \circ \mathcal{S}) \otimes \mathrm{id}_K^{\otimes i-1}\bigr) \circ f,\\
Z_X(\tau_{d_i})(f) &= \bigl(\mathrm{id}_K^{\otimes g-i} \otimes (\mathcal{O} \circ (\mathcal{T} \otimes \mathcal{T})) \otimes \mathrm{id}_K^{\otimes i-2}\bigr) \circ f \:\:\: \text{for } i \geq 2, \\
Z_X(\tau_{e_i})(f) &= \left(\mathrm{id}_K^{\otimes g-i} \otimes \bigl((\mathcal{T} \otimes \theta_{K^{\otimes i-1}}) \circ \mathcal{Q}_{K^{\otimes i-1}}\bigr) \right) \circ f \:\:\: \text{for } i \geq 2.
\end{split}
\end{equation}
Recall that the curves $a_i, b_i, d_i, e_i$ are represented in Figure \ref{figureCourbesCanoniques}. Since these Dehn twists are a generating set, we have an operator $Z_X(f)$ for all $f \in \mathrm{MCG}(\Sigma_g^{\mathrm{o}})$. If moreover we take $X = \mathbf{1\!\!\!1}$, the unit object of $\mathcal{C}$, then this defines a projective representation $Z_{\mathbf{1\!\!\!1}} : \mathrm{MCG}(\Sigma_g) \to \mathrm{PGL}(V_{\mathbf{1\!\!\!1}})$ of the mapping class group of $\Sigma_g$.
\smallskip\\
\indent Now, let us explicit the above formulas in the case of $\mathcal{C} = \mathrm{mod}_l(H)$. Recall from section \ref{modH} that the category $\mathrm{mod}_l(H)$ has braiding $c_{X,Y} : X \otimes Y \to Y \otimes X$ and twist $\theta_X : X \to X$ given by
\[ c_{X,Y}(x \otimes y) = b_i \cdot y \otimes a_i \cdot x, \:\:\:\: \theta_X(x) = v^{-1} \cdot x \]
and that the action on the dual module $X^*$ is $h \cdot \varphi = \varphi(S(h) \cdot ?)$ for all $\varphi \in X^*, h \in H$.
\smallskip\\
\indent It is well-known (and not difficult to see) that $K$ is $H^*$ endowed with the coadjoint action:
\[ \forall \, h \in H, \:\forall \, \varphi \in K, \:\: h \varphi = \varphi(S(h')?h'') \]
and that the dinatural transformation of $K$ is
\[ i_X(\psi \otimes x) = \psi(? \cdot x) \in K. \]
Note that $\psi(? \cdot x)$ is just a matrix coefficient of the module $X$. The dinatural family $d_{X,Y}$ of \eqref{dinatProduit} is
\[ d_{X,Y}(\varphi \otimes x \otimes \psi \otimes y) = \psi(S(b_i)?b_j \cdot y) \varphi(?a_ja_i \cdot x) \]
where in the right of the equality it is the usual product in $H^*$: $\langle fg, h\rangle = f(h')g(h'')$. To compute the product $m_K$ in $K$ explicitly, observe that $i_{H_{\mathrm{reg}}}(\varphi \otimes 1) = \varphi$, where $H_{\mathrm{reg}}$ is the regular representation of $H$. Thus
\begin{align*}
m_K(\varphi \otimes \psi) &= m_K \circ (i_{H_{\mathrm{reg}}} \otimes i_{H_{\mathrm{reg}}})\!\left( \varphi \otimes 1 \otimes \psi \otimes 1 \right) = d_{H_{\mathrm{reg}}, H_{\mathrm{reg}}}(\varphi \otimes 1 \otimes \psi \otimes 1)\\
&= \psi(S(b_i)?b_j) \varphi(?a_ja_i) = \varphi(a_j ? a_i) \psi(S(b_i)b_j ?)
\end{align*}
where we used $R \Delta = \Delta^{\mathrm{op}} R$ for the last equality. Moreover, the unit element of $K$ is $1_K = \varepsilon$, the counit of $H$. We record the following lemma, already given in \cite{lyu95b}.

\begin{lemma}\label{cointegraleBilatere}
Let $\mu^r \in H^*$ be the right integral on $H$ (unique up to scalar). Then $\mu^r$ is the two-sided cointegral in $K$ (unique up to scalar):
\[ \forall\, \varphi \in K, \:\:\: m_K(\mu^r \otimes \varphi) = m_K(\varphi \otimes \mu^r) = \varepsilon_K(\varphi)\mu^r \]
where $\varepsilon_K(\varphi) = \varphi(1)$.
\end{lemma}
\begin{proof}
Using \eqref{quasiCyclic} and the basic properties of $R$, we get
\[ m_K(\mu^r \otimes \varphi) =  \mu^r(a_j ? a_i) \varphi(S(b_i)b_j ?) = \mu^r(S^2(a_i)a_j ?) \varphi(S(b_i)b_j ?) = \mu^r \varphi = \varphi(1)\mu^r. \]
Similarly:
\[ m_K(\varphi \otimes \mu^r) = \mu^r(S(b_i)?b_j) \varphi(?a_ja_i) = \mu^r(S^2(b_jS^{-1}(b_i))?) \varphi(?a_ja_i) = \mu^r \varphi = \varphi(1) \mu^r. \]
\end{proof}

\indent The dinatural families of \eqref{dinatRep} are
\begin{equation*}
\begin{split}
&\alpha_X(\varphi \otimes x) = \varphi(v^{-1}? \cdot x), \\
&\beta_{X,Y}(\varphi \otimes x \otimes \psi \otimes y) = \varphi(?b_j a_i \cdot x) \otimes \psi(S(a_jb_i)? \cdot y) = \varphi(?v'^{-1}v \cdot x) \otimes \psi(S(v''^{-1})v? \cdot y), \\
&\gamma_X^Y(\varphi \otimes x \otimes y) = \varphi(?b_j a_i \cdot x) \otimes a_j b_i \cdot y = \varphi(?v'^{-1}v \cdot x) \otimes v''^{-1}v \cdot y
\end{split}
\end{equation*}
where we used \eqref{ribbon}. It follows that the morphisms defined in \eqref{morphismesRep} are
\begin{equation*}
\begin{split}
&\mathcal{T}(\varphi) = \mathcal{T} \circ i_{H_{\mathrm{reg}}}(\varphi \otimes 1) = \alpha_{H_{\mathrm{reg}}}(\varphi \otimes 1) = \varphi(v^{-1}?), \\
&\mathcal{O}(\varphi \otimes \psi) = \mathcal{O} \circ (i_{H_{\mathrm{reg}}} \otimes i_{H_{\mathrm{reg}}})(\varphi \otimes 1 \otimes \psi \otimes 1) = \beta_{H_{\mathrm{reg}}, H_{\mathrm{reg}}}(\varphi \otimes 1 \otimes \psi \otimes 1) = \varphi(?v'^{-1}v) \otimes \psi(S(v''^{-1})v?), \\
&\mathcal{Q}_Y(\varphi \otimes y) = \mathcal{Q}_Y \circ (i_{H_{\mathrm{reg}}} \otimes \mathrm{id}_Y)(\varphi \otimes 1 \otimes y) = \gamma_{H_{\mathrm{reg}}}^Y(\varphi \otimes 1 \otimes y) = \varphi(?v'^{-1}v) \otimes v''^{-1}v \cdot y.
\end{split}
\end{equation*}
\noindent In view of \eqref{repLyubashenkoC}, note that $(\mathcal{T} \otimes \theta_Y) \circ \mathcal{Q}_Y(\varphi \otimes y) = \varphi(? v'^{-1}) \otimes v''^{-1} \cdot y$. Finally, thanks to Lemma \ref{cointegraleBilatere}, the morphism $\mathcal{S}$ is
\[ \mathcal{S}(\varphi) = \varphi\!\left(v'^{-1}v\right) \mu^r\!\left(S(v''^{-1})v?\right) = \varphi\!\left(S^{-1}(v''^{-1})v\right) \mu^r\!\left(v'^{-1}v?\right) \]
where the second equality is due to $v'^{-1} \otimes S(v''^{-1}) = S^{-1}(v''^{-1}) \otimes v'^{-1}$ (which follows from $S(v^{-1}) = v^{-1}$). Moreover, we will need the following lemma to prove the equivalence of the representations.
\begin{lemma}\label{lemmaST}
Let $\rho$ be the representation of $\mathcal{L}_{1,0}(H)$ on $H^*$, then the following formulas hold:
\begin{equation*}
\begin{split}
&\mathcal{T} = \rho(v_A^{-1}) = (v^{-1})_*, \:\:\:\: \mathcal{S} = \mu^l(v^{-1}) g^{-1}_* \circ \rho(v_A^2v_B) \circ g_*,\\
&\mathcal{S}^{-1} \circ \mathcal{T} \circ \mathcal{S} = (g^{-1} v)_* \circ \rho(v_B^{-1}) \circ (g v^{-1})_*,
\end{split}
\end{equation*}
where $h_*(\varphi) = \varphi(?h)$ for all $h \in H$ and $\varphi \in H^*$.
\end{lemma}
\begin{proof}
The formula for $\mathcal{T}$ is obvious. Propositions \ref{actionAB} and \ref{propVIntegrale} give $\rho(v_B)$ and then we compute using \eqref{muLmuRgCarre} and \eqref{integraleShifte}:
\begin{align*}
\rho(v_B)(\varphi) &= v_B \triangleright \varphi = \mu^l(v^{-1})^{-1}\left( \mu^l(g^{-1}v^{-1}?)\varphi^v \right)^{v^{-1}} = \mu^l(v^{-1})^{-1}\left( \mu^r(gv^{-1}?)\varphi^v \right)^{v^{-1}}\\
&= \mu^l(v^{-1})^{-1} \mu^r\!\left(v'^{-1}?gv^{-1}\right) \varphi\!\left(S^{-1}(v''^{-1})g^{-1}v\right)\\
&= \mu^l(v^{-1})^{-1} \left( gv^{-2} \right)_*\!\left(\mu^r(v v'^{-1}?)\right) \left\langle g^{-1}_*(\varphi), S^{-1}(v''^{-1})v \right\rangle\\
&=  \mu^l(v^{-1})^{-1} \left( gv^{-2} \right)_* \circ \mathcal{S} \circ g^{-1}_*(\varphi) = \mu^l(v^{-1})^{-1} \rho(v_A^{-2}) \circ g_* \circ \mathcal{S} \circ g^{-1}_*(\varphi)
\end{align*}
where $\varphi^h = \varphi(h?)$ for $h \in H$. The last claimed formula follows from $\mathcal{S} = \mu^l(v^{-1}) (g^{-1}v)_* \circ \rho(v_Av_Bv_A) \circ (gv^{-1})_*$ and the fact that  $v_A, v_B \in \mathcal{L}_{1,0}(H)$ satisfy the braid relation $v_A v_B v_A = v_B v_A v_B$ (see Proposition \ref{braidV}).
\end{proof}

For the representation space, we take $X = H_{\mathrm{reg}}$, so that $V_X = \mathrm{Hom}_H(H_{\mathrm{reg}}, K^{\otimes g}) \cong K^{\otimes g}$. Then by the previous formulas, we get the Lyubashenko projective representation of $\mathrm{MCG}(\Sigma_g^{\mathrm{o}})$ \eqref{repLyubashenkoC} applied to $\mathrm{mod}_l(H)$:
\begin{equation}\label{formulesLyubashenko}
\begin{split}
Z_{H_{\mathrm{reg}}}(\tau_{a_i})(\varphi_1 \otimes \ldots \otimes \varphi_g) &= \varphi_1 \otimes \ldots \otimes \varphi_{g-i+1}(v^{-1}?) \otimes \ldots \otimes \varphi_g, \\
Z_{H_{\mathrm{reg}}}(\tau_{b_i})(\varphi_1 \otimes \ldots \otimes \varphi_g) &= \varphi_1 \otimes \ldots \otimes (g^{-1}v)_* \circ \rho(v_B^{-1}) \circ (gv^{-1})_*(\varphi_{g-i+1}) \otimes \ldots \otimes \varphi_g, \\
Z_{H_{\mathrm{reg}}}(\tau_{d_i})(\varphi_1 \otimes \ldots \otimes \varphi_g) &=  \varphi_1 \otimes \ldots \otimes \varphi_{g-i+1}\!\left(? v'^{-1}\right) \otimes \varphi_{g-i+2}\!\left(S(v''^{-1})?\right) \otimes \ldots \otimes \varphi_g, \\
Z_{H_{\mathrm{reg}}}(\tau_{e_i})(\varphi_1 \otimes \ldots \otimes \varphi_g) &= \varphi_1 \otimes \ldots \otimes \varphi_{g-i} \otimes \varphi_{g-i+1}\!\left(? v^{(1)-1}\right) \otimes \varphi_{g-i+2}\!\left( S(v^{(2)-1}) ? v^{(3)-1} \right) \otimes \ldots \\
&\:\:\:\:\: \otimes \varphi_g\!\left( S(v^{(2i-2)-1}) ? v^{(2i-1)-1} \right),
\end{split}
\end{equation}
with $i \geq 2$ for the two last formulas. If we take $X = \mathbb{C}$, we get 
\[ V_{\mathbb{C}} = \mathrm{Hom}_H(\mathbb{C}, K^{\otimes g}) = (K^{\otimes g})^{\mathrm{inv}} = \left\{ f \in K^{\otimes g} \left| \, \forall\, h \in H, \:\: h \cdot f = \varepsilon(h)f\right.\right\} \]
where by definition of the action of $H$ on $K$, the action of $H$ on $K^{\otimes g}$ is
\begin{equation}\label{actionHKg}
h \cdot \varphi_1 \otimes \ldots \otimes \varphi_g = \varphi_1\!\left(S(h^{(1)}) ? h^{(2)}\right) \otimes \ldots \otimes \varphi_g\!\left(S(h^{(2g-1)}) ? h^{(2g)}\right).
\end{equation}
Then $Z_{\mathbb{C}}$ is a projective representation of $\mathrm{MCG}(\Sigma_g)$ (note that $Z_{\mathbb{C}}$ is just  $Z_{H_{\mathrm{reg}}}$ restricted to $(K^{\otimes g})^{\mathrm{inv}}$).

\medskip

\indent To conclude this section, we explain how to see $\mathcal{L}_{0,1}(H)$ as a coend. Recall the algebra $\mathcal{F}_{0,1}(H) \cong \mathcal{L}_{0,1}(H)$ from Remark \ref{F01}. We put a left $H$-module structure on it by letting $h \cdot \varphi = \varphi \cdot S^{-1}(h) = \varphi\!\left( S^{-1}(h'') ? h' \right)$. Since $h \cdot (\varphi \ast \psi) = (h'' \cdot \varphi) \ast (h' \cdot \psi)$, $\mathcal{F}_{0,1}(H)$ is an algebra in $\mathrm{mod}_l(H^{\mathrm{cop}})$, where $H^{\mathrm{cop}}$ is $H$ with opposite coproduct. Moreover, in $H^{\mathrm{cop}}$, we replace $\Delta$ by $\Delta^{\mathrm{op}}$, $R$ by $R'$ and $S$ by $S^{-1}$ so that the formulas for the product and the $H$-action in the coend of $\mathrm{mod}_l(H^{\mathrm{cop}})$ are exactly those of $\mathcal{F}_{0,1}(H)$. We state this as a proposition.
\begin{proposition}\label{L01Coend}
 It holds:
\[ \mathcal{L}_{0,1}(H) \cong \mathcal{F}_{0,1}(H) = \int^{X \in \mathrm{mod}_l(H^{\mathrm{cop}})} X^* \otimes X. \]
\end{proposition}

\subsection{Equivalence of the representations}
\indent Recall the map $F : H^* \to H^*$ 
\[ F(\varphi) = \varphi\!\left(a_i ? b_i\right), \:\:\:\:\: F^{-1}(\varphi) = \varphi\!\left(S^{-1}(a_i) ? b_i\right) \] 
(already used in the proof of Theorem \ref{formulesExplicites}) and let $\sigma : (H^*)^{\otimes g} \to (H^*)^{\otimes g}$ be the permutation
\[ \sigma(\varphi_1 \otimes \varphi_2 \otimes \ldots \otimes \varphi_{g-1} \otimes \varphi_g) = \varphi_g \otimes \varphi_{g-1} \otimes \ldots \otimes \varphi_2 \otimes \varphi_1. \]
It satisfies $\sigma^{-1} = \sigma$.
\begin{theorem}\label{thmEquivalenceReps}
The representation of Theorem \ref{thmRepMCG} and the Lyubashenko representation of $\mathrm{MCG}(\Sigma_g^{\mathrm{o}})$ and $\mathrm{MCG}(\Sigma_g)$ are equivalent. More precisely:\\
1) The isomorphism of vector spaces
\[ \fonc{(F \circ S)^{\otimes g} \circ \sigma}{K^{\otimes g}}{(H^*)^{\otimes g}}{\varphi_1 \otimes \ldots \otimes \varphi_g}{\varphi_g\!\left(b_i S(?) a_i\right) \otimes \ldots \otimes \varphi_1\!\left(b_i S(?) a_i\right)} \]
is an intertwiner between the two representations:
\[ \left[(F \circ S)^{\otimes g} \circ \sigma \right] \circ Z_{H_{\mathrm{reg}}}(f) = \theta_g^{\mathrm{o}}(f) \circ \left[(F \circ S)^{\otimes g} \circ \sigma \right]. \]
2) The isomorphism of vector spaces
\[ (F \circ S)^{\otimes g} \circ \sigma : (K^{\otimes g})^{\mathrm{inv}} \rightarrow \mathrm{Inv}\!\left((H^*)^{\otimes g}\right) \]
is an intertwiner between the two representations:
\[ \left[(F \circ S)^{\otimes g} \circ \sigma \right] \circ Z_{\mathbb{C}}(f) = \theta_g(f) \circ \left[(F \circ S)^{\otimes g} \circ \sigma \right]. \]
\end{theorem}
\begin{proof}
1) We show that this isomorphism intertwines the formulas of Theorem \ref{formulesExplicites} and of \eqref{formulesLyubashenko}. Thanks to the properties of $v$ \eqref{ribbon}, it is clear that $(F \circ S)^{\otimes g} \circ \sigma \circ Z_{H_{\mathrm{reg}}}(\tau_{a_i}) = \theta_g^{\mathrm{o}}(\tau_{a_i}) \circ (F \circ S)^{\otimes g} \circ \sigma$. Next, thanks to \eqref{actionvAvB}, \eqref{muLmuRgCarre} and \eqref{integraleShifte}, we have
\[ \theta_1^{\mathrm{o}}(\tau_b)(\varphi) = v_B^{-1} \triangleright \varphi = \mu^l(v)^{-1} \mu^r\!\left(gv^{-1} v'?\right) \varphi\!\left( vS^{-1}(gv'') \right). \]
Hence, for $\varphi \in H^*$,
\begin{align*}
\theta_1^{\mathrm{o}}(\tau_b) \circ (F \circ S)(\varphi) &= \mu^l(v)^{-1} \mu^r\!\left(gv^{-1}v'?\right) \varphi\!\left( v b_i g v'' a_i \right) = \mu^l(v)^{-1} \mu^r\!\left(g\overline{Y}_j?\right) \varphi\!\left( v^2 b_i g \overline{X}_j a_i \right)\\
&= \mu^l(v)^{-1} \mu^r\!\left(g S(a_j)S^{-1}(b_k)?\right) \varphi\!\left( v^2 g S^{-2}(b_i) b_j a_k a_i \right) = (\star)
\end{align*}
with $\overline{X}_i \otimes \overline{Y}_i = (RR')^{-1}$. We have a Yang-Baxter relation
\begin{align*}
S(a_j)S^{-1}(b_k) \otimes S^{-2}(b_i)b_j \otimes a_k a_i &= a_jS^{-1}(b_k) \otimes S^{-1}\!\left(b_j S^{-1}(b_i)\right) \otimes a_k a_i = (\mathrm{id} \otimes S^{-1} \otimes \mathrm{id})(R_{12} R_{31}^{-1} R_{32}^{-1})\\
&= (\mathrm{id} \otimes S^{-1} \otimes \mathrm{id})( R_{32}^{-1} R_{31}^{-1} R_{12}) = S^{-1}(b_k)S(a_j) \otimes b_j S^{-2}(b_i) \otimes a_i a_k
\end{align*}
which allows us to continue the computation:
\begin{align*}
(\star) &= \mu^l(v)^{-1} \mu^r\!\left(g S^{-1}(b_k)S(a_j)?\right) \varphi\!\left( v^2 g b_j S^{-2}(b_i) a_i a_k \right) = \mu^l(v)^{-1} \mu^r\!\left(g S^{-1}(b_k)S(a_j)?\right) \varphi\!\left( v S^2(b_j) a_k \right)\\
&= \mu^l(v)^{-1} \mu^r\!\left(g S^{-1}(a_j b_k)?\right) \varphi\!\left( v b_j a_k \right) = \mu^l(v)^{-1} \mu^r\!\left(g v S^{-1}(v''^{-1})?\right) \varphi\!\left( v^2 v'^{-1} \right).
\end{align*}
We used \eqref{elementDrinfeld} and \eqref{ribbon}. On the other hand, we compute
\begin{align*}
(F \circ S) \circ Z_{H_{\mathrm{reg}}}(\tau_b)(\varphi) &= (F \circ S) \circ (\mathcal{S}^{-1} \circ \mathcal{T} \circ \mathcal{S})(\varphi) = (F \circ S) \circ (g^{-1} v)_* \circ \rho(v_B^{-1}) \circ (g v^{-1})_*(\varphi)\\
&= F \circ S\!\left( \mu^l(v)^{-1} \mu^r\!\left( v'?\right) \varphi\!\left( S^{-1}(v'') \right) \right) = \mu^l(v)^{-1} \mu^r\!\left( v' b_i S(?) a_i\right) \varphi\!\left( S^{-1}(v'') \right)\\
&= \mu^l(v)^{-1} \mu^r\!\left( v S^2(a_i) \overline{Y}_j b_i S(?)\right) \varphi\!\left( v S^{-1}(\overline{X}_i) \right)\\
&= \mu^l(v)^{-1} \mu^r\!\left( v S^2(a_i) S(a_j) S^{-1}(b_k) b_i S(?)\right) \varphi\!\left( v S^{-1}(b_ja_k) \right) = (\star\star)
\end{align*}
where $\rho$ is the representation of $\mathcal{L}_{g,0}(H)$ on $(H^*)^{\otimes g}$. We used Lemma \ref{lemmaST}, \eqref{quasiCyclic} and \eqref{ribbon}. As previously, we have a Yang-Baxter relation
\[ S^2(a_i) S(a_j) \otimes S^{-1}(b_k) b_i \otimes b_ja_k = S(a_j)S^2(a_i) \otimes b_i S^{-1}(b_k) \otimes a_k b_j \]
which allows us to continue the computation:
\begin{align*}
(\star \star) &= \mu^l(v)^{-1} \mu^r\!\left( v S(a_j) S^2(a_i) b_i S^{-1}(b_k) S(?)\right) \varphi\!\left( v S^{-1}(a_k b_j) \right)\\
& = \mu^l(v)^{-1} \mu^r\!\left(  S(a_j) g S^{-1}(b_k) S(?)\right) \varphi\!\left( v S^{-1}(a_k b_j) \right) = \mu^l(v)^{-1} \mu^r\!\left( g a_j b_k S(?)\right) \varphi\!\left( v b_j a_k \right)\\
&= \mu^l(v)^{-1} \mu^r\!\left( g v v''^{-1} S(?)\right) \varphi\!\left( v^2 v'^{-1} \right) = \mu^l(v)^{-1} \mu^r \circ S\!\left( ? S^{-1}(v''^{-1})vg^{-1}\right) \varphi\!\left( v^2 v'^{-1} \right)\\
&= \mu^l(v)^{-1} \mu^l\!\left( ? S^{-1}(v''^{-1})vg^{-1}\right) \varphi\!\left( v^2 v'^{-1} \right) = \mu^l(v)^{-1} \mu^r\!\left( gv S^{-1}(v''^{-1})? \right) \varphi\!\left( v^2 v'^{-1} \right).
\end{align*}
We used \eqref{elementDrinfeld} to simplify $S^2(a_i) b_i = S(S^{-1}(b_i)S(a_i)) = gv^{-1}$ and the properties of $\mu^l$ and $\mu^r$ recorded in section \ref{rappelHopf}. Hence, it holds $\theta_1^{\mathrm{o}}(\tau_b) \circ (F \circ S) = (F \circ S) \circ Z_{H_{\mathrm{reg}}}(\tau_b)$, which clearly implies that $\theta_g^{\mathrm{o}}(\tau_{b_i}) \circ (F \circ S)^{\otimes g} \circ \sigma = (F \circ S)^{\otimes g} \circ \sigma \circ Z_{H_{\mathrm{reg}}}(\tau_{b_i})$. Let us now proceed with $\tau_{d_i}$ ($i \geq 2$):
\begin{align*}
 &(F \circ S)^{\otimes g} \circ \sigma \circ Z_{H_{\mathrm{reg}}}(\tau_{d_i}) \circ \sigma \circ (S^{-1} \circ F^{-1})^{\otimes g}\!\left(\varphi_1 \otimes \ldots \varphi_g\right)\\
&=(F \circ S)^{\otimes g} \circ \sigma \circ Z_{H_{\mathrm{reg}}}(\tau_{d_i})\!\left( \varphi_g\!\left(S^{-1}(a_j) S^{-1}(?) b_j\right) \otimes \ldots \otimes \varphi_1\!\left(S^{-1}(a_j) S^{-1}(?) b_j\right) \right)\\
&= (F \circ S)^{\otimes g} \circ \sigma \!\left( \varphi_g\!\left(S^{-1}(a_j) S^{-1}(?) b_j\right) \otimes \ldots \otimes \varphi_i\!\left(S^{-1}(a_j) S^{-1}(v'^{-1})S^{-1}(?) b_j\right) \otimes  \varphi_{i-1}\!\left(S^{-1}(a_j) S^{-1}(?) v''^{-1} b_j\right) \right.\\
&\:\:\:\:\:\:\:\:\:\:\:\:\:\:\:\:\:\:\:\:\:\:\:\:\:\:\:\:\:\:\:\:\: \left. \otimes \ldots \otimes \varphi_1\!\left(S^{-1}(a_j) S^{-1}(?) b_j\right) \right)\\
&=  \varphi_1 \otimes \ldots \otimes \varphi_{i-1}\!\left(S^{-1}(a_j)a_k ? b_k v''^{-1} b_j \right) \otimes \varphi_i\!\left(S^{-1}(a_j) S^{-1}(v'^{-1})a_k ? b_k b_j \right) \otimes \ldots \otimes \varphi_g\\
&= \theta_g^{\mathrm{o}}(\tau_{d_i})\bigl(\varphi_1 \otimes \ldots \varphi_g\bigr).
\end{align*}
Finally, for $\tau_{e_i}$ ($i \geq 2$):
\begin{align*}
 &(F \circ S)^{\otimes g} \circ \sigma \circ Z_{H_{\mathrm{reg}}}(\tau_{e_i}) \circ \sigma \circ (S^{-1} \circ F^{-1})^{\otimes g}\!\left(\varphi_1 \otimes \ldots \varphi_g\right)\\
&= (F \circ S)^{\otimes g} \circ \sigma \circ Z_{H_{\mathrm{reg}}}(\tau_{e_i})\!\left( \varphi_g\!\left(S^{-1}(a_j) S^{-1}(?) b_j \right) \otimes \ldots \otimes \varphi_1\!\left(S^{-1}(a_j) S^{-1}(?) b_j\right) \right)\\
&= (F \circ S)^{\otimes g} \circ \sigma\!\left( \varphi_g\!\left(S^{-1}(a_j) S^{-1}(?) b_j\right) \otimes \ldots \otimes \varphi_{i+1}\!\left(S^{-1}(a_j) S^{-1}(?) b_j\right) \otimes \varphi_{i}\!\left(S^{-1}(a_j)S^{-1}(v^{(1)-1}) S^{-1}(?) b_j\right) \right. \\
& \:\:\:\:\:\:\:\:\:\:\:\:\:\:\:\:\: \left.  \otimes \, \varphi_{i-1}\!\left( S^{-1}(a_j) S^{-1}(v^{(3)-1}) S^{-1}(?) v^{(2)-1} b_j \right) \otimes \ldots \otimes \varphi_1\!\left( S^{-1}(a_j) S^{-1}(v^{(2i-1)-1})  S^{-1}(?) v^{(2i-2)-1} b_j \right) \right)\\
&= \varphi_1\!\left( S^{-1}(a_j) S^{-1}(v^{(2i-1)-1}) a_k ? b_k v^{(2i-2)-1} b_j \right) \otimes \ldots \otimes \varphi_{i-1}\!\left( S^{-1}(a_j) S^{-1}(v^{(3)-1}) a_k ? b_k v^{(2)-1} b_j \right) \\
& \:\:\:\: \otimes \varphi_{i}\!\left(S^{-1}(a_j)S^{-1}(v^{(1)-1}) a_k ? b_k b_j\right) \otimes \varphi_{i+1} \otimes \ldots \otimes \varphi_g\\
&= \varphi_1\!\left( S^{-1}(v^{(2i-2)-1})  ? v^{(2i-1)-1} \right) \otimes \ldots \otimes \varphi_{i-1}\!\left( S^{-1}(v^{(2)-1}) ? v^{(3)-1} \right) \otimes \varphi_{i}\!\left(S^{-1}(a_j)S^{-1}(v^{(1)-1}) a_k ? b_k b_j\right) \otimes \varphi_{i+1}\\
&\:\:\:\: \otimes \ldots \otimes \varphi_g\\
&= \theta_g^{\mathrm{o}}(\tau_{e_i})\bigl(\varphi_1 \otimes \ldots \otimes \varphi_g\bigr).
\end{align*}
We used $\Delta^{\mathrm{op}}R = R\Delta$ for the last equality.\\
\noindent 2) It is not difficult to see that $(F \circ S)^{\otimes g} \circ \sigma : K^{\otimes g} \to (H^*)^{\otimes g}$ is a morphism of $H$-modules, where $K^{\otimes g}$ is endowed with the action \eqref{actionHKg} and $(H^*)^{\otimes g}$ is endowed with the action \eqref{actionH} (with $n=0$). Hence, the restriction of $(F \circ S)^{\otimes g} \circ \sigma$ to $(K^{\otimes g})^{\mathrm{inv}}$ indeed takes values in $\mathrm{Inv}\bigl((H^*)^{\otimes g}\bigr)$. Since 
\[ Z_{\mathbb{C}}(f) = \left(Z_{H_{\mathrm{reg}}}(f)\right)_{\bigl| (K^{\otimes g})^{\mathrm{inv}}\bigr.} \:\:\: \text{ and } \:\:\: \theta_g(f) = \theta_g^{\mathrm{o}}(f)_{\bigl| \mathrm{Inv}\left((H^*)^{\otimes g}\right)\bigr.}, \]
the result follows from the first part of the theorem.
\end{proof}

\chapter{Graphical calculus and relation to skein theory}\label{chapitreGraphiqueSkein}
\indent The main topic of this last chapter is to explain how to compute in $\mathcal{L}_{g,n}(H)$ in a graphical way and to define the Wilson loop map $W$. For that, we introduce a graphical element associated to a matrix with coefficients in $\mathcal{L}_{g,n}(H)$ and combine it with the ones recalled in section \ref{modH}. Hence, this can be seen as an extension of the Reshetikhin-Turaev functor (section \ref{modH}), even though the map ``evaluation of a diagram'' is no longer a functor (the evaluation of a diagram is an element of $\mathcal{L}_{g,n}(H) \otimes V$, where $V$ is a $H$-module). Note that such a graphical calculus can be used to compute in any algebra defined by means of matrix coefficients associated to a braided Hopf algebra (like $H$, $\mathcal{O}(H)$, $\mathcal{H}(\mathcal{O}(H))$).

\smallskip

\indent We first reformulate graphically the defining relations of $\mathcal{L}_{g,n}(H)$. Then we graphically define the Wilson loop map $W$, which to an oriented, colored and framed link in $\Sigma_{g,n}^{\mathrm{o}}$ associates an element of $\mathcal{L}_{g,n}^{\mathrm{inv}}(H)$. This definition and the resulting properties are equivalent, but maybe simpler, to the Wilson loops of \cite{BR2}, whose definition was based on chord diagrams and did not used a graphical presentation of $\mathcal{L}_{g,n}(H)$, and to the Wilson loops of \cite{BFKB, BFKB2}, whose formal definition in the setting of $\mathcal{F}_{g,n}(H)$ (functions on connections, see section \ref{lienAvecLesLGFT}) was based on their notion of multitangles. Here we choose the canonical thickened graph of Figure \ref{surfaceAvecMatrices}; with this choice, the definition of the Wilson loop map becomes obvious and natural-looking, thanks to the use of the graphical calculus. 

\smallskip

\indent With the gauge algebra $H = \bar U_q$, we use the Wilson loop map and its particular properties in this setting to obtain representations of skein algebras. Note that, as in the case of the representations of mapping class groups, the restriction to $\mathrm{Inv}(V)$ allows us to glue back the disc $D$ and to obtain a representation of $\mathcal{S}_q(\Sigma_{g,0})$ (and not just of $\mathcal{S}_q(\Sigma_{g,0}^{\mathrm{o}})$).

\smallskip

\indent The main results of this chapter are
\begin{itemize}
\item The definition of the Wilson loop map $W$ (Definition \ref{defDefWilson}) and its natural and expected properties (Theorem \ref{wilsonStack}, Propositions \ref{wilsonInv}, \ref{orientationEtWilson}, \ref{propWilsonSimple}), which indicate that the definition of $W$ is the good one.
\item The representation of the skein algebra $\mathcal{S}_q(\Sigma_{g,0})$ on $\mathrm{Inv}\bigl( (\bar U_q)^{\otimes g} \bigr)$ (Theorem \ref{theoRepSkein}).
\item The explicit study of the representation of $\mathcal{S}_q(\Sigma_1)$ on $\mathrm{SLF}(\bar U_q)$ (Propositions \ref{structureRepSqS1} and \ref{propLienSkeinEtL10}).
\end{itemize}

\smallskip

\indent As previously, we assume that $H$ is a finite dimensional factorizable ribbon Hopf algebra, even though these assumptions can be weakened for this chapter.

\section{Diagrammatic description of $\mathcal{L}_{g,n}(H)$}

\indent Let $I$ be a $H$-module. We have an isomorphism of vector spaces :
\begin{equation*}
\begin{array}{rcl} 
\mathcal{L}_{g,n}(H) \otimes I \otimes I^* & \overset{\sim}{\rightarrow} & \mathcal{L}_{g,n}(H) \otimes \mathrm{End}_{\mathbb{C}}(I)\\
x \otimes u \otimes \varphi & \mapsto & x \otimes \left(y \mapsto \varphi(y)u\right).
\end{array} 
\end{equation*}
Let us choose a basis $(v_i)$ of $I$ and let $(v^j)$ be its dual basis. Then $\mathcal{L}_{g,n}(H) \otimes \mathrm{End}_{\mathbb{C}}(I) \cong \mathcal{L}_{g,n}(H) \otimes \mathrm{Mat}_{\dim(I)}(\mathbb{C}) = \mathrm{Mat}_{\dim(I)}\!\left(\mathcal{L}_{g,n}(H)\right)$ and the inverse of the above isomorphism is
\begin{equation}\label{decMatrice}
\begin{array}{rcl} 
\mathrm{Mat}_{\dim(I)}\!\left(\mathcal{L}_{g,n}(H)\right) & \overset{\sim}{\rightarrow} & \mathcal{L}_{g,n}(H) \otimes I \otimes I^* \\
M & \mapsto & M^i_j  \otimes v_i \otimes v^j.
\end{array} 
\end{equation}
In this chapter, we systematically identify a matrix $M \in \mathrm{Mat}_{\dim(I)}\!\left(\mathcal{L}_{g,n}(H)\right)$ with $M^i_j \otimes v_i \otimes v^j$, written more shortly $M^i_j \, v_i \otimes v^j$.

\smallskip

\indent We denote by $\overset{I}{X}$ an element of $\mathcal{L}_{g,n}(H) \otimes \mathrm{End}_{\mathbb{C}}(I) = \mathrm{Mat}_{\dim(I)}\!\left(\mathcal{L}_{g,n}(H)\right)$. In general, we will restrict $\overset{I}{X}$ to be a product of the matrices of generators $\overset{I}{A}(i), \overset{I}{B}(j), \overset{I}{M}(k)$ up to some normalization $\overset{I}{v}{^r}$, namely:
\begin{equation}\label{matriceGenerale}
\overset{I}{v}{^r}\overset{I}{A}(i_1)^{l_1} \overset{I}{B}(j_1)^{m_1} \overset{I}{M}(k_1)^{n_1} \ldots \overset{I}{A}(i_s)^{l_s} \overset{I}{B}(j_s)^{m_s} \overset{I}{M}(k_s)^{n_s} \in \mathrm{Mat}_{\dim(I)}\!\left(\mathcal{L}_{g,n}(H)\right)
\end{equation}
with $r, l_{\alpha}, m_{\alpha}, n_{\alpha} \in \mathbb{Z}$ and $1 \leq i_{\alpha}, j_{\alpha} \leq g, \: g+1 \leq k_{\alpha} \leq g+n$; for instance $\overset{I}{M}(3)^{-2}\overset{I}{B}(1){^{-1}}\overset{I}{A}(2)$. Using the identification \eqref{decMatrice}, we represent graphically $\overset{I}{X} = \overset{I}{X}{^i_j} \, v_i \otimes v^j$ by the following diagram:

\begin{equation}\label{anse}
\begingroup%
  \makeatletter%
  \providecommand\color[2][]{%
    \errmessage{(Inkscape) Color is used for the text in Inkscape, but the package 'color.sty' is not loaded}%
    \renewcommand\color[2][]{}%
  }%
  \providecommand\transparent[1]{%
    \errmessage{(Inkscape) Transparency is used (non-zero) for the text in Inkscape, but the package 'transparent.sty' is not loaded}%
    \renewcommand\transparent[1]{}%
  }%
  \providecommand\rotatebox[2]{#2}%
  \newcommand*\fsize{\dimexpr\f@size pt\relax}%
  \newcommand*\lineheight[1]{\fontsize{\fsize}{#1\fsize}\selectfont}%
  \ifx\svgwidth\undefined%
    \setlength{\unitlength}{148.68101591bp}%
    \ifx\svgscale\undefined%
      \relax%
    \else%
      \setlength{\unitlength}{\unitlength * \real{\svgscale}}%
    \fi%
  \else%
    \setlength{\unitlength}{\svgwidth}%
  \fi%
  \global\let\svgwidth\undefined%
  \global\let\svgscale\undefined%
  \makeatother%
  \begin{picture}(1,0.41355276)%
    \lineheight{1}%
    \setlength\tabcolsep{0pt}%
    \put(0,0){\includegraphics[width=\unitlength,page=1]{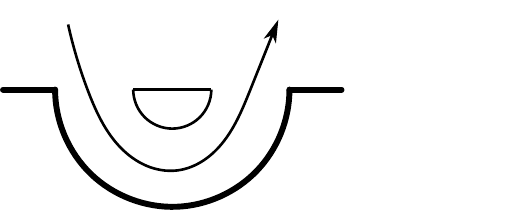}}%
    \put(0.06826181,0.37682666){\color[rgb]{0,0,0}\makebox(0,0)[lt]{\lineheight{1.25}\smash{\begin{tabular}[t]{l}$I$\end{tabular}}}}%
    \put(0.52430761,0.010939){\color[rgb]{0,0,0}\makebox(0,0)[lt]{\lineheight{1.25}\smash{\begin{tabular}[t]{l}$\overset{I}{X}$\end{tabular}}}}%
    \put(0,0){\includegraphics[width=\unitlength,page=2]{anse.pdf}}%
  \end{picture}%
\endgroup%

\end{equation}

\noindent The module $I$ colors the strand while the matrix $\overset{I}{X}$ colors the handle. Mimicking \eqref{dualitePrime}, we define a graphical element corresponding to the negative orientation of the strand:
\begin{equation}\label{sensOppose}
\begingroup%
  \makeatletter%
  \providecommand\color[2][]{%
    \errmessage{(Inkscape) Color is used for the text in Inkscape, but the package 'color.sty' is not loaded}%
    \renewcommand\color[2][]{}%
  }%
  \providecommand\transparent[1]{%
    \errmessage{(Inkscape) Transparency is used (non-zero) for the text in Inkscape, but the package 'transparent.sty' is not loaded}%
    \renewcommand\transparent[1]{}%
  }%
  \providecommand\rotatebox[2]{#2}%
  \newcommand*\fsize{\dimexpr\f@size pt\relax}%
  \newcommand*\lineheight[1]{\fontsize{\fsize}{#1\fsize}\selectfont}%
  \ifx\svgwidth\undefined%
    \setlength{\unitlength}{285.39483121bp}%
    \ifx\svgscale\undefined%
      \relax%
    \else%
      \setlength{\unitlength}{\unitlength * \real{\svgscale}}%
    \fi%
  \else%
    \setlength{\unitlength}{\svgwidth}%
  \fi%
  \global\let\svgwidth\undefined%
  \global\let\svgscale\undefined%
  \makeatother%
  \begin{picture}(1,0.28793971)%
    \lineheight{1}%
    \setlength\tabcolsep{0pt}%
    \put(0,0){\includegraphics[width=\unitlength,page=1]{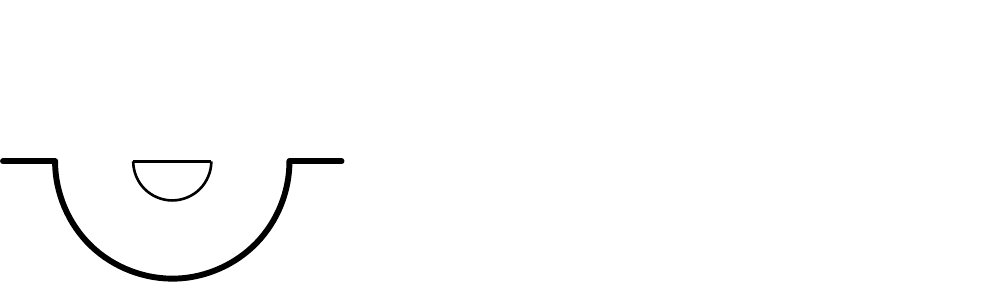}}%
    \put(0.27314648,0.00569885){\color[rgb]{0,0,0}\makebox(0,0)[lt]{\lineheight{1.25}\smash{\begin{tabular}[t]{l}$\overset{I}{X}$\end{tabular}}}}%
    \put(0,0){\includegraphics[width=\unitlength,page=2]{sens_oppose.pdf}}%
    \put(0.71456752,0.00569885){\color[rgb]{0,0,0}\makebox(0,0)[lt]{\lineheight{1.25}\smash{\begin{tabular}[t]{l}$\overset{I^*}{X}$\end{tabular}}}}%
    \put(0,0){\includegraphics[width=\unitlength,page=3]{sens_oppose.pdf}}%
    \put(0.37787999,0.11677759){\color[rgb]{0,0,0}\makebox(0,0)[lt]{\lineheight{1.25}\smash{\begin{tabular}[t]{l}=\end{tabular}}}}%
    \put(0,0){\includegraphics[width=\unitlength,page=4]{sens_oppose.pdf}}%
    \put(0.66750775,0.19272084){\color[rgb]{0,0,0}\makebox(0,0)[lt]{\lineheight{1.25}\smash{\begin{tabular}[t]{l}$e_I$\end{tabular}}}}%
    \put(0,0){\includegraphics[width=\unitlength,page=5]{sens_oppose.pdf}}%
    \put(0.50148027,0.19176672){\color[rgb]{0,0,0}\makebox(0,0)[lt]{\lineheight{1.25}\smash{\begin{tabular}[t]{l}$\mathrm{id}_{I^*}$\end{tabular}}}}%
    \put(0.4960813,0.25591023){\color[rgb]{0,0,0}\makebox(0,0)[lt]{\lineheight{1.25}\smash{\begin{tabular}[t]{l}$I$\end{tabular}}}}%
    \put(0.69708759,0.25786879){\color[rgb]{0,0,0}\makebox(0,0)[lt]{\lineheight{1.25}\smash{\begin{tabular}[t]{l}$I$\end{tabular}}}}%
    \put(0.48669005,0.13817477){\color[rgb]{0,0,0}\makebox(0,0)[lt]{\lineheight{1.25}\smash{\begin{tabular}[t]{l}$I^*$\end{tabular}}}}%
    \put(0,0){\includegraphics[width=\unitlength,page=6]{sens_oppose.pdf}}%
    \put(0.26232042,0.16461482){\color[rgb]{0,0,0}\makebox(0,0)[lt]{\lineheight{1.25}\smash{\begin{tabular}[t]{l}$I$\end{tabular}}}}%
  \end{picture}%
\endgroup%

\end{equation}

\noindent where $e_I: I^{**} \to I$ is the isomorphism \eqref{identificationBidual}. Let us explain \eqref{sensOppose}. To define the graphical element on the left, we put a ribbon graph atop the one defined in \eqref{anse}. This ribbon graph represents a morphism $I^* \otimes I^{**} \to I^* \otimes I$ in $\mathrm{mod}_l(H)$ (see section \ref{modH}), which can be applied to $\overset{I^*}{X}{^i_j} \, v^i \otimes \langle ?, v_j \rangle$ and thus gives a well-defined element in $\mathcal{L}_{g,n}(H) \otimes I^* \otimes I$. For further use we record that, due to \eqref{propCoherenceBidual} and \eqref{naturalite} below, the converse of \eqref{sensOppose} is
\begin{equation}\label{sensOpposeBis}
\begingroup%
  \makeatletter%
  \providecommand\color[2][]{%
    \errmessage{(Inkscape) Color is used for the text in Inkscape, but the package 'color.sty' is not loaded}%
    \renewcommand\color[2][]{}%
  }%
  \providecommand\transparent[1]{%
    \errmessage{(Inkscape) Transparency is used (non-zero) for the text in Inkscape, but the package 'transparent.sty' is not loaded}%
    \renewcommand\transparent[1]{}%
  }%
  \providecommand\rotatebox[2]{#2}%
  \newcommand*\fsize{\dimexpr\f@size pt\relax}%
  \newcommand*\lineheight[1]{\fontsize{\fsize}{#1\fsize}\selectfont}%
  \ifx\svgwidth\undefined%
    \setlength{\unitlength}{285.39482813bp}%
    \ifx\svgscale\undefined%
      \relax%
    \else%
      \setlength{\unitlength}{\unitlength * \real{\svgscale}}%
    \fi%
  \else%
    \setlength{\unitlength}{\svgwidth}%
  \fi%
  \global\let\svgwidth\undefined%
  \global\let\svgscale\undefined%
  \makeatother%
  \begin{picture}(1,0.28793973)%
    \lineheight{1}%
    \setlength\tabcolsep{0pt}%
    \put(0,0){\includegraphics[width=\unitlength,page=1]{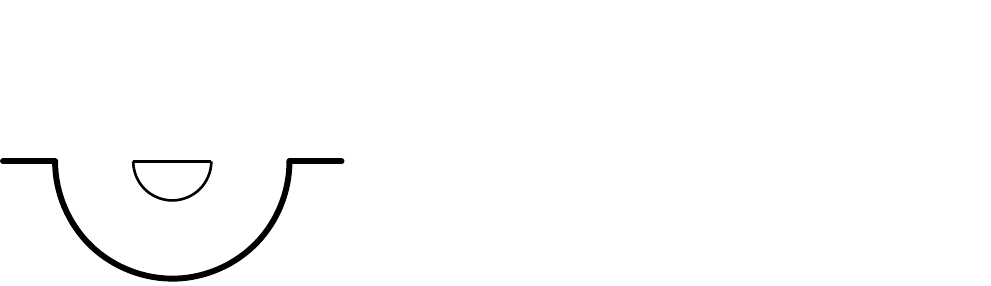}}%
    \put(0.27314648,0.00569885){\color[rgb]{0,0,0}\makebox(0,0)[lt]{\lineheight{1.25}\smash{\begin{tabular}[t]{l}$\overset{I}{X}$\end{tabular}}}}%
    \put(0,0){\includegraphics[width=\unitlength,page=2]{sens_oppose_bis.pdf}}%
    \put(0.71456752,0.00569885){\color[rgb]{0,0,0}\makebox(0,0)[lt]{\lineheight{1.25}\smash{\begin{tabular}[t]{l}$\overset{I^*}{X}$\end{tabular}}}}%
    \put(0,0){\includegraphics[width=\unitlength,page=3]{sens_oppose_bis.pdf}}%
    \put(0.37787998,0.11677761){\color[rgb]{0,0,0}\makebox(0,0)[lt]{\lineheight{1.25}\smash{\begin{tabular}[t]{l}=\end{tabular}}}}%
    \put(0,0){\includegraphics[width=\unitlength,page=4]{sens_oppose_bis.pdf}}%
    \put(0.50983132,0.19272086){\color[rgb]{0,0,0}\makebox(0,0)[lt]{\lineheight{1.25}\smash{\begin{tabular}[t]{l}$e_I$\end{tabular}}}}%
    \put(0,0){\includegraphics[width=\unitlength,page=5]{sens_oppose_bis.pdf}}%
    \put(0.65915654,0.19176672){\color[rgb]{0,0,0}\makebox(0,0)[lt]{\lineheight{1.25}\smash{\begin{tabular}[t]{l}$\mathrm{id}_{I^*}$\end{tabular}}}}%
    \put(0.6537576,0.25591023){\color[rgb]{0,0,0}\makebox(0,0)[lt]{\lineheight{1.25}\smash{\begin{tabular}[t]{l}$I$\end{tabular}}}}%
    \put(0.53941116,0.25786881){\color[rgb]{0,0,0}\makebox(0,0)[lt]{\lineheight{1.25}\smash{\begin{tabular}[t]{l}$I$\end{tabular}}}}%
    \put(0,0){\includegraphics[width=\unitlength,page=6]{sens_oppose_bis.pdf}}%
    \put(0.06175625,0.16049307){\color[rgb]{0,0,0}\makebox(0,0)[lt]{\lineheight{1.25}\smash{\begin{tabular}[t]{l}$I$\end{tabular}}}}%
    \put(0,0){\includegraphics[width=\unitlength,page=7]{sens_oppose_bis.pdf}}%
    \put(0.69296005,0.13671711){\color[rgb]{0,0,0}\makebox(0,0)[lt]{\lineheight{1.25}\smash{\begin{tabular}[t]{l}$I^*$\end{tabular}}}}%
  \end{picture}%
\endgroup%

\end{equation}

\indent The tensor product of two matrices $\overset{I}{X}, \overset{J}{Y}$, defined by
\[ \overset{I}{X} \otimes \overset{J}{Y} = \overset{I}{X}{^i_j} \overset{J}{Y}{^k_l} \, v_i \otimes v^j \otimes w_k \otimes w^l \]
(where $(v_i)$ is a basis of $I$, $(w_k)$ is a basis of $J$ and $(v^j), (w^l)$ are their respective dual bases), is represented by the gluing of the corresponding graphical elements:
\begin{equation}\label{tensorProd}
\begingroup%
  \makeatletter%
  \providecommand\color[2][]{%
    \errmessage{(Inkscape) Color is used for the text in Inkscape, but the package 'color.sty' is not loaded}%
    \renewcommand\color[2][]{}%
  }%
  \providecommand\transparent[1]{%
    \errmessage{(Inkscape) Transparency is used (non-zero) for the text in Inkscape, but the package 'transparent.sty' is not loaded}%
    \renewcommand\transparent[1]{}%
  }%
  \providecommand\rotatebox[2]{#2}%
  \newcommand*\fsize{\dimexpr\f@size pt\relax}%
  \newcommand*\lineheight[1]{\fontsize{\fsize}{#1\fsize}\selectfont}%
  \ifx\svgwidth\undefined%
    \setlength{\unitlength}{245.58551116bp}%
    \ifx\svgscale\undefined%
      \relax%
    \else%
      \setlength{\unitlength}{\unitlength * \real{\svgscale}}%
    \fi%
  \else%
    \setlength{\unitlength}{\svgwidth}%
  \fi%
  \global\let\svgwidth\undefined%
  \global\let\svgscale\undefined%
  \makeatother%
  \begin{picture}(1,0.25187776)%
    \lineheight{1}%
    \setlength\tabcolsep{0pt}%
    \put(0,0){\includegraphics[width=\unitlength,page=1]{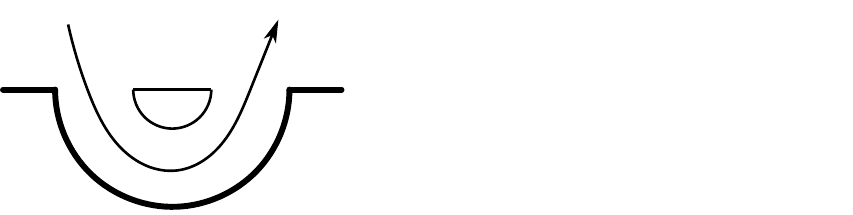}}%
    \put(0.05040422,0.22964325){\color[rgb]{0,0,0}\makebox(0,0)[lt]{\lineheight{1.25}\smash{\begin{tabular}[t]{l}$I$\end{tabular}}}}%
    \put(0.31379241,0.00812957){\color[rgb]{0,0,0}\makebox(0,0)[lt]{\lineheight{1.25}\smash{\begin{tabular}[t]{l}$\overset{I}{X}$\end{tabular}}}}%
    \put(0,0){\includegraphics[width=\unitlength,page=2]{produit_Kronecker.pdf}}%
    \put(0.71443375,0.00812957){\color[rgb]{0,0,0}\makebox(0,0)[lt]{\lineheight{1.25}\smash{\begin{tabular}[t]{l}$\overset{J}{Y}$\end{tabular}}}}%
    \put(0,0){\includegraphics[width=\unitlength,page=3]{produit_Kronecker.pdf}}%
    \put(0.44796585,0.22889765){\color[rgb]{0,0,0}\makebox(0,0)[lt]{\lineheight{1.25}\smash{\begin{tabular}[t]{l}$J$\end{tabular}}}}%
  \end{picture}%
\endgroup%

\end{equation}

\begin{figure}[h]
\centering
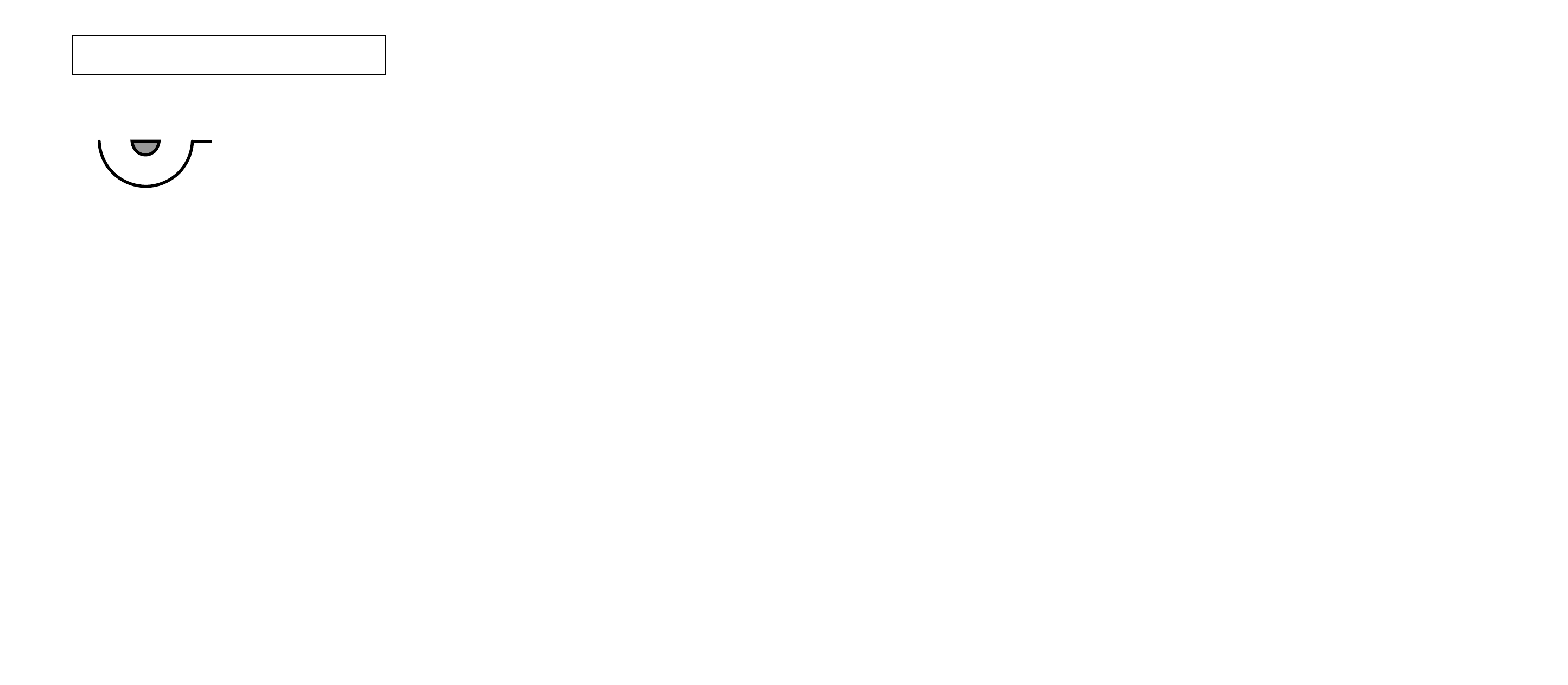
\caption{The evaluation map $\widetilde{F}_{\mathrm{RT}}$ (the double arrows in the first equality mean any orientation).}
\label{figureDefEvaluation}
\end{figure}

\begin{definition}\label{defEvaluation}
A diagram is obtained by gluing (as in \eqref{tensorProd}) several copies of the handle diagrams introduced in \eqref{anse} and \eqref{sensOppose}, and by putting atop an oriented and colored ribbon graph $G$ (see section \ref{modH}). The evaluation of a diagram, depicted in Figure \ref{figureDefEvaluation}, is a map $\widetilde{F}_{\mathrm{RT}}$ which consists of applying $F_{\mathrm{RT}}(G)$ to the matrices associated to the handle diagrams introduced previously, where $F_{\mathrm{RT}}$ is the Reshetikhin-Turaev functor (see section \ref{modH}).  The evaluation of a diagram is an element of $\mathcal{L}_{g,n}(H) \otimes J_1 \otimes \ldots \otimes J_l$, where $J_1, \ldots, J_l$ are $H$-modules.
\end{definition}
\noindent In the sequel, we always identify a diagram with its evaluation through $\widetilde{F}_{\mathrm{RT}}$.

\smallskip

\indent For instance, consider the following diagram:
\begin{center}
\begingroup%
  \makeatletter%
  \providecommand\color[2][]{%
    \errmessage{(Inkscape) Color is used for the text in Inkscape, but the package 'color.sty' is not loaded}%
    \renewcommand\color[2][]{}%
  }%
  \providecommand\transparent[1]{%
    \errmessage{(Inkscape) Transparency is used (non-zero) for the text in Inkscape, but the package 'transparent.sty' is not loaded}%
    \renewcommand\transparent[1]{}%
  }%
  \providecommand\rotatebox[2]{#2}%
  \newcommand*\fsize{\dimexpr\f@size pt\relax}%
  \newcommand*\lineheight[1]{\fontsize{\fsize}{#1\fsize}\selectfont}%
  \ifx\svgwidth\undefined%
    \setlength{\unitlength}{245.58549787bp}%
    \ifx\svgscale\undefined%
      \relax%
    \else%
      \setlength{\unitlength}{\unitlength * \real{\svgscale}}%
    \fi%
  \else%
    \setlength{\unitlength}{\svgwidth}%
  \fi%
  \global\let\svgwidth\undefined%
  \global\let\svgscale\undefined%
  \makeatother%
  \begin{picture}(1,0.24220106)%
    \lineheight{1}%
    \setlength\tabcolsep{0pt}%
    \put(0,0){\includegraphics[width=\unitlength,page=1]{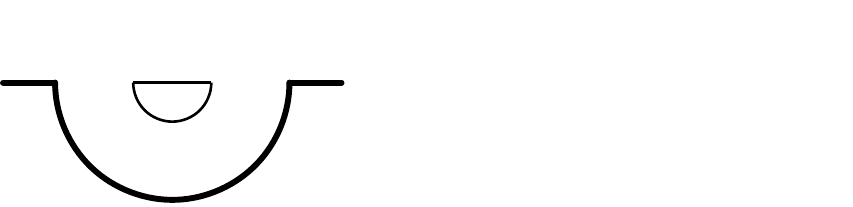}}%
    \put(0.04314224,0.21996655){\color[rgb]{0,0,0}\makebox(0,0)[lt]{\lineheight{1.25}\smash{\begin{tabular}[t]{l}$I$\end{tabular}}}}%
    \put(0.31379245,0.00662263){\color[rgb]{0,0,0}\makebox(0,0)[lt]{\lineheight{1.25}\smash{\begin{tabular}[t]{l}$\overset{I}{X}$\end{tabular}}}}%
    \put(0,0){\includegraphics[width=\unitlength,page=2]{produit_matriciel.pdf}}%
    \put(0.71443374,0.00662263){\color[rgb]{0,0,0}\makebox(0,0)[lt]{\lineheight{1.25}\smash{\begin{tabular}[t]{l}$\overset{I}{Y}$\end{tabular}}}}%
    \put(0,0){\includegraphics[width=\unitlength,page=3]{produit_matriciel.pdf}}%
  \end{picture}%
\endgroup%

\end{center}
It is evaluated as
\[ \overset{I}{X}{^i_j} \overset{I}{Y}{^k_l} \, \mathrm{id}_I \otimes d_I \otimes \mathrm{id}_{I^*}\!\left(v_i \otimes v^j \otimes v_k \otimes v^l\right) = \overset{I}{X}{^i_j} \overset{I}{Y}{^j_l} \, v_i \otimes v^l =   (\overset{I}{X} \overset{I}{Y}){^i_l} \, v_i \otimes v^l. \]
Hence, we see that this diagram represents the matrix product $\overset{I}{X} \overset{I}{Y}$. Similarly, the diagram
\begin{center}
\begingroup%
  \makeatletter%
  \providecommand\color[2][]{%
    \errmessage{(Inkscape) Color is used for the text in Inkscape, but the package 'color.sty' is not loaded}%
    \renewcommand\color[2][]{}%
  }%
  \providecommand\transparent[1]{%
    \errmessage{(Inkscape) Transparency is used (non-zero) for the text in Inkscape, but the package 'transparent.sty' is not loaded}%
    \renewcommand\transparent[1]{}%
  }%
  \providecommand\rotatebox[2]{#2}%
  \newcommand*\fsize{\dimexpr\f@size pt\relax}%
  \newcommand*\lineheight[1]{\fontsize{\fsize}{#1\fsize}\selectfont}%
  \ifx\svgwidth\undefined%
    \setlength{\unitlength}{148.68101445bp}%
    \ifx\svgscale\undefined%
      \relax%
    \else%
      \setlength{\unitlength}{\unitlength * \real{\svgscale}}%
    \fi%
  \else%
    \setlength{\unitlength}{\svgwidth}%
  \fi%
  \global\let\svgwidth\undefined%
  \global\let\svgscale\undefined%
  \makeatother%
  \begin{picture}(1,0.42182373)%
    \lineheight{1}%
    \setlength\tabcolsep{0pt}%
    \put(0,0){\includegraphics[width=\unitlength,page=1]{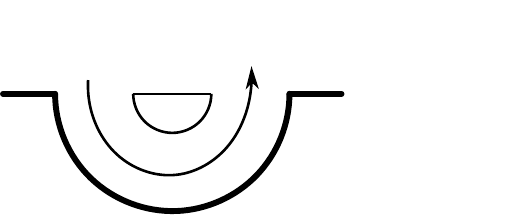}}%
    \put(0.11324335,0.28536387){\color[rgb]{0,0,0}\makebox(0,0)[lt]{\lineheight{1.25}\smash{\begin{tabular}[t]{l}$I$\end{tabular}}}}%
    \put(0.5243076,0.010939){\color[rgb]{0,0,0}\makebox(0,0)[lt]{\lineheight{1.25}\smash{\begin{tabular}[t]{l}$\overset{I}{X}$\end{tabular}}}}%
    \put(0,0){\includegraphics[width=\unitlength,page=2]{trace_quantique.pdf}}%
  \end{picture}%
\endgroup%

\end{center}
is evaluated as 
\[ \overset{I}{X}{^i_j} \, d'_I(v_i \otimes v^j) = \overset{I}{X}{^i_j} \, v^j(gv_i) = \overset{I}{X}{^i_j} \, \overset{I}{g}^j_i = \mathrm{tr}\bigl(\overset{I}{g} \overset{I}{X}\bigr) = \mathrm{tr}_q\bigl(\overset{I}{X}\bigr), \]
and we see that it represents the quantum trace of $\overset{I}{X}$.

\indent Recall that the matrices $\overset{I}{A}(i), \overset{I}{B}(j), \overset{I}{M}(k)$ commute with the morphisms (naturality, see \eqref{naturaliteLgn}). Hence, this is also true for matrices of the form \eqref{matriceGenerale}. Namely we have $f \overset{I}{X} = \overset{J}{X} f$, where $f : I \to J$ is a morphism and we identify $f$ with its matrix. Let us see the diagrammatic description of this fact. We have
\begin{align*}
\overset{I}{X}{^i_j} \, f \otimes \mathrm{id}_{I^*}(v_i \otimes v^j)  &= \overset{I}{X}{^i_j} \, f_i^kv_k \otimes v^j = \bigl(f\overset{I}{X}\bigr)^k_j \, v_k \otimes v^j = \bigl(\overset{J}{X}f\bigr)^k_j \, w_k \otimes w^j\\
& = \overset{J}{X}{^k_l} \, w_k \otimes f^l_jw^j = \overset{J}{X}{^k_l} \, \mathrm{id}_I \otimes f^*(w_k \otimes w^l)
\end{align*}
where $f^* : J^* \to I^*$ is the transpose of $f$. Thus we get the first diagram below. The second diagram is a consequence of the first thanks to \eqref{sensOppose} and the equality $f \circ e_I = e_J \circ f^{**}$.
\begin{equation}\label{naturalite}
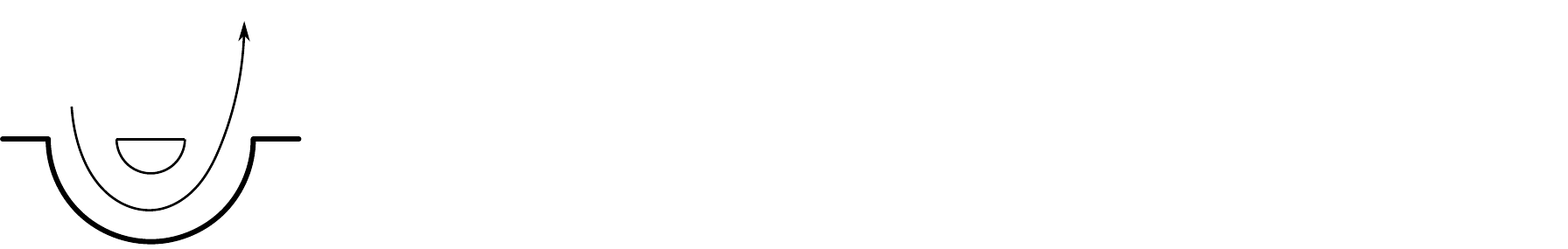
\end{equation}

\indent Let us now write the defining relations of $\mathcal{L}_{g,n}(H)$ in a diagrammatic form. Note that for $a,b \in H$, it holds
\begin{equation}\label{trucEvidentdAlgebreLineaire}
\bigl(\overset{I}{a} \overset{I}{X} \overset{I}{b}\bigr)^i_j \, v_i \otimes v^j = \overset{I}{X}{^i_j} \, av_i \otimes S^{-1}(b)v^j.
\end{equation}
We will use this fact several times in the sequel.

\smallskip

\noindent \textbullet ~~ {\em Fusion relation of  $\mathcal{L}_{0,1}(H)$: } Recall that 
\[ \overset{I \otimes J}{M} = \overset{I}{M}_1 (\overset{IJ}{R'})_{12} \overset{J}{M}_2 (\overset{IJ}{R'}){^{-1}_{12}} = \left(\overset{I}{M} \, \overset{I}{b_i} \, \overset{I}{b_j}\right)_1 \left(\overset{J}{a_i} \, \overset{J}{M} \, \overset{J}{S(a_j)}\right)_2. \]
Hence, we have:
\begin{align*}
\overset{I \otimes J}{M}{^{km}_{ln}} \, v_k \otimes w_m \otimes w^n \otimes v^l &= \left(\overset{I}{M} \, \overset{I}{b_i} \, \overset{I}{b_j}\right)^k_l \left(\overset{J}{a_i} \, \overset{J}{M} \, \overset{J}{S(a_j)}\right)^m_n \, v_k \otimes w_m \otimes w^n \otimes v^l \\
&= \overset{I}{M}{^k_l} \overset{J}{M}{^m_n} \, v_k \otimes a_iw_m \otimes a_j w^n \otimes S^{-1}(b_i b_j) v^l\\
&= \overset{I}{M}{^k_l} \overset{J}{M}{^m_n} \, v_k \otimes S(a_i) w_m \otimes S(a_j) w^n \otimes b_j b_i v^l\\
&= \overset{I}{M}{^k_l} \overset{J}{M}{^m_n} \, \mathrm{id}_I \otimes \mathrm{id}_J \otimes c_{J^*, I^*}^{-1}\!\left(v_k \otimes S(a_i) w_m \otimes b_i v^l \otimes w^n \right)\\
&= \overset{I}{M}{^k_l} \overset{J}{M}{^m_n} \, \left(\mathrm{id}_I \otimes \mathrm{id}_J \otimes c_{J^*, I^*}^{-1} \right) \circ \left( \mathrm{id}_I \otimes c_{J,I^*}^{-1} \otimes \mathrm{id}_{J^*}\right)\!\left(v_k \otimes v^l \otimes w_m \otimes w^n \right).\\
\end{align*}
We thus obtain the diagrammatic identity below:
\begin{equation}\label{relationFusion}
\begingroup%
  \makeatletter%
  \providecommand\color[2][]{%
    \errmessage{(Inkscape) Color is used for the text in Inkscape, but the package 'color.sty' is not loaded}%
    \renewcommand\color[2][]{}%
  }%
  \providecommand\transparent[1]{%
    \errmessage{(Inkscape) Transparency is used (non-zero) for the text in Inkscape, but the package 'transparent.sty' is not loaded}%
    \renewcommand\transparent[1]{}%
  }%
  \providecommand\rotatebox[2]{#2}%
  \newcommand*\fsize{\dimexpr\f@size pt\relax}%
  \newcommand*\lineheight[1]{\fontsize{\fsize}{#1\fsize}\selectfont}%
  \ifx\svgwidth\undefined%
    \setlength{\unitlength}{505.84185204bp}%
    \ifx\svgscale\undefined%
      \relax%
    \else%
      \setlength{\unitlength}{\unitlength * \real{\svgscale}}%
    \fi%
  \else%
    \setlength{\unitlength}{\svgwidth}%
  \fi%
  \global\let\svgwidth\undefined%
  \global\let\svgscale\undefined%
  \makeatother%
  \begin{picture}(1,0.24480261)%
    \lineheight{1}%
    \setlength\tabcolsep{0pt}%
    \put(0,0){\includegraphics[width=\unitlength,page=1]{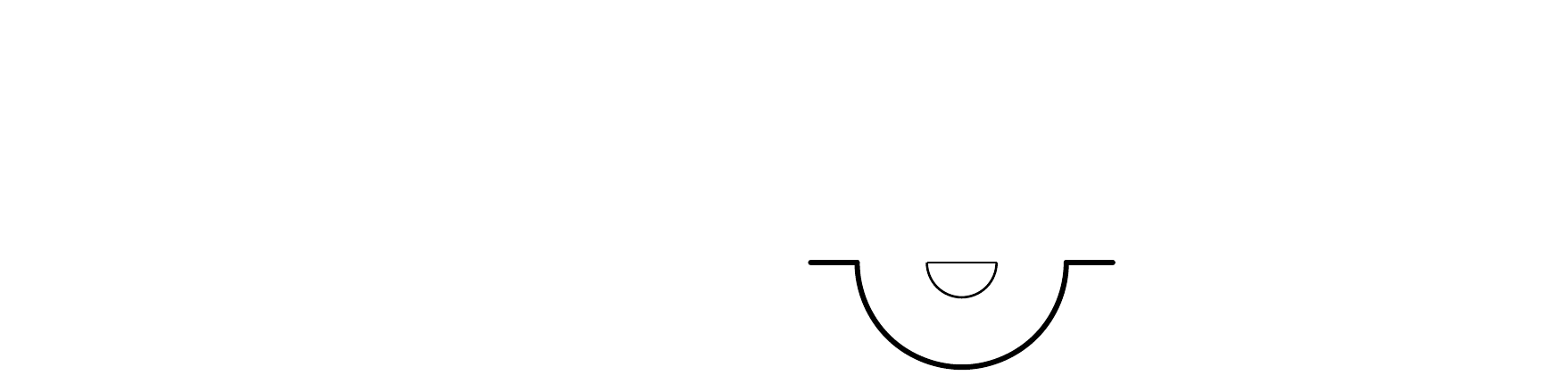}}%
    \put(0.52910019,0.23244731){\color[rgb]{0,0,0}\makebox(0,0)[lt]{\lineheight{1.25}\smash{\begin{tabular}[t]{l}$I$\end{tabular}}}}%
    \put(0.6625636,0.00333888){\color[rgb]{0,0,0}\makebox(0,0)[lt]{\lineheight{1.25}\smash{\begin{tabular}[t]{l}$\overset{I}{M}$\end{tabular}}}}%
    \put(0,0){\includegraphics[width=\unitlength,page=2]{relation_fusion.pdf}}%
    \put(0.8573616,0.00417054){\color[rgb]{0,0,0}\makebox(0,0)[lt]{\lineheight{1.25}\smash{\begin{tabular}[t]{l}$\overset{J}{M}$\end{tabular}}}}%
    \put(0,0){\includegraphics[width=\unitlength,page=3]{relation_fusion.pdf}}%
    \put(-0.00146332,0.11381494){\color[rgb]{0,0,0}\makebox(0,0)[lt]{\lineheight{1.25}\smash{\begin{tabular}[t]{l}$I \otimes J$\end{tabular}}}}%
    \put(0.17243856,0.0039469){\color[rgb]{0,0,0}\makebox(0,0)[lt]{\lineheight{1.25}\smash{\begin{tabular}[t]{l}$\overset{I \otimes J}{M}$\end{tabular}}}}%
    \put(0,0){\includegraphics[width=\unitlength,page=4]{relation_fusion.pdf}}%
    \put(0.28715026,0.11392188){\color[rgb]{0,0,0}\makebox(0,0)[lt]{\lineheight{1.25}\smash{\begin{tabular}[t]{l}$I$\end{tabular}}}}%
    \put(0.42844705,0.00394694){\color[rgb]{0,0,0}\makebox(0,0)[lt]{\lineheight{1.25}\smash{\begin{tabular}[t]{l}$\overset{I \otimes J}{M}$\end{tabular}}}}%
    \put(0,0){\includegraphics[width=\unitlength,page=5]{relation_fusion.pdf}}%
    \put(0.31061585,0.11509596){\color[rgb]{0,0,0}\makebox(0,0)[lt]{\lineheight{1.25}\smash{\begin{tabular}[t]{l}$J$\end{tabular}}}}%
    \put(0,0){\includegraphics[width=\unitlength,page=6]{relation_fusion.pdf}}%
    \put(0.63347318,0.23400778){\color[rgb]{0,0,0}\makebox(0,0)[lt]{\lineheight{1.25}\smash{\begin{tabular}[t]{l}$J$\end{tabular}}}}%
    \put(0.2334006,0.07086557){\color[rgb]{0,0,0}\makebox(0,0)[lt]{\lineheight{1.25}\smash{\begin{tabular}[t]{l}=\end{tabular}}}}%
    \put(0.48338631,0.07148884){\color[rgb]{0,0,0}\makebox(0,0)[lt]{\lineheight{1.25}\smash{\begin{tabular}[t]{l}=\end{tabular}}}}%
  \end{picture}%
\endgroup%

\end{equation}

\begin{figure}[!h]
\centering
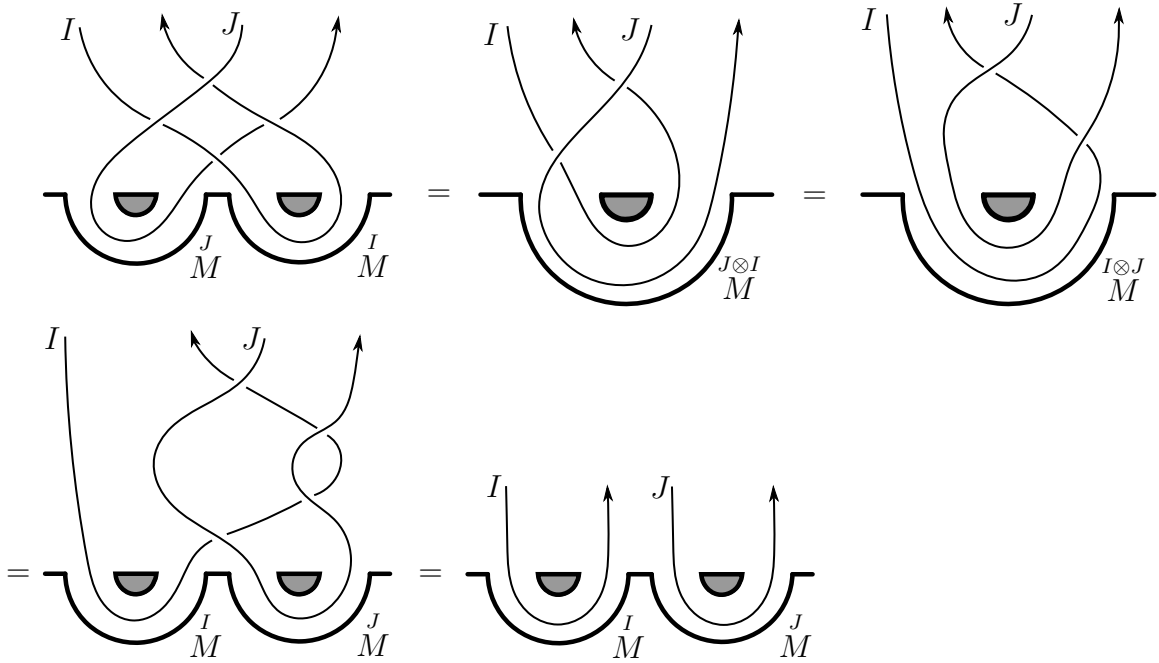
\caption{Proof of the reflection equation in $\mathcal{L}_{0,1}(H)$.}
\label{preuveReflexion}
\end{figure}

\noindent \textbullet ~~ {\em Reflection equation:} The reflection equation in $\mathcal{L}_{0,1}(H)$ is the following exchange relation :
\[ \overset{IJ}{R}_{12} \overset{I}{M}_1 (\overset{IJ}{R'})_{12} \overset{J}{M}_2 = \overset{J}{M}_2 \overset{IJ}{R}_{12} \overset{I}{M}_1 (\overset{IJ}{R'})_{12}. \]
The graphical representation of this equation is depicted as follows:
\begin{equation}\label{equationReflexion}
\begingroup%
  \makeatletter%
  \providecommand\color[2][]{%
    \errmessage{(Inkscape) Color is used for the text in Inkscape, but the package 'color.sty' is not loaded}%
    \renewcommand\color[2][]{}%
  }%
  \providecommand\transparent[1]{%
    \errmessage{(Inkscape) Transparency is used (non-zero) for the text in Inkscape, but the package 'transparent.sty' is not loaded}%
    \renewcommand\transparent[1]{}%
  }%
  \providecommand\rotatebox[2]{#2}%
  \newcommand*\fsize{\dimexpr\f@size pt\relax}%
  \newcommand*\lineheight[1]{\fontsize{\fsize}{#1\fsize}\selectfont}%
  \ifx\svgwidth\undefined%
    \setlength{\unitlength}{474.76305933bp}%
    \ifx\svgscale\undefined%
      \relax%
    \else%
      \setlength{\unitlength}{\unitlength * \real{\svgscale}}%
    \fi%
  \else%
    \setlength{\unitlength}{\svgwidth}%
  \fi%
  \global\let\svgwidth\undefined%
  \global\let\svgscale\undefined%
  \makeatother%
  \begin{picture}(1,0.2700562)%
    \lineheight{1}%
    \setlength\tabcolsep{0pt}%
    \put(0,0){\includegraphics[width=\unitlength,page=1]{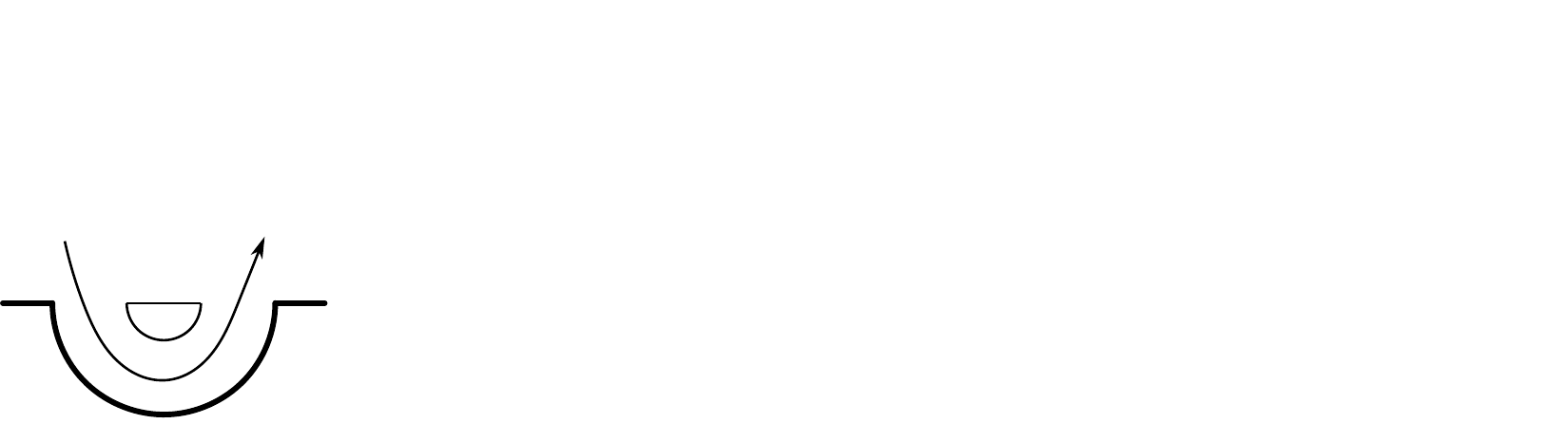}}%
    \put(0.0260731,0.11948785){\color[rgb]{0,0,0}\makebox(0,0)[lt]{\lineheight{1.25}\smash{\begin{tabular}[t]{l}$I$\end{tabular}}}}%
    \put(0.16231857,0.00490331){\color[rgb]{0,0,0}\makebox(0,0)[lt]{\lineheight{1.25}\smash{\begin{tabular}[t]{l}$\overset{I}{M}$\end{tabular}}}}%
    \put(0,0){\includegraphics[width=\unitlength,page=2]{equation_reflexion.pdf}}%
    \put(0.36956243,0.00490331){\color[rgb]{0,0,0}\makebox(0,0)[lt]{\lineheight{1.25}\smash{\begin{tabular}[t]{l}$\overset{J}{M}$\end{tabular}}}}%
    \put(0,0){\includegraphics[width=\unitlength,page=3]{equation_reflexion.pdf}}%
    \put(0.23172386,0.11910221){\color[rgb]{0,0,0}\makebox(0,0)[lt]{\lineheight{1.25}\smash{\begin{tabular}[t]{l}$J$\end{tabular}}}}%
    \put(0.43559719,0.07207022){\color[rgb]{0,0,0}\makebox(0,0)[lt]{\lineheight{1.25}\smash{\begin{tabular}[t]{l}=\end{tabular}}}}%
    \put(0.63862125,0.00420527){\color[rgb]{0,0,0}\makebox(0,0)[lt]{\lineheight{1.25}\smash{\begin{tabular}[t]{l}$\overset{J}{M}$\end{tabular}}}}%
    \put(0.84802425,0.0044642){\color[rgb]{0,0,0}\makebox(0,0)[lt]{\lineheight{1.25}\smash{\begin{tabular}[t]{l}$\overset{I}{M}$\end{tabular}}}}%
    \put(0,0){\includegraphics[width=\unitlength,page=4]{equation_reflexion.pdf}}%
    \put(0.51327501,0.25694817){\color[rgb]{0,0,0}\makebox(0,0)[lt]{\lineheight{1.25}\smash{\begin{tabular}[t]{l}$I$\end{tabular}}}}%
    \put(0,0){\includegraphics[width=\unitlength,page=5]{equation_reflexion.pdf}}%
    \put(0.71826158,0.25855473){\color[rgb]{0,0,0}\makebox(0,0)[lt]{\lineheight{1.25}\smash{\begin{tabular}[t]{l}$J$\end{tabular}}}}%
    \put(0,0){\includegraphics[width=\unitlength,page=6]{equation_reflexion.pdf}}%
  \end{picture}%
\endgroup%

\end{equation}
\noindent A diagrammatic proof of this relation is shown in Figure \ref{preuveReflexion}; this is simply a graphical reformulation of the proof of Proposition \ref{EqRef}. For the second equality, we used naturality \eqref{naturalite} and the fact that $c_{J,I}^* = c_{J^*, I^*}$.

If we plug the inverse of the tangle (which is a braid) and we exchange $I$ and $J$ in the reflection equation above, we get the following relation, which we also call reflection equation:
\begin{equation}\label{equationReflexionBis}
\begingroup%
  \makeatletter%
  \providecommand\color[2][]{%
    \errmessage{(Inkscape) Color is used for the text in Inkscape, but the package 'color.sty' is not loaded}%
    \renewcommand\color[2][]{}%
  }%
  \providecommand\transparent[1]{%
    \errmessage{(Inkscape) Transparency is used (non-zero) for the text in Inkscape, but the package 'transparent.sty' is not loaded}%
    \renewcommand\transparent[1]{}%
  }%
  \providecommand\rotatebox[2]{#2}%
  \newcommand*\fsize{\dimexpr\f@size pt\relax}%
  \newcommand*\lineheight[1]{\fontsize{\fsize}{#1\fsize}\selectfont}%
  \ifx\svgwidth\undefined%
    \setlength{\unitlength}{474.76305703bp}%
    \ifx\svgscale\undefined%
      \relax%
    \else%
      \setlength{\unitlength}{\unitlength * \real{\svgscale}}%
    \fi%
  \else%
    \setlength{\unitlength}{\svgwidth}%
  \fi%
  \global\let\svgwidth\undefined%
  \global\let\svgscale\undefined%
  \makeatother%
  \begin{picture}(1,0.27005623)%
    \lineheight{1}%
    \setlength\tabcolsep{0pt}%
    \put(0,0){\includegraphics[width=\unitlength,page=1]{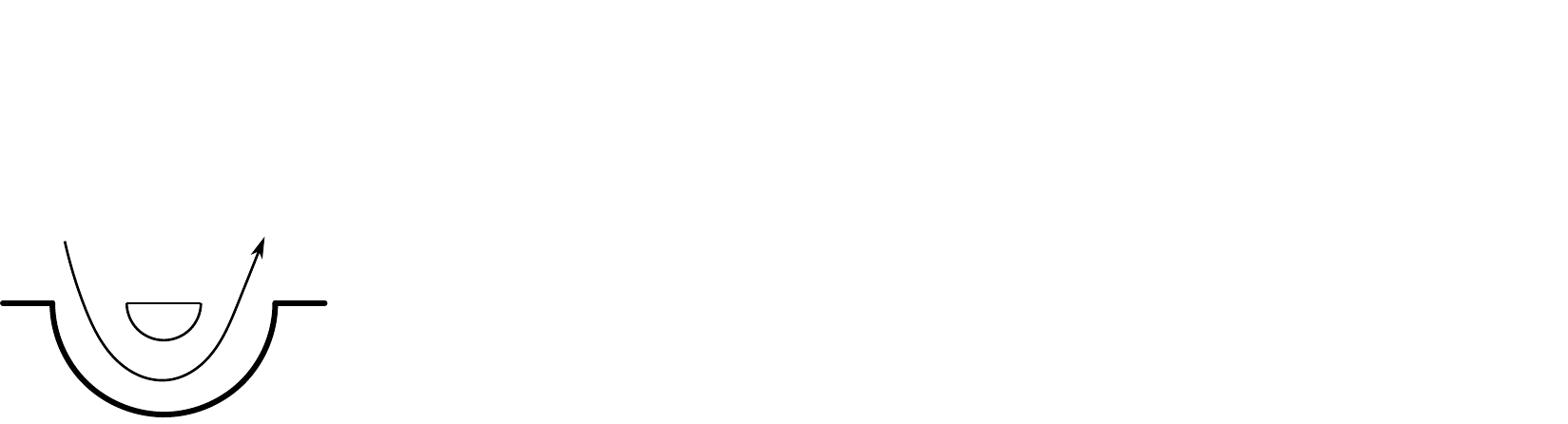}}%
    \put(0.0260731,0.11948787){\color[rgb]{0,0,0}\makebox(0,0)[lt]{\lineheight{1.25}\smash{\begin{tabular}[t]{l}$I$\end{tabular}}}}%
    \put(0.16231857,0.00490326){\color[rgb]{0,0,0}\makebox(0,0)[lt]{\lineheight{1.25}\smash{\begin{tabular}[t]{l}$\overset{I}{M}$\end{tabular}}}}%
    \put(0,0){\includegraphics[width=\unitlength,page=2]{equation_reflexion_bis.pdf}}%
    \put(0.36956247,0.00490326){\color[rgb]{0,0,0}\makebox(0,0)[lt]{\lineheight{1.25}\smash{\begin{tabular}[t]{l}$\overset{J}{M}$\end{tabular}}}}%
    \put(0,0){\includegraphics[width=\unitlength,page=3]{equation_reflexion_bis.pdf}}%
    \put(0.23172386,0.11910223){\color[rgb]{0,0,0}\makebox(0,0)[lt]{\lineheight{1.25}\smash{\begin{tabular}[t]{l}$J$\end{tabular}}}}%
    \put(0.43559714,0.07207019){\color[rgb]{0,0,0}\makebox(0,0)[lt]{\lineheight{1.25}\smash{\begin{tabular}[t]{l}=\end{tabular}}}}%
    \put(0.63862125,0.00420527){\color[rgb]{0,0,0}\makebox(0,0)[lt]{\lineheight{1.25}\smash{\begin{tabular}[t]{l}$\overset{J}{M}$\end{tabular}}}}%
    \put(0.84802425,0.00446422){\color[rgb]{0,0,0}\makebox(0,0)[lt]{\lineheight{1.25}\smash{\begin{tabular}[t]{l}$\overset{I}{M}$\end{tabular}}}}%
    \put(0,0){\includegraphics[width=\unitlength,page=4]{equation_reflexion_bis.pdf}}%
    \put(0.51327496,0.25694821){\color[rgb]{0,0,0}\makebox(0,0)[lt]{\lineheight{1.25}\smash{\begin{tabular}[t]{l}$I$\end{tabular}}}}%
    \put(0,0){\includegraphics[width=\unitlength,page=5]{equation_reflexion_bis.pdf}}%
    \put(0.71826158,0.25855476){\color[rgb]{0,0,0}\makebox(0,0)[lt]{\lineheight{1.25}\smash{\begin{tabular}[t]{l}$J$\end{tabular}}}}%
    \put(0,0){\includegraphics[width=\unitlength,page=6]{equation_reflexion_bis.pdf}}%
  \end{picture}%
\endgroup%

\end{equation}

\noindent \textbullet ~~ {\em Negative orientation and inverse of $\overset{I}{M}$: } Recall the algebra $\mathcal{F}_{0,1}(H) \cong \mathcal{L}_{0,1}(H)$ of Remark \ref{F01}, which is $H^*$ endowed with the product $\varphi \ast \psi = \varphi\!\left( ? b_j S(b_i) \right) \psi\!\left(a_i ? a_j\right)$. We identify $\mathcal{L}_{0,1}(H)$ and $\mathcal{F}_{0,1}(H)$ by $\overset{I}{M} \mapsto \overset{I}{T}$. In other words, we consider $\overset{I}{M}$ as a matrix whose coefficients are linear forms on $H$; the evaluation of $\overset{I}{M}$ on $x \in H$ is obviously defined by $\overset{I}{M}(x){^i_j} = \overset{I}{M}{^i_j}(x) = \overset{I}{x}{^i_j}$.

\begin{lemma}\label{matriceInverse}
Under the above identification, it holds
\[ \overset{I}{M}{^{-1}}(x) = \left( u^{-1} S(b_i) S(x) a_i \right)^I \]
where $u = gv$ is the Drinfeld element \eqref{elementDrinfeld}.
\end{lemma}
\begin{proof}
\begin{align*}
\biggl( \overset{I}{M}{^i_j} \ast (\overset{I}{M}{^{-1}}){^j_k} \biggr)\!(x) &= \left\langle \overset{I}{M}{^i_j}, x' b_m S(b_l) \right\rangle \left\langle (\overset{I}{M}{^{-1}}){^j_k}, a_l x'' a_m \right\rangle\\
&= \biggl( \overset{I}{x'} \overset{I}{b_m} \overset{I}{b_l} \biggr)^i_j \biggl( \overset{I}{u^{-1}} \overset{I}{S(b_n)} \overset{I}{S(a_m)} \overset{I}{S(x'')} \overset{I}{a_l} \overset{I}{a_n} \biggr)^j_k\\
&= \biggl( \overset{I}{x'} \overset{I}{b_m} \overset{I}{u^{-1}} \overset{I}{S^2(b_l S^{-1}(b_n))} \overset{I}{S(a_m)} \overset{I}{S(x'')} \overset{I}{a_l} \overset{I}{a_n}  \biggr)^i_k\\
& = \biggl( \overset{I}{x'} \overset{I}{u^{-1}} \overset{I}{S^2(b_m)} \overset{I}{S(a_m)} \overset{I}{S(x'')}  \biggr)^i_k = \biggl( \overset{I}{x'} \overset{I}{S(x'')}  \biggr)^i_k = \varepsilon(x) \delta^i_k.
\end{align*}
We used that $S^2(b_m)S(a_m) = S(b_m)a_m = u$.
\end{proof}

\noindent Still under the identification $\mathcal{L}_{0,1}(H) = \mathcal{F}_{0,1}(H)$ and by \eqref{antipodeT}, we have $\overset{I^*}{M}(x) = {^t}\overset{I}{S(x)}$. It follows that the formula of Lemma \ref{matriceInverse} can be rewritten as 
$\overset{I}{M}{^{-1}}(x) = \left( u^{-1} S(b_i)\right)^I {^t}\overset{I^*}{M} \overset{I}{a_i}$, and finally:
\[ \overset{I^*}{M} = \exposantGauche{t}{\biggl( \overset{I}{u} \overset{I}{b_i} \overset{I}{M}{^{-1}} \overset{I}{S^2(a_i)} \biggr)}. \]
Let us represent this formula graphically:
\begin{align*}
\overset{I^*}{M}{^i_j} \, \mathrm{id}_{I^*} \otimes e_I\bigl(v^i \otimes \langle ?, v_j \rangle \bigr) &= \overset{I^*}{M}{^i_j} \, v^i \otimes g^{-1}v_j = \biggl(\overset{I}{u} \overset{I}{b_k} \overset{I}{M}{^{-1}} \overset{I}{S^2(a_k)}\biggr)^j_i \, v^i \otimes g^{-1}v_j\\
&= (\overset{I}{M}{^{-1}})^j_i \, S(a_k)v^i \otimes vb_k v_j = (\overset{I}{M}{^{-1}})^j_i \, (\mathrm{id}_{I^*} \otimes \theta_I^{-1}) \circ c_{I^*, I}^{-1}(v_j \otimes v^i).
\end{align*}
By definition of the value of a negatively oriented strand in a handle \eqref{sensOppose}, we get:
\begin{equation}\label{dessinMatriceInverse}
\begingroup%
  \makeatletter%
  \providecommand\color[2][]{%
    \errmessage{(Inkscape) Color is used for the text in Inkscape, but the package 'color.sty' is not loaded}%
    \renewcommand\color[2][]{}%
  }%
  \providecommand\transparent[1]{%
    \errmessage{(Inkscape) Transparency is used (non-zero) for the text in Inkscape, but the package 'transparent.sty' is not loaded}%
    \renewcommand\transparent[1]{}%
  }%
  \providecommand\rotatebox[2]{#2}%
  \newcommand*\fsize{\dimexpr\f@size pt\relax}%
  \newcommand*\lineheight[1]{\fontsize{\fsize}{#1\fsize}\selectfont}%
  \ifx\svgwidth\undefined%
    \setlength{\unitlength}{307.01495918bp}%
    \ifx\svgscale\undefined%
      \relax%
    \else%
      \setlength{\unitlength}{\unitlength * \real{\svgscale}}%
    \fi%
  \else%
    \setlength{\unitlength}{\svgwidth}%
  \fi%
  \global\let\svgwidth\undefined%
  \global\let\svgscale\undefined%
  \makeatother%
  \begin{picture}(1,0.26757672)%
    \lineheight{1}%
    \setlength\tabcolsep{0pt}%
    \put(0,0){\includegraphics[width=\unitlength,page=1]{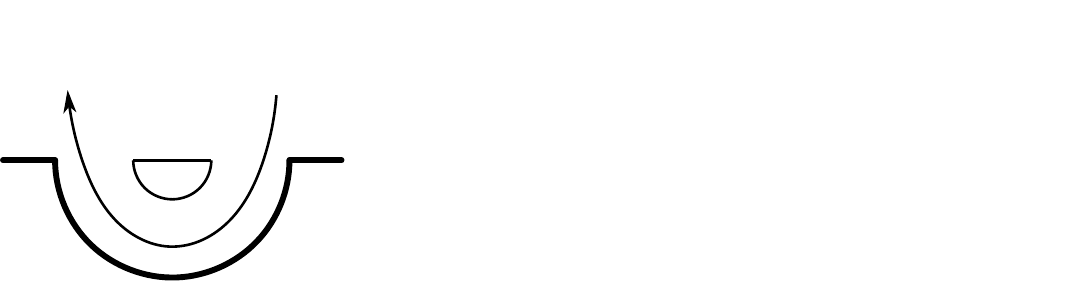}}%
    \put(0.23237783,0.17814955){\color[rgb]{0,0,0}\makebox(0,0)[lt]{\lineheight{1.25}\smash{\begin{tabular}[t]{l}$I$\end{tabular}}}}%
    \put(0.25391134,0.00626829){\color[rgb]{0,0,0}\makebox(0,0)[lt]{\lineheight{1.25}\smash{\begin{tabular}[t]{l}$\overset{I}{M}$\end{tabular}}}}%
    \put(0,0){\includegraphics[width=\unitlength,page=2]{matrice_inverse.pdf}}%
    \put(0.35409346,0.11003823){\color[rgb]{0,0,0}\makebox(0,0)[lt]{\lineheight{1.25}\smash{\begin{tabular}[t]{l}=\end{tabular}}}}%
    \put(0,0){\includegraphics[width=\unitlength,page=3]{matrice_inverse.pdf}}%
    \put(0.60056912,0.24979102){\color[rgb]{0,0,0}\makebox(0,0)[lt]{\lineheight{1.25}\smash{\begin{tabular}[t]{l}$I$\end{tabular}}}}%
    \put(0.65989041,0.00529753){\color[rgb]{0,0,0}\makebox(0,0)[lt]{\lineheight{1.25}\smash{\begin{tabular}[t]{l}$\overset{I}{M}{^{-1}}$\end{tabular}}}}%
    \put(0,0){\includegraphics[width=\unitlength,page=4]{matrice_inverse.pdf}}%
  \end{picture}%
\endgroup%

\end{equation}
Note that since this formula is true in $\mathcal{L}_{0,1}(H)$, it will be true for any matrix $\overset{I}{X}$ with coefficients in $\mathcal{L}_{g,n}(H)$ which satisfies the fusion relation. We also mention that \eqref{dessinMatriceInverse} implies the following relation, which will be used later:
\begin{equation}\label{traceInverseTraceDual}
\mathrm{tr}_q\bigl(\overset{I}{M}{^{-1}}\bigr) = \mathrm{tr}_q\bigl(\overset{I^*}{M}\bigr).
\end{equation}
Indeed, thanks to \eqref{sensOppose} and \eqref{dualitePrime}, we have:
\begin{center}
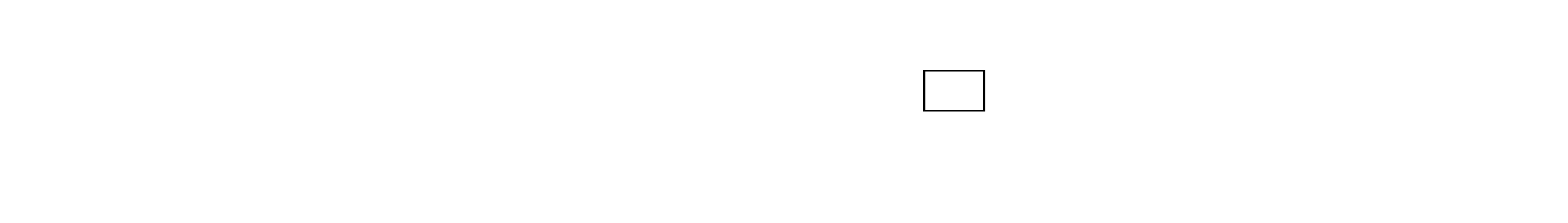
\end{center}

\noindent \textbullet ~~ {\em Exchange relation of $\mathcal{L}_{1,0}(H)$: } Recall from \eqref{echangeL10Inverse} that this relation can be written
\[ \overset{I}{A}_1 \overset{J}{B}_2 = \overset{J}{(a_i)}_2\overset{IJ}{(R')}_{12}\overset{J}{B}_2\overset{IJ}{R}_{12}\overset{I}{A}_1\overset{IJ}{(R')}_{12}\overset{I}{S(b_i)}_1 = \biggl( \overset{J}{a_i} \, \overset{J}{a_j} \, \overset{J}{B} \, \overset{J}{b_k} \, \overset{J}{a_l} \biggr)_{\!\!2} \, \biggl( \overset{I}{b_j} \, \overset{I}{a_k} \, \overset{I}{A} \, \overset{I}{b_l} \, \overset{I}{S(b_i)} \biggr)_{\!\!1}. \]
Hence we have:
\begin{align*}
& \: \overset{I}{A}{^m_n} \overset{J}{B}{^o_p} \, v_m \otimes v^n \otimes w_o \otimes w^p \\
=& \: \biggl( \overset{J}{a_i} \, \overset{J}{a_j} \, \overset{J}{B} \, \overset{J}{b_k} \, \overset{J}{a_l} \biggr)^o_p \, \biggl( \overset{I}{b_j} \, \overset{I}{a_k} \, \overset{I}{A} \, \overset{I}{b_l} \, \overset{I}{S(b_i)} \biggr)^m_n \, v_m \otimes v^n \otimes w_o \otimes w^p\\
=& \: \overset{J}{B}{^o_p} \overset{I}{A}{^m_n} \, b_j a_k v_m \otimes S^{-1}\!\left(b_lS(b_i)\right)v^n \otimes a_ia_j w_o \otimes S^{-1}(b_k a_l) w^p\\
=& \: \overset{J}{B}{^o_p} \overset{I}{A}{^m_n} \, b_jS(a_k)v_m \otimes b_i b_lv^n \otimes a_ia_j w_o \otimes a_l b_k w^p\\
=& \: \overset{J}{B}{^o_p} \overset{I}{A}{^m_n} \, \mathrm{id}_I \otimes c_{J,I^*} \otimes \mathrm{id}_{J^*}\bigl(b_jS(a_k)v_m \otimes a_j w_o \otimes b_lv^n \otimes a_l b_k w^p\bigr)\\
=& \: \overset{J}{B}{^o_p} \overset{I}{A}{^m_n} \, (\mathrm{id}_I \otimes c_{J,I^*} \otimes \mathrm{id}_{J^*}) \circ (c_{J,I} \otimes c_{J^*, I^*})\bigl(w_o \otimes S(a_k)v_m \otimes b_k w^p \otimes v^n\bigr)\\
=& \: \overset{J}{B}{^o_p} \overset{I}{A}{^m_n} \, (\mathrm{id}_I \otimes c_{J,I^*} \otimes \mathrm{id}_{J^*}) \circ (c_{J,I} \otimes c_{J^*, I^*}) \circ (\mathrm{id}_J \otimes c_{I,J^*}^{-1} \otimes \mathrm{id}_{I^*})\bigl(w_o \otimes w^p \otimes v_m \otimes v^n\bigr)
\end{align*}
and this yields the diagrammatic identity below.
\begin{equation}\label{dessinEchangeL10}
\begingroup%
  \makeatletter%
  \providecommand\color[2][]{%
    \errmessage{(Inkscape) Color is used for the text in Inkscape, but the package 'color.sty' is not loaded}%
    \renewcommand\color[2][]{}%
  }%
  \providecommand\transparent[1]{%
    \errmessage{(Inkscape) Transparency is used (non-zero) for the text in Inkscape, but the package 'transparent.sty' is not loaded}%
    \renewcommand\transparent[1]{}%
  }%
  \providecommand\rotatebox[2]{#2}%
  \newcommand*\fsize{\dimexpr\f@size pt\relax}%
  \newcommand*\lineheight[1]{\fontsize{\fsize}{#1\fsize}\selectfont}%
  \ifx\svgwidth\undefined%
    \setlength{\unitlength}{468.26884423bp}%
    \ifx\svgscale\undefined%
      \relax%
    \else%
      \setlength{\unitlength}{\unitlength * \real{\svgscale}}%
    \fi%
  \else%
    \setlength{\unitlength}{\svgwidth}%
  \fi%
  \global\let\svgwidth\undefined%
  \global\let\svgscale\undefined%
  \makeatother%
  \begin{picture}(1,0.27309375)%
    \lineheight{1}%
    \setlength\tabcolsep{0pt}%
    \put(0,0){\includegraphics[width=\unitlength,page=1]{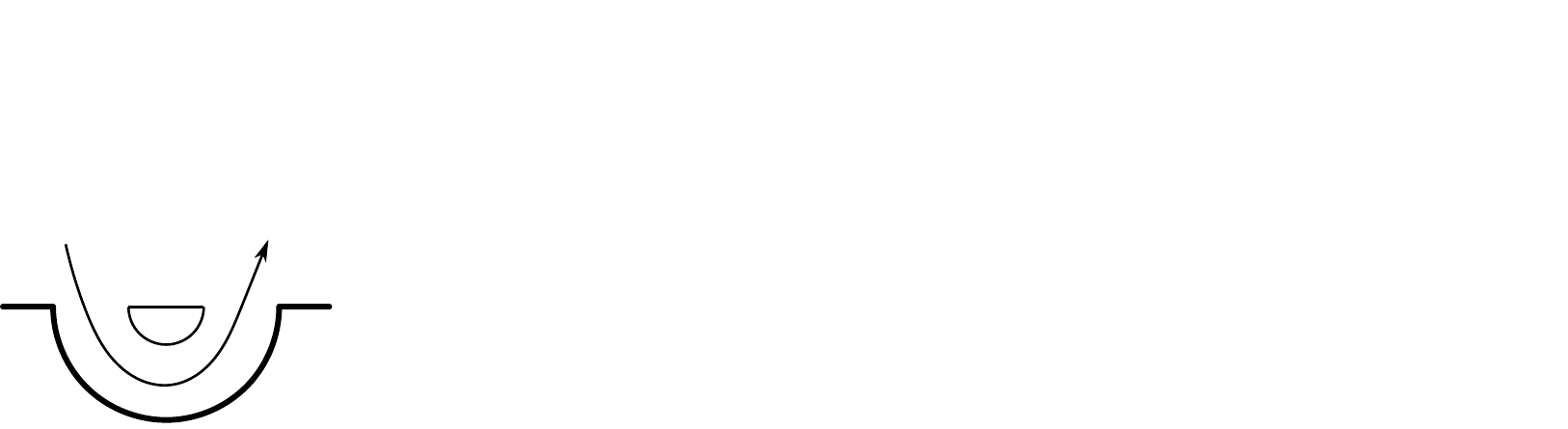}}%
    \put(0.0264347,0.12068439){\color[rgb]{0,0,0}\makebox(0,0)[lt]{\lineheight{1.25}\smash{\begin{tabular}[t]{l}$I$\end{tabular}}}}%
    \put(0.1645697,0.00451084){\color[rgb]{0,0,0}\makebox(0,0)[lt]{\lineheight{1.25}\smash{\begin{tabular}[t]{l}$\overset{I}{A}$\end{tabular}}}}%
    \put(0,0){\includegraphics[width=\unitlength,page=2]{echange_L10.pdf}}%
    \put(0.37468773,0.00451084){\color[rgb]{0,0,0}\makebox(0,0)[lt]{\lineheight{1.25}\smash{\begin{tabular}[t]{l}$\overset{J}{B}$\end{tabular}}}}%
    \put(0,0){\includegraphics[width=\unitlength,page=3]{echange_L10.pdf}}%
    \put(0.23493754,0.1202934){\color[rgb]{0,0,0}\makebox(0,0)[lt]{\lineheight{1.25}\smash{\begin{tabular}[t]{l}$J$\end{tabular}}}}%
    \put(0,0){\includegraphics[width=\unitlength,page=4]{echange_L10.pdf}}%
    \put(0.5142077,0.25980394){\color[rgb]{0,0,0}\makebox(0,0)[lt]{\lineheight{1.25}\smash{\begin{tabular}[t]{l}$I$\end{tabular}}}}%
    \put(0.63887723,0.00426359){\color[rgb]{0,0,0}\makebox(0,0)[lt]{\lineheight{1.25}\smash{\begin{tabular}[t]{l}$\overset{J}{B}$\end{tabular}}}}%
    \put(0,0){\includegraphics[width=\unitlength,page=5]{echange_L10.pdf}}%
    \put(0.84899539,0.00426359){\color[rgb]{0,0,0}\makebox(0,0)[lt]{\lineheight{1.25}\smash{\begin{tabular}[t]{l}$\overset{I}{A}$\end{tabular}}}}%
    \put(0,0){\includegraphics[width=\unitlength,page=6]{echange_L10.pdf}}%
    \put(0.72203715,0.26143277){\color[rgb]{0,0,0}\makebox(0,0)[lt]{\lineheight{1.25}\smash{\begin{tabular}[t]{l}$J$\end{tabular}}}}%
    \put(0,0){\includegraphics[width=\unitlength,page=7]{echange_L10.pdf}}%
    \put(0.43793607,0.07264657){\color[rgb]{0,0,0}\makebox(0,0)[lt]{\lineheight{1.25}\smash{\begin{tabular}[t]{l}=\end{tabular}}}}%
  \end{picture}%
\endgroup%

\end{equation}
Remark that, as above, the tangle appearing in \eqref{dessinEchangeL10} is in fact a braid and it can be inverted in order to exchange $\overset{J}{B}$ and $\overset{I}{A}$.

\noindent \textbullet ~~ {\em Exchange relation of $\mathcal{L}_{g,n}(H)$: } Recall from \eqref{eqEchangeLgn} that for $\alpha < \beta$
\[\overset{J}{V}(\beta)_2 \, \overset{I}{U}(\alpha)_1 = (\overset{I}{a_i})_1 \overset{IJ}{R}_{12} \overset{I}{U}(\alpha)_1 \overset{IJ}{R}{^{-1}_{12}}\overset{J}{V}(\beta)_2 \overset{IJ}{R}_{12} \overset{J}{S(b_i)}_2 = \biggl( \overset{I}{a_i} \overset{I}{a_j} \overset{I}{U}(\alpha) \overset{I}{S(a_k)} \overset{I}{a_l} \biggr)_{\!\! 1} \, \biggl( \overset{J}{b_j} \, \overset{J}{b_k} \, \overset{J}{V}(\beta) \, \overset{J}{b_l} \, \overset{J}{S(b_i)} \biggr)_{\!\!2}. \]
where $U, V$ are $A$ or $B$.
\noindent Hence we have:
\begin{align*}
&\: \overset{J}{V}(\beta)^m_n \overset{I}{U}(\alpha)^o_p \, v_m \otimes v^n \otimes w_o \otimes w^p\\
=& \: \biggl( \overset{I}{a_i} \overset{I}{a_j} \overset{I}{U}(\alpha) \overset{I}{S(a_k)} \overset{I}{a_l} \biggr)^o_p \, \biggl( \overset{J}{b_j} \, \overset{J}{b_k} \, \overset{J}{V}(\beta) \, \overset{J}{b_l} \, \overset{J}{S(b_i)} \biggr)^m_n \, v_m \otimes v^n \otimes w_o \otimes w^p\\
=& \: \overset{I}{U}(\alpha)^o_p \overset{J}{V}(\beta)^m_n\, b_j b_k v_m \otimes S^{-1}(b_l S(b_i))v^n \otimes a_i a_j w_o \otimes S^{-1}(S(a_k)a_l) w^p\\
=& \: \overset{J}{V}(\beta)^o_p \overset{I}{U}(\alpha)^m_n\, b_j b_k v_m \otimes b_i b_l v^n \otimes a_i a_j w_o \otimes a_l a_k w^p\\
=& \: \overset{J}{V}(\beta)^o_p \overset{I}{U}(\alpha)^m_n\, \mathrm{id}_I \otimes c_{J,I^*} \otimes \mathrm{id}_{J^*} \bigl(b_j b_k v_m \otimes a_j w_o \otimes b_l v^n \otimes a_l a_k w^p\bigr)\\
=& \: \overset{J}{V}(\beta)^o_p \overset{I}{U}(\alpha)^m_n\, (\mathrm{id}_I \otimes c_{J,I^*} \otimes \mathrm{id}_{J^*}) \circ (c_{J,I} \otimes c_{J^*,I^*}) \bigl(w_o \otimes b_k v_m \otimes a_k w^p \otimes v^n\bigr)\\
=& \: \overset{J}{V}(\beta)^o_p \overset{I}{U}(\alpha)^m_n\, (\mathrm{id}_I \otimes c_{J,I^*} \otimes \mathrm{id}_{J^*}) \circ (c_{J,I} \otimes c_{J^*,I^*}) \circ (\mathrm{id}_J \otimes c_{J^*,I} \otimes \mathrm{id}_{I^*}) \bigl(w_o \otimes w^p \otimes v_m \otimes v^n\bigr)
\end{align*}
and this yields the diagrammatic identity below:

\begin{equation}\label{dessinEchangeLgn}
\begingroup%
  \makeatletter%
  \providecommand\color[2][]{%
    \errmessage{(Inkscape) Color is used for the text in Inkscape, but the package 'color.sty' is not loaded}%
    \renewcommand\color[2][]{}%
  }%
  \providecommand\transparent[1]{%
    \errmessage{(Inkscape) Transparency is used (non-zero) for the text in Inkscape, but the package 'transparent.sty' is not loaded}%
    \renewcommand\transparent[1]{}%
  }%
  \providecommand\rotatebox[2]{#2}%
  \newcommand*\fsize{\dimexpr\f@size pt\relax}%
  \newcommand*\lineheight[1]{\fontsize{\fsize}{#1\fsize}\selectfont}%
  \ifx\svgwidth\undefined%
    \setlength{\unitlength}{495.34809045bp}%
    \ifx\svgscale\undefined%
      \relax%
    \else%
      \setlength{\unitlength}{\unitlength * \real{\svgscale}}%
    \fi%
  \else%
    \setlength{\unitlength}{\svgwidth}%
  \fi%
  \global\let\svgwidth\undefined%
  \global\let\svgscale\undefined%
  \makeatother%
  \begin{picture}(1,0.2583218)%
    \lineheight{1}%
    \setlength\tabcolsep{0pt}%
    \put(0,0){\includegraphics[width=\unitlength,page=1]{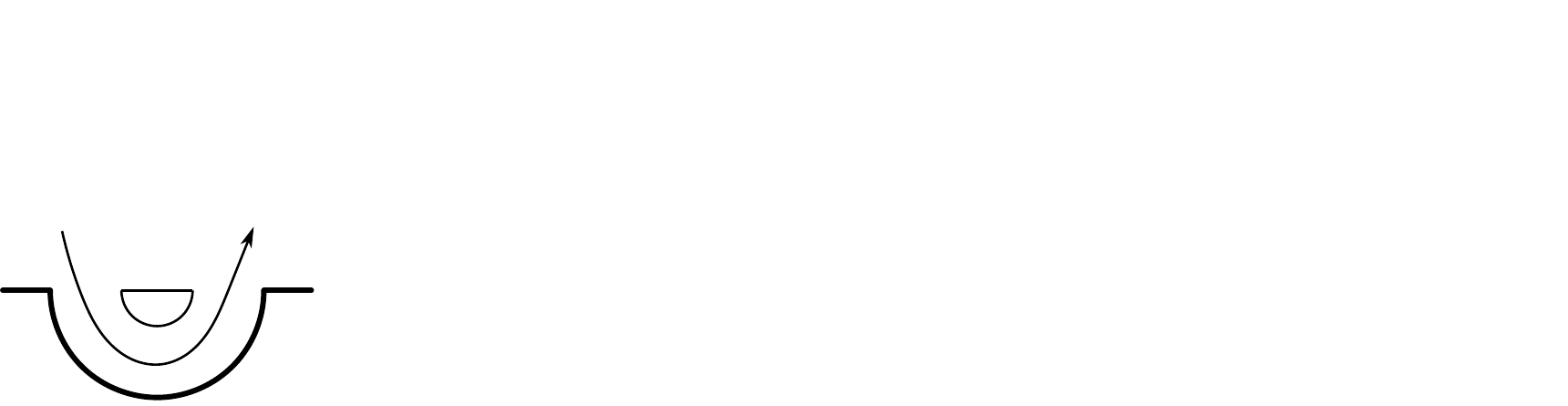}}%
    \put(0.02498959,0.11424421){\color[rgb]{0,0,0}\makebox(0,0)[lt]{\lineheight{1.25}\smash{\begin{tabular}[t]{l}$I$\end{tabular}}}}%
    \put(0.15557314,0.00442154){\color[rgb]{0,0,0}\makebox(0,0)[lt]{\lineheight{1.25}\smash{\begin{tabular}[t]{l}$\overset{J}{V}(\beta)$\end{tabular}}}}%
    \put(0,0){\includegraphics[width=\unitlength,page=2]{echange_Lgn.pdf}}%
    \put(0.35420467,0.00442154){\color[rgb]{0,0,0}\makebox(0,0)[lt]{\lineheight{1.25}\smash{\begin{tabular}[t]{l}$\overset{I}{U}(\alpha)$\end{tabular}}}}%
    \put(0,0){\includegraphics[width=\unitlength,page=3]{echange_Lgn.pdf}}%
    \put(0.22209417,0.1138746){\color[rgb]{0,0,0}\makebox(0,0)[lt]{\lineheight{1.25}\smash{\begin{tabular}[t]{l}$J$\end{tabular}}}}%
    \put(0,0){\includegraphics[width=\unitlength,page=4]{echange_Lgn.pdf}}%
    \put(0.48609749,0.24575849){\color[rgb]{0,0,0}\makebox(0,0)[lt]{\lineheight{1.25}\smash{\begin{tabular}[t]{l}$I$\end{tabular}}}}%
    \put(0.60395165,0.0041878){\color[rgb]{0,0,0}\makebox(0,0)[lt]{\lineheight{1.25}\smash{\begin{tabular}[t]{l}$\overset{I}{U}(\alpha)$\end{tabular}}}}%
    \put(0,0){\includegraphics[width=\unitlength,page=5]{echange_Lgn.pdf}}%
    \put(0.80258334,0.0041878){\color[rgb]{0,0,0}\makebox(0,0)[lt]{\lineheight{1.25}\smash{\begin{tabular}[t]{l}$\overset{J}{V}(\beta)$\end{tabular}}}}%
    \put(0,0){\includegraphics[width=\unitlength,page=6]{echange_Lgn.pdf}}%
    \put(0.68256546,0.24729829){\color[rgb]{0,0,0}\makebox(0,0)[lt]{\lineheight{1.25}\smash{\begin{tabular}[t]{l}$J$\end{tabular}}}}%
    \put(0.41573488,0.06830278){\color[rgb]{0,0,0}\makebox(0,0)[lt]{\lineheight{1.25}\smash{\begin{tabular}[t]{l}=\end{tabular}}}}%
    \put(0,0){\includegraphics[width=\unitlength,page=7]{echange_Lgn.pdf}}%
  \end{picture}%
\endgroup%

\end{equation}

\section{The Wilson loop map}
\indent In what follows, we will consider framed links which are oriented and colored, up to isotopy (equivalently, oriented and colored ribbons up to isotopy). By colored we mean that any connected component of the link is labelled by a $H$-module. We denote by $\mathcal{R}^{\mathrm{OC}}_{g,n}$ the set of isotopy classes of oriented, framed and colored links in $\Sigma_{g,n}^{\mathrm{o}} \times [0, 1]$, and by $\mathbb{C}\mathcal{R}^{\mathrm{OC}}_{g,n}$ the $\mathbb{C}$-vector space whose basis is $\mathcal{R}^{\mathrm{OC}}_{g,n}$ (formal linear combinations of elements of $\mathcal{R}^{\mathrm{OC}}_{g,n}$).

\smallskip

\indent Recall the view of $\Sigma_{g,n}^{\mathrm{o}}$ depicted in Figure \ref{surfaceAvecMatrices} and assume that it represents $\Sigma_{g,n}^{\mathrm{o}} \times \{ 0 \} \subset \Sigma_{g,n}^{\mathrm{o}} \times [0,1]$ (thickened surface). If we have a (framed oriented) link $L \in \Sigma_{g,n}^{\mathrm{o}} \times [0,1]$, we may assume up to isotopy that each of the thickened handles simply contains a bunch parallel arcs (\textit{i.e.} it does not contains cups, caps or crossings) and that the thickened rectangle contains a $(m,0)$-tangle (with $m$ even) projecting onto $\Sigma_{g,n}^{\mathrm{o}} \times \{ 0 \}$, as follows:
\begin{center}
\begingroup%
  \makeatletter%
  \providecommand\color[2][]{%
    \errmessage{(Inkscape) Color is used for the text in Inkscape, but the package 'color.sty' is not loaded}%
    \renewcommand\color[2][]{}%
  }%
  \providecommand\transparent[1]{%
    \errmessage{(Inkscape) Transparency is used (non-zero) for the text in Inkscape, but the package 'transparent.sty' is not loaded}%
    \renewcommand\transparent[1]{}%
  }%
  \providecommand\rotatebox[2]{#2}%
  \newcommand*\fsize{\dimexpr\f@size pt\relax}%
  \newcommand*\lineheight[1]{\fontsize{\fsize}{#1\fsize}\selectfont}%
  \ifx\svgwidth\undefined%
    \setlength{\unitlength}{321.23092933bp}%
    \ifx\svgscale\undefined%
      \relax%
    \else%
      \setlength{\unitlength}{\unitlength * \real{\svgscale}}%
    \fi%
  \else%
    \setlength{\unitlength}{\svgwidth}%
  \fi%
  \global\let\svgwidth\undefined%
  \global\let\svgscale\undefined%
  \makeatother%
  \begin{picture}(1,0.28915527)%
    \lineheight{1}%
    \setlength\tabcolsep{0pt}%
    \put(0,0){\includegraphics[width=\unitlength,page=1]{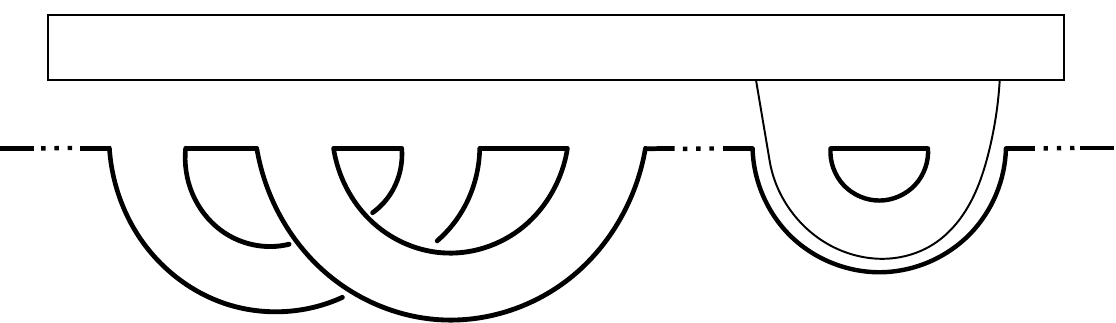}}%
    \put(0.37437431,0.23570034){\color[rgb]{0,0,0}\makebox(0,0)[lt]{\lineheight{1.25}\smash{\begin{tabular}[t]{l}$(m,0)$-tangle $T$\end{tabular}}}}%
    \put(0,0){\includegraphics[width=\unitlength,page=2]{formeStandard.pdf}}%
    \put(-0.07751617,-0.5182153){\color[rgb]{0,0,0}\makebox(0,0)[lt]{\begin{minipage}{0.04602984\unitlength}\raggedright \end{minipage}}}%
  \end{picture}%
\endgroup%

\end{center}
\noindent Hence we see that $L$ can be represented by a (non-unique) $(m,0)$-tangle $T$. The non-unicity comes from the fact that we can drag crossings, cups and caps of $T$ along the handles and obtain another $T'$ which also represents the link $L$.

\begin{definition}\label{defDefWilson}
Let $L$ be an oriented and colored framed link represented by a $(m,0)$-tangle as explained above. The Wilson loop around $L$ is an element $W(L) \in \mathcal{L}_{g,n}(H)$ defined as the evaluation of the diagram at the bottom of Figure \ref{defWilson} (recall Definition \ref{defEvaluation}). We extend $W$ to $\mathbb{C}\mathcal{R}^{\mathrm{OC}}_{g,n}$ by linearity and this gives a map $W : \mathbb{C}\mathcal{R}^{\mathrm{OC}}_{g,n} \to \mathcal{L}_{g,n}(H)$, which we call the Wilson loop map. 
\begin{figure}[!h]
\centering
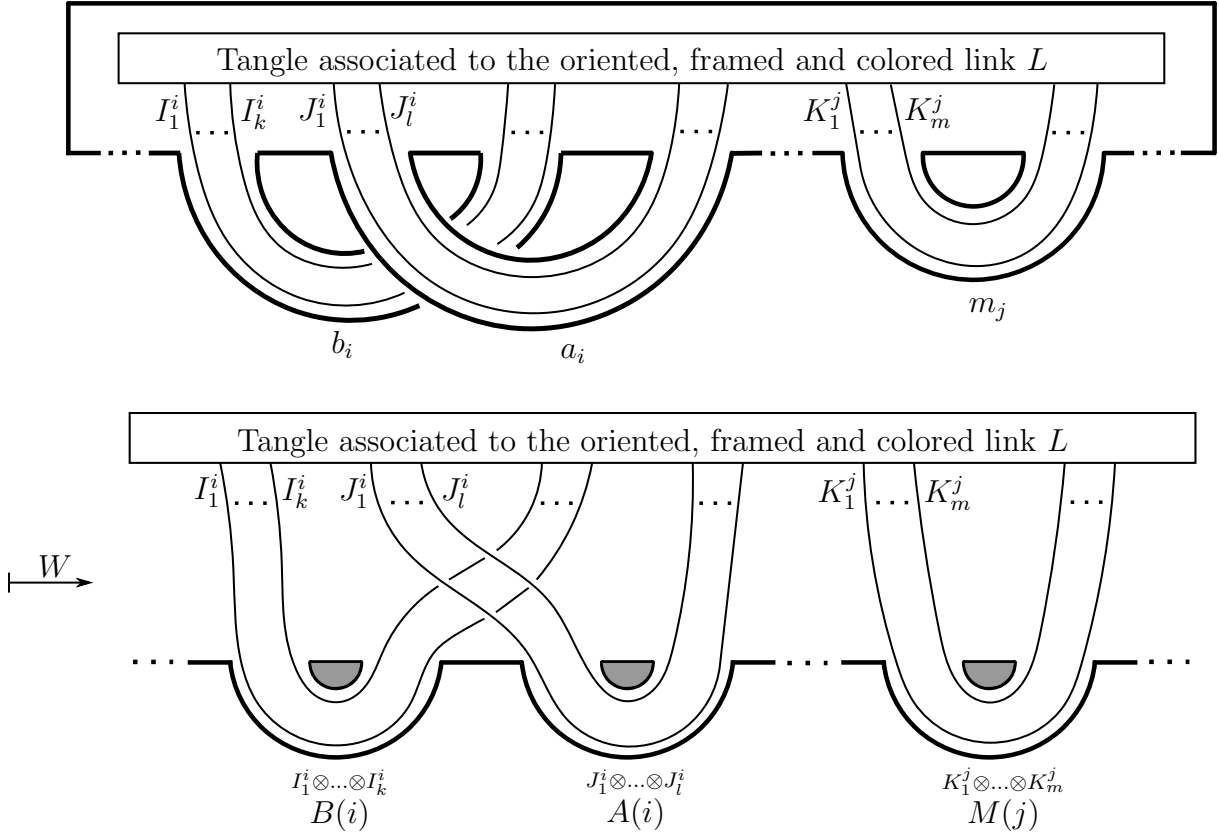
\caption{Definition of the Wilson loop map $W$.}
\label{defWilson}
\end{figure}
\end{definition}
\noindent Observe that since a tangle $T$ representing the link $L$ is of type $(m,0)$ (no outgoing strands), $F_{\mathrm{RT}}(T)$ is a morphism with values in $\mathbb{C}$. Hence, by Definition \ref{defEvaluation}, the evaluation of the diagram in Figure \ref{defWilson} is indeed an element of $\mathcal{L}_{g,n}(H) \otimes \mathbb{C} = \mathcal{L}_{g,n}(H)$.

\smallskip

\begin{figure}[!h]
\centering
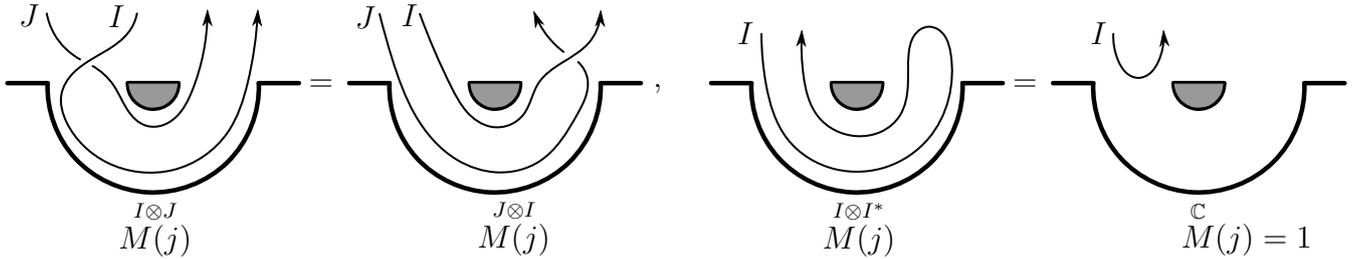
\caption{Examples of consequences of the naturality.}
\label{WilsonBienDef}
\end{figure}

\indent We note that $W(L)$ does not depend on the choice of a tangle $T$ representing $L$. First, since the Reshetikhin-Turaev functor is an isotopy invariant, the evaluation of the diagram in Figure \ref{defWilson} depends only on the isotopy class of $T$. Moreover, if we drag certain crossings, cups and caps along the handles in order to obtain another tangle $T'$ representing $L$, then this does not change the value of $W(L)$ thanks to naturality \eqref{naturalite}, see Figure \ref{WilsonBienDef}.

\indent We stress that the diagrammatic rules introduced above allow us to compute the value of $W$ in a purely graphical way. In order to clarify the definition of $W$, we will compute two examples below. Recall from \eqref{WAWBWBA} the notations $\overset{I}{W}_{\! A} = \mathrm{tr}_q(\overset{I}{A})$ and $\overset{I}{W}_{\! B} = \mathrm{tr}_q(\overset{I}{B})$, which corresponds to Wilson loops around the loops $a$ and $b$. More generally if $\overset{I}{X}$ is a matrix of the form \eqref{matriceGenerale} we let $\overset{I}{W}_{\! X} = \mathrm{tr}_q(\overset{I}{X})$, according to \eqref{defWexplicite} and \eqref{notationEmbed}.

\begin{exemple}
In Figure \ref{exempleWilson1}, we compute the value of the Wilson loop around the simple closed curve $b^{-1}a \subset \Sigma_{1,0}^{\mathrm{o}} \times \{0\}$. The result is not surprising: this is simply the quantum trace of the lift of the simple closed curve $b^{-1}a$, that is the quantum trace of the holonomy of the closed curve $b^{-1}a$. This property is always true, see Proposition \ref{propWilsonSimple}.
\begin{figure}
\centering
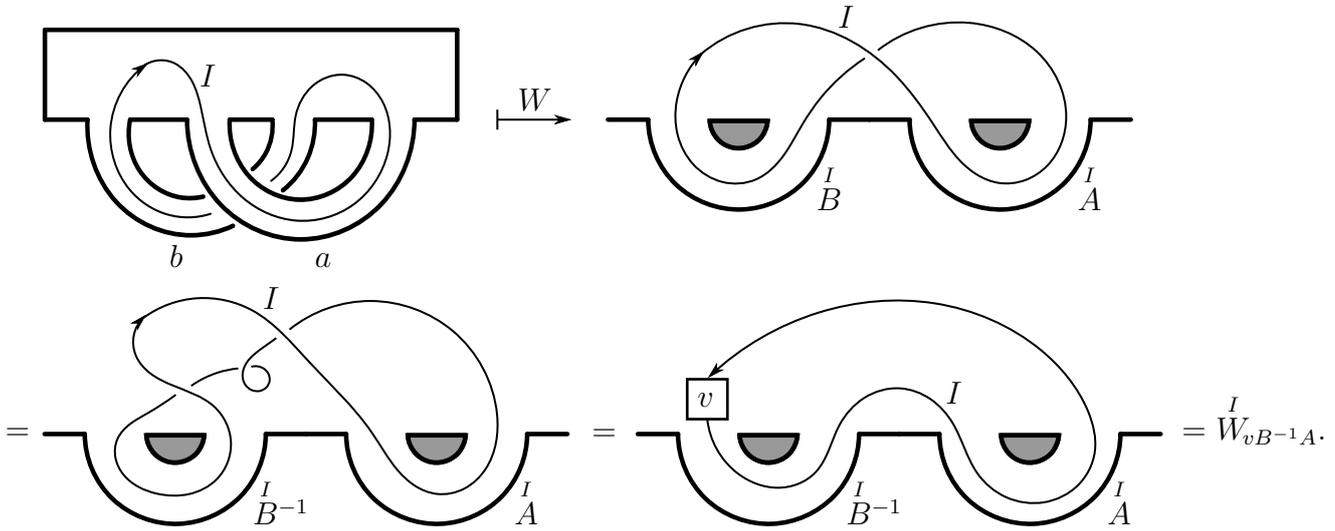
\caption{Example of computation of a Wilson loop.}
\label{exempleWilson1}
\end{figure}
\finEx
\end{exemple}

\indent We now state the properties of $W$; they are all natural-looking. The first property is that Wilson loops are invariant elements.
\begin{proposition}\label{wilsonInv}
For any $L \in \mathbb{C}\mathcal{R}^{\mathrm{OC}}_{g,n}$, it holds $W(L) \in \mathcal{L}_{g,n}^{\mathrm{inv}}(H)$.
\end{proposition}
\begin{proof}
We give a detailed proof although it is rather obvious. Keep Figure \ref{defWilson} in mind. The link $L$ is represented by some $(m,0)$-tangle $T$. We can assume that the orientations of the strands in the handles are all positively oriented. Indeed, if in the original link a strand enters an handle with the negative orientation then we apply \eqref{sensOppose} and we push the coupons $\mathrm{id}_{I^*}$ and $e_I$ near the tangle $T$, obtaining a ribbon graph denoted $T'$. Moreover, we can assume that there is only one strand passing in each handle. Indeed, if there are several strands, we use the coupons \eqref{couponsId} below and we push them near $T'$, obtaining a new ribbon graph denoted $T''$.
\begin{equation}\label{couponsId}
\begingroup%
  \makeatletter%
  \providecommand\color[2][]{%
    \errmessage{(Inkscape) Color is used for the text in Inkscape, but the package 'color.sty' is not loaded}%
    \renewcommand\color[2][]{}%
  }%
  \providecommand\transparent[1]{%
    \errmessage{(Inkscape) Transparency is used (non-zero) for the text in Inkscape, but the package 'transparent.sty' is not loaded}%
    \renewcommand\transparent[1]{}%
  }%
  \providecommand\rotatebox[2]{#2}%
  \newcommand*\fsize{\dimexpr\f@size pt\relax}%
  \newcommand*\lineheight[1]{\fontsize{\fsize}{#1\fsize}\selectfont}%
  \ifx\svgwidth\undefined%
    \setlength{\unitlength}{228.77554018bp}%
    \ifx\svgscale\undefined%
      \relax%
    \else%
      \setlength{\unitlength}{\unitlength * \real{\svgscale}}%
    \fi%
  \else%
    \setlength{\unitlength}{\svgwidth}%
  \fi%
  \global\let\svgwidth\undefined%
  \global\let\svgscale\undefined%
  \makeatother%
  \begin{picture}(1,0.25223868)%
    \lineheight{1}%
    \setlength\tabcolsep{0pt}%
    \put(0,0){\includegraphics[width=\unitlength,page=1]{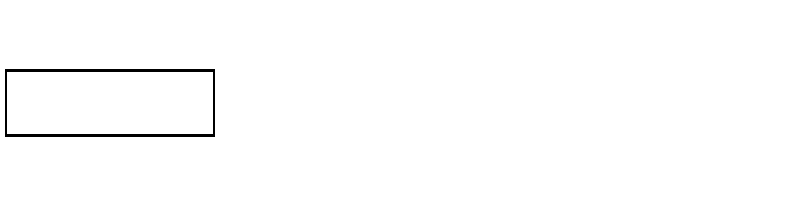}}%
    \put(0.11705895,0.10688433){\color[rgb]{0,0,0}\makebox(0,0)[lt]{\lineheight{1.25}\smash{\begin{tabular}[t]{l}$\mathrm{id}$\end{tabular}}}}%
    \put(0,0){\includegraphics[width=\unitlength,page=2]{couponIdBrinSimple.pdf}}%
    \put(-0.00134444,0.20507548){\color[rgb]{0,0,0}\makebox(0,0)[lt]{\lineheight{1.25}\smash{\begin{tabular}[t]{l}$V_1$\end{tabular}}}}%
    \put(0.23465952,0.20472729){\color[rgb]{0,0,0}\makebox(0,0)[lt]{\lineheight{1.25}\smash{\begin{tabular}[t]{l}$V_l$\end{tabular}}}}%
    \put(0.15657812,0.02723877){\color[rgb]{0,0,0}\makebox(0,0)[lt]{\lineheight{1.25}\smash{\begin{tabular}[t]{l}$V_1 \otimes \ldots \otimes V_l$\end{tabular}}}}%
    \put(0,0){\includegraphics[width=\unitlength,page=3]{couponIdBrinSimple.pdf}}%
    \put(0.704046,0.10590987){\color[rgb]{0,0,0}\makebox(0,0)[lt]{\lineheight{1.25}\smash{\begin{tabular}[t]{l}$\mathrm{id}$\end{tabular}}}}%
    \put(0,0){\includegraphics[width=\unitlength,page=4]{couponIdBrinSimple.pdf}}%
    \put(0.58794614,0.20507548){\color[rgb]{0,0,0}\makebox(0,0)[lt]{\lineheight{1.25}\smash{\begin{tabular}[t]{l}$V_1$\end{tabular}}}}%
    \put(0.82475742,0.20472729){\color[rgb]{0,0,0}\makebox(0,0)[lt]{\lineheight{1.25}\smash{\begin{tabular}[t]{l}$V_l$\end{tabular}}}}%
    \put(0.7411637,0.02586075){\color[rgb]{0,0,0}\makebox(0,0)[lt]{\lineheight{1.25}\smash{\begin{tabular}[t]{l}$V_1 \otimes \ldots \otimes V_l$\end{tabular}}}}%
    \put(0,0){\includegraphics[width=\unitlength,page=5]{couponIdBrinSimple.pdf}}%
  \end{picture}%
\endgroup%

\end{equation}
The diagram of $W(L)$ now looks as follows:
\begin{center}
\begingroup%
  \makeatletter%
  \providecommand\color[2][]{%
    \errmessage{(Inkscape) Color is used for the text in Inkscape, but the package 'color.sty' is not loaded}%
    \renewcommand\color[2][]{}%
  }%
  \providecommand\transparent[1]{%
    \errmessage{(Inkscape) Transparency is used (non-zero) for the text in Inkscape, but the package 'transparent.sty' is not loaded}%
    \renewcommand\transparent[1]{}%
  }%
  \providecommand\rotatebox[2]{#2}%
  \newcommand*\fsize{\dimexpr\f@size pt\relax}%
  \newcommand*\lineheight[1]{\fontsize{\fsize}{#1\fsize}\selectfont}%
  \ifx\svgwidth\undefined%
    \setlength{\unitlength}{541.11151717bp}%
    \ifx\svgscale\undefined%
      \relax%
    \else%
      \setlength{\unitlength}{\unitlength * \real{\svgscale}}%
    \fi%
  \else%
    \setlength{\unitlength}{\svgwidth}%
  \fi%
  \global\let\svgwidth\undefined%
  \global\let\svgscale\undefined%
  \makeatother%
  \begin{picture}(1,0.20563362)%
    \lineheight{1}%
    \setlength\tabcolsep{0pt}%
    \put(0,0){\includegraphics[width=\unitlength,page=1]{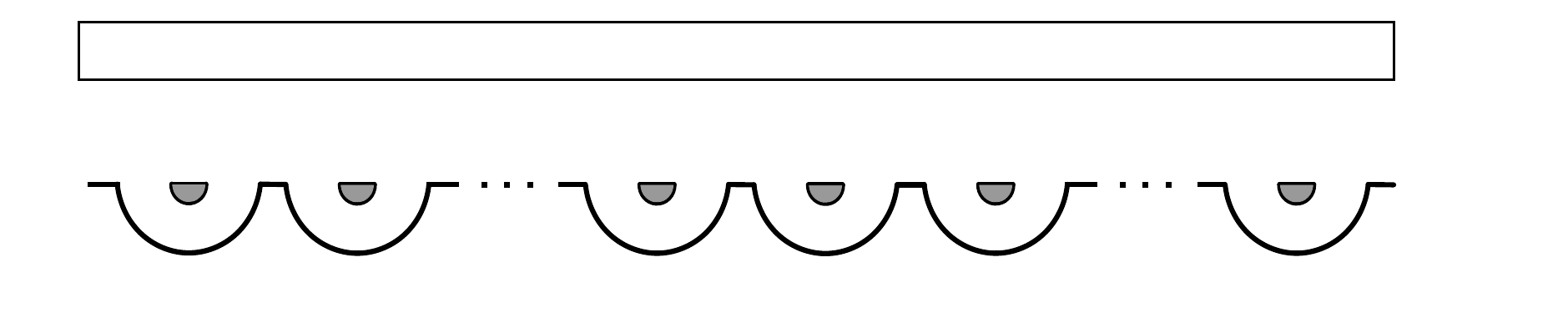}}%
    \put(0.47146724,0.16553898){\color[rgb]{0,0,0}\makebox(0,0)[lt]{\lineheight{1.25}\smash{\begin{tabular}[t]{l}$T''$\end{tabular}}}}%
    \put(0,0){\includegraphics[width=\unitlength,page=2]{tanglePreuveInv.pdf}}%
    \put(-0.00118481,0.1468969){\color[rgb]{0,0,0}\makebox(0,0)[lt]{\lineheight{1.25}\smash{\begin{tabular}[t]{l}$T_c$\end{tabular}}}}%
    \put(0.10507705,0.00407029){\color[rgb]{0,0,0}\makebox(0,0)[lt]{\lineheight{1.25}\smash{\begin{tabular}[t]{l}$\overset{I(1)}{B(1)}$\end{tabular}}}}%
    \put(0.21345864,0.00368778){\color[rgb]{0,0,0}\makebox(0,0)[lt]{\lineheight{1.25}\smash{\begin{tabular}[t]{l}$\overset{J(1)}{A(1)}$\end{tabular}}}}%
    \put(0.39988341,0.00358544){\color[rgb]{0,0,0}\makebox(0,0)[lt]{\lineheight{1.25}\smash{\begin{tabular}[t]{l}$\overset{I(g)}{B(g)}$\end{tabular}}}}%
    \put(0.50826501,0.00320292){\color[rgb]{0,0,0}\makebox(0,0)[lt]{\lineheight{1.25}\smash{\begin{tabular}[t]{l}$\overset{J(g)}{A(g)}$\end{tabular}}}}%
    \put(0.60303645,0.00315882){\color[rgb]{0,0,0}\makebox(0,0)[lt]{\lineheight{1.25}\smash{\begin{tabular}[t]{l}$\overset{K(g+1)}{M(g+1)}$\end{tabular}}}}%
    \put(0.78402223,0.00254659){\color[rgb]{0,0,0}\makebox(0,0)[lt]{\lineheight{1.25}\smash{\begin{tabular}[t]{l}$\overset{K(g+n)}{M(g+n)}$\end{tabular}}}}%
    \put(0.05079786,0.13855439){\color[rgb]{0,0,0}\makebox(0,0)[lt]{\lineheight{1.25}\smash{\begin{tabular}[t]{l}$_{I(1)}$\end{tabular}}}}%
    \put(0.10544494,0.13851027){\color[rgb]{0,0,0}\makebox(0,0)[lt]{\lineheight{1.25}\smash{\begin{tabular}[t]{l}$_{J(1)}$\end{tabular}}}}%
    \put(0.35266517,0.13809512){\color[rgb]{0,0,0}\makebox(0,0)[lt]{\lineheight{1.25}\smash{\begin{tabular}[t]{l}$_{I(g)}$\end{tabular}}}}%
    \put(0.40731222,0.13805101){\color[rgb]{0,0,0}\makebox(0,0)[lt]{\lineheight{1.25}\smash{\begin{tabular}[t]{l}$_{J(g)}$\end{tabular}}}}%
    \put(0.60040604,0.13829318){\color[rgb]{0,0,0}\makebox(0,0)[lt]{\lineheight{1.25}\smash{\begin{tabular}[t]{l}$_{K(g+1)}$\end{tabular}}}}%
    \put(0.73140578,0.13800193){\color[rgb]{0,0,0}\makebox(0,0)[lt]{\lineheight{1.25}\smash{\begin{tabular}[t]{l}$_{K(g+n)}$\end{tabular}}}}%
  \end{picture}%
\endgroup%

\end{center}
where $I(i), J(i), K(j)$ are tensor products as in \eqref{couponsId}. Let $\bigl(u(\alpha)_i\bigr)$ (resp. $\bigl(v(\alpha)_i\bigr)$, $\bigl(w(\beta)_i\bigr)$) be a basis of $I(\alpha)$ (resp. $J(\alpha)$, $K(\beta)$), then
\begin{align*}
&W(L) \cdot h \\
&= \prod_{\alpha=1}^g \bigl(\overset{I(\alpha)}{h^{(4\alpha-3)}} \overset{I(\alpha)}{B}\!_{\!}(\alpha) \overset{I(\alpha)}{S(h^{(4\alpha-2)})} \bigr)^{i_{\alpha}}_{j_{\alpha}} \bigl(\overset{J(\alpha)}{h^{(4\alpha-1)}} \overset{J(\alpha)}{A}\!\!(\alpha) \overset{J(\alpha)}{S(h^{(4\alpha)})} \bigr)^{k_{\alpha}}_{l_{\alpha}} \prod_{\beta=g+1}^{g+n} \bigl( \overset{K(\beta)}{h^{(4g+2\beta-1)}} \overset{K(\beta)}{M}\!\!(\beta) \overset{K(\beta)}{S(h^{(4g+2\beta)})} \bigr)^{m_{\beta}}_{p_{\beta}} \\
& \quad F_{\mathrm{RT}}(T_c)\biggl(\bigotimes_{\alpha = 1}^g u(\alpha)_{i_{\alpha}} \otimes u(\alpha)^{j_{\alpha}} \otimes v(\alpha)_{k_{\alpha}} \otimes v(\alpha)^{l_{\alpha}} \otimes \bigotimes_{\beta = g+1}^{g+n} w(\beta)_{m_{\beta}} \otimes w(\beta)^{p_{\beta}} \biggr)\\
&= \prod_{\alpha=1}^g \overset{I(\alpha)}{B}\!_{\!}(\alpha){^{i_{\alpha}}_{j_{\alpha}}} \overset{J(\alpha)}{A}\!\!(\alpha){^{k_{\alpha}}_{l_{\alpha}}} \prod_{\beta=g+1}^{g+n} \overset{K(\beta)}{M}\!\!(\beta){^{m_{\beta}}_{p_{\beta}}}
 F_{\mathrm{RT}}(T_c)\biggl(\bigotimes_{\alpha = 1}^g h^{(4\alpha-3)} u(\alpha)_{i_{\alpha}} \otimes h^{(4\alpha-2)} u(\alpha)^{j_{\alpha}}\\
& \quad \otimes h^{(4\alpha-1)} v(\alpha)_{k_{\alpha}} \otimes h^{(4\alpha)} v(\alpha)^{l_{\alpha}} \otimes \bigotimes_{\beta = g+1}^{g+n} h^{(4g+2\beta-1)} w(\beta)_{m_{\beta}} \otimes h^{(4g+2\beta)}w(\beta)^{p_{\beta}} \biggr) = \varepsilon(h) W(L).
\end{align*}
We simply used \eqref{actionHLgn}, \eqref{trucEvidentdAlgebreLineaire} and the fact that $F_{\mathrm{RT}}(T_c)$ is a morphism $\bigotimes_{\alpha=1}^g I(\alpha) \otimes I(\alpha)^* \otimes J(\alpha) \otimes J(\alpha)^* \otimes \bigotimes_{\beta = g+1}^{g+n} K(\beta) \otimes K(\beta)^* \to \mathbb{C}$.
\end{proof}

\indent Let us see the behaviour of $W$ under change of orientation. This is the generalization of the corresponding fact for the Reshetikhin-Turaev functor (see \textit{e.g.} \cite[Lemma 3.18]{KM}). Let $L$ be an oriented framed link; we have $L = (L_1, I_1) \sqcup \ldots \sqcup (L_k, I_k)$ where the $L_i$'s are the connected components of $L$ and the $I_i$'s are their colors. We denote by $L_i^{-1}$ the oriented framed curve whose orientation is opposite to that of $L_i$.
\begin{proposition}\label{orientationEtWilson}
With the notation above, it holds
\[ W\bigl((L^{-1}_1, I_1) \sqcup \ldots \sqcup (L_k, I_k)\bigr) = W\bigl( (L_1, I^*_1) \sqcup \ldots \sqcup (L_k, I_k) \bigr). \]
\end{proposition}
\begin{proof}
We can assume that $L = L_1$. The result follows from the application of the local equalities \eqref{dualitePrime}, \eqref{sensOppose} and \eqref{sensOpposeBis} on $W(L_1)$ together with the fact that the diagram $W(L_1)$ contains an equal number of cups and caps (a part of strand crossing a handle is considered as a cup).
\end{proof}

\indent Now we show that $W$ is compatible with the stack product.
\begin{definition}\label{stackProduct}
Let $L_1, L_2 \in \mathcal{R}^{\mathrm{OC}}_{g,n}$ and let $L^-_1 \in \Sigma_{g,n}^{\mathrm{o}} \times [0, \frac{1}{2}[$ be isotopic to $L_1$ and $L^+_2 \in \Sigma_{g,n}^{\mathrm{o}} \times ]\frac{1}{2}, 1]$ be isotopic to $L_2$. The stack product of $L_1$ and $L_2$ is $L_1 \ast L_2 = L^-_1 \cup L^+_2 \in \mathcal{R}^{\mathrm{OC}}_{g,n}$. We extend the stack product bilinearly to $\mathbb{C}\mathcal{R}^{\mathrm{OC}}_{g,n}$.
\end{definition}

\begin{exemple}
In Figure \ref{exempleWilson2}, we compute the value of $W$ on the stack product $a \ast b$. We get that the Wilson loop around the stack product of the two links is the product of the Wilson loops around each of these links: $W(a \ast b) = W(a)W(b)$. This is a general fact, as we shall see now.
\begin{figure}
\centering
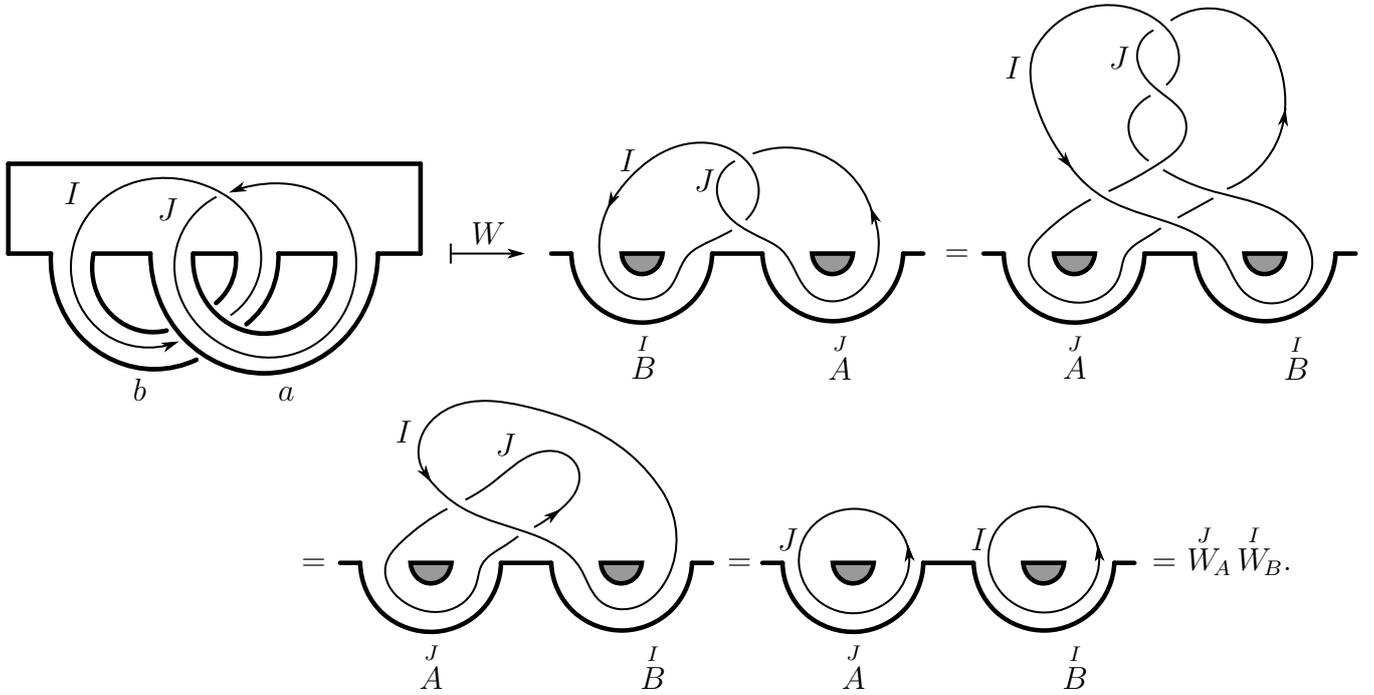
\caption{Example of computation of a Wilson loop.}
\label{exempleWilson2}
\end{figure}
\finEx
\end{exemple}

\begin{figure}[p]
\centering
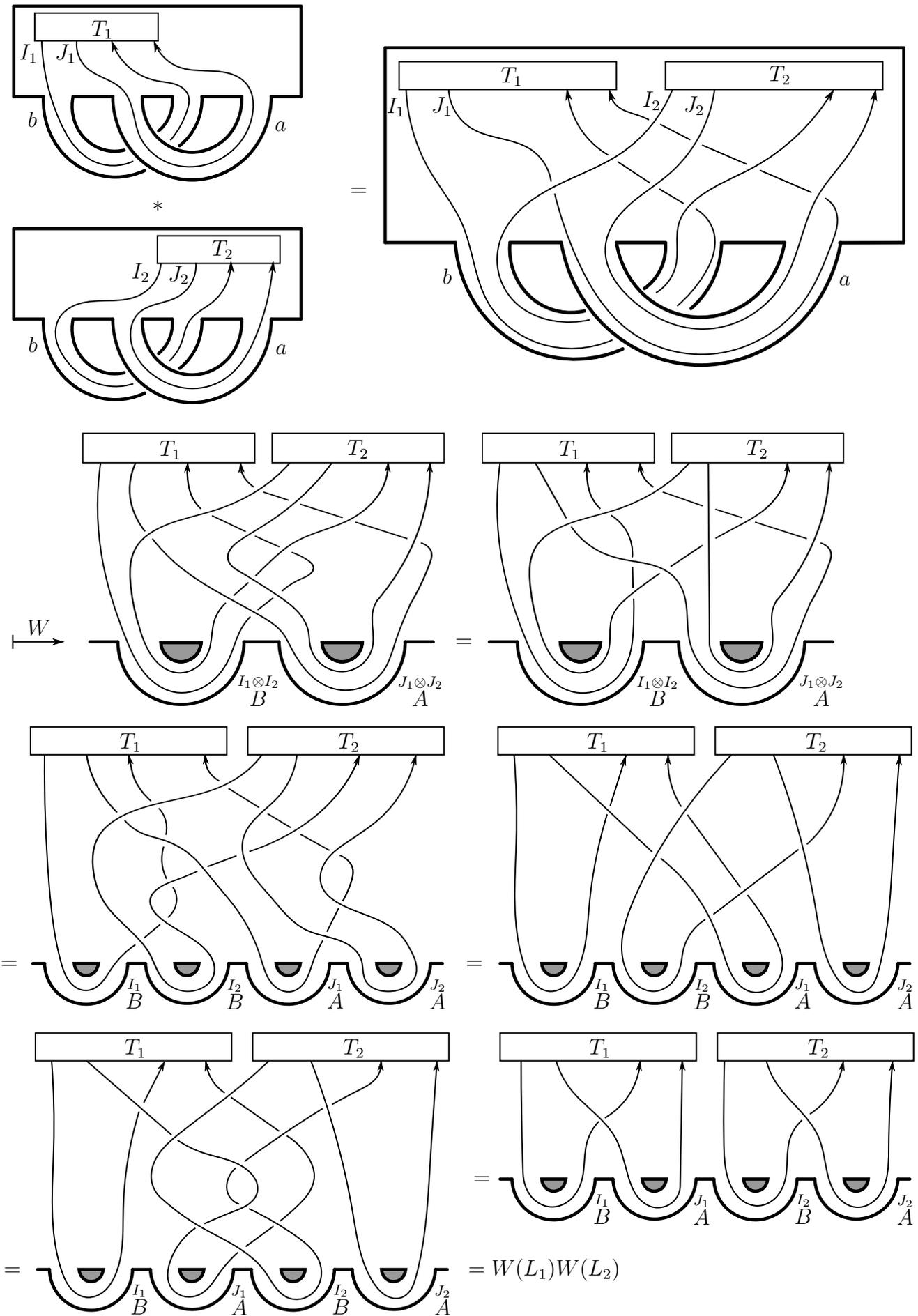
\caption{Compatibility of $W$ with the stack product.}
\label{preuveWilson}
\end{figure}

\begin{theorem}\label{wilsonStack}
The map $W : \mathbb{C}\mathcal{R}^{\mathrm{OC}}_{g,n} \to \mathcal{L}_{g,n}(H)$ is a morphism of algebras:
\[ W(L_1 \ast L_2) = W(L_1)W(L_2). \]
\end{theorem}
\begin{proof}
The proof is purely diagrammatic. For instance, the proof for $g=1, n=0$ is depicted in Figure \ref{preuveWilson}. We used the definition of $W$ (Figure \ref{defWilson}), the fusion relation \eqref{relationFusion}, the exchange relation \eqref{dessinEchangeL10} and obvious topological simplifications. The proof for the general case is similar: use the definition of $W$, then use the fusion relation and simplify the diagram, and finally use the exchange relations and simplify the diagram. As explained in the proof of Proposition \ref{wilsonInv}, we can assume that all the strands are positively oriented when they go through the handles and that there is only one strand in each handle.
\end{proof}

\smallskip

\indent Finally, we show that the Wilson loop around a simple closed curve (simple loop or oriented circle\footnote{Recall that according to the terminology fixed in section \ref{MCGSigmaG}, a circle is a simple closed curve up to free homotopy (not necessarily oriented, unless stated) while a loop (element of the fundamental group) is a based oriented curve up to fixed-basepoint homotopy (not necessarily simple, unless stated).}) is simply the quantum trace of its holonomy. We restrict to $\Sigma_g^{\mathrm{o}}$, where $\Sigma_g = \Sigma_{g,0}$; however, Proposition \ref{propWilsonSimple} and its corollaries are undoubtedly true for $\Sigma_{g,n}^{\mathrm{o}}$.
\\
\indent An oriented circle $\gamma \subset \Sigma_g^{\mathrm{o}}$ can be seen as an oriented framed link lying in $\Sigma_g^{\mathrm{o}} \times \{0\}$ and thus $W(\gamma)$ makes sense. Let $x \in \pi_1(\Sigma_{g}^{\mathrm{o}})$ be a simple loop colored by some $H$-module $I$ and let $[x]$ be its free homotopy class. We define
\[ W(x) = W\bigl([x]\bigr)\]
Recall the lift $\overset{I}{\widetilde{x}}$ of a simple loop $x$ from Definition \ref{defLiftLoop} (express $x$ in terms of the generators $b_i,a_i$ of $\pi_1(\Sigma_g^{\mathrm{o}})$ and replace the loops $b_i,a_i$ by the matrices $\overset{I}{B}(i), \overset{I}{A}(i)$ up to the normalization $\overset{I}{v}{^{N(x)}}$, where $N(x)$ is defined in section \ref{sectionNorma}) and the lift $\widetilde{f}$ of a homeomorphism $f \in \mathrm{MCG}(\Sigma_{g}^{\mathrm{o}})$ from Definition \ref{liftHomeo} (it satisfies $\widetilde{f}\bigl(\overset{I}{\widetilde{x}}\bigr) = \overset{I}{\widetilde{f(x)}}$ for any simple loop $x \in \pi_1(\Sigma_g^{\mathrm{o}})$).

\begin{lemma}\label{MCGcommuteW}
Let $f \in \mathrm{MCG}(\Sigma_{g}^{\mathrm{o}})$ and let $\gamma \subset \Sigma_{g}^{\mathrm{o}} \times \{0\}$ be an oriented circle (colored by a $H$-module), then
\[ W\bigl( f(\gamma) \bigr) = \widetilde{f}\bigl( W(\gamma) \bigr). \]
\end{lemma}
\begin{proof}
We can assume that $f$ is one of the Humphries generators. The proof is purely diagrammatic. Let $I$ be the color of $\gamma$. We represent $\gamma$ by a tangle $T$ which does not contains crossings. As explained in the proof of Proposition \ref{wilsonInv}, we can assume that all the strands are positively oriented when they go through the handles and that there is only one strand in each handle. Here since there is just one color, a strand going through a handle will be colored by $I(\epsilon) = I^{\epsilon_1} \otimes \ldots \otimes I^{\epsilon_k}$, where $\epsilon = (\epsilon_1, \ldots, \epsilon_k)$ is a sequence of $\pm$ signs and $I^+ = I, I^- = I^*$. 
Now, if $f = \tau_{a_1}$, we can restrict to $\Sigma_1^{\mathrm{o}}$ and we perform the graphical computation represented in Figure \ref{twistAetMCG}. We used the fusion relation and the reflection equation \eqref{equationReflexionBis}. The equalities for the others Humphries generators are shown similarly (but the diagrammatic computations are more cumbersome): for $\tau_{b_i}$ we can also restrict to $\Sigma_1^{\mathrm{o}}$ and for $\tau_{d_i}, \tau_{e_2}$ we can restrict to $\Sigma_2^{\mathrm{o}}$.
\begin{figure}[p]
\centering
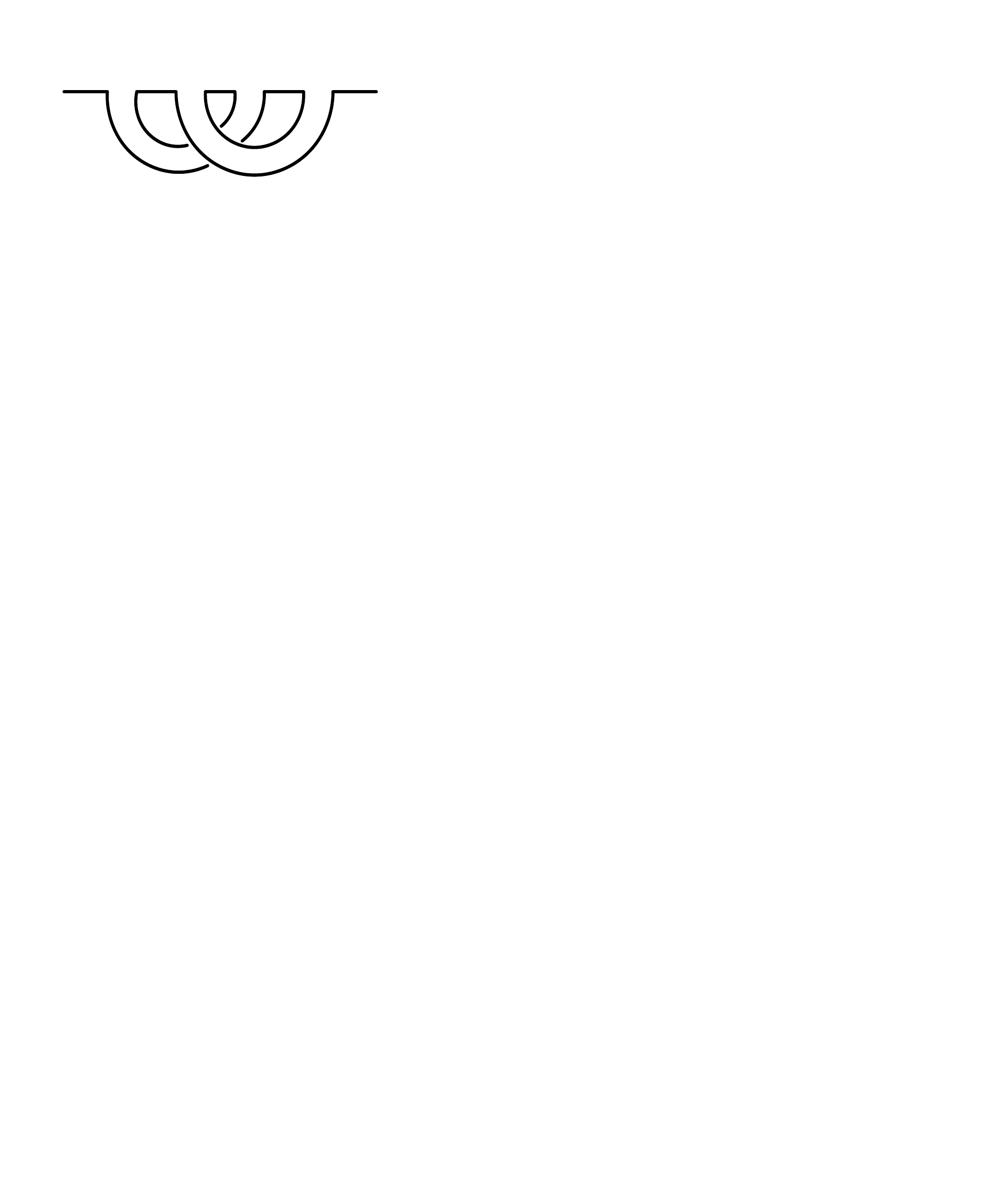
\caption{Proof of the equality $W \circ \tau_a(\gamma) = \widetilde{\tau_a} \circ W(\gamma)$.}
\label{twistAetMCG}
\end{figure}
\end{proof}

\begin{proposition}\label{propWilsonSimple}
Let $x \in \pi_1(\Sigma_{g}^{\mathrm{o}})$ be a simple loop colored by a $H$-module $I$, then
\[ W(x) = \mathrm{tr}_q\bigl(\overset{I}{\widetilde{x}}\bigr). \]
\end{proposition}
\begin{proof}
Assume first that $x$ is positively oriented (recall Definition \ref{defPositivelyOriented}). Then by Lemma \ref{transfoLoop} there exists a homeomorphism $f$ such that $x = f(y)$ where $y$ is either $a_1$ or $s_1 \ldots$ or $s_g$ (recall \eqref{courbesGenHumphries}). One can check by direct computation that the result is true for these particular loops. Note that $[f(y)] = f([y])$ and hence, due to Lemmas \ref{MCGcommuteW} and \ref{pasSurprenant}, we get:
\[ W(x) = W\bigl(f(y)\bigr) = \widetilde{f}\bigl(W(y)\bigr) = \widetilde{f}\biggl(\mathrm{tr}_q\bigl(\overset{I}{\widetilde{y}}\bigr)\biggr) = \mathrm{tr}_q\biggl(\widetilde{f}\bigl(\overset{I}{\widetilde{y}}\bigr)\biggr) = \mathrm{tr}_q\bigl(\overset{I}{\widetilde{f(y)}}\bigr) = \mathrm{tr}_q\bigl( \overset{I}{\widetilde{x}}\bigr). \]
Assume now that $x$ is negatively oriented. Thanks to Proposition \ref{orientationEtWilson}, we have $W(x) = W\bigl((x^{-1})^*\bigr)$ where $(x^{-1})^*$ is $x^{-1}$ colored by $I^*$. Moreover, since $x^{-1}$ is positively oriented, we have a morphism $j_{\widetilde{x^{-1}}} : \mathcal{L}_{0,1}(H) \to \mathcal{L}_{g,0}(H)$ defined by $\overset{J}{M} \mapsto \overset{J}{\widetilde{x^{-1}}}$ (see \eqref{morphismeEmbed}). Hence applying $j_{\widetilde{x^{-1}}}$ to \eqref{traceInverseTraceDual} yields $\mathrm{tr}_q\biggl(\overset{I^*}{\widetilde{x^{-1}}}\biggr) = \mathrm{tr}_q\biggl(\bigl(\overset{I}{\widetilde{x^{-1}}}\bigr)^{-1}\biggr)$. Thus, using the formula just established for positively oriented circles and the definition of the lift of a negatively oriented loop, we get
\[ W(x) = W\bigl((x^{-1})^*\bigr) = \mathrm{tr}_q\biggl( \overset{I^*}{\widetilde{x^{-1}}} \biggr) = \mathrm{tr}_q\biggl( \bigl(\overset{I}{\widetilde{x^{-1}}}\bigr)^{-1} \biggr) = \mathrm{tr}_q\bigl( \overset{I}{\widetilde{x}} \bigr), \]
as desired.
\end{proof}

\begin{remark}\label{pourquoiDefLift}
The second part of the above proof (and more precisely the formula \eqref{traceInverseTraceDual}) reveals why we have been forced to distinguish the positively oriented case from the negatively oriented case when we have defined the lift of a simple loop (Definition \ref{defLiftLoop}). Indeed, if we use the same formula to define the lifts of the positively oriented simple loops and the negatively oriented simple loops, then we get that the lift of $a_1^{-1}$ (for instance) is $\overset{I}{v} \overset{I}{A}(1)^{-1}$; but $W(a_1^{-1}) = \mathrm{tr}_q\bigl( \overset{I}{A}(1)^{-1}\bigr)$ and thus with this definition $W(a_1^{-1})$ would not have been the trace of the holonomy of $a_1^{-1}$.
\finEx
\end{remark}

\begin{corollary}\label{indepLoop}
If $x, y \in \pi_1(\Sigma_{g,0}^{\mathrm{o}})$ are simple loops colored by a $H$-module $I$ and such that $[x] = [y]$, it holds:
\[ \mathrm{tr}_q( \overset{I}{\widetilde{x}}) = \mathrm{tr}_q( \overset{I}{\widetilde{y}}). \]
\end{corollary}
\begin{proof}
$\mathrm{tr}_q( \overset{I}{\widetilde{x}}) = W(x) = W\bigl([x]\bigr) = W\bigl([y]\bigr) = W(y) = \mathrm{tr}_q( \overset{I}{\widetilde{y}}).$
\end{proof}

\begin{corollary}
Let $\gamma \subset \Sigma_{g,0}^{\mathrm{o}} \times \{0\}$ be an oriented circle colored by a $H$-module $I$, then
\[ W(\gamma) = \mathrm{tr}_q\bigl(\overset{I}{\widetilde{\gamma}}\bigr). \]
Here $\overset{I}{\widetilde{\gamma}}$ is a lift of $\gamma$, defined by $\overset{I}{\widetilde{\gamma}} = \overset{I}{\widetilde{y}}$, where $y \in \pi_1(\Sigma_{g}^{\mathrm{o}})$ is a simple loop such that $[y] = \gamma$.
\end{corollary}

\begin{figure}[h]
\centering
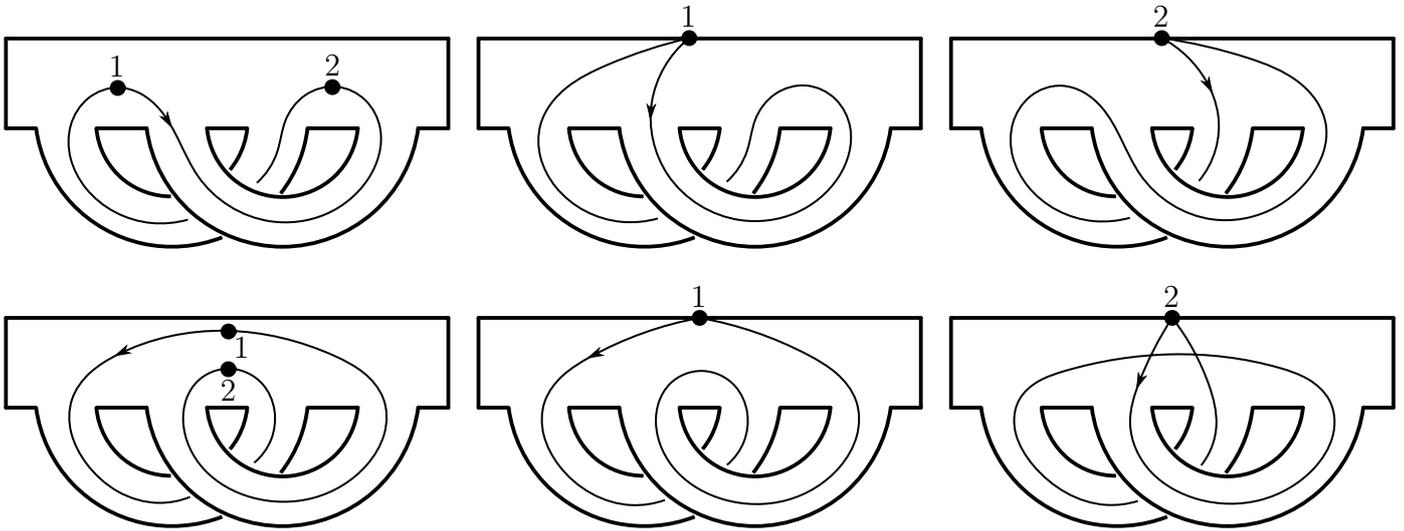
\caption{Basepoints on the oriented circles $[b^{-1}a]$ and $[ba]$.}
\label{changementPointBase}
\end{figure}

\indent To conclude this section we mention that, thanks to Corollary \ref{indepLoop}, different choices of basepoints on an oriented circle imply equalities between traces. For instance in $\mathcal{L}_{1,0}(H)$, the equality
\[ \mathrm{tr}_q\bigl( \overset{I}{v} \overset{I}{B}{^{-1}} \overset{I}{A} \bigr) = \mathrm{tr}_q\bigl( \overset{I}{v}{^{-1}} \overset{I}{A} \overset{I}{B}{^{-1}}\bigr) \]
follows from choosing two basepoints on the circle $[b^{-1}a] = [ab^{-1}]$, which are depicted at the top of Figure \ref{changementPointBase}. Reversing the orientation, we similarly obtain $\mathrm{tr}_q\bigl( \overset{I}{v} \overset{I}{B} \overset{I}{A}{^{-1}}\bigr) = \mathrm{tr}_q\bigl( \overset{I}{v}{^{-1}}\overset{I}{A}{^{-1}}\overset{I}{B}\bigr)$. In contrast,
\[ \mathrm{tr}_q\bigl( \overset{I}{v}{^m} \overset{I}{B} \overset{I}{A}\bigr) \neq \mathrm{tr}_q\bigl( \overset{I}{v}{^n}\overset{I}{A}\overset{I}{B}\bigr).\]
for all $m,n \in \mathbb{Z}$. This is due to the fact that if we choose the basepoint numbered $2$ on $[ba]$ at the bottom Figure \ref{changementPointBase}, we do not get a simple loop and the previous results do not apply. For instance, if $H = \bar U_q(\mathfrak{sl}_2)$ and $I = \mathcal{X}^+(2)$ (see Chapter \ref{chapitreUqSl2}), $\overset{I}{v}$ is just a scalar and a computation reveals that
\begin{align*}
\mathrm{tr}_q\bigl(\overset{I}{A}\overset{I}{B}\bigr) &= -qa_1a_2 - qb_1c_2 - q^{-1}c_1b_2 - q^{-1}d_1d_2\\
\mathrm{tr}_q\bigl(\overset{I}{B} \overset{I}{A}\bigr) &= - a_1a_2 + (1- q^{-2})a_1d_2 - q^2b_1c_2 - c_1b_2 + (1-q^{-2})d_1a_2 + (-1 + q^{-2} - q^{-4})d_1d_2.
\end{align*}
These two elements are not proportional since the monomials \eqref{PBWL10} are linearly independent. 
Hence such relations between traces in $\mathcal{L}_{g,0}(H)$ follow from the geometry of simple closed curves on $\Sigma_g^{\mathrm{o}}$ and are not simply algebraic coincidences.

\section{Graphical calculus when $H = \bar U_q(\mathfrak{sl}_2)$}\label{calculGraphiqueSl2}
\indent We take $H = \bar U_q = \bar U_q(\mathfrak{sl}_2)$. Moreover from now on, and otherwise indicated explicitly, {\em we assume that all the strands are colored  by the fundamental representation $\mathcal{X}^+(2)$}.\footnote{Recall from section \ref{braidedExtension} that the $R$-matrix belongs to an extension of $\bar U_q$ by a square root of $K$. In order to evaluate a crossing between two strands colored by $\mathcal{X}^+(2)$, we obviously define the action of $K^{1/2}$ on $\mathcal{X}^+(2) = \mathbb{C}v_0 \oplus \mathbb{C}v_1$ by $K^{1/2}v_0 = q^{1/2}v_0, \: K^{1/2}v_1 = q^{-1/2}v_1$.} In this case, the diagrammatic calculus (and \textit{a fortiori} the Wilson loop map) explained in the previous sections satisfies the Jones skein relation:
\begin{equation}\label{jones}
\begingroup%
  \makeatletter%
  \providecommand\color[2][]{%
    \errmessage{(Inkscape) Color is used for the text in Inkscape, but the package 'color.sty' is not loaded}%
    \renewcommand\color[2][]{}%
  }%
  \providecommand\transparent[1]{%
    \errmessage{(Inkscape) Transparency is used (non-zero) for the text in Inkscape, but the package 'transparent.sty' is not loaded}%
    \renewcommand\transparent[1]{}%
  }%
  \providecommand\rotatebox[2]{#2}%
  \newcommand*\fsize{\dimexpr\f@size pt\relax}%
  \newcommand*\lineheight[1]{\fontsize{\fsize}{#1\fsize}\selectfont}%
  \ifx\svgwidth\undefined%
    \setlength{\unitlength}{305.48680869bp}%
    \ifx\svgscale\undefined%
      \relax%
    \else%
      \setlength{\unitlength}{\unitlength * \real{\svgscale}}%
    \fi%
  \else%
    \setlength{\unitlength}{\svgwidth}%
  \fi%
  \global\let\svgwidth\undefined%
  \global\let\svgscale\undefined%
  \makeatother%
  \begin{picture}(1,0.14485981)%
    \lineheight{1}%
    \setlength\tabcolsep{0pt}%
    \put(0,0){\includegraphics[width=\unitlength,page=1]{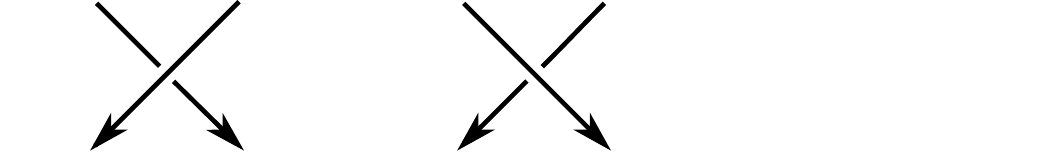}}%
    \put(-0.0013802,0.07491234){\color[rgb]{0,0,0}\makebox(0,0)[lt]{\lineheight{1.25}\smash{\begin{tabular}[t]{l}$q^{1/2}$\end{tabular}}}}%
    \put(0.26891972,0.07537867){\color[rgb]{0,0,0}\makebox(0,0)[lt]{\lineheight{1.25}\smash{\begin{tabular}[t]{l}$- \:\:\:\: q^{-1/2}$\end{tabular}}}}%
    \put(0,0){\includegraphics[width=\unitlength,page=2]{skein_jones.pdf}}%
    \put(0.58599193,0.07509718){\color[rgb]{0,0,0}\makebox(0,0)[lt]{\lineheight{1.25}\smash{\begin{tabular}[t]{l}$\:\:\: = \:\:\: (q - q^{-1})$\end{tabular}}}}%
    \put(0,0){\includegraphics[width=\unitlength,page=3]{skein_jones.pdf}}%
  \end{picture}%
\endgroup%

\end{equation}
This relation means that if $d_+, d_-$ and $d_{||}$ are three diagrams which are equal except in a small disk $D$, such that $d_+$ looks like a positive crossing in $D$, $d_-$ looks like a negative crossing in $D$ and $d_{||}$ looks like two vertical strands in $D$, then it holds $q^{1/2}d_+ - q^{-1/2}d_- = \hat q d_{||}$.

\smallskip

\indent An important fact is that $\mathcal{X}^+(2)$ is self-dual. Using this and the same kind of argument that in \cite[Lem 3.18]{KM}\footnote{Be aware that the algebra $\mathbf{U}_q$ of  \cite{KM} is not exactly $\bar U_q$ since their $K$ is a square root of our $K$, and that they choose $g = K^2$ for the pivotal element ($K$ in our notations). This choice of $g$ changes some signs between their formulas and ours, see for instance \cite[Th 4.3]{KM}.}, we are going to show that the Wilson loop map $W$ does not depend of the orientation of the link. 

\smallskip

\indent Let $\{v_0, v_1\}$ be the canonical basis of $\mathcal{X}^+(2)$ and $\{v^0, v^1\}$ be its dual basis, then

\begin{equation}\label{isoDual}
\begin{array}{rcrcl} 
D& : & \mathcal{X}^+(2)^* & \overset{\sim}{\rightarrow} & \mathcal{X}^+(2)\\
& & v^0 & \mapsto & -q v_1\\
& & v^1 & \mapsto & v_0
\end{array}
\end{equation}
is an isomorphism of $\bar U_q$-modules. We denote $e : \mathcal{X}^+(2)^{**} \overset{\sim}{\rightarrow} \mathcal{X}^+(2)$ the isomorphism with the bidual (see \eqref{identificationBidual}).

\begin{lemma}\label{lemmaDualite}
It holds $e \circ D^* = D$ and it follows that 
\begin{center}
\begingroup%
  \makeatletter%
  \providecommand\color[2][]{%
    \errmessage{(Inkscape) Color is used for the text in Inkscape, but the package 'color.sty' is not loaded}%
    \renewcommand\color[2][]{}%
  }%
  \providecommand\transparent[1]{%
    \errmessage{(Inkscape) Transparency is used (non-zero) for the text in Inkscape, but the package 'transparent.sty' is not loaded}%
    \renewcommand\transparent[1]{}%
  }%
  \providecommand\rotatebox[2]{#2}%
  \newcommand*\fsize{\dimexpr\f@size pt\relax}%
  \newcommand*\lineheight[1]{\fontsize{\fsize}{#1\fsize}\selectfont}%
  \ifx\svgwidth\undefined%
    \setlength{\unitlength}{501.21685591bp}%
    \ifx\svgscale\undefined%
      \relax%
    \else%
      \setlength{\unitlength}{\unitlength * \real{\svgscale}}%
    \fi%
  \else%
    \setlength{\unitlength}{\svgwidth}%
  \fi%
  \global\let\svgwidth\undefined%
  \global\let\svgscale\undefined%
  \makeatother%
  \begin{picture}(1,0.144249)%
    \lineheight{1}%
    \setlength\tabcolsep{0pt}%
    \put(0,0){\includegraphics[width=\unitlength,page=1]{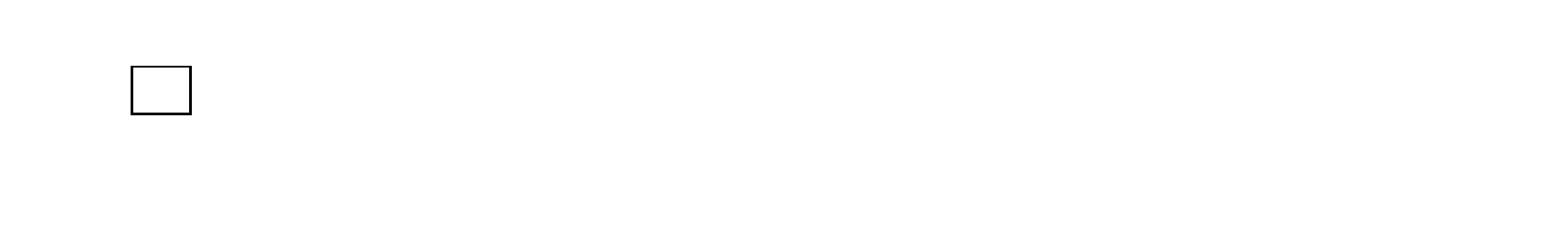}}%
    \put(0.09219551,0.07868311){\color[rgb]{0,0,0}\makebox(0,0)[lt]{\lineheight{1.25}\smash{\begin{tabular}[t]{l}$D$\end{tabular}}}}%
    \put(0,0){\includegraphics[width=\unitlength,page=2]{changeOrientation.pdf}}%
    \put(0.11690171,0.00268484){\color[rgb]{0,0,0}\makebox(0,0)[lt]{\lineheight{1.25}\smash{\begin{tabular}[t]{l}$\overset{\mathcal{X}^+(2)}{X}$\end{tabular}}}}%
    \put(0,0){\includegraphics[width=\unitlength,page=3]{changeOrientation.pdf}}%
    \put(0.3113426,0.00171056){\color[rgb]{0,0,0}\makebox(0,0)[lt]{\lineheight{1.25}\smash{\begin{tabular}[t]{l}$\overset{\mathcal{X}^+(2)}{X}$\end{tabular}}}}%
    \put(0,0){\includegraphics[width=\unitlength,page=4]{changeOrientation.pdf}}%
    \put(0.15915328,0.0506542){\color[rgb]{0,0,0}\makebox(0,0)[lt]{\lineheight{1.25}\smash{\begin{tabular}[t]{l}$=$\end{tabular}}}}%
    \put(0,0){\includegraphics[width=\unitlength,page=5]{changeOrientation.pdf}}%
    \put(0.21599266,0.07915929){\color[rgb]{0,0,0}\makebox(0,0)[lt]{\lineheight{1.25}\smash{\begin{tabular}[t]{l}$D$\end{tabular}}}}%
    \put(0,0){\includegraphics[width=\unitlength,page=6]{changeOrientation.pdf}}%
    \put(0.52291347,0.06621919){\color[rgb]{0,0,0}\makebox(0,0)[lt]{\lineheight{1.25}\smash{\begin{tabular}[t]{l}$=$\end{tabular}}}}%
    \put(0,0){\includegraphics[width=\unitlength,page=7]{changeOrientation.pdf}}%
    \put(0.47263923,0.0796064){\color[rgb]{0,0,0}\makebox(0,0)[lt]{\lineheight{1.25}\smash{\begin{tabular}[t]{l}$D$\end{tabular}}}}%
    \put(0.57061409,0.07922912){\color[rgb]{0,0,0}\makebox(0,0)[lt]{\lineheight{1.25}\smash{\begin{tabular}[t]{l}$D$\end{tabular}}}}%
    \put(0,0){\includegraphics[width=\unitlength,page=8]{changeOrientation.pdf}}%
    \put(0.85946508,0.08005134){\color[rgb]{0,0,0}\makebox(0,0)[lt]{\lineheight{1.25}\smash{\begin{tabular}[t]{l}$=$\end{tabular}}}}%
    \put(0.81026169,0.06631904){\color[rgb]{0,0,0}\makebox(0,0)[lt]{\lineheight{1.25}\smash{\begin{tabular}[t]{l}$D$\end{tabular}}}}%
    \put(0.90712936,0.06631904){\color[rgb]{0,0,0}\makebox(0,0)[lt]{\lineheight{1.25}\smash{\begin{tabular}[t]{l}$D$\end{tabular}}}}%
  \end{picture}%
\endgroup%

\end{center}
\end{lemma}
\begin{proof}
One checks easily that 
\[ D^*(v^0) = \langle ?, v_1 \rangle, \:\:\:\: D^*(v^1) = -q \langle ?, v_0 \rangle. \]
Hence $e \circ D^*(v^0) = e(\langle ?, v_1 \rangle) = K^{p-1}v_1 = -qv_1 = D(v^0)$ and similarly for $v^1$. Using this and \eqref{sensOppose}, we get:
\begin{center}
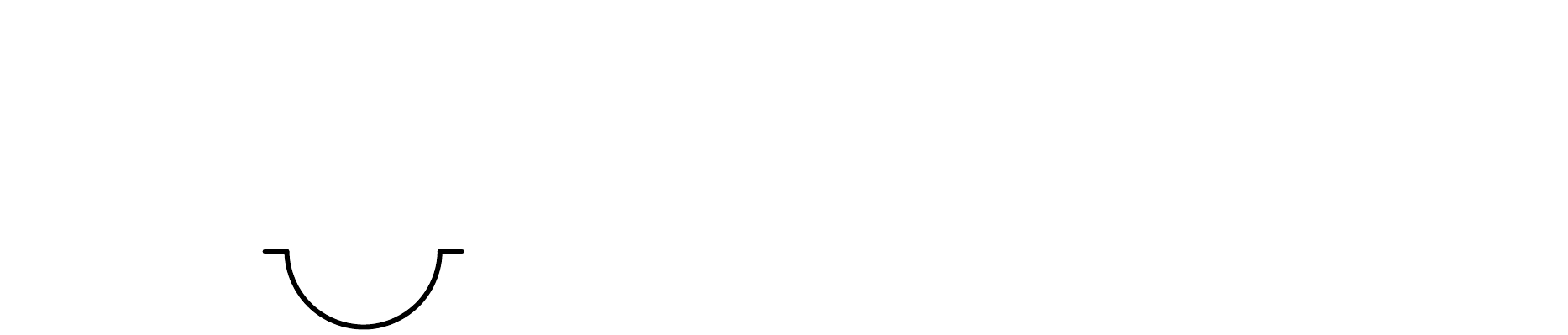
\end{center}
and the proofs of the others equalities are similar.
\end{proof}

\begin{remark}
If we choose $K$ instead of $K^{p+1}$ for the pivotal element $g$ of $\bar U_q$, then a minus sign appears in the equalities of Lemma \ref{lemmaDualite}. In this case, the value of $W(L)$ may depend of the orientation of $L$ up to a sign.
\finEx
\end{remark}

\begin{corollary}\label{indepOrientation}
When $H = \bar U_q$, $g = K^{p+1}$ and all the strands of $L$ are colored by $\mathcal{X}^+(2)$, $W(L)$ does not depend of the orientation of the strands of $L$.
\end{corollary}
\begin{proof}
We can assume that $L$ simply contains one connected component. Apply $W$ to $L$. Then in the diagram representing $W(L)$, introduce two coupons $D$ and $D^{-1}$ in the strand, according to the orientation of $L$. Between these coupons, the orientation of the strand is reversed. Using the previous lemma, we move $D$ along $L$: at each time it passes through a cup, a cap or a crossing, it changes the orientation. At the end, the coupon $D$ arrives on the other side of $D^{-1}$ and they both collapse, leaving $W(L)$ with the opposite orientation. See the figure above, where we have introduced $D$ and $D^{-1}$ near a cap in $W(L)$.
\begin{center}
\begingroup%
  \makeatletter%
  \providecommand\color[2][]{%
    \errmessage{(Inkscape) Color is used for the text in Inkscape, but the package 'color.sty' is not loaded}%
    \renewcommand\color[2][]{}%
  }%
  \providecommand\transparent[1]{%
    \errmessage{(Inkscape) Transparency is used (non-zero) for the text in Inkscape, but the package 'transparent.sty' is not loaded}%
    \renewcommand\transparent[1]{}%
  }%
  \providecommand\rotatebox[2]{#2}%
  \newcommand*\fsize{\dimexpr\f@size pt\relax}%
  \newcommand*\lineheight[1]{\fontsize{\fsize}{#1\fsize}\selectfont}%
  \ifx\svgwidth\undefined%
    \setlength{\unitlength}{275.03233925bp}%
    \ifx\svgscale\undefined%
      \relax%
    \else%
      \setlength{\unitlength}{\unitlength * \real{\svgscale}}%
    \fi%
  \else%
    \setlength{\unitlength}{\svgwidth}%
  \fi%
  \global\let\svgwidth\undefined%
  \global\let\svgscale\undefined%
  \makeatother%
  \begin{picture}(1,0.40463611)%
    \lineheight{1}%
    \setlength\tabcolsep{0pt}%
    \put(0.49503619,0.2068748){\color[rgb]{0,0,0}\makebox(0,0)[lt]{\lineheight{1.25}\smash{\begin{tabular}[t]{l}=\end{tabular}}}}%
    \put(0,0){\includegraphics[width=\unitlength,page=1]{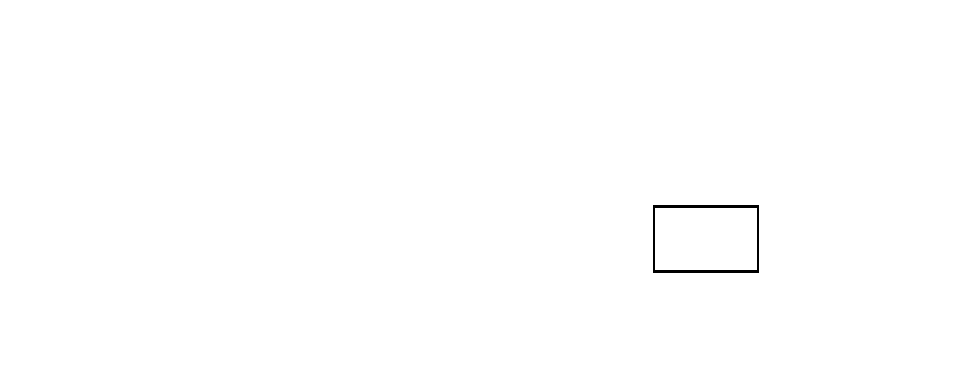}}%
    \put(0.72399504,0.13802084){\color[rgb]{0,0,0}\makebox(0,0)[lt]{\lineheight{1.25}\smash{\begin{tabular}[t]{l}$D$\end{tabular}}}}%
    \put(0,0){\includegraphics[width=\unitlength,page=2]{couponDD.pdf}}%
    \put(0.70221532,0.28499843){\color[rgb]{0,0,0}\makebox(0,0)[lt]{\lineheight{1.25}\smash{\begin{tabular}[t]{l}$D^{-1}$\end{tabular}}}}%
    \put(0,0){\includegraphics[width=\unitlength,page=3]{couponDD.pdf}}%
    \put(0.81214976,0.20660264){\color[rgb]{0,0,0}\makebox(0,0)[lt]{\lineheight{1.25}\smash{\begin{tabular}[t]{l}=\end{tabular}}}}%
    \put(0,0){\includegraphics[width=\unitlength,page=4]{couponDD.pdf}}%
    \put(0.37180919,0.14085811){\color[rgb]{0,0,0}\makebox(0,0)[lt]{\lineheight{1.25}\smash{\begin{tabular}[t]{l}$D^{-1}$\end{tabular}}}}%
    \put(0,0){\includegraphics[width=\unitlength,page=5]{couponDD.pdf}}%
    \put(0.39254115,0.29127463){\color[rgb]{0,0,0}\makebox(0,0)[lt]{\lineheight{1.25}\smash{\begin{tabular}[t]{l}$D$\end{tabular}}}}%
    \put(0,0){\includegraphics[width=\unitlength,page=6]{couponDD.pdf}}%
    \put(0.16361577,0.20793758){\color[rgb]{0,0,0}\makebox(0,0)[lt]{\lineheight{1.25}\smash{\begin{tabular}[t]{l}=\end{tabular}}}}%
    \put(0,0){\includegraphics[width=\unitlength,page=7]{couponDD.pdf}}%
  \end{picture}%
\endgroup%

\end{center}
\end{proof}

Hence, we can define the value of $W$ on a non-oriented link: this is just the value of $W$ on $L$ with an arbitrary orientation (and colored by $\mathcal{X}^+(2)$). Equivalently, we can define directly a non-oriented diagrammatic calculus for $\mathcal{L}_{g,n}(\bar U_q)$, with unoriented cups, caps, crossings and handles:
\begin{center}
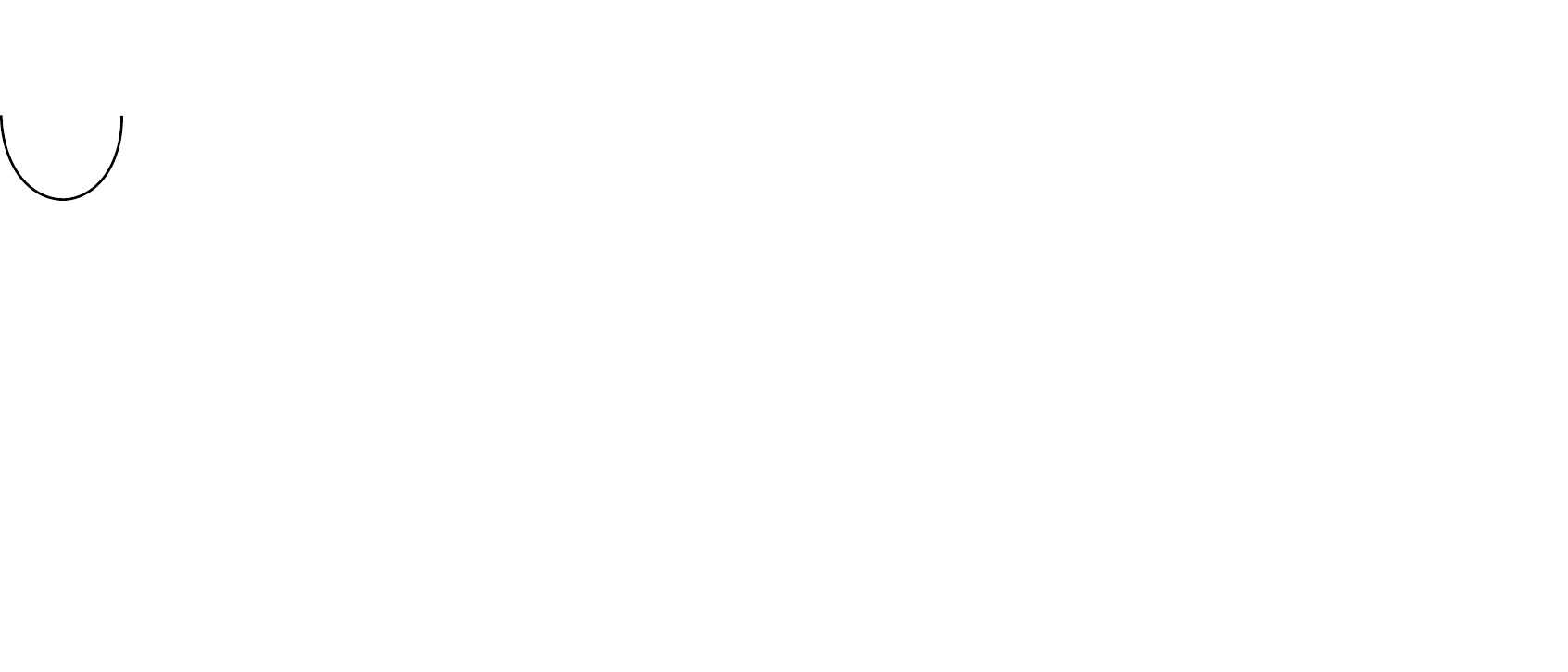
\end{center}
Note that until now we simply write $X$ for the matrices $\overset{\mathcal{X}^+(2)}{X}$ labelling the handles. As in section \ref{sectionL01Uq}, we denote by $a,b,c,d$ the coefficients of the matrix $M$, which generate $\mathcal{L}_{0,1}(\bar U_q)$:
\[ M = 
\left(
\begin{array}{cc}
a & b\\
c & d
\end{array}
\right) \in \mathrm{Mat}_2\!\left(\mathcal{L}_{0,1}(\bar U_q)\right).
 \]
The explicit values of the unoriented graphical elements are:
\[ \cup(1) = v_1 \otimes v_0 - q \, v_0 \otimes v_1, \:\:\:\:\:\:\:\:
\begin{array}{ll}
\cap(v_0 \otimes v_0) = 0, & \cap(v_0 \otimes v_1) = 1,\\
\cap(v_1 \otimes v_0) = -q^{-1}, & \cap(v_1 \otimes v_1) = 0,
\end{array}
\]
\begin{center}
\begingroup%
  \makeatletter%
  \providecommand\color[2][]{%
    \errmessage{(Inkscape) Color is used for the text in Inkscape, but the package 'color.sty' is not loaded}%
    \renewcommand\color[2][]{}%
  }%
  \providecommand\transparent[1]{%
    \errmessage{(Inkscape) Transparency is used (non-zero) for the text in Inkscape, but the package 'transparent.sty' is not loaded}%
    \renewcommand\transparent[1]{}%
  }%
  \providecommand\rotatebox[2]{#2}%
  \newcommand*\fsize{\dimexpr\f@size pt\relax}%
  \newcommand*\lineheight[1]{\fontsize{\fsize}{#1\fsize}\selectfont}%
  \ifx\svgwidth\undefined%
    \setlength{\unitlength}{1127.756506bp}%
    \ifx\svgscale\undefined%
      \relax%
    \else%
      \setlength{\unitlength}{\unitlength * \real{\svgscale}}%
    \fi%
  \else%
    \setlength{\unitlength}{\svgwidth}%
  \fi%
  \global\let\svgwidth\undefined%
  \global\let\svgscale\undefined%
  \makeatother%
  \begin{picture}(1,0.09405225)%
    \lineheight{1}%
    \setlength\tabcolsep{0pt}%
    \put(0,0){\includegraphics[width=\unitlength,page=1]{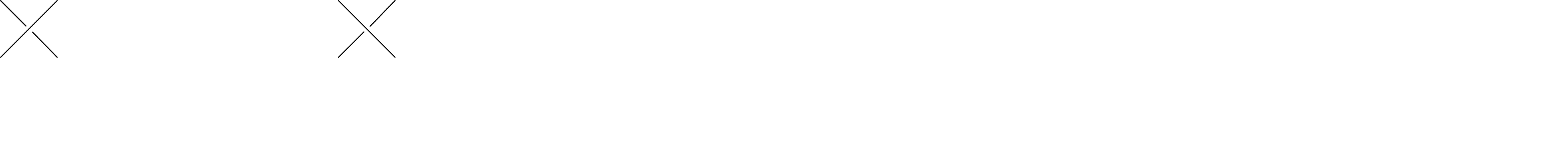}}%
    \put(0.0425307,0.07252502){\color[rgb]{0,0,0}\makebox(0,0)[lt]{\lineheight{1.25}\smash{\begin{tabular}[t]{l}$=q^{-1/2}\left(\begin{array}{cccc} q & 0 & 0 & 0 \\0 & 0 & 1 & 0\\0 & 1 & \hat q & 0\\0 & 0 & 0 & q\end{array}\right),$\end{tabular}}}}%
    \put(0.25903948,0.0723261){\color[rgb]{0,0,0}\makebox(0,0)[lt]{\lineheight{1.25}\smash{\begin{tabular}[t]{l}$=q^{1/2}\left(\begin{array}{cccc} q^{-1} & 0 & 0 & 0 \\0 &  -\hat q & 1 & 0\\0 & 1 & 0 & 0\\0 & 0 & 0 & q^{-1}\end{array}\right),$\end{tabular}}}}%
    \put(0.12600817,0.00129968){\color[rgb]{0,0,0}\makebox(0,0)[lt]{\lineheight{1.25}\smash{\begin{tabular}[t]{l}$M$\end{tabular}}}}%
    \put(0,0){\includegraphics[width=\unitlength,page=2]{valeurRmatrice.pdf}}%
    \put(0.14027338,0.02435856){\color[rgb]{0,0,0}\makebox(0,0)[lt]{\lineheight{1.25}\smash{\begin{tabular}[t]{l}$= -qa \, v_0 \otimes v_1  + b \, v_0 \otimes v_0 - qc \, v_1 \otimes v_1 + d \, v_1 \otimes v_0$.\end{tabular}}}}%
  \end{picture}%
\endgroup%

\end{center}
A simple computation yields the Kauffman skein relation:
\begin{equation}\label{relationKauffman}
\begingroup%
  \makeatletter%
  \providecommand\color[2][]{%
    \errmessage{(Inkscape) Color is used for the text in Inkscape, but the package 'color.sty' is not loaded}%
    \renewcommand\color[2][]{}%
  }%
  \providecommand\transparent[1]{%
    \errmessage{(Inkscape) Transparency is used (non-zero) for the text in Inkscape, but the package 'transparent.sty' is not loaded}%
    \renewcommand\transparent[1]{}%
  }%
  \providecommand\rotatebox[2]{#2}%
  \newcommand*\fsize{\dimexpr\f@size pt\relax}%
  \newcommand*\lineheight[1]{\fontsize{\fsize}{#1\fsize}\selectfont}%
  \ifx\svgwidth\undefined%
    \setlength{\unitlength}{252.13898009bp}%
    \ifx\svgscale\undefined%
      \relax%
    \else%
      \setlength{\unitlength}{\unitlength * \real{\svgscale}}%
    \fi%
  \else%
    \setlength{\unitlength}{\svgwidth}%
  \fi%
  \global\let\svgwidth\undefined%
  \global\let\svgscale\undefined%
  \makeatother%
  \begin{picture}(1,0.16834092)%
    \lineheight{1}%
    \setlength\tabcolsep{0pt}%
    \put(0,0){\includegraphics[width=\unitlength,page=1]{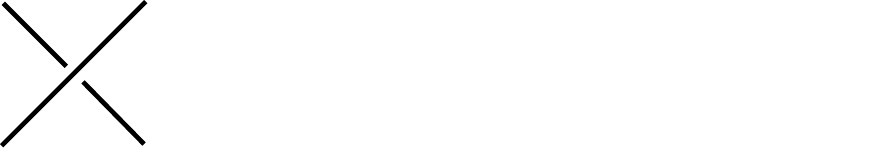}}%
    \put(0.22206066,0.08447801){\color[rgb]{0,0,0}\makebox(0,0)[lt]{\lineheight{1.25}\smash{\begin{tabular}[t]{l}$= \:\:\:q^{1/2}$\end{tabular}}}}%
    \put(0.63281951,0.08415885){\color[rgb]{0,0,0}\makebox(0,0)[lt]{\lineheight{1.25}\smash{\begin{tabular}[t]{l}$+ \:\:\:\: q^{-1/2}$\end{tabular}}}}%
    \put(0,0){\includegraphics[width=\unitlength,page=2]{skein_Kauffman.pdf}}%
  \end{picture}%
\endgroup%

\end{equation}
Alternatively, one can derive this relation from the Jones skein relation \eqref{jones} and the independence of the orientation. Relation \eqref{relationKauffman} means that if $d_+, d_{||}$ and $d_{=}$ are three diagrams which are equal except in a small disk $D$, such that $d_+$ looks like a positive crossing in $D$, $d_{||}$ looks like two vertical strands in $D$ and $d_{=}$ looks like two horizontal strands in $D$, it holds $d_+ = q^{1/2}d_{||} + q^{-1/2}d_{=}$. 
It is also useful to record the relation for negative crossings, for twists and for contractible circles:
\begin{center}
\begingroup%
  \makeatletter%
  \providecommand\color[2][]{%
    \errmessage{(Inkscape) Color is used for the text in Inkscape, but the package 'color.sty' is not loaded}%
    \renewcommand\color[2][]{}%
  }%
  \providecommand\transparent[1]{%
    \errmessage{(Inkscape) Transparency is used (non-zero) for the text in Inkscape, but the package 'transparent.sty' is not loaded}%
    \renewcommand\transparent[1]{}%
  }%
  \providecommand\rotatebox[2]{#2}%
  \newcommand*\fsize{\dimexpr\f@size pt\relax}%
  \newcommand*\lineheight[1]{\fontsize{\fsize}{#1\fsize}\selectfont}%
  \ifx\svgwidth\undefined%
    \setlength{\unitlength}{570.46577488bp}%
    \ifx\svgscale\undefined%
      \relax%
    \else%
      \setlength{\unitlength}{\unitlength * \real{\svgscale}}%
    \fi%
  \else%
    \setlength{\unitlength}{\svgwidth}%
  \fi%
  \global\let\svgwidth\undefined%
  \global\let\svgscale\undefined%
  \makeatother%
  \begin{picture}(1,0.21502647)%
    \lineheight{1}%
    \setlength\tabcolsep{0pt}%
    \put(0,0){\includegraphics[width=\unitlength,page=1]{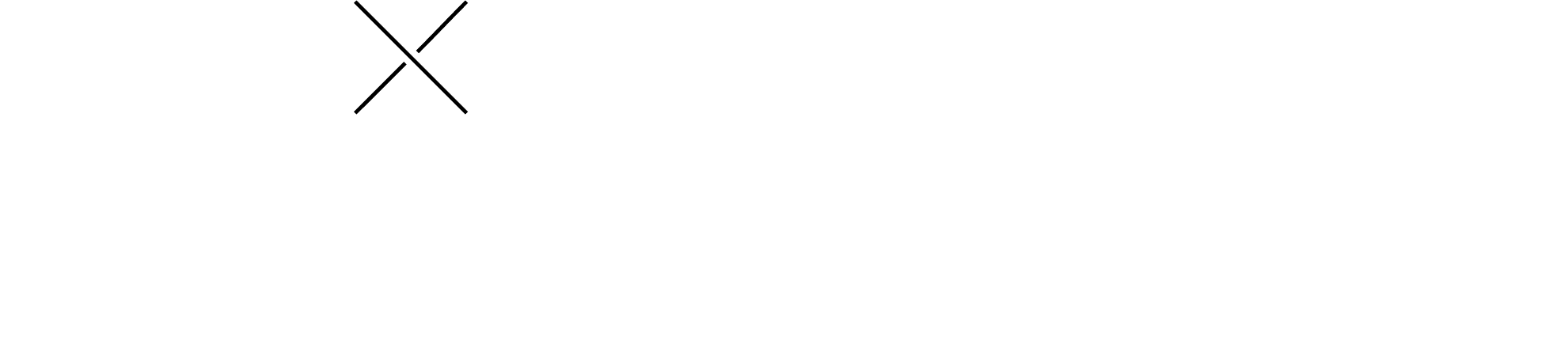}}%
    \put(0.32428341,0.17854343){\color[rgb]{0,0,0}\makebox(0,0)[lt]{\lineheight{1.25}\smash{\begin{tabular}[t]{l}$= \:\:\:q^{-1/2}$\end{tabular}}}}%
    \put(0.50583388,0.17840237){\color[rgb]{0,0,0}\makebox(0,0)[lt]{\lineheight{1.25}\smash{\begin{tabular}[t]{l}$+ \:\:\:\: q^{1/2}$\end{tabular}}}}%
    \put(0,0){\includegraphics[width=\unitlength,page=2]{skein_Kauffman_inverse.pdf}}%
    \put(0.06498782,0.04475902){\color[rgb]{0,0,0}\makebox(0,0)[lt]{\lineheight{1.25}\smash{\begin{tabular}[t]{l}=\end{tabular}}}}%
    \put(0.17080381,0.04525543){\color[rgb]{0,0,0}\makebox(0,0)[lt]{\lineheight{1.25}\smash{\begin{tabular}[t]{l}$= \:\:\:-q^{3/2}$\end{tabular}}}}%
    \put(0,0){\includegraphics[width=\unitlength,page=3]{skein_Kauffman_inverse.pdf}}%
    \put(0.4473898,0.04896038){\color[rgb]{0,0,0}\makebox(0,0)[lt]{\lineheight{1.25}\smash{\begin{tabular}[t]{l}=\end{tabular}}}}%
    \put(0.55199493,0.04765269){\color[rgb]{0,0,0}\makebox(0,0)[lt]{\lineheight{1.25}\smash{\begin{tabular}[t]{l}$= \:\:\:-q^{-3/2}$\end{tabular}}}}%
    \put(0,0){\includegraphics[width=\unitlength,page=4]{skein_Kauffman_inverse.pdf}}%
    \put(0.81486461,0.04655856){\color[rgb]{0,0,0}\makebox(0,0)[lt]{\lineheight{1.25}\smash{\begin{tabular}[t]{l}$= \:- [2] \: \varnothing$\end{tabular}}}}%
  \end{picture}%
\endgroup%

\end{center}
where $\varnothing$ is the empty diagram. The Kauffman skein relation allows us to resolve all the crossings in the diagrammatic calculus and is very useful to derive identities in $\mathcal{L}_{g,n}(\bar U_q)$, like the following proposition.

\begin{proposition}
It holds:
\[ M^2 + q^{-1} W_{\! M} M + q^{-2}\mathbb{I}_2 = 0 \]
where we recall that $W_M = \mathrm{tr}_q(M) = \mathrm{tr}(\overset{\mathcal{X}^+(2)}{K^{p+1}}M)$ and $\mathbb{I}_2$ is the $2 \times 2$ identity matrix. We call this relation the quantum Cayley-Hamilton identity.
\end{proposition}
\begin{proof}
Using \eqref{dessinMatriceInverse} and the Kauffman skein relation, we get:
\begin{center}
\begingroup%
  \makeatletter%
  \providecommand\color[2][]{%
    \errmessage{(Inkscape) Color is used for the text in Inkscape, but the package 'color.sty' is not loaded}%
    \renewcommand\color[2][]{}%
  }%
  \providecommand\transparent[1]{%
    \errmessage{(Inkscape) Transparency is used (non-zero) for the text in Inkscape, but the package 'transparent.sty' is not loaded}%
    \renewcommand\transparent[1]{}%
  }%
  \providecommand\rotatebox[2]{#2}%
  \newcommand*\fsize{\dimexpr\f@size pt\relax}%
  \newcommand*\lineheight[1]{\fontsize{\fsize}{#1\fsize}\selectfont}%
  \ifx\svgwidth\undefined%
    \setlength{\unitlength}{522.20162539bp}%
    \ifx\svgscale\undefined%
      \relax%
    \else%
      \setlength{\unitlength}{\unitlength * \real{\svgscale}}%
    \fi%
  \else%
    \setlength{\unitlength}{\svgwidth}%
  \fi%
  \global\let\svgwidth\undefined%
  \global\let\svgscale\undefined%
  \makeatother%
  \begin{picture}(1,0.14537564)%
    \lineheight{1}%
    \setlength\tabcolsep{0pt}%
    \put(0.10897372,0.00311455){\color[rgb]{0,0,0}\makebox(0,0)[lt]{\lineheight{1.25}\smash{\begin{tabular}[t]{l}$M^{-1}$\end{tabular}}}}%
    \put(0,0){\includegraphics[width=\unitlength,page=1]{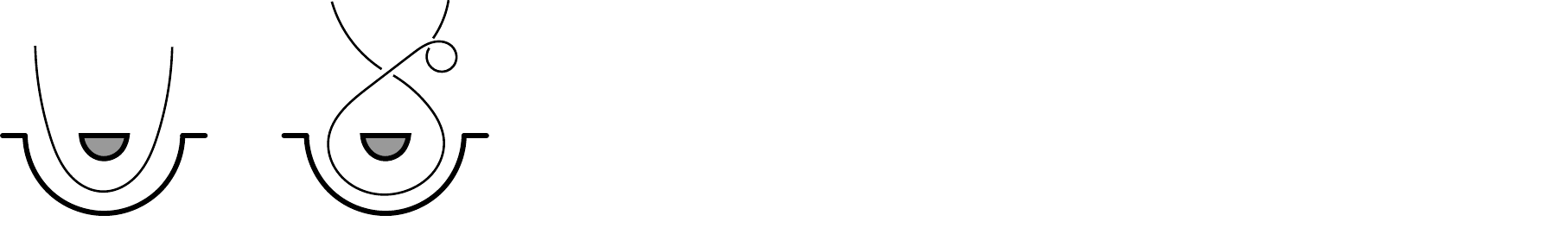}}%
    \put(0.14810555,0.05241269){\color[rgb]{0,0,0}\makebox(0,0)[lt]{\lineheight{1.25}\smash{\begin{tabular}[t]{l}=\end{tabular}}}}%
    \put(0,0){\includegraphics[width=\unitlength,page=2]{calcul_inverse_M.pdf}}%
    \put(0.3616491,0.05383037){\color[rgb]{0,0,0}\makebox(0,0)[lt]{\lineheight{1.25}\smash{\begin{tabular}[t]{l}$-q^{3/2}$\end{tabular}}}}%
    \put(0.33135091,0.05244323){\color[rgb]{0,0,0}\makebox(0,0)[lt]{\lineheight{1.25}\smash{\begin{tabular}[t]{l}=\end{tabular}}}}%
    \put(0,0){\includegraphics[width=\unitlength,page=3]{calcul_inverse_M.pdf}}%
    \put(0.60580764,0.05383049){\color[rgb]{0,0,0}\makebox(0,0)[lt]{\lineheight{1.25}\smash{\begin{tabular}[t]{l}$-q^{2}$\end{tabular}}}}%
    \put(0.57550941,0.05244335){\color[rgb]{0,0,0}\makebox(0,0)[lt]{\lineheight{1.25}\smash{\begin{tabular}[t]{l}=\end{tabular}}}}%
    \put(0,0){\includegraphics[width=\unitlength,page=4]{calcul_inverse_M.pdf}}%
    \put(0.80087236,0.05212266){\color[rgb]{0,0,0}\makebox(0,0)[lt]{\lineheight{1.25}\smash{\begin{tabular}[t]{l}$- \:\:\: q$\end{tabular}}}}%
    \put(0,0){\includegraphics[width=\unitlength,page=5]{calcul_inverse_M.pdf}}%
    \put(0.2884384,0.00629356){\color[rgb]{0,0,0}\makebox(0,0)[lt]{\lineheight{1.25}\smash{\begin{tabular}[t]{l}$M$\end{tabular}}}}%
    \put(0.53469357,0.00747042){\color[rgb]{0,0,0}\makebox(0,0)[lt]{\lineheight{1.25}\smash{\begin{tabular}[t]{l}$M$\end{tabular}}}}%
    \put(0.75657764,0.00826518){\color[rgb]{0,0,0}\makebox(0,0)[lt]{\lineheight{1.25}\smash{\begin{tabular}[t]{l}$M$\end{tabular}}}}%
    \put(0.96845285,0.00861422){\color[rgb]{0,0,0}\makebox(0,0)[lt]{\lineheight{1.25}\smash{\begin{tabular}[t]{l}$M$\end{tabular}}}}%
  \end{picture}%
\endgroup%

\end{center}
which means that $M^{-1} = -q^2 M - q W_{\!M} \mathbb{I}_2$. The desired identity immediately follows.
\end{proof}

\noindent The proof of the previous proposition also shows that 
$ M^{-1} = 
\left(
\begin{array}{cc}
d & -q^2 b\\
-q^2 c & q^2 a -q \hat q d
\end{array}
\right).
$
Moreover, it is easy to show by induction that
\begin{equation}\label{puissanceM}
 M^n = (-1)^{n+1}q^{-n+1}R_n(W_{\! M})M + (-1)^{n+1}q^{-n}R_{n-1}(W_{\! M})\mathbb{I}_2 
\end{equation}
where the polynomials $R_n$ are defined by $R_0(X) = 0, R_1(X) = 1$ and $R_{n+1}(X) = X R_n(X) - R_{n-1}(X)$ for $n \geq 1$ (Chebychev polynomials of the second kind).

\smallskip

We now focus on $\Sigma_{1,0}^{\mathrm{o}}$, and we want to prove the following proposition. In fact, it follows from the property that $\mathcal{S}_q(\Sigma_1)$ is generated by the isotopy classes $[b], [a]$. However, we want to show it on the $\mathcal{L}_{1,0}(\bar U_q)$ side, with the graphical calculus. Moreover, enhancing slightly the proof of Lemma \ref{diagrammeReduit}, we can show a more general result: if $X$ is any product of $B^{\pm 1}, A^{\pm 1}$, then $W_{\! X} \in \langle W_{\! A}, W_{\! B} \rangle$.

\begin{proposition}\label{propWAWB}
For any framed link $L \subset \Sigma_{1,0}^{\mathrm{o}} \times [0, 1]$ (whose all strands are colored by $\mathcal{X}^+(2)$), we have $W(L) \in \mathbb{C}\langle W_{\! A}, W_{\! B} \rangle$. In other words, every $W(L)$ can be written as a (non-commutative) polynomial in $W_{\! A}$ and $W_{\! B}$.
\end{proposition}

We need two lemmas.

\begin{lemma}\label{lemmaWAWB}
We have:
\begin{align*}
&1. \:\:\: AW_{\! A} = W_{\! A}A, \:\:\: BW_{\! B} = W_{\! B}B, \:\:\: BAW_{\! BA} = W_{\! BA}BA, \\
&2. \:\:\: A W_{\! B} = q^{-1} W_{\! B} A - q \hat q BA, \\
&3. \:\:\: B W_{\! A} = q W_{\! A} B + q^2 \hat q BA, \\
&4. \:\:\: BA W_{\! A} = q^{-1} W_{\! A} BA - q^{-2} \hat q B, \\
&5. \:\:\: BA W_{\! B} = q W_{\! B} BA + q^{-1} \hat q A, \\
&6. \:\:\: W_{\! BA} = q^{-2}\hat q^{-1} W_{\! B} W_{\! A} - q^{-1}\hat q^{-1} W_{\! A} W_{\! B} \:\: \text{(in particular, }W_{\! BA} \in \mathbb{C}\langle W_{\! A}, W_{\! B} \rangle).
\end{align*}
It follows that, if $P \in \mathbb{C}\langle x_1, x_2 \rangle$ is a (non-commutative) polynomial, then there exist $Q, R, S \in \mathbb{C}\langle x_1, x_2 \rangle$ such that:
\[ B P(W_{\! A}, W_{\!B}) = Q(W_{\! A}, W_{\!B})A + R(W_{\! A}, W_{\!B})B + S(W_{\! A}, W_{\!B})BA. \]
This is also true if we replace $B$ by $A$ or $BA$ but we will not need it.
\end{lemma}
\begin{proof}
Since $W_{\! M}$ is central in $\mathcal{L}_{0,1}(H)$, we have $MW_{\! M} = W_{\! M}M$. But $A, B$ and\footnote{Recall that in this section, all is evaluated in $\mathcal{X}^+(2)$ and thus $v$ is identified with the scalar $-q^{-3/2}$, by \eqref{valueVRep}.} $v^{-1}BA = -q^{3/2}BA$ satisfy the fusion relation of $\mathcal{L}_{0,1}(H)$, thus we can apply the morphisms $j_A, j_B, j_{v^{-1}BA}$ defined by $j_A(M) = A, j_B(M) = B, j_{v^{-1}BA}(M) = v^{-1}BA$ and we get the three equalities of 1. Next, using the exchange relation between $A$ and $B$ \eqref{dessinEchangeL10}, we prove relations 2 and 3 diagramatically, see Figures \ref{preuveRel2} and \ref{preuveRel3}. Relations 4 and 5 are immediate consequences. For instance:
\[ BA W_{\! A} = B W_{\! A} A = q W_{\! A} BA + q^2 \hat q BA^2 = q W_{\! A} BA - q \hat q BAW_{\! A} - \hat q B = q W_{\! A} BA - q^2 BAW_{\! A} + BAW_{\! A} - \hat q B. \]
\begin{figure}[h]
\centering
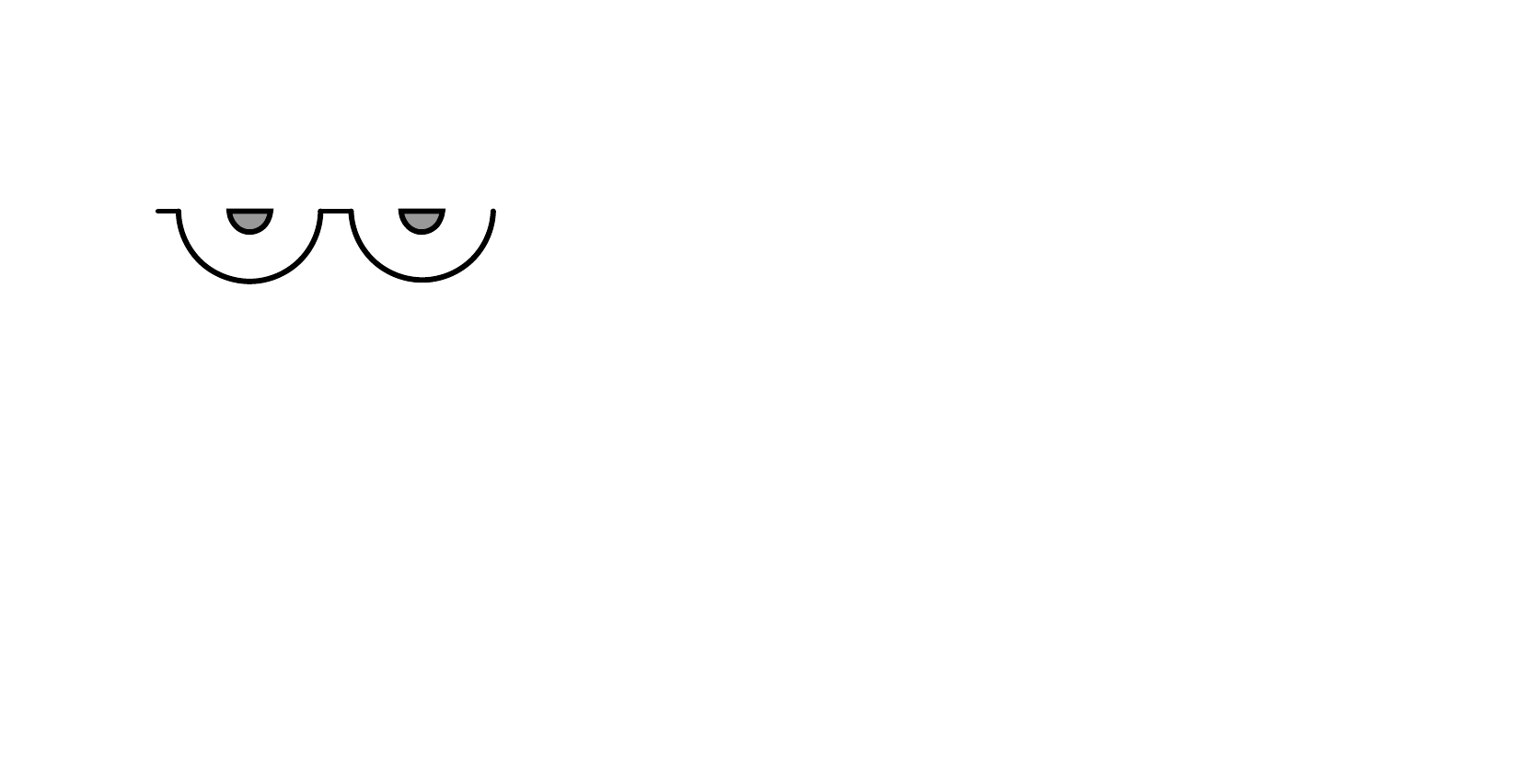
\caption{Proof of relation 2.}
\label{preuveRel2}
\end{figure}
\begin{figure}[h]
\centering
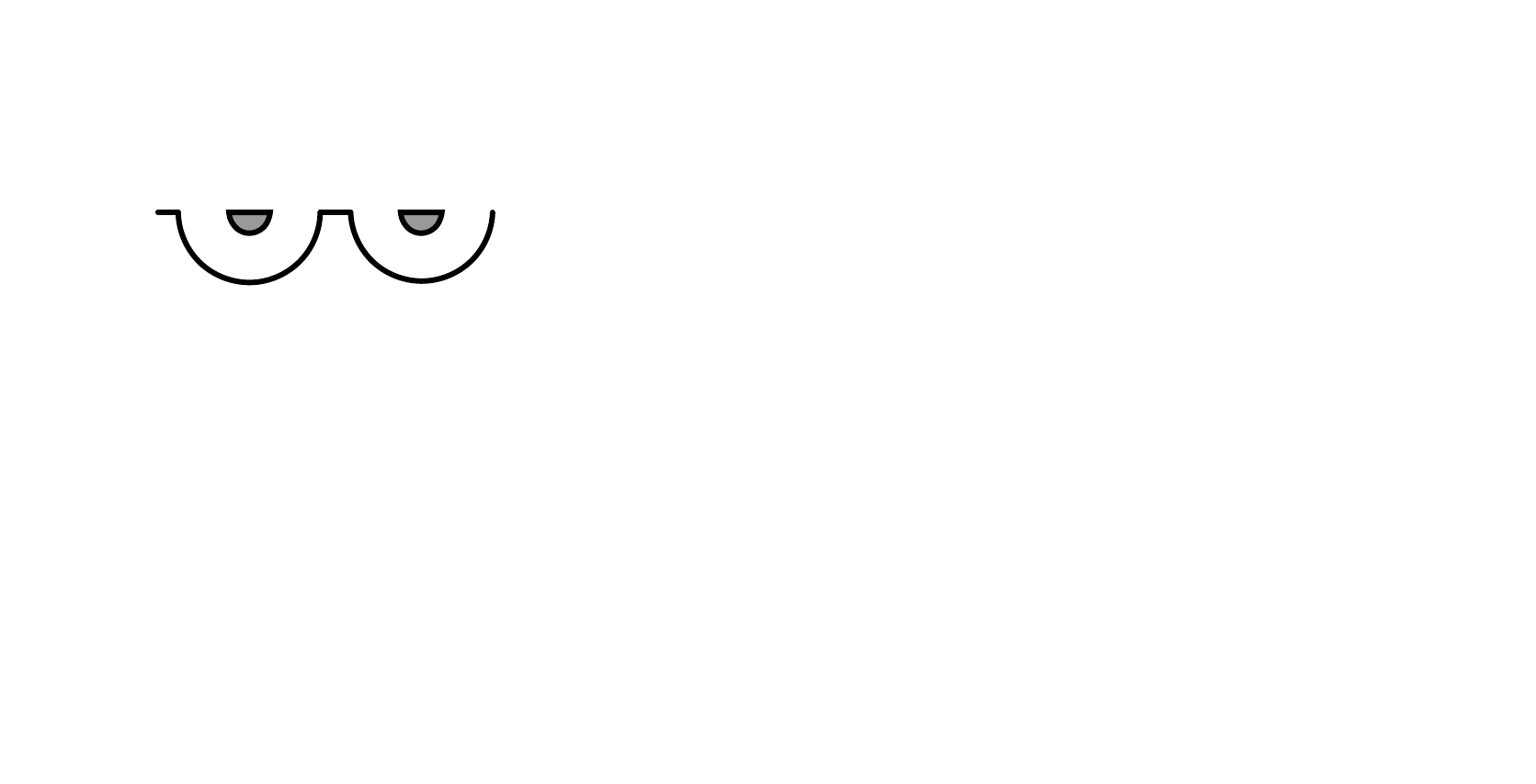
\caption{Proof of relation 3.}
\label{preuveRel3}
\end{figure}
We used the quantum Cayley-Hamilton identity. It follows that $0 = q W_{\! A} BA - q^2 BAW_{\! A} - \hat q B$, as desired. To show relation 6, simply apply $\mathrm{tr}_q$ to relation 2. 
For the second part of the lemma, we observe that it is a consequence of a more general fact, namely: if $P \in \mathbb{C}\langle x_1, x_2 \rangle$ is a (non-commutative) polynomial, then there exist $Q_i,R_i,S_i \in \mathbb{C}\langle x_1, x_2 \rangle$, $i = 1,2,3$, such that:
\begin{align*}
A P(W_{\! A}, W_{\!B}) &= Q_1(W_{\! A}, W_{\!B})A + R_1(W_{\! A}, W_{\!B})B + S_1(W_{\! A}, W_{\!B})BA, \\
B P(W_{\! A}, W_{\!B}) &= Q_2(W_{\! A}, W_{\!B})A + R_2(W_{\! A}, W_{\!B})B + S_2(W_{\! A}, W_{\!B})BA, \\
BA P(W_{\! A}, W_{\!B}) &= Q_3(W_{\! A}, W_{\!B})A + R_3(W_{\! A}, W_{\!B})B + S_3(W_{\! A}, W_{\!B})BA.
\end{align*}
Indeed, we can assume that $P$ is a monomial and show this set of three equalities by induction on the length of $P$ (for instance, $P = x_1x_2x_1^2$ has length $4$) thanks to the previous commutation relations.
\end{proof}

\begin{lemma}\label{diagrammeReduit}
For any a framed link $L \subset \Sigma_{1,0}^{\mathrm{o}} \times [0, 1]$ (whose all strands are colored by $\mathcal{X}^+(2)$), $W(L)$ is a linear combination of elements of the form
\begin{center}
\begingroup%
  \makeatletter%
  \providecommand\color[2][]{%
    \errmessage{(Inkscape) Color is used for the text in Inkscape, but the package 'color.sty' is not loaded}%
    \renewcommand\color[2][]{}%
  }%
  \providecommand\transparent[1]{%
    \errmessage{(Inkscape) Transparency is used (non-zero) for the text in Inkscape, but the package 'transparent.sty' is not loaded}%
    \renewcommand\transparent[1]{}%
  }%
  \providecommand\rotatebox[2]{#2}%
  \newcommand*\fsize{\dimexpr\f@size pt\relax}%
  \newcommand*\lineheight[1]{\fontsize{\fsize}{#1\fsize}\selectfont}%
  \ifx\svgwidth\undefined%
    \setlength{\unitlength}{486.61889561bp}%
    \ifx\svgscale\undefined%
      \relax%
    \else%
      \setlength{\unitlength}{\unitlength * \real{\svgscale}}%
    \fi%
  \else%
    \setlength{\unitlength}{\svgwidth}%
  \fi%
  \global\let\svgwidth\undefined%
  \global\let\svgscale\undefined%
  \makeatother%
  \begin{picture}(1,0.29007597)%
    \lineheight{1}%
    \setlength\tabcolsep{0pt}%
    \put(0,0){\includegraphics[width=\unitlength,page=1]{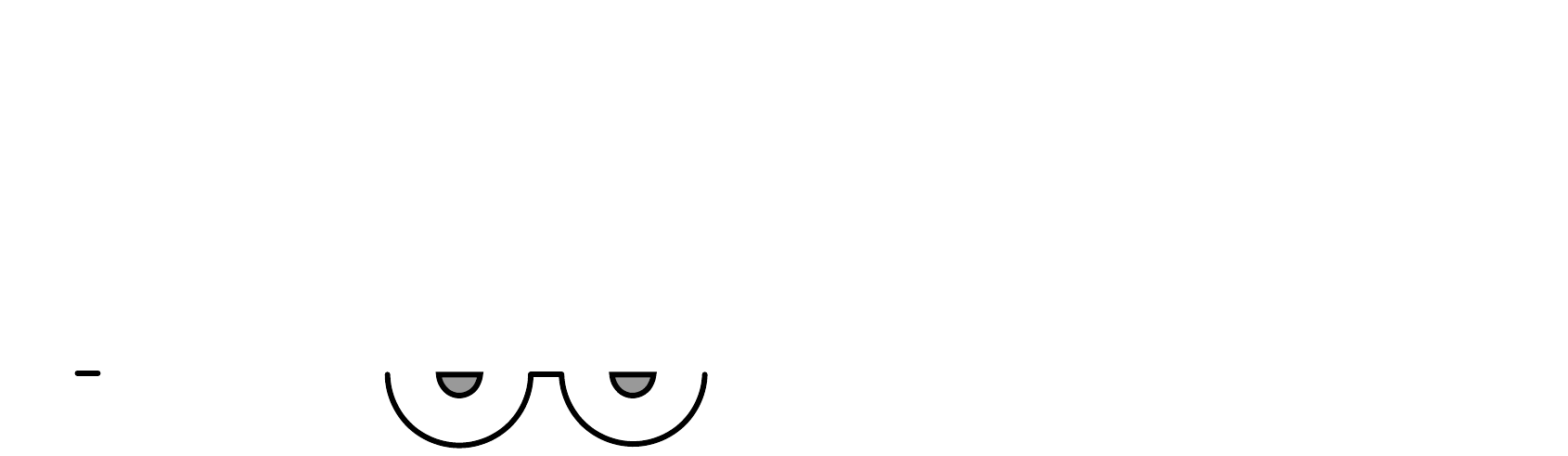}}%
    \put(0.32538474,0.00301206){\color[rgb]{0,0,0}\makebox(0,0)[lt]{\lineheight{1.25}\smash{\begin{tabular}[t]{l}$B$\end{tabular}}}}%
    \put(0.43830279,0.00482156){\color[rgb]{0,0,0}\makebox(0,0)[lt]{\lineheight{1.25}\smash{\begin{tabular}[t]{l}$B$\end{tabular}}}}%
    \put(0,0){\includegraphics[width=\unitlength,page=2]{elements_reduits.pdf}}%
    \put(0.55108255,0.0053328){\color[rgb]{0,0,0}\makebox(0,0)[lt]{\lineheight{1.25}\smash{\begin{tabular}[t]{l}$A$\end{tabular}}}}%
    \put(0.66592707,0.00819767){\color[rgb]{0,0,0}\makebox(0,0)[lt]{\lineheight{1.25}\smash{\begin{tabular}[t]{l}$A$\end{tabular}}}}%
    \put(0,0){\includegraphics[width=\unitlength,page=3]{elements_reduits.pdf}}%
    \put(0.14836254,0.00885557){\color[rgb]{0,0,0}\makebox(0,0)[lt]{\lineheight{1.25}\smash{\begin{tabular}[t]{l}$B$\end{tabular}}}}%
    \put(0,0){\includegraphics[width=\unitlength,page=4]{elements_reduits.pdf}}%
    \put(0.85181284,0.01129953){\color[rgb]{0,0,0}\makebox(0,0)[lt]{\lineheight{1.25}\smash{\begin{tabular}[t]{l}$A$\end{tabular}}}}%
    \put(0,0){\includegraphics[width=\unitlength,page=5]{elements_reduits.pdf}}%
    \put(-0.00136865,0.04584016){\color[rgb]{0,0,0}\makebox(0,0)[lt]{\lineheight{1.25}\smash{\begin{tabular}[t]{l}$W_{\!B}^m$\end{tabular}}}}%
    \put(0.88151519,0.04400205){\color[rgb]{0,0,0}\makebox(0,0)[lt]{\lineheight{1.25}\smash{\begin{tabular}[t]{l}$W_{\!A}^n$\end{tabular}}}}%
  \end{picture}%
\endgroup%

\end{center}
namely $W_{\! B}^m \mathrm{tr}_q\!\left( B \ldots \mathrm{tr}_q\!\left(B \mathrm{tr}_q(BA) A\right) \ldots  A \right) W_{\! A}^n$.
\end{lemma}
\begin{proof}
Apply the fusion relation to the two handles ($B$ and $A$) of $W(L)$. Then resolve all the crossings with the Kauffman relation. We get a diagram without crossings, with a lot of handles labelled $B$ which contain only one strand at the left and a lot of handles labelled $A$ which contain only one strand at the right. In general, after this transformation,  the diagram will contain pieces which for instance look like
\begin{equation}\label{anseMultiplesU}
\begingroup%
  \makeatletter%
  \providecommand\color[2][]{%
    \errmessage{(Inkscape) Color is used for the text in Inkscape, but the package 'color.sty' is not loaded}%
    \renewcommand\color[2][]{}%
  }%
  \providecommand\transparent[1]{%
    \errmessage{(Inkscape) Transparency is used (non-zero) for the text in Inkscape, but the package 'transparent.sty' is not loaded}%
    \renewcommand\transparent[1]{}%
  }%
  \providecommand\rotatebox[2]{#2}%
  \newcommand*\fsize{\dimexpr\f@size pt\relax}%
  \newcommand*\lineheight[1]{\fontsize{\fsize}{#1\fsize}\selectfont}%
  \ifx\svgwidth\undefined%
    \setlength{\unitlength}{318.61406274bp}%
    \ifx\svgscale\undefined%
      \relax%
    \else%
      \setlength{\unitlength}{\unitlength * \real{\svgscale}}%
    \fi%
  \else%
    \setlength{\unitlength}{\svgwidth}%
  \fi%
  \global\let\svgwidth\undefined%
  \global\let\svgscale\undefined%
  \makeatother%
  \begin{picture}(1,0.22967715)%
    \lineheight{1}%
    \setlength\tabcolsep{0pt}%
    \put(0,0){\includegraphics[width=\unitlength,page=1]{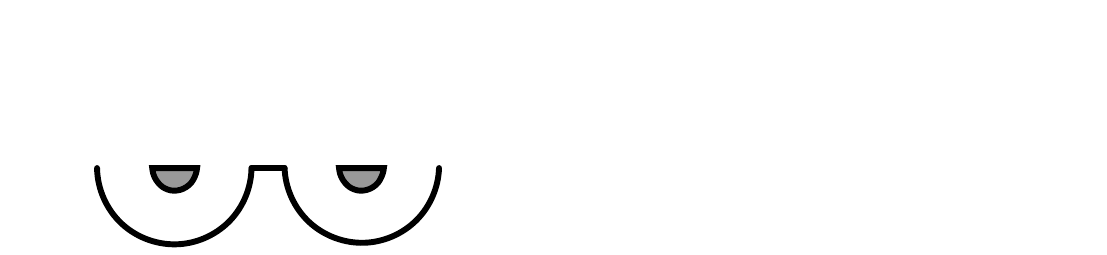}}%
    \put(0.21205417,0.00460032){\color[rgb]{0,0,0}\makebox(0,0)[lt]{\lineheight{1.25}\smash{\begin{tabular}[t]{l}$U$\end{tabular}}}}%
    \put(0.38451375,0.0052649){\color[rgb]{0,0,0}\makebox(0,0)[lt]{\lineheight{1.25}\smash{\begin{tabular}[t]{l}$U$\end{tabular}}}}%
    \put(0,0){\includegraphics[width=\unitlength,page=2]{exemple_preuve_WA_WB.pdf}}%
    \put(0.55676218,0.00744509){\color[rgb]{0,0,0}\makebox(0,0)[lt]{\lineheight{1.25}\smash{\begin{tabular}[t]{l}$U$\end{tabular}}}}%
    \put(0.73146441,0.0118206){\color[rgb]{0,0,0}\makebox(0,0)[lt]{\lineheight{1.25}\smash{\begin{tabular}[t]{l}$U$\end{tabular}}}}%
    \put(0,0){\includegraphics[width=\unitlength,page=3]{exemple_preuve_WA_WB.pdf}}%
    \put(0.89572269,0.00612263){\color[rgb]{0,0,0}\makebox(0,0)[lt]{\lineheight{1.25}\smash{\begin{tabular}[t]{l}$U$\end{tabular}}}}%
    \put(0,0){\includegraphics[width=\unitlength,page=4]{exemple_preuve_WA_WB.pdf}}%
  \end{picture}%
\endgroup%

\end{equation}
where $U$ is $B$ or $A$. But these elements can be transformed in a polynomial in $W_{\! U}$ thanks to \eqref{puissanceM}. For instance, with the piece of diagram below:
\begin{align*}
\mathrm{tr}_q\!\left( U \mathrm{tr}_q(U^2) U^2 \right) &= \mathrm{tr}_q\!\left( U \mathrm{tr}_q(-q^{-1}W_{\! U} U - q^{-2} \mathbb{I}_2) U^2 \right) = (-q^{-1}W_{\! U}^2 + q^{-2}\hat q) \mathrm{tr}_q\!\left(U^3\right)\\
& = (-q^{-1}W_{\! U}^2 + q^{-2}\hat q) \mathrm{tr}_q\!\left(q^{-2}(W_{\! U}^2-1)U + q^{-3}W_{\! U} \mathbb{I}_2\right)\\
&= (-q^{-1}W_{\! U}^2 + q^{-2}\hat q)(q^{-2}(W_{\! U}^2-1)W_{\! U} - \hat qq^{-3}W_{\! U}).
\end{align*}
Then, since $W_{\! U} U = U W_{\! U}$, we can drag the powers of $W_B$ on the left and the powers of $W_A$ on the right:
\begin{center}
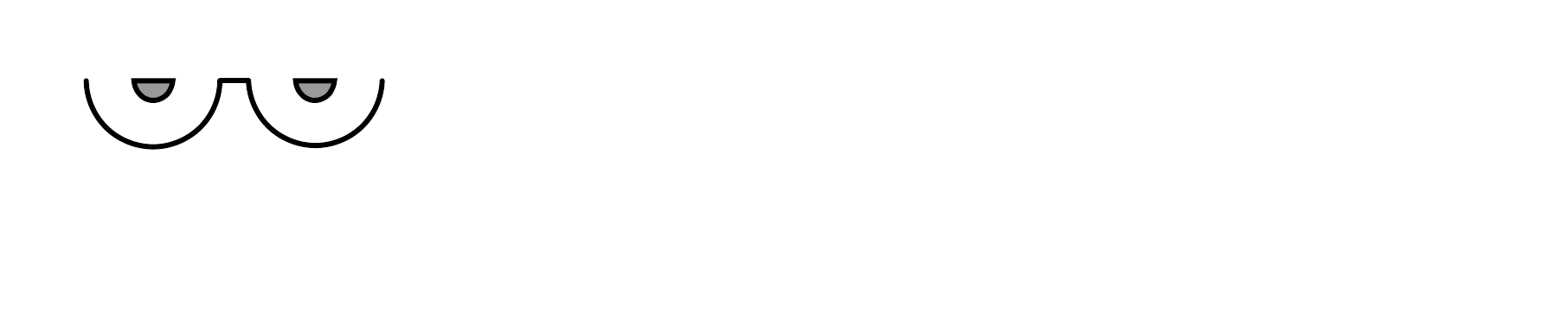
\end{center}
The remaining strands are of the following form: from left to right, they meet several {\em consecutive} handles labelled $B$ and then they meet several {\em consecutive} handles labelled $A$ (note that thanks to the above transformations, the strands cannot meet several handles with the same label ($B$ or $A$) which are not consecutive). Again, this is resolved thanks to \eqref{puissanceM}. For instance, with $B$ and two consecutive handles:
\begin{center}
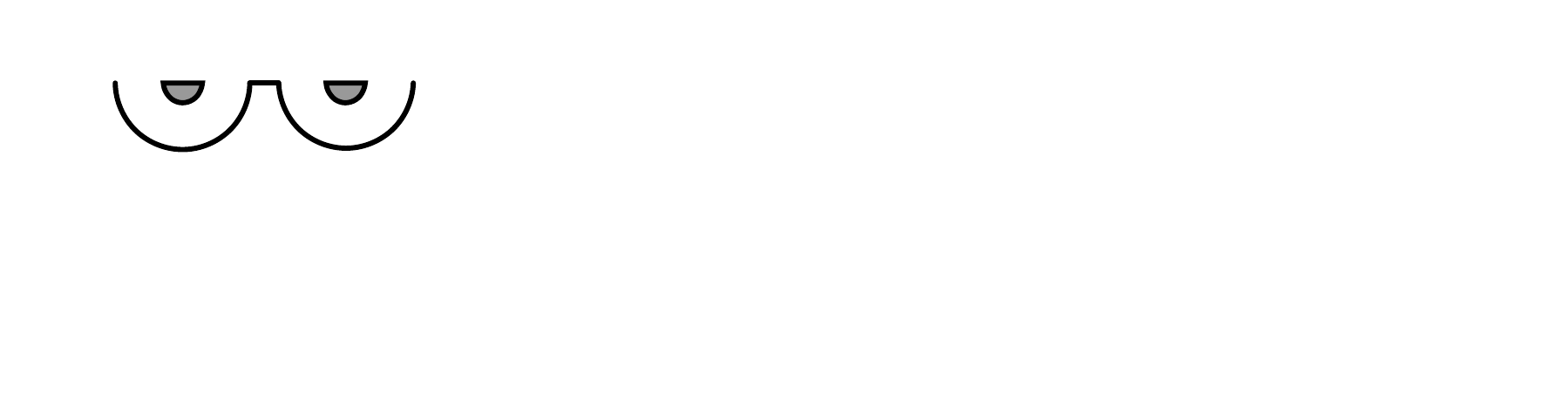
\end{center}
where $U$ is $B$ or $A$. As above, we drag the circle ($=W_{\! B}$) in the first term to the left. We know that the strand $\gamma$ in the second term does not meet any handle labelled by $B$ and using isotopy it is evaluated as $\mathrm{tr}_q(A^i)$ for some $i$; we transform this into a polynomial in $W_{\! A}$ thanks to \eqref{puissanceM} and we drag the result to the right. This gives the desired form.
\end{proof}

\begin{proof}[Proof of Proposition \ref{propWAWB}] It suffices to show that the result is true for the elements of Lemma \ref{diagrammeReduit}. These elements form a sequence which can be defined by induction:
\[ t_1 = W_{\! BA}, \:\:\:\: t_{n+1} = \mathrm{tr}_q\!\left(Bt_nA\right). \]
Thanks to Lemma \ref{lemmaWAWB}, $t_1 \in \mathbb{C}\langle W_{\! A}, W_{\! B} \rangle$. Let us assume that $t_n = P(W_{\! A}, W_{\! B})$ for some non-commutative polynomial $P \in \mathbb{C}\langle x_1, x_2 \rangle$, then due to Lemma \ref{lemmaWAWB}:
\begin{align*}
t_{n+1} &= \mathrm{tr}_q\!\left(BP(W_{\! A}, W_{\! B})A\right) = \mathrm{tr}_q\!\left(Q(W_{\! A}, W_{\!B})A^2 + R(W_{\! A}, W_{\!B})BA + S(W_{\! A}, W_{\!B})BA^2\right)\\
&= \mathrm{tr}_q\!\left(-q^{-1}Q(W_{\! A}, W_{\!B})W_{\!A}A -q^{-2}Q(W_{\! A}, W_{\!B})\mathbb{I}_2  + R(W_{\! A}, W_{\!B})BA - q^{-1}S(W_{\! A}, W_{\!B})BAW_{\! A}\right.\\
&\:\:\:\:\:\:\:\:\:\:\:\: \left.- q^{-2}S(W_{\! A}, W_{\!B})B\right)\\
&= -q^{-1}Q(W_{\! A}, W_{\!B})W_{\!A}^2 + q^{-1}[2]Q(W_{\! A}, W_{\!B}) + R(W_{\! A}, W_{\!B})W_{\! BA} - q^{-1}S(W_{\! A}, W_{\!B})W_{\! BA}W_{\! A}\\ 
&\:\:\:\:\:- q^{-2}S(W_{\! A}, W_{\!B})W_{\! B}\\
&= T(W_{\! A}, W_{\! B}),
\end{align*}
for some $T\in \mathbb{C}\langle x_1, x_2 \rangle$, again due to the fact that $W_{\! BA} \in \mathbb{C}\langle W_{\! A}, W_{\! B} \rangle$. This proves the result by induction.
\end{proof}

\section{Representation of the skein algebra at roots of unity}\label{sectionRepSkeinAlg}
\indent Recall that $q = e^{i\pi/p}$ is a $2p$-th root of unity. Let $\Sigma$ be a compact oriented surface. We denote by $\mathcal{R}(\Sigma)$ the set of isotopy classes of framed links (\textit{i.e.} ribbons) in $\Sigma \times [0, 1]$, and by $\mathbb{C}\mathcal{R}(\Sigma)$ the $\mathbb{C}$-vector space whose basis is $\mathcal{R}(\Sigma)$ (formal linear combinations of elements of $\mathcal{R}(\Sigma)$). Recall that the stack product $\ast$ endows $\mathbb{C}\mathcal{R}(\Sigma)$ with an associative algebra structure. 
Let $\mathcal{K}(\Sigma)$ be the ideal generated by the Kauffman skein relation \eqref{relationKauffman}. The Kauffman skein algebra of $\Sigma$, denoted by $\mathcal{S}_q(\Sigma)$, is $\mathbb{C}\mathcal{R}(\Sigma)/\mathcal{K}(\Sigma)$.

\smallskip

\indent Thanks to Corollary \ref{indepOrientation}, we have a map $W : \mathbb{C}\mathcal{R}(\Sigma_{g,n}^{\mathrm{o}}) \to \mathcal{L}_{g,n}(\bar U_q)$ (where we implicitly color all the elements of $\mathcal{R}(\Sigma_{g,n}^{\mathrm{o}})$ by $\mathcal{X}^+(2)$). Since for $H = \bar U_q$ the diagrammatic calculus satisfies the Kauffman skein relation, we get a morphism of algebras $W : \mathcal{S}_q(\Sigma_{g,n}^{\mathrm{o}}) \to \mathcal{L}_{g,n}(\bar U_q)$, and more precisely $W : \mathcal{S}_q(\Sigma_{g,n}^{\mathrm{o}}) \to \mathcal{L}_{g,n}^{\mathrm{inv}}(\bar U_q)$. This provides representations of the skein algebra of $\Sigma_{g,n}^{\mathrm{o}}$ for all $g,n$. Analogously to the representation of the mapping class group, we can glue back the disk $D$ by passing to the invariants of the representation, and for $n=0$ we will see that this indeed provides a representation of $\mathcal{S}_q(\Sigma_{g})$ (the result is probably true for any $\Sigma_{g,n}$). 

\smallskip

\indent Let $j : \Sigma_{g}^{\mathrm{o}} \to \Sigma_{g}$ be the canonical embedding and let $j = j \times \mathrm{id} : \Sigma_{g}^{\mathrm{o}} \times [0,1] \to \Sigma_{g} \times [0,1]$ be the corresponding embedding. The corresponding map $j : \mathcal{R}(\Sigma_{g}^{\mathrm{o}}) \to \mathcal{R}(\Sigma_{g})$ is surjective, as well as the corresponding morphism $j : \mathcal{S}_q(\Sigma_{g}^{\mathrm{o}}) \to \mathcal{S}_q(\Sigma_{g})$. If $L \in \mathcal{S}_q(\Sigma_{g})$, we denote by $L^{\mathrm{o}}$ any element of $\mathcal{S}_q(\Sigma_{g}^{\mathrm{o}})$ such that $j(L^{\mathrm{o}}) = L$. Since $j$ is a morphism, it holds $(L_1 \ast L_2)^{\mathrm{o}} = L_1^{\mathrm{o}} \ast L_2^{\mathrm{o}}$.

\begin{theorem}\label{theoRepSkein}
1. Let $\rho : \mathcal{L}_{g,n}(\bar U_q) \to \mathrm{End}_{\mathbb{C}}(V)$ be a representation (with $V = (\bar U_q^*)^{\otimes g} \otimes I_1 \otimes \ldots \otimes I_n$, where $I_1, \ldots, I_n$ are representations of $\bar U_q$). The map
\[ \fleche{\mathcal{S}_q(\Sigma_{g,n}^{\mathrm{o}})}{\mathrm{End}_{\mathbb{C}}(V)}{L}{\rho(W(L))} \]
is a representation of $\mathcal{S}_q(\Sigma_{g,n}^{\mathrm{o}})$.
\\2. Assume $n=0$ and let $\rho_{\mathrm{inv}}$ be the representation of $\mathcal{L}_{g,0}^{\mathrm{inv}}(\bar U_q)$ on $\mathrm{Inv}\bigl( (\bar U_q^*)^{\otimes g} \bigr)$. The map
\[ \fleche{\mathcal{S}_q(\Sigma_{g})}{\mathrm{End}_{\mathbb{C}}\!\left(\mathrm{Inv}\bigl( (\bar U_q^*)^{\otimes g} \bigr)\right)}{L}{\rho_{\mathrm{inv}}\!\left(W(L^{\mathrm{o}})\right)} \]
is well-defined and is a representation of $\mathcal{S}_q(\Sigma_{g})$.
\end{theorem}
\begin{proof}
1. is obvious.\\
2. Since $j$ is a morphism, it holds $(L_1 \ast L_2)^{\mathrm{o}} = L_1^{\mathrm{o}} \ast L_2^{\mathrm{o}}$. It remains to show that it is well-defined. Thanks to the Kauffman skein relation, we can assume that $L$ is a simple closed curve in $\Sigma_{g} \times \{0\}$. Hence, since $\pi_1(\Sigma_{g}) = \pi_1(\Sigma_{g}^{\mathrm{o}})/\langle c_{g,0} \rangle$ (where $c_{g,0}$ is the boundary curve), it suffices to check that $\rho_{\mathrm{inv}}(W(\gamma)) = \rho_{\mathrm{inv}}(W(\gamma_c))$, where $\gamma, \gamma_c \subset \Sigma_g^{\mathrm{o}} \times \{0\}$  are the simple closed curves depicted below:
\begin{center}
\begingroup%
  \makeatletter%
  \providecommand\color[2][]{%
    \errmessage{(Inkscape) Color is used for the text in Inkscape, but the package 'color.sty' is not loaded}%
    \renewcommand\color[2][]{}%
  }%
  \providecommand\transparent[1]{%
    \errmessage{(Inkscape) Transparency is used (non-zero) for the text in Inkscape, but the package 'transparent.sty' is not loaded}%
    \renewcommand\transparent[1]{}%
  }%
  \providecommand\rotatebox[2]{#2}%
  \newcommand*\fsize{\dimexpr\f@size pt\relax}%
  \newcommand*\lineheight[1]{\fontsize{\fsize}{#1\fsize}\selectfont}%
  \ifx\svgwidth\undefined%
    \setlength{\unitlength}{316.01829055bp}%
    \ifx\svgscale\undefined%
      \relax%
    \else%
      \setlength{\unitlength}{\unitlength * \real{\svgscale}}%
    \fi%
  \else%
    \setlength{\unitlength}{\svgwidth}%
  \fi%
  \global\let\svgwidth\undefined%
  \global\let\svgscale\undefined%
  \makeatother%
  \begin{picture}(1,0.15719691)%
    \lineheight{1}%
    \setlength\tabcolsep{0pt}%
    \put(0,0){\includegraphics[width=\unitlength,page=1]{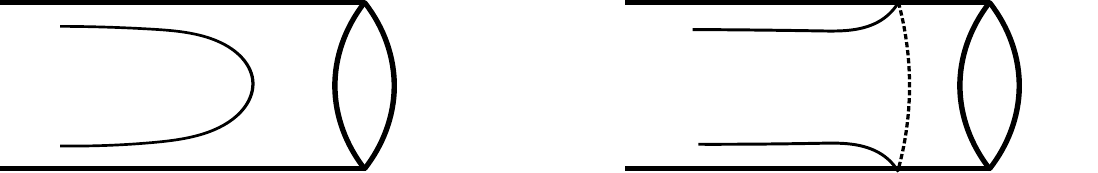}}%
    \put(0.22092436,0.11967175){\color[rgb]{0,0,0}\makebox(0,0)[lt]{\lineheight{1.25}\smash{\begin{tabular}[t]{l}$\gamma$\end{tabular}}}}%
    \put(0.77943887,0.10149035){\color[rgb]{0,0,0}\makebox(0,0)[lt]{\lineheight{1.25}\smash{\begin{tabular}[t]{l}$\gamma_c$\end{tabular}}}}%
    \put(0,0){\includegraphics[width=\unitlength,page=2]{surfacePreuveSkein.pdf}}%
  \end{picture}%
\endgroup%

\end{center}
These pictures represent a neighborhood of the boundary (see Figure \ref{figureSurface}). Take a basepoint on each circle and endow it with the positive orientation as follows:
\begin{center}
\begingroup%
  \makeatletter%
  \providecommand\color[2][]{%
    \errmessage{(Inkscape) Color is used for the text in Inkscape, but the package 'color.sty' is not loaded}%
    \renewcommand\color[2][]{}%
  }%
  \providecommand\transparent[1]{%
    \errmessage{(Inkscape) Transparency is used (non-zero) for the text in Inkscape, but the package 'transparent.sty' is not loaded}%
    \renewcommand\transparent[1]{}%
  }%
  \providecommand\rotatebox[2]{#2}%
  \newcommand*\fsize{\dimexpr\f@size pt\relax}%
  \newcommand*\lineheight[1]{\fontsize{\fsize}{#1\fsize}\selectfont}%
  \ifx\svgwidth\undefined%
    \setlength{\unitlength}{294.33269813bp}%
    \ifx\svgscale\undefined%
      \relax%
    \else%
      \setlength{\unitlength}{\unitlength * \real{\svgscale}}%
    \fi%
  \else%
    \setlength{\unitlength}{\svgwidth}%
  \fi%
  \global\let\svgwidth\undefined%
  \global\let\svgscale\undefined%
  \makeatother%
  \begin{picture}(1,0.16826678)%
    \lineheight{1}%
    \setlength\tabcolsep{0pt}%
    \put(0,0){\includegraphics[width=\unitlength,page=1]{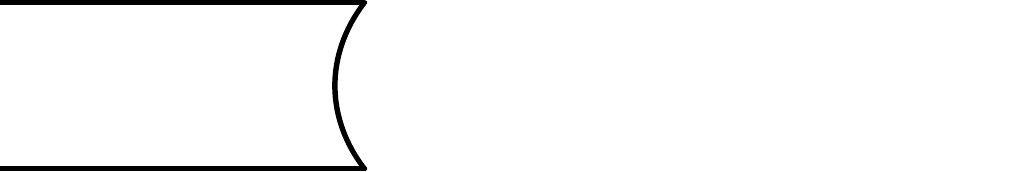}}%
    \put(0.20748183,0.01779206){\color[rgb]{0,0,0}\makebox(0,0)[lt]{\lineheight{1.25}\smash{\begin{tabular}[t]{l}$x$\end{tabular}}}}%
    \put(0,0){\includegraphics[width=\unitlength,page=2]{surfacePreuveSkeinPointBase.pdf}}%
    \put(0.86996362,0.03758763){\color[rgb]{0,0,0}\makebox(0,0)[lt]{\lineheight{1.25}\smash{\begin{tabular}[t]{l}$x_c$\end{tabular}}}}%
    \put(0,0){\includegraphics[width=\unitlength,page=3]{surfacePreuveSkeinPointBase.pdf}}%
  \end{picture}%
\endgroup%

\end{center}
We get two positively oriented simple loops $x, x_c=x^{-1}c_{g,0} \in \pi_1(\Sigma_g^{\mathrm{o}})$ such that $[x] = \gamma, [x_c] = \gamma_c$. It holds\footnote{In general it is of course not true that $\widetilde{xy} = \widetilde{x}\widetilde{y}$.} $\widetilde{x_c} = \widetilde{x^{-1}}C_{g,0}$, with $C_{g,0} = \widetilde{c_{g,0}}$ (the lifts are implicitly considered in the fundamental representation $\mathcal{X}^+(2)$). Indeed, to join the loops $x^{-1}$ and $c_{g,0}$ we must necessarily add a cap going from right to left and thus $N(x_c) = N(x^{-1}) + N(c_{g,0}) - 1$ (see section \ref{sectionNorma}). Moreover, due to \eqref{NLoopNFree} and \eqref{NInvFree},
\[ N(x^{-1}) = N(\gamma^{-1}) = -N(\gamma) = -N(x) + 1. \]
This yields 
\[
\begin{split}
\widetilde{x_c} &= \mathrm{ev}_{\mathcal{X}^+(2)}\bigl( v^{N(x^{-1}c_{g,0})}x^{-1}c_{g,0} \bigr) = \mathrm{ev}_{\mathcal{X}^+(2)}\bigl( v^{N(x^{-1}) + N(c_{g,0}) - 1}x^{-1}c_{g,0} \bigr) = \mathrm{ev}_{\mathcal{X}^+(2)}\bigl( v^{-N(x) + N(c_{g,0})}x^{-1}c_{g,0} \bigr)\\
& = \mathrm{ev}_{\mathcal{X}^+(2)}\bigl( v^{-N(x)}x^{-1} \bigr)\mathrm{ev}_{\mathcal{X}^+(2)}\bigl( v^{N(c_{g,0})}c_{g,0} \bigr) = \mathrm{ev}_{\mathcal{X}^+(2)}\bigl( v^{N(x)}x \bigr)^{-1} \mathrm{ev}_{\mathcal{X}^+(2)}\bigl( v^{N(c_{g,0})}c_{g,0} \bigr) = \widetilde{x^{-1}} C_{g,0}
\end{split}
\]
Hence, by Proposition \ref{propWilsonSimple} we have
\begin{align*}
\rho_{\mathrm{inv}}\bigl(W(\gamma_c)\bigr) &= \rho_{\mathrm{inv}}\bigl(W(x^{-1}c_{g,0})\bigr) = \rho_{\mathrm{inv}}\bigl(\mathrm{tr}_q(\widetilde{x^{-1}c_{g,0}})\bigr) = \rho_{\mathrm{inv}}\bigl(\mathrm{tr}_q(\widetilde{x^{-1}} C_{g,0})\bigr) = \rho_{\mathrm{inv}}\bigl(\mathrm{tr}_q(\widetilde{x}^{-1})\bigl)\\
& = \rho_{\mathrm{inv}}\bigl(W(x^{-1})\bigr) = \rho_{\mathrm{inv}}\bigl(W(x)\bigr) = \rho_{\mathrm{inv}}\bigl(W(\gamma)\bigr).
\end{align*}
We used that $\rho_{\mathrm{inv}}(C_{g,0}) = \mathbb{I}_2$ (by definition) and that $W$ does not depend of the basepoint (always true) nor of the orientation (when $H = \bar U_q$).
\end{proof}

\section{Explicit study of the representation of $\mathcal{S}_q(\Sigma_1)$}\label{repSqSigma1Uq}
\indent By Theorem \ref{theoRepSkein}, we have a representation of $\mathcal{S}_q(\Sigma_1)$ on $\mathrm{Inv}(\bar U_q^*) = \mathrm{SLF}(\bar U_q)$ and by Proposition \ref{propWAWB}, $W\!\left(\mathcal{S}_q(\Sigma_1^{\mathrm{o}})\right) = \langle W_{\! A}, W_{\! B} \rangle$. Hence, to study the representation of $\mathcal{S}_q(\Sigma_1)$ on $\mathrm{SLF}(\bar U_q)$, it suffices to study the action of the operators $\rho_{\mathrm{inv}}(W_{\! A})$ and $\rho_{\mathrm{inv}}(W_{\! B})$. For this, we use again the GTA basis.

\subsection{Structure of the representation}\label{structRepSqS1}
Recall that we denote by $\triangleright$ the representation of $\mathcal{L}_{1,0}(\bar U_q)$ on $\bar U_q^*$. Also recall the formulas (which are consequences of Proposition \ref{actionAB} and of \eqref{DWChi}):
\[ W_{\! A} \triangleright \varphi = -\hat q^2 \varphi^C, \:\:\:\:\:\: W_{\! B} \triangleright \varphi = \left( \chi^+_2 \varphi^v \right)^{v^{-1}}, \]
where $\varphi^z = \varphi(z?)$ and 
\begin{equation*}
C = FE + \frac{qK + q^{-1}K^{-1}}{(q - q^{-1})^2} = \sum_{j=0}^p \frac{q^j + q^{-j}}{ (q - q^{-1})^2 } e_j + \sum_{k = 1}^{p-1} (w^+_k + w^-_k) \in \mathcal{Z}(\bar U_q)
\end{equation*}
is the Casimir element, with its decomposition in the canonical basis of the center (see Definition \ref{defBaseCanoCentre}) . We denote by 
\begin{equation}\label{symetrieV}
v_s = v_{\mathcal{X}^+(s)} = v_{\mathcal{X}^-(p-s)} = (-1)^{s-1}q^{\frac{-(s^2-1)}{2}}
\end{equation}
the scalar corresponding to the action of $v$ on the simple module $\mathcal{X}^+(s)$ or $\mathcal{X}^-(p-s)$\footnote{The symmetry property expressed in the second equality of \eqref{symetrieV} is true for any central element of $\bar U_q$, see \eqref{symetrieCentraux} 
} ($v_0$ being $v_{\mathcal{X}^-(p)}$).

\smallskip

Let us compute the actions of $W_{\! A}$ and $W_{\! B}$ on the GTA basis. First, using the expression of $C$ in the canonical basis of $\mathcal{Z}(\bar U_q)$ above and the formulas \eqref{actionCentreSLF} for the action of $\mathcal{Z}(\bar U_q)$ on $\mathrm{SLF}(\bar U_q)$, we get:
\begin{equation}\label{actionWAWB}
W_{\! A} \triangleright \chi^{\epsilon}_s = -\epsilon (q^s + q^{-s})\chi^{\epsilon}_s, \:\:\:\:\:\: W_{\! A} \triangleright G_s = -(q^s + q^{-s})G_s - \hat q^2(\chi^+_s + \chi^-_{p-s}). 
\end{equation}
To compute the action of $W_{\! B}$, we must use the multiplication rules in the GTA basis (Theorem \ref{ProduitArike}), the expressions of $v$ and $v^{-1}$ in the canonical basis of $\mathcal{Z}(\bar U_q)$ \eqref{rubanCentre} and the the formulas \eqref{actionCentreSLF} for the action of $\mathcal{Z}(\bar U_q)$ on $\mathrm{SLF}(\bar U_q)$. 
If $1 < s < p$:
\[ W_{\! B} \triangleright \chi^{\epsilon}_s = v_{\mathcal{X}^{\epsilon}(s)}(\chi^+_2 \chi^+_s)^{v^{-1}} = v_{\mathcal{X}^{\epsilon}(s)}(\chi^{\epsilon}_{s-1} + \chi^{\epsilon}_{s+1})^{v^{-1}} = \frac{v_{\mathcal{X}^{\epsilon}(s)}}{v_{\mathcal{X}^{\epsilon}(s-1)}}\chi^{\epsilon}_{s-1} + \frac{v_{\mathcal{X}^{\epsilon}(s)}}{v_{\mathcal{X}^{\epsilon}(s+1)}}\chi^{\epsilon}_{s+1}. \]
For $s=1$:
\[ W_{\! B} \triangleright \chi^{\epsilon}_1 = v_{\mathcal{X}^{\epsilon}(1)}(\chi^+_2 \chi^{\epsilon}_1)^{v^{-1}} = v_{\mathcal{X}^{\epsilon}(1)}(\chi^{\epsilon}_2)^{v^{-1}}  = \frac{v_{\mathcal{X}^{\epsilon}(1)}}{v_{\mathcal{X}^{\epsilon}(2)}}\chi^{\epsilon}_2. \]
And for $s=p$:
\begin{align*}
 W_{\! B} \triangleright \chi^{\epsilon}_p &= v_{\mathcal{X}^{\epsilon}(p)}(\chi^+_2 \chi^{\epsilon}_p)^{v^{-1}} = v_{\mathcal{X}^{\epsilon}(p)}(2\chi^{\epsilon}_{p-1} + 2\chi^{-\epsilon}_{1})^{v^{-1}} = 2\frac{v_{\mathcal{X}^{\epsilon}(p)}}{v_{\mathcal{X}^{\epsilon}(p-1)}}\chi^{\epsilon}_{p-1} + 2\frac{v_{\mathcal{X}^{\epsilon}(p)}}{v_{\mathcal{X}^{-\epsilon}(1)}}\chi^{-\epsilon}_1\\
&= 2\frac{v_{\mathcal{X}^{\epsilon}(p)}}{v_{\mathcal{X}^{\epsilon}(p-1)}}(\chi^{\epsilon}_{p-1} + \chi^{-\epsilon}_1),
\end{align*}
thanks to the \eqref{symetrieV}.

\smallskip

\indent Let 
\[ \mathcal{P} = \mathrm{vect}\!\left( \chi^+_s + \chi^-_{p-s}, \chi^{\pm}_p \right)_{1 \leq s \leq p-1} = \mathrm{vect}\!\left(\chi^P\right)_{P \in \mathrm{Proj}_{\bar U_q}} \]
be the subspace generated by the characters of the projective $\bar U_q$-modules. Introduce notation for the basis elements of $\mathcal{P}$:
\[ X_0 = \chi^-_p, \:\:\:\:\: X_s = \chi^+_s + \chi^-_{p-s} \:\: \text{ for } 1 \leq s \leq p-1 , \:\:\:\:\: X_p = \chi^+_p. \]
The formulas above and \eqref{symetrieV} give
\begin{eqnarray*}
W_{\! A} \triangleright X_s &=& -(q^s + q^{-s})X_s\:\:\:\:\:\:\:\:\:  \text{for } 0 \leq s \leq p,\\
W_{\! B} \triangleright X_0 &=&  2\frac{v_0}{v_1}X_1 = -2q^{\frac{1}{2}}X_1 , \\
W_{\! B} \triangleright X_s &=&  \frac{v_s}{v_{s-1}}X_{s-1} + \frac{v_s}{v_{s+1}}X_{s+1} = -q^{-s + \frac{1}{2}}X_{s-1} - q^{s + \frac{1}{2}}X_{s+1} \:\:\:\:\:\:\:\:\: \text{for } 1 \leq s \leq p-1.\\
W_{\! B} \triangleright X_p &=&  2\frac{v_p}{v_{p-1}}X_{p-1} = 2q^{\frac{1}{2}}X_{p-1}.
\end{eqnarray*}
In particular, $\mathcal{P}$ is a submodule of $\mathrm{SLF}(\bar U_q)$ under the action of $\mathcal{S}_q(\Sigma_1)$.

\smallskip

\indent Now we compute the action on $G_s$. First, note that
\begin{align*}
G_t^v &= v_t G_t + \hat q v_t \frac{p}{[t]}\chi^+_t - \hat q v_t \frac{t}{[t]}X_t,\\
G_t^{v^{-1}} &= v_t^{-1} G_t - \hat q v_t^{-1} \frac{p}{[t]}\chi^+_t + \hat q v_t^{-1} \frac{t}{[t]}X_t.
\end{align*}
Hence:
\begin{align*}
&W_{\! B} \triangleright G_s = \left( \chi^+_2 G_s^v \right)^{v^{-1}} = v_s\left( \chi^+_2 G_s \right)^{v^{-1}} + \hat q v_s\frac{p}{[s]}\left( \chi^+_2 \chi^+_s \right)^{v^{-1}} - \hat q v_s\frac{s}{[s]}\left( \chi^+_2 X_s\right)^{v^{-1}}\\
&= v_s\frac{[s-1]}{[s]}G_{s-1}^{v^{-1}} + v_s\frac{[s+1]}{[s]}G_{s+1}^{v^{-1}} + \hat q \frac{v_s}{v_{s-1}} \frac{p}{[s]} \chi^+_{s-1} + \hat q \frac{ v_s }{ v_{s+1} } \frac{p}{[s]} \chi^+_{s+1}
 - \hat q \frac{ v_s }{ v_{s-1} } \frac{s}{[s]} X_{s-1} - \hat q \frac{ v_s }{ v_{s+1} } \frac{s}{[s]} X_{s+1}\\
&= \frac{ v_s }{ v_{s-1} } \frac{[s-1]}{[s]} G_{s-1} + \frac{ v_s }{ v_{s+1} } \frac{[s-1]}{[s]} G_{s+1} - \hat q \frac{v_s}{v_{s-1}} \frac{p}{[s]} \chi^+_{s-1} + \hat q \frac{v_s}{v_{s-1}} \frac{s-1}{[s]} X_{s-1} - \hat q \frac{v_s}{v_{s+1}} \frac{p}{[s]} \chi^+_{s+1}\\
&\:\:\:\:\: + \hat q \frac{v_s}{v_{s+1}} \frac{s+1}{[s]} X_{s+1} + \hat q \frac{v_s}{v_{s-1}} \frac{p}{[s]} \chi^+_{s-1} + \hat q \frac{v_s}{v_{s+1}} \frac{p}{[s]} \chi^+_{s+1} - \hat q \frac{v_s}{v_{s-1}} \frac{s}{[s]}X_{s-1} - \hat q \frac{v_s}{v_{s+1}} \frac{s}{[s]} X_{s+1}\\
&=  \frac{ v_s }{ v_{s-1} } \frac{[s-1]}{[s]} G_{s-1} + \frac{ v_s }{ v_{s+1} } \frac{[s+1]}{[s]} G_{s+1} - \hat q \frac{v_s}{v_{s-1}} \frac{1}{[s]} X_{s-1} + \hat q \frac{v_s}{v_{s+1}} \frac{1}{[s]} X_{s+1}\\
&= -q^{-s + \frac{1}{2}}\frac{[s-1]}{[s]} G_{s-1} - q^{s + \frac{1}{2}}\frac{[s+1]}{[s]} G_{s+1} + \hat q q^{-s + \frac{1}{2}} \frac{1}{[s]} X_{s-1} - \hat q q^{s + \frac{1}{2}} \frac{1}{[s]} X_{s+1}.
\end{align*}
for all $1 \leq s \leq p-1$ (with the convention that for $s=1$ and $s=p-1$ the undefined terms are $0$).

\smallskip

\indent Consider the following subspaces of $\mathrm{SLF}(\bar U_q)$:
\[ \mathcal{U} = \mathrm{vect}\!\left( \chi^+_s \right)_{1 \leq s \leq p-1}, \:\:\:\:\:\:\:\:\:\: \mathcal{V} = \mathrm{vect}\!\left( G_s \right)_{1 \leq s \leq p-1}. \]
The formulas above reveal that the structure of $\mathrm{SLF}(\bar U_q)$ under the action of $\mathcal{S}_q(\Sigma_1)$ has the following shape:
\[
\xymatrix{
\mathcal{U}  \ar[rd]_{W_{\! B}} & &  \mathcal{V} \ar[ld]^{W_{\! A}, W_{\! B}}\\  
 &\mathcal{P} &
}
\]
Let us check that this gives rise to a composition series:
\[ J_1 = \mathcal{P} \: \subset \: J_2 = J_1 \oplus \mathcal{U} \: \subset \: J_3 = J_2 \oplus \mathcal{V} = \mathrm{SLF}(\bar U_q). \]
%

\begin{lemma}\label{lemmeEsWsPourIndec}
Recall Definition \ref{defBaseCanoCentre} and notation \eqref{notationImbed}. Then $(e_s)_A \in \mathbb{C}\langle W_{\! A} \rangle$ for all $0 \leq s \leq p$ and $(w^+_t + w^-_t)_A \in \mathbb{C}\langle W_{\! A} \rangle$ for all $1 \leq t \leq p-1$; in particular these elements belong to $W\bigl( \mathcal{S}_q(\Sigma_1) \bigr)$.
\end{lemma}
\begin{proof}
Recall from \eqref{engendreParCasimir} that the central elements $e_s$ and $w^+_t + w^-_{p-t}$ belong to the subalgebra of $\mathcal{Z}(\bar U_q)$ generated by the Casimir element $C$. Now, in \eqref{WCasimir}, we computed that
\[ \overset{\mathcal{X}^+(2)}{W_{\! M}} = \mathrm{tr}\bigl( \overset{\mathcal{X}^+(2)}{K^{p+1}} \overset{\mathcal{X}^+(2)}{M} \bigr) = - \hat q^2 C \]
with $\hat q = q - q^{-1}$, under the identification $\Psi_{0,1}$ between $\mathcal{L}_{0,1}(\bar U_q)$ and $\bar U_q$. It follows that the elements $e_s, w^+_t + w^-_t$ can be written as polynomials of $\overset{\mathcal{X}^+(2)}{W_{\! M}}$. Hence, applying the morphism $j_A : \mathcal{L}_{0,1}(\bar U_q) \to \mathcal{L}_{1,0}(\bar U_q)$ defined by $\overset{I}{M} \mapsto \overset{I}{A}$, we get that the elements $(e_s)_A = j_A(e_s), (w^+_t + w^-_t)_A = j_A(w^+_t + w^-_t)$ can be written as polynomials in $W_{\! A} = \overset{\mathcal{X}^+(2)}{W_{\! A}} = j_A\bigl(\overset{\mathcal{X}^+(2)}{W_{\! M}}\bigr)$.
\end{proof}

\begin{proposition}\label{structureRepSqS1}
$J_1 \subset J_2 \subset J_3$ is a composition series of $\mathrm{SLF}(\bar U_q)$ under the action of $\mathcal{S}_q(\Sigma_1)$. More precisely, the structure of the representation is schematized by the following diagram:
\[
\xymatrix{
\mathcal{U}  \ar[rd]_{W_{\! B}} & &  \mathcal{V} \ar[ld]^{W_{\! A}, W_{\! B}}\\  
 &\mathcal{P} &
}
\]
Moreover, this representation is indecomposable.
\end{proposition}
\begin{proof}
We will use the element $W(b^{-1}a) = W_{\! vB^{-1}A}$ which by Lemma \ref{actionBmoinsUnA} implements the multiplication by $\chi^+_2$:
\[ W_{\! vB^{-1}A} \triangleright \varphi = \chi^+_2 \varphi. \]
\textbullet ~~ {\em $J_1$ is simple:} Let $0 \neq S \subset J_1$ be a submodule, and let $0 \neq \psi = x_0 X_0 + \ldots  + x_p X_p \in S$. Thanks to Lemma \ref{lemmeEsWsPourIndec} we can use the elements $(e_s)_A$. Note that \eqref{actionCentreSLF} gives $(e_i)_A \triangleright X_j = \delta_{i,j}X_j$, and thus
\[ \forall \, 0 \leq j \leq p, \:\:\:\:\: (e_j)_A \triangleright \psi = \psi(e_j ?) = x_jX_j \in S. \]
Since $x \neq 0$ one of the $x_j$, say $x_s$, is not $0$. Then $X_s \in S$, and using $Y$ we get 
\begin{equation*}
\begin{split}
&(e_{s-1})_A W_{\! vB^{-1}A} \triangleright X_s = (e_{s-1})_A \triangleright (X_{s-1} + X_{s+1}) = X_{s-1} \in S,\\
&(e_{s+1})_A W_{\! vB^{-1}A} \triangleright X_s = (e_{s+1})_A \triangleright (X_{s-1} + X_{s+1}) = X_{s+1} \in S.
\end{split}
\end{equation*}
Continuing like this, we get step by step that all the $X_j$'s belong to $S$, and hence $S = J_1$ as desired.
\smallskip

\noindent \textbullet ~~ {\em $J_2/J_1$ is simple:} Let $\overline{\chi}^+_s = \chi^+_s + J_1$ for $1 \leq s \leq p-1$; these elements form a basis of $J_2/J_1$. We have
\begin{equation*}
\begin{split}
&(e_i)_A \triangleright \overline{\chi}^+_j = \delta_{i,j}\overline{\chi}^+_j, \\
W_{\! vB^{-1}A} \triangleright \overline{\chi}^+_1 = \overline{\chi}^+_2, \:\:\:\:\:\:\: &W_{\! vB^{-1}A} \triangleright \overline{\chi}^+_s = \overline{\chi}^+_{s-1} + \overline{\chi}^+_{s+1} \: \text{ for } 2 \leq s \leq p-2, \:\:\:\:\:\:\: W_{\! vB^{-1}A} \triangleright \overline{\chi}^+_{p-1} = \overline{\chi}^+_{p-2}.
\end{split}
\end{equation*}
The same reasoning as for $J_1$ gives the result.

\smallskip

\noindent \textbullet ~~ {\em $J_3/J_2$ is simple:} Let $\overline{G}_s = G_s + J_2$ for $1 \leq s \leq p-1$; these elements form a basis of $J_3/J_2$. We have
\begin{align*}
&(e_i)_A \triangleright \overline{G}_j = \delta_{i,j}\overline{G}_j, \\
&W_{\! vB^{-1}A} \triangleright \overline{G}_1 = [2]\overline{G}_2,\\
&W_{\! vB^{-1}A} \triangleright \overline{G}_s = \frac{[s-1]}{[s]}\overline{G}_{s-1} + \frac{[s+1]}{[s]}\overline{G}_{s+1} \: \text{ for } 2 \leq s \leq p-2, \\
&W_{\! vB^{-1}A} \triangleright \overline{G}_{p-1} = [2]\overline{G}_{p-2}.
\end{align*}
The same reasoning as for $J_1$ gives the result.

\smallskip

\indent For the last claim, write $\mathrm{SLF}(\bar U_q) = U_1 \oplus U_2$. At least one of the two subspaces $U_1$ or $U_2$ contains an element of the form $\varphi = G_1 + \sum_{i \neq 1} \lambda_i G_i + \sum_{j, \epsilon} \beta^{\epsilon}_j\chi^{\epsilon}_j$. Assume for instance that it is $U_1$. Then, thanks to Lemma \ref{lemmeEsWsPourIndec} we can use the elements $(w_s^+ + w^-_s)_A$ and \eqref{actionCentreSLF} yields
\[ (w^+_1 + w^-_1)_A \triangleright \varphi = \varphi\bigl( (w_1^+ + w^-_1) ? \bigr) = G_1\bigl( (w_1^+ + w^-_1) ? \bigr) = \chi^+_1 + \chi^-_{p-1} \in \mathcal{P} \cap U_1. \]
It follows that $\mathcal{P} \subset U_1$. Now, let $\psi = \sum_{i} \eta_i G_i + \sum_{j, \epsilon} \gamma^{\epsilon}_j\chi^{\epsilon}_j \in U_2$; then 
\begin{equation*}
\begin{split}
&\forall \, 1 \leq i \leq p-1, \:\:\:\: (w^+_i + w^-_i)_A \triangleright \psi = \eta_i (\chi^+_i + \chi^-_{p-i}) \in \mathcal{P} \cap U_2 = \{ 0 \} \:\:\: \text{ and thus } \eta_i = 0,\\
&\forall \, 0 \leq j \leq p, \:\:\:\: (e_j)_A \triangleright \psi = \eta_j G_j + \gamma^+_j\chi^+_j + \gamma^-_{p-j}\chi^-_{p-j} = \gamma^+_j\chi^+_j + \gamma^-_{p-j}\chi^-_{p-j} \in U_2
\end{split}
\end{equation*}
with the convention that $\chi^+_0 = \chi^-_0 = 0$. Now for any $j$, we have
\begin{equation*}
\begin{split}
&(e_{j-1})_AW_{vB^{-1}A} \triangleright (\gamma^+_j\chi^+_j + \gamma^-_{p-j}\chi^-_{p-j}) = \gamma^+_j\chi^+_{j-1} + \gamma^-_{p-j}\chi^-_{p-j+1} \in U_2,\\
&(e_{j+1})_AW_{vB^{-1}A} \triangleright (\gamma^+_j\chi^+_j + \gamma^-_{p-j}\chi^-_{p-j}) = \gamma^+_j\chi^+_{j+1} + \gamma^-_{p-j}\chi^-_{p-j-1} \in U_2.
\end{split}
\end{equation*}
Continuing as long as necessary to apply $(e_{j \pm k})_AW_{vB^{-1}A}$, we will get that $\gamma^+_j\chi^+_p, \gamma^-_j\chi^-_p \in U_2$. But $\chi^+_p, \chi^-_p \in \mathcal{P} \subset U_1$ and it follows that $\gamma^+_j=0, \gamma^-_j=0$. Hence $\psi = 0$ and $U_2 = \{ 0 \}$, as desired.
\end{proof}

\begin{remark}
The first claim of Proposition \ref{propStructureSLFSousL10inv} is a consequence of the previous proposition. Indeed, $\mathrm{SLF}(\bar U_q)$ is indecomposable under the action of $W\bigl(\mathcal{S}_q(\Sigma_1)\bigr) = \mathbb{C}\langle W_{\! A}, W_{\! B} \rangle \subset \mathcal{L}_{1,0}^{\mathrm{inv}}(\bar U_q)$, and hence it is indecomposable under the action of the whole algebra $\mathcal{L}_{1,0}^{\mathrm{inv}}(\bar U_q)$.
\end{remark}

\subsection{Relationship with the skein representation}
\indent Consider a handlebody $H_g \subset \mathbb{R}^3,$ such that $\partial H_g =  \Sigma_g \subset \mathbb{R}^3$, for instance $H_g = \Sigma_{0,g+1} \times [0,1]$. Let $d$ be the Euclidean distance on $\mathbb{R}^3$, let  $\varepsilon > 0$ and define
\begin{equation}\label{decoupageToreSolide}
H_g^{\leq\varepsilon/2} = \bigl\{ x \in H_g \, \bigl| \, d(x, \partial H_g) \leq \frac{\varepsilon}{2} \bigr. \bigr\}, \:\:\:\:\:\:\: H_g^{\geq\varepsilon} = \bigl\{ x \in H_g \, \bigl| \, d(x, \partial H_g) \geq \varepsilon \bigr. \bigr\}.
\end{equation}
Take $\varepsilon$ sufficiently small, so that $H_g^{\leq\varepsilon/2}$ is diffeomorphic to $\Sigma_g \times [0,1]$ and $H_g^{\geq\varepsilon}$ is diffeomorphic to $H_g$. This dichotomy gives a representation $\rho$ of the skein algebra $\mathcal{S}_q(\Sigma_g) = \mathcal{S}_q(H_g^{\leq\varepsilon/2})$ on the skein module $\mathcal{S}_q(H_g) = \mathcal{S}_q(H_g^{\geq\varepsilon})$, defined by:
\begin{equation}\label{repSkeinHg}
\rho(L_1)(L_2) = \langle L_1 \cup L_2 \rangle, 
\end{equation}
where $\langle L \rangle$ is the value of $L$ in $\mathcal{S}_q(H_g)$. In practice, this just means that we put the link $L_1 \subset \Sigma_g \times [0,1]$ very close to $\partial H_g$ and the link $L_2 \subset H_g$ very close to the core of $H_g$.

\smallskip

\begin{figure}[!h]
\centering
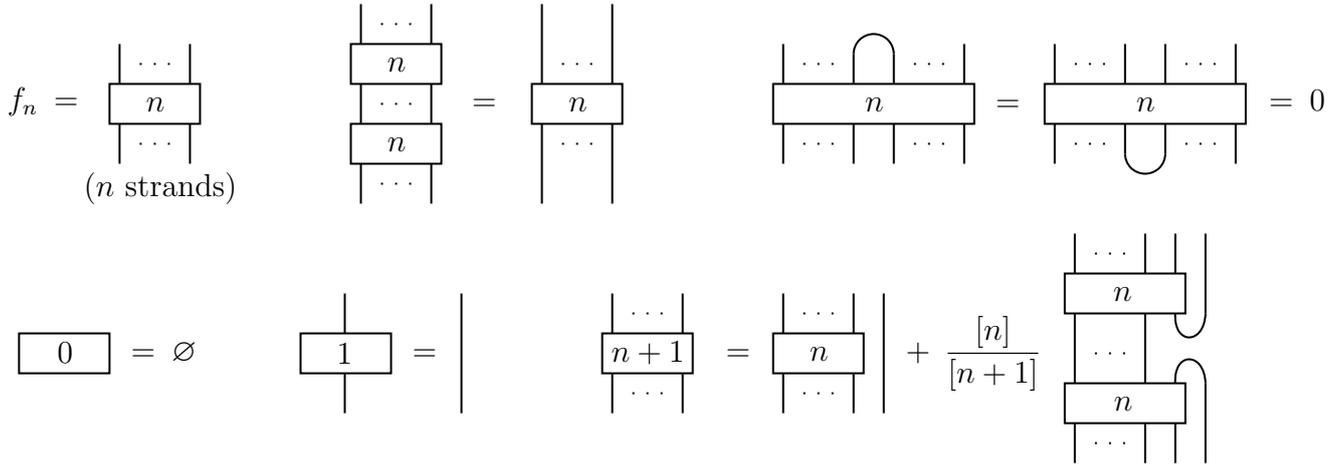
\caption{Jones-Wenzl idempotents.}
\label{JonesWenzl}
\end{figure}

Recall, for $0 \leq n \leq p-1$ (where $q^{2p}=1$), the $n$-th Jones-Wenzl idempotent $f_n$. This is an element of the Temperley-Lieb algebra on $n$ strands $\mathrm{TL}_{q,n} = \mathrm{End}_{\bar U_q}\!\bigl(\mathcal{X}^+(2)^{\otimes n}\bigr)$ (see \textit{e.g.} \cite{CFS} and the references therein)\footnote{Note that for $n \geq p$ the identification between the Temperley-Lieb algebra on $n$ strands and the centralizer of $\bar U_q$ on $\mathcal{X}^+(2)^{\otimes n}$ is not true: instead we have a {\em strict} embedding $\mathrm{TL}_{q,n} \hookrightarrow \mathrm{End}_{\bar U_q}\!\bigl(\mathcal{X}^+(2)^{\otimes n}\bigr)$.}. The properties of these elements $f_n$ are listed in Figure \ref{JonesWenzl}. Note that $\mathcal{X}^+(2)^{\otimes n} \cong V \oplus \mathcal{X}^+(n+1)$, where $V$ does not contain $\mathcal{X}^+(n+1)$ as a direct summand; then $f_n \in \mathrm{End}_{\bar U_q}\!\bigl(\mathcal{X}^+(2)^{\otimes n}\bigr)$ is the unique (up to scalar) morphism which factorizes through $\mathcal{X}^+(n+1)$: $\mathcal{X}^+(2)^{\otimes n} \to \mathcal{X}^+(n+1) \to \mathcal{X}^+(2)^{\otimes n}$.

\smallskip

\indent For each $n$, we have a closure map $\mathrm{cl} : \mathrm{TL}_{q,n} \to \mathcal{S}_q(H_1)$ (note that $H_1$ is a thickened annulus $ S^1 \times [0,1] \times [0,1]$):
\begin{center}
\begingroup%
  \makeatletter%
  \providecommand\color[2][]{%
    \errmessage{(Inkscape) Color is used for the text in Inkscape, but the package 'color.sty' is not loaded}%
    \renewcommand\color[2][]{}%
  }%
  \providecommand\transparent[1]{%
    \errmessage{(Inkscape) Transparency is used (non-zero) for the text in Inkscape, but the package 'transparent.sty' is not loaded}%
    \renewcommand\transparent[1]{}%
  }%
  \providecommand\rotatebox[2]{#2}%
  \newcommand*\fsize{\dimexpr\f@size pt\relax}%
  \newcommand*\lineheight[1]{\fontsize{\fsize}{#1\fsize}\selectfont}%
  \ifx\svgwidth\undefined%
    \setlength{\unitlength}{181.91044253bp}%
    \ifx\svgscale\undefined%
      \relax%
    \else%
      \setlength{\unitlength}{\unitlength * \real{\svgscale}}%
    \fi%
  \else%
    \setlength{\unitlength}{\svgwidth}%
  \fi%
  \global\let\svgwidth\undefined%
  \global\let\svgscale\undefined%
  \makeatother%
  \begin{picture}(1,0.38655259)%
    \lineheight{1}%
    \setlength\tabcolsep{0pt}%
    \put(0,0){\includegraphics[width=\unitlength,page=1]{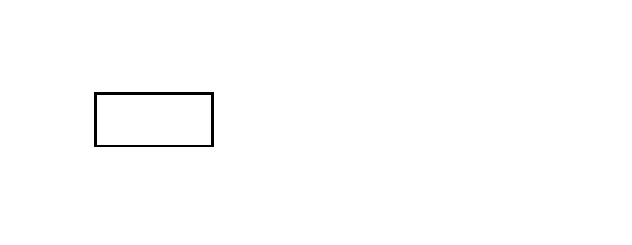}}%
    \put(0.22606178,0.1820538){\color[rgb]{0,0,0}\makebox(0,0)[lt]{\lineheight{1.25}\smash{\begin{tabular}[t]{l}$D$\end{tabular}}}}%
    \put(0,0){\includegraphics[width=\unitlength,page=2]{closure.pdf}}%
    \put(0.09814349,0.00647782){\color[rgb]{0,0,0}\makebox(0,0)[lt]{\lineheight{1.25}\smash{\begin{tabular}[t]{l}($n$ strands)\end{tabular}}}}%
    \put(0,0){\includegraphics[width=\unitlength,page=3]{closure.pdf}}%
    \put(-0.00243187,0.18526966){\color[rgb]{0,0,0}\makebox(0,0)[lt]{\lineheight{1.25}\smash{\begin{tabular}[t]{l}cl\end{tabular}}}}%
    \put(0,0){\includegraphics[width=\unitlength,page=4]{closure.pdf}}%
    \put(0.87446381,0.17740565){\color[rgb]{0,0,0}\makebox(0,0)[lt]{\lineheight{1.25}\smash{\begin{tabular}[t]{l}$D$\end{tabular}}}}%
    \put(0,0){\includegraphics[width=\unitlength,page=5]{closure.pdf}}%
    \put(0.47007172,0.18325509){\color[rgb]{0,0,0}\makebox(0,0)[lt]{\lineheight{1.25}\smash{\begin{tabular}[t]{l}$=$\end{tabular}}}}%
  \end{picture}%
\endgroup%

\end{center}

\smallskip

\indent Since $q$ is a $2p$-th root of unity, we can consider the reduced skein module $\mathcal{S}^{\mathrm{red}}_q(H_g)$, which is a quotient of $\mathcal{S}_q(H_g)$ by relations involving the Jones-Wenzl idempotents, see \textit{e.g.} \cite{costantino}. In particular, any diagram containing $f_{p-1}$ is null in the reduced skein module. Here we consider the case $g=1$, and a basis of $\mathcal{S}^{\mathrm{red}}_q(H_1)$ is given by the (classes of the) closures $\mathrm{cl}(f_n)$, with $0 \leq n \leq p-2$.

\smallskip

\indent Let us study the representation of $\mathcal{S}_q(\Sigma_1)$ on $\mathcal{S}_q^{\mathrm{red}}(H_1)$. It is known that $\mathcal{S}_q(\Sigma_1)$ is generated by the circles $[a]$ and $[b]$; by definition of the representation, their actions are given by
\begin{center}
\begingroup%
  \makeatletter%
  \providecommand\color[2][]{%
    \errmessage{(Inkscape) Color is used for the text in Inkscape, but the package 'color.sty' is not loaded}%
    \renewcommand\color[2][]{}%
  }%
  \providecommand\transparent[1]{%
    \errmessage{(Inkscape) Transparency is used (non-zero) for the text in Inkscape, but the package 'transparent.sty' is not loaded}%
    \renewcommand\transparent[1]{}%
  }%
  \providecommand\rotatebox[2]{#2}%
  \newcommand*\fsize{\dimexpr\f@size pt\relax}%
  \newcommand*\lineheight[1]{\fontsize{\fsize}{#1\fsize}\selectfont}%
  \ifx\svgwidth\undefined%
    \setlength{\unitlength}{454.06618026bp}%
    \ifx\svgscale\undefined%
      \relax%
    \else%
      \setlength{\unitlength}{\unitlength * \real{\svgscale}}%
    \fi%
  \else%
    \setlength{\unitlength}{\svgwidth}%
  \fi%
  \global\let\svgwidth\undefined%
  \global\let\svgscale\undefined%
  \makeatother%
  \begin{picture}(1,0.19697485)%
    \lineheight{1}%
    \setlength\tabcolsep{0pt}%
    \put(-0.00146572,0.09682656){\color[rgb]{0,0,0}\makebox(0,0)[lt]{\lineheight{1.25}\smash{\begin{tabular}[t]{l}$\rho(a)$\end{tabular}}}}%
    \put(0,0){\includegraphics[width=\unitlength,page=1]{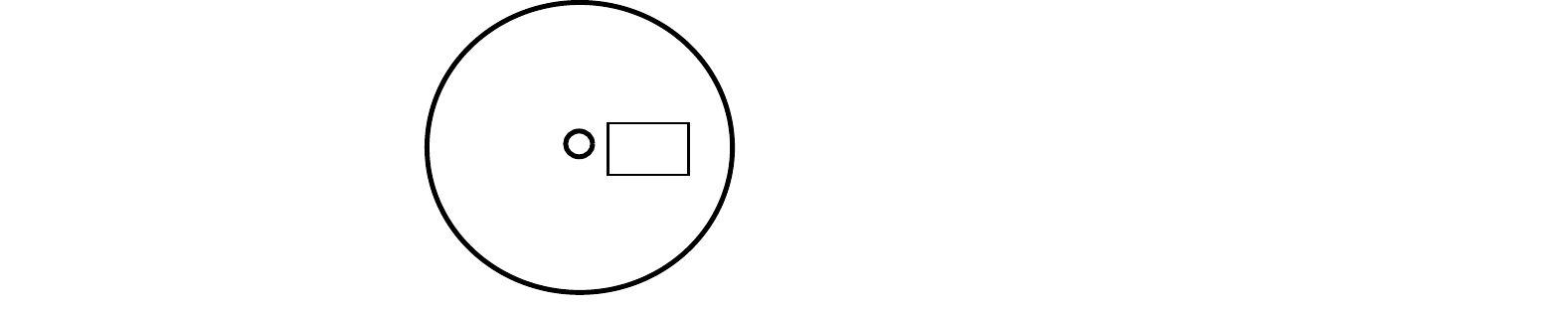}}%
    \put(0.40467072,0.09771989){\color[rgb]{0,0,0}\makebox(0,0)[lt]{\lineheight{1.25}\smash{\begin{tabular}[t]{l}$n$\end{tabular}}}}%
    \put(0,0){\includegraphics[width=\unitlength,page=2]{actionSkein.pdf}}%
    \put(0.23190911,0.09783952){\color[rgb]{0,0,0}\makebox(0,0)[lt]{\lineheight{1.25}\smash{\begin{tabular}[t]{l}$=$\end{tabular}}}}%
    \put(0,0){\includegraphics[width=\unitlength,page=3]{actionSkein.pdf}}%
    \put(0.93588671,0.08900477){\color[rgb]{0,0,0}\makebox(0,0)[lt]{\lineheight{1.25}\smash{\begin{tabular}[t]{l}$n$\end{tabular}}}}%
    \put(0,0){\includegraphics[width=\unitlength,page=4]{actionSkein.pdf}}%
    \put(0.16926048,0.0969217){\color[rgb]{0,0,0}\makebox(0,0)[lt]{\lineheight{1.25}\smash{\begin{tabular}[t]{l}$n$\end{tabular}}}}%
    \put(0,0){\includegraphics[width=\unitlength,page=5]{actionSkein.pdf}}%
    \put(0.52977926,0.09240465){\color[rgb]{0,0,0}\makebox(0,0)[lt]{\lineheight{1.25}\smash{\begin{tabular}[t]{l}$\rho(b)$\end{tabular}}}}%
    \put(0.76315412,0.09341761){\color[rgb]{0,0,0}\makebox(0,0)[lt]{\lineheight{1.25}\smash{\begin{tabular}[t]{l}$=$\end{tabular}}}}%
    \put(0,0){\includegraphics[width=\unitlength,page=6]{actionSkein.pdf}}%
    \put(0.70050548,0.09249978){\color[rgb]{0,0,0}\makebox(0,0)[lt]{\lineheight{1.25}\smash{\begin{tabular}[t]{l}$n$\end{tabular}}}}%
    \put(0,0){\includegraphics[width=\unitlength,page=7]{actionSkein.pdf}}%
  \end{picture}%
\endgroup%

\end{center}
where we denote $\rho(a), \rho(b)$ instead of $\rho([a]), \rho([b])$.

\begin{lemma}\label{lemmeActionSkein}
In $\mathcal{S}_q^{\mathrm{red}}(H_1)$ it holds
\begin{align*}
\rho(a)\bigl( \mathrm{cl}(f_n) \bigr) &= -(q^{n+1} + q^{-(n+1)})\mathrm{cl}(f_n),\\
\rho(b)\bigl( \mathrm{cl}(f_0) \bigr) = \mathrm{cl}(f_1), \:\:\:\:\: \rho(b)\bigl( \mathrm{cl}(f_n) \bigr) &= \mathrm{cl}(f_{n-1}) + \mathrm{cl}(f_{n+1}), \:\:\:\:\: \rho(b)\bigl( \mathrm{cl}(f_{p-2}) \bigr) = \mathrm{cl}(f_{p-3}).\\
\end{align*}
\end{lemma}
\begin{proof}
Observe that, for $n \geq 2$,
\begin{center}
\begingroup%
  \makeatletter%
  \providecommand\color[2][]{%
    \errmessage{(Inkscape) Color is used for the text in Inkscape, but the package 'color.sty' is not loaded}%
    \renewcommand\color[2][]{}%
  }%
  \providecommand\transparent[1]{%
    \errmessage{(Inkscape) Transparency is used (non-zero) for the text in Inkscape, but the package 'transparent.sty' is not loaded}%
    \renewcommand\transparent[1]{}%
  }%
  \providecommand\rotatebox[2]{#2}%
  \newcommand*\fsize{\dimexpr\f@size pt\relax}%
  \newcommand*\lineheight[1]{\fontsize{\fsize}{#1\fsize}\selectfont}%
  \ifx\svgwidth\undefined%
    \setlength{\unitlength}{725.92621739bp}%
    \ifx\svgscale\undefined%
      \relax%
    \else%
      \setlength{\unitlength}{\unitlength * \real{\svgscale}}%
    \fi%
  \else%
    \setlength{\unitlength}{\svgwidth}%
  \fi%
  \global\let\svgwidth\undefined%
  \global\let\svgscale\undefined%
  \makeatother%
  \begin{picture}(1,0.14117926)%
    \lineheight{1}%
    \setlength\tabcolsep{0pt}%
    \put(0,0){\includegraphics[width=\unitlength,page=1]{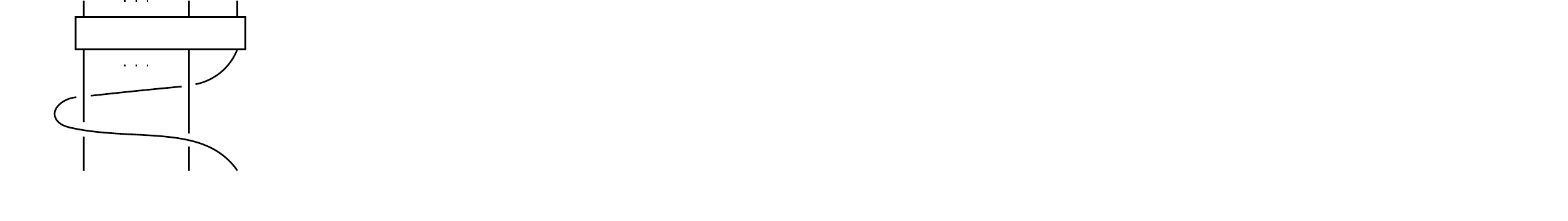}}%
    \put(0.0969382,0.11648139){\color[rgb]{0,0,0}\makebox(0,0)[lt]{\lineheight{1.25}\smash{\begin{tabular}[t]{l}$n$\end{tabular}}}}%
    \put(-0.00114198,0.07575292){\color[rgb]{0,0,0}\makebox(0,0)[lt]{\lineheight{1.25}\smash{\begin{tabular}[t]{l}cl\end{tabular}}}}%
    \put(0,0){\includegraphics[width=\unitlength,page=2]{cl_JW_braiding.pdf}}%
    \put(0.36556053,0.11648139){\color[rgb]{0,0,0}\makebox(0,0)[lt]{\lineheight{1.25}\smash{\begin{tabular}[t]{l}$n$\end{tabular}}}}%
    \put(0.20826078,0.07415036){\color[rgb]{0,0,0}\makebox(0,0)[lt]{\lineheight{1.25}\smash{\begin{tabular}[t]{l}$= \:\: q^{-1/2}\mathrm{cl}$\end{tabular}}}}%
    \put(0,0){\includegraphics[width=\unitlength,page=3]{cl_JW_braiding.pdf}}%
    \put(0.62329778,0.11611241){\color[rgb]{0,0,0}\makebox(0,0)[lt]{\lineheight{1.25}\smash{\begin{tabular}[t]{l}$n$\end{tabular}}}}%
    \put(0.47632965,0.07378135){\color[rgb]{0,0,0}\makebox(0,0)[lt]{\lineheight{1.25}\smash{\begin{tabular}[t]{l}$= \:\: q^{-1}\mathrm{cl}$\end{tabular}}}}%
    \put(0,0){\includegraphics[width=\unitlength,page=4]{cl_JW_braiding.pdf}}%
    \put(0.47932632,0.00488122){\color[rgb]{0,0,0}\makebox(0,0)[lt]{\lineheight{1.25}\smash{\begin{tabular}[t]{l}$= \:\: \ldots \:\: =  \: q^{-n+1}\mathrm{cl}(f_n)$\end{tabular}}}}%
  \end{picture}%
\endgroup%

\end{center}
We used the fact that the composition of a cup or a cap with the Jones-Wenzl idempotents is $0$ and the cyclicity of $\mathrm{cl}$.
Now, assume by induction that $\rho(a)\bigl( \mathrm{cl}(f_n) \bigr) = \lambda_n \mathrm{cl}(f_n)$ for a family of scalars $\lambda_n$. It is easily checked that $\lambda_0 = -(q + q^{-1}), \: \lambda_1 = - (q^2 + q^{-2})$. Applying the Kauffman skein relation twice in $\mathrm{TL}_{q,n}$, we obtain
\begin{center}
\begingroup%
  \makeatletter%
  \providecommand\color[2][]{%
    \errmessage{(Inkscape) Color is used for the text in Inkscape, but the package 'color.sty' is not loaded}%
    \renewcommand\color[2][]{}%
  }%
  \providecommand\transparent[1]{%
    \errmessage{(Inkscape) Transparency is used (non-zero) for the text in Inkscape, but the package 'transparent.sty' is not loaded}%
    \renewcommand\transparent[1]{}%
  }%
  \providecommand\rotatebox[2]{#2}%
  \newcommand*\fsize{\dimexpr\f@size pt\relax}%
  \newcommand*\lineheight[1]{\fontsize{\fsize}{#1\fsize}\selectfont}%
  \ifx\svgwidth\undefined%
    \setlength{\unitlength}{514.61310072bp}%
    \ifx\svgscale\undefined%
      \relax%
    \else%
      \setlength{\unitlength}{\unitlength * \real{\svgscale}}%
    \fi%
  \else%
    \setlength{\unitlength}{\svgwidth}%
  \fi%
  \global\let\svgwidth\undefined%
  \global\let\svgscale\undefined%
  \makeatother%
  \begin{picture}(1,0.18480462)%
    \lineheight{1}%
    \setlength\tabcolsep{0pt}%
    \put(0,0){\includegraphics[width=\unitlength,page=1]{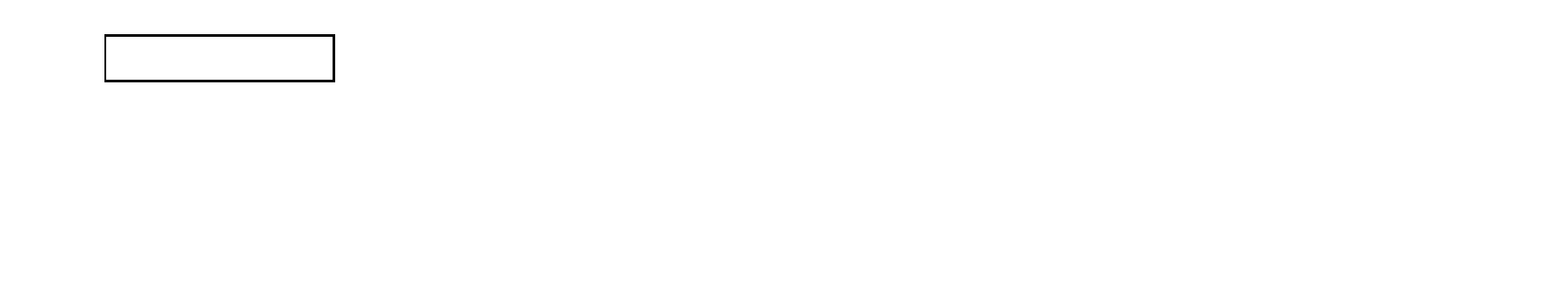}}%
    \put(0.11571332,0.14046989){\color[rgb]{0,0,0}\makebox(0,0)[lt]{\lineheight{1.25}\smash{\begin{tabular}[t]{l}$n+1$\end{tabular}}}}%
    \put(0,0){\includegraphics[width=\unitlength,page=2]{JW_entoure.pdf}}%
    \put(0.28485027,0.10539688){\color[rgb]{0,0,0}\makebox(0,0)[lt]{\lineheight{1.25}\smash{\begin{tabular}[t]{l}$= \:\:\: q$\end{tabular}}}}%
    \put(0.62089611,0.10370184){\color[rgb]{0,0,0}\makebox(0,0)[lt]{\lineheight{1.25}\smash{\begin{tabular}[t]{l}$+ \:\: (1-q^{-2})$\end{tabular}}}}%
    \put(0,0){\includegraphics[width=\unitlength,page=3]{JW_entoure.pdf}}%
    \put(0.46041194,0.14083792){\color[rgb]{0,0,0}\makebox(0,0)[lt]{\lineheight{1.25}\smash{\begin{tabular}[t]{l}$n+1$\end{tabular}}}}%
    \put(0,0){\includegraphics[width=\unitlength,page=4]{JW_entoure.pdf}}%
    \put(0.86193461,0.14959714){\color[rgb]{0,0,0}\makebox(0,0)[lt]{\lineheight{1.25}\smash{\begin{tabular}[t]{l}$n+1$\end{tabular}}}}%
    \put(0,0){\includegraphics[width=\unitlength,page=5]{JW_entoure.pdf}}%
    \put(-0.00139823,0.10293032){\color[rgb]{0,0,0}\makebox(0,0)[lt]{\lineheight{1.25}\smash{\begin{tabular}[t]{l}cl\end{tabular}}}}%
    \put(0,0){\includegraphics[width=\unitlength,page=6]{JW_entoure.pdf}}%
    \put(0.34109208,0.10293031){\color[rgb]{0,0,0}\makebox(0,0)[lt]{\lineheight{1.25}\smash{\begin{tabular}[t]{l}cl\end{tabular}}}}%
    \put(0,0){\includegraphics[width=\unitlength,page=7]{JW_entoure.pdf}}%
    \put(0.74187852,0.1029303){\color[rgb]{0,0,0}\makebox(0,0)[lt]{\lineheight{1.25}\smash{\begin{tabular}[t]{l}cl\end{tabular}}}}%
    \put(0,0){\includegraphics[width=\unitlength,page=8]{JW_entoure.pdf}}%
    \put(0.28492424,0.0036744){\color[rgb]{0,0,0}\makebox(0,0)[lt]{\lineheight{1.25}\smash{\begin{tabular}[t]{l}$= \bigl( q\lambda_n + (1-q^{-2})q^{-n}\bigr)\mathrm{cl}(f_{n+1})$\end{tabular}}}}%
  \end{picture}%
\endgroup%

\end{center}
For the second equality we used the recurrence formula for the Jones-Wenzl idempotents together with the induction hypothesis. Hence $\lambda_{n+1} = q\lambda_n + (1-q^{-2})q^{-n}$, and it follows that $\lambda_n = -(q^{n+1} - q^{-(n+1)})$, as desired. To compute the action of $\rho(b)$, note first that thanks to the reccurence formula, we have
\begin{center}
\begingroup%
  \makeatletter%
  \providecommand\color[2][]{%
    \errmessage{(Inkscape) Color is used for the text in Inkscape, but the package 'color.sty' is not loaded}%
    \renewcommand\color[2][]{}%
  }%
  \providecommand\transparent[1]{%
    \errmessage{(Inkscape) Transparency is used (non-zero) for the text in Inkscape, but the package 'transparent.sty' is not loaded}%
    \renewcommand\transparent[1]{}%
  }%
  \providecommand\rotatebox[2]{#2}%
  \newcommand*\fsize{\dimexpr\f@size pt\relax}%
  \newcommand*\lineheight[1]{\fontsize{\fsize}{#1\fsize}\selectfont}%
  \ifx\svgwidth\undefined%
    \setlength{\unitlength}{285.19048096bp}%
    \ifx\svgscale\undefined%
      \relax%
    \else%
      \setlength{\unitlength}{\unitlength * \real{\svgscale}}%
    \fi%
  \else%
    \setlength{\unitlength}{\svgwidth}%
  \fi%
  \global\let\svgwidth\undefined%
  \global\let\svgscale\undefined%
  \makeatother%
  \begin{picture}(1,0.14001573)%
    \lineheight{1}%
    \setlength\tabcolsep{0pt}%
    \put(0,0){\includegraphics[width=\unitlength,page=1]{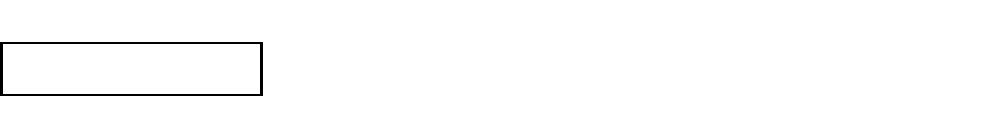}}%
    \put(0.11166472,0.05894996){\color[rgb]{0,0,0}\makebox(0,0)[lt]{\lineheight{1.25}\smash{\begin{tabular}[t]{l}$n$\end{tabular}}}}%
    \put(0,0){\includegraphics[width=\unitlength,page=2]{trR_JW.pdf}}%
    \put(0.32921563,0.06134621){\color[rgb]{0,0,0}\makebox(0,0)[lt]{\lineheight{1.25}\smash{\begin{tabular}[t]{l}$\displaystyle= \:\: -\frac{[n+1]}{[n]} f_{n-1}$\end{tabular}}}}%
  \end{picture}%
\endgroup%

\end{center}
Hence, using again the recurrence formula of the $f_n$'s and the cyclicity of $\mathrm{cl}$, we get
\begin{center}
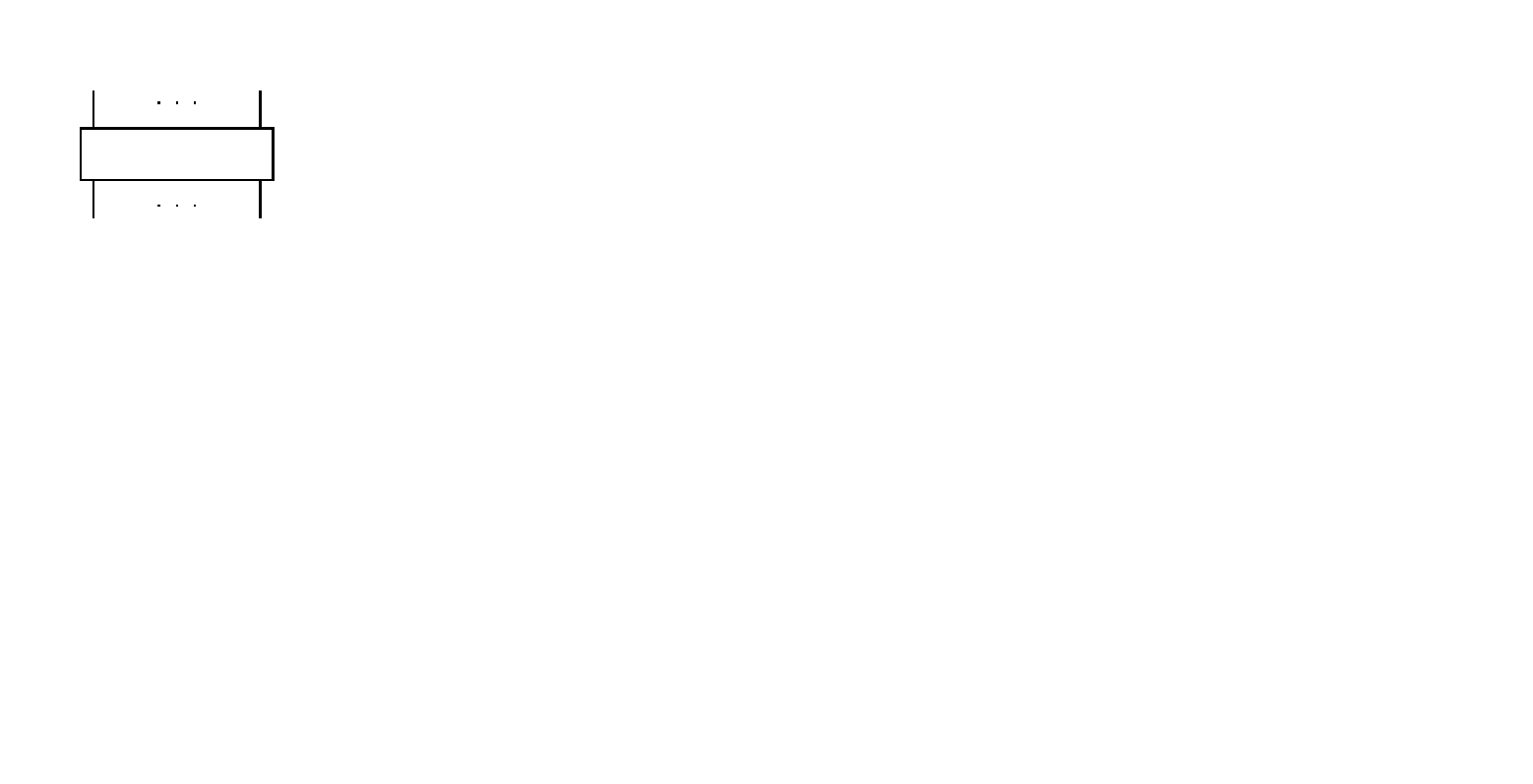
\end{center}
The case $n=0$ is obvious and the case $n=p-2$ follows from the previous equality and the fact that $f_{p-1} = 0$ in $\mathcal{S}_q^{\mathrm{red}}(H_1)$.
\end{proof}

Let $\overline{\mathcal{U}} = J_2/J_1 = \mathrm{vect}\bigl( \overline{\chi}^+_s \bigr)_{1 \leq s \leq p-1}$, where $\overline{\chi}^+_s$ is the class of $\chi^+_s$ modulo $\mathcal{P}$ (see Proposition \ref{structureRepSqS1}).

\begin{proposition}\label{propLienSkeinEtL10}
The $\mathcal{S}_q(\Sigma_1)$-modules $\mathcal{S}_q^{\mathrm{red}}(H_1)$ and $\overline{\mathcal{U}}$ are isomorphic. The isomorphism is given by
\[ \fonc{F}{\mathcal{S}_q^{\mathrm{red}}(H_1)}{\overline{\mathcal{U}}}{\mathrm{cl}(f_n)}{v_A^{-1} \triangleright \overline{\chi}^+_{n+1} = v_{n+1}^{-1} \overline{\chi}^+_{n+1}} \]
where $v_s$ is defined in \eqref{symetrieV}.
\end{proposition}
\begin{proof}
Recall that
\begin{equation*}
\begin{split}
&W(a) \triangleright \overline{\chi}^+_s = -(q^s + q^{-s})\overline{\chi}^+_s, \\
W(b) \triangleright \overline{\chi}^+_1 = \frac{v_1}{v_2} \overline{\chi}^+_2, \:\:\:\:\:\:\: &W(b) \triangleright \overline{\chi}^+_s = \frac{v_s}{v_{s-1}} \overline{\chi}^+_{s-1} +  \frac{v_s}{v_{s+1}}\overline{\chi}^+_{s+1}, \:\:\:\:\:\:\: W(b) \triangleright \overline{\chi}^+_{p-1} = \frac{v_{p-1}}{v_{p-2}} \overline{\chi}^+_{p-2}.
\end{split}
\end{equation*}
The result follows by comparison with the formulas of Lemma \ref{lemmeActionSkein}.
\end{proof}

\indent The representation $\mathcal{S}_q^{\mathrm{red}}(H_1)$ can be described in terms of the Reshetikhin-Turaev topological quantum field theory (RT TQFT for short). Recall the cobordism category $\mathcal{C}$ for this TQFT: the objects are closed oriented surfaces and the morphisms $\mathrm{Hom}_{\mathcal{C}}(S_1, S_2)$ are pairs $(M,L)$, where $M$ is a compact oriented $3$-manifold endowed with a $p_1$-structure such that $\partial M = (-S_1) \sqcup S_2$ and $L \subset M$ is a framed link up to isotopy. For more informations, see \cite{RT2, BHMV}, the lecture notes \cite{costantino} and also \cite{BW}, where the representation of $\mathcal{S}_q(\Sigma_g)$ given by the RT TQFT is shown to be irreducible. The $3$-manifolds $H_1^{\leq\varepsilon/2} \cong \Sigma_1 \times [0,1], H_1^{\geq\varepsilon} \cong H_1$ of \eqref{decoupageToreSolide} with framed links inside them are cobordisms:
\[ (\Sigma_1 \times [0,1], L_1) \in \mathrm{Hom}_{\mathcal{C}}(\Sigma_1, \Sigma_1), \:\:\:\:\:\:\:  (H_1, L_2) \in \mathrm{Hom}_{\mathcal{C}}(\varnothing, \Sigma_1). \]
The functor $Z : \mathcal{C} \to \mathrm{Vect}_{\mathbb{C}}$ of the RT TQFT (which depends on the primitive root of unity $q$), gives  linear maps
\[ Z(\Sigma_1 \times [0,1], L_1) : Z(\Sigma_1) \to Z(\Sigma_1), \:\:\:\:\:\:\: Z(H_1, L_2) : Z(\varnothing) = \mathbb{C} \to Z(\Sigma_1). \]
Hence $Z(H_1, L_2)$ is just a choice of an element of $Z(\Sigma_1)$. As recalled in \cite[Lemma 5]{BW}, every element of $Z(\Sigma_1)$ can be written $Z(H_1, L_2)$ for some link $L_2 \subset H_1$. Moreover, the maps $Z(M, \cdot)$ satisfy the Kauffman skein relation:
\[ Z(M, L_+) = q^{1/2}Z(M,L_{||}) + q^{-1/2}Z(M, L_=) \]
where $L_+, L_{||}, L_= \subset M$ are identical except in a little ball in which they look like in \eqref{relationKauffman}. It follows that  $Z(\Sigma_1 \times [0,1], \cdot) : \mathcal{S}_q(\Sigma_1) \to \mathrm{End}_{\mathbb{C}}\bigl( Z(\Sigma_1) \bigr)$ furnishes a representation of $\mathcal{S}_q(\Sigma_1)$ on $Z(\Sigma_1)$:
\begin{equation}\label{repTQFT}
L_1 \cdot Z(H_1, L_2) = Z\bigl( (\Sigma_1 \times [0,1], L_1) \circ (H_1, L_2) \bigr) \in \mathrm{Hom}_{\mathcal{C}}(Z(\varnothing), Z(\Sigma_1)) = Z(\Sigma_1).  
\end{equation}
The map $Z_{H_1} : \mathcal{S}_q(H_1) \to Z(\Sigma_1)$ defined by $Z_{H_1}(L) = Z(H_1, L)$ is surjective, and more precisely it gives rise to an isomorphism $Z_{H_1}^{\mathrm{red}} : \mathcal{S}^{\mathrm{red}}_q(H_1) \to Z(\Sigma_1)$. Moroever, $Z_{H_1}^{\mathrm{red}}$ is an isomorphism of $\mathcal{S}_q(\Sigma_1)$-modules between $\mathcal{S}^{\mathrm{red}}_q(H_1)$ endowed with the action \eqref{repSkeinHg} and $Z(\Sigma_1)$ endowed with the action  \eqref{repTQFT}.

\smallskip

\indent We have seen that, for $\Sigma_1$, the representation of Theorem \ref{theoRepSkein} ``contains'' the natural skein representation \eqref{repSkeinHg} on $\mathcal{S}_q^{\mathrm{red}}(H_1)$ (or equivalently \eqref{repTQFT}), in the sense that it is the composition factor $J_2/J_1$. By \cite[Theorem 7]{BW}, we know that the representation \eqref{repTQFT} is irreducible in any genus. Hence we conjecture the following.

\begin{conjecture}\label{conjectureSkein}
The representation \eqref{repSkeinHg} of $\mathcal{S}_q(\Sigma_g)$ on $\mathcal{S}_q^{\mathrm{red}}(H_g)$ is a composition factor $J_{i+1}/J_i$ of the representation of Theorem \ref{theoRepSkein} (namely the representation induced by the Wilson loop map $W : \mathcal{S}_q(\Sigma_g) \to \mathcal{L}_{g,0}(\bar U_q)$ and the representation of $\mathcal{L}_{g,0}(\bar U_q)$ on $(\bar U_q^*)^{\otimes g}$).
\end{conjecture}

\end{document}